\documentclass[10pt,final,mytheorems]{book}
\usepackage[latin1]{inputenc}
\pdfoutput=1

\usepackage{bm}
\usepackage{epsfig}
\usepackage{makeidx}

\usepackage{theorem}
\usepackage{amsmath}
\usepackage{amssymb}
\usepackage{amscd}
\usepackage[mathscr]{euscript}
\usepackage{rotating}
\usepackage[all]{xy}
\usepackage{chngcntr}
\usepackage{hyperref}

\theorembodyfont{\sl}

\newtheorem{theorem}[subsection]{Theorem}
\newtheorem{lemma}[subsection]{Lemma}
\newtheorem{proposition}[subsection]{Proposition}
\newtheorem{corollary}[subsection]{Corollary}

\newtheorem{subtheoreme}[subsection]{Theorem}
\newtheorem{sublemme}[subsection]{Lemma}
\newtheorem{subproposition}[subsection]{Proposition}
\newtheorem{subcorollaire}[subsection]{Corollary}

\newtheorem{thA}{Theorem A}

\newtheorem{thB}{Theorem B}

\theorembodyfont{\rmfamily}

\newtheorem{notation}[subsection]{Notation}
\newtheorem{remarks}[subsection]{Remarks}
\newtheorem{remark}[subsection]{Remark}
\newtheorem{definition}[subsection]{Definition}
\newtheorem{example}[subsection]{Example}

\newtheorem{subremarques}[subsection]{Remarks}
\newtheorem{subremarque}[subsection]{Remark}
\newtheorem{subdefinition}[subsection]{Definition}

\newenvironment{proof}{\vspace{0.3cm}\emph{Proof.}}
{\hfill $\square$ \vspace{0.3cm}}
\newenvironment{prooft}{\vspace{0.3cm}\emph{Proof of the theorem.}}
{\hfill $\square$ \vspace{0.3cm}}
\newenvironment{proofp}{\vspace{0.3cm}\emph{Proof of the proposition.}}
{\hfill $\square$ \vspace{0.3cm}}
\newenvironment{proofpn}[1]{\vspace{0.3cm}\emph{Proof of
proposition #1.}} {\hfill $\square$ \vspace{0.3cm}}

\newenvironment{bulletlist}
  {\itemize\let\origitem\item
   \renewcommand{\item}[1][default]
   {\origitem[\csname seamus##1\endcsname]}}
  {\enditemize}

\newcommand\C{\mathbb{C}}
\newcommand\Nat{\mathbb{N}}
\newcommand\R{\mathbb{R}}
\newcommand\Z{\mathbb{Z}}
\newcommand\Q{\mathbb {Q}}
\newcommand\Fi{\mathbb{F}}
\newcommand\Af{\mathbb{A}_f}
\newcommand\Of{\mathcal{O}}
\newcommand\Ade{\mathbb{A}}

\DeclareMathOperator\Hom{\mathrm{Hom}}

\newcommand\SD{\mathbb{S}}
\newcommand\T{{\bf T}}
\newcommand\G{{\bf G}}
\renewcommand\H{{\bf H}}
\newcommand\GU{{\bf GU}}
\newcommand\U{{\bf U}}

\newcommand\SU{{\bf SU}}
\newcommand\GSp{{\bf GSp}}

\newcommand\GSO{{\bf GSO}}
\newcommand\GO{{\bf GO}}

\newcommand\GL{{\bf GL}}
\newcommand\Pa{{\bf P}}
\newcommand\QP{{\bf Q}}
\newcommand\RP{{\bf R}}
\newcommand\B{{\bf B}}
\newcommand\Gr{\mathbb{G}}
\newcommand\N{{\bf N}}
\newcommand\M{{\bf M}}
\newcommand\Le{{\bf L}}
\newcommand\Se{{\bf S}}
\newcommand\A{{\bf A}}
\newcommand\Ar{\mathrm{A}}
\newcommand\K{\mathrm {K}}
\newcommand\Hr{\mathrm{H}}
\newcommand\PGU{{\bf PGU}}
\newcommand\PGL{{\bf PGL}}
\newcommand\PGSp{{\bf PGSp}}
\newcommand\PGSO{{\bf PGSO}}

\newcommand\Pro{{\bf P}}

\newcommand\F{{\mathcal F}}
\newcommand\Gf{{\mathcal G}}

\newcommand\Dgoth{{\mathfrak D}}

\newcommand\ggoth{{\mathfrak g}}
\newcommand\hgoth{{\mathfrak h}}
\newcommand\tgoth{{\mathfrak t}}
\newcommand\Kgoth{{\mathfrak K}}
\newcommand\agoth{{\mathfrak a}}
\newcommand\Sgoth{{\mathfrak S}}

\DeclareMathOperator{\Ad}{Ad}
\DeclareMathOperator{\Aut}{Aut}
\DeclareMathOperator{\Cent}{Cent}
\DeclareMathOperator{\Gal}{Gal}
\DeclareMathOperator{\Ho}{H}
\DeclareMathOperator{\Int}{Int}
\DeclareMathOperator{\Ker}{Ker}
\DeclareMathOperator{\Lie}{Lie}
\DeclareMathOperator{\Nor}{Nor}
\DeclareMathOperator{\Ob}{Ob}
\DeclareMathOperator{\Out}{Out}
\DeclareMathOperator{\Tr}{Tr}
\DeclareMathOperator{\vol}{vol}
\newcommand\Bo{\mathcal{B}}
\newcommand\Ch{\mathcal{C}}
\newcommand\Diag{{\bf D}}

\newcommand\Ell{\mathcal{E}}

\newcommand\Ens{\mathcal{E}}
\newcommand\Gcal{\mathcal{G}}
\newcommand\Hecke{{\mathcal H}}
\newcommand\Hcal{{\mathcal H}}
\newcommand\Ka{{\mathfrak K}}
\newcommand\Kens{{\mathcal K}}
\newcommand\Levi{{\mathcal L}}
\newcommand\Mcal{{\mathcal M}}

\newcommand\Ncal{{\mathcal N}}
\newcommand\Norme{{\mathcal N}}
\newcommand\nugras{{\bm\nu}}

\DeclareMathOperator{\Spec}{Spec}
\newcommand{\quash}[1]{}

\newcommand\as{{\underline{a}}}
\newcommand\es{{\underline{e}}}

\newcommand\sgras{{\bm s}}
\newcommand\Mod{{\mathcal M}}
\newcommand\Par{{\mathcal P}}

\newcommand\Rens{{\mathcal R}}
\newcommand\sous{\setminus}
\newcommand\til{\widetilde}
\newcommand\X{{\mathcal X}}
\newcommand\Y{{\mathcal Y}}
\newcommand\fl{\longrightarrow}
\newcommand\fle{\longmapsto}
\newcommand\iso{\stackrel {\sim} {\fl}}

\newcommand\ungras{1\mkern -5mu\mathrm{l}}
\newcommand\ddotsinv{\begin{turn}{45}\large\ldots\end{turn}}

\newcommand{\LG}{{}^L G} 
\newcommand{\LR}{{}^L R}
\newcommand{\LT}{{}^L T} 
\newcommand{\LH}{{}^L H}
\newcommand{\Ltheta}{{}^L \theta}

\DeclareMathOperator{\inv}{inv}

\DeclareMathOperator{\Res}{Res}
\DeclareMathOperator{\cok}{cok} 
\DeclareMathOperator{\un}{un} 
\DeclareMathOperator{\cont}{cont} 
\DeclareMathOperator{\sgn}{sgn}

\makeindex

\begin{document}

\frontmatter

\title{On the cohomology of certain non-compact Shimura varieties}

\author{Sophie Morel \\ with an appendix by Robert Kottwitz}

% \makehalftitle
\maketitle

\tableofcontents

%\begin{thepreface}

The goal of this text is to calculate the trace of a Hecke correspondence
composed with a (big enough) power of the Frobenius automorphism at a good
place on the intersection cohomology of the Baily-Borel compactification
of certain Shimura varieties, and then to stabilize the result for the
Shimura varieties associated to unitary groups over $\Q$.

The main result is theorem \ref{th:stab_FT_IC_GL_n}. It expresses the
above trace in terms of the twisted trace formula on products of general
linear groups, for well-chosen test functions.

Here are two applications of this result. The first (corollary
\ref{cor:calcul_L_GL_n}) is about the calculation of the $L$-function of
the intersection complex of the Baily-Borel compactification :

\begin{thA} Let $E$ be a quadratic imaginary extension of $\Q$,
$\G=\GU(p,q)$ one of the unitary groups defined by using $E$ (cf
\ref{groupes1}), $\K$ a neat open compact subgroup of $\G(\Af)$,
$M^{\K}(\G,\X)$ the associated Shimura variety (cf \ref{groupes1} and
\ref{points_fixes1}) and $V$ an irreducible algebraic representation of $\G$.
Denote by $IC^{\K}V$ the intersection complex of the Baily-Borel
compactification of $M^{\K}(\G,\X)$ with coefficients in $V$.
Let $\Ell_G$ be the set of elliptic endoscopic groups $\G(\U^*(n_1)\times
\U^*(n_2))$ of $\G$, where $n_1,n_2\in\Nat$ are such that $n_1+n_2=p+q$
and that
$n_2$ is even. For every $\H\in\Ell_G$, let $\Pi_H$ be the set of equivalence
classes of automorphic representations of $\H(\Ade_E)$. 

Assume that $\K$ is small enough.
Then there exist, for every $\H\in\Ell_G$, an explicit finite set $R_H$ of
algebraic representations of ${}^L\H_E$ and a family of complex numbers
$(c_H(\pi_H,r_H))_{\pi_H\in\Pi_H,r_H\in R_H}$, almost all zero,
such that, for every finite place $\wp$ of $E$ above a prime number where $\K$
is hyperspecial,
\[\log L_\wp(s,IC^{\K}V)=\sum_{\H\in\Ell_G}
\sum_{\pi_H\in\Pi_H}\sum_{r_H\in R_H}
c_H(\pi_H,r_H)\log L(s-\frac{d}{2},\pi_{H,\wp},r_H),\]
where $d=pq$ is the dimension of $M^{\K}(\G,\X)$.

\end{thA} 

See the statement of corollary \ref{cor:calcul_L_GL_n} for more details.
The second application is corollary \ref{cor:GL_n_E}. We give a simplified
statement of this corollary and refer to \ref{GL_n_applications4} for
the definitions :

\begin{thB} Let $n$ be a positive integer that is not dividible by $4$ and
$E$ an imaginary quadratic extension of $\Q$. Denote by $\theta$ the
automorphism $g\fle{}^t\overline{g}^{-1}$ of $R_{E/\Q}\GL_{n,E}$.
Let $\pi$ be a $\theta$-stable cuspidal automorphic representation of
$\GL_n(\Ade_E)$ that is regular algebraic.
Let $S$ be the union of the set of prime numbers that ramify in $E$ and of the
set of prime numbers under finite places of $E$ where
$\pi$ is ramified.
Then there exists a number field $K$, a positive integer $N$ and, for
every finite place $\lambda$ of $K$, a continuous
finite-dimensional representation $\sigma_\lambda$ of
$\Gal(\overline{\Q}/E)$ with coefficients in $K_\lambda$, such that :
\begin{itemize}
\item[(i)] The representation $\sigma_\lambda$ is unramified outside of
$S\cup\{\ell\}$, where $\ell$ is the prime number under $\lambda$, and
pure of weight $1-n$.
\item[(ii)] For every place $\wp$ of $E$ above a prime number $p\not\in S$,
for every finite place $\lambda\not|p$ of $K$,
\[\log L_\wp(s,\sigma_\lambda)=N\log L(s+\frac{n-1}{2},\pi_{\wp}).\]

\end{itemize}

In particular, $\pi$ satisfies the Ramanujan-Petersson conjecture at every
finite unramified place.

\end{thB}

There is also a result for $n$ dividible by $4$, but it is weaker and
its statement is longer. See also chapter \ref{applications} for applications
of the stabilized fixed point formula (corollary \ref{cor:stab_FT_IC})
that do not use base change to $\GL_n$.

The method used in this text is the one developed by Ihara, Langlands and
Kottwitz : comparison of the fixed point formula of Grothendieck-Lefschetz
and of the trace formula of Arthur-Selberg. In the case where the Shimura
variety is compact and where the group has no endoscopy, this method is
explained in the article \cite{K-LAR} of Kottwitz. Using base change
from unitary groups to $\GL_n$, Clozel deduced from this a version of
corollary B with supplementary conditions on the automorphic representation
at a finite place (cf \cite{Cl-RGRA} and the article \cite{CL} of Clozel and
Labesse).
The case of a compact Shimura variety (more generally, of the cohomology
with compact support) and of a group that might have non-trivial endoscopy
is treated by Kottwitz in the articles \cite{K-PSSV} and \cite{K-SVLR}, modulo
the fundamental lemma. For unitary groups, the fundamental lemma
(and the twisted
version that is used in the stabilization of the fixed point formula)
is now known thanks to the work of Laumon-Ngo (\cite{LN}), Hales
(\cite{Ha}) and Waldspurger (\cite{Wa1}, \cite{Wa2}, \cite{Wa3}); note that
the fundamental lemma is even known in general thanks to the recent article
of Ngo (\cite{Ng}).

The case of $\GL_2$ over $\Q$ (ie of modular curves) has been treated in the
book \cite{DK}, and the case of $\GL_2$ over a totally real number field in
the article \cite{BL} of Brylinski and Labesse. In these two cases, the
Shimura variety is non-compact (but its Baily-Borel compactification is
not too complicated) and the group has no endoscopy.

One of the simplest cases where the Shimura variety is non-compact and
the group has non-trivial endoscopy is that of the unitary group
$\GU(2,1)$. This case has been studied in the book \cite{LR}. For groups
of semi-simple rank $1$, Rapoport proved in \cite{Ra} the (unstabilized)
fixed point formula in the case where the Hecke correspondance is trivial.
The stabilized fixed point formula for the symplectic groups $\GSp_{2n}$ is
proved in \cite{M3} (and some applications to $\GSp_4$ and $\GSp_6$, similar
to corollary A, are given). In \cite{Lau} and \cite{Lau2},
Laumon obtained results similar
to corollary B for the groups $\GSp_4$, by using the cohomology with
compact support (instead of the intersection cohomology).

Finally, note that Shin obtained recently results analogous to
corollary B, and also results about ramified places, by using the
cohomology of Igusa varieties attached to compact unitary Shimura
varieties (cf \cite{Shi1}, \cite{Shi2}, \cite{Shi3}). This builds on
previous work of Harris and Taylor (\cite{HT}).

We give a quick description of the different chapters.

Chapter \ref{points_fixes} contains ``known facts'' about the fixed point
formula. When the Shimura variety is associated to a unitary group over $\Q$
and the Hecke correspondence is trivial, the fixed point formula has been
proved in \cite{M} (theorem 5.3.3.1).
The article \cite{M2} contains the theoretical tools needed to treat the
case of non-trivial Hecke correspondences for Siegel modular varieties
(proposition 5.1.5 and theorem 5.2.2), but does not finish the calculation.
We generalize here the results of \cite{M2} under certain conditions on the
group (that are satisfied by unitary groups over $\Q$ and by symplectic
groups), then use them to calculate the trace on the intersection cohomology
of a Hecke correspondence twisted by a high enough power of Frobenius 
in the case when the Shimura variety and the boundary strata of its Baily-Borel
compactification are of the type considered by Kottwitz in his article
\cite{K-PSSV} (ie PEL of type $A$ or $C$).
The result is given in theorem \ref{th:points_fixes_moi}.

Chapters \ref{groupes} to \ref{stabilisation_K} treat the stabilization of
the fixed point formula for unitary groups over $\Q$.
We prove conjecture (10.1) of the article
\cite{K-SVLR} (corollary \ref{cor:stab_FT_IC}).
Kottwitz stabilized the elliptic part of the fixed point formula in
\cite{K-SVLR}, and the method of this book to stabilize the terms coming from
the boundary is inspired by his method.
The most complicated calculations are at the infinite place, where we need
to show a formula for the values at certain elements of the stable
characters of discrete series (proposition
\ref{prop:calcul_Phi_M}). This formula looks a little like the formulas
established by Goresky, Kottwitz et MacPherson (\cite{GKM} 
theorems 5.1 et 5.2), though it has less terms. This is special to unitary
groups : the analogous formula for symplectic groups (cf section 4 of
\cite{M3}) is much more complicated, and more different from the formulas
of \cite{GKM}.

In chapter \ref{groupes}, we define the unitary groups over $\Q$ that we will
study, as well as their Shimura data, and we recall some facts about their
endoscopy.

Chapter \ref{serie_discrete} contains the calculations at the infinite place.

Chapter \ref{partie_en_p} contains explicit calculations, at an unramified
place of the group, of the Satake transform, the base change map, the
transfer map and the twisted transfer map, and a compatibility result
for the twisted transfer and constant term maps.

In chapter \ref{FT_stable_geometrique}, we recall the stabilization by
Kottwitz of the geometric side of the invariant trace formula when the test
function is stable cuspidal at infinity (cf \cite{K-NP}). This stabilization
relies on the calculation by Arthur of the geometric side of the invariant
trace formula for a function that is stable cuspidal at infinity
(cf \cite{A-L2} formula (3.5) and theorem 6.1), and uses only the fundamental
lemma (and not the weighted fundamental lemma). Unfortunately, this result
is unpublished. 
Chapter \ref{FT_stable_geometrique} also contains the normalization of the
Haar measures and of the transfer factors, the statement of the fundamental
lemmas that we use and a summary of the results that are known about these
fundamental lemmas.

In chapter \ref{stabilisation_K}, we put the results of chapters
\ref{groupes}, \ref{serie_discrete} and \ref{partie_en_p} together and
stabilize the fixed point formula.

Chapter \ref{applications} gives applications of the stabilized fixed point
formula that do not use base change to $\GL_n$. First, in section
\ref{applications0}, we show how to make the results of this chapter
formally independent from Kottwitz's unpublished article \cite{K-NP}
(this is merely a formal game, because of course a large part of this book was
inspired by \cite{K-NP} in the first place).
In \ref{applications1}, we express the logarithm of the $L$-function of the
intersection complex at a finite place above a big enough prime number as
a sum (a priori with complex coefficients) of logarithms of local $L$-functions
of automorphic representations of unitary groups.
We also give, in section \ref{applications2}, an application to the
Ramanujan-Petersson
conjecture (at unramified places) for certain discrete automorphic
representations of unitary groups.

Chapter \ref{GL_n_applications} gives applications of the stabilized
fixed point formula that use base change to $\GL_n$. In section
\ref{GL_n_applications1}, we recall some facts about non-connected groups.
In sections \ref{GL_n_applications2} and \ref{GL_n_applications3}, we study
the twisted trace formula for certain test functions. We give applications
of this in \ref{GL_n_applications4}; in particular, we obtain another
formula for the $L$-function of the intersection complex, this time in
terms of local $L$-functions of automorphic
representations of general linear groups.
The simple twisted trace formula proved in this chapter implies some
weak base change results; these have been worked out in section
\ref{GL_n_applications5}.

In chapter \ref{lemme_fondamental}, we prove the particular case of
the twisted fundamental lemma that is used in the stabilization of the
fixed point formula in the article \cite{K-SVLR} of
Kottwitz, the article \cite{Lau} of Laumon and chapter \ref{stabilisation_K}.
The methods of this chapter are not new, and no attempt was made to obtain
the most general result possible.
Waldspurger showed in \cite{Wa3} that the twisted fundamental lemma
for the unit of the Hecke algebra is a consequence of the ordinary fundamental
lemma (and, in the general case, of the non-standard fundamental lemma).
We show that, in the particular case that we need, the twisted fundamental
lemma for the unit of the Hecke algebra implies the twisted fundamental
lemma for all the functions of the Hecke algebra. The method is the same as
in the article \cite{Ha} of Hales (ie it is the method inspired by the
article \cite{Cl-LF} of Clozel, and by the remark of the referee of this
article).

The appendix by Robert Kottwitz contains a comparison theorem between the
twisted transfer factors of \cite{KS} and of \cite{K-SVLR}. This result is
needed to use the twisted fundamental lemma in the stabilization of the
fixed point formula.

It is a great pleasure for me to thank Robert Kottwitz and Gérard Laumon.
Robert Kottwitz very kindly allowed me to read his unpublished manuscript
\cite{K-NP}, that has been extremely helpful
to me in writing this text. He also
helped me fix a problem in the proof of proposition
\ref{prop:FT_invariante_tordue}, pointed out several mistakes in chapter
\ref{lemme_fondamental} and accepted to write his proof of the comparison
of twisted transfer factors as an appendix of this book.
Gérard Laumon suggested that I study the intersection cohomology of
non-compact Shimura varieties and has spent countless hours patiently
explainig the subject to me. I also thank Jean-Loup Waldspurger for
sending me a complete version of his manuscript \cite{Wa3} on twisted
endoscopy before it was published.

I am grateful to the other mathematicians who have answered my questions or
pointed out simpler arguments to me, in particular
Pierre-Henri Chaudouard, Laurent Fargues,
Günter Harder, Colette Moeglin, Bao Chau Ngo,
Sug Woo Shin and Marie-France Vignéras (I am especially grateful to Sug Woo
Shin for repeatedly correcting my misconceptions about the spectral side of
the twisted trace formula).

Finally, I would like to express my gratitude to the anonymous referee for
finding several mistakes and inaccuracies in the first version of this text.

This text was written entirely while I was a Clay Research Fellow of the
Clay Mathematics Institute, and worked as a member
at the Institute for Advanced Study in Princeton. Moreover, I have been
partially supported by the National Science Foundation under agreements
number DMS-0111298 and DMS-0635607.

%\end{thepreface}

\mainmatter

\chapter{The fixed point formula}
\label{points_fixes}

\section{Shimura varieties}
\label{points_fixes1}

The reference for this section is \cite{P2} \S3.

Let $\SD=R_{\C/\R}\Gr_{m,\R}$\index{SD@$\SD$}.
Identify $\SD(\C)=(\C\otimes_\R\C)^
\times$ and $\C^\times\times\C^\times$ using the morphism
$a\otimes 1+b\otimes i\fle (a+ib,
a-ib)$, and write $\mu_0:\Gr_{m,\C}\fl\SD_\C$\index{$\mu_0$}
for the morphism $z\fle (z,1)$.

The definition of (pure) Shimura data that will be used here is that of
\cite{P2} (3.1),
up to condition (3.1.4).
So a pure Shimura datum\index{pure Shimura datum}\index{Shimura datum}
is a triple
$(\G,\X,h)$ (that will often be written simply $(\G,\X)$), where $\G$
is a connected reductive linear algebraic group over
$\Q$, $\X$ is a set with a transitive action of
$\G(\R)$, and $h:\X\fl\Hom(\SD,\G_\R)$ is a
$\G(\R)$-equivariant morphism, satisfying conditions (3.1.1), (3.1.2), (3.1.3)
and (3.1.5) of \cite{P2}, but not necessarily condition (3.1.4) (ie the group
$\G^{ad}$ may have a simple factor of compact type defined over $\Q$).

Let $(\G,\X,h)$ be a Shimura datum.
The field of definition $F$ of the conjugacy class of cocharacters
$h_x\circ\mu_0:\Gr_{m,\C}\fl\G_\C$, $x\in\X$, is called the
\emph{reflex field}\index{reflex field}
\index{F@$F$\quad reflex field}
of the datum.
If $\K$ is an open compact subgroup of $\G(\Af)$, there is an associated
Shimura variety\index{Shimura variety}
$M^{\K}(\G,\X)$, that is a quasi-projective algebraic
variety over $F$ satisfying
\[M^{\K}(\G,\X)(\C)=\G(\Q)\sous(\X\times\G(\Af)/\K).\]
\index{MKGX@$M^{\K}(\G,\X)$\quad Shimura variety}
If moreover $\K$ is \emph{neat}\index{neat}
(cf \cite{P1} 0.6), 
then $M^{\K}(\G,\X)$ is smooth over $F$.
Let $M(\G,\X)$ be the inverse limit of the $M^{\K}(\G,\X)$, taken over
the set of open compact subgroups $\K$ of $\G(\Af)$.

Let $g,g'\in\G(\Af)$, and let $\K,\K'$ be open compact subgroups of $\G(\Af)$
such that $\K'\subset g\K g^{-1}$. Then there is a finite morphism
\[T_g:M^{\K'}(\G,\X)\fl M^{\K}(\G,\X),\]
\index{Tg@$T_g$\quad Hecke operator}
that is given on complex points by
\[\left\{\begin{array}{rcl}\G(\Q)\sous(\X\times\G(\Af)/\K') & \fl &
\G(\Q)\sous(\X\times\G(\Af)/\K) \\
\G(\Q)(x,h\K') & \fle & \G(\Q)(x,hg\K)\end{array}\right.\]
If $\K$ is neat, then the morphism $T_g$ is étale.

Fix $\K$. The Shimura variety $M^{\K}(\G,\X)$ is not projective over $F$
in general, but it has a compactification
$j:M^{\K}(\G,\X)\fl M^{\K}(\G,\X)^*$, the Satake-Baily-Borel (or Baily-Borel,
or minimal Satake, or minimal) compactification,\index{Baily-Borel
compactification} such that
$M^{\K}(\G,\X)^*$ is a normal projective variety over $F$ and $M^{\K}(\G,\X)$
is open dense in $M^{\K}(\G,\X)^*$.
Note that $M^{\K}(\G,\X)^*$ is not smooth in
general (even when $\K$ is neat).
The set of complex points of $M^{\K}(\G,\X)^*$ is
\[M^{\K}(\G,\X)^*(\C)=\G(\Q)\sous(\X^*\times\G(\Af)/\K),\]
where $\X^*$ is a topological space having $\X$ as an open dense subset
and such that the $\G(\Q)$-action on $\X$ extends to a continuous
$\G(\Q)$-action on $\X^*$.
As a set, $\X^*$ is the disjoint union of $\X$ and of boundary components
$\X_P$ indexed by the set of admissible parabolic subgroups of $\G$
(a parabolic subgroup of $\G$ is called \emph{admissible}
\index{admissible parabolic subgroup}
if it is not
equal to $\G$ and if its image in every simple factor $\G'$ of $\G^{ad}$
is equal to $\G'$ or to a maximal parabolic subgroup of $\G'$,
cf \cite{P1} 4.5).
If $\Pa$ is an admissible parabolic subgroup of $\G$, then
$\Pa(\Q)=Stab_{\G(\Q)}(\X_P)$; the $\Pa(\Q)$-action on $\X_P$ extends to a
transitive $\Pa(\R)$-action, and the unipotent radical of $\Pa$ acts
trivially on $\X_P$.

For every $g,\K,\K'$ as above, there is a finite morphism
$\overline{T}_g:M^{\K'}(\G,\X)^*\fl M^{\K}(\G,\X)^*$ extending the
morphism $T_g$.

From now on, we will assume that $\G$ satisfies the following condition :
Let $\Pa$ be an admissible parabolic subgroup of $\G$, let $\N_P$ be its
unipotent radical, $\U_P$ the center of $\N_P$ and $\M_P=\Pa/\N_P$ the Levi
quotient. Then there exists two connected reductive subgroups $\Le_P$ and
$\G_P$ of $\M_P$ such that :
\begin{bulletlist}
\item $\M_P$ is the direct product of $\Le_P$ and $\G_P$;
\item $\G_P$ contains $\G_1$, where $\G_1$ is the normal subgroup of
$\M_P$ defined by Pink in \cite{P2} (3.6) (on page 220), and the quotient
$\G_P/\G_1Z(\G_P)$ is $\R$-anisotropic;
\item $\Le_P\subset Cent_{\M_P}(\U_P)\subset Z(\M_P)\Le_P$;
\item $\G_P(\R)$ acts transitively on $\X_P$, and $\Le_P(\R)$ acts
trivially on $\X_P$;
\item for every neat open compact subgroup $\K_M$ of $\M_P(\Af)$,
$\K_M\cap\Le_P(\Q)=\K_M\cap Cent_{\M_P(\Q)}(\X_P)$.

\end{bulletlist}
Denote by $\QP_P$ the inverse image of $\G_P$ in $\Pa$.

\begin{remark}
If $\G$ satisfies this condition, then, for every admissible parabolic
subgroup $\Pa$ of $\G$, the group $\G_P$ satisfies the same condition.

\end{remark}

\begin{example} Any interior form of the general symplectic group
$\GSp_{2n}$ or of the quasi-split unitary group $\GU^*(n)$ defined in
\ref{groupes1} satisfies the condition.

\end{example}

The boundary of $\M^{\K}(\G,\X)^*$ has a natural stratification (this
stratification exists in general, but its description is a little
simpler when $\G$ satisfies the above condition). Let $\Pa$ be an
admissible parabolic subgroup of $\G$. Pink has defined a morphism
$\X_P\fl\Hom(\SD,\G_{P,\R})$ (\cite{P2} (3.6.1)) such that $(\G_P,\X_P)$
is a Shimura datum, and the reflex field of $(\G_P,\X_P)$ is $F$.
Let $g\in\G(\Af)$. 
Let
$\Hr_P=g\K g^{-1}\cap\Pa(\Q)\QP_P(\Af)$, %\index{HrP@$\Hr_P$} 
$\Hr_L=g\K g^{-1}\cap\Le_P(\Q)\N_P(\Af)$, %\index{HrL@$\Hr_L$} 
$\K_Q=g\K g^{-1}\cap\QP_P(\Af)$, %\index{KQ@$\K_Q$}
$\K_N=g\K g^{-1}\cap\N_P(\Af)$. %\index{KN@$\K_N$}
Then (cf \cite{P2} (3.7)) there is a morphism, finite over its image,
\[M^{\K_Q/\K_N}(\G_P,\X_P)\fl M^{\K}(\G,\X)^*-M^{\K}(\G,\X).\]
The group $\Hr_P$ acts on the right on
$M^{\K_Q/\K_N}(\G_P,\X_P)$, and this action factors through the finite
group $\Hr_P/\Hr_L\K_Q$. Denote by $i_{P,g}$\index{iPg@$i_{P,g}$}
the locally closed immersion
\[M^{\K_Q/\K_N}(\G_P,\X_P)/\Hr_P\fl M^{\K}(\G,\X)^*.\]
This immersion extends to a finite morphism
\[\overline{i}_{P,g}:M^{\K_Q/\K_N}(\G_P,\X_P)^*/\Hr_P\fl M^{\K}(\G,\X)^*\]
(this morphism is not a closed immersion in general).
The boundary of $M^{\K}(\G,\X)^*$ is the union of the images of the
morphisms $i_{P,g}$, for $\Pa$ an admissible parabolic subgroup of
$\G$ and $g\in\G(\Af)$.
If $\Pa'$ is another admissible parabolic subgroup of $\G$ and
$g'\in\G(\Af)$, then the images of the immersions $i_{P,g}$ and
$i_{P',g'}$ are equal if and only if there exists $\gamma\in\G(\Q)$ such that
$\Pa'=\gamma\Pa\gamma^{-1}$ and $\Pa(\Q)\QP_P(\Af)g\K=\Pa(\Q)\QP_P(\Af)
\gamma^{-1}g'\K$; if there is no such $\gamma$, then these images are
disjoint.
If $\K$ is neat, then $\K_Q/\K_N$ is also neat and the action of
$\Hr_P/\Hr_L\K_N$ on $M^{\K_Q/\K_N}(\G_P,\X_P)$ is free (so $M^{\K_Q/\K_N}
(\G_P,\X_P)/\Hr_P$ is smooth).

The images of the morphisms $i_{P,g}$, $g\in\G(\Af)$, are the
\emph{boundary strata}\index{boundary stratum}
of $M^{\K}(\G,\X)^*$ associated to $\Pa$.

To simplify notations, assume from now on that $\G^{ad}$ is simple.
Fix a minimal parabolic subgroup $\Pa_0$ of $\G$. A parabolic subgroup
of $\G$ is called \emph{standard}\index{standard parabolic subgroup}
if it contains $\Pa_0$. Let
$\Pa_1,\dots,\Pa_n$ be the maximal standard parabolic subgroups of $\G$,
with the numbering satisfying : 
$r\leq s$ if and only if $\U_{P_r}\subset\U_{P_s}$ (cf \cite{GHM} (22.3)).
Write $\N_r=\N_{P_r}$, $\G_r=\G_{P_r}$, $\Le_r=\Le_{P_r}$, $i_{r,g}=i_{P_r,
g}$, etc. 

Let $\Pa$ be a standard parabolic subgroup of $\G$. Write $\Pa=\Pa_{n_1}\cap
\dots\cap\Pa_{n_r}$, with $n_1<\dots<n_r$. 
The Levi quotient $\M_P=\Pa/\N_P$ is the direct product of $\G_{n_r}$ and
of a Levi subgroup $\Le_P$ of $\Le_{n_r}$.
Let $\Ch_P$\index{CP:@$\Ch_P$}
be the set of $n$-uples $(X_1,\dots,X_r)$, where :
\begin{bulletlist}
\item $X_1$ is a boundary stratum of $M^{\K}(\G,\X)^*$ associated
to $\Pa_{n_1}$;
\item for every $i\in\{1,\dots,r-1\}$, $X_{i+1}$ is a boundary stratum of
$X_i$ associated to the maximal parabolic subgroup
$(\Pa_{n_{i+1}}\cap\QP_{n_i})/\N_{n_i}$ of $\G_{n_i}$.

\end{bulletlist}
Let $\Ch^1_P$ be the quotient of $\G(\Af)\times\QP_{n_1}(\Af)\times\dots\times
\QP_{n_{r-1}}(\Af)$ by the following equivalence relation :
$(g_1,\dots,g_r)$ is equivalent to $(g'_1,\dots,g'_r)$ if and only if, for
every $i\in\{1,\dots,r\}$,
\[(\Pa_{n_1}\cap\dots\cap\Pa_{n_i})(\Q)\QP_{n_i}(\Af)g_i\dots g_1\K=
(\Pa_{n_1}\cap\dots\cap\Pa_{n_i})(\Q)\QP_{n_i}(\Af)g'_i\dots g'_1\K.\]

\begin{proposition}
\label{prop:description_chaines}
\begin{itemize}
\item[(i)] The map $\G(\Af)\fl\Ch^1_P$ that sends $g$ to the class of
$(g,1,\dots,1)$ induces a bijection $\Pa(\Q)\QP_{n_r}(\Af)\sous\G(\Af)
/\K\iso\Ch^1_P$.

\item[(ii)] Define a map $\varphi':\Ch^1_P\fl\Ch_P$ in the following way :
Let $(g_1,\dots,g_r)\in\G(\Af)\times\QP_{n_1}(\Af)\times
\dots\times\QP_{n_{r-1}}(\Af)$. For every $i\in\{1,\dots,r\}$, write 
\[\Hr_i=(g_i\dots g_1)\K (g_i\dots g_1)^{-1}\cap(\Pa_{n_1}\cap\dots\cap
\Pa_{n_i})(\Q)\QP_{n_i}(\Af)\]
and let $\K_i$ be the image of $\Hr_i\cap\QP_{n_i}(\Af)$ by the obvious
morphism $\QP_{n_i}(\Af)\fl\G_{n_i}(\Af)$. 
Then $\varphi'$ sends the class of $(g_1,\dots,g_r)$ to the $n$-uple
$(X_1,\dots,X_r)$, where $X_1=Im(i_{n_1,g_1})=M^{\K_1}(\G_{n_1},\X_{n_1})/
\Hr_1$ and, for every
$i\in\{1,\dots,r-1\}$, $X_{i+1}$ is the boundary stratum of
$X_i=M^{\K_i}(\G_{n_i},\X_{n_i})/\Hr_i$ image of the morphism $i_{P,g}$,
with $\Pa=(\Pa_{n_{i+1}}\cap\QP_{n_i})/\N_{n_i}$ (a maximal parabolic subgroup
of $\G_{n_i}$) and $g=g_{i+1}\N_{n_i}(\Af)\in\G_{n_i}(\Af)$.

Then this map $\Ch^1_P\fl\Ch_P$ is well-defined and bijective.

\end{itemize}
\end{proposition}

The proposition gives a bijection
$\varphi_P:\Pa(\Q)\QP_{n_r}(\Af)\sous\G(\Af)/\K
\iso\Ch_P$. On the other hands, there is a map from $\Ch_P$ to the set of
boundary strata of $M^{\K}(\G,\X)^*$ associated to $\Pa_{n_r}$, defined
by sending $(X_1,\dots,X_r)$ to the image of $X_r$ in $M^{\K}(\G,\X)^*$.
After identifying $\Ch_P$ to $\Pa(\Q)\QP_{n_r}(\Af)\sous\G(\Af)/\K$
using $\varphi_P$ and the second set to
$\Pa_{n_r}(\Q)\QP_{n_r}(\Af)\sous\G(\Af)/\K$ using
$g\fle Im(i_{n_r,g})$, this map becomes the obvious projection
$\Pa(\Q)\QP_{n_r}(\Af)\sous\G(\Af)/\K\fl\Pa_{n_r}(\Q)\QP_{n_r}(\Af)\sous
\G(\Af)/\K$.

\begin{proof}\begin{itemize}
\item[(i)] As $\QP_{n_r}\subset\QP_{n_{r-1}}\subset\dots\subset\QP_{n_1}$,
it is easy to see that, in the definition of $\Ch^1_P$, $(g_1,\dots,g_r)$
is equivalent to $(g'_1,\dots,g'_r)$ if and only if
\[(\Pa_{n_1}\cap\dots\cap\Pa_{n_r})(\Q)\QP_{n_r}(\Af)g_r\dots g_1\K=(\Pa_{n_1}
\cap\dots\cap\Pa_{n_r})(\Q)\QP_{n_r}(\Af)g'_r\dots g'_1\K.\]
The results now follows from the fact that
$\Pa=\Pa_{n_1}\cap\dots\cap\Pa_{n_r}$.

\item[(ii)] We first check that $\varphi'$ is well-defined. Let
$i\in\{1,\dots,r-1\}$. If $X_i=M^{\K_i}(\G_{n_i},\X_{n_i})/\Hr_i$ and
$X_{i+1}$ is the boundary stratum $Im(i_{P,g})$ of $X_i$, with $\Pa$ and $g$
as in the proposition, then $X_{i+1}=M^{\K'}(\G_{n_
{i+1}},\X_{n_{i+1}})/\Hr'$, where $\Hr'=g_{i+1}\Hr_i g_{i+1}^{-1}\cap
\Pa_{n_{i+1}}(\Q)\QP_{n_{i+1}}(\Af)$ and $\K'$ is the image of $\Hr'\cap\QP_
{n_{i+1}}(\Af)$ in $\G_{n_{i+1}}(\Af)$. As $g_{i+1}\in\QP_{n_i}(\Af)$,
\[\Hr'=(g_{i+1}\dots g_1)\K(g_{i+1}\dots g_1)^{-1}\cap (\Pa_{n_1}\cap\dots\cap
\Pa_{n_i})(\Q)\QP_{n_i}(\Af)\cap\Pa_{n_{i+1}}(\Q)\QP_{n_{i+1}}(\Af).\]
On the other hand, it is easy to see that
\[(\Pa_{n_1}\cap\dots\cap\Pa_{n_i})(\Q)\QP_{n_i}(\Af)\cap\Pa_{n_{i+1}}(\Q)
\QP_{n_{i+1}}(\Af)=(\Pa_{n_1}\cap\dots\cap\Pa_{n_{i+1}})(\Q)\QP_{n_{i+1}}(\Af).
\]
Hence $\Hr'=\Hr_{i+1}$, and $X_{i+1}=M^{\K_{i+1}}(\G_{n_{i+1}},\X_{n_{i+1}})/
\Hr_{i+1}$. It is also clear that the $n$-uple $(X_1,\dots,X_r)$ defined
in the proposition doesn't change if $(g_1,\dots,g_r)$ is replaced by an
equivalent $r$-uple.

It is clear that $\varphi'$ is surjective. We want to show that it is
injective.
Let $c,c'\in\Ch^1_P$; write $(X_1,\dots,X_r)=\varphi'(c)$ and
$(X'_1,\dots,X'_r)=\varphi'(c')$, and suppose that $(X_1,\dots,X_r)=
(X'_1,\dots,X'_r)$. Fix representatives $(g_1,\dots,g_n)$ and $(g'_1,
\dots,g'_n)$ of $c$ and $c'$. As before, write, for every
$i\in\{1,\dots,n\}$,
\[\Hr_i=(g_i\dots g_1)\K(g_i\dots g_1)^{-1}\cap(\Pa_{n_1}\cap\dots\cap\Pa_
{n_i})(\Q)\QP_{n_i}(\Af)\]
\[\Hr'_i=(g'_i\dots g'_1)\K(g'_i\dots g'_1)^{-1}\cap(\Pa_{n_1}\cap\dots\cap\Pa_
{n_i})(\Q)\QP_{n_i}(\Af).\]
Then the equality $X_1=X'_1$ implies that $\Pa_{n_1}(\Q)\QP_{n_1}(\Af)g_1\K=
\Pa_{n_1}(\Q)\QP_{n_1}(\Af)g'_1\K$ and, for every $i\in\{1,\dots,r-1\}$, the
equality $X_{i+1}=X'_{i+1}$ implies that
\[\Pa_{n_{i+1}}(\Q)\QP_{n_{i+1}}(\Af)g_{i+1}\Hr_i(g_i\dots g_1)=
\Pa_{n_{i+1}}(\Q)\QP_{n_{i+1}}(\Af)g'_{i+1}\Hr'_i(g'_i\dots g'_1).\]
So $(g_1,\dots,g_r)$ and $(g'_1,\dots,g'_r)$ are equivalent, and $c=c'$.

\end{itemize}
\end{proof}

 % variétés de Shimura
\section{Local systems and Pink's theorem}
\label{points_fixes2}

Fix a number field $K$. If $\G$ is a linear algebraic group over $\Q$,
let $Rep_\G$\index{RepG@$Rep_\G$}
be the category of algebraic representations of $\G$ defined over
$K$.
Fix a prime number $\ell$ and a place $\lambda$ of $K$ over $\ell$.

Let $\M$ be a connected reductive group over $\Q$,
$\Le$ and $\G$ connected reductive subgroups of $\M$ such that $\M$ is
the direct product of $\Le$ and $\G$, and $(\G,\X)$ a Shimura datum.
Extend the $\G(\Af)$-action on $M(\G,\X)$ to a $\M(\Af)$-action by the
obvious map $\M(\Af)\fl\G(\Af)$ (so $\Le(\Af)$ acts trivially). Let $\K_M$
be a neat open compact subgroup of $\M(\Af)$.
Write $\Hr=\K_M\cap\Le(\Q)
\G(\Af)$, $\Hr_L=\K_M\cap\Le(\Q)$ (an arithmetic subgroup of
$\Le(\Q)$) and $\K=\K_M\cap\G(\Af)$. The group $\Hr$ acts on the Shimura
variety $M^{\K}(\G,\X)$, and the quotient $M^{\K}(\G,\X)/\Hr$ is equal to
$M^{\Hr/\Hr_L}(\G,\X)$ ($\Hr/\Hr_L$ is a neat open compact subgroup of
$\G(\Af)$).

\begin{remark} It is possible to generalize the morphisms $T_g$ of
\ref{points_fixes1} : If $m\in\Le(\Q)\G(\Af)$ and $\K'_M$ is an open
compact subgroup of $\M(\Af)$ such that
$\K'_M\cap\Le(\Q)\G(\Af)\subset m\Hr m^{-1}$, then there is a morphism
\[T_m:M(\G,\X)/\Hr'\fl M(\G,\X)/\Hr,\]
where $\Hr'=\K'_M\cap\Le(\Q)\G(\Af)$ and $\Hr=\K_M\cap\Le(\Q)\G(\Af)$.
This morphism is simply the one induced by the injection $\Hr'\fl\Hr$,
$h\fle mhm^{-1}$ (equivalently, it is induced by the endomorphism
$x\fle xm$ of $M(\G,\X)$).

\end{remark}

There is an additive triangulated functor $V\fle\F^{\Hr/\Hr_L}R\Gamma(\Hr_L,V)$
\index{FHrHrL@$\F^{\Hr/\Hr_L}$}
from the category $D^b(Rep_\M)$ to the category of $\lambda$-adic complexes
on $M^{\K}(\G,\X)/\Hr$,\footnote{Here, and in the rest of the book, the
notation $R\Gamma$ will be used to denote the right derived functor of the
functor $\H^0$.}
\index{RGamma@$R\Gamma$}
constructed using the functors $\mu_{\Gamma,\varphi}$
of Pink (cf \cite{P1} (1.10)) for the profinite étale (and Galois of
group $\Hr/\Hr_L$) covering $M(\G,\X)\fl M^{\K}(\G,\X)/\Hr$ and the properties
of the arithmetic subgroups of $\Le(\Q)$. This construction is explained
in \cite{M} 2.1.4. For every $V\in\Ob D^b(Rep_\M)$ and $k\in\Z$,
$\Ho^k\F^{\Hr/\Hr_L}R\Gamma(\Hr_L,V)$ is a lisse $\lambda$-adic sheaf on
$M^{\K}(\G,\X)/\Hr$, whose fiber is (noncanonically) isomorphic to
\[\bigoplus_{i+j=k}\Ho^i(\Hr_L,\Ho^jV).\]

\begin{remark} If $\Gamma$ is a neat arithmetic subgroup of $\Le(\Q)$
(eg $\Gamma=\Hr_L$), then it is possible to compute $R\Gamma(\Gamma,V)$ 
in the category $D^b(Rep_\G)$, because $\Gamma$ is of type FL (cf \cite{BuW}, 
theorem 3.14). 

\end{remark}

We will now state a theorem of Pink about the direct image of the complexes
$\F^{\Hr/\Hr_L}R\Gamma(\Hr_L,V)$
by the open immersion $j:M^{\K}(\G,\X)/\Hr\fl M^{\K}(\G,\X)^*/\Hr$.
Let $\Pa$ be an admissible parabolique subgroup of $\G$ and $g\in\G(\Af)$.
Write
\[\Hr_P=g\Hr g^{-1}\cap\Le(\Q)\Pa(\Q)\QP_P(\Af),\]
\[\Hr_{P,L}=g\Hr g^{-1}\cap\Le(\Q)\Le_P(\Q)\N_P(\Af),\]
\[\K_N=g\Hr g^{-1}\cap\N_P(\Af),\]
\[\K_G=(g\Hr g^{-1}\cap\QP_P(\Af))/(g\Hr g^{-1}\cap\N_P(\Af)),\]
and $i=i_{P,g}:M^{\K_G}(\G_P,\X_P)/\Hr_P\fl M^{\K}(\G,\X)^*/\Hr$.

Then theorem 4.2.1 of \cite{P2} implies the following results (cf \cite{M}
2.2) :

\begin{theorem}\label{th:Pink}\index{Pink's theorem}
For every $V\in\Ob D^b(Rep_\M)$, there are canonical isomorphisms
\[\begin{array}{rcl}\displaystyle{i^*Rj_*\F^{\Hr/\Hr_L}R\Gamma(\Hr_L,V)} & 
\simeq & \displaystyle{\F^{\Hr_P/\Hr_{P,L}}R\Gamma(\Hr_{P,L},V)} \\
& \simeq & \displaystyle{\F^{\Hr_P/\Hr_{P,L}}R\Gamma(\Hr_{P,L}/\K_N,R\Gamma
(Lie(\N_P),V)).}\end{array}\]

\end{theorem}

The last isomorphism uses van Est's theorem, as stated (and proved) in
\cite{GKM} \S24.

\vspace{.5cm}

We will also use local systems on locally symmetric spaces that are
not necessarily Hermitian. We will need the following notation. Let $\G$ be
a connected reductive group over $\Q$. Fix a maximal compact subgroup
$\K_\infty$ of $\G(\R)$.
Let $\A_G$\index{AG@$\A_G$}
be the maximal ($\Q$-)split torus of te center of $\G$,
$\X=\G(\R)/\K_\infty\A_G(\R)^0$ and
$q(\G)=\dim(\X)/2\in\frac{1}{2}\Z$. Write
\[M^{\K}(\G,\X)(\C)=\G(\Q)\sous (\X\times\G(\Af)/\K)\]
(even though $(\G,\X)$ is not a Shimura datum in general, and
$\M^{\K}(\G,\X)(\C)$ is not always the set of complex points of an algebraic
variety).
If $\K$ is small enough (eg neat), this quotient is a real analytic
variety. There are morphisms $T_g$ ($g\in\G(\Af)$) defined exactly as in
\ref{points_fixes1}.

Let $V\in\Ob Rep_\G$. Let $\F^{\K}V$\index{FKV@$\F^{\K}V$}
be the sheaf of local sections of the morphism
\[\G(\Q)\sous (V\times\X\times\G(\Af)/\K)\fl\G(\Q)\sous(\X\times\G(\Af)/\K)\]
(where $\G(\Q)$ acts on $V\times\X\times\G(\Af)/\K$ by $(\gamma,(v,x,g\K))\fle
(\gamma.v,\gamma.x,\gamma g\K)$). As suggested by the notation, there is a
connection between this sheaf and the local systems defined above :
if $(\G,\X)$ is a Shimura datum, then $\F^{\K}V\otimes K_\lambda$ is the
inverse image on $M^{\K}(\G,\X)(\C)$ of the $\lambda$-adic sheaf
$\F^{\K}V$ on $M^{\K}(\G,\X)$ (cf \cite{L1} p 38 or \cite{M} 2.1.4.1)

Let $\Gamma$ be a neat arithmetic subgroup of $\G(\Q)$. Then the quotient
$\Gamma\sous\X$ is a real analytic variety. For every $V\in\Ob Rep_\G$,
let $\F^{\Gamma}V$ be the sheaf of local sections of the morphism
\[\Gamma\sous (V\times\X)\fl\Gamma\sous\X\]
(where $\Gamma$ acts on $V\times\X$ by $(\gamma,(v,x))\fle (\gamma.v,\gamma.x)$
). 

Let $\K$ be a neat open compact subgroup of $\G(\Af)$, and let
$(g_i)_{i\in I}$ be a system of representatives of the double quotient
$\G(\Q)\sous\G(\Af)/\K$. For every $i\in I$, let $\Gamma_i=g_i\K g_i^{-1}
\cap\G(\Q)$. Then the $\Gamma_i$ are neat arithmetic subgroups of
$\G(\Q)$, 
\[M^{\K}(\G,\X)(\C)=\coprod_{i\in I}\Gamma_i\sous\X\]
and, for every $V\in\Ob Rep_\G$,
\[\F^{\K}V=\bigoplus_{i\in I}\F^{\Gamma_i}V.\]

 % systèmes locaux et théorème de Pink
\section{Integral models}
\label{points_fixes3}

Notations are as in \ref{points_fixes1}. Let
$(\G,\X)$ be a Shimura datum such that $\G^{ad}$ is simple and that the
maximal parabolic subgroups of $\G$ satisfy the condition
of \ref{points_fixes1}. The goal here is to show that there exist integral
models
\index{integral model of a Shimura variety}
(ie models over a suitable localization of $\Of_F$) of the varieties and
sheaves of \ref{points_fixes1} and \ref{points_fixes2} such that Pink's
theorem is still true. The exact conditions that we want these models
to satisfy are given more precisely below (conditions (1)-(7)).

Fix a minimal parabolic subgroup $\Pa_0$ of $\G$,
and let $(\G_1,\X_1),\dots,(\G_n,\X_n)$ be the Shimura data associated
to the standard maximal parabolic subgroups of $\G$.
We will also write $(\G_0,\X_0)=(\G,\X)$.
Note that, for every $i\in\{0,\dots,n\}$, $\Pa_0$ determines a minimal
parabolic subgroup of $\G_i$.
It is clear that, for every $i\in\{0,\dots,n\}$, the Shimura data
associated to the maximal parabolic subgroups of $\G_i$ are the $(\G_j,\X_j)$,
with $i+1\leq j\leq n$.

Remember that $F$ is the reflex field of $(\G,\X)$. It is also the reflex
field of all the $(\G_i,\X_i)$ (\cite{P1} 12.1 and 11.2(c)).
Let $\overline{\Q}$ be the algebraic closure of $\Q$ in $\C$;
as $F$ is by definition a subfield of $\C$, it is included in $\overline{\Q}$.
For every prime number $p$, fix an algebraic closure $\overline{\Q}_p$
of $\Q_p$ and an injection $F\subset\overline{\Q}_p$.

Fix a point $x_0$ of $\X$, and let $h_0:\SD\fl\G_\R$ be the morphism
corresponding to $x_0$.
Let $w$ be the composition of $h_0$ and of the injection
$\Gr_{m,\R}\subset\SD$.
\index{w@$w$\quad weight morphism}
Then $w$ is independant of the choice of $h_0$ and it
is defined over $\Q$ (cf \cite{P2} 5.4). An algebraic representation
$\rho:\G\fl\GL(V)$ of $\G$ is said to be \emph{pure of weight $m$}
if $\rho\circ w$ is the multiplication by the character $\lambda\fle\lambda^m$
of $\Gr_m$ (note that the sign convention here is not the same as in
\cite{P2} 5.4).
\index{weight of an algebraic representation}

Consider the following data :
\begin{bulletlist}
\item for every $i\in\{0,\dots,n\}$, a set $\Kens_i$ of neat open compact
subgroups of $\G_i(\Af)$, stable by $\G(\Af)$-conjugacy;
\item for every $i\in\{0,\dots,n\}$, a subset $A_i$ of $\G_i(\Af)$ such
that $1\in A_i$;
\item for every $i\in\{0,\dots,n\}$, a full abelian subcategory
$\Rens_i$ of $Rep_{\G_i}$, stable by taking direct factors

\end{bulletlist}
satisfying the following conditions :
Let $i,j\in\{0,\dots,n\}$ be such that $j>i$, and $\K\in\Kens_i$.
Let $\Pa$ be the standard maximal parabolic subgroup of $\G_i$ associated to
$(\G_j,\X_j)$. Then :
\begin{itemize}
\item[(a)] For every $g\in\G_i(\Af)$,
\[(g\K g^{-1}\cap\QP_P(\Af))/(g\K g^{-1}\cap\N_P(\Af))\in\Kens_j,\]
and, for every $g\in\G_i(\Af)$ and every standard parabolic subgroup $\Pa'$
of $\G_i$ such that $\QP_P\subset\Pa'\subset\Pa$,
\[(g\K g^{-1}\cap\Pa'(\Q)\N_{P'}(\Af)\QP_P(\Af))/(g\K g^{-1}\cap
\Le_{P'}(\Q)\N_{P'}(\Af))\in\Kens_j\]
\[(g\K g^{-1}\cap\Pa'(\Af))/(g\K g^{-1}\cap\Le_{P'}(\Af)\N_{P'}(\Af))\in
\Kens_j.\]
\item[(b)] Let $g\in A_i$ and $\K'\in\Kens_i$ be such that $\K'\subset 
g\K g^{-1}$. Let $h\in\Pa(\Q)\QP_P(\Af)\sous\G(\Af)/\K$ and $h'\in
\Pa(\Q)\QP_P(\Af)\sous\G(\Af)/\K'$ be such that $\Pa(\Q)\QP_P(\Af)h\K=
\Pa(\Q)\QP_P(\Af)h'g\K$. Then there exists $p\in\Le_P(\Q)$ and $q\in\QP_P(
\Af)$ such that $pqh\K=h'g\K$ and that the image of $q$ in $\G_j(\Af)=
\QP_P(\Af)/\N_P(\Af)$ is in $A_j$.
\item[(c)] For every $g\in\G_i(\Af)$ and $V\in\Ob\Rens_i$,
\[R\Gamma(\Gamma_L,R\Gamma(Lie(\N_P),V))\in\Ob D^b(\Rens_j),\]
where
\[\Gamma_L=(g\K g^{-1}\cap\Pa(\Q)\QP_P(\Af))/(g\K g^{-1}\cap\QP_P(\Af)).\]

\end{itemize}

Let $\Sigma$ be a finite set of prime numbers such that the groups
$\G_0,\dots,\G_n$ are unramified outside $\Sigma$. For every
$p\not\in\Sigma$, fix $\Z_p$-models of these groups.
Note
\[\Ade_\Sigma=\prod_{p\in\Sigma}\Q_p.\]
Fix $\ell\in\Sigma$ and a place $\lambda$ of $K$ above $\ell$, and 
consider the following conditions on $\Sigma$ :
\begin{itemize}
\item[(1)] For every $i\in\{0,\dots,n\}$, $A_i\subset\G_i(\Ade_
{\Sigma})$ and every $\G(\Af)$-conjugacy class in $\Kens_i$ has a
representative of the form $\K_\Sigma\K^\Sigma$, with
$\K_\Sigma\subset\G_i(\Ade_{\Sigma})$ and $\K^\Sigma=\prod\limits_
{p\not\in\Sigma}\G_i(\Z_p)$.
\item[(2)] For every $i\in\{0,\dots,n\}$ and $\K\in\Kens_i$, there exists
a smooth quasi-projective scheme $\Mod^{\K}(\G_i,\X_i)$ on $\Spec(\Of_F
[1/\Sigma])$ whose generic fiber is $M^{\K}(\G_i,\X_i)$.
\item[(3)] For every $i\in\{0,\dots,n\}$ and $\K\in\Kens_i$, there exists
a normal scheme $\Mod^{\K}(\G_i,\X_i)^*$, projective over
$\Spec(\Of_F[1/\Sigma])$, containing $\Mod^{\K}(\G_i,\X_i)$ as a dense open
subscheme and with generic fiber $M^{\K}(\G_i,\X_i)^*$. Moreover, the morphisms
$i_{P,g}$ (resp. $\overline{i}_{P,g}$) of \ref{points_fixes1} extend to locally
closed immersions (resp. finite morphisms) between the models over
$\Spec(\Of_F[1/\Sigma])$, and the boundary of $\Mod^{\K}(\G_i,\X_i)^*-\Mod^{\K}
(\G_i,\X_i)$ is still the disjoint union of the images of the immersions
$i_{P,g}$.
\item[(4)] For every $i\in\{0,\dots,n\}$, $g\in A_i$ and $\K,\K'\in\Kens_i$
such that $\K'\subset g\K g^{-1}$, the morphism $\overline{T}_g:M^{\K'}(\G_i,
\X_i)^*\fl M^{\K}(\G_i,\X_i)^*$ extends to a finite morphism
$\Mod^{\K'}(\G_i,\X_i)^*\fl\Mod^{\K}(\G_i,\X_i)^*$, that will still be denoted
by $\overline{T}_g$,
whose restriction to the strata of $\Mod^{\K'}(\G_i,\X_i)^*$ (including the
open stratum $\Mod^{\K'}(\G_i,\X_i)$) is étale.
\item[(5)] For every $i\in\{0,\dots,n\}$ and $\K\in\Kens_i$, there exists a
functor $\F^{\K}$ from $\Rens_i$ to the category of lisse $\lambda$-adic
sheaves on $\Mod^{\K}(\G_i,\X_i)$ that, after passing to the special fiber,
is isomorphic to the functor $\F^{\K}$ of \ref{points_fixes2}.
\item[(6)] For every $i\in\{0,\dots,n\}$, $\K\in\Kens_i$ and $V\in\Ob\Rens_i$,
the isomorphisms of Pink's theorem (\ref{th:Pink}) extend to isomorphisms
between $\lambda$-adic complexes on the $\Spec(\Of_F[1/\Sigma])$-models.
\item[(7)] For every $i\in\{0,\dots,n\}$, $\K\in\Kens_i$ and $V\in\Ob\Rens_i$,
the sheaf $\F^{\K}V$ on $\Mod^{\K}(\G_i,\X_i)$ is mixed (\cite{D2} 1.2.2).
If moreover $V$ is pure of weight $m$, then $\F^{\K}V$ is pure
of weight $-m$.

\end{itemize}

The fact that suitable integral models exist for PEL Shimura varieties
has been proved by Kai-Wen Lan, who constructed the toroidal
and minimal compactifications of the integral models.

\begin{proposition}\label{prop:modeles_entiers} Suppose that the Shimura
datum $(\G,\X)$ is of the type considered in \cite{K-PSSV} \S5; more precisely,
we suppose fixed data as in \cite{Lan} 1.2. Let $\Sigma$ be a finite set of
prime numbers that contains all bad primes (in the sense of \cite{Lan}
1.4.1.1).
For every $i\in\{0,\dots,n\}$,
let $A_i=\G_i(\Ade_\Sigma)$ and let $\Kens_i$ be the
union of the $\G_i(\Af)$-conjugacy classes of neat open compact subgroups of
the form $\K_\Sigma\K^\Sigma$ with $\K_\Sigma\subset\prod\limits_{p\in\Sigma}
\G_i(\Z_p)$ and
$\K^\Sigma=\prod\limits_{p\not\in\Sigma}\G_i(\Z_p)$.

Then the set $\Sigma$ satisfies conditions (1)-(7), and moreover the
schemes $\Mod^{\K}(\G_i,\X_i)$ of (2) are the schemes representing
the moduli problem of \cite{Lan} 1.4.

\end{proposition}

\begin{proof} This is just putting together Lan's and Pink's results.
Condition (1) is automatic. Condition (2) (in the more precise form given
in the proposition) is a consequence of theorem 1.4.1.12 of \cite{Lan}.
Conditions (3) and (4) are implied by theorem 7.2.4.1 and proposition
7.2.5.1 of \cite{Lan}. The construction of the sheaves in condition (5)
is the same as in \cite{P2} 5.1, once the integral models of condition (2)
are known to exist. In \cite{P2} 4.9, Pink observed that the proof of his
theorem extends to integral models if toroidal compactifications and a
minimal compactification of the integral model satisfying the properties
of section 3 of \cite{P2} have been constructed. This has been done by Lan
(see, in addition to the results cited above, theorem 6.4.1.1 and propositions
6.4.2.3, 6.4.2.9 and 6.4.3.4 of \cite{Lan}), so condition (6) is also
satisfied.
In the PEL case, $\G^{ad}$ is automatically if abelian type in the sense
of \cite{P2} 5.6.2 (cf \cite{K-PSSV} \S5). So
$\G_i^{ad}$ is of
abelian type for all $i$, and condition (7) is implied by proposition
5.6.2 in \cite{P2}.

\end{proof}

\begin{remark}\label{rq:ensemble_Sigma}
\index{bad primes}
Let $(\G,\X)$ be one of the Shimura
data defined in \ref{groupes1}, and let $\K$ be a neat open compact subgroup
of $\G(\Af)$. Then there exist a finite set $S$ of primes such that
$\K=\K_S\prod_{p\not\in S}\G(\Z_p)$, with $\K_S\subset\prod_{v\in S}\G(\Q_v)$
(and with the $\Z$-structure on $\G$ defined in remark
\ref{rq:groupes_sur_Z}). Let $\Sigma$ be the union of $S$ and of all
prime numbers that are ramified in $E$. Then $\Sigma$ contains all bad
primes, so proposition \ref{prop:modeles_entiers} above applies to
$\Sigma$.

\end{remark}

\begin{remark}\label{rq:modeles_canoniques}
The convention we use here for the action of the Galois group on the
canonical model is that of Pink (\cite{P2} 5.5), that is different from
the convention of Deligne (in \cite{D}) and hence also from the convention
of Kottwitz (in \cite{K-PSSV}); 
so what Kottwitz calls canonical model of the Shimura variety associated
to the Shimura datum
$(\G,\X,h^{-1})$ is here the canonical model of the Shimura variety associated
to the Shimura datum $(\G,\X,h)$.

\end{remark}

Let us indicate another way to find integral models when the Shimura datum is
not necessarily PEL. The
problem with this approach is that the set $\Sigma$ of ``bad'' primes is
unknown.

\begin{proposition}\label{prop:modeles_entiers_bis}
Let $\Kens_i$ and $A_i$ be as above (and satisfying conditions (a) and
(b)). Suppose that, for every $i\in\{0,\dots,n\}$,
$\Kens_i$ is finite modulo $\G_i(\Af)$-conjugacy and $A_i$ is finite.
If $\G^{ad}$ is of abelian type (in the sense of \cite{P2} 5.6.2),
then here exists a finite set $\Sigma$ of
prime numbers satisfying conditions (1)-(7), with
$\Rens_i=Rep_{\G_i}$ for every $i\in\{0,\dots,n\}$.

In general, there exists a finite set $\Sigma$ of prime numbers satisfying
conditions (1)-(6), with $\Rens_i=Rep_{\G_i}$ for every
$i\in\{0,\dots,n\}$. Let $\Rens'_i$, $0\leq i\leq n$, be full subcategories
of the $Rep_{\G_i}$, stable by taking direct factors and by isomorphism,
containing the trivial representation, satisfying condition (c) and minimal
for all these properties (this determines the $\Rens'_i$).
Then there exists $\Sigma'\supset\Sigma$ finite such that
$\Sigma'$, the $\Kens_i$ and the $\Rens'_i$ satisfy condition (7).

\end{proposition}

This proposition will typically be applied to the following situation :
$g\in\G(\Af)$ and
$\K,\K'$ are neat open compact subgroups of $\G(\Af)$ such that
$\K'\subset\K\cap g\K g^{-1}$, and we want to study the Hecke correspondence
$(T_g,T_1):M^{\K'}(\G,\X)^*\fl (M^{\K}(\G,\X)^*)^2$. In order to reduce this
situation modulo $p$, choose sets
$\Kens_i$ such that $\K,\K'\in\Kens_0$ and that condition (a) is satisfied,
and minimal for these properties,
sets $A_i$ such that $1,g\in A_0$ and that condition (b) is satisfied,
and minimal for these properties; take $\Rens_i=Rep_{\G_i}$ if $\G^{ad}$
is of abelian type and $\Rens_i$ equal to the $\Rens'_i$ defined in the
proposition in the other cases; fix $\Sigma$ such that conditions (1)-(7)
are satisfied, and reduce modulo $p\not\in\Sigma$.

\begin{proof} First we show that, in the general case, there is a finite
set $\Sigma$ of prime numbers satisfying conditions (1)-(6), with
$\Rens_i=Rep_{\G_i}$. It is obviously possible to find $\Sigma$ satisfying
conditions (1)-(4). Proposition 3.6 of \cite{W} implies that we can find
$\Sigma$ satisfying conditions (1)-(5).
To show that there exists $\Sigma$ satisfying conditions (1)-(6),
reason as in the proof of proposition 3.7 of \cite{W}, using the generic
base change theorem of Deligne (cf SGA
4 1/2 [Th. finitude] théorème 1.9).
As in the proof of proposition \ref{prop:modeles_entiers}, if $\G^{ad}$ is
of abelian type, then condition (7) is true by proposition 5.6.2 of
\cite{P2}. In the general case, let $\Rens'_i$ be defined as in
the statement of the proposition. Condition (7) for these subcategories
is a consequence of proposition 5.6.1 of \cite{P2} (reason as in the second
proof of \cite{P2} 5.6.6).

\end{proof}

\begin{remark} Note that it is clear from the proof that, after replacing
$\Sigma$ by a bigger finite set, we can choose the
integral models $\Mod^{\K}(\G_i,\X_i)$ to be any integral models specified
before (as long as they satisfy the conditions of (2)).

\end{remark}

When we later talk about reducing
Shimura varieties modulo $p$, we will always implicitely fix $\Sigma$ as
in proposition \ref{prop:modeles_entiers} (or proposition
\ref{prop:modeles_entiers_bis})
and take $p\not\in\Sigma$.
The prime number
$\ell$ will be chosen among elements of $\Sigma$ (or added to $\Sigma$).
\index{reduction modulo $p$ (for Shimura varieties)}

 % modèles entiers 
\section{Weighted cohomology complexes and intersection complex}
\label{points_fixes4}

Let $(\G,\X)$ be a Shimura datum and $\K$ be a neat open compact subgroup
of $\G(\Af)$. Assume that $\G$ satisfies the conditions of \ref{points_fixes}
and that $\G^{ad}$ is simple.
Fix a minimal parabolic subgroup $\Pa_0$ of $\G$ and maximal standard
parabolic subgroups
$\Pa_1,\dots,\Pa_n$ as before proposition \ref{prop:description_chaines}.
Fix prime numbers $p$ and $\ell$ as in the end of \ref{points_fixes3},
and a place $\lambda$ of $K$ above $\ell$.
In this section, we will write $\M^{\K}(\G,\X)$, etc, for the reduction
modulo $p$ of the varieties of \ref{points_fixes1}.

Write $M_0=M^{\K}(\G,\X)$ and $d=\dim M_0$, and, for every $r\in\{1,\dots,n\}$,
denote by $M_r$ the union of the boundary strata of $M^{\K}(\G,\X)^*$
associated to $\Pa_r$, by $d_r$ the dimension of $M_r$ and by $i_r$ the
inclusion of $M_r$ in $M^{\K}(\G,\X)^*$. 
Then $(M_0,\dots,M_n)$ is a stratification of
$M^{\K}(\G,\X)^*$ in the sense of \cite{M2} 3.3.1. Hence, for every
$\as=(a_0,\dots,a_n)\in (\Z\cup\{\pm\infty\})^{n+1}$, the functors
$w_{\leq\as}$ and $w_{>\as}$ of \cite{M2} 3.3.2 are defined
(on the category $D^b_m(M^{\K}(\G,\X)^*,K_\lambda)$ of mixed $\lambda$-adic
complexes on $M^{\K}(\G,\X)^*$).
We will recall the definition of the intersection complex and of the
weighted cohomology complexes.
Remember that $j$ is the open immersion $M^{\K}(\G,\X)\fl
M^{\K}(\G,\X)^*$.

\begin{remark} We will need to use the fact that the sheaves
$\F^{\K}V$ are mixed with known weights.
So we fix categories
$\Rens_0,\dots,\Rens_n$ as in \ref{points_fixes3}, satisfying conditions
(c) and (7) of \ref{points_fixes3}. 
If $\G^{ad}$ is of abelian type, we can simply take $\Rens_0=Rep_\G$.

\end{remark}

\begin{definition}\begin{itemize}
\item[(i)] \index{intersection complex}
Let $V\in\Ob Rep_\G$. The intersection complex on
$M^{\K}(\G,\X)^*$ with coefficients in $V$ is the complex
\[IC^{\K}V=(j_{!*}(\F^{\K}V[d]))[-d].\]
\item[(ii)] \index{weighted cohomology complex}
(cf \cite{M2} 4.1.3) Let $t_1,\dots,t_n\in\Z\cup\{\pm\infty\}$. 
For every $r\in\{1,\dots,n\}$, write $a_r=-t_r+d_r$. Define an additive
triangulated functor
\[W^{\geq t_1,\dots,\geq t_n}:D^b(\Rens_0)\fl D^b_m(M^{\K}(\G,\X)^*,
K_\lambda)\]
in the following way : for every $m\in\Z$, if $V\in\Ob D^b(\Rens_0)$
is such that all $\Ho^iV$, $i\in\Z$, are pure of weight $m$, then
\[W^{\geq t_1,\dots,\geq t_n}V=w_{\leq (-m+d,-m+a_1,\dots,-m+a_n)}Rj_*\F^{\K}V.
\]

\end{itemize}
\end{definition}

Proposition 4.1.5 of \cite{M2} admits the following obvious generalization :

\begin{proposition}\label{prop:IC_est_pondere} Let $t_1,\dots,t_n\in\Z$ be
such that, for every $r\in\{1,\dots,n\}$, $d_r-d\leq t_r\leq 1+d_r
-d$. Then, for every $V\in\Ob\Rens_0$, there is a canonical isomorphism
\[IC^{\K}V\simeq W^{\geq t_1,\dots,\geq t_n}V.\]

\end{proposition}

We now want to calculate the restriction to boundary strata of the weighted
cohomology complexes.
The following theorem is a consequence of propositions 3.3.4 and 3.4.2 of 
\cite{M2}.

\begin{theorem}\label{th:formule_w}
Let $\as=(a_0,\dots,a_n)\in (\Z\cup\{\pm\infty\})^{n+1}$.
Then, for every $L\in\Ob D^b_m(M^{\K}(\G,\X),K_\lambda)$ such that all
perverse cohomology sheaves of $L$ are pure of weight $a_0$, there is an
equality of classes in the Grothendieck group of $D^b_m(M^{\K}(\G,\X)^*,
K_\lambda)$) :
\[[w_{\leq\as}Rj_*L]=\sum_{1\leq n_1<\dots<n_r\leq n}(-1)^r[i_{n_r!}w_{\leq 
a_{n_r}}i_{n_r}^!\dots i_{n_1!}w_{\leq a_{n_1}}i_{n_1}^!j_!L].\]

\end{theorem}

Therefore it is enough to calculate the restriction to boundary strata of
the complexes
$i_{n_r!}w_{\leq a_{n_r}}i_{n_r}^!\dots i_{n_1!}w_{\leq a_{n_1}}
i_{n_1}^!j_!\F^{\K}V$, $1\leq n_1<\dots<n_r\leq n$. The following proposition
generalizes proposition 4.2.3 of \cite{M2} and proposition 5.2.3 of 
\cite{M} :

\begin{proposition}\label{prop:restriction_ponderes_strates} Let
$n_1,\dots,n_r\in\{1,\dots,n\}$ be such that $n_1<\dots<n_r$, $a_1,\dots,a_r\in
\Z\cup\{\pm\infty\}$, $V\in\Ob D^b(\Rens_0)$ and $g\in\G(\Af)$. Write
$\Pa=\Pa_{n_1}\cap\dots\cap\Pa_{n_r}$; remember that, in \ref{points_fixes1},
before proposition \ref{prop:description_chaines}, we constructed a set
$\Ch_P\simeq\Pa(\Q)\QP_{n_r}(\Af)
\sous\G(\Af)/\K$ and a map from this set to the set of boundary strata
of $M^{\K}(\G,\X)^*$ associated to $\Pa_{n_r}$.
For every $i\in\{1,\dots,r\}$, let $w_i:\Gr_m\fl\G_{n_i}$ be the cocharacter
associated to the Shimura datum $(\G_{n_i},\X_{n_i})$ as in
\ref{points_fixes3}; the image of $w_i$ is contained in the center of
$\G_{n_i}$, and $w_i$ can be seen as a cocharacter of $\M_P$.
For every $i\in\{1,\dots,r\}$, write $t_i=-a_i+d_{n_i}$. Let
\[L=i_{n_r,g}^*Ri_{n_r*}w_{>a_r}i_{n_r}^*\dots Ri_{n_1*}w_{>a_1}i_{n_1}^*
Rj_*\F^{\K}V.\]

Then there is a canonical isomorphism
\[L\simeq\bigoplus_C T_{C*}L_C,\]
where the direct sum is over the set of $C=(X_1,\dots,X_r)\in\Ch_P$ that
are sent to the stratum $Im(i_{n_r,g})$, $T_C$ is the obvious morphism
$X_r\fl Im(i_{n_r},g)$ (a finite étale morphism) and $L_C$ is an
$\lambda$-adic complex on $X_r$ such that, if
$h\in\G(\Af)$ is a representative of $C$, there is an isomorphism
\[L_C\simeq\F^{\Hr/\Hr_L}R\Gamma(\Hr_L/\K_N,R\Gamma(Lie(\N_P),V)_{<t_1,\dots,
<t_r}),\]
where $\Hr=h\K h^{-1}\cap\Pa(\Q)\QP_{n_r}(\Af)$, $\Hr_L=h\K h^{-1}\cap\Pa(\Q)
\N_{n_r}(\Af)\cap\Le_{n_r}(\Q)\N_{n_r}(\Af)$, $\K_N=h\K h^{-1}\cap\N_P(\Q)
\N_{n_r}(\Af)$ and, for every
$i\in\{1,\dots,r\}$, the subscript $>t_i$ means that the complex
$R\Gamma(Lie(\N_P),V)$ of representations of $\M_P$ is truncated by
the weights of $w_i(\Gr_m)$ (cf \cite{M2} 4.1.1).

\end{proposition}

Remember that the Levi quotient $\M_P$ is the direct product of $\G_{n_r}$
and of a Levi subgroup $\Le_P$ of $\Le_{n_r}$. Write $\Gamma_L=
\Hr_L/\K_N$ and $X_L=\Le_P(\R)/\K_{L,\infty}\A_{L_P}(\R)^0$, where 
$\K_{L,\infty}$ is a maximal compact subgroup of $\Le_P(\R)$ and $\A_{L_P}$
is, as in \ref{points_fixes2}, the maximal split subtorus of the center
of $\Le_P$; also remember that $q(\Le_P)=\dim(X_L)/2$.
Then $\Gamma_L$ is a neat arithmetic subgroup of $\Le_P(\Q)$, and, for every
$W\in\Ob D^b(Rep_{\Le_P})$,
\[R\Gamma(\Gamma_L,W)=R\Gamma(\Gamma_L\sous X_L,\F^{\Gamma_L}W).\]
Write
\[R\Gamma_c(\Gamma_L,W)=R\Gamma_c(\Gamma_L\sous X_L,\F^{\Gamma_L}W).\]
If $W\in\Ob D^b(Rep_{\M_P})$, then
this complex can be seen as an object of $D^b(Rep_{\G_{n_r}})$, because
it is the dual of $R\Gamma(\Gamma_L,W^*)[\dim(X_L)]$ (where $W^*$ is the
contragredient of $W$).
Define in the same way a complex $R\Gamma_c(\K_L,W)$ for $\K_L$ a neat open
compact subgroup of $\Le_P(\Af)$ and $W\in\Ob D^b(Rep_{\Le_P})$.

\begin{corollary}\label{cor:restriction_ponderes_strates} Write
\[M=i_{n_r,g}^*i_{n_r!}w_{\leq a_r}i_{n_r}^!\dots i_{n_1!}w_{\leq a_1}i_{n_1}^!
j_!\F^{\K}V.\]
Then there is a canonical isomorphism
\[M\simeq \bigoplus_C T_{C*}M_C,\]
where the sum is as in the proposition above and, for every $C=(X_1,\dots,X_r)
\in\Ch_P$ that is sent to the stratum $Im(i_{n_r,g})$, 
$M_C$ is an $\lambda$-adic complex on $X_r$ such that,
if $h$ is a representative of $C$, then
there is an isomorphism (with the notations of the proposition)
\[M_C\simeq\F^{\Hr/\Hr_L}R\Gamma_c(\Hr_L/\K_N,R\Gamma(Lie(\N_P),V)_{\geq t_1,
\dots,\geq t_r})[-\dim(\A_{M_P}/\A_G)].\]

\end{corollary} 

\begin{proof}
Let $V^*$ be the contragredient of $V$. The complex dual to $M$ is :
\[\begin{array}{rcl}D(M) & = & i_{n_r,g}^*Ri_{n_r*}w_{\geq -a_r}i_{n_r}^*
\dots Ri_{n_1*}w_{\geq -a_1}i_{n_1}^*Rj_*D(\F^{\K}V) \\
& = & i_{n_r,g}^*Ri_{n_r*}w_{\geq -a_r}i_{n_r}^*\dots Ri_{n_1*}w_{\geq -a_1}
i_{n_1}^*Rj_*(\F^{\K}V^*[2d](d)) \\
& = & (i_{n_r,g}^*Ri_{n_r*}w_{\geq 2d-a_r}i_{n_r}^*\dots Ri_{n_1*}w_{\geq 2d-
a_1}i_{n_1}^*Rj_*\F^{\K}V^*)[2d](d).\end{array}\]
For every $i\in\{1,\dots,r\}$, let $s_i=-(2d-a_i-1)+d_{n_i}=1-t_i-2(d-d_
{n_i})$. By proposition \ref{prop:restriction_ponderes_strates},
\[D(M)\simeq\bigoplus_C T_{C*}M'_C,\]
with
\[M'_C\simeq\F^{\Hr/\Hr_L}R\Gamma(\Hr_L/\K_N,R\Gamma(Lie(\N_P),V^*)_{<s_1,\dots
,<s_r})[2d](d).\]
Take $M_C=D(M'_C)$. It remains to prove the formula for $M_C$.

Let $m=\dim(\N_P)$. By lemma (10.9) of \cite{GHM}, 
\[R\Gamma(Lie(\N_P),V)_{\geq t_1,\dots,\geq t_r}\simeq R\Hom(R\Gamma(Lie(\N_P),
V^*)_{<s_1,\dots,<s_r},\Ho^{m}(Lie(\N_P),\Q))[-m],\]
and $\Ho^m(Lie(\N_P),\Q)$ is the character
$\gamma\fle\det(\Ad(\gamma),Lie(\N_P))^{-1}$ of $\M_P$ (only the case of groups
$\G$ with anisotropic center is treated in \cite{GHM}, but the general case
is similar).
In particular, $\Hr_L/\K_N$ acts trivially on $\Ho^m(Lie(\N_P),\Q)$, and
the group $w_{r}(\Gr_m)$ acts by the character
$\lambda\fle\lambda^{2(d-d_{n_r})}$ ($w_r$ is defined as in proposition
\ref{prop:restriction_ponderes_strates}).
Hence
\[M_C\simeq\F^{\Hr/\Hr_L}R\Gamma_c(\Hr_L/\K_N,R\Gamma(Lie(\N_P),V)_{\geq t_1,
\dots,\geq t_r})[a],\]
with
\[a=2d_{n_r}+m+2q(\Le_P)-2d=2q(\G_{n_r})+2q(\Le_P)+\dim(\N_P)-2q(\G)
=-\dim(\A_{M_P}/\A_G).\]

\end{proof}

\begin{proofpn}{\ref{prop:restriction_ponderes_strates}}
Let $C=(X_1,\dots,X_r)\in\Ch_P$. Let $I_1$ be the locally closed immersion
$X_1\fl M^{\K}(\G,\X)$ and, for every
$m\in\{1,\dots,r-1\}$, denote by $j_m$ the open immersion $X_m\fl X_m^*$
and by $I_{m+1}$ the locally closed immersion $X_{m+1}\fl X_m^*$ (where $X_m^*$
is the Baily-Borel compactification of $X_m$).
Define a complex $L_C$ on $X_r$ by :
\[L_C=w_{>a_r}I_r^*Rj_{r-1*}w_{>a_{r-1}}I_{r-1}^*\dots w_{>a_1}I_1^*Rj_*\F^{\K}
V.\]

Let us show by induction on $r$ that $L$ is isomorphic to the direct sum of
the $T_{C*}L_C$, for $C\in\Ch_P$ that is sent to the stratum
$Y:=Im(i_{n_r,g})$. The statement is obvious if $r=1$.
Suppose that $r\geq 2$ and that the statement is true for $r-1$.
Let $Y_1,\dots,Y_m$ be the boundary strata of $M^{\K}(\G,\X)^*$
associated to $\Pa_{n_{r-1}}$ whose adherence contains $Y$.
For every $i\in\{1,\dots,m\}$, let $u_i:Y_i\fl M^{\K}(\G,\X)^*$ be the
inclusion and let
\[L_i=u_i^*Ri_{n_{r-1}*}w_{>a_{r-1}}i_{n_{r-1}}^*\dots Ri_{n_1*}w_{>a_1}
i_{n_1}^*Rj_*\F^{\K}V.\]
It is obvious that
\[L=\bigoplus_{i=1}^mi_{n_r,g}^*Ru_{i*}w_{>a_r}L_i.\]
Write $\Pa'=\Pa_{n_1}\cap\dots\cap\Pa_{n_{r-1}}$.
Let $i\in\{1,\dots,m\}$. By the induction hypothesis, $L_i$ is isomorphic to
the direct sum of the $T_{C'*}L_{C'}$ over the set of $C'\in\Ch_{P'}$
that are sent to $Y_i$, where $L_{C'}$ is defined in the same way as
$L_C$.
Fix $C'=(X_1,\dots,X_{r-1})$ that is sent to $Y_i$; let us calculate
$i_{n_r,g}^*Ru_{i*}w_{>a_r}T_{C'*}L_{C'}$.
There is a commutative diagram, with squares cartesien up to nilpotent
elements :
\[\xymatrix{Y'\ar[r]^-{I'}\ar[d]_T & X_{r-1}^*\ar[d]_{\overline{T}_{C'}} & 
X_{r-1}\ar[l]_-{j_{r-1}}\ar[d]_{T_{C'}} \\
Y\ar[r] & \overline{Y}_i & Y_i\ar[l]}\]
where $Y'$ is a disjoint union of boundary strata of $X_{r-1}^*$ associated
to the parabolic subgroup $(\Pa_{n_r}\cap\QP_{n_{r-1}})/\N_{r_{n-1}}$.
Moreover, the vertical arrows are finite maps, and the maps
$T$ and $T_{C'}$ are étales. By the proper base change isomorphism and
the fact the functors $w_{>a}$ commute with taking the direct image by
a finite étale morphism, there is an isomorphism :
\[i_{n_r,g}^*Ru_{i*}w_{>a_r}T_{C'*}L_{C'}=T_*w_{>a_r}{I'}^*Rj_{r-1*}L_{C'}.\]
The right hand side is the direct sum of the complexes
\[(T\circ I_r)_*w_{>a_r}I_r^*Rj_{r-1*}L_{C'}=T_{C*}L_C,\]
for $I_r:X_r\fl X_{r-1}^*$ in the set of boundary strata of $X_{r-1}^*$
included in $Y'$ and for $C=(X_1,\dots,X_r)$. This calculations clearly
imply the statement that we were trying to prove.

It remains to prove the formula for $L_C$ given in the proposition.
Again, use induction on $r$. If $r=1$, the formula for $L_C$ is a direct
consequence of Pink's theorem (\ref{th:Pink}) and of lemma 4.1.2 of \cite{M2}.
Supppose that $r\geq 2$ and that the result is known for $r-1$. Let
$C=(X_1,\dots,X_r)\in\Ch_P$, and let $h\in\G(\Af)$ be a representative of $C$.
Write $\Pa'=\Pa_{r_1}\cap\dots\cap
\Pa_{r_{n-1}}$, $C'=(X_1,\dots,X_{r-1})$,
\[\Hr=h\K h^{-1}\cap\Pa(\Q)\QP_{n_r}(\Af),\]
\[\Hr_L=h\K h^{-1}\cap\Pa(\Q)\N_{n_r}(\Af)\cap\Le_{n_r}(\Q)\N_{n_r}(\Af)=
\Hr\cap\Le_{n_r}(\Q)\N_{n_r}(\Af),\]
\[\K_N=h\K h^{-1}\cap\N_P(\Q)\N_{n_r}(\Af),\] 
\[\Hr'=h\K h^{-1}\cap\Pa'(\Q)\QP_{n_{r-1}}(\Af),\]
\[\Hr'_L=h\K h^{-1}\cap\Pa'(\Q)\N_{n_{r-1}}(\Af)\cap\Le_{n_{r-1}}(\Q)\N_{n_
{r-1}}(\Af)=\Hr'\cap\Le_{n_{r-1}}(\Q)\N_{n_{r-1}}(\Af),\]
\[\K'_N=h\K h^{-1}\cap\N_{P'}(\Q)\N_{n_{r-1}}(\Af).\] 
By the induction hypothesis, there is a canonical isomorphism
\[L_{C'}\simeq\F^{\Hr'/\Hr'_L}R\Gamma(\Hr'_L/\K'_N,R\Gamma(Lie(\N_{P'}),
V)_{<t_1,\dots,<t_{r-1}}).\]
Applying Pink's theorem, we get a canonical isomorphism
\[L_C\simeq w_{>a_r}\F^{\Hr/\Hr_L}R\Gamma(\Hr_L/\Hr'_L,R\Gamma(\Hr'_L/\K'_N,
R\Gamma(Lie(\N_{n_{r-1}}),V)_{<t_1,\dots,<t_{r-1}}))).\]
There are canonical isomorphisms
\[\begin{array}{rcl}R\Gamma(\Hr_L/\Hr'_L,R\Gamma(\Hr'_L/\K'_N,-)) & \simeq & 
R\Gamma(\Hr_L/\K_N,R\Gamma(\K_N/\K'_N,-)) \\
& \simeq & R\Gamma(\Hr_L/\K_N,R\Gamma(Lie(\N_{n_r}/\N_{n_{r-1}}),-))
\end{array}\]
(the last isomorphism comes from van Est's theorem, cf \cite{GKM} \S24).
On the other hand, for every $i\in\{1,\dots,r-1\}$, the image of the
cocharacter $w_i:\Gr_m\fl\G_{n_i}$ is contained in the center of $\G_{n_i}$,
hence it commutes with $\G_{n_{r-1}}$. This implies that
\[R\Gamma(Lie(\N_{n_r}/\N_{n_{r-1}}),R\Gamma(Lie(\N_{n_{r-1}}),V)_{<t_1,\dots,
t_{r-1}})=R\Gamma(Lie(\N_{n_r}),V)_{<t_1,\dots,<t_{r-1}},\]
so that
\[L_C\simeq w_{>a_r}\F^{\Hr/\Hr_L}R\Gamma(\Hr_L/\K_N,R\Gamma(Lie(\N_{n_r}),V)_
{<t_1,\dots,<t_{r-1}}).\]
To finish the proof, it suffices to apply lemma 4.1.2 of \cite{M2}
and to notice that the image of
$w_r:\Gr_m\fl\G_{n_r}$ commutes with $\Le_{n_r}(\Q)$, hence also with its
subgroup $\Hr_L/\K_N$.

\end{proofpn}

 % complexes pondérés et complexe d'intersection
\section{Cohomological correspondences}
\label{points_fixes5}

\begin{notation} Let $(T_1,T_2):X'\fl X_1\times X_2$ be a correspondence
of separated schemes of finite type over a finite field, and let
$c:T_1^*L_1\fl T_2^!L_2$ be a cohomological correspondence with support
in $(T_1,T_2)$. Denote by $\Phi$ the absolute Frobenius morphism of $X_1$.
\index{$\Phi$\quad absolute Frobenius morphism}
For every $j\in\Nat$, we will write $\Phi^j c$ for the cohomological
correspondence with support in $(\Phi^j\circ T_1,T_2)$ defined as the
following composition of maps :
\[(\Phi^j\circ T_1)^*L_1=T_1^*\Phi^{j*}L_1\simeq T_1^*L_1\stackrel{c}{\fl}
T_2^!L_2.\]

\end{notation}

\vspace{.5cm}

First we will define Hecke correspondences on the complexes of
\ref{points_fixes2}.
\index{Hecke correspondence}
Fix $\M$, $\Le$ and $(\G,\X)$ as in
\ref{points_fixes2}. Let $m_1,m_2\in\Le(\Q)\G(\Af)$ and
$\K'_M,\K_M^{(1)},\K_M^{(2)}$ be neat open compact subgroups of
$\M(\Af)$ such that $\Hr'\subset m_1\Hr^{(1)}m_1^{-1}\cap m_2
\Hr^{(2)}m_2^{-1}$, where $\Hr'=\K'_M\cap\Le(\Q)\G(\Af)$ and
$\Hr^{(i)}=\K_M^{(i)}\cap\Le(\Q)\G(\Af)$.
This gives two finite étale morphisms
$T_{m_i}:M(\G,\X)/\Hr'\fl M(\G,\X)/\Hr^{(i)}$, $i=1,2$. Write
$\Hr_L^{(i)}=\Hr^{(i)}\cap\Le(\Q)$ and $\Hr'_L=\Hr'\cap\Le(\Q)$. Let $V\in\Ob 
Rep_\M$. For $i=1,2$, write
\[L_i=\F^{\Hr^{(i)}/\Hr_L^{(i)}}R\Gamma(\Hr_L^{(i)},V).\]
By \cite{P2} 1.11.5, there are canonical isomorphisms
\[T_{m_i}^*L_i\simeq\F^{\Hr'/\Hr'_L}\theta_i^*R\Gamma(\Hr^{(i)}_L,V)\]
where $\theta_i^*R\Gamma(\Hr^{(i)}_L,V)$ is the inverse image by the morphism
$\Hr'/\Hr'_L\fl\Hr^{(i)}/\Hr^{(i)}_L$, $h\fle m_i^{-1}hm_i$, of the complex of
$\Hr^{(i)}/\Hr_L^{(i)}$-modules $R\Gamma(\Hr^{(i)}_L,V)$. 
Using the injections
$\Hr'_L\fl\Hr_L^{(i)}$, $h\fle m_i^{-1}hm_i$, we get an adjonction morphism
$\theta_1^*R\Gamma(\Hr_L^{(1)},V)\stackrel{adj}{\fl}R\Gamma(\Hr'_L,V)$ and a
trace morphism $R\Gamma(\Hr'_L,V)\stackrel{Tr}{\fl}\theta_2^*R\Gamma
(\Hr_L^{(2)},V)$ (this last morphism
exists because the index of $\Hr'_L$ in $\Hr_L^{(2)}$ is finite);
these morphisms are $\Hr'/\Hr'_L$-equivariant. The Hecke correspondence
\[c_{m_1,m_2}:T_{m_1}^*L_1\fl T_{m_2}^!L_2=T_{m_2}^*L_2\]
is the map
\[T_{m_1}^*L_1\simeq\F^{\Hr'/\Hr'_L}\theta_1^*R\Gamma(\Hr^{(1)}_L,V)\stackrel
{adj}{\fl}\F^{\Hr'/\Hr'_L}R\Gamma(\Hr'_L,V)\stackrel{Tr}{\fl}\F^{\Hr'/\Hr'_L}
\theta_2^*R\Gamma(\Hr_L^{(2)},V)\simeq T_{m_2}^*L_2.\]
Note that, if $\Le=\{1\}$, then this correspondence is an isomorphism.

\begin{remarks}\begin{itemize}
\item[(1)] Assume that $\K'_M\subset m_1\K_M^{(1)}m_1^{-1}\cap
m_2\K_M^{(2)}m_2^{-1}$, and write $\K'_L=\K'_M\cap\Le(\Af)$ and
$\K_L^{(i)}=\K_M^{(i)}\cap\Le(\Af)$. Using the methods of
\cite{M} 2.1.4 (and the fact that, for every open compact subgroup $\K_L$
of $\Le(\Af)$, $R\Gamma(\K_L,V)=\bigoplus_{i\in I}R\Gamma
(g_i\K_L g_i^{-1}\cap\Le(\Q),V)$, where $(g_i)_{i\in I}$ is a system of
representatives of $\Le(\Q)\sous\Le(\Af)/\K_L$), it is possible to construct
complexes $M_i=\F^{\K_M^{(i)}/\K_L^{(i)}}R\Gamma(\K_L^{(i)},V)$ and
$\F^{\K'_M/\K'_L}R\Gamma(\K'_L,V)$. There is a correspondence
\[(T_{m_1},T_{m_2}):M^{\K'_M/\K'_L}(\G,\X)\fl M^{\K_M^{(1)}/\K_L^{(1)}}(\G,
\X)\times M^{\K_M^{(2)}/\K_L^{(2)}}(\G,\X),\]
and a cohomological correspondence, constructed as above,
\[c_{m_1,m_2}:T_{m_1}^*M_1\fl T_{m_2}^!M_2.\]

\item[(2)] There are analogous correspondences, constructed by replacing
$R\Gamma(\Hr_L^{(i)},V)$ and $R\Gamma(\Hr'_L,V)$ (resp. $R\Gamma(\K_L^{(i)},V)$
and $R\Gamma(\K_L',V)$) with $R\Gamma_c(\Hr_L^{(i)},V)$ and $R\Gamma_c(\Hr'_L,
V)$ (resp. $R\Gamma_c(\K_L^{(i)},V)$ and $R\Gamma_c(\K'_L,V)$). We will
still use the notation $c_{m_1,m_2}$ for these correspondences.

\end{itemize}
\end{remarks}

\vspace{.5cm}

Use the notations of \ref{points_fixes4}, and fix $g\in\G(\Af)$ and a second
open compact subgroup $\K'$ of $\G(\Af)$, such that
$\K'\subset\K\cap g\K g^{-1}$. Fix prime numbers $p$ and $\ell$ as in
the end of \ref{points_fixes3}.
In particular, it is assumed that $g\in\G(\Af^p)$ and that $\K$ (resp. $\K'$)
is of the form $\K^p\G(\Z_p)$ (resp. ${\K'}^p\G(\Z_p)$), with $\K^p\subset
\G(\Af^p)$ (resp. ${\K'}^p\subset\G(\Af^p)$) and $\G(\Z_p)$ a hyperspecial
maximal compact subgroup of $\G(\Q_p)$.
As in \ref{points_fixes4},
we will use the notations $M^{\K}(\G,\X)$, etc, for the reductions modulo $p$
of the varieties of \ref{points_fixes1}.

Let $\Phi$ be the absolute Frobenius morphism of $M^{\K}(\G,\X)^*$.
For every $V\in\Ob D^b(Rep_\G)$ and $j\in\Z$, let
$u_j:(\Phi^j T_g)^*\F^{\K}V\fl T_1^!\F^{\K}V$ be the cohomological
correspondence $\Phi^j c_{g,1}$ on $\F^{\K}V$ (with support in
$(\Phi^jT_g,T_1)$).
\index{uj@$u_j$}

Let $V\in\Ob D^b(\Rens_0)$.
By \cite{M2} 5.1.2 et 5.1.3 :
\begin{bulletlist}
\item for every $t_1,\dots,t_n\in\Z\cup\{\pm\infty\}$, the
correspondence $u_j$ extends in a unique way to a correspondence
\[\overline{u}_j:(\Phi^j\overline{T}_g)^*W^{\geq t_1,\dots,\geq t_n}V\fl
\overline{T}_1^!W^{\geq t_1,\dots,\geq t_n}V;\]
\item for every $n_1,\dots,n_r\in\{1,\dots,n\}$ such that $n_1<\dots
<n_r$ and every $a_1,\dots,a_r\in\Z\cup\{\pm\infty\}$, the correspondence
$u_j$ gives in a natural way a cohomological correspondence on
$i_{n_r!}w_{\leq a_r}i_{n_r}^!\dots i_{n_1!}w_{\leq a_1}i_{n_1}^!j_!\F^{\K}V$
with support in $(\Phi^j\overline{T}_g,\overline{T}_1)$; write
$i_{n_r!}w_{\leq a_r}i_{n_r}^!\dots i_{n_1!}w_{\leq a_1}i_{n_1}^!j_!u_j$
for this correspondence.

\end{bulletlist}
Moreover, there is an analog of theorem \ref{th:formule_w} for cohomological
correspondences (cf \cite{M2} 5.1.5). The goal of this section is to calculate
the correspondences
$i_{n_r!}w_{\leq a_r}i_{n_r}^!\dots i_{n_1!}w_{\leq a_1}i_{n_1}^!j_!u_j$.

Fix $n_1,\dots,n_r\in\{1,\dots,n\}$ such that $n_1<\dots<n_r$ and
$a_1,\dots,a_r\in\Z\cup\{\pm\infty\}$, and write
\[L=i_{n_r!}w_{\leq a_r}i_{n_r}^!\dots i_{n_1!}w_{\leq a_1}i_{n_1}^!j_!\F^{\K}
V\]
\[u=i_{n_r!}w_{\leq a_r}i_{n_r}^!\dots i_{n_1!}w_{\leq a_1}i_{n_1}^!j_!u_j.\]
Use the notations of corollary
\ref{cor:restriction_ponderes_strates}. By this corollary, there is an
isomorphism
\[L\simeq\bigoplus_{C\in\Ch_P}(i_CT_C)_!L_C,\]
where, for every $C=(X_1,\dots,X_r)\in\Ch_P$, $i_C$ is the inclusion in
$M^{\K}(\G,\X)^*$ of the boundary stratum image of $X_r$ (ie of the stratum
$Im(i_{n_r,h})$, if $h\in\G(\Af)$ is a representative of $C$).
Hence the correspondence $u$ can be seen as a matrix
$(u_{C_1,C_2})_{C_1,C_2\in\Ch_P}$,
and we want to calculate the entries of this matrix.

Let $\Ch'_P$ be the analog of the set $\Ch_P$ obtained when $\K$ is replaced
with $\K'$. The morphisms $\overline{T}_g$, $\overline{T}_1$ define maps
$T_g,T_1:\Ch'_P\fl\Ch_P$, and these maps correspond via the bijections
$\Ch_P\simeq\Pa(\Q)\QP_{n_r}(\Af)\sous\G(\Af)/\K$ and
$\Ch'_P\simeq\Pa(\Q)\QP_{n_r}(\Af)\sous\G(\Af)/\K'$
of proposition \ref{prop:description_chaines} to the maps induced by
$h\fle hg$ and $h\fle h$.

Let $C_1=(X_1^{(1)},\dots,X_r^{(1)}),C_2=(X_1^{(2)},\dots,X_r^
{(2)})\in\Ch_P$, and choose representatives $h_1,h_2\in\G(\Af)$ of $C_1$ and
$C_2$.
Let $C'=(X_1',\dots,X_r')\in\Ch'_P$ be such that $T_g(C')=C_1$
and $T_1(C')=C_2$.
Fix a representative $h'\in\G(\Af)$ of $C'$. There exist
$q_1,q_2\in\Pa(\Q)\QP_{n_r}(\Af)$ such that $q_1h'\in h_1g\K$ and
$q_2h'\in h_2\K$.
Let $\overline{q}_1,\overline{q}_2$ be the images of $q_1,q_2$ in 
$\Le_{n_r}(\Q)\G_{n_r}(\Af)$. The following diagrams are commutative :
\[\xymatrix@C=40pt{X'_r\ar[r]^-{i_{C'}T_{C'}}\ar[d]_{T_{\overline{q}_1}} & 
M^{\K'}(\G,\X)^*\ar[d]_{\overline{T}_g} \\ X_r^{(1)}\ar[r]^-{i_{C_1}T_{C^{
(1)}}} & M^{\K}(\G,\X)^*}\qquad
\xymatrix@C=40pt{X'_r\ar[r]^-{i_{C'}T_{C'}}\ar[d]_{T_{\overline{q}_2}} & 
M^{\K'}(\G,\X)^*\ar[d]_{\overline{T}_1} \\ X_r^{(2)}\ar[r]^-{i_{C_2}T_{C^
{(2)}}} & M^{\K}(\G,\X)^*}\]
By corollary \ref{cor:restriction_ponderes_strates}, there are isomorphisms
\[L_{C_1}\simeq\F^{\Hr^{(1)}/\Hr^{(1)}_L}R\Gamma_c(\Hr^{(1)}_L/
\K^{(1)}_N,R\Gamma(Lie(\N_{n_r}),V)_{\geq t_1,\dots,\geq t_r})[a]\]
and
\[L_{C_2}\simeq\F^{\Hr^{(2)}/\Hr^{(2)}_L}R\Gamma_c(\Hr^{(2)}_L/\K^{(2)}_N,
R\Gamma(Lie(\N_{n_r}),V)_{\geq t_1,\dots,\geq t_r})[a],\]
where $t_1,\dots,t_r$ are defined as in proposition
\ref{prop:restriction_ponderes_strates},
$a=-\dim(\A_{M_P}/\A_G)$, 
$\Hr^{(i)}=h_i\K h_i^{-1}\cap\Pa(\Q)\QP_{n_r}(\Af)$,
$\Hr_L^{(i)}=\Hr^{(i)}\cap\Le_{n_r}(\Q)\N_{n_r}(\Af)$ and
$\K_N^{(i)}=\Hr^{(i)}\cap\N_{n_r}(\Af)$. We get a cohomological correspondence
\[\Phi^jc_{\overline{q}_1,\overline{q}_2}:
(\Phi^jT_{\overline{q}_1})^*L_{C_1}\fl T_{\overline{q}_2}^!L_{C_2}.\]
Define a cohomological correspondace
\[u_{C'}:(\Phi^j\overline{T}_g)^*(i_{C_1}T_{C_1})_!L_{C_1}\fl\overline{T}
_1^!(i_{C_2}T_{C_2})_!L_{C_2}\]
by taking the direct image with compact support
of the previous correspondence by
$(i_{C_1}T_{C_1},i_{C_2}T_{C_2})$ (the direct image of a correspondence by a
proper morphism is defined in SGA 5 III 3.3; the direct image by a locally
closed immersion is defined in \cite{M2} 5.1.1 (following \cite{F} 1.3.1),
and the direct image with
compact support is defined by duality).
Finally, write
\[N_{C'}=[\K_N^{(2)}:h_2\K'h_2^{-1}\cap\N_{n_r}(\Af)].\]

\begin{proposition}\label{prop:restriction_corr_ponderees_strates} The
coefficient $u_{C_1,C_2}$ in the above matrix is equal to
\[\sum_{C'}N_{C'}u_{C'},\]
where the sum is taken over the set of $C'\in\Ch'_P$ such that
$T_g(C')=C_1$ and $T_1(C')=C_2$.

\end{proposition}

This proposition generalizes (the dual version of) theorem 5.2.2 of \cite{M2}
and can be proved exactly in the same way (by induction on $r$, as in the
proof of proposition \ref{prop:restriction_ponderes_strates}). 
The proof of theorem 5.2.2 of \cite{M2} uses proposition 2.2.3 of 
\cite{M2} (via the proof of corollary 5.2.4), but this proposition
is simply a reformulation of proposition 4.8.5 de \cite{P2}, and it is true
as well for the Shimura varieties considered here.

 % correspondances cohomologiques
\section{The fixed point formulas of Kottwitz and Goresky-Kottwitz-MacPherson}
\label{points_fixes6}

In this section, we recall two results about the fixed points of
Hecke correspondences, that will be used in \ref{points_fixes7}.

\begin{theorem}\label{th:points_fixes_Kottwitz}(\cite{K-PSSV} 19.6)
\index{counting point formula of Kottwitz}
Notations are as in \ref{points_fixes5}. Assume that the Shimura datum
$(\G,\X)$ is of the type considered in \cite{K-PSSV} \S5, and that we are not
in case (D) of that article (ie that $\G$ is not an orthogonal group).
Fix an algebraic closure $\Fi$ of $\Fi_p$.
Let $V\in\Ob Rep_\G$.
For every $j\geq 1$, denote by $T(j,g)$ the sum over the set of fixed points
in $M^{\K'}(\G,\X)(\Fi)$ of the correspondence
$(\Phi^j\circ T_g,T_1)$ of the naive local terms (cf \cite{P3} 1.5)
of the cohomological correspondence $u_j$ on $\F^{\K}V$
defined in \ref{points_fixes5}.
Then
\[T(j,g)=\sum_{(\gamma_0;\gamma,\delta)\in C_{\G,j}}c(\gamma_0;\gamma,\delta)
O_\gamma(f^p)TO_\delta(\phi_j^{\G})\Tr(\gamma_0,V).\]

\end{theorem}

Let us explain briefly the notations (see \cite{K-SVLR} \S2 and 3 for more
detailed explanations).

The function $f^p\in C_c^\infty(\G(\Af^p))$ is defined by the formula
\[f^p=\frac{\ungras_{g\K^p}}{\vol({\K'}^p)}.\]
For every $\gamma\in\G(\Af^p)$, write
\[O_\gamma(f^p)=\int_{\G(\Af^p)_\gamma\sous\G(\Af^p)}f^p(x^{-1}\gamma x)
d\overline{x},\]
where $\G(\Af^p)_\gamma$ is the centralizer of $\gamma$ in
$\G(\Af^p)$.
\index{Ogamma@$O_\gamma$\quad orbital integral}

Remember that we fixed an injection $F\subset\overline{\Q}_p$; this determines
a place $\wp$ of $F$ over $p$.
Let $\Q_p^{nr}$ be the maximal unramified extension of $\Q_p$ in
$\overline{\Q}_p$, $L$ be the unramified extension of degree $j$ of $F_\wp$
in $\overline{\Q}_p$,
$r=[L:\Q_p]$, $\varpi_L$ be a uniformizer of $L$ and
$\sigma\in\Gal(\Q_p^{nr}/\Q_p)$ be the element lifting the arithmetic Frobenius
morphism of $\Gal(\Fi/\Fi_p)$. Let $\delta\in\G(L)$. Define the norm
$N\delta$ of $\delta$ by
\[N\delta=\delta\sigma(\delta)\dots\sigma^{r-1}(\delta)\in\G(L).\]
\index{Ndelta@$N\delta$\quad naïve norm}
The $\sigma$-centralizer of $\delta$ in $\G(L)$ is by definition
\[\G(L)^\sigma_\delta=\{x\in\G(L)|x\delta=\delta\sigma(x)\}.\]
\index{$\sigma$-centralizer}
We say that $\delta'\in\G(L)$ is $\sigma$-conjugate to $\delta$ in $\G(L)$
if there exists $x\in\G(L)$ such that $\delta'=x^{-1}\delta\sigma(x)$. 
\index{$\sigma$-conjugacy}

By definition of the reflex field $F$, the conjugacy class of
cocharacters
$h_x\circ\mu_0:\Gr_{m,\C}\fl\G_\C$, $x\in\X$, of \ref{points_fixes1} 
is defined over $F$. Choose an element $\mu$ in this conjugacy class
that factors through a maximal split torus of $\G$ over $\Of_L$
(cf \cite{K-SVLR} \S3 p173), and write
\[\phi_j^{\G}=\ungras_{\G(\Of_L)\mu(\varpi_L^{-1})\G(\Of_L)}\in\Hecke(\G(L),
\G(\Of_L)).\]
\index{$\phi_j^\G$}
($\Hecke(\G(L),\G(\Of_L))$ is the Hecke algebra of functions with compact
support on $\G(L)$ that are bi-invariant by $\G(\Of_L)$.)
For every $\delta\in\G(L)$ and $\phi\in C_c^\infty(\G(L))$, write
\[TO_\delta(\phi)=\int_{\G(L)_\delta^\sigma\sous\G(L)}\phi(y^{-1}\delta
\sigma(y))d\overline{y}.\]
\index{TOdelta@$TO_\delta$\quad twisted orbital integral}
Let $\widehat{\T}$ be a maximal torus of $\widehat{\G}$. The conjugacy class
of cocharacters $h_x\circ\mu_0$, $x\in\X$, corresponds to
a Weyl group orbit of characters of $\widehat{\T}$; denote by
$\mu_1$ the restriction to $Z(\widehat{\G})$ of any of these characters
(this does not depend on the choices).

\index{CGj@$C_{\G,j}$\quad set of Kottwitz triples}
It remains to define the set $C_{\G,j}$ indexing the sum of the theorem
and the coefficients $c(\gamma_0;\gamma,\delta)$. Consider the set of triples
$(\gamma_0;\gamma,\delta)\in\G(\Q)\times\G(\Af^p)\times\G(L)$ satisfying
the following conditions (we will later write (C) for the list of these
conditions) :
\begin{bulletlist}
\item $\gamma_0$ is semi-simple and elliptic in $\G(\R)$ 
(in there exists an elliptic maximal torus $\T$ of $\G_\R$ such that
$\gamma_0\in\T(\R)$).
\item For every place $v\not=p,\infty$ of $\Q$, $\gamma_v$
(the local component of $\gamma$ at $v$) is
$\G(\overline{\Q}_v)$-conjugate to $\gamma_0$.
\item $N\delta$ and $\gamma_0$ are $\G(\overline{\Q}_p)$-conjugate.
\item The image of the $\sigma$-conjugacy class of $\delta$ by the map
$B(\G_{\Q_p})\fl X^*(Z(\widehat{\G})^{\Gal(\overline{\Q}_p/\Q_p)})$ of
\cite{K-SVLR} 6.1 is the restriction of $-\mu_1$ to
$Z(\widehat{\G})^{\Gal(\overline{\Q}_p/\Q_p)}$.

\end{bulletlist}
Two triples $(\gamma_0;\gamma,\delta)$ and $(\gamma_0';\gamma',
\delta')$ are called equivalent if $\gamma_0$ and $\gamma'_0$ are
$\G(\overline{\Q})$-conjugate, $\gamma$ and $\gamma'$ are
$\G(\Af^p)$-conjugate, and $\delta$ and $\delta'$ are $\sigma$-conjugate
in $\G(L)$.

Let $(\gamma_0;\gamma,\delta)$ be a triple satisfying conditions (C).
Let $I_0$ be the centralizer of $\gamma_0$ in $\G$.
\index{I0@$I_0$\quad centralizer of $\gamma_0$}
There is a canonical morphism $Z(\widehat{\G})\fl Z(\widehat{I}_0)$,
and the exact sequence
\[1\fl Z(\widehat{\G})\fl Z(\widehat{I}_0)\fl Z(\widehat{I}_0)/Z(\widehat
{\G})\fl 1\]
induces a morphism
\[\pi_0((Z(\widehat{I}_0)/Z(\widehat{\G}))^{\Gal(\overline{\Q}/\Q)})\fl\Ho^1
(\Q,Z(\widehat{\G})).\]
Denote by $\Ka(I_0/\Q)$ the inverse image by this morphism of the subgroup
\[\Ker^1(\Q,Z(\widehat{\G})):=\Ker(\Ho^1(\Q,Z(\widehat{\G}))\fl\prod_{v\ 
place\ of\ \Q}\Ho^1(\Q_v,Z(\widehat{\G}))).\]
\index{KIQ@$\Ka(I_0/\Q)$}
In \cite{K-SVLR} \S2, Kottwitz defines an element $\alpha(\gamma_0;\gamma,
\delta) \in\Ka(I_0/\Q)^D$ (where, for every group $A$,
$A^D=\Hom(A,\C^\times)$); this element depends only on the equivalence
class of $(\gamma_0;\gamma,\delta)$.
\index{$\alpha(\gamma_0;\gamma,\delta)$}
For every place $v\not=p,\infty$ of $\Q$, denote by $I(v)$ the centralizer
of $\gamma_v$ in $\G_{\Q_v}$; as $\gamma_0$ and $\gamma_v$ are
$\G(\overline{\Q}_v)$-conjugate, the group $I(v)$ is an inner form of $I_0$
over $\Q_v$.
On the other hand, there exists a $\Q_p$-group $I(p)$ such that
$I(p)(\Q_p)=\G(L)_\delta^\sigma$, and this group is an inner form of $I_0$
over $\Q_p$.
\index{Ip@$I(p)$}
There is a similar object for the infinite place : in
the beginning of \cite{K-SVLR} \S3, Kottwitz defines an inner form
$I(\infty)$ of $I_0$; $I(\infty)$ is an algebraic group over $\R$, anisotropic
modulo $\A_G$.
\index{Iinf@$I(\infty)$}
Kottwitz shows that, if $\alpha(\gamma_0;\gamma,\delta)=1$, then there exists
an inner form $I$ of $I_0$ over $\Q$ such that, for every place $v$ of $\Q$,
$I_{\Q_v}$ and $I(v)$ are isomorphic (Kottwitz's statement is more precise,
cf \cite{K-SVLR} p 171-172).
\index{I@$I$\quad ``centralizer'' of a Kottwitz triple}

The set $C_{\G,j}$ indexing the sum of the theorem is the set of
equivalence classes of triples $(\gamma_0;\gamma,\delta)$ satisfying
conditions (C) and such that $\alpha(\gamma_0;\gamma,\delta)=1$.
For every $(\gamma_0,\gamma,\delta)$ in $C_{\G,j}$, let
\[c(\gamma_0;\gamma,\delta)=\vol(I(\Q)\sous I(\Af))|\Ker(\Ker^1(\Q,I_0)\fl
\Ker^1(\Q,\G))|.\]
\index{cgamma@$c(\gamma_0;\gamma,\delta)$}

Finally, the Haar measures are normalized as in \cite{K-SVLR} \S3 :
Take on $\G(\Af^p)$ (resp. $\G(\Q_p)$, resp. $\G(L)$) the Haar measure
such that the volume of $\K^p$ (resp. $\G(\Z_p)$, resp. $\G(\Of_L)$) is equal
to $1$. Take on $I(\Af^p)$ (resp. $I(\Q_p)$) a Haar measure such that the
volume
of every open compact subgroup is a rational number, and use inner twistings
to transport these measures to $\G(\Af^p)_\gamma$ and
$\G(L)_\delta^\sigma$.

\begin{remark} If $\K'=\K\cap g\K g^{-1}$, we may replace
$f^p$ with the function
\[\frac{\ungras_{\K^pg\K^p}}{\vol(\K^p)}\in\Hecke(\G(\Af^p),\K^p):=
C_c^\infty(\K^p\sous\G(\Af^p)/\K^p)\]
(cf \cite{K-PSSV} \S16 p 432).

\end{remark}

\begin{remark} There are two differences between the formula given here and
formula (19.6) of \cite{K-PSSV} :
\begin{itemize}
\item[(1)] Kottwitz considers the correspondence $(T_g,\Phi^j\circ T_1)$
(and not $(\Phi^j\circ T_g,T_1)$) and does not define the naive local term in
the same way as Pink (cf \cite{K-PSSV} \S16 p 433). But is is easy to see
(by comparing the definitions of the naive local terms and composing
Kottwitz's correspondence by $T_{g^{-1}}$) that the number
$T(j,f)$ of \cite{K-PSSV} (19.6) is equal to $T(j,g^{-1})$.
This explains that the function of $C_c^{\infty}(\G(\Af^p))$ appearing in
theorem \ref{th:points_fixes_Kottwitz} is
$\vol({\K'}^p)^{-1}\ungras_{g\K^p}$, instead of the function
$\widetilde{f}^p=\vol({\K'}^p)^{-1}\ungras_{\K^pg^{-1}}$ of \cite{K-PSSV} \S16
p 432. (Kottwitz also takes systematically $\K'=\K\cap g\K g^{-1}$, but
his result generalizes immediately to the case where $\K'$ is of finite index
in $\K\cap g\K g^{-1}$).
\item[(2)] Below formula (19.6) of \cite{K-PSSV}, Kottwitz notes that
this formula is true for the canonical model of a Shimura variety associated
to the datum $(\G,\X,h^{-1})$ (and not $(\G,\X,h)$). The normalization
of the global class field isomorphism used in \cite{K-SVLR}, \cite{K-PSSV}
and here are the same (it is also the normalization of \cite{D} 0.8 and
\cite{P2} 5.5). However, the convention for the action of the Galois
group on the special points of the canonical model that is used here
is the convention of \cite{P2} 5.5, and it differs (by a sign)
from the convention of \cite{D} 2.2.4 (because the reciprocity morphism
of \cite{P2} 5.5 is the inverse of the reciprocity morphism of
\cite{D} 2.2.3). As Kottwitz uses Deligne's conventions, what he calls
canonical model of a Shimura variety associated to the datum
$(\G,\X,h^{-1})$ is what is called here canonical model of a Shimura
variety associated to the datum $(\G,\X,h)$.

\end{itemize}
\end{remark}

\begin{remark}\label{rq:Kottwitz_en_mieux}
Actually, Kottwitz proves a stronger result in
\cite{K-PSSV} \S19 :
For every $\gamma\in\G(\Af^p)$, let $N(\gamma)$ be the number of fixed
points $x'$ in $M^{\K'}(\G,\X)(\Fi)$ that can be represented by an element
$\til{x}$ of $M(\G,\X)(\Fi)$ such that there exists $k\in\K$ and
$g\in\G(\Af)$ with $\Phi^j(\til{x})g=
\til{x}k$ and $gk^{-1}$ $\G(\Af^p)$-conjugate to $\gamma$
(this condition depends only on $x'$, and not on the choice of $\til{x}$).
Then
\[N(\gamma)=\sum_\delta c(\gamma_0;\gamma,\delta)O_\gamma(f^p)TO_\delta(\phi_j
^{\G}),\]
where the sum is taken over the set of $\sigma$-conjugacy classes of
$\delta\in\G(L)$ such that there exists $\gamma_0\in\G(\Q)$ such that
the triple $(\gamma_0;\gamma,\delta)$ is in $C_{\G,j}$ (if such a $\gamma_0$
exists, it is unique up to $\G(\overline{\Q})$-conjugacy, because, for every
place $v\not=p,\infty$ of $\Q$, it is conjugate under
$\G(\overline{\Q}_v)$ to the component at $v$ of $\gamma$).
Moreover, if $x'$ is a fixed point contributing to $N(\gamma)$, then
the naive local term at $x'$ is
$\Tr(\gamma_\ell,V)$ (where $\gamma_\ell$ is the $\ell$-adic component of
$\gamma$).

\end{remark}

\begin{remark}\label{rq:Kottwitz_en_encore_plus_mieux} Some of the
Shimura varieties that will be used later are not of the type considered
in \cite{K-PSSV} \S5, so we will need another generalization of Kottwitz's
result, in a very particular (and easy) case.
Let $(\G,\X,h)$ be a Shimura datum (in the sense of \ref{points_fixes1})
such that $\G$ is a torus.
Let $\Y$ be the image of $\X$ by the morphism $h:\X\fl\Hom(\SD,\G)$
($\Y$ is a point because $\G$ is commutative, but the cardinality of
$\X$ can be greater than $1$ in general; remember that the morphism $h$ is
assumed to have finite fibers, but that it is not assumed to be injective).
Let $\G(\R)^+$ be the subgroup of $\G(\R)$ stabilizing a connected component of
$\X$ (this group does not depend on the choice of the connected component)
and $\G(\Q)^+=\G(\Q)\cap\G(\R)^+$. The results of theorem
\ref{th:points_fixes_Kottwitz} and of remark \ref{rq:Kottwitz_en_mieux}
are true for the Shimura datum $(\G,\Y)$ (in this case, they are a
consequence of the description of the action of the Galois group on the
special points of the canonical model, cf \cite{P2} 5.5).
For the Shimura datum $(\G,\X)$, these results are also true if the
following changes are made :
\begin{itemize}
\item[-] multiply the formula giving the trace in theorem
\ref{th:points_fixes_Kottwitz} and the formula giving the number of fixed
points in remark \ref{rq:Kottwitz_en_mieux} by $|\X|$;
\item[-] replace $C_{\G,j}$ with the subset of triples
$(\gamma_0;\gamma,\delta)\in C_{\G,j}$ such that $\gamma_0\in\G(\Q)^+$.

\end{itemize}
This fact is also an easy consequence of \cite{P2} 5.5.

\end{remark}

\vspace{.5cm}

The fixed point formula of Goresky, Kottwitz and MacPherson applies to
a different situation, that of the end of
\ref{points_fixes2}. Use the notations introduced there.
Let $V\in\Ob Rep_\G$, $g\in\G(\Af)$, and let $\K,\K'$ be neat open compact
subgroups of $\G(\Af)$ such that $\K'\subset\K\cap g\K g^{-1}$. This gives
two finite étale morphisms $T_g,T_1:M^{\K'}(\G,\X)(\C)\fl M^{\K}(\G,\X)(\C)$.
Define a cohomological correspondence
\[u_g:T_g^*\F^{\K}V\iso T_1^!\F^{\K}V\]
as in the beginning of \ref{points_fixes5}.
The following theorem is a particular case of theorem 7.14.B of \cite{GKM}
(cf \cite{GKM} (7.17)).

\begin{theorem}\label{th:points_fixes_GKM}
\index{fixed point formula of Goresky-Kottwitz-MacPherson}
The trace of the correspondence
$u_g$ on the cohomology with compact support
$R\Gamma_c(M^{\K}(\G,\X)(\C),\F^{\K}V)$ is equal to
\[\sum_\M(-1)^{\dim(\A_M/\A_G)}(n_M^G)^{-1}\sum_\gamma\iota^M(\gamma)
^{-1}\chi(\M_\gamma)O_\gamma(f^{\infty}_\M)|D_M^G(\gamma)|^{1/2}Tr(\gamma,V),\]
where the first sum is taken over the set of $\G(\Q)$-conjugacy classes
of cuspidal
Levi subgroups $\M$ of $\G$ and, for every $\M$, the second sum is taken
over the set $\gamma$ of semi-simple $\M(\Q)$-conjugacy classes that
are elliptic in $\M(\R)$.

\end{theorem}

Let us explain the notations.
\begin{bulletlist}
\item $\displaystyle{f^{\infty}=\frac{\ungras_{g\K}}{\vol(\K')}
\in C_c^{\infty}(\G(\Af))}$, and $f^{\infty}_\M$ is the constant term
of $f^{\infty}$ at $\M$ (cf \cite{GKM} (7.13.2)).
\item Let $\M$ be a Levi subgroup of $\G$.
Let $\A_M$ be the maximal ($\Q$-)split subtorus of the center of $\M$ and
\[n_M^G=|\Nor_\G(\M)(\Q)/\M(\Q)|.\]
\index{nMG@$n_M^G$}
$\M$ is called \emph{cuspidal} if the group $\M_\R$ has a maximal
($\R$-)torus $\T$ such that $\T/\A_{M,\R}$ is anisotropic.
\index{cuspidal group}
\item Let $\M$ be a Levi subgroup of $\G$ and $\gamma\in
\M(\Q)$. Let $\M^\gamma$ be the centralizer of $\gamma$ in $\M$,
$\M_\gamma=(\M^\gamma)^0$,
\index{Ggamma@$\G_\gamma$\quad connected centralizer} 
\[\iota^M(\gamma)=|\M^\gamma(\Q)/\M_\gamma(\Q)|\]
and
\[D_M^G(\gamma)=\det(1-\Ad(\gamma),Lie(\G)/Lie(\M)).\]
\index{DMG@$D_M^G$\quad partial Weyl denominator}
\item $\chi(\M_\gamma)$ is the Euler characteristic of
$\M_\gamma$, cf \cite{GKM} (7.10).
\index{$\chi(\G)$\quad Euler characteristic of a group}

\end{bulletlist}

\begin{remark} According to \cite{GKM} 7.14.B, the formula of the theorem
should give $Tr(\gamma,V^*)$ (or $Tr(\gamma^{-1},V)$) and not $Tr(\gamma,V)$.
The difference between the formula given here and that of \cite{GKM}
comes from the fact that \cite{GKM} uses a different convention to
define the trace of $u_g$ (cf \cite{GKM} (7.7)); the convention used here
is that of SGA 5 III and of \cite{P1}.

\end{remark}

 % rappels des formules des points fixes de Kottwitz
                      % et Goresky-Kottwitz-MacPherson
\section{The fixed point formula}
\label{points_fixes7}

Use the notations introduced before proposition
\ref{prop:restriction_corr_ponderees_strates} and the notations of
\ref{points_fixes6}. Assume that the Shimura data
$(\G,\X)$ and $(\G_i,\X_i)$, $1\leq i\leq n-1$, are of the type considered
 \cite{K-PSSV} \S5, with case (D) excluded.
(In particular, $\G^{ad}$ is of abelian type, so we can take $\Rens_0=Rep_\G$,
ie choose any $V\in\Ob D^b(Rep_\G)$.)
Assume moreover that $(\G_n,\X_n)$ is of the type considered in
\cite{K-PSSV} \S5 (case (D) excluded) or that
$\G_n$ is a torus.

We want to calculate the trace of the cohomological correspondence
\[\overline{u}_j:(\Phi^j\overline{T}_g)^*W^{\geq t_1,\dots,\geq t_n}V\fl
\overline{T}_1^!W^{\geq t_1,\dots,\geq t_n}V.\]
Assume that $w(\Gr_m)$ acts on the
$\Ho^iV$, $i\in\Z$, by $t\fle t^m$, for a certain $m\in\Z$
(where $w:\Gr_m\fl\G$ is the cocharacter of \ref{points_fixes3}).

Let
\[f^{\infty,p}=\vol({\K'}^p)^{-1}\ungras_{g\K^p}.\]
Let $\Pa$ be a standard parabolic subgroup of $\G$. Write
$\Pa=\Pa_{n_1}\cap\dots\cap\Pa_{n_r}$, with $n_1<\dots<n_r$. Let
\begin{flushleft}$\displaystyle{T_P=m_P\sum_\Le(-1)
^{\dim(\A_{L}/\A_{L_P})}(n_L^{L_P})^{-1}\sum_{\gamma_L}\iota^L(\gamma_L)
^{-1}\chi(\Le_{\gamma_L})|D_L^{L_P}(\gamma_L)|^{1/2}
}$\end{flushleft}
\begin{center}$\displaystyle{\sum_{(\gamma_0;\gamma,\delta)\in C'_{\G_{n_r},j}}
c(\gamma_0;\gamma,\delta)
O_{\gamma_L\gamma}(f^{\infty,p}_{\Le\G_{n_r}})O_{\gamma_L}(\ungras_{\Le(\Z_p)})
}$\end{center}
\begin{flushright}$\displaystyle{\delta_{P(\Q_p)}^{1/2}(\gamma_0)TO_\delta(
\phi_j^{\G_{n_r}})
\delta_{P(\R)}^{1/2}(\gamma_L\gamma_0)\Tr(\gamma_L\gamma_0,R\Gamma(Lie(\N_P),
V)_{\geq t_{n_1}+m,\dots,\geq t_{n_r}+m}),}$\end{flushright}
where the first sum is taken over the set
of $\Le_P(\Q)$-conjugacy classes of cuspidal Levi subgroups $\Le$ of
$\Le_P$, the second sum is taken over the set of
semi-simple conjugacy classes $\gamma_L\in\Le(\Q)$ that are elliptic
in $\Le(\R)$, and :
\begin{itemize}
\item[-] $\Le(\Z_p)$ is a hyperspecial maximal compact subgroup of
$\Le(\Q_p)$;
\item[-] $m_P=1$ if $n_r<n$ or if $(\G_n,\X_n)$ is of the type considered
in \cite{K-PSSV} \S5, and $m_P=|\X_{n_r}|$ if $n_r=n$ and $\G_{n_r}$ is
a torus;
\item[-] $C'_{\G_{n_r},j}=C_{\G_{n_r},j}$ if $n_r<n$ or if $(\G_n,\X_n)$
is of the type considered in \cite{K-PSSV} \S5, and, if $\G_n$ is a torus,
$C'_{\G_n,j}$ is the subset of $C_{\G_n,j}$ defined in remark
\ref{rq:Kottwitz_en_encore_plus_mieux}.

\end{itemize}
Write also
\[T_G=\sum_{(\gamma_0;\gamma,\delta)\in C_{\G,j}}c(\gamma_0;\gamma,\delta)
O_\gamma(f^{\infty,p})TO_\delta(\phi_j^{\G})\Tr(\gamma_0,V).\]

\begin{theorem}\label{th:points_fixes_moi} If $j$ is positive and
big enough, then
\[\Tr(\overline{u}_j,R\Gamma(M^{\K}(\G,\X)^*_{\Fi},(W^{\geq t_1,\dots,\geq t_n}
V)_\Fi))=T_G+\sum_\Pa T_P,\]
where the sum is taken over the set of standard parabolic subgroups of $\G$.
Moreover, if $g=1$ and $\K=\K'$, then this formula is true for every
$j\in\Nat^*$.

\end{theorem}

\begin{proof} For every $i\in\{1,\dots,n\}$, let $a_i=-t_i-m+\dim(M_i)$.
For every standard parabolic subgroup $\Pa=\Pa_{n_1}\cap\dots\cap
\Pa_{n_r}$, with $n_1<\dots<n_r$, let
\[T'_P=(-1)^r\Tr(i_{n_r!}w_{\leq a_{n_r}}i_{n_r}^!\dots i_{n_1!}w_{\leq
a_{n_1}}i_{n_1}^!\overline{u}_j).\]
Let
\[T'_G=\Tr(\overline{u}_j,R\Gamma(M^{\K}(\G,\X)^*_\Fi,(j_!\F^{\K}V)_\Fi)).\]
Then, by the dual of proposition 5.1.5 of \cite{M2} and by the definition
of $W^{\geq t_1,\dots,\geq t_n}V$,
\[\Tr(\overline{u}_j,R\Gamma(M^{\K}(\G,\X)^*_{\Fi},(W^{\geq t_1,\dots,\geq t_n}
V)_\Fi))=T'_G+\sum_\Pa T'_P,\]
where the sum is taken over the set of standard parabolic subgroups of $\G$.
So we want to show that $T'_G=T_G$ and $T'_P=T_P$. Fix $\Pa\not=\G$ 
(and $n_1,\dots,n_r$). 
It is easy to see that
\[\dim(\A_{M_P}/\A_G)=r.\]
Let $h\in\G(\Af^p)$. Write
\[\K_{N,h}=h\K h^{-1}\cap\N(\Af)\]
\[\K_{P,h}=h\K h^{-1}\cap\Pa(\Af)\]
\[\K_{M,h}=\K_{P,h}/\K_{N,h}\]
\[\K_{L,h}=(h\K h^{-1}\cap\Le_P(\Af)\N_P(\Af))/\K_{N,h}\]
\[\Hr_h=h\K h^{-1}\cap\Pa(\Q)\QP_{n_r}(\Af)\]
\[\Hr_{L,h}=h\K h^{-1}\cap\Le_P(\Q)\N_P(\Af).\]
Define in the same way groups $\K'_{N,h}$, etc, by replacing $\K$ with
$\K'$. If there exists $q\in\Pa(\Q)\QP_{n_r}(\Af)$ such that $qh\K=hg\K$, let
$\overline{q}$ be the image of $q$ in $\M_P(\Af)$, and let
$u_h$ be the cohomological correspondence on
$\F^{\Hr_h/\Hr_{L,h}}R\Gamma_c(\Hr_{L,h},R\Gamma(Lie(\N_{n_r}),V)_
{\geq t_{n_1},\dots,\geq t_{n_r}})[a]$ with support in
$(\Phi^j T_{\overline{q}},T_1)$ equal to
$\Phi^j c_{\overline{q},1}$ (we may assume that $q\in\Pa(\Af^p)$, 
hence that $\overline{q}\in\M_P(\Af^p)$). This correspondence is called
$u_{C'}$ in \ref{points_fixes5}, where $C'$ is the image of $h$ in
$\Ch'_P$. If there is no such $q\in\Pa(\Q)\QP_{n_r}(\Af)$, take
$u_h=0$. Similarly, if there exists $q\in\Pa(\Af)$ such that $qh\K=hg\K$,
let $\overline{q}$ be the image of $q$ in $\M_P(\Af)$, and let
$v_h$ be the cohomological correspondence on
$\F^{\K_{M,h}/\K_{L,h}}R\Gamma_c
(\K_{L,h},R\Gamma(Lie(\N_P),V)_{\geq t_{n_1},\dots,\geq t_{n_r}})[a]$ with
support in $(\Phi^jT_{\overline{q}},T_1)$ equal to $\Phi^j c_
{\overline{q},1}$ (we may assume that $q\in\Pa(\Af^p)$). 
If there is no such $q\in\Pa(\Af)$, take $v_h=0$.
Finally, let $N_h=[\K_{N,h}:\K'_{N,h}]$.

Let $h\in\G(\Af^p)$ be such that there exists $q\in\Pa(\Af)$ with $qh\K=hg\K$.
By proposition \ref{prop:flip} below,
\[\Tr(v_h)=\sum_{h'}\Tr(u_{h'}),\]
where the sum is taken over a system of representatives $h'\in\G(\Af^p)$ of
the double classes in
$\Pa(\Q)\QP_{n_r}(\Af)\sous\G(\Af)/\K'$ that are sent to the class of $h$ in
$\Pa(\Af)\sous\G(\Af)/K'$ (apply proposition \ref{prop:flip}
with $\M=\M_P$,
$\K_M=\K_{M,h}$, $m$ equal to the image of $q$ in $\M_P(\Af)$).
On the other hand, by proposition 
\ref{prop:restriction_corr_ponderees_strates},
\[T'_P=(-1)^r\sum_h N_h \Tr(u_h),\]
where the sum is taken over a system of representatives $h\in\G(\Af^p)$
of the double classes in $\Pa(\Q)\QP_{n_r}(\Af)\sous\G(\Af)/\K'$. Hence
\[T'_P=(-1)^r\sum_h N_h \Tr(v_h),\]
where the sum is taken over a system of representatives $h\in\G(\Af^p)$
of the double classes in $\Pa(\Af)\sous\G(\Af)/\K'$.

Let $h\in\G(\Af^p)$. Assume that there exists $q\in\Pa(\Af^p)$ such that
$qh\K=hg\K$. Let $\overline{q}$ be the image of $q$ in $\M_P(\Af^p)$. Write
$\overline{q}=q_Lq_H$, with $q_L\in\Le_P(\Af^p)$ and $q_H\in\G_{n_r}
(\Af^p)$. Let
\[f^{\infty,p}_{G,h}=\vol(\K'_{M,h}/\K'_{L,h})^{-1}\ungras_{q_H(\K_{M,h}/\K_
{L,h})}.\]
Notice that $\K'_{L,h}\subset\K_{L,h}\cap q_L\K_{L,h} q_L^{-1}$. Let
$u_{q_L}$ be the endomorphism of $R\Gamma_c(\K_{L,h},R\Gamma(Lie(\N_P),V)_{\geq
t_{n_1},\dots,\geq t_{n_r}})$ induced by the cohomological correspondence
$c_{q_L,1}$.

To calulate the trace of $v_h$, we will use Deligne's conjecture, that has
been proved by Pink (cf \cite{P3}) assuming some hypotheses (that are
satisfied here), and in general by
Fujiwara (\cite{F}) and Varshavsky (\cite{V}). This conjecture (that should
now be called theorem) says that, if $j$ is big enough, then the fixed
points of the correspondence between schemes underlying $v_h$ are all
isolated, and that the trace of $v_h$ is the sum over these fixed points
of the naive local terms.
By theorem \ref{th:points_fixes_Kottwitz} and remarks
\ref{rq:Kottwitz_en_mieux} et \ref{rq:Kottwitz_en_encore_plus_mieux},
if $j$ is big enough, then
\begin{flushleft}$\displaystyle{\Tr(v_h)=(-1)^rm_P
\sum_{(\gamma_0;\gamma,\delta)
\in C'_{\G_{n_r},j}}c(\gamma_0;\gamma,\delta)O_\gamma(f^{\infty,p}_{G,h})
TO_\delta(\phi_j^{\G_{n_r}})}$\end{flushleft}
\begin{flushright}$\displaystyle{\Tr(u_{q_L}\gamma_0,R\Gamma_c(\K_{L,h},R\Gamma
(Lie(\N_P),V)_{\geq t_{n_1},\dots,\geq t_{n_r}})).}$\end{flushright}
Let
\[f^\infty_{L_P,h}=\vol(\K'_{L,h})^{-1}\ungras_{q_L\K_{L,h}}.\]
Then
\[f^\infty_{L_P,h}=\ungras_{\Le_P(\Z_p)}f^{\infty,p}_{L_P,h},\]
with $f^{\infty,p}_{L_P,h}\in C_c^\infty(\Le_P(\Af^p))$.
By theorem \ref{th:points_fixes_GKM}, for every $\gamma_0\in\G_{n_r}(\Q)$,
\begin{flushleft}$\displaystyle{\Tr(u_{q_L}\gamma_0,R\Gamma_c(\K_{L,h},R\Gamma
(Lie(\N_P),V)_{\geq t_{n_1},\dots,\geq t_{n_r}}))=\sum_\Le(-1)^{\dim(\A_L/\A_
{L_P})}(n_L^{L_P})^{-1}}$\end{flushleft}
\begin{flushright}$\displaystyle{\sum_{\gamma_L}\iota^L(\gamma_L)^{-1}\chi(\Le_
{\gamma_L})|D_L^{L_P}(\gamma_L)|^{1/2}
O_{\gamma_L}((f^\infty_{L_P,h})_L)\Tr(\gamma_L\gamma_0,R\Gamma(Lie
(\N_P),V)_{\geq t_{n_1},\dots,\geq t_{n_r}}),}$\end{flushright}
where the first sum is taken over the set of conjugacy classes
of cuspidal Levi subgroups $\Le$ of $\Le_P$ and the second sum is taken
over the set of semi-simple conjugacy classes $\gamma_L$ of $\Le(\Q)$ that
are elliptic in $\Le(\R)$.
To show that $T'_P=T_P$, it is enough to show that,
for every Levi subgroup $\Le$ of $\Le_P$, for every $\gamma_L\in\Le(\Q)$
and every $(\gamma_0;\gamma,\delta)\in C_{\G_{n_r},j}$, 
\[\sum_h N_hO_{\gamma_L}((f^\infty_{L_P,h})_L)O_\gamma(f^{\infty,p}_{G,h})=
O_{\gamma_L\gamma}(f^{\infty,p}_{\Le\G_{n_r}})\delta_{P(\Q_p)}^{1/2}(\gamma_L
\gamma_0)O_{\gamma_L}(\ungras_{\Le(\Z_p)})\delta_{P(\R)}^{1/2}(\gamma_0),\]
where the sum is taken over a system of representatives
$h\in\G(\Af^p)$ of the double classes in $\Pa(\Af)\sous\G(\Af)/\K'$
(with $f^\infty_{L_P,h}=0$ and
$f^{\infty,p}_{G,h}=0$ if there is no $q\in\Pa(\Af)$ such that
$qh\K=hg\K$).

Fix a parabolic subgroup $\RP$ of $\Le_P$ with Levi subgroup $\Le$,
and let $\Pa'=\RP\G_{n_r}\N_P$ (a parabolic subgroup of $\G$ with Levi subgroup
$\Le\G_{n_r}$). 
Fix a system of representatives $(h_i)_{i\in I}$ in
$\G(\Af^p)$ of $\Pa(\Af)\sous\G(\Af)/\K'$. For every $i\in I$, fix a system
of representatives
$(m_{ij})_{j\in J_i}$ in $\Le_P(\Af^p)$ of
$\RP(\Af)\sous\Le_P(\Af)/\K'_{L,h_i}$. Then $(m_{ij}h_i)_{i,j}$ is a system
of representatives of $\Pa'(\Af)\sous\G(\Af)/\K'$. By lemma
\ref{lemme:terme_constant_O} below,
\[O_{\gamma_L\gamma}(f^{\infty,p}_{\Le\G_{n_r}})=\delta_{P'(\Af^p)}^{1/2}
(\gamma_L\gamma)\sum_{i,j}r(m_{ij}h_i)O_{\gamma_L\gamma}(f_{P',m_{ij}h_i}),\]
where
\[r(m_{ij}h_i)=[(m_{ij}h_i)\K (m_{ij}h_i)^{-1}\cap\N_{P'}(\Af):
(m_{ij}h_i)\K'(m_{ij}h_i)^{-1}\cap\N_{P'}(\Af)]\]
and $f_{P',m_{ij}h_i}$ is equal to the product of
\[\vol(((m_{ij}h_i)\K'(m_{ij}h_i)^{-1}\cap\Pa'(\Af))/((m_{ij}h_i)\K' (m_{ij}
h_i)^{-1}\cap\N_{P'}(\Af)))^{-1}\]
and of the characteristic function of the image in
$(\Le\G_{n_r})(\Af^p)=\M_{P'}(\Af^p)$ of
$(m_{ij}h_i)g\K(m_{ij}h_i)^{-1}\cap\Pa'(\Af^p)$.
Note that
\[r(m_{ij}h_i)=N_{h_i}r'(m_{ij}),\]
where
\[r'(m_{ij})=[m_{ij}\K_{L,h_i}m_{ij}^{-1}\cap\N_R(\Af):m_{ij}\K'_{L,h_i}
m_{ij}^{-1}\cap\N_R(\Af)],\]
that
\[\delta_{P'(\Af^p)}(\gamma_L\gamma)=\delta_{R(\Af^p)}(\gamma_L)\delta_{P(
\Af^p)}(\gamma_L\gamma),\]
and that
\[f_{P',m_{ij}h_i}=f_{R,m_{ij}}f^{\infty,p}_{G,h_i},\]
where $f_{R,m_{ij}}$ is the product of
\[\vol((m_{ij}\K'_{L,h_i}m_{ij}^{-1}\cap\RP(\Af))/(m_{ij}\K'_{L,h_i}m_{ij}^{-1}
\cap\N_R(\Af)))^{-1}\]
and of the characteristic function of the image in $\Le(\Af^p)=\M_R(\Af^p)$
of $(m_{ij}h_i)g\K(m_{ij}h_i)^{-1}\cap\RP(\Af)\N_P(\Af)$. By applying lemma
\ref{lemme:terme_constant_O} again, we find, for every $i\in I$,
\[\sum_{j\in J_i}r'(m_{ij})O_{\gamma_L}(f_{R,m_{ij}})=\delta_{R(\Af^p)}
^{-1/2}(\gamma_L)O_{\gamma_L}((f^{\infty,p}_{L_P,h_i})_\Le).\]
Finally,
\[\sum_{i\in I}N_{h_i}O_{\gamma_L}((f^\infty_
{L_P,h_i})_\Le)O_\gamma(f^{\infty,p}_{G,h_i})=O_{\gamma_L}
((\ungras_{\Le_P(\Z_p)})_\Le)\sum_{i\in I}N_{h_i}O_{\gamma_L}((f^{\infty,p}_
{L_P,h_i})_\Le)O_\gamma(f^{\infty,p}_{G,h_i})\]
\begin{flushright}$\displaystyle{=O_{\gamma_L}((\ungras_{\Le_P(\Z_p)})_\Le)
\sum_{i\in I}N_{h_i}O_\gamma(f^{\infty,p}_{G,h_i})\delta_{R(\Af^p)}^{1/2}
(\gamma_L)\sum_{j\in J_i}r'(m_{ij})O_{\gamma_L}(f_{R,m_{ij}})}$\end{flushright}
\begin{flushright}$\displaystyle{=O_{\gamma_L}((\ungras_{\Le_P(\Z_p)})_\Le)
\delta_{R(\Af^p)}^{1/2}(\gamma_L)\sum_{i\in I}\sum_{j\in J_i}r(m_{ij}h_i)
O_{\gamma_L\gamma}(f_{P',m_{ij}h_i})}$\end{flushright}
\begin{flushright}$\displaystyle{=O_{\gamma_L}((\ungras_{\Le_P(\Z_p)})_\Le)
\delta_{R(\Af^p)}^{1/2}(\gamma_L)\delta_{P'(\Af^p)}^{-1/2}(\gamma_L\gamma)
O_{\gamma_L\gamma}(f^{\infty,p}_{\Le\G_{n_r}})}$\end{flushright}
\begin{flushright}$\displaystyle{=O_{\gamma_L}((\ungras_{\Le_P(\Z_p)})_\Le)
\delta_{P(\Af^p)}^{-1/2}(\gamma_L\gamma)
O_{\gamma_L\gamma}(f^{\infty,p}_{\Le\G_{n_r}}).}$\end{flushright}
To finish the proof, it suffices to notice that
$(\ungras_{\Le_P(\Z_p)})_\Le=\ungras_{\Le(\Z_p)}$, that
$\delta_{P(\Af^p)}^{-1/2}(\gamma_L\gamma)=\delta_{P(\Af^p)}
^{-1/2}(\gamma_L\gamma_0)$, that, as $\gamma_L\gamma_0\in\M_P(\Q)$, the
product formula gives
\[\delta_{P(\Af^p)}^{-1/2}(\gamma_L\gamma_0)=\delta_{P(\Q_p)}^{1/2}(\gamma_L
\gamma_0)\delta_{P(\R)}^{1/2}(\gamma_L\gamma_0)\]
and that
\[\delta_{P(\Q_p)}(\gamma_L\gamma_0)=\delta_{P(\Q_p)}(\gamma_L)\delta_{P(\Q_p)}
(\gamma_0)=\delta_{P(\Q_p)}(\gamma_0)\]
if $O_{\gamma_L}(\ungras_{\Le(\Z_p)})\not=0$ (because this implies that
$\gamma_L$ is conjugate in $\Le(\Q_p)$ to an element of $\Le(\Z_p)$).

If $j$ is big enough, we can calculate $T'_G$ using theorem
\ref{th:points_fixes_Kottwitz} and Deligne's conjecture.
It is obvious
$T'_G=T_G$.

If $g=1$ and $\K=\K'$, then $\overline{u}_j$ is simply the cohomological
correspondence induced by $\Phi^j$.
In this case, we can calculate the trace of $\overline{u}_j$, for every
$j\in\Nat^*$, using
Grothendieck's trace formula (cf
SGA 4 1/2 [Rapport] 3.2).

\end{proof}

\begin{proposition}\label{prop:flip}
Let $\M$, $\Le$ and $(\G,\X)$ be as in \ref{points_fixes2}. 
Let $m\in\M(\Af)$ and let $\K'_M$, $\K_M$ be neat open compact subgroups
of $\M(\Af)$ such that $\K'_M\subset\K_M\cap m\K_M m
^{-1}$. Let $\K_L=\K_M\cap\Le(\Af)$ and $\K=\K_M/\K_L$.
Consider a system of representatives $(m_i)_{i\in I}$ of the set
of double classes $c\in\Le(\Q)\G(\Af)\sous\M(\Af)/\K'_M$ such that
$cm\K_M=c\K_M$. For every $i\in I$, fix $l_i\in\Le(\Q)$ and 
$g_i\in\G(\Af)$ such that $l_ig_im_i\in m_im\K_M$.
Assume that the Shimura varieties and the morphisms that we get from the
above data have good reduction modulo $p$ (in particular,
$\K_M$ and $\K'_M$ are hyperspecial at $p$, and $m,m_i\in\M(\Af^p)$,
$g_i\in\G(\Af^p)$). Let
$\Fi_q$ be the field of definition of these varieties and $\Fi$ be an algebraic
closure of $\Fi_q$.

For every $i\in I$, let $\Hr_i=m_i\K_M m_i^{-1}\cap
\Le(\Q)\G(\Af)$, $\Hr_{i,L}=\Hr_i\cap\Le(\Q)$ and $\K_i=\Hr_i/\Hr_{i,L}$.
Fixe $V\in\Ob Rep_\G$. Let
\[L=\F^{\K}R\Gamma(\K_L,V)\]
\[L_i=\F^{\K_i}R\Gamma(\Hr_{i,L},V)\]
\[M=\F^{\K}R\Gamma_c(\K_L,V)\]
\[M_i=\F^{\K}R\Gamma_c(\Hr_{i,L},V).\]
Then, for every $\sigma\in\Gal(\Fi/\Fi_q)$,
\begin{itemize}
\item[(1)] $\displaystyle{\sum_{i\in I}\Tr(\sigma c_{l_ig_i,1},R\Gamma(M^{\K_i}
(\G,\X)_\Fi,L_{i,\Fi}))=\Tr(\sigma c_{m,1},R\Gamma(M^{\K}(\G,\X)_\Fi,L_\Fi)).}$
\item[(2)] $\displaystyle{\sum_{i\in I}\Tr(\sigma c_{l_ig_i,1},R\Gamma_c(M^
{\K_i}(\G,\X)_\Fi,M_{i,\Fi}))=\Tr(\sigma c_{m,1},R\Gamma_c(M^{\K}(\G,\X)_\Fi,
M_\Fi)).}$

\end{itemize}
\end{proposition}

\begin{proof} Write $m=lg$, with $l\in\Le(\Af)$ and $g\in\G(\Af)$.
We may assume that $m_i\in\Le(\Af)$, hence $g_i=g$, for every $i\in I$.
Let $\K^0=\Hr_i\cap\G(\Af)=m\K_M m^{-1}\cap\G(\Af)$.

Point (1) implies point (2) by duality.

Let us prove (1). Let $c_m$ be the endomorphism of $R\Gamma(\K_M,V)$ equal to
\[R\Gamma(\K_M,V)\fl R\Gamma(\K'_M,V)\stackrel{Tr}{\fl} R\Gamma(\K_M,V),\]
where the first map is induced by the injection $\K'_M\fl\K_M$, $k\fle
m^{-1}km$, and the second map is the trace morphism associated to the
injection $\K'_M\subset\K_M$. Define in the same way, for every $i\in I$, an 
endomorphism $c_{l_ig_i}$ of $R\Gamma(\Hr_i,V)$. Then
\[R\Gamma(\K_M,V)\simeq\bigoplus_{i\in I}R\Gamma(\Hr_i,V)\]
and $c_m=\bigoplus\limits_{i\in I}c_{l_ig_i}$, so it is enough to show that
this decomposition is $\Gal(\Fi/\Fi_q)$-equivariant.
Let $\sigma\in\Gal(\Fi/\Fi_q)$. Then $\sigma$ induces an endomorphism
of $R\Gamma(\K^0,V)=R\Gamma(M^{\K^0}(\G,\X)_\Fi,\F^{\K^0}V_\Fi)$, that will
still be denoted by $\sigma$, and, by the lemma below, the endomorphism of
$R\Gamma(\K_M,V)$ (resp. $R\Gamma(\Hr_i,V)$) induced by $\sigma$ is
\[R\Gamma(\K_M/(\K_M\cap\Le(\Af)),\sigma)\]
\[(resp. \quad R\Gamma(\Hr_i/(\Hr_i\cap\Le(\Q)),\sigma)).\]
This finishes the proof.

\end{proof}

\begin{lemma} Let $\M$, $\Le$ and $(\G,\X)$ be as in the proposition above.
Let $\K_M$ be a neat open compact subgroup of $\M(\Af)$. Let
$\K_L=\K_M\cap\M(\Af)$, $\K_G=\K_M\cap\G(\Af)$, $\Hr=\K_M\cap
\Le(\Q)\G(\Af)$, $\Hr_L=\K_M\cap\Le(\Q)$, $\K=\K_M/\K_L$ and $\K'=\Hr/\Hr_L$.
Let $V\in\Ob Rep_\M$ and $\sigma\in\Gal(\Fi/\Fi_q)$. 
The element $\sigma$ induces an
endomorphism of $R\Gamma(\K_G,V)=R\Gamma(M^{\K_G}(\G,X)_\Fi,\F^{\K_G}V_\Fi
)$ (resp. $R\Gamma(\K_M,V)=R\Gamma(M^{\K}(\G,\X)_\Fi,\F^{\K}R\Gamma(\K_L,V)_
\Fi)$, resp. $R\Gamma(\Hr,V)=R\Gamma(M^{\K'}(\G,\X)_\Fi,\F^{\K'}R\Gamma(\Hr_L,
V)_\Fi)$), that will be denoted by
$\varphi_0$ (resp. $\varphi$, resp. $\varphi'$).
Then
\[\varphi=R\Gamma(\K_M/\K_G,\varphi_0)\]
and
\[\varphi'=R\Gamma(\Hr/\K_G,\varphi_0).\]

\end{lemma}

\begin{proof} The two equalities are proved in the same way. Let us prove
the first one. Let $Y=M^{\K_G}(\G,\X)$, $X=M^{\K}(\G,\X)$, let $f:Y\fl X$
be the (finite étale) morphism $T_1$ and
$L=\F^{\K}R\Gamma(\K_L,V)$. Then, $f^*L=\F^{\K_G}R\Gamma(\K_L,V)$
by \cite{P1} (1.11.5), and $L$ is canonically a direct factor of
$f_*f^*L$ because $f$ is finite étale, so it is enough to show that
the endomorphism of
\[R\Gamma(Y_\Fi,f^*L)=R\Gamma(\K_G,R\Gamma(\K_L,V))=R\Gamma(\K_L,R\Gamma(
\K_G,V))\]
induced by $\sigma$ is equal to $R\Gamma(\K_L,\varphi_0)$. The complex
$M=\F^{\K_G}V$ on $Y$ is a complex of $\K_L$-sheaves in the sense of
\cite{P2} (1.2), and $R\Gamma(\K_L,M)=f^*L$ by \cite{P2}
(1.9.3). To conclude, apply \cite{P2} (1.6.4).

\end{proof}

The following lemma of \cite{GKM} is used in the proof of theorem
\ref{th:points_fixes_moi}. Let $\G$ be a connected reductive group over $\Q$,
$\M$ a Levi subgroup of $\G$ and $\Pa$ a parabolic subgroup of $\G$ with Levi
subgroup $\M$. Let $\N$ be the unipotent radical of $\Pa$.
If $f\in C_c^{\infty}(\G(\Af))$, the constant term
$f_M\in C_c^\infty(\M(\Af))$ of $f$ at $\M$ is defined in
\cite{GKM} (7.13) (the fonction $f_M$ depends on the choice of $\Pa$, but its
orbital integrals do not depend on that choice).
\index{constant term}
For every $g\in\M(\Af)$, let
\[\delta_{P(\Af)}(g)=|\det(\Ad(g),Lie(\N)\otimes\Af)|_{\Af}.\]
\index{$\delta_P$}

Let $g\in\G(\Af)$ and let $\K',\K$ be open compact subgroups of $\G(\Af)$
such that $\K'\subset g\K g^{-1}$. For every $h\in\G(\Af)$, let
$\K_{M}(h)$ be the image in $\M(\Af)$ of $hg\K h^{-1}\cap\Pa(\Af)$,
\[f_{P,h}=\vol((h\K' h^{-1}\cap\Pa(\Af))/(h\K' h^{-1}\cap\N(\Af)))^{-1}
\ungras_{\K_{M}(h)}\in C_c^{\infty}(\M(\Af)),\]
and
\[r(h)=[h\K h^{-1}\cap\N(\Af):h\K' h^{-1}\cap\N(\Af)].\]
(Note that, if there is no element $q\in\Pa(\Af)$ such that
$qh\K=hg\K$, then $\K_M(h)$ is empty, hence $f_{P,h}=0$.)
Let
\[f=\vol(\K')^{-1}\ungras_{g\K}\]
and
\[f_P=\sum_h r(h)f_{P,h},\]
where the sum is taken over a system of representatives of the
double quotient $\Pa(\Af)\sous\G(\Af)/\K'$. 

\begin{lemma}\label{lemme:terme_constant_O}(\cite{GKM} 7.13.A) 
The functions $f_M$ and $\delta_{P(\Af)}^{1/2}f_P$ have the same orbital
integrals.

\end{lemma}

In \cite{GKM}, the ``$g$'' is on the right of the ``$\K$'' (and not on the
left), and $\delta_{P(\Af)}^{-1/2}$ appears in the formula instead of
$\delta_{P(\Af)}^{1/2}$, but is is easy to see that their proof adapts to
the case considered here. There are obvious variants of this lemma
obtained by replacing $\Af$ with $\Af^p$ or $\Q_p$, where $p$ is a prime
number.

\begin{remark}\label{rq:support_terme_constant} The above lemma implies
in particular that the function $\gamma\fle O_\gamma(f_\M)$ on
$\M(\Af)$ has a support contained in a set of the form
$\bigcup\limits_{m\in\M(\Af)}mX m^{-1}$, where $X$ is a compact subset of
$\M(\Af)$), because the support of $\gamma\fle O_\gamma(f_\M)$
is contained in the union of the conjugates of $\K_M(h)$,
for $h$ in a system of representatives of the finite set
$\Pa(\Af)\sous\G(\Af)/\K'$.
Moreover, if $g=1$, then we may assume that $X$ is a finite union of
compact subgroups of $\M(\Af)$, that are neat of $\K$ is neat
(because the $\K_M(h)$ are subgroups of $\M(\Af)$ in that case).

\end{remark}

 % formule des points fixes

\chapter{The groups}
\label{groupes}

In the next chapters, we will apply the fixed point formula to certain
unitary groups over $\Q$.
The goal of this chapter is to define these unitary groups and their Shimura
data, and to recall the description of their parabolic subgroups and of
their endoscopic groups.

\section{Definition of the groups and of the Shimura data}
\label{groupes1}

For $n\in\Nat^*$, write
\[I=I_n=\left(\begin{array}{ccc}1 & & 0 \\ & \ddots & \\ 0 & & 1\end{array}
\right)\in\GL_n(\Z)\]
and
\[A_n=\left(\begin{array}{ccc}0 & & 1 \\ & \begin{turn}{45}\large\ldots
\end{turn}& \\ 1 & & 0\end{array}\right)\in\GL_n(\Z).\]

Let $E=\Q[\sqrt{-b}]$ ($b\in\Nat^*$ square-free) be an imaginary quadratic
extension of $\Q$.
\index{E@$E$\quad imaginary quadratic extension of $\Q$}
The nontrivial automorphism of $E$ will be denoted
by \raisebox{5pt}{$\overline{\hspace{0.2cm}}$}.
\index{\raisebox{5pt}{$\overline{\hspace{0.2cm}}$}\quad non-trivial
automorphism of $E$}
Fix once and for all an injection
$E\subset\overline{\Q}\subset\C$, and an injection $\overline{\Q}\subset
\overline{\Q}_p$ for every prime number $p$.

Let $n\in\Nat^*$ and let $J\in\GL_n(\Q)$ be a symmetric matrix.
Define an algebraic group $\GU(J)$ over $\Q$ by :
\[\GU(J)(A)=\{g\in\GL_n(E\otimes_\Q A)|g^* Jg=c(g)J, c(g)\in A^\times\},\]
for every $\Q$-algebra $A$ (for $g\in \GL_n(E\otimes_\Q A)$, we write
$g^*={}^t \overline{g}$). The group $\GU(J)$ comes with two morphisms
of algebraic groups over $\Q$ :
\[c:\GU(J)\fl\Gr_m\mbox{ et }\det:\GU(J)\fl R_{E/\Q}\Gr_m.\]
Let $\U(J)=\Ker(c)$ and $\SU(J)=\Ker(c)\cap\Ker(\det)$.
\index{unitary group}
\index{GUJ@$\GU(J)$\quad general unitary group}
\index{UJ@$\U(J)$\quad unitary group}
\index{SUJ@$\SU(J)$\quad special unitary group}
\index{c@$c$\quad multiplier morphism}

The group $\SU(J)$ is the derived group of $\GU(J)$ and $\U(J)$.
The groups $\GU(J)$ and $\U(J)$ are connected reductive, and the group
$\SU(J)$ is semi-simple and simply connected.

Let $p,q\in\Nat$ be such that $n:=p+q\geq 1$. 
Let
\[J=J_{p,q}:=\left(\begin{array}{cc}I_p & 0 \\ 0 & -I_q\end{array}\right),\]
and set $\GU(p,q)=\GU(J)$, $\U(p,q)=\U(J)$ and $\SU(p,q)=\SU(J)$.
If $q=0$, we also write $\GU(p)=\GU(p,q)$, etc.
These groups are quasi-split over $\Q$ if and only if $|p-q|\leq 1$.
The semi-simple $\Q$-rank and the semi-simple $\R$-rank of $\GU(p,q)$ are
both equal to $\min(p,q)$.
\index{GUpq@$\GU(p,q)$\quad general unitary group}
\index{Upq@$\U(p,q)$\quad unitary group}
\index{SUpq@$\SU(p,q)$\quad special unitary group}

Let $n\in\Nat^*$. Let
\[\GU^*(n)=\left\{\begin{array}{ll}\GU(n/2,n/2) & \mbox{ if }n\mbox{ is even}
\\ \GU((n+1)/2,(n-1)/2) & \mbox{ if }n\mbox{ is odd}\end{array}\right..\]
The group $\GU^*(n)$ is the quasi-split inner form of $\GU(J)$,
for every symmetric $J\in\GL_n(\Q)$. Write $\U^*(n)=\Ker(c:\GU^*(n)\fl\Gr_m)$
and $\SU^*(n)=\Ker(\det:\U^*(n)\fl R_{E/\Q}\Gr_m)$.
\index{GUn@$\GU^*(n)$\quad quasi-split general unitary group}
\index{Un@$\U^*(n)$\quad quasi-split unitary group}
\index{SUn@$\SU^*(n)$\quad quasi-split special unitary group}

Finally, let $\GU^*(0)=\GU(0,0)=\Gr_m$, and
$(c:\GU^*(0)\fl\Gr_m)=id$.

Let $n_1,\dots,n_r\in\Nat$ and let $J_1\in\GL_{n_1}(\Q),\dots,J_r\in\GL_{n_r}
(\Q)$ be symmetric matrices.
Write
\[\G(\U(J_1)\times\dots\times\U(J_r))=\{(g_1,\dots,g_r)\in\GU(J_1)\times
\dots\times\GU(J_r)|
c(g_1)=\dots=c(g_r)\}.\]
Similarly, write
\[\G(\U^*(n_1)\times\dots\times\U^*(n_r))=\{(g_1,\dots,g_r)\in\GU^*(n_1)\times
\dots\times\GU^*(n_r)|
c(g_1)=\dots=c(g_r)\}.\]
\index{GUp1q1@$\G(\U(p_1,q_1)\times\dots\times\U(p_r,q_r))$}
\index{GUn1@$\G(\U^*(n_1)\times\dots\times\U^*(n_r))$}
\index{GUp1q1@$\G(\U(J_1)\times\dots\times\U(J_r))$}

\begin{remark}\label{rq:groupes_sur_Z}
If the matrix $J$ is in $\GL_n(\Z)$, then there is an obvious way to extend
$\GU(J)$ to a group scheme $\Gf$ over $\Z$ : for every $\Z$-algebra $A$, set
\[\Gf(A)=\{g\in\GL_n(A\otimes_\Z\Of_E)|g^*Jg=c(g)J,c(g)\in A^\times\}.\]
If $\ell$ is a prime number unramified in $E$, then $\Gf_{\Fi_\ell}$ is a
connected reductive algebraic group over $\Fi_\ell$.

In particular, this construction applies to the groups
$\GU(p,q)$ and $\GU^*(n)$.

\end{remark}

\vspace{1cm}

We now define the Shimura data.
Let as before $\SD=R_{\C/\R}\Gr_m$.

Let $p,q\in\Nat$ be such that $n:=p+q\geq 1$, and let $\G=\GU(p,q)$.
If $p\not=q$ (resp. $p=q$), let $\X$ be the set of $q$-dimensional
subspaces of $\C^n$ on which the Hermitian form
$(v,w)\fle {}^t\overline{v}J_{p,q}w$ is negative definite
(resp. positive or negative definite).
Let $x_0\in\X$ be the subspace of $\C^n$ generated by the $q$ vectors
$e_{n+1-q},\dots,e_n$, where $(e_1,\dots,e_n)$ is the canonical
basis of $\C^n$.

The group $\G(\R)$ acts on $\X$ via the injection $\G(\R)\subset
\GL_n(\R\otimes_\Q E)\simeq\GL_n(\C)$, and this action is transitive.
Define a $\G(\R)$-equivariant morphism $h:\X\fl\Hom(\SD,\G_\R)$ by
\[ h_0=h(x_0)=\left\{\begin{array}{ccl} \SD & \fl & \quad\quad\G_\R \\
                             z & \fle  &  
\left({\begin{array}{cc} zI_p & 0 \\
\\ 0 & \overline{z}I_q\end{array}}\right)\end{array}\right..\]
Then $(\G,\X,h)$ is a Shimura datum in the sense of \cite{P1} 2.1.

The group $\SD(\C)=(\C\otimes_\R\C)^\times$ is isomorphic to $\C^\times\times
\C^\times$ by the morphism $a\otimes 1+b\otimes i\fle (a+ib,a-ib)$. Let
$r:\Gr_{m,\C}\fl\SD_\C$ be the morphism $z\fle (z,1)$, and let $\mu=h_0\circ
r:\Gr_{m,\C}\fl\G_\C$. 

Identify $\G_E$ with a subgroup of $\GL_{n,E}\times\GL_{n,E}$ 
by the isomorphism $(R_{E/\Q}\GL_{n,\Q})_E\simeq\GL_{n,E}\times
\GL_{n,E}$
that sends $X\otimes 1+Y\otimes \sqrt{-b}$ to $(X+\sqrt{-b}Y,X-\sqrt{-b}Y)$.
Then, for every $z\in (R_{E/\Q}\Gr_{m,\Q})$,
\[\mu(z)=\left(\begin{array}{cc}(z,1)I_p & 0 \\ 0 & (1,z) I_q\end{array}
\right).\]

\begin{notation}\label{notation:mu_p}
Let $p'\in\{1,\dots,n\}$. 
Define a cocharacter $\mu_{p'}:\Gr_{m,E}\fl\G_E$ by :
\[\mu_{p'}(z)=
\left(\begin{array}{cc}(z,1)I_{p'} \\ 0 &  (1,z)I_{n-p'}\end{array}\right).\]
\index{$\mu_p$}

\end{notation}

\section{Parabolic subgroups}
\label{groupes2}

Let $\G$ be a connected reductive algebraic group over $\Q$.
Fix a minimal parabolic subgroup $\Pa_0$ of $\G$. Remember that a parabolic
subgroup of $\G$ is called \emph{standard} if it contains $\Pa_0$.
Fix a Levi subgroup $\M_0$ of $\Pa_0$. Then a Levi subgroup $\M$ of $\G$
will be called \emph{standard} if $\M$ is a Levi subgroup of a standard
parabolic subgroup and $\M\supset\M_0$.
\index{standard Levi subgroup}
Any parabolic subgroup of $\G$ is $\G(\Q)$-conjugate to a unique standard
parabolic subgroup, so it is enough to describe the standard parabolic
subgroups.

Let $p,q\in\Nat$ be such that $n:=p+q\geq 1$. We are interested in the
parabolic subgroups of $\GU(p,q)$. As $\GU(p,q)=\GU(q,p)$, we may assume
that $p\geq q$. Then the matrix $J_{p,q}$ is $\GL_n(\Q)$-conjugate to
\[A_{p,q}:=\left(\begin{array}{ccc}0 & 0 & A_q \\ 0 & I_{p-q} & 0 \\ A_q & 0 & 
0\end{array}\right),\]
so $\GU(p,q)$ is isomorphic to the unitary group $\G:=\GU(A_{p,q})$, and it
is enough to describe the parabolic subgroups of $\G$.
A maximal torus of $\G$ is the diagonal torus
\[\T=\left\{\left(\begin{array}{ccc}\lambda_1 &  & 0 \\  & \ddots &  
\\ 0 &  & \lambda_n\end{array}\right),\lambda_1,\dots,\lambda_n\in R_{E/\Q}
\Gr_m,\lambda_1\overline{\lambda}_n=\dots=\lambda_q\overline{\lambda}_{p+1}=
\lambda_{q+1}\overline{\lambda}_{q+1}=\dots=\lambda_p\overline{\lambda}_p\in
\Gr_m\right\}.\]
The maximal split subtorus of $\T$ is
\[\Se=\left\{\lambda\left(\begin{array}{ccccccc}\lambda_1 &  &  &  &  &  
& 0 \\  & \ddots &  &  &  &  &  \\
 &  & \lambda_q &  &  &  &  \\
 &  &  &  I_{p-q} &  &  &  \\
 &  &  &  & \lambda_q^{-1} &  &  \\
 &  &  &  &  & \ddots &  \\
0 &  &  &  &  &  & \lambda_1^{-1} \end{array}\right),\lambda,\lambda_1,
\dots,\lambda_q\in\Gr_m\right\}\]
if $p>q$, and
\[\Se=\left\{\left(\begin{array}{cccccc}\lambda\lambda_1 &  &  &  &  & 0 \\
 & \ddots &  &  &  &  \\
 &  & \lambda\lambda_q &  &  &  \\
 &  &  & \lambda_q^{-1} &  &  \\
 &  &  &  & \ddots &  \\
0 &  &  &  &  & \lambda_1^{-1}\end{array}\right),
\lambda,\lambda_1,\dots,\lambda_q\in\Gr_m\right\}\]
if $p=q$.

A minimal parabolic subgroup of $\G$ containing $\Se$ is
\[\Pa_0=\left\{\left(\begin{array}{ccc}A &  & * \\
 & B &  \\
0 &  & C\end{array}\right),A,C\in R_{E/\Q}\B_q,B\in R_{E/\Q}\GL_{p-q}\right\}
\cap\G,\]
where $\B_q\subset\GL_q$ is the subgroup of upper triangular matrices.

The standard parabolic subgroups of $\G$ are indexed by the subsets of
$\{1,\dots,q\}$ in the following way.

Let $S\subset\{1,\dots,q\}$. Write $S=\{r_1,r_1+r_2,\dots,r_1+\dots+r_m\}$
with $r_1,\dots,r_m\in\Nat^*$, and let $r=r_1+\dots+r_m$. The standard
parabolic subgroup $\Pa_S$ corresponding to $S$ is the intersection of $\G$
and of the group
\[\left(\begin{array}{ccccccc}R_{E/\Q}\GL_{r_1} &  &  &  &  &  & * \\
 & \ddots &  &  &  &  &  \\
 &  & R_{E/\Q}\GL_{r_m} &  &  &  &  \\
 &  &  & \GU(A_{p-r,q-r}) &  &  &  \\
 &  &  &  & R_{E/\Q}\GL_{r_m} &  &  \\
 &  &  &  &  & \ddots &  \\
0 &  &  &  &  &  & R_{E/\Q}\GL_{r_1}\end{array}\right).\]
In particular, the standard maximal parabolic subgroups of $\G$ are the
\[\Pa_r:=\Pa_{\{r\}}=\left(\begin{array}{ccc}R_{E/\Q}\GL_r &  & * \\
 & \GU(A_{p-r,q-r}) &  \\
0 &  & R_{E/\Q}\GL_r\end{array}\right)\cap\G\]
for $r\in\{1,\dots,q\}$, and $\Pa_S=\bigcap\limits_{r\in S}\Pa_r$.
Note that $\Pa_0=\Pa_{\{1,\dots,q\}}$.
\index{PS@$\Pa_S$\quad standard parabolic subgroup}
\index{Pr@$\Pa_r$\quad standard maximal parabolic subgroup}

Let $\N_S$ (or $\N_{P_S}$) be the unipotent radical of $\Pa_S$, $\M_S$
(or $\M_{P_S}$) the obvious Levi
subgroup (of block diagonal matrices) and $\A_{M_S}$ the maximal split
subtorus of the center of $\M_S$.
\index{NS@$\N_S$\quad unipotent radical of $\Pa_S$}
\index{MS@$\M_S$\quad Levi quotient of $\Pa_S$}
Write as before $S=\{r_1,\dots,r_1+\dots+r_m\}$ and $r=r_1+\dots+r_m$. Then
there is an isomorphism
\[\begin{array}{ccl}\M_S & \iso & R_{E/\Q}\GL_{r_1}\times\dots\times
 R_{E/\Q}\GL_{r_m}\times\GU(p-r,q-r) \\
diag(g_1,\dots,g_m,g,h_m,\dots,h_1) & \fle & (c(g)^{-1}g_1,\dots,c(g)^{-1}g_m,
g)\end{array}.\]
The inverse image by this isomorphism of $R_{E/\Q}\GL_{r_1}\times\dots\times
R_{E/\Q}\GL_{r_m}$ is called \emph{linear part} of $\M_S$ and denoted by
$\Le_S$ (or $\Le_{P_S}$). The inverse image of $\GU(p-r,q-r)$ is called
\emph{Hermitian part} of $\M_S$ and denoted by $\G_r$ (or $\G_{P_S}$).
Note that the maximal parabolic subgroups of $\G$ satisfy the condition of
\ref{points_fixes1}.
\index{LS@$\Le_S$\quad linear part of $\M_S$}
\index{Gr@$\G_r$\quad Hermitian part of $\M_S$}

\section{Endoscopic groups}
\label{groupes3}

In this section, we want to study the elliptic endoscopic triples for the
groups $\G$ defined in \ref{groupes1}. It is enough to consider the quasi-split
forms.
We will use the definition of elliptic endoscopic triples and of
isomorphisms of endoscopic triples given in \cite{K-STF:CTT} 7.4 et 7.5.
\index{elliptic endoscopic triple}

Let $n_1,\dots,n_r\in\Nat^*$ and $\G=\G(\U^*(n_1)\times\dots\times
\U^*(n_r))$; here we use the Hermitian forms $A_{p,q}$ of \ref{groupes2}
to define $\G$.
We first calculate the dual group $\widehat{\G}$ of $\G$. As $\G$ splits
over $E$, the action of $\Gal(\overline{\Q}/\Q)$ on
$\widehat{\G}$ factors through $\Gal(E/\Q)$. Let $\tau$ be the nontrivial
element of $\Gal(E/\Q)$.
\index{dual group}
\index{G@$\widehat{\G}$\quad dual group of $\G$}
\index{$\tau$\quad non-trivial element of $\Gal(E/\Q)$}

Let $\varphi$ be the isomorphism from $\G_E\subset\GL_{n_1,E\otimes E}\times
\dots\times\GL_{n_r,E\otimes E}$ to $\Gr_{m,E}\times\GL_{n_1,E}\times\dots
\times\GL_{n_r,E}$ that sends $g=(X_1\otimes 1+Y_1\otimes\sqrt{-b},\dots,
X_r\otimes 1+Y_r\otimes\sqrt{-b})\in\G_E$ to $(c(g),X_1+\sqrt{-b}Y_1,\dots,
X_r+\sqrt{-b}Y_r)$. 
Let $\T$ be the diagonal
torus of $\G$ (a maximal torus of $\G$) and $\B$ be the subgroup of upper
triangular matrices in $\G$ (this is a Borel subgroup of $\G$ because of the
choice of the Hermitian form).
There is a canonical isomorphism
\begin{flushleft}$\displaystyle{\T=\{((\lambda_{1,1},\dots,\lambda_{1,n_1}),
\dots,(\lambda_{r,1},\dots,\lambda_{r,n_r}))\in R_{E/\Q}\Gr_m
^{n_1+\dots+n_r}|}$
\end{flushleft}
\begin{flushright}$\displaystyle{\exists\lambda\in\Gr_m,\forall
i\in\{1,\dots,r\},\forall j\in\{1,\dots,n_i\},\lambda_{i,j}\overline{\lambda}_
{i,n_i+1-j}=\lambda\}.}$\end{flushright}
The restriction of $\varphi$ to $\T_E$ induces an isomorphism
\[\T_E\iso\Gr_{m,E}\times\Gr_{m,E}^{n_1}\times\dots\times\Gr_{m,E}^
{n_r}.\]
For every $i\in\{1,\dots,r\}$ and $j\in\{1,\dots,n_i\}$, let $e_{i,j}$ be
the character of $\T$ defined by
\[e_{i,j}(\varphi^{-1}((\lambda,(\lambda_{1,1},\dots,\lambda_{1,n_1}),\dots,
(\lambda_{r,1},\dots,\lambda_{r,n_r}))))=\lambda_{i,j}.\]
Then the group of characters of $\T$ is
\[X^*(\T)=\Z c\oplus\bigoplus_{i=1}^r\bigoplus_{j=1}^{n_i}\Z e_{i,j},\]
and $\Gal(E/\Q)$ acts on $X^*(\T)$ by
\[\tau(c)=c\]
\[\tau(e_{i,j})=c-e_{i,n_i+1-j}.\]

Hence the dual torus of $\T$ is
\[\widehat{\T}=\C^\times\times (\C^\times)^{n_1}\times\dots\times (\C^\times)
^{n_r},\]
with the action of $\Gal(E/\Q)$ given by
\[\tau((\lambda,(\lambda_{i,j})_{1\leq i\leq r,1\leq j\leq n_i}))=
(\lambda\prod_{i,j}\lambda_{i,j},(\lambda_{i,n_i+1-j}^{-1})_{1\leq i\leq r,
1\leq j\leq n_i}).\]
\index{T@$\widehat{\T}$\quad dual group of $\T$}

The set of roots of $\T$ in $Lie(\G)$ is
\[\Phi=\Phi(\T,\G)=\{e_{i,j}-e_{i,j'},1\leq i\leq r,1\leq j,j'\leq n_i,j
\not=j'\}.\]
The subset of simple roots determined by $\B$ is
\[\Delta=\{\alpha_{i,j}=e_{i,j+1}-e_{i,j},1\leq i\leq r,1\leq j\leq n_i-1\}.\]
The group $\Gal(E/\Q)$ acts on $\Delta$ by :
\[\tau(\alpha_{i,j})=\alpha_{i,n_i-j}.\]

For every $n\in\Nat^*$, let $\Phi_n\in\GL_n(\Z)$ be the matrix with
entries
\[(\Phi_n)_{ij}=(-1)^{i+1}\delta_{i,n+1-j}.\]
\index{$\Phi_n$}
The dual group of $\G$ is
\[\widehat{\G}=\C^\times\times\GL_{n_1}(\C)\times\dots\times\GL_{n_r}(\C),\]
with $\widehat{\T}$ immersed diagonally.
The action of $\Gal(E/\Q)$ that respects the obvious splitting is :
\[\tau((\lambda,g_1,\dots,g_r))=(\lambda\det(g_1)\dots\det(g_r),\Phi_{n_1}^{-1}
({}^tg_1)^{-1}\Phi_{n_1},\dots,\Phi_{n_r}^{-1}({}^tg_r)^{-1}\Phi_{n_r}).\]

\begin{proposition}\label{prop:groupes_endoscopiques}
For every $i\in\{1,\dots,r\}$, let $n_i^+,n_i^-\in\Nat$ be
such that $n_i=n_i^++n_i^-$. Suppose that $n_1^-+\dots+n_r^-$ is even.
Set
\[s=(1,diag(\overbrace{1,\dots,1}^{n_1^+},\overbrace{-1,\dots,-1}^{n_1^-}),
\dots,diag(\overbrace{1,\dots,1}^{n_r^+},\overbrace{-1,\dots,-1}^{n_r^-}))
\in\widehat{\G}\]
\[\H=\G(\U^*(n_1^+)\times\U^*(n_1^-)\times\dots\times\U^*(n_r^+)\times
\U^*(n_r^-))\]
and define
\begin{flushleft}$\displaystyle{\eta_0:\widehat{\H}=\C^\times\times\GL_{n_1^+}
(\C)
\times\GL_{n_1^-}(\C)\times\dots\times\GL_{n_r^+}(\C)\times\GL_{n_r^-}(\C)}$
\end{flushleft}
\begin{flushright}$\displaystyle{\fl\widehat{\G}=\C^\times\times\GL_{n_1}
(\C)\times\dots\times\GL_{n_r}(\C)}$\end{flushright}
by
\[\eta_0((\lambda,g_1^+,g_1^-,\dots,g_r^+,g_r^-))=
\left(\lambda,\left(\begin{array}{cc}g_1^+ & 0 \\ 0 & g_1^- \end{array}\right),
\dots,\left(\begin{array}{cc}g_r^+ & 0 \\ 0 & g_r^-\end{array}\right)\right).\]

Then $(\H,s,\eta_0)$ is an elliptic endoscopic triple for $\G$.
The group $\Lambda(\H,s,\eta_0)$ of \cite{K-STF:CTT} 7.5 is isomorphic to
$(\Z/2\Z)^I$, where $I=\{i\in\{1,\dots,r\}|n_i^+=n_i^-\}$. 
\index{$\Lambda(\H,s,\eta_0)$}

Moreover, the elliptic endoscopic triples for $\G$ determined by
$((n_1^+,n_1^-),\dots,(n_r^+,n_r^-))$ and $((m_1^+,m_1^-),\dots,(m_r^+,m_r^-))$
are isomorphic of and only if, for every $i\in\{1,\dots,r\}$,
$(n_i^+,n_i^-)=(m_i^+,m_i^-)$ or $(n_i^+,n_i^-)=(m_i^-,m_i^+)$.

Finally, every elliptic endoscopic triple for $\G$ is isomorphic to one
of the triples defined above.

\end{proposition}

Note that an elliptic endoscopic triple $(\H,s,\eta_0)$ is uniquely
determined by $s$ and that, for every elliptic endoscopic triple
$(\H,s,\eta_0)$, the group $\H_\R$ has an elliptic maximal torus.

\begin{proof} Let $(\H,s,\eta_0)$ be determined by $((n_1^+,n_1^-),\dots,
(n_r^+,n_r^-))$ as above. To show that $(\H,s,\eta_0)$ is an endoscopic triple
of $\G$, we have to check conditions (7.4.1)-(7.4.3) of
\cite{K-STF:CTT}. Conditions (7.4.1) and (7.4.2) are obviously satisfied, and
condition (7.4.3) is a consequence of the fact that
$s\in Z(\widehat{\H})^{\Gal(E/\Q)}$.
(Note that the condition ``$n_1^-+\dots+n_r^-$ even'' is necessary
for $s\in Z(\widehat{\H})$ to be fixed by $\Gal(E/\Q)$.)

We next show that $(\H,s,\eta_0)$ is elliptic.
The center of $\widehat{\H}$ is
\[Z(\widehat{\H})=\{(\lambda,\lambda_1^+I_{n_1^+},\lambda_1^-I_{n_1^-},\dots,
\lambda_r^+I_{n_r^+},\lambda_r^-I_{n_r^-}),\lambda,\lambda_1^+,\lambda_1^-,
\dots,\lambda_r^+,\lambda_r^-\in\C^\times\},\]
with the action of $\Gal(E/\Q)$ given by
\begin{flushleft}$\displaystyle{\tau((\lambda,\lambda_1^+I_{n_1^+},\lambda_1^-
I_{n_1^-},\dots,\lambda_r^+I_{n_r}^+,\lambda_r^-I_{n_r^-}))=}$\end{flushleft}
\begin{flushright}$\displaystyle{(\lambda(\lambda_1^+)^{n_1^+}(\lambda_1^-)
^{n_1^-}\dots (\lambda_r^+)^{n_r^+}(\lambda_r^-)^{n_r^-},(\lambda_1^+)^{-1}
I_{n_1^+},(\lambda_1^-)^{-1}I_{n_1^-},\dots,(\lambda_r^+)^{-1}I_{n_r^+},
(\lambda_r^-)^{-1}I_{n_r^-}).}$\end{flushright}
Hence $\left(Z(\widehat{\H})^{\Gal(E/\Q)}\right)^0=\C^\times\times\{1\}
\subset Z(\widehat{\G})$, and $(\H,s,\eta_0)$ is elliptic.

We want to calculate the group of outer automorphisms of $(\H,s,\eta_0)$.
It is the same to calculate the group of outer automorphisms of the endoscopic
data $(s,\rho)$ associated to $(\H,s,\eta_0)$ (cf \cite{K-STF:CTT} 7.2
and 7.6). Let
\[I=\{i\in\{1,\dots,r\}|n_i^+=n_i^-\}.\]
Let $g\in\widehat{\G}$ be such that $\Int(g)(\eta_0(\widehat{\H}))=\eta_0(
\widehat{\H})$. 

Let $a,b\in\Nat$ be such that $a+b=n>0$, and 
\[\G'=\left(\begin{array}{cc}\GL_a & 0 \\ 0 & \GL_b\end{array}\right)\subset
\GL_n.\]
If $a\not=b$, then the normalizer of $\G'$ in $\GL_n$ is $\G'$. If $a=b$,
then the normalizer of $\G'$ in $\GL_n$ is generated by $\G'$ and by
\[I_{a,b}:=\left(\begin{array}{cc}0 & I_a \\
I_b & 0\end{array}\right).\]

By applying this remark to $\eta_0(\widehat{\H})\subset\widehat{\G}$, we
find that $g$ is in the subgroup of $\widehat{\G}$ generated by
$\eta_0(\widehat{\H})$ and by the elements $(1,\dots,1,I_{n_i^+,n_i^-},1,
\dots,1)$, $i\in I$. It is easy to see that all the elements of this group
define automorphisms of $(s,\rho)$. Hence
\[\Lambda(\H,s,\eta_0)=\Lambda(s,\rho)=\Aut(s,\rho)/\Int(\widehat{\H})
\simeq (\Z/2\Z)^I.\]

The statement about isomorphisms between the endoscopic triples defined
in the proposition is obvious.

Let $(\H,s,\eta_0)$ be an elliptic endoscopic triple for $\G$.
We want to show that $(\H,s,\eta_0)$ is isomorphic to one of the triples
defined above. We may assume (without changing the isomorphism class
of $(\H,s,\eta_0)$) that $s\in\widehat{\T}$. We know that
$\Ker^1(\Q,\G)=\{1\}$ by lemma
\ref{lemme:Tamagawa} below, so condition (7.4.3) of \cite{K-STF:CTT} implies
that the image of $s$ in
$\pi_0\left((Z(\widehat{\H})/Z(\widehat{\G}))^{\Gal(\overline{\Q}/\Q)}\right)$ 
comes from an element of $Z(\widehat{\H})^{\Gal(\overline{\Q}/\Q)}$. As
$(\H,s,\eta_0)$ is elliptic,
\[\pi_0\left((Z(\widehat{\H})/Z(\widehat{\G}))^{\Gal(\overline{\Q}/\Q)}\right)=
(Z(\widehat{\H})/Z(\widehat{\G}))^{\Gal(\overline{\Q}/\Q)},\]
so the image of $s$ in $Z(\widehat{\H})/Z(\widehat{\G})$ comes from an
element of $Z(\widehat{\H})^{\Gal(\overline{\Q}/\Q)}$. We may assume
that is $s$ is fixed by $\Gal(\overline{\Q}/\Q)$ (because replacing $s$ by
a $Z(\widehat{\G})$-translate does not change the isomorphism class of
$(\H,s,\eta_0)$).

Let us first suppose that $r=1$. Write $n=n_1$. We may assume that
\[s=(1,\left(\begin{array}{ccc}\lambda_1 I_{m_1} & 0 & 0 \\
0 & \ddots & 0 \\ 0 & 0 & \lambda_t I_{m_t}\end{array}\right)),\]
with $\lambda_1,\dots,\lambda_t\in\C^\times$, $\lambda_i\not=\lambda_j$
if $i\not =j$ and $m_1,\dots,m_t\in\Nat^*$ such that $m_1+\dots m_t=n$.
Then $\widehat{\H}=\Cent_{\widehat{\G}}(s)\simeq\C^\times\times\GL_{m_1}(\C)
\times\dots\times\GL_{m_t}(\C)$ and $Z(\widehat{\H})\simeq\C^\times\times
(\C^\times)^t$. 

As $(\H,s,\eta_0)$ is elliptic, we must have $\left(Z(\widehat{\H})^{
\Gal(\overline{\Q}/\Q)}\right)^0\subset Z(\widehat{\G})^{\Gal(
\overline{\Q}/\Q)}\subset\C^\times\times\{\pm I_n\}$.
The only way $Z(\widehat{\H})^{\Gal(\overline{\Q}/
\Q)}/Z(\widehat{\G})^{\Gal(\overline{\Q}/\Q)}$ can be finite is if
\[Z(\widehat{\H})^{\Gal(\overline{\Q}/\Q)}\subset\C^\times\times\{\pm 1\}^t.\]
But $s=(1,\lambda_1,\dots,\lambda_t)\in Z(\widehat{\H})^{\Gal(
\overline{\Q}/\Q)}$ and the $\lambda_i$ are pairwise distinct,
so $t\leq 2$. If $t=1$, then $s\in Z(\widehat{\G})$ and
$(\H,s,\eta_0)$ is isomorphic to the trivial endoscopic triple $(\G,1,id)$. 

Suppose that $t=2$. We may assume that $\lambda_1=1$ and $\lambda_2=-1$. 
By condition (7.1.1) of \cite{K-STF:CTT},
\[\tau((\lambda,\lambda_1,\lambda_2))=(\lambda\lambda_1^{m_1}
\lambda_2^{m_2},\lambda_{w(1)}^{-1},\lambda_{w(2)}^{-1}),\]
for a permutation $w\in{\mathfrak S}_2$.
In particular, $(-1)^{m_2}=1$, so $m_2$ is even.

It remains to determine the morphism
$\rho:\Gal(\overline{\Q}/\Q)\fl\Out(\widehat{\H})$ associated to
$(\H,s,\eta_0)$ in \cite{K-STF:CTT} 7.6.
As the derived group of $\widehat{\G}$ is simply connected and $\G$
splits over $E$, $\H$ also splits over $E$ (cf definition 1.8.1
in \cite{Ng}).
So the action of $\Gal(\overline{\Q}/\Q)$ on $\widehat{\H}$ factors
through $\Gal(E/\Q)$, and in particular $\rho$ factors through $\Gal(E/\Q)$.
By condition (7.4.2) of \cite{K-STF:CTT}, there exists
$g_\tau\in\widehat{\G}$ such that $(g_\tau,\tau)$ normalizes
$\widehat{\H}$ in $\widehat{\G}\rtimes\Gal(\overline{\Q}/\Q)$ and that
$\rho(\tau)=\Int((g_\tau,\tau))$ in $\Out(\widehat{\H})$.
Hence $g_\tau=gw_0$, with $g\in\Nor_{\widehat{\G}}
(\widehat{\H})$ and $w_0=\left(\begin{array}{cc}0 & I_{m_1} \\ I_{m_2} & 0
\end{array}\right)$.

Suppose first that $m_1\not=m_2$. Then $\Nor_{\widehat{\G}}
(\widehat{\H})=\widehat{\H}$, so $\rho(\tau)=\Int((w_0,\tau))$.
It is now clear that $(\H,s,\eta_0)$ is isomorphic to one of the triples
defined above.

Suppose that $m_1=m_2$. Then $\Nor_{\widehat{\G}}(\widehat{\H})$
is the subgroup of $\widehat{\G}$ generated by $\widehat{\H}$ and $w_0$.
Hence $\rho(\tau)=\Int((1,\tau))$ or $\Int((w_0,\tau))$.
If $\rho(\tau)=\Int((1,\tau))$, then
\[Z(\widehat{\H})^{\Gal(E/\Q)}\simeq\{(\lambda,\lambda_1,\lambda_2)\in 
(\C^\times)^3|(\lambda_1\lambda_2)^{m_1}=1\mbox{ and }\lambda_1=\lambda_2^
{-1}\}=\{(\lambda,\lambda_1,\lambda_1^{-1}),\lambda,\lambda_1\in\C^\times\},\]
and $s$ is not in the image of $Z(\widehat{\H})^{\Gal(E/\Q)}$. Hence
$\rho(\tau)=\Int((w_0,\tau))$, and $(\H,s,\eta_0)$ is isomorphic to one of
the triples defined above.

If $r>1$, the reasoning is the same (but with more complicated notations).

\end{proof}

Fix
$n_1^+,n_1^-,\dots,n_r^+,n_r^-\in\Nat$ such that $n_i^++n_i^-=n_i$ for every
$i\in\{1,\dots,r\}$ and that $n_1^-+\dots+n_r^-$ is even.
Let $(\H,s,\eta_0)$ be the elliptic endoscopic triple for $\G$ associated
to this data as in proposition \ref{prop:groupes_endoscopiques}.
The derived group of $\G$ is simply connected, so, by proposition 1
of \cite{L2}, there exists a $L$-morphism
$\eta:{}^L\H:=\widehat{\H}\rtimes W_\Q\fl{}^L\G:=\widehat{\G}\rtimes W_\Q$
extending $\eta_0:\widehat{\H}\fl\widehat{\G}$. We want to give an explicit
formula for such a $\eta$.

For every place $v$ of $\Q$, we fixed an injection $\overline{\Q}\subset
\overline{\Q}_v$; this gives a morphism $\Gal(\overline{\Q}_v/
\Q_v)\fl\Gal(\overline{\Q}/\Q)$, and we fix a morphism $W_{\Q_v}\fl W_\Q$
above this morphism of Galois groups.

Let $\omega_{E/\Q}:\Ade^\times/\Q^\times\fl\{\pm 1\}$ be the quadratic
character of $E/\Q$. (Note that, for every prime number $p$ unramified in $E$,
the character $\omega_{E/\Q}$ is unramified at $p$.)

The following proposition is the adaptation to unitary similitude groups of
\cite{Ro2} 1.2 and is easy to prove.

\begin{proposition}\label{prop:prolongement_eta_0} Let $\mu:W_E\fl
\C^\times$ be the character corresponding by the class field isomorphism
$W_E^{ab}\simeq\Ade_E^\times/E^\times$ to a character extending
$\omega_{E/\Q}$. We may, and will, assume that $\mu$ is unitary.
\footnote{In fact, in this case, we may even assume that $\mu$ is of finite
order, but we will not need this fact.}
Let $c\in W_\Q$ be an element lifting the nontrivial element of
$\Gal(E/\Q)$. Define a morphism
$\varphi:W_\Q\fl{}^L\G$ in the following way :
\begin{bulletlist}
\item $\varphi(c)=(A,c)$, where
\[A=(1,\left(\left(\begin{array}{cc}\Phi_{n_1^+} & 0 \\ 0 & (-1)^{n_1^+}
\Phi_{n_1^-}
\end{array}\right)\Phi_{n_1}^{-1},\dots,\left(\begin{array}{cc}\Phi_{n_r^+}
& 0 \\
0 & (-1)^{n_r^+}\Phi_{n_r^-}\end{array}\right)\Phi_{n_r}^{-1}\right));\]
\item on $W_E$, $\varphi$ is given by
\[\varphi_{|W_E}=(1,\left(\left(\begin{array}{cc}\mu^{n_1^-} I_{n_1^+} 
& 0 \\
0 & \mu^{-n_1^+} I_{n_1^-}\end{array}\right),\dots,\left(\begin{array}{cc}
\mu^{n_r^-} I_{n_r^+} & 0 \\ 0 & \mu^{-n_r^+} I_{n_r^-}\end{array}\right)
\right),id).\]

\end{bulletlist}
Then $\varphi$ is well-defined, and $\eta:{}^L\H\fl{}^L\G$, $(h,w)\fle(\eta_0
(h),1)\varphi(w)$, is a $L$-morphism extending $\eta_0$.

\end{proposition}

For every place $v$ of $\Q$, let $\varphi_v$ be the composition of $\varphi$
and of the morphism $W_{\Q_v}\fl W_\Q$. We have the following consequences
of the properties of $\varphi$ in the proposition :

Let $p$ be a prime number unramified in $E$, and fix $\sigma\in W_{\Q_p}$
lifting the arithmetic Frobenius. Set $r=1$ if $p$ splits totally in $E$,
and $r=2$ if $p$ is inert in $E$. Then
\quash{$\varphi_p$ is unramified (so $\eta$
is unramified at $p$), and}
\[\varphi_p(\sigma^r)=(1,(I_{n_1},\dots,I_{n_r}),\sigma^r).\]
\quash{
\[B=\left(\left(\begin{array}{cc}(-1)^{n_1^-}I_{n_1^+} &
0 \\
0 & (-1)^{n_1^+} I_{n_1^-}\end{array}\right),\dots,\left(\begin{array}{cc}
(-1)^{n_r^-} I_{n_r^+} & 0 \\ 0 & (-1)^{n_r^+} I_{n_r^-}\end{array}\right)
\right).\]}

On the other hand, there exists an odd integer $C\in\Z$ such that,
for every $z\in\C^\times=W_\C$,
\[\varphi_\infty(z)=((1,(B_1(z),\dots,B_r(z))),z),\]
with
\[B_i(z)=\left(\begin{array}{cc}z^{Cn_i^-/2}\overline{z}^{-Cn_i^-/2}I_{n_i^+} 
& 0 \\ 0 & z^{-Cn_i^+/2}\overline{z}^{Cn_i^+/2}I_{n_i^-}\end{array}\right).\]

\vspace{.5cm}

We finish this section by a calculation of Tamagawa numbers.
\index{$\tau(\G)$\quad Tamagawa number of $\G$}

\begin{lemma}\label{lemme:Tamagawa}\begin{itemize}
\item[(i)] Let $n_1,\dots,n_r\in\Nat^*$ and $\G=\G(\U^*(n_1)\times
\dots\times\U^*(n_r))$. Then $\Ker^1(\Q,\G)=\{1\}$, and 
$Z(\widehat{\G})^{\Gal(E/\Q)}\simeq\C^\times\times\{(\epsilon_1,\dots,
\epsilon_r)\in\{\pm 1\}^r|\epsilon_1^{n_1}\dots\epsilon_r^{n_r}=1\}$.
Hence the Tamagawa number of $\G$ is
\[\tau(\G)=\left\{\begin{array}{ll}2^r & \mbox{ if all the }n_i\mbox{ are
even} \\
2^{r-1} & \mbox{ otherwise}\end{array}\right..\]
\item[(ii)] Let $F$ be a finite extension of $\Q$ and
$\Le=R_{F/\Q}\GL_{n,F}$,
with $n\in\Nat^*$. Then $\tau(\Le)=1$.

\end{itemize}
\end{lemma}

\begin{proof}
Remember that, by \cite{K-STF:CTT} 4.2.2 and 5.1.1, \cite{K-TN} and
\cite{C},
for every connected reductive algebraic group $\G$ on $\Q$,
\[\tau(\G)=|\pi_0(Z(\widehat{\G})^{\Gal(\overline{\Q}/\Q)})|.|\Ker^1(\Q,\G)|
^{-1}.\]

\begin{itemize}
\item[(i)] It is enough to prove the first two statements.
By \cite{K-PSSV} \S7, the canonical morphism
\[\Ker^1(\Q,Z(\G))\fl\Ker^1(\Q,\G)\]
is an isomorphism. The center of $\G$ is $\{(\lambda_1,\dots,\lambda_r
)\in (R_{E/\Q}\Gr_m)^r|\lambda_1\overline{\lambda}_1=\dots=\lambda_r\overline{
\lambda}_r\}$, so it is isomorphic to $R_{E/\Q}\Gr_m\times\U(1)^{r-1}$
(by the map $(\lambda_1,\dots,\lambda_r)\fle(\lambda_1,\lambda_2\lambda_1^{-1},
\dots,\lambda_r\lambda_1^{-1})$).
As
\[\Ho^1(\Q,R_{E/\Q}\Gr_m)=\Ho^1(E,\Gr_m)=\{1\},\]
it remains to show that $\Ker^1(\Q,\U(1))=\{1\}$. Let
$c:\Gal(\overline{\Q}/\Q)\fl\U(1)(\overline{\Q})$ be a $1$-cocyle
representing an element of $\Ker^1(\Q,\U(1))$. Note that $\U(1)(\overline{\Q})
\simeq\overline{\Q}^\times$, that $\Gal(\overline{\Q}/\Q)$ acts
on $\U(1)(\overline{\Q})$ via its quotient $\Gal(E/\Q)$, and
that $\tau\in\Gal(E/\Q)$ acts by $t\fle t^{-1}$. In particular, the restriction
of $c$ to $\Gal(\overline{\Q}/E)$ is a group morphism $\Gal(\overline{\Q}/E)
\fl\overline{\Q}^\times$. As this restriction is locally trivial, the
{\v C}ebotarev
density theorem implies that $c(\Gal(\overline{\Q}/E))=1$. So we can
see $c$ as a $1$-cocycle $\Gal(E/\Q)\fl\U(1)(\overline{\Q})$. As
$\Gal(E/\Q)\simeq\Gal(\C/\R)$ and $c$ is locally a coboundary, this
implies that $c$ is a coboundary.

By the description of $\widehat{\G}$ given above,
\[Z(\widehat{\G})=\{(\lambda,\lambda_1 I_{n_1},\dots,\lambda_r I_{n_r}),
\lambda,\lambda_1,\dots,\lambda_r\in\C^\times\},\]
with the action of $\Gal(E/\Q)$ given by
\[\tau((\lambda,\lambda_1I_{n_1},\dots,\lambda_rI_{n_r}))=
(\lambda\lambda_1^{n_1}\dots,\lambda_r^{n_r},\lambda_1^{-1}I_{n_1},\dots,
\lambda_r^{-1}I_{n_r}).\]
The second statement is now clear.

\item[(ii)] It suffices to show that $\Ker^1(\Q,\Le)=\{1\}$ and that
$Z(\widehat{\Le})^{\Gal(F/\Q)}$ is connected. The first equality comes
from the fact that
\[\Ho^1(\Q,\Le)=\Ho^1(F,\GL_n)=\{1\}.\]
On the other hand, $\widehat{\Le}=\GL_n(\C)^{[F:\Q]}$, with the obvious action
of $\Gal(F/\Q)$, so
$Z(\widehat{\Le})^{\Gal(F/\Q)}\simeq\C^\times$ is connected.

\end{itemize}
\end{proof}

\section{Levi subgroups and endoscopic groups}
\label{groupes4}

In this section, we recall some notions defined in section 7 of
\cite{K-NP}. Notations and definitions are as in section 7 of
\cite{K-STF:CTT}.

Let $\G$ be a connected reductive group on a local or global field $F$.
Let $\Ell(\G)$ be the set of isomorphism classes of elliptic endoscopic
triples for $\G$ (in the sense of \cite{K-STF:CTT} 7.4)
\index{EG@$\Ell(\G)$}
and $\Levi(\G)$
\index{LG@$\Levi(\G)$}
be the set of $\G(F)$-conjugacy classes of Levi subgroups of $\G$.
Let $\M$ be a Levi subgroup of $\G$. There is a canonical
$\Gal(\overline{F}/F)$-equivariant embedding
$Z(\widehat{\G})\fl Z(\widehat{\M})$.

\begin{definition}(\cite{K-NP} 7.1) An endoscopic $\G$-triple for $\M$ is
an endoscopic triple $(\M',s_M,\eta_{M,0})$ for $\M$ such that :
\begin{itemize}
\item[(i)] the image of $s_M$ in $Z(\widehat{\M'})/Z(\widehat{\G})$ is fixed
by $\Gal(\overline{F}/F)$;
\item[(ii)] the image of $s_M$ in $\Ho^1(F,Z(\widehat{\G}))$ (via the
morphism $\pi_0((Z(\widehat{\M'})/Z(\widehat{\G}))^{\Gal(\overline{F}/F)})
\fl\Ho^1(F,Z(\widehat{\G}))$ of \cite{K-STF:CTT} 7.1) is trivial if $F$ is
local, and in $\Ker^1(F,Z(\widehat{\G}))$ if $F$ is global.
\index{endoscopic $\G$-triple}

\end{itemize}
The $\G$-triple $(\M',s_M,\eta_{M,0})$ is called elliptic if it is
elliptic as an endoscopic triple for $\M$.

Let $(\M'_1,s_1,\eta_{1,0})$ and $(\M'_2,s_2,\eta_{2,0})$ be endoscopic
$\G$-triples for $\M$. An isomorphism of endoscopic $\G$-triples from
$(\M'_1,s_1,\eta_{1,0})$ to $(\M'_2,s_2,\eta_{2,0})$ is an isomorphism
$\alpha:\M'_1\fl\M'_2$ of endoscopic triples for $\M$ (in the sense of
\cite{K-STF:CTT} 7.5) such that the images of $s_1$ and $\widehat{\alpha}(s_2)$
in $Z(\widehat{\M'_1})/Z(\widehat{\G})$ are equal.

\end{definition}

Let $(\M',s_M,\eta_{M,0})$ be an endoscopic $\G$-triple for $\M$.
Then there is an isomorphism class of endoscopic triples for $\G$ associated
to $(\M',s_M,\eta_{M,0})$ in the following way (cf \cite{K-NP} 3.7 et 7.4) :
There is a canonical $\widehat{\G}$-conjugacy class of embeddings
${}^L\M\fl{}^L\G$; fix an element in this class, and use it to see
${}^L\M$ as a subgroup of ${}^L\G$. Define a subgroup $\Mcal$ of ${}^L\M$
as follows : an element $x\in{}^L\M$ is in $\Mcal$ if and only if there exists
$y\in{}^L\M'$ such that the images of $x$ and $y$ by the projections
${}^L\M\fl W_F$ and ${}^L\M'\fl W_F$ are the same
and that $\Int(x)\circ\eta_0=\eta_0\circ\Int(y)$. Then the restriction to
$\Mcal$ of the projection ${}^L\M\fl W_F$ is surjective, and $\Mcal\cap
\widehat{\M}=\eta_0(\widehat{\M'})$. Moreover, $\Mcal$ is a closed subgroup of
${}^L\M$.
\footnote{\cite{K-NP} 3.4 : Let $K$ be a finite extension of $F$ over which
$\M'$ and $\M$ split. Define a subgroup $\Mcal_K$ of
$\widehat{\M}\rtimes\Gal(K/F)$ in the same way as $\Mcal$.
This subgroup is obviously closed, and $\Mcal$ is the inverse image of
$\Mcal_K$.}
Set $\widehat{\H}=\Cent_{\widehat{\G}}(s_M)^0$, and $\Hcal=\Mcal
\widehat{\H}$. Then $\Hcal$ is a closed subgroup of ${}^L\G$,
the restriction to $\Hcal$ of the projection ${}^L\G\fl W_F$ is surjective,
and $\Hcal\cap\widehat{\G}=\widehat{\H}$. Hence $\Hcal$ induces a morphism
$\rho:W_F\fl\Out(\widehat{\H})$. Moreover, there exists a finite extension
$K$ of $F$ and a closed subgroup $\Hcal_K$ of $\widehat{\G}\rtimes\Gal(
K/F)$ such that $\Hcal$ is the inverse image of $\Hcal_K$.
\footnote{\cite{K-NP} 3.5 : Let $K'$ be a finite extension of $F$ over which
$\G$ splits. The group $\widehat{\H}$ is of finite index in its normalizer
in $\widehat{\G}$, so the group $\Hcal$ is of finite index in its normalizer
$\Ncal$ in ${}^L\G$. Hence the intersection of $\Hcal$ with the subgroup
$\widehat{\H}\times W_{K'}$ of $\Ncal$ is a closed subgroup of finite index
of $\widehat{\H}\times W_{K'}$; so it is also an open subgroup.
Hence $\Hcal$ contains an open subgroup of $\widehat{\H}\times W_{K'}$,
ie it contains a subgroup $\widehat{\H}\times W_K$, with $K$ a finite
extension of $K'$. The sought-for group $\Hcal_K$ is
$\widehat{\H}\times W_K$.}
Hence $\rho$ factors through $W_F\fl\Gal(K/F)$, and $\rho$ can be seen as
a morphism $\Gal(\overline{F}/F)\fl\Out(\widehat{\H})$. It is easy to see
that $(s_M\mod Z(\widehat{\G}),\rho)$ is an endoscopic datum for $\G$
(in the sense of \cite{K-STF:CTT}), and that its isomorphism class depends
only on the isomorphism class of $(\M',s_M,\eta_{M,0})$. We associate
to $(\M',s_M,\eta_{M,0})$ the isomorphism class of endoscopic triples
for $\G$ corresponding to
$(s_M\mod Z(\widehat{\G}),\rho)$ (cf \cite{K-STF:CTT} 7.6).

Let $\Ell_\G(\M)$ be the set of isomorphism classes of endoscopic $\G$-triples
$(\M',s_M,\eta_{M,0})$
\index{EGM@$\Ell_\G(\M)$}
for $\M$ such that the isomorphism class of endoscopic triples for $\G$
associated to $(\M',s_M,\eta_{M,0})$ is elliptic.
There are obvious maps
$\Ell_\G(\M)\fl\Ell(\M)$ and $\Ell_\G(\M)\fl\Ell(\G)$.
For every endoscopic $\G$-triple $(\M',s_M,\eta_{M,0})$ for $\M$, let
$\Aut_G(\M',s_M,\eta_{M,0})$ be the group of $\G$-automorphisms
of $(\M',s_M,\eta_{M,0})$ and $\Lambda_G(\M',s_M,\eta_{M,0})=
\Aut_G(\M',s_M,\eta_{M,0})/\M'_{ad}(F)$ be the group of outer
$\G$-automorphisms; if $\M=\G$, we will omit the subscript $\G$.
\index{$\Lambda_G(\M',s_M,\eta_{M,0})$}

Remember that we write $n_M^G=|\Nor_\G(\M)(F)/\M(F)|$ (cf
\ref{points_fixes6}).
\index{nMG@$n_M^G$}

Lemma \ref{lemme:Levi_endoscopiques} below is a particular case of
lemma 7.2 of \cite{K-NP}. As \cite{K-NP} is (as yet) unpublished, we prove
lemma \ref{lemme:Levi_endoscopiques} by a direct calculation.
Assume that $\G$ is one of the unitary groups of
\ref{groupes1} (and that $F=\Q$). If $(\M',s_M,\eta_{M,0})\in\Ell_\G(\M)$
and if $(\H,s,\eta_0)$ is its image in $\Ell(\G)$, then it is easy to see
that $\M'$ determines a $\H(\Q)$-conjugacy class of Levi subgroups of $\H$.
\footnote{In the case of unitary groups, this can be seen simply by
writing explicit formulas for $\H$, $\M$ and $\M'$. Actually, this fact is
true in greater generality and proved in \cite{K-NP} 7.4 (but we will not need
this here) : with notations as above, the group $\Mcal$ is a Levi
subgroup of $\Hcal$ (for a suitable definition of ``Levi subgroup'' in that
context), and gives a conjugacy class of Levi subgroups of $\H$ because
$\H$ is quasi-split.}

\begin{lemma}\label{lemme:Levi_endoscopiques} Assume that $\G$ is
quasi-split. Let $\varphi:\coprod\limits
_{(\H,s,\eta_0)\in\Ell(\G)}\Levi(\H)\fl\C$. Then
\begin{flushleft}$\displaystyle{\sum_{(\H,s,\eta_0)\in\Ell(\G)}|\Lambda(\H,s,
\eta_0)|^{-1}\sum_{\M_H\in\Levi(\H)}(n_{M_H}^H)^{-1}\varphi(\H,\M_H)}$
\end{flushleft}
\begin{flushright}$\displaystyle{=\sum_{\M\in\Levi(\G)}(n_M^G)^{-1}
\sum_{(\M',s_M,\eta_{M,0})\in\Ell_\G(\M)}|\Lambda_G(\M',s_M,\eta_{M,0})|^{-1}
\varphi(\H,\M_H),}$\end{flushright}
where, in the second sum, $(\H,s,\eta_0)$ is the image of
$(\M',s_M,\eta_{M,0})$ in $\Ell(\G)$ and $\M_H$ is the element of
$\Levi(\H)$ associated to $(\M',s_M,\eta_{M,0})$.

\end{lemma}

(As $\M'$ and $\M_H$ are isomorphic, we will sometimes write $\M'$ instead
of $\M_H$.)

We will use this lemma only for functions $\varphi_H$ that vanish when
their second argument is not a cuspidal Levi subgroup (see theorem
\ref{th:points_fixes_GKM} for the definition of a cuspidal Levi subgroup).
In that case, the lemma is an easy consequence of lemma
\ref{lemme:Levi_endoscopiques_GU} below, that is proved in the same way
as proposition \ref{prop:groupes_endoscopiques}.

In the next lemma, we consider only the case of the group $\GU^*(n)$ in
order to simplify the notations. The case of
$\G(\U^*(n_1)\times\dots\times\U^*(n_r))$
is similar.

\begin{lemma}\label{lemme:Levi_endoscopiques_GU} Let $n\in\Nat^*$ and
$\G=\GU^*(n)$. Let $\T$ be the diagonal torus of $\G$, and identify
$\widehat{\T}$ with $\C^\times\times(\C^\times)^n$ as in \ref{groupes3}.
Let $\M$ be a cuspidal Levi subgroup of $\G$. Then $\M$ is isomorphic
to $(R_{E/\Q}\Gr_m)^r\times\GU^*(m)$, with $r,m\in\Nat$ such that $n=m+2r$.
Let $\T_M$ be the diagonal torus of $\M$. The dual group $\widehat{\M}$
is isomorphic to the Levi subgroup
\[\C^\times\times\left(\begin{array}{ccc}\begin{array}{ccc}* & & 0 \\
& \ddots & \\ 0 & & *\end{array} & & 0 \\
 & \GL_m(\C) & \\
0 & & \begin{array}{ccc}* & & 0 \\ & \ddots & \\ 0 & & *\end{array}\end{array}
\right)\]
(with blocks of size $r,m,r$) of $\widehat{\G}$. Fix an isomorphism
$\widehat{\T}_M\simeq\widehat{\T}$ compatible with this identification.

Then an element $(\M',s_M,\eta_{M,0})$ of $\Ell_\G(\M)$ is uniquely
determined by $s_M$. If we assume (as we may) that
$s_M\in\widehat{\T}_M\simeq\widehat{\T}$, then 
$s_M\in Z(\widehat{\G})(\{1\}\times\{\pm 1\}^n)^{\Gal(\overline{\Q}/\Q)}$.
For every $A\subset\{1,\dots,r\}$ and $m_1,m_2\in\Nat$ such that $m=m_1+
m_2$ and that $m_2$ is even, set
\[s_{A,m_1,m_2}=(s_1,\dots,s_r,\overbrace{1,\dots,1}^{m_1},\overbrace{-1,\dots,
-1}^{m_2},s_r,\dots,s_1),\]
with $s_i=-1$ if $i\in A$ and $s_i=1$ if $i\not\in A$.
If $r<n/2$, then the set of $(1,s_{A,m_1,m_2})$ is a system of representatives
of the set of equivalence classes of possible $s_M$.
If $r=n/2$ (so $m=0$), then every $s_M$ is equivalent to a $(1,s_{A,0,0})$,
and $(1,s_{A,0,0})$ and $(1,s_{A',0,0})$ are equivalent
if and only if $\{1,\dots,r\}=A\sqcup A'$.

Let $s_M=(1,(s_1,\dots,s_n))\in(\{1\}\times\{\pm 1\}^n)^{\Gal(\overline{\Q}/
\Q)}$. Let $(\M',s_M,\eta_{M,0})$ be the element of $\Ell_\G(\M)$ associated
to $s_M$, and $(\H,s,\eta_0)$ be its image in $\Ell(\G)$. Let
$n_1=|\{i\in\{1,\dots,n\}|s_i=1\}$, $n_2=n-n_1$,
$m_1=|\{i\in\{r+1,\dots,r+m\}|s_i=1\}|$, $m_2=m-m_1$, $r_1=
(n_1-m_1)/2$, $r_2=(n_2-m_2)/2$ ($r_1$ and $r_2$ are integers
by the condition on $s_M$).
Then $\H=\G(\U^*(n_1)\times
\U^*(n_2))$, $\M'=(R_{E/\Q}\Gr_m)^r\times\G(\U^*(m_1)\times\U^*(m_2))$, and
$n_{M'}^H=2^r(r_1)!(r_2)!$. Moreover, $|\Lambda_G(\M',s_M,\eta_{M,0})|$
is equal to $1$ if $\M\not=\G$.

\end{lemma}

We end this section by recalling a result of \cite{K-NP} 7.3.
Assume again that $\G$ is any connected reductive group on a local or
global field $F$. Let $\M$ be a Levi subgroup of $\G$.

\begin{definition} Let $\gamma\in\M(F)$ be semi-simple.
An endoscopic $\G$-quadruple for $(\M,\gamma)$ is a quadruple
$(\M',s_M,\eta_{M,0},\gamma')$, where $(\M',s_M,\eta_{M,0})$ is an
endoscopic $\G$-triple for $\M$ and $\gamma'\in\M'(F)$ is a semi-simple
$(\M,\M')$-regular element such that $\gamma$ is an image of $\gamma'$
(the unexplained expressions in this sentence are defined in
\cite{K-STF:EST} 3).
\index{endoscopic $\G$-quadruple}
An isomorphism of endoscopic $\G$-quadruples $\alpha:(\M'_1,s_{M,1},\eta_
{M,0,1},\gamma'_1)\fl(\M'_2,s_{M,2},\eta_{M,0,2},\gamma'_2)$ is an
isomorphism of endoscopic $\G$-triples $\alpha:\M'_1\fl\M'_2$ such that
$\alpha(\gamma'_1)$ and $\gamma'_2$ are stably conjugate.

\end{definition}

Let $I$ be a connected reductive subgroup of $\G$ that contains a maximal
torus of $\G$.
There is a canonical
$\Gal(\overline{F}/F)$-equivariant inclusion $Z(\widehat{\G})\subset
Z(\widehat{I})$. Let $\Kgoth_\G(I/F)$ be the set of elements in
$(Z(\widehat{I})/Z(\widehat{\G}))^{\Gal(\overline{F}/F)}$ whose image by
the morphism $(Z(\widehat{I})/Z(\widehat{\G}))^{\Gal(\overline{F}/F)}\fl
\H^1(F,Z(\widehat{\G}))$ (coming from the exact sequence
$1\fl Z(\widehat{\G})\fl Z(\widehat{I})\fl Z(\widehat{I})/Z(\widehat{\G})\fl
1$)
is trivial if $F$ is local and locally trivial
if $F$ is global.
\index{KGIF@$\Kgoth_\G(I/F)$}
\footnote{This definition is coherent with the definition of
$\Kgoth(I_0/\Q)$ in \ref{points_fixes6} : in \ref{points_fixes6}, $I_0$
is the centralizer of a semi-simple elliptic element, so
$(Z(\widehat{I_0})/Z(\widehat{\G}))^{\Gal(\overline{\Q}/\Q)}$
is finite.}
If $I$ is included in $\M$, there is an obvious morphism
$\Kgoth_\G(I/F)\fl\Kgoth_\M(I/F)$.

Fix $\gamma\in\M(F)$ semi-simple, and let $I=\Cent_\M(\gamma)^0$.
Let $(\M',s_M,\eta_{M,0},\gamma')$ be an endoscopic $\G$-quadruple for
$(\M,\gamma)$. Let $I'=\Cent_{\M'}(\gamma')^0$. As $\gamma'$ is
$(\M,\M')$-regular, $I'$ is an inner form of $I$ (cf \cite{K-STF:EST} 3),
so there is a canonical isomorphism $Z(\widehat{I})\simeq
Z(\widehat{I'})$. Let $\kappa(\M',s_M,\eta_{M,0},\gamma')$ be the image of
$s_M$ by the morphism $Z(\widehat{\M'})\subset Z(\widehat{I'})\simeq
Z(\widehat{I})$.

\begin{lemma}\label{lemme:param_groupes_endoscopiques} The map
$(\M',s_M,\eta_{M,0},\gamma')\fle\kappa(\M',s_M,\eta_{M,0},\gamma')$ induces
a bijection from the set of isomorphism classes of endoscopic
$\G$-quadruples for $(\M,\gamma)$ to $\Kgoth_\G(I/F)$. Moreover,
the automorphisms of endoscopic $\G$-quadruples for $(\M,\gamma)$
are all inner.

\end{lemma}

This lemma is lemma 7.1 of \cite{K-NP}. It is a generalization of lemma
9.7 of \cite{K-STF:EST} and can be proved in the same way.

\chapter{Discrete series}
\label{serie_discrete}

\section{Notations}
\label{serie_discrete1}

Let $\G$ be a connected reductive algebraic group over $\R$.
In this chapter, we form the $L$-groups with the Weil group $W_\R$.
\index{WR@$W_\R$\quad Weil group of $\R$}
\index{WC@$W_\C$\quad Weil group of $\C$}
Remember that $W_\R=W_\C\sqcup W_\C\tau$, with $W_\C=\C^\times$, $\tau^2=-1\in
\C^\times$ 
and, for every $z\in\C^\times$, $\tau z\tau^{-1}=\overline{z}$, and that
$W_\R$ acts on $\widehat{\G}$ via its quotient  $\Gal(\C/\R)\simeq W_\R/W_\C$.
Let $\Pi(\G(\R))$ (resp. $\Pi_{temp}(\G(\R))$) be the set of equivalence
classes of irreductible (resp. irreducible and tempered)
admissible representations of $\G(\R)$.
\index{$\Pi(\G(\R))$}
\index{$\Pi_{temp}(\G(\R))$}
For every $\pi\in\Pi(\G(\R))$, let $\Theta_\pi$ be the Harish-Chandra
character of $\pi$ (seen as a real analytic function on the set
$\G_{reg}(\R)$ of regular elements of $\G(\R)$).
\index{$\Theta_\pi$\quad Harish-Chandra character of $\pi$}
\index{Greg@$\G_{reg}$\quad set of regular elements}

Assume that $\G(\R)$ has a discrete series. Let $\A_G$ be the maximal
($\R$-)split torus in the center of $\G$ and $\overline{\G}$ be an inner
form of $\G$ such that $\overline{\G}/\A_G$ is $\R$-anisotropic.
\index{G@$\overline{\G}$}
Write $q(G)=\dim(X)/2$, where $X$ is the symmetric space of $\G(\R)$.
\index{qG@$q(\G)$}
Let $\Pi_{disc}(\G(\R))\subset\Pi(\G(\R))$
be the set of equivalence classes of representations in the discrete series.
\index{$\Pi_{disc}(\G(\R))$}

The set $\Pi_{disc}(\G(\R))$ is the disjoint union of finite subsets called
$L$-packets; $L$-packets all have the same number of elements and are
parametrized by equivalence classes of elliptic Langlands parameters
$\varphi:W_\R\fl{}^L\G$, or,
equivalently, by isomorphism classes of irreducible representations $E$
of $\overline{\G}(\R)$. Let $\Pi(\varphi)$ (resp. $\Pi(E)$) be the $L$-packet
associated to the parameter $\varphi$ (resp. to the representation $E$), and
let $d(\G)$ be the cardinality of a $L$-packet of $\Pi_{disc}(\G)$.
\index{discrete series $L$-packet}
\index{$\Pi(\varphi)$}
\index{$\Pi(E)$}
\index{dG@$d(\G)$}

If $\pi\in\Pi_{disc}(\G(\R))$, we will write $f_\pi$ for a
pseudo-coefficient of $\pi$ (cf \cite{CD}).
\index{fpi@$f_\pi$\quad pseudo-coefficient of $\pi$}

For every elliptic Langlands parameter $\varphi:W_\R\fl{}^L\G$, write
\[S\Theta_\varphi=\sum_{\pi\in\Pi(\varphi)}\Theta_\pi.\]
\index{STheta@$S\Theta_\varphi$}

We are going to calculate the integer $d(\G)$ for the unitary groups of
\ref{groupes1}. The following definition will be useful (this notion
already appeared in theorem \ref{th:points_fixes_GKM} and in section
\ref{groupes4}).

\begin{definition}\label{def:groupe_cuspidal}
Let $\G$ be a connected reductive group over $\Q$. Denote by $\A_G$ the maximal
$\Q$-split torus in the center of $\G$. $\G$ is called \emph{cuspidal}
if the group $(\G/\A_G)_\R$ has a maximal $\R$-torus that is $\R$-anisotropic.

\end{definition}
\index{cuspidal group}

Let
$p,q\in\Nat$ be such that $p\geq q$ and $p+q\geq 1$. Let $\G=\GU(p,q)$, and
use the Hermitian form $A_{p,q}$ of \ref{groupes2} to define $\G$.
Then $\A_G=\Gr_mI_{p+q}$.
Let $\T$ be the diagonal maximal torus of $\G$.
Let
\[\T_{ell}=\left\{\left(\begin{array}{ccc}
\begin{array}{ccc}a_1 & & 0 \\ & \ddots & \\ 0 & & a_q\end{array} 
& \begin{array}{ccc}&& \\ & 0 & \\ && \end{array} 
& \begin{array}{ccc}0 & & b_1 \\ & \ddotsinv & \\ b_q & & 0 \end{array} \\
\begin{array}{ccc}&& \\ & 0 & \\ && \end{array} 
& \begin{array}{ccc}c_1 & & 0 \\ & \ddots & \\ 0 & & c_{p-q}  \end{array}
& \begin{array}{ccc}&& \\ & 0 & \\ &&\end{array} \\
\begin{array}{ccc}0 & & b_q \\ & \ddotsinv & \\ b_1 & & 0 \end{array}
& \begin{array}{ccc}&& \\ & 0 & \\ &&\end{array}
& \begin{array}{ccc}a_q & & 0 \\ & \ddots & \\ 0 & & a_1\end{array}
\end{array}\right)\right\},\]
where $a_i,b_i,c_i\in R_{E/\Q}\Gr_m$ are such that :
\[\left\{\begin{array}{ll}a_i\overline{b}_i+\overline{a}_ib_i=0\mbox{ for }
1\leq i\leq q \\
a_1\overline{a}_1+b_1\overline{b}_1=\dots=a_q\overline{a}_q+b_q\overline{b}_q=
c_1\overline{c}_1=\dots=c_{p-q}\overline{c}_{p-q}\end{array}\right..\]
Then $\T_{ell}$ is a maximal torus of $\G$, and $\T_{ell}/\A_G$ is
$\R$-anisotropic. So $\G$ is cuspidal.
Write
\[u_G=\frac{1}{\sqrt{2}}\left(\begin{array}{ccc}I_q & 0 & -(i\otimes i)J_q \\
0 & (\sqrt{2}\otimes 1)I_{p-q} & 0 \\ (i\otimes i)J_q & 0 & I_q\end{array}
\right)\in\SU(p,q)(\C).\]
Conjugacy by $u_G^{-1}$ is an isomorphism $\alpha:\T_{ell,\C}\iso
\T_\C$. Use $\alpha$ to identify $\widehat{\T}_{ell}$ and
$\widehat{\T}=\C^\times\times(\C^\times)^{p+q}$. Then the action of
$\Gal(\C/\R)=\{1,\tau\}$ on $\widehat{\T}_{ell}$ is given by :
\[\tau((\lambda,(\lambda_1,\dots,\lambda_{p+q})))=(\lambda\lambda_1\dots
\lambda_{p+q},(\lambda_1^{-1},\dots,\lambda_{p+q}^{-1})).\]
Let $\Omega_\G=W(\T_{ell}(\C),\G(\C))$ and $\Omega_{\G(\R)}=W(\T_{ell}(\R),
\G(\R))$ be the Weyl groups of $\T_{ell}$ over $\C$ and $\R$. The group
$\Omega_\G\simeq W(\T(\C),\G(\C))\simeq\Sgoth_{p+q}$ acts on
$\T(\C)$ by permuting the diagonal entries.
The subgroup $\Omega_{\G(\R)}$ of $\Omega_{\G}$ is the group
$\Sgoth_p\times\Sgoth_q$ if $p\not=q$, and the union of
$\Sgoth_q\times\Sgoth_q$ and of the set of permutations that send
$\{1,\dots,q\}$ to $\{q+1,\dots,n\}$ if $p=q$.
Hence
\[d(\G)=\left\{\begin{array}{ll}\displaystyle{\frac{(p+q)!}{p!q!}} &
\mbox{ if }p\not=q\\
\displaystyle{\frac{(2q)!}{2(q!)^2}} & \mbox{ if }p=q\end{array}\right..\]

\begin{remark} The torus $\T_{ell}$ is isomorphic to
$\G(\U(1)^{p+q})$ by the morphism
\begin{flushleft}$\displaystyle{\beta_\G:\left(\begin{array}{ccc}
diag(a_1,\dots,a_q) & 0 & diag(b_1,\dots,b_q)J_q \\
0 & diag(c_1,\dots,c_{p-q}) & 0 \\
J_q diag (b_1,\dots,b_q) & 0 & diag(a_q,\dots,a_1)
\end{array}\right)}$\end{flushleft}
\begin{flushright}$\displaystyle{\fle (a_1-b_1,\dots,a_q-b_q,c_1,\dots,c_{p-q},
a_q+b_q,\dots,a_1+b_1).}$\end{flushright}

\end{remark}

These constructions have obvious generalizations to the groups
$\G(\U(p_1,q_1)\times\dots\times\U(p_r,q_r))$. (In particular, these groups
are also cuspidal.)

\section{The fonctions $\Phi_M(\gamma,\Theta)$}
\label{serie_discrete2}

In this section, we recall a construction of Arthur and Shelstad.

Let $\G$ be a connected reductive group on $\R$. A \emph{virtual character}
$\Theta$ on $\G(\R)$ is a linear combination with coefficients in $\Z$
of functions $\Theta_\pi$, $\pi\in\Pi(\G(\R))$. The virual character $\Theta$
is called \emph{stable} if
$\Theta(\gamma)=\Theta(\gamma')$ for every $\gamma,\gamma'\in\G_{reg}(\R)$
that are stably conjugate.
\index{virtual character}
\index{stable virtual character}

Let $\T$ be a maximal torus of $\G$. Let $\A$ be the maximal split
subtorus of $\T$ and $\M=\Cent_\G(\A)$ (a Levi subgroup of $\G$).
For every $\gamma\in\M(\R)$, set
\[D_M^G(\gamma)=\det(1-\Ad(\gamma),Lie(\G)/Lie(\M)).\]
\index{DMG@$D_M^G$\quad partial Weyl denominator}

\begin{lemma}\label{lemme:Phi_M}
(\cite{A-L2} 4.1, \cite{GKM} 4.1) Let $\Theta$ be a stable virtual character
on $\G(\R)$. Then the function
\[\gamma\fle |D_M^G(\gamma)|_\R^{1/2}\Theta(\gamma)\]
on $\T_{reg}(\R)$ extends to a continuous function on $\T(\R)$, that will be
denoted by $\Phi_M(.,\Theta)$ or $\Phi_M^G(.,\Theta)$.

\end{lemma}
\index{$\Phi_M^G$}

We will often see $\Phi_M(.,\Theta)$ as a function on $\M(\R)$,
defined as follows : if $\gamma\in\M(\R)$ is $\M(\R)$-conjugate to a
$\gamma'\in\T(\R)$, set $\Phi_M(\gamma,\Theta)=\Phi_M(\gamma',\Theta)$; if
there is no element of $\T(\R)$ conjugate to $\gamma\in\M(\R)$, set
$\Phi_M(\gamma,\Theta)=0$.

\begin{remark}\label{rq:Phi_M}
The function $\Phi_M(.,\Theta)$ on $\M(\R)$ is invariant by conjugacy by
$\Nor_\G(\M)(\R)$ (because $\Theta$ and $D_M^G$ are).

\end{remark}

\section{Transfer}
\label{serie_discrete3}

We first recall some definitions from \cite{K-SVLR} \S7.

Let $\G$ be a connected reductive algebraic group over $\Q$. For every
maximal torus $\T$ of $\G$, let $\Bo_G(\T)$ be the set of Borel subgroups
\index{BGT@$\Bo_G(\T)$}
of $\G_\C$ containing $\T$. Assume that $\G$ has a maximal torus
$\T_G$ such that $(\T_G/\A_G)_\R$ is anisotropic, and let
$\overline{\G}$ be an inner form of $\G$ over $\R$ such that
$\overline{\G}/\A_{G,\R}$ is anisotropic.
Write $\Omega_G=W(\T_G(\C),\G(\C))$.
\index{$\Omega_G$}
Let $\varphi:W_\R\fl{}^L\G$ be an elliptic Langlands parameter.

Let $(\H,s,\eta_0)$ be an elliptic endoscopic triple for $\G$. Choose a
$L$-morphism $\eta:{}^L\H\fl{}^L\G$ extending $\eta_0:\widehat{\H}\fl
\widehat{\G}$ (we assume that such a $\eta$ exists),
and let $\Phi_H(\varphi)$ be the set of equivalence classes
of Langlands parameters $\varphi_H:W_\R\fl{}^L\H$ such that
$\eta\circ\varphi_H$ and $\varphi$ are equivalent.
\index{$\Phi_H(\varphi)$}
Assume that the torus $\T_G$ comes from a maximal torus $\T_H$ of $\H$, and fix
an admissible isomorphism $j:\T_H\iso\T_G$. Write
$\Omega_H=W(\T_H(\C),\H(\C))$. Then
$j_*(\Phi(\T_H,\H))\subset\Phi(\T_G,\G)$, so $j$ induces a map
$j^*:\Bo_G(\T_G)\fl\Bo_H(\T_H)$ and an injective morphism
$\Omega_H\fl\Omega_G$; we use this morphism to see $\Omega_H$ as a subgroup
of $\Omega_G$.

Let $\B\in\Bo_G(\T_G)$, and let $\B_H=j^*(\B)$. Set
\[\begin{array}{rcl}\Omega_* & = & \{\omega\in\Omega_G|j^*(\omega(\B))=\B_H\}\\
& = & \{\omega\in\Omega_G|\omega^{-1}(j_*(\Phi(\T_H,\B_H)))\subset\Phi(\T_G,
\B)\}.\end{array}\]
\index{$\Omega_*$}
Then, for every $\omega\in\Omega_G$, there exists a unique pair
$(\omega_H,\omega_*)\in\Omega_H\times\Omega_*$ such that $\omega=\omega_H
\omega_*$. Moreover, there is a bijection $\Phi_H(\varphi)\iso\Omega_*$
defined as follows :
if $\varphi_H\in\Phi_H(\varphi)$, send it to the unique $\omega_*(\varphi_H)
\in\Omega_*$ such that $(\omega_*(\varphi_H)^{-1}\circ j,\B,\B_H)$ is aligned
with $\varphi_H$ (in the sense of \cite{K-SVLR} \S7 p 184).

The Borel subgroup $\B$ also defines a $L$-morphism
$\eta_B:{}^L\T_G\fl{}^L\G$, unique up to $\widehat{\G}$-conjugacy
(cf \cite{K-SVLR} p 183).
\index{$\eta_B$}

We will use the normalization of the transfer factors of \cite{K-SVLR} \S7,
that we recall in the next definition.

\begin{definition}
For every $\gamma_H\in\T_H(\R)$, set (notations are as above) :
\[\Delta_{j,B}(\gamma_H,\gamma)=(-1)^{q(\G)+q(\H)}\chi_B(\gamma)\prod_
{\alpha\in\Phi(\T_G,\B)-j_*(\Phi(\T_H,\B_H))}(1-\alpha(\gamma^{-1})),\]
where $\gamma=j(\gamma_H)$ and $\chi_B$ is the quasi-character of
$\T_G(\R)$ associated to the $1$-cocyle $a:W_\R\fl\widehat{\T}_G$ such that
$\eta\circ\eta_{B_H}\circ\widehat{j}$ and $\eta_B.a$ are conjugate under
$\widehat{\G}$.

\end{definition}
\index{$\Delta_{j,B}$\quad real transfer factor}

\begin{remark}\label{rq:facteur_transfert}\begin{itemize}
\item[(1)] Let $\varphi_H\in\Phi_H(\varphi)$ be
such that $\omega_*(\varphi_H)=1$.
After replacing $\varphi$ (resp. $\varphi_H$) by a $\widehat{\G}$-conjugate
(resp. a $\widehat{\H}$-conjugate), we can write
$\varphi=\eta_B\circ\varphi_B$ (resp. $\varphi_H=\eta_{B_H}\circ
\varphi_{B_H}$), where $\varphi_B$ (resp. $\varphi_{B_H}$) is a Langlands
parameter for $\T_G$ (resp. $\T_H$).
Let $\chi_{\varphi,B}$ (resp. $\chi_{\varphi_H,B_H}$) be the quasi-character
of $\T_G(\R)$ (resp. $\T_H(\R)$) associated to $\varphi_B$
(resp. $\varphi_{B_H}$). Then $\chi_B=\chi_{\varphi,B}(\chi_{\varphi_H,B_H}
\circ j^{-1})^{-1}$.
\index{$\chi_{\varphi,B}$}
\item[(2)] Let $\omega\in\Omega_G$. Write $\omega=\omega_H\omega_*$, with
$\omega_H\in\Omega_H$ and $\omega_*\in\Omega_*$. Then $\Delta_{j,\omega(B)}=
\det(\omega_*)\Delta_{j,B}$, where $\det(\omega_*)=\det(\omega_*,X^*(\T_G))$.

\end{itemize}
\end{remark}

Let $p,q\in\Nat$ be such that $p\geq q$ and that $n:=p+q\geq 1$.
Fix $n_1,n_2\in\Nat^*$ such that $n_2$ is even and
$n_1+n_2=n$. Let $\G$ be the group $\GU(p,q)$ and $(\H,s,\eta_0)$ be the
elliptic endoscopic triple for $\G$ associated to $(n_1,n_2)$ as in
proposition \ref{prop:groupes_endoscopiques}. In section
\ref{serie_discrete1}, we defined elliptic maximal tori
$\T_G=\T_{G,ell}$ and $\T_H=\T_{H,ell}$ of $\G$ and $\H$,
and isomorphisms $\beta_G:\T_G\iso\G(\U(1)^n)$ and $\beta_H:\T_H\iso\G(
\U(1)^n)$. Take $j=\beta_G^{-1}\circ\beta_H:\T_H\iso\T_G$. (It is easy to
see that this is an admissible isomorphism.)
We also defined
$u_G\in\G(\C)$ such that $\Int(u_G^{-1})$ sends $\T_{G,\C}$ to the diagonal
torus $\T$ of $\G_\C$. This defines an isomorphism (not compatible with
Galois actions in general) $\widehat{\T}_G\simeq\widehat{\T}$. The
composition of this isomorphism and of the embedding $\widehat{\T}\subset
\widehat{\G}$ defined in \ref{groupes3} gives an embedding
$\widehat{\T}_G\subset\widehat{\G}$. Conjugacy by $u_G$ also gives an
isomorphism $\Omega_G\simeq\Sgoth_n$. Via this isomorphism,
$\Omega_H=\Sgoth_{n_1}\times\Sgoth_{n_2}$, embedded in $\Sgoth_n$ in the
obvious way.

\begin{remark}\label{rq:le_bon_Borel}
It is easy to give simple descriptions of the subset $\Omega_*$ of
$\Omega_G$ and of the bijection $\Phi_H(\varphi)\iso\Omega_*$ for a
particular choice of $\B\in\Bo_G(\T_G)$.
Let
\[\B=\Int(u_G)\left(\begin{array}{ccc}* & & * \\ & \ddots & \\ 0 & & *
\end{array}\right).\]
Then
\[\Omega_*=\{\sigma\in\Sgoth_n|\sigma^{-1}_{|\{1,\dots,n_1\}}\mbox{ and }
\sigma^{-1}_{|\{n_1+1,\dots,n\}}\mbox{ are non-decreasing}\}.\]
As $W_\C$ is commutative, we may assume after replacing $\eta$ by a
$\widehat{\G}$-conjugate that $\eta$ sends $\{1\}\times W_\C\subset{}^L\H$
to $\widehat{\T}_G\times W_\C\subset{}^L\G$. As moreover $W_\C$ acts
trivially on $\widehat{\H}$, 
\[\eta((1,z))=((z^a\overline{z}^{b},\left(\begin{array}{cc}z^{a_1}\overline{z}
^{b_1}I_{n_1} & 0 \\ 0 & z^{a_2}\overline{z}^{b_2}I_{n_2}\end{array}\right))
,z),\quad z\in W_\C,\]
with $a,b,a_1,a_2,b_1,b_2\in\C$ such that $a-b,a_1-b_1,a_2-b_2\in\Z$.
Let $\varphi:W_\R\fl {}^L\G$ be an elliptic Langlands parameter.
We may assume that $\varphi_{|W_\C}$ is
\[z\fle ((z^\lambda\overline{z}^\mu,\left(\begin{array}{ccc}z^{\lambda_1}
\overline{z}^{\mu_1}& & 0 \\ & \ddots & \\ 0 & & z^{\lambda_n}\overline{z}^
{\mu_n}\end{array}\right)),z),\] 
with $\lambda,\mu,\lambda_1,\dots,\lambda_n,\mu_1,\dots,\mu_n\in\C$ such that
$\lambda-\mu\in\Z$, $\lambda_i-\mu_i\in\Z$ for every $i\in\{1,\dots,n\}$ and 
that the $\lambda_i$ are pairwise distinct.
Then there is a commutative diagram
\[\xymatrix{\Phi_H(\varphi)\ar[rr]^-{\sim} & & \Omega_*\ar[ld]^-\sim \\
& \{I\subset\{1,\dots,n\},|I|=n_1\}\ar[lu]^-\sim & }\]
where :
\begin{itemize}
\item[$*$] the horizontal arrow is $\varphi_H\fle\omega_*(\varphi_H)$;
\item[$*$] the arrow on the right is $\omega_*\fle\omega_*^{-1}
(\{1,\dots,n_1\})$;
\item[$*$] if $I\subset\{1,\dots,n\}$ has $n_1$ elements, write
$I=\{i_1,\dots,i_{n_1}\}$ and $\{1,\dots,n\}-I=\{j_1,\dots,j_{n_2}\}$ with
$i_1<\dots<i_{n_1}$ and $j_1<\dots<j_{n_2}$, and associate to $I$ the
unique $\varphi_H\in\Phi_H(\varphi)$ such that, for $z\in W_\C$,
\begin{flushleft}$\displaystyle{\varphi_H(z)=((z^{\lambda-a}\overline{z}^
{\mu-b},\left(\begin{array}{ccc}
z^{\lambda_{i_1}-a_1}
\overline{z}^{\mu_{i_1}-b_1} & & 0 \\ & \ddots & \\ 0 & & z^{\lambda_{i_{n_1}
-a_1}}\overline{z}^{\mu_{i_{n_1}-b_1}}\end{array}\right),}$\end{flushleft}
\begin{flushright}$\displaystyle{\left(\begin{array}{ccc}z^{\lambda_{j_1}-
a_2}\overline{z}^{\mu_{j_1}-b_2} & & 0 \\ & \ddots & \\ 0 & & z^{
\lambda_{j_{n_2}-a_2}}
\overline{z}^{\mu_{j_{n_2}-b_2}}\end{array}\right)),z).}$\end{flushright}

\end{itemize}

\end{remark}

\vspace{.5cm}

Remember that a Levi subgroup of $\G$ or $\H$ is called \emph{standard}
if it is a Levi subgroup of a standard parabolic subgroup and contains
the diagonal torus.
Let $\M$ be a cuspidal standard Levi subgroup of $\G$,
and let $r\in\{1,\dots,q\}$ be such that $\M=\M_{\{1,\dots,r\}}\simeq
(R_{E/\Q}\Gr_m)^r\times\GU(p-r,q-r)$. Let $(\M',s_M,\eta_{M,0})$ be an
element of $\Ell_\G(\M)$ (cf \ref{groupes4}) whose image in $\Ell(\G)$ is
$(\H,s,\eta_0)$.
There is a conjugacy class of Levi subgroups of $\H$ associated to
$(\M',s_M,\eta_{M,0})$; let $\M_H$ be the standard Levi subgroup in this
class.
If $m_1,m_2,r_1,r_2$ are defined as in lemma
\ref{lemme:Levi_endoscopiques_GU}, then
\[\M_H=\H\cap\left(\left(\begin{array}{ccc}
\T_{r_1} & 0 & 0 \\
0 & \GU^*(m_1) & 0 \\
0 & 0 & \T_{r_1}\end{array}\right),\left(\begin{array}{ccc}
\T_{r_2} & 0 & 0 \\
0 & \GU^*(m_2) & 0 \\
0 & 0 & \T_{r_2}\end{array}\right)\right),\]
with
\[\T_{r_1}=\left(\begin{array}{ccc}* &  & 0 \\  & \ddots &  \\ 0 &  & *
\end{array}\right)\subset R_{E/\Q}\GL_{r_1}\quad\mbox{and}\quad\T_{r_2}=
\left(\begin{array}{ccc}* & & 0 \\ & \ddots & \\ 0 & & *\end{array}\right)
\subset R_{E/\Q}\GL_{r_2}.\]
Hence $\M_H\simeq\G(\U^*(m_1)\times\U^*(m_2))\times(R_{E/\Q}\Gr_m)^{r_1+r_2}$.
Set
\[\T_{M_H}=\T_{\G(\U^*(m_1)\times\U^*(m_2)),ell}\times(R_{E/\Q}\Gr_m)
^{r_1+r_2}\]
\[\T_M=\T_{\GU(p-r,q-r),ell}\times(R_{E/\Q}\Gr_m)^r.\]
Then $\T_{M_H}$ (resp. $\T_M)$ is an elliptic maximal torus of
$\M_H$ (resp. $\M$).
We have isomorphisms
\begin{flushleft}$\displaystyle{\beta_{\G(\U^*(m_1)\times\U^*(m_2))}\times id
:\T_{M_H}=\T_{\G(\U^*(m_1)\times\U^*(m_2)),ell}\times(R_{E/\Q}\Gr_m)^r}$
\end{flushleft}
\begin{flushright}$\displaystyle{\iso\G(\U(1)^{m_1+m_2})\times(R_{E/\Q}
\Gr_m)^r}$\end{flushright}
\begin{flushleft}$\displaystyle{\beta_{\GU(p-r,q-r)}^{-1}\times id :
\G(\U(1)^{p+q-2r})\times(R_{E/\Q}\Gr_m)^{r}}$\end{flushleft}
\begin{flushright}$\displaystyle{\iso\T_{\GU(p-r,q-r),ell}\times(R_{E/\Q}\Gr_m)
^r=\T_M.}$\end{flushright}
Let
\[j_M:\T_{M_H}\iso\T_M\]
be the composition of these isomorphisms (note that $p+q-2r=m_1+m_2$). 

As before, the isomorphism $j_M$ is admissible and induces maps
$j_{M*}:\Phi(\T_{M_H},\H)\fl\Phi(\T_M,\G)$
and $j_M^*:\Bo_G(\T_M)\fl\Bo_H(\T_{M_H})$ (and similar maps if we replace
$\G$ by $\M$ and $\H$ by $\M_H$). It is easy to see that all the real roots of
$\Phi(\T_M,\G)$ (resp. $\Phi(\T_M,\M)$) are in
$j_{M*}(\Phi(\T_{M_H},\H))$ (resp. $j_{M*}(\Phi(\T_{M_H},\M_H))$).

Define an element $u\in\GU(p-r,q-r)(\C)$ in the same way as the elemnt
$u_G$ of \ref{serie_discrete1}, such that $\Int(u^{-1})$ sends
$\T_{\GU(p-r,q-r),ell,\C}$ to the diagonal torus of $\GU(p-r,q-r)_\C$. Let
$u_M=diag(I_r,u,I_r)u_G^{-1}\in\G(\C)$. Then $\Int(u_M^{-1})$ sends
$\T_{M,\C}$ onto $\T_{G,\C}$. Similarly, we get
$u_{M_H}\in\H(\C)$ such that $\Int(u_{M_H}^{-1})$ sends $\T_{M_H}$ to
$\T_{H,\C}$.
The following diagram is commutative :
\[\xymatrix@C=40pt{\T_{M,\C}\ar[r]^-{\Int(u_M^{-1})}\ar[d]_{j_M} & \T_{G,\C}
\ar[d]^j \\
\T_{M_H,\C}\ar[r]_-{\Int(u_{M_H}^{-1})} & \T_{H,\C}}\]

Use conjugacy by $u_M$ (resp. $u_{M_H}$) to identify
$\Omega_G$ (resp. $\Omega_H$) and $W(\T_M(\C),\G(\C))$ (resp.
$W(\T_{M_H}(\C),\H(\C))$). If $\B\in\Bo_G(\T_M)$, we use
$\Int(u_M^{-1})(\B)\in\Bo_G(\T_G)$ to define (as before) a subset
$\Omega_*$ of $\Omega_G$ and a bijection $\Phi_H(\varphi)\iso
\Omega_*$.

By \cite{K-NP} p 23, the morphism $\eta$ determines a $L$-morphism
$\eta_M:{}^L\M_H={}^L\M'\fl{}^L\M$, unique up to $\widehat{\M}$-conjugacy
and extending $\eta_{M,0}$. We use this morphism $\eta_M$ to define transfer
factors $\Delta_{j_M,\B_M}$, for every $\B_M\in\Bo_M(\T_M)$ : if
$\gamma_H\in\T_{M_H}(\R)$, set
\[\Delta_{j_M,\B_M}(\gamma_H,\gamma)=(-1)^{q(\G)+q(\H)}\chi_{\B_M}(\gamma)
\prod_{\alpha\in\Phi(\T_M,\B_M)-j_{M*}(\Phi(\T_{M_H},\B_{M_H}))}
(1-\alpha(\gamma^{-1}))\]
(note the sign),
where $\gamma=j_M(\gamma_H)$, $\B_{M_H}=j_M^*(\B_M)$ and $\chi_{B_M}$ is the
quasi-character of $\T_M(\R)$ associated to the $1$-cocyle
$a_M:W_\R\fl\widehat{\T}_M$ such that
$\eta_M\circ\eta_{B_{M_H}}\circ\widehat{j}_M$ and $\eta_{B_M}.a_M$ are
$\widehat{\M}$-conjugate.
\index{$\Delta_{j_M,\B_M}$\quad real transfer factor}

The next proposition is a generalization of the calculations of
\cite{K-SVLR} p 186.

\begin{proposition}\label{prop:transfert_caracteres}
Fix $\B\in\Bo_G(\T_M)$ (that determines $\Omega_*$
and $\Phi_H(\varphi)\iso\Omega_*$), et and let $\B_M=\B\cap\M$.
Let $\gamma_H\in\T_{M_H}(\R)$ and $\gamma=j_M(\gamma_H)$. 
Then
\[\Delta_{j_M,B_M}(\gamma_H,
\gamma)\Phi_M(\gamma^{-1},S\Theta_\varphi)=\sum_{\varphi_H\in\Phi_H
(\varphi)}\det(\omega_*(\varphi_H))\Phi_{M_H}(\gamma_H^{-1},S\Theta_
{\varphi_H}).\]

\end{proposition}

\begin{proof} Both sides of the equality that we want to prove depend on
the choice of the $L$-morphism $\eta:{}^L\H\fl{}^L\G$ extending
$\eta_0:\widehat{\H}\fl\widehat{\G}$. Let $\eta'$ be another such
$L$-morphism. Then the difference between $\eta'$ and $\eta$ is given by
an element of $\Ho^1(W_\R,Z(\widehat{\H}))$; let $\chi$ be the corresponding
quasi-character of $\H(\R)$. Then, if we replace $\eta$ by $\eta'$, the
transfer factor $\Delta_{j_M,\B_M}$ is multiplied by $\chi$ and the stable
characters $S\Theta_{\varphi_H}$, $\varphi_H\in
\Phi_H(\varphi)$, are multiplied by $\chi^{-1}$; hence both sides of the
equality are multiplied by $\chi(\gamma_H)$. It is therefore enough to
prove the proposition for a particular choice of $\eta$.

We choose $\eta$ such that
\[\eta((1,\tau))=((1,\left(\begin{array}{cc}0 & I_{n_1} \\ (-1)^{n_1}I_{n_2} 
& 0 \end{array}\right)),\tau)\]
and that, for every $z\in\C^\times=W_\C$,
\[\eta((1,z))=((1,\left(\begin{array}{cc}
z^{n_2/2}\overline{z}^{-n_2/2}I_{n_1} & 0 \\
0 & z^{-n_1/2}\overline{z}^{n_1/2}I_{n_2}\end{array}\right)),z).\]

We first recall the formulas for $\Phi_M(.,S\Theta_
\varphi)$ and $\Phi_{M_H}(.,S\Theta_{\varphi_H})$. The reference for this
is \cite{A-L2} p 272-274.
\footnote{Note that there is a mistake in this reference. Namely, with the
notations used below, the formula of \cite{A-L2} is correct for
elements in the image of the map $Z\times\tgoth_M(\R)\fl\T_M(\R),
(z,X)\fle z\exp(X)$, but it is not true in general, as claimed in \cite{A-L2},
that the stable discrete series characters vanish outside of the image of this
map. This is not a problem here because the exponential map $\tgoth_M(\R)\fl
\T_M(\R)$ is surjective unless $\M$ is a torus, and, if $\M$ is a torus,
elements in $\M(\R)$ that are not in the image of the exponential map
are also not in $Z(\G)(\R)\G_{der}(\R)$, so all discrete series characters
vanish on these elements.
A formula that is correct in the general case can be found
in section 4 of \cite{GKM}.}

Let $V_\varphi$ be the irreducible representation of $\overline{\G}(\R)$
corresponding to $\varphi$ and $\xi_\varphi$ be the quasi-character
by which $\A_G(\R)^0$ acts on $V_\varphi$. Let $\B_0=\Int(u_M^{-1})(\B)$,
$Z$ be the maximal compact subgroup of the center of $\G(\R)$ and
$\tgoth_G=\Lie(\T_G)$.
Define functions $\rho_G$ and $\Delta_G$ on
$\tgoth_G(\R)$ by :
\[\rho_G=\frac{1}{2}\sum_{\alpha\in\Phi(\T_G,\B_0)}\alpha\]
\[\Delta_G=\prod_{\alpha\in\Phi(\T_G,\B_0)}(e^{\alpha/2}-e^{-\alpha/2})\]
(we use the same notations for characters on $\T_G$ and the linear forms
on $\tgoth_G$ that are defined by differentiating these characters).
Notice that $\T_G(\R)=Z\exp(\tgoth_G(\R))$ (this is a general fact; here,
$\T_G(\R)$ is even equal to $\exp(\tgoth_G(\R))$).
The representation $V_\varphi$ corresponds to a pair $(\zeta_\varphi,
\lambda_\varphi)$, where $\zeta_\varphi$ is a quasi-character of $Z$ and
$\lambda_\varphi$ is a linear form on $\tgoth_G(\C)$, such that :
\begin{bulletlist}
\item $\lambda_\varphi$ is regular dominant;
\item the morphism $Z\times\tgoth_G(\R)\fl\C^\times$, $(z,X)
\fle\zeta_\varphi(z)e^{(\lambda_\varphi-\rho_G)(X)}$, factors through the
surjective morphism $Z\times\tgoth_G(\R)\fl\T_G(\R)$, $(z,X)\fle
z\exp(X)$, and defines a quasi-character on $\T_G(\R)$, whose restriction
to $\A_G(\R)^0$ is $\xi_\varphi$.

\end{bulletlist}
Note that the quasi-character $z\exp(X)\fle\zeta_\varphi(z)e^{
(\lambda_\varphi-\rho_G)(X)}$ on $\T_G(\R)$ is equal to the quasi-character
$\chi_{\varphi,B_0}$ defined in remark \ref{rq:facteur_transfert} (1).
Remember the Weyl character formula :
if $\gamma\in\T_{G,reg}(\R)$ is such that $\gamma=z\exp(X)$, with $z\in Z$ and
$X\in\tgoth_G(\R)$, then
\[Tr(\gamma,V_\varphi)=(-1)^{q(\G)}S\Theta_{\varphi}(\gamma)=
\Delta_G(X)^{-1}\zeta_\varphi(z)\sum_{\omega\in\Omega_G}\det(\omega)
e^{(\omega\lambda_\varphi)(X)}.\]

Let $R$ be a root system whose Weyl group $W(R)$ contains $-1$. Then, to
every pair $(Q^+,R^+)$ such that $R^+\subset R$ and
$Q^+\subset R^{\vee}$ are positive root systems, we can associate an integer
$\overline{c}(Q^+,R^+)$. The definition of $\overline{c}(Q^+,R^+)$ is
recalled in \cite{A-L2} p 273.

Let $R$ be the set of real roots in $\Phi(\T_M,\G)$, 
$R^+=R\cap\Phi(\T_M,\B)$ and $\tgoth_M=\Lie(\T_M)$. If $X$ is a regular
element of $\tgoth_M(\R)$, let 
\[R_X^+=\{\alpha\in R|\alpha(X)>0\},\]
\[\varepsilon_R(X)=(-1)^{|R_X^+\cap(-R^+)|}.\]
If $\nu$ is a linear form on $\tgoth_M(\C)$ such that
$\nu(\alpha^{\vee})\not=0$ for every $\alpha^{\vee}\in R^{\vee}$, let
\[Q_\nu^+=\{\alpha^{\vee}\in R^{\vee}|\nu(\alpha^{\vee})>0\}.\]
Define a function $\Delta_M$ on $\tgoth_M(\R)$ by
\[\Delta_M=\prod_{\alpha\in\Phi(\T_M,\B_M)}(e^{\alpha/2}-e^{-\alpha/2}).\]
As $\Int(u_M)$ sends $\T_G(\C)$ onto $\T_M(\C)$, $\Ad(u_M)$ 
defines an isomorphism $\tgoth_M(\C)^*\iso\tgoth_G(\C)^*$.

Let $\gamma\in\T_{M,reg}(\R)$. If there exist $z\in Z$ and $X\in\tgoth_M(\R)$
such that $\gamma=z\exp(X)$, then (formula (4.8) of \cite{A-L2})
\[\Phi_M(\gamma,S\Theta_\varphi)=(-1)^{q(\G)}
\Delta_M(X)^{-1}\varepsilon_R(X)\zeta_\varphi
(z)\sum_{\omega\in\Omega_G}\det(\omega)\overline{c}(Q^+_{\Ad(u_M)\omega
\lambda},R^+_X)e^{(\Ad(u_M)\omega\lambda)(X)}.\]

There are similar objects, defined by replacing $\G$ by $\H$, etc, and
similar formulas for the functions $\Phi_{M_H}(.,S\Theta_
{\varphi_H})$.

Let $\gamma_H\in\T_{M_H}(\R)$ and $\gamma=j_M(\gamma_H)$; we want to prove
the equality of the proposition. We may assume that $\gamma$ is regular in
$\G$ (because the set of $\gamma_H$ such that $j_M(\gamma_H)$ is $\G$-regular
is dense in $\T_{M_H}(\R)$).
Note that, as $\T_M$ are $\T_{M_H}$ are both
isomorphic to $R_{E/\Q}\Gr_m^r\times\G(\U(1)^{m_1+m_2})$,
the exponential maps $\tgoth_M(\R)\fl\T_M(\R)$ and $\tgoth_{M_H}(\R)
\fl\T_{M_H}(\R)$ are surjective unless $\M$ is a torus (ie $m_1+m_2=0$).
If $\M$ is a torus (so $\M_H$ is also a torus) and $\gamma_H$ is not in the
image of the exponential map, then $c(\gamma_H)=c(\gamma)<0$, so
$\gamma\not\in Z(\G)(\R)\G_{der}(\R)$ and $\gamma_H\not\in Z(\H)(\R)
\H_{der}(\R)$, and all discrete series characters vanish on $\gamma$ and
$\gamma_H$; so the equality of the proposition is obvious.

These remarks show that we may assume that there exists
$X_H\in\tgoth_{M_H}(\R)$ such that $\gamma_H^{-1}=\exp(X_H)$. Then
$\gamma^{-1}=\exp(X)$, with $X=j_M(X_H)\in\tgoth_M(\R)$. As all the real roots
of $\Phi(\T_M,\G)$ are in $j_{M*}(\Phi(\T_{M_H},\H))$, $R$ is equal
to $j_{M*}(R_H)$. Hence
$R_X^+=j_{M*}(R_{H,X_H}^+)$ and $\varepsilon_R(X)=\varepsilon_{R_H}(X_H)$.
On the other hand, using the description of the bijection $\Phi_H(\varphi)
\iso\Omega_*$ given above and the choice of $\eta$,
it is easy to see that, for every $\varphi_H\in
\Phi_H(\varphi)$, $\zeta_{\varphi_H}=\zeta_\varphi$ and $\lambda_{\varphi_H}=
\omega_*(\varphi_H)(\lambda_\varphi)\circ j_M+\rho_H-\rho_G\circ j_M$. As
$\rho_H-\rho_G\circ j_M$ is $\Omega_H$-invariant and vanishes on the elements
of $R_H^{\vee}$, this implies that, for every $\varphi_H\in\Phi_H(
\varphi)$ and $\omega_H\in\Omega_H$, $Q^+_{\Ad(u_M^{-1})\omega_H
\omega_*(\varphi_H)\lambda_\varphi}=j_{M*}(Q^+_{\Ad(u_{M_H}^{-1})\omega_H
\lambda_{\varphi_H}})$. So we get :
\begin{flushleft}$\displaystyle{(-1)^{q(\H)}
\sum_{\varphi_H\in\Phi_H(\varphi)}\det(\omega_*
(\varphi_H))\Phi_{M_H}(\gamma_H^{-1},S\Theta_{\varphi_H})}$\end{flushleft}
\begin{center}$\displaystyle{=\Delta_{M_H}(X_H)^{-1}
\varepsilon_R(X)\sum_{\omega_*\in\Omega_*}\det(\omega_*)\sum_
{\omega_H\in\Omega_H}\det(\omega_H)\overline{c}(Q^+_{\Ad(u_M^{-1})\omega_H
\omega_*\lambda_\varphi},R_X^+)}$\end{center}
\begin{flushright}$\displaystyle{e^{\Ad(u_M^{-1})(\omega_H\omega_*\lambda_
\varphi+\rho_H\circ j_M^{-1}-\rho_G)(X)}}$\end{flushright}
\begin{flushright}$\displaystyle{=(-1)^{q(\G)}\Delta_{M_H}(X_H)^{-1}\Delta_M(X)
e^{\rho_H(X_H)-\rho_G(X)}\Phi_M(\gamma^{-1},S\Theta_\varphi).}$\end{flushright}

To finish the proof, it is enough to show that
\[(-1)^{q(\G)+q(\H)}\Delta_{j_M,\B_M}(\gamma_H,\gamma)=\Delta_{M_H}(X_H)^{-1}
\Delta_M(X)e^{\rho_H(X_H)-\rho_G(X)}.\]
Let $\Phi^+=\Phi(\T_M,\B)$, $\Phi_H^+=j_{M*}(\Phi(\T_{M_H},j_M^*(\B)))$,
$\Phi_M^+=\Phi(\T_M,\B_M)$, $\Phi_{M_H}^+=j_{M*}(\Phi(\T_{M_H},j_M^*(\B_M)))$.
Then
\[\begin{array}{rcl}\Delta_{M_H}(X_H)^{-1}\Delta_M(X)e^{\rho_H(X_H)-\rho_G(X)}
& = & \displaystyle{
\prod_{\alpha\in\Phi_M^+-\Phi_{M_H}^+}(e^{\alpha(X)/2}-e^{-\alpha(X)/2})
\prod_{\alpha\in\Phi^+-\Phi_H^+}e^{-\alpha(X)/2}} \\
& = & \displaystyle{\prod_{\alpha\in\Phi_M^+-\Phi_{M_H}^+}(1-\alpha(\gamma^
{-1}))
\prod_{\alpha\in\Phi^+-(\Phi_H^+\cup\Phi_M^+)}e^{-\alpha(X)/2}}.\end{array}\]
So it is enough to show that
\[\chi_{B_M}(\gamma)=\prod_{\alpha\in\Phi^+-(\Phi_H^+\cup\Phi_M^+)}
e^{-\alpha(X)/2}.\]
Remember that $\chi_{B_M}$ is the quasi-character of $\T_M(\R)$
corresponding to the $1$-cocyle $a_M:W_\R\fl\widehat{\T}_M$ such that
$\eta_M\circ\eta_{B_{M_H}}\circ\widehat{j}_M$ and $\eta_{B_M}.a_M$ are
$\widehat{\M}$-conjugate, where $\B_{M_H}=j_M^*(\B_M)$. So the equality above
is an easy consequence of the definitions of $\eta_{B_M}$ and $\eta_{B_{M_H}}$
and of the choice of $\eta$.

\end{proof}

\begin{remark}\label{rq:M_M_H_regulier} As $\Delta_{j_M,B_M}(\gamma_H,
\gamma)=0$ if $\gamma_H$ is not $(\M,\M_H)$-regular, the right hand side
of the equality of the proposition is non-zero only if
$\gamma_H$ is $(\M,\M_H)$-regular.

\end{remark}

\section{Calculation of certain $\Phi_M(\gamma,\Theta)$}
\label{serie_discrete4}

As before, let $\G=\GU(p,q)$, with $p,q\in\Nat$ such that $p\geq q$ and
$n:=p+q\geq 1$.
Fix $s\in\{1,\dots,q\}$, and set $S=\{1,\dots,s\}$ and $\M=\M_S$
(with notations as in \ref{groupes2}).
The goal of this section is to calculate
$\Phi_M(\gamma,\Theta)$, for $\Theta$ the character of a $L$-packet of the
discrete series of $\G(\R)$ associated to an algebraic representation of
$\G_\C$ and certain $\gamma$ in $\M(\R)$.

The linear part of $\M$ (resp. $\M_s$) is $\Le_S=(R_{E/\Q}\Gr_m)^s$ 
(resp. $\Le_s=R_{E/\Q}\GL_s$), and its Hermitian part is
$\G_s=\GU(p-s,q-s)$. The group $\Le_S$ is a minimal Levi subgroup of
$\Le_s$. The Weyl group $W(\Le_S(\Q),\Le_s(\Q))$ is obviously isomorphic to
$\Sgoth_s$, and we identify these groups; we extend the action of
$\Sgoth_s$ on $\Le_S$ to an action on $\M=\Le_S\times\G_s$, by declaring
that the action is trivial on $\G_s$.

For every $r\in\{1,\dots,q\}$, let $t_r=r(r-n)$.

\begin{proposition}\label{prop:calcul_Phi_M}
Let $E$ be an irreducible algebraic representation of $\G_\C$. Let
$m\in\Z$ be such that the central torus $\Gr_{m,\C}I_n$ of $\G_\C$ acts
on $E$ by multiplication by the character $z\fle z^m$. For every
$r\in\{1,\dots,q\}$, let $t'_r=t_r+m$.
Choose an elliptic Langlands parameter $\varphi:W_\R\fl{}^L\G$ corresponding
to $E$ (seen as an irreducible representation of $\GU(n)(\R)\subset
\GU(n)(\C)\simeq\G(\C)$), and let
$\Theta=(-1)^{q(\G)}S\Theta_\varphi$.
Let $\gamma\in\M(\R)$ be semi-simple elliptic. Write $\gamma=\gamma_l
\gamma_h$, with $\gamma_l=(\lambda_1,\dots,\lambda_s)\in(\C^\times)^s=
\Le_S(\R)$ and $\gamma_h\in\G_s(\R)$.
Then $c(\gamma)=c(\gamma_h)>0$ unless $\M$ is a torus (ie $s=q$). If
$c(\gamma)<0$, then $\Phi_M(\gamma,\Theta)=0$. If $c(\gamma)>0$, then :
\begin{itemize}
\item[(i)] If $c(\gamma)|\lambda_r|^2\geq 1$ for every $r\in\{1,\dots,s\}$, 
then $\Phi_{M}(\gamma,\Theta)$ is equal to
\begin{flushleft}$\displaystyle{2^s\sum_{{S'\subset S}\atop{S'\ni s}}(-1)^
{\dim(\Ar_{M}/\Ar_{M_{S'}})}|W(\Le_S(\Q),\Le_{S'}(\Q))|^{-1}}$\end{flushleft}
\begin{flushright}$\displaystyle{\sum_{\sigma\in
\Sgoth_s}|D_{\M}^{\M_{S'}}(\sigma\gamma)|_\R^{1/2}\delta_{P_{S'}(\R)}
^{1/2}(\sigma\gamma)Tr(\sigma\gamma,R\Gamma(Lie(\N_{S'}),E)_{<t'_r,r\in S'}).}$
\end{flushright}

\item[(ii)] If $0<c(\gamma)|\lambda_r|^2\leq 1$ for every $r\in\{1,\dots,s\}$, 
then $\Phi_{M}(\gamma,\Theta)$ is equal to
\begin{flushleft}$\displaystyle{(-1)^s2^s\sum_{{S'\subset S}\atop{S'\ni s}}(-1)
^{\dim(\Ar_{M}/\Ar_{M_{S'}})}|W(\Le_S(\Q),\Le_{S'}(\Q))|^{-1}}$\end{flushleft}
\begin{flushright}$\displaystyle{\sum_{\sigma\in\Sgoth_s}|D_{\M}^
{\M_{S'}}(\sigma\gamma)|_\R^{1/2}\delta_{P_{S'}(\R)}
^{1/2}(\sigma\gamma)Tr(\sigma\gamma,R\Gamma(Lie(\N_{S'}),E)_{>t'_r,r\in S'}).}$
\end{flushright}

\end{itemize}
\end{proposition}

The notations $R\Gamma(Lie(\N_{S'}),E)_{<t'_r,r\in S'}$ and
$R\Gamma(Lie(\N_{S'}),E)_{>t'_r,r\in S'}$ are those of proposition
\ref{prop:restriction_ponderes_strates}.

\begin{proof} 
Let $\gamma\in\M(\R)$ be semi-simple elliptic.
Use the notations of \ref{serie_discrete3}, in particular of the proof of
proposition \ref{prop:transfert_caracteres}. As $\gamma$ is elliptic in
$\M(\R)$, we may assume that $\gamma\in\T_M(\R)$. The fact that $\Phi_M
(\gamma,\Theta)=0$ if $c(\gamma)<0$ (and that this can happen only if $\M$
is a torus) has already been noted in the proof of proposition
\ref{prop:transfert_caracteres}. So we may assume that $c(\gamma)>0$.

The proofs of (i) and (ii) are similar. Let us prove (ii). Assume that
$c(\gamma)|\lambda_r|^2\leq 1$ for every $r\in S$.
As both sides of the equality we want to prove are continuous functions of
$\gamma$, we may assume that $c(\gamma)|\lambda_r|^2<1$ for every $r\in S$
and that $\gamma$ is regular in $\G$. Let $X$ be an element of
$\tgoth_M(\R)$ such that $\gamma=\exp(X)$ (remember that, as the torus
$\T_M$ is isomorphic to $(R_{E/\Q}\Gr_m)^s\times\G(\U(1)^{n-2s})$,
such a $X$ exists if and only if $c(\gamma)>0$).
Choose an element $\B$ of $\Bo_G(\T_M)$ such that $\B\subset\Pa_S$.
There is a pair $(\zeta,\lambda_G)$ associated to $E$ as in
\ref{serie_discrete3} ($\zeta$ is a quasi-character of $Z$, and $\lambda_G
\in\tgoth_G(\C)^*$). Write
$\lambda=\Ad(u_M)(\lambda_G)\in\tgoth_M(\C)^*$ and $\rho_B=\frac{1}{2}\sum
\limits_{\alpha\in\Phi(\T_M,\B)}\alpha$. Then $\lambda-\rho_B$ is the highest
weight of $E$ relative to $(\T_M,\B)$. 

Let $S'\subset S$ be such that $s\in S'$. We use Kostant's theorem 
(see eg \cite{GHM} \S11) to calculate the trace of $\gamma$ on
$R\Gamma(Lie(\N_{S'}),E)_{>t'_r,r\in S'}$.
Let
$\Omega=W(\T_M(\C),\G(\C))$, $\ell$ be the length function on $\Omega$,
$\Omega_{S'}=W(\T_M(\C),\M_{S'}(\C))$ and $\Phi^+=\Phi(\T_M,\B)$. For every
$\omega\in\Omega$, let $\Phi^+(\omega)=\{\alpha\in\Phi^+|\omega^{-1}\alpha
\in-\Phi^+\}$. Then $\Omega'_{S'}:=\{\omega\in\Omega|\Phi^+(\omega)\subset
\Phi(\T_M,\N_{S'})\}$ is a system of representatives of
$\Omega_{S'}\sous\Omega$. Kostant's theorem says that, for every $k\in\Nat$,
\[\Ho^k(Lie(\N_{S'}),E)\simeq\bigoplus_{\omega\in\Omega'_{S'},\ell(\omega)=k}
V_{\omega(\lambda)-\rho_B},\]
where, for every $\omega\in\Omega$, $V_{\omega(\lambda)-\rho_B}$ is the
algebraic representation of $\M_{S',\C}$ with highest weight
$\omega(\lambda)-\rho_B$ (relative to $(\T_M,\B\cap\M_{S',\C})$). 

For every $r\in\{1,\dots,s\}$, let
\[\varpi_r:\Gr_m\fl\T_M,\lambda\fl\left(\begin{array}{ccc}\lambda I_r & & 0 \\
& I_{n-2r} & \\ 0 & & \lambda^{-1}I_r\end{array}\right);\]
we use the same notation for the morphism $Lie(\Gr_m)\fl\tgoth_M$ obtained
by differentiating $\varpi_r$.
Let $k\in\Nat$. By definition of the truncation,
\[\Ho^k(Lie(\N_{S'}),E)_{>t'_r,r\in S'}\simeq\bigoplus_{\omega}
V_{\omega(\lambda)-\rho_B},\]
where the sum is taken over the set of $\omega\in\Omega'_{S'}$
of length $k$ and such that, for every $r\in S'$,
$\langle\omega(\lambda)-\rho_B,\varpi_r\rangle>t_r$. As $t_r=
\langle-\rho_B,\varpi_r\rangle$ for every $r\in\{1,\dots,s\}$, the last
condition on $\omega$ is equivalent to : 
$\langle\omega(\lambda),\varpi_r\rangle>0$, for every $r\in S'$. 

On the other hand, by the Weyl character formula, for every
$\omega\in\Omega'_{S'}$,
\[Tr(\gamma,V_{\omega(\lambda)-\rho_B})=\Delta_{M_{S'}}(X)^{-1}\sum_{\omega_M
\in\Omega_{S'}}\det(\omega_M)e^{(\omega_M(\omega(\lambda)-\rho_B+\rho_{S'}))(X)
},\]
where $\rho_{S'}=\frac{1}{2}\sum\limits_{\alpha\in\Phi(\T_M,\B\cap\M_{S',\C})}
\alpha$. As $\rho_{S'}-\rho_B$ is invariant by $\Omega_{S'}$ and 
\[e^{(\rho_{S'}-\rho_B)(X)}=\delta_{P_{S'}(\R)}^{-1/2}(\gamma),\]
this formula becomes
\[Tr(\gamma,V_{\omega(\lambda)-\rho_B})=\Delta_{M_{S'}}(X)^{-1}\delta_{P_{S'}
(\R)}^{-1/2}(\gamma)\sum_{\omega_M\in\Omega_{S'}}\det(\omega_M)e^{(\omega_M
\omega(\lambda))(X)}.\]
Hence
\begin{flushleft}$\displaystyle{Tr(\gamma,R\Gamma(Lie(\N_{S'}),E)_{>t'_r,r\in 
S'})}$\end{flushleft}
\begin{flushright}$\displaystyle{=\Delta_{M_{S'}}(X)^{-1}\delta_{P_{S'}(\R)}
^{-1/2}(\gamma)\sum_{\omega_M\in\Omega_{S'}}\sum_{\omega}\det(\omega_M\omega)
e^{(\omega_M\omega(\lambda))(X)},}$\end{flushright}
where the second sum is taken over the set of $\omega\in\Omega'_{S'}$
such that $\langle\omega(\lambda),\varpi_r\rangle>0$ for every $r\in S'$.
As the $\varpi_r$, $r\in S'$, are invariant by $\Omega_{S'}$,
for every $\omega_M\in\Omega_{S'}$, $\omega\in\Omega'_{S'}$ and $r\in S'$,
$\langle\omega(\lambda),\varpi_r\rangle=\langle\omega_M\omega(\lambda),
\varpi_r\rangle$. Hence :
\[Tr(\gamma,R\Gamma(Lie(\N_{S'}),E)_{>t'_r,r\in S'})=\Delta_{M_{S'}}(X)^{-1}
\delta_{P_{S'}(\R)}^{-1/2}(\gamma)\sum_\omega\det(\omega)e^{(\omega(\lambda))
(X)},\]
where the sum is taken over the set of $\omega\in\Omega$ such that
$\langle\omega(\lambda),\varpi_r\rangle>0$ for every $r\in S'$.

Moreover,
\[|D_{\M}^{\M_{S'}}(\gamma)|^{1/2}=|\Delta_{M_{S'}}(X)||\Delta_M
(X)|^{-1}.\]
But all the roots of $\T_M$ in $Lie(\M_{S'})/Lie(\M)$ are complex, so :
\[\Delta_{M_{S'}}(X)\Delta_M(X)^{-1}=\prod_{{\alpha\in\Phi(\T_M,
Lie(\M_{S'})/Lie(\M))}\atop{\alpha>0}}(e^{\alpha(X)/2}-e^{-\alpha(X)/2})\in
\R^+,\]
and
\[|D_{\M}^{\M_{S'}}(\gamma)|^{1/2}=\Delta_{M_{S'}}(X)\Delta_{M}(X)
^{-1}.\]
Finally,
\[|D_{\M}^{\M_{S'}}(\gamma)|^{1/2}\delta_{P_{S'}(\R)}^{1/2}(\gamma)Tr(\gamma,
R\Gamma(Lie(\N_{S'}),E)_{>t'_r,r\in S'})=\Delta_M(X)^{-1}\sum_\omega\det
(\omega)e^{(\omega(\lambda))(X)},\]
where the sum is taken as before on the set of $\omega\in\Omega$ such that
$\langle\omega(\lambda),\varpi_r\rangle>0$ for every $r\in S'$.

The action of the group $\Sgoth_s$ on $\T_M(\C)$ gives an injective morphism
$\Sgoth_s\fl\Omega$. Use this morphism to see $\Sgoth_s$ as a subgroup of 
$\Omega$. For every $\sigma\in\Sgoth_s$, $\det(\sigma)=1$, and the function
$\Delta_M$ is invariant by $\Sgoth_s$. Hence :
\begin{flushleft}$\displaystyle{\sum_{\sigma\in\Sgoth_s}|D_{\M}^{\M_{S'}}(
\sigma\gamma)|^{1/2}\delta_{P_{S'}(\R)}(\sigma\gamma)^{1/2}Tr(\sigma\gamma,
R\Gamma(Lie(\N_{S'}),E)_{>t'_r,r\in S'})}$\end{flushleft}
\begin{flushright}$\displaystyle{=\Delta_M(X)^{-1}\sum_{\omega\in\Omega}
\det(\omega)e^{(\omega(\lambda))(X)}|\{\sigma\in\Sgoth_s|\langle\sigma\omega
(\lambda),\varpi_r\rangle>0\mbox{ for every }r\in S'\}|.}$\end{flushright}

We now use the formula of \cite{A-L2} \ 272-274
(recalled in the proof of proposition \ref{prop:transfert_caracteres})
to calculate $\Phi_M(\gamma,\Theta)$. Let $R$ be the set of real roots
in $\Phi(\T_M,\G)$. For every $r\in\{1,\dots,s\}$, let
\[\alpha_r:\T_M\simeq(R_{E/\Q}\Gr_m)^s\times\G(\U(1)^{n-2s})\fl\Gr_m,
((\lambda_1,\dots,\lambda_s),g)\fle c(g)\lambda_r\overline{\lambda}_r.\]
Then $R=\{\pm\alpha_1,\dots,\pm\alpha_s\}$, 
$R^+:=R\cap\Phi^+=\{\alpha_1,\dots,\alpha_s\}$, and, for every
$r\in\{1,\dots,s\}$, the coroot $\alpha_r^{\vee}$ is the morphism
\[\Gr_m\fl\T_M,\lambda\fle ((\underbrace{1,\dots,1}_{r-1},\lambda,\underbrace{
1,\dots,1}_{s-r}),1).\]
Note that $\varpi_r=\alpha_1^\vee+\dots+\alpha_r^\vee$.
As $c(\gamma)|\lambda_r|^2\in ]0,1[$ for every $r\in\{1,\dots,s\}$,
$R_X^+=\{-\alpha_1,\dots,-\alpha_s\}=-R^+$ and $\varepsilon_R(X)=(-1)^s$.
Let $Q^+$ be a positive root system in $R^\vee$. If $Q^+\not=\{
\alpha_1^\vee,\dots,\alpha_s^\vee\}$, then
$\overline{c}(Q^+,R_X^+)=0$ by property (ii) of the function
$\overline{c}$ of \cite{A-L2} p 273. Suppose that $Q^+=\{\alpha_1^\vee,
\dots,\alpha_s^\vee\}$. Note that $R$ is the product of the root systems
$\{\pm\alpha_r\}$, $1\leq r\leq s$. Hence
\[\overline{c}(Q^+,R_X^+)=\prod_{r=1}^s\overline{c}(\{\alpha_r^\vee\},
\{-\alpha_r\}).\]
But it is easy to see that $\overline{c}(\{\alpha_r^\vee\},\{-\alpha_r\})=2$ 
for every $r\in\{1,\dots,s\}$ (by property (iii) of \cite{A-L2} p 273).
Hence $\overline{c}(Q^+,R_X^+)=2^s$. Finally, we find
\[\Phi_M(\gamma,\Theta)=(-1)^s2^s\Delta_M(X)^{-1}\sum_\omega\det(\omega)
e^{(\omega(\lambda))(X)},\]
where the sum is taken over the set of $\omega\in\Omega$ such that
$\langle\omega(\lambda),\alpha_r^\vee\rangle>0$ for every $r\in\{1,\dots,s\}$.

To finish the proof, it is enough to show that, if $\omega\in\Omega$ is fixed,
then
\[\sum_{{S'\subset S}\atop{s\in S'}}(-1)^{|S|-|S'|}|W(\Le_S(\Q),\Le_{S'}(\Q))|
^{-1}|\{\sigma\in\Sgoth_s|\langle\sigma\omega(\lambda),\varpi_r\rangle>0
\mbox{ for every }r\in S'\}|\]
is equal to $1$ if $\langle\omega(\lambda),\alpha_r^\vee\rangle>0$ for every
$r\in\{1,\dots,s\}$ and to $0$ otherwise.
This is proved in lemma \ref{lemme_combinatoire_idiot} below.

\end{proof}

Let $n\in\Nat^*$.
Let $S\subset\{1,\dots,n\}$. If $\lambda=(\lambda_1,\dots,\lambda_n)\in
\R^n$, we say that $\lambda>_S0$ if, for every $r\in S$, $\lambda_1+\dots+
\lambda_r>0$, and we write
\[\Sgoth_S(\lambda)=\{\sigma\in\Sgoth_n|\sigma(\lambda)>_S0\}.\]
If $S=\{r_1,\dots,r_k\}$ with $r_1<\dots<r_k$, write
\[w_S=r_1!\prod_{i=1}^{k-1}(r_{i+1}-r_i)!.\]

\begin{lemma}\label{lemme_combinatoire_idiot} Let $\lambda=(\lambda_1,
\dots,\lambda_n)\in\R^n$. Then
\[\sum_{{S\subset\{1,\dots,n\}}\atop{S\ni n}
}(-1)^{|S|}w_S^{-1}|\Sgoth_S(\lambda)|=\left\{
\begin{array}{ll}(-1)^n & \mbox{ if }\lambda_r>0\mbox{ for every }r\in\{1,
\dots,n\} \\ 0 & \mbox{ otherwise}\end{array}\right..\]

\end{lemma}

\begin{proof} First we reformulate the problem.
Let $\lambda=(\lambda_1,\dots,\lambda_n)\in\R^n$. Write $\lambda>0$ if
$\lambda_1>0,\lambda_1+\lambda_2>0,\dots,\lambda_1+\dots+\lambda_n>0$. 
For every $I\subset\{1,\dots,n\}$, let $s_I(\lambda)=\sum\limits_{i\in I}
\lambda_i$. Let $\Par_{ord}(n)$ be the set of ordered partitions of
$\{1,\dots,n\}$. For every $p=(I_1,\dots,I_k)\in\Par_{ord}(n)$, set
$|p|=k$ and
$\lambda_p=(s_{I_1}(\lambda),\dots,s_{I_k}(\lambda))\in\R^k$. Let
\[\Par_{ord}(\lambda)=\{p\in\Par_{ord}(n)|\lambda_p>0\}.\]
Then it is obvious that :
\[\sum_{{S\subset\{1,\dots,n\}}\atop{S\ni n}}(-1)^{|S|}w_S^{-1}|\Sgoth_S(
\lambda)|=\sum_{p\in\Par_{ord}(\lambda)}(-1)^{|p|}.\]

We show the lemma by induction on the pair $(n,|\Par_{ord}(\lambda)|)$
(we use the lexicographical ordering). If $n=1$ or if $\Par_{ord}(\lambda)=
\varnothing$
(ie $\lambda_1+\dots+\lambda_n\leq 0$), the result is obvious.
Assume that $n\geq 2$, that $\Par_{ord}(\lambda)\not=\varnothing$ and
that the result is known for :
\begin{itemize}
\item[-] all the elements of $\R^m$ if $1\leq m<n$;
\item[-] all the $\lambda'\in\R^n$ such that $|\Par_{ord}(\lambda')|
<|\Par_{ord}(\lambda)|$.

\end{itemize}

Let $\Ens$ be the set of $I\subset\{1,\dots,n\}$ such that there exists
$p=(I_1,\dots,I_k)\in\Par_{ord}(\lambda)$ with $I_1=I$; $\Ens$ is non-empty
because $\Par_{ord}(\lambda)$ is non-empty. Let $\varepsilon$ be the minimum
of the $s_I(\lambda)/|I|$, for $I\in\Ens$. Let $I\in\Ens$ be an element with
minimal cardinality among the elements $J$ of $\Ens$ such that
$s_J(\lambda)=\varepsilon|J|$. Define $\lambda'=(\lambda'_1,\dots,
\lambda'_n)\in\R^n$ by :
\[\lambda'_i=\left\{\begin{array}{ll}\lambda_i & \mbox{ if }i\not\in I \\
\lambda_i-\varepsilon & \mbox{ if }i\in I\end{array}\right..\]
Let $\Par'$ be the set of $p=(I_1,\dots,I_k)\in\Par_{ord}(\lambda)$ such that
there exists $r\in\{1,\dots,k\}$ with $I=I_1\cup\dots\cup I_r$. Then
$\Par'\not=\varnothing$.

It is obvious that $\Par_{ord}(\lambda')\subset\Par_{ord}(\lambda)-\Par'$
(because $s_I(\lambda')=0$). Let us show that $\Par_{ord}(\lambda')=
\Par_{ord}(\lambda)-\Par'$. Let $p=(I_1,\dots,I_k)\in\Par_{ord}(\lambda)-
\Par'$, and let us show that $p\in\Par_{ord}(\lambda')$. It is enough to
show that $s_{I_1}(\lambda')>0$, because $(I_1\cup\dots\cup I_r,I_{r+1},
\dots,I_k)\in\Par_{ord}(\lambda)-\Par'$ for every $r\in\{1,\dots,k\}$.
By definition of $\lambda'$, $s_{I_1}(\lambda')=s_{I_1}(\lambda)-
\varepsilon|I\cap I_1|$. If $s_{I_1}(\lambda)>\varepsilon|I_1|$, then
$s_{I_1}(\lambda')>\varepsilon (|I_1|-|I\cap I_1|)\geq 0$. If $s_{I_1}(\lambda)
=\varepsilon|I_1|$, then $|I_1|\geq |I|$ by definition of $I$ and
$I_1\not=I$ because $p\not\in\Par'$, so $I_1\not\subset I$, and
$s_{I_1}(\lambda')=\varepsilon(|I_1|-|I\cap I_1|)>0$.

As $s_I(\lambda')=0$, there exists $i\in I$ such that $\lambda'_i\leq 0$.
By the induction hypothesis, $\sum\limits_{p\in\Par_{ord}(\lambda')}
(-1)^{|p|}=0$. Hence $\sum\limits_{p\in\Par_{ord}(\lambda)}(-1)^{|p|}=
\sum\limits_{p\in\Par'}(-1)^{|p|}$. As the equality of the lemma does not
change if the $\lambda_i$ are permuted, we may assume that there exists
$m\in\{1,\dots,n\}$ such that $I=\{1,\dots,m\}$. Assume first that
$m<n$. Let $\mu=(\lambda_1,\dots,\lambda_m)$ and $\nu=(\lambda_{m+1},
\dots,\lambda_n)$. Identify $\{m+1,\dots,n\}$ to $\{1,\dots,n-m\}$ by
the map $k\fle k-m$, and define a map
$\varphi:\Par'\fl\Par_{ord}(m)\times\Par_{ord}(n-m)$ as follows :
if $p=(I_1,\dots,I_k)\in\Par'$ and if $r\in\{1,\dots,k\}$ is such that
$I_1\cup\dots\cup I_r=I$, set $\varphi(p)=((I_1,\dots,I_r),
(I_{r+1},\dots,I_k))$. The map $\varphi$ is clearly injective.
Let us show that the image of $\varphi$ is $\Par_{ord}(\mu)\times\Par_{ord}
(\nu)$.
The inclusion $\Par_{ord}(\mu)\times\Par_{ord}(\nu)\subset\varphi(\Par')$ is
obvious. Let $p=(I_1,\dots,I_k)\in\Par'$, and let $r\in\{1,\dots,k\}$ be
such that $I=I_1\cup\dots\cup I_r$. We want to show that
$\varphi(p)\in\Par_{ord}(\mu)\times\Par_{ord}(\nu)$, ie that, for every
$s\in\{r+1,\dots,k\}$, $s_{I_{r+1}\cup\dots\cup I_s}(\lambda)>0$.
After replacing $p$ by
$(I_1\cup\dots\cup I_r,I_{r+1}\cup\dots\cup I_s,I_{s+1},\dots,I_k)$, we
may assume that $r=1$ (hence $I=I_1$) and $s=2$. Then
\[s_{I_2}(\lambda)=s_{I\cup I_2}(\lambda)-s_{I}(\lambda)=s_{I\cup I_2}(\lambda)
-\varepsilon|I|\geq \varepsilon|I\cup I_2|-\varepsilon|I|>0.\]
Finally :
\[\sum_{p\in\Par'}(-1)^{|p|}=\left(\sum_{p\in\Par_{ord}(\mu)}(-1)^{|p|}\right)
\left(\sum_{p\in\Par_{ord}(\nu)}(-1)^{|p|}\right).\]
Hence the conclusion of the lemma is a consequence of the induction
hypothesis, applied to $\mu$ and $\nu$.

We still have to treat the case $I=\{1,\dots,n\}$. Let us show that there
is no partition $\{I_1,I_2\}$ of $\{1,\dots,n\}$ such that
$s_{I_1}(\lambda)>0$ and $s_{I_2}(\lambda)>0$ (in particular, there exists
at least one $i$ such that $\lambda_i\leq 0$). If such a partition existed,
then, by definition of $I$, we would have inequalities
$s_{I_1}(\lambda)>\varepsilon|I_1|$ and $s_{I_2}(\lambda)>\varepsilon|I_2|$,
hence $s_I(\lambda)>\varepsilon|I|$; but that is impossible.
Let $\Par(n)$ be the set of (unordered) partitions of $\{1,\dots,n\}$.
For every $q=\{I_\alpha,\alpha\in A\}\in\Par(n)$, write $|q|=|A|$.
Let $q=\{I_\alpha,\alpha\in A\}\in\Par(n)$.
By lemma \ref{sous_lemme_combinatoire_idiot}, applied to $(s_{I_{\alpha_1}}
(\lambda),\dots,s_{I_{\alpha_k}}(\lambda))$ for a numbering
$(\alpha_1,\dots,\alpha_k)$ of $A$ (the choice of numbering is unimportant),
there are exactly $(|q|-1)!$ way to order $q$ in order to get an element of
$\Par_{ord}(\lambda)$. Hence
\[\sum_{p\in\Par_{ord}(\lambda)}(-1)^{|p|}=\sum_{q\in\Par(n)}(-1)^
{|q|}(|q|-1)!.\]
If $q'$ is a partition of $\{1,\dots,n-1\}$, we can associate to it a partition
$q$ of $\{1,\dots,n\}$ in one of the following ways :
\begin{itemize}
\item[(i)] Adding $n$ to one of the sets of $q'$. There are $|q'|$ ways of
doing this, and we get $|q|=|q'|$.
\item[(ii)] Adding to $q'$ the set $\{n\}$. There is only one way of doing
this, and we get $|q|=|q'|+1$.

\end{itemize}
We get every partition of $\{1,\dots,n\}$ in this way, and we get it only once.
Hence (remember that $n\geq 2$) :
\[\sum_{q\in\Par(n)}(-1)^{|q|}(|q|-1)!=\sum_{q'\in\Par(n-1)}(-1)^{|q'|}
|q'|(|q'|-1)!+\sum_{q'\in\Par(n-1)}(-1)^{|q'|+1}(|q'|)!=0.\]

\end{proof}

In the lemma below, $\Sgoth_n$ acts on $\R^n$ in the usual way (permuting
the coordinates).

\begin{lemma}\label{sous_lemme_combinatoire_idiot} Let $\lambda=(\lambda_1,
\dots,\lambda_n)\in\R^n$. Assume that $\lambda_1+\dots+\lambda_n>0$ and that
there is no partition $\{I_1,I_2\}$ of $\{1,\dots,n\}$ such that
$s_{I_1}(\lambda)>0$ and $s_{I_2}(\lambda)>0$. Then
\[|\{\sigma\in\Sgoth_n|\sigma(\lambda)>0\}|=(n-1)!.\]

\end{lemma}

\begin{proof} Let $\Sgoth(\lambda)=\{\sigma\in\Sgoth_n|\sigma(\lambda)>0
\}$. Let $\tau\in\Sgoth_n$ be the permutation that sends an element
$i$ of $\{1,\dots,n-1\}$ to $i+1$, and sends $n$ to $1$.
Let us show that there exists a unique $k\in\{1,\dots,n\}$ such that
$\tau^k\in\Sgoth(\lambda)$. Let $s=\min\{\lambda_1+\dots\lambda_l,1\leq l\leq
n\}$. Let $k$ be the biggest element of $\{1,\dots,n\}$ such that
$\lambda_1+\dots+\lambda_k=s$. If $l\in\{k+1,\dots,n\}$, then $\lambda_1+\dots+
\lambda_l>\lambda_1+\dots+\lambda_k$, hence $\lambda_{k+1}+\dots+\lambda_l>0$.
If $l\in\{1,\dots,k\}$, then 
\[\begin{array}{rcl}\lambda_{k+1}+\dots+\lambda_n+\lambda_1+\dots+\lambda_l &
= & (\lambda_1+\dots+\lambda_n)-(\lambda_1+\dots+\lambda_k)+(\lambda_1+\dots+
\lambda_l) \\
 & > & -(\lambda_1+\dots+\lambda_k)+(\lambda_1+\dots+\lambda_l) \\
& \geq & 0.\end{array}\]
This proves that $\tau^k(\lambda)=(\lambda_{k+1},\dots,
\lambda_n,\lambda_1,\dots,\lambda_k)>0$.
Suppose that there exists $k,l\in\{1,\dots,n\}$ such that
$k<l$, $\tau^k(\lambda)>0$ and $\tau^l(\lambda)>0$. Let
$I_1=\{k+1,\dots,l\}$ and $I_2=\{1,\dots,n\}-I_1$. Then $s_{I_1}(\lambda)=
\lambda_{k+1}+\dots+\lambda_l>0$ because $\tau^k(\lambda)>0$, and
$s_{I_2}(\lambda)=\lambda_{l+1}+\dots+\lambda_n+\lambda_1+\dots+\lambda_k>0$
because $\tau^l(\lambda)>0$. This contradicts the assumption on $\lambda$.

Applying the above reasoning to $\sigma(\lambda)$, for $\sigma\in\Sgoth_n$,
we see that $\Sgoth_n$ is the disjoint union of the subsets
$\tau^k\Sgoth(\lambda)$, $1\leq k\leq n$. This implies the conclusion of
the lemma.

\end{proof}

\chapter{Orbital integrals at $p$}
\label{partie_en_p}

\section{A Satake transform calculation (after Kottwitz)}
\label{partie_en_p1}

\begin{lemma}\label{def_r_mu}(cf \cite{K-SVTOI} 2.1.2, \cite{K-SVLR} p 193)
Let $F$ be a local or global field and $\G$ be a connected reductive
algebraic group over $F$.
For every cocharacter $\mu:\Gr_{m,F}\fl\G$, there exists a representation
$r_\mu$ of ${}^L\G(=\widehat{\G}\rtimes W_F$), unique up to isomorphism,
satisfying the following conditions :
\begin{itemize}
\item[(a)] The restriction of $r_\mu$ to $\widehat{\G}$ is
irreducible algebraic of highest weight $\mu$.
\item[(b)] For every $\Gal(\overline{F}/F)$-fixed splitting of
$\widehat{\G}$, the group $W_F$, embedded in ${}^L\G$ by the section
associated to the splitting, acts trivially on the highest weight
subspace of $r_\mu$ (determined by the same splitting).
\index{rmu@$r_\mu$}

\end{itemize}
\end{lemma}

Let $p$ be a prime number, $\overline{\Q}_p$ be an algebraic closure of
$\Q_p$, $\Q_p^{ur}$ be the maximal unramified extension of $\Q_p$ in
$\overline{\Q}_p$, $F\subset\overline{\Q}_p$ be a finite unramified extension
of $\Q_p$, $W_F$ be the Weyl group of $F$,
$\varpi_F$ be a uniformizer of $F$, $n=[F:\Q_p]$. The cardinality of
the residual field of $F$ is $p^n$.

Let $\G$ be a connected reductive algebraic group over $F$, and assume
that $\G$ is unramified. Fix a hyperspecial maximal compact subgroup
$\K$ of $\G(F)$, and let
$\Hecke=\Hecke(\G(F),\K):=C_c^\infty(\K\sous\G(F)/\K)$
be the associated Hecke algebra.
For every cocharacter $\mu:\Gr_{m,F}\fl\G$ of $\G$, let
\[f_\mu=\frac{\ungras_{\K\mu(\varpi_F^{-1})\K}}{\vol(\K)}\in\Hecke.\]
\index{fmu@$f_\mu$}

In this section, the $L$-group of $\G$ will be
${}^L\G=\widehat{\G}\rtimes W_F$. Let $\varphi\fle
\pi_\varphi$ be the bijection between the set of equivalence classes
of admissible unramified morphisms $\varphi:W_F\fl{}^L\G$ and the set of
isomorphism classes of spherical representations of $\G(F)$.

\begin{theorem}\label{th:transformee_de_Satake}(\cite{K-SVTOI} 2.1.3)
Let $\mu:\Gr_{m,F}\fl\G$ be such that the weights of the representation
$\Ad\circ\mu:\Gr_{m,F}\fl Lie(\G_{\overline{\Q}_p})$ are in
$\{-1,0,1\}$. Fix a maximal torus $\T$ of $\G$ such that $\mu$ factors
through $\T$, and an ordering on the roots of $\T$ in $\G$ such that
$\mu$ is dominant. Let $\rho$ be half the sum of the positive roots.

Then, for every admissible unramified morphism
$\varphi:W_F\fl{}^L\G$, 
\[\Tr(\pi_\varphi(f_\mu))=p^{n<\rho,\mu>}\Tr(r_{-\mu}(\varphi(\Phi_F))),\]
where $r_{-\mu}$ is the representation defined in lemma \ref{def_r_mu}
above and $\Phi_F\in W(\Q_p^{ur}/F)$ is the geometric Frobenius.

\end{theorem}

\begin{remark} In 2.1.3 of \cite{K-SVTOI}, the arithmetic Frobenius
is used instead of the geometric Frobenius. The difference comes
from the fact that we use here the other normalization of the
class field isomorhism (cf \cite{K-SVLR} p 193).

\end{remark}

\section{Explicit calculations for unitary groups}
\label{partie_en_p2}

This section contains explicit descriptions of the Satake isomorphism,
the base change map, the transfer map and the twisted transfer (or unstable
base change) map for the spherical Hecke algebras of
the unitary groups of \ref{groupes1}. These calculations will be useful
when proving proposition \ref{prop:identite_en_p} and in  the applications
of chapter \ref{applications} and section \ref{GL_n_applications4}.

Let $p$ be a prime number, and let $\overline{\Q}_p$ and $\Q_p^{ur}$
be as in \ref{partie_en_p1}.
Remember that, if $\G$ and $\H$ are unramified
groups over $\Q_p$ and if $\eta:{}^L\H:=\widehat{\H}\rtimes W_{\Q_p}\fl
{}^L\G:=\widehat{\G}\rtimes W_{\Q_p}$ is an unramified $L$-morphism
(ie a $L$-morphism that comes by inflation from a morphism
$\widehat{\H}\rtimes W(\Q_p^{ur}/\Q_p)\fl\widehat{\G}\rtimes
W(\Q_p^{ur}/\Q_p)$), then it induces a morphism of algebras
\index{unramified $L$-morphism}
$b_\eta:\Hecke_\H\fl\Hecke_\G$, where $\Hecke_\H$ (resp. $\Hecke_\G$) is
the spherical Hecke algebra of $\H$ (resp. $\G$). This construction is
recalled in more detail in \ref{lemme_fondamental2} (just before lemma
\ref{lemme:DL2}).
\index{beta@$b_\eta$}

Let $E=\Q[\sqrt{-b}]$
be an imaginary quadratic extension where $p$ is unramified,
and fix a place $\wp$ of $E$ above $p$ (ie an embedding $E\subset\Q_p^{ur}$).
Let $n_1,\dots,n_r\in\Nat^ *$ and let $J_1\in\GL_{n_1}(\Z),\dots,J_r\in
\GL_{n_r}(\Z)$ be symmetric matrices that are antidiagonal (ie in the subset
$\left(\begin{array}{ccc}0 & & * \\ & \begin{turn}{45}\large\ldots
\end{turn}& \\ * & & 0\end{array}\right)$). Let
$q_i$ be the floor (integral part) of $n_i/2$, $n=n_1+\dots+n_r$ and
$\G=\G(\U(J_1)\times\dots\times\U(J_r))$.
Then the group $\G$ is unramified over $\Q_p$. Hence $\G$ extends
to a reductive
group scheme over $\Z_p$ (ie to a group scheme over $\Z_p$ with connected
geometric fibers whose special fiber is a reductive group over $\Fi_p$);
we gave an example of such a group scheme in remark
\ref{rq:groupes_sur_Z}.
We will still denote this group scheme by $\G$.
In this section, the $L$-group of $\G$ will be the $L$-group over $\Q_p$,
ie ${}^L\G=\G\rtimes W_{\Q_p}$.

\subsection*{Satake isomorphism}
\index{Satake isomorphism}

Let $L\subset\Q_p^{ur}$ be an unramified extension of $\Q_p$.
Let $\K_L=\G(\Of_L)$; it is a hyperspecial maximal
compact subgroup of $\G(L)$. We calculate the Satake
isomorphism for $\Hecke(\G(L),\K_L)$.

Suppose first that $\G$ splits over $L$, ie $L\supset E$
(if $L=\Q_p$, this means that
$p$ splits in $E$). Then $\G_L\simeq\Gr_{m,L}\times\GL_{n_1,L}
\times\dots\times\GL_{n_r,L}$.
For every $i\in\{1,\dots,r\}$, let $\T_i$ be the diagonal torus of $\GU(J_i)$.
Identify $\T_{i,L}$ with the torus $\Gr_{m,L}\times\Gr_{m,L}^{n_i}$ by the
isomorphism
\[g=diag(\lambda_1,\dots,\lambda_{n_i})\fle (c(g),(u(\lambda_1),\dots,
u(\lambda_{n_i}))),\]
where $u$ is the morphism $L\otimes_{\Q_p}E\fl L,x\otimes 1+y\otimes\sqrt{-b}
\fle x+y\sqrt{-b}$.
A split maximal torus of $\G_L$ is the diagonal torus
$\T_{G,L}$, where
\[\T_G=\{(g_1,\dots,g_r)\in\T_1\times\dots\times\T_r|c(g_1)=\dots=c(g_r)\}.\]
The above isomorphisms give an isomorphism $\T_{G,L}\simeq\Gr_{m,L}
\times\Gr_{m,L}^n$.
Let $\Omega_G(L)=W(\T_G(L),\G(L))$ be the relative Weyl group of
$\T_{G,L}$ (as $\G$ splits over $L$, this group is actually equal
to the absolute Weyl group). Then
$\Omega_G(L)\simeq\Sgoth_{n_1}\times\dots\times\Sgoth_{n_r}$.
The Satake isomorphism is an isomorphism
\[\Hecke(\G(L),\K_L)\iso\C[X_*(\T_G)]^{\Omega_G(L)}.\]

There is an isomorphism
\[\C[X_*(\T_G)]\simeq\C[X^{\pm 1},X_{i,j}^{\pm 1},1\leq i\leq r,1\leq j
\leq n_i]\]
induced by the isomorphism $\T_{G,L}\simeq\Gr_{m,L}\times\Gr_{m,L}^n$ defined
above. Explicitely
\begin{bulletlist}
\item $X$ corresponds to the cocharacter
\[\lambda\fle\left(\frac{\lambda+1}{2}\otimes 1+\frac{1-\lambda}{2\sqrt{-b}}
\otimes\sqrt{-b}\right)(I_{n_1},\dots,I_{n_r}).\]
\item Let $i\in\{1,\dots,r\}$ and $s\in\{1,\dots,q_i\}\cup\{n_i+1-
q_i,\dots,n_i\}$. Then $X_{i,s}$ corresponds to the cocharacter
\[\lambda\fle (I_1,\dots,I_{i-1},diag(a_1(\lambda),\dots,a_{n_i}(\lambda)),
I_{i+1},\dots,I_r),\]
with :
\[a_j(\lambda)=\left\{\begin{array}{ll}\frac{\lambda+1}{2}\otimes 1+\frac{
\lambda-1}{2\sqrt{-b}}\otimes\sqrt{-b} & \mbox{ if }j=s \\
\frac{\lambda^{-1}+1}{2}\otimes 1+\frac{1-\lambda^{-1}}{2\sqrt{-b}}\otimes
\sqrt{-b} & \mbox{ if }s=n+1-j \\
1 & \mbox{ otherwise}\end{array}\right..\]
\item If $i\in\{1,\dots,r\}$ is such that $n_i$ is odd,
then $X_{i,(n_i+1)/2}$ corresponds to the cocharacter
\[\lambda\fle(I_1,\dots,I_{i-1},\left(\begin{array}{ccc}I_{q_i} & & 0 \\
& \frac{\lambda+\lambda^{-1}}{2}\otimes 1+\frac{\lambda-\lambda^{-1}}{2
\sqrt{-b}}\otimes\sqrt{-b} & \\
0 & & I_{q_i}\end{array}\right),I_{i+1},\dots,I_r).\]

\end{bulletlist}

We get an isomorphism
\[\C[X_*(\T_G)]^{\Omega_G(L)}\simeq\C[X^{\pm 1}]\otimes\C[X_{i,j}^{\pm 1},
1\leq i\leq r,1\leq j\leq n_i]^{\Sgoth_{n_1}\times\dots\times\Sgoth_{n_r}},\]
where $\Sgoth_{n_i}$ acts by permutations on $X_{i,1},\dots,X_{i,n_i}$ 
and trivially on the $X_{i',j}$ if $i'\not=i$.

Suppose now that $\G$ does not split over $L$ (this implies that
$p$ is inert in $E$).
For every $i\in\{1,\dots,r\}$, a maximal split torus of
$\GU(J_i)_{L}$ is $\Se_{i,L}$, where
\[\Se_i=\{diag(\lambda\lambda_1,\dots,\lambda\lambda_{q_i},\lambda_{q_i}^{-1},
\dots,\lambda_1^{-1}),\lambda,\lambda_1,\dots,\lambda_{q_i}
\in\Gr_{m,\Q_p}\}\simeq\Gr_{m,\Q_p}^{q_i+1}\]
if $n_i$ is even, and
\[\Se_i=\{diag(\lambda\lambda_1,\dots,\lambda\lambda_{q_i},\lambda,\lambda
\lambda_{q_i}^{-1},\dots,\lambda\lambda_1^{-1}),\lambda,\lambda_1,\dots,
\lambda_{q_i}\in\Gr_{m,\Q_p}\}\simeq\Gr_{m,\Q_p}^{q_i+1}\]
if $n_i$ is odd.
A maximal split torus of $\G_{L}$ is $\Se_{G,L}$, where
\[\Se_G=\{(g_1,\dots,g_r)\in\Se_1\times\dots\times\Se_r|c(g_1)=\dots=c(g_r)\}^0
\simeq\Gr_{m,\Q_p}^{q_1+\dots+q_r+1}.\]

Let $\Omega_G(L)=W(\Se_G(L),\G(L))$ be the relative Weyl group of
$\Se_G(L)$. Then
$\Omega_G(L)\simeq\Omega_1\times\dots\times\Omega_r$, where
$\Omega_i$ is the subgroup of $\Sgoth_{n_i}$ generated by the transposition
$(1,n_i)$ and by the image of the morphism
\[\Sgoth_{q_i}\fl\Sgoth_{n_i},\quad \sigma\fle\left(\tau:j\fle\left\{
\begin{array}{ll}\sigma(j) & \mbox{ if }1\leq j\leq q_i \\
j & \mbox{ if }q_i+1\leq j\leq n_i-q_i \\
n_i+1-\sigma(n_i+1-j) & \mbox{ if }n_i+1-q_i\leq j\leq n_i\end{array}\right.
\right).\]
Hence $\Omega_i$ is isomorphic to the semi-direct product
$\{\pm 1\}^{q_i}\rtimes\Sgoth_{q_i}$, where $\Sgoth_{q_i}$ acts on
$\{\pm 1\}^{q_i}$ by $(\sigma,(
\varepsilon_1,\dots,\varepsilon_{q_i}))\fle (\varepsilon_{\sigma^{-1}(1)},
\dots,\varepsilon_{\sigma^{-1}(q_i)})$. 
The Satake isomorphism is an isomorphism
\[\Hecke(\G(L),\K_L)\iso\C[X_*(\Se_G)]^{\Omega_G(L)}.\]

Assume that $n_i$ is even. Then there is an isomorphism
\[\C[X_*(\Se_i)]\simeq\C[{X'_i}^{\pm 1},X_{i,1}^{\pm 1},\dots,X_{i,q_i}^
{\pm 1}]\]
that sends $X'_i$ to the cocharacter
\[\lambda\fle\left(\begin{array}{cc}\lambda I_{q_i} & 0 \\ 0 & I_{q_i}
\end{array}\right)\]
and $X_{i,s}$, $1\leq s\leq q_i$, to the cocharacter
\[\lambda\fle diag(a_1(\lambda),\dots,a_{n_i}(\lambda)),\]
with
\[a_j(\lambda)=\left\{\begin{array}{ll}\lambda & \mbox{ if }j=s \\
\lambda^{-1} & \mbox{ if }j=n_i+1-s \\
1 & \mbox{ otherwise}\end{array}\right..\]
Hence we get an isomorphism
\[\C[X_*(\Se_i)]^{\Omega_i}\simeq\C[{X'_i}^{\pm 1},X_{i,1}^{\pm 1},\dots,X_{i,
q_i}^{\pm 1}]^{\{\pm 1\}^{q_i}\rtimes\Sgoth_{q_i}},\]
where $\Sgoth_{q_i}$ acts by permutations on $X_{i,1},\dots,X_{i,q_i}$ and
trivially on $X'_i$, and $\{\pm 1\}^{q_i}$ acts by $((\varepsilon_1,
\dots,\varepsilon_{q_i}),X_{i,j})\fle X_{i,j}^{\varepsilon_j}$ and
$((\varepsilon_1,\dots,\varepsilon_{q_i}),X'_i)\fle X'_i\displaystyle{\prod_
{j\ tq\ \varepsilon_j=-1}X_{i,j}^{-1}}$.
Note that the ($\Omega_i$-invariant) cocharacter
$\lambda\fle\lambda
I_{n_i}$ corresponds to $X_i:={X_i'}^2X_{i,1}^{-1}\dots X_{i,q_i}^{-1}$.

Assume that $n_i$ is odd. Then there is an isomorphism
\[\C[X_*(\Se_i)]\simeq\C[X_i^{\pm 1},X_{i,1}^{\pm 1},\dots,X_{i,q_i}^{\pm 1}]\]
that sends $X_i$ to the cocharacter $\lambda\fle\lambda I_{n_i}$
and $X_{i,s}$, $1\leq s\leq q_i$, to the cocharacter defined by the same
formula as when $n_i$ is even.
Hence we get an isomorphism
\[\C[X_*(\Se_i)]^{\Omega_i}\simeq\C[X_i^{\pm 1},X_{i,1}^{\pm 1},
\dots,X_{i,q_i}^{\pm 1}]^{\{\pm 1\}^{q_i}\rtimes\Sgoth_{q_i}},\]
where $\{\pm 1\}^{q_i}\rtimes\Sgoth_{q_i}$ acts as before on
$X_{i,1},\dots,X_{i,q_i}$ (and trivially on $X_i$).

Finally, we get, if all the $n_i$ are even,
\[\C[X_*(\Se_G)]^{\Omega_G(L)}\simeq\C[(X_1'\dots X_r')^{\pm 1},X_{i,j}^{
\pm 1},1\leq i\leq r,1\leq j\leq q_i]^{\Omega_1\times\dots\times\Omega_r},\]
and, if at least one of the $n_i$ is odd,
\[\C[X_*(\Se_G)]^{\Omega_G(L)}\simeq\C[(X_1\dots X_r)^{\pm 1},X_{i,j}^
{\pm 1},1\leq i\leq r,1\leq j\leq q_i]^{\Omega_1\times\dots\times\Omega_r}.\]
Let $X'=X'_1\dots X'_r$ and $X=X_1\dots X_r$.

In order to unify notations later, write, for every
$i\in\{1,\dots,r\}$ and $j\in\{n_i+1-q_i,\dots,n_i\}$,
$X_{i,j}=X_{i,n_i+1-j}^{-1}$ and, for every $i\in\{1,\dots,r\}$ such that
$n_i$ is odd, $X_{i,\frac{n_i+1}{2}}=1$.

Note that $X$ does not stand for the same cocharacter if $\G$ splits over
$L$ or does not (neither do the $X_{i,j}$, but this is more obvious). Let
$\nu:\Gr_{m,L}\fl\G_{L}$ be the cocharacter corresponding to $X$.
If $\G$ does not split over $L$,
then $\nu$ is defined over $\Q$ and $c(\nu(\lambda))=
\lambda^2$ for every $\lambda\in\Gr_m$. If $\G$ splits over $L$, then $\nu$
is defined over $E$ (and is not defined over $\Q$) and $c(\nu(\lambda))=
\lambda$ for every $\lambda\in\Gr_{m,E}$.

We end this subsection with an explicit version of the result of
\ref{partie_en_p1}. Assume that $\G$ splits over $L$ (ie that
$L$ contains $E_\wp$), and fix a uniformizer $\varpi_L$ of $L$.
Set $d=[L:\Q_p]$.
Let $s_1,\dots,s_r\in\Nat$ be such that $s_i\leq n_i$. For every
$i\in\{1,\dots,r\}$, there is a cocharacter
$\mu_{s_i}:\Gr_{m,E}\fl\GU^*(n_i)_E$, defined in \ref{notation:mu_p}.
Let $\mu=(\mu_{s_1},\dots,\mu_{s_r}):\Gr_{m,E}\fl\G_E$ and
\[\phi=\frac{\ungras_{\K_L\mu(\varpi_L^{-1})\K_L}}{vol(\K_L)}\in\Hecke(\G(L),
\K_L)\]
(with the notations of \ref{partie_en_p1}, $\phi=f_\mu$).
Let $r_{-\mu}$ be the representation of ${}^L\G_{E_\wp}$ associated to $-\mu$
as in lemma \ref{def_r_mu}, and $\Phi\in W(\Q_p^{ur}/\Q_p)$ be the
geometric Frobenius (so $\Phi^d$ is a generator of $W(\Q_p^{ur}/L)$).

\begin{proposition}\label{prop:transfo_Satake_CB}
For every admissible unramified morphism
$\varphi:W_{L}\fl{}^L\G$,
\[\Tr(\pi_\varphi(\phi))=p^{d(s_1(n_1-s_1)/2+\dots+s_r(n_r-s_r)/2)}
\Tr(r_{-\mu}(\varphi(\Phi^d))).\]
In other words, the Satake transform of $\phi$ is
\[p^{d(s_1(n_1-s_1)/2+\dots+s_r(n_r-s_r)/2)}X^{-1}\sum_{{I_1\subset\{1,\dots,
n_1\}}\atop{|I_1|=s_1}}\dots\sum_{{I_r\subset\{1,\dots,n_r\}}\atop{|I_r|=
s_r}}\prod_{i=1}^r\prod_{j\in I_i}X_{i,j}^{-1}.\]

\end{proposition}

\begin{proof} To deduce the first formula from theorem
\ref{th:transformee_de_Satake} (ie theorem 2.1.3 of
\cite{K-SVTOI}), it is enough to show that
\[<\rho,\mu>=s_1(n_1-s_1)/2+\dots+s_r(n_r-s_r)/2,\]
where $\rho$ is half the sum of the roots of $\T_G$ in the standard Borel
subgroup of $\G$ (ie the group of upper triangular matrices).
This is an easy consequence of the definition of $\mu$.

In the reformulation below formula (2.3.4) of \cite{K-SVTOI},
theorem 2.1.3 of \cite{K-SVTOI} says that the Satake transform of $\phi$ is
\[p^{d(s_1(n_1-s_1)/2+\dots+s_r(n_r-s_r)/2)}\sum_{\nu\in\Omega_G(L)(-\mu)}\nu.
\]
To prove the second formula, it is therefore enough to notice that
$-\mu\in X_*(T_G)$ corresponds by the isomorphism
$\C[X_*(\T_G)]\simeq\C[X^{\pm 1}]\otimes\C[X_{i,j}^{\pm 1}]$ to
$X^{-1}\prod_{i=1}^r\limits(X_{i,1}\dots X_{i,s_i})^{-1}$.

\end{proof}

\subsection*{The base change map}

In this subsection, $L$ is still an unramified extension of $\Q_p$.
Write $\K_0=\G(\Z_p)$ and $d=[L:\Q_p]$. If $L$ contains $E_\wp$, write
$a=[L:E_\wp]$.
Let $R=R_{L/\Q_p}\G_L$. Then there is a ``diagonal'' $L$-morphism
$\eta:{}^L\G\fl{}^LR$ (cf example \ref{ex:groupes_non_connexes}). It induces
a morphism
$b_\eta:\Hecke(\G(L),\K_L)\fl\Hecke(\G(\Q_p),\K_0)$, called \emph{base
change map} (or \emph{stable base change map}).
\index{base change map}
\index{stable base change map}
We want to calculate this morphism. To avoid confusion,
when writing the Satake isomorphism for $\Hecke(\G(L),\K_L)$, we will use
the letter $Z$ (instead of $X$) for the indeterminates (and we will still
use $X$ when writing the Satake isomorphism for
$\Hecke(\G(\Q_p),\K_0)$).

Assume first that $\G$ does not split over $L$ (so that $p$ is inert in $E$).
Then the base change morphism corresponds by the Satake isomorphisms to the
morphism induced by
\[\begin{array}{rcl}\C[Z^{\pm 1}]\otimes\C[Z_{i,j}^{\pm 1},1\leq i\leq r,
1\leq j\leq q_i] & \fl
& \C[X^{\pm 1},X_{i,j}^{\pm 1},1\leq i\leq r,1\leq j\leq q_i]\\
Z & \fle & X^d \\
Z_{i,j} & \fle & X_{i,j}^d\end{array}\]

Assume that $\G$ splits over $L$ but not over $\Q_p$. Then $L\supset E_\wp$,
$p$ is inert in $E$ and $d=2a$.
The base change morphism corresponds by the Satake
isomorphisms to the morphism induced by
\[\begin{array}{rcl}\C[Z^{\pm 1}]\otimes\C[Z_{i,j}^{\pm 1},1\leq i\leq r,
1\leq j\leq n_i] & \fl
& \C[X^{\pm 1},X_{i,j}^{\pm 1},1\leq i\leq r,1\leq j\leq q_i]\\
Z & \fle & X^a \\
Z_{i,j} & \fle & \left\{\begin{array}{ll}X_{i,j}^a & \mbox{ if }1\leq j\leq q_i
 \\ 1 & \mbox{ if }q_i+1\leq j\leq n_i-q_i \\ X_{i,n_i+1-j}^{-a} & \mbox{ if }
n_i+1-q_i\leq j\leq n_i\end{array}\right.\end{array}\]

Assume that $\G$ splits over $L$ and $\Q_p$. Then $p$ splits in $E$ and
$L\supset E_\wp=\Q_p$, so $d=a$.
The base change morphism corresponds by the Satake isomorphisms to the morphism
induced by
\[\begin{array}{rcl}\C[Z^{\pm 1}]\otimes\C[Z_{i,j}^{\pm 1},1\leq i\leq r,
1\leq j\leq n_i]
& \fl & \C[X^{\pm 1}]\otimes\C[X_{i,j}^{\pm 1},1\leq i\leq r,1\leq j\leq n_i]\\
Z & \fle & X^a \\
Z_{i,j} & \fle & X_{i,j}^a\end{array}\]

Notice that, with the conventions of the previous subsection in the case when
$\G$ does not split over $\Q_p$, the base change morphism is given by the same
formulas in the last two cases (this was the point of the conventions).

\begin{remark}\label{rq:image_BC} Assume that $p$ is inert in $E$ and that
$L=E_\wp=E_p$. Then the image of the base change morphism is
$\C[(X_1\dots X_r)^{\pm 1},X_{i,j}^{\pm 1},1\leq i\leq r,1\leq j\leq q_i]^
{\Omega_1\times\dots\times\Omega_r}$. In particular, the base change morphism
is surjective if and only if one of the $n_i$ is odd.

\end{remark}

\subsection*{The transfer map}

In this subsection and the next, we consider, to simplify notations, the group
$\G=\GU(J)$ with $J\in\GL_n(\Z)$ symmetric and antidiagonal,
but all results extend in an obvious way to the groups considered
before.

Let $n_1,n_2\in\Nat$ be such that $n_2$ is even and $n=n_1+n_2$. 
Let $(\H,s,\eta_0)$ be the elliptic endoscopic triple of $\G$ associated to
$(n_1,n_2)$ as in proposition \ref{prop:groupes_endoscopiques}
(note that this endoscopic triple is not always elliptic over $\Q_p$).
Let $q$ (resp. $q_1$, $q_2$) be the integral part of
$n/2$ (resp. $n_1/2$, $n_2/2$).
The group $\H$ is unramified over $\Q_p$. We will write
$\H$ for the group scheme over $\Z_p$ extending $\H$ that
is defined in remark \ref{rq:groupes_sur_Z} and  $\K_{H,0}$ for $\H(\Z_p)$.

Any unramified $L$-morphism $\eta:{}^L\H\fl{}^L\G$ extending $\eta_0$
induces a morphism $b_\eta:\Hecke(\G(\Q_p),\K_0)\fl\Hecke(\H(\Q_p),
\K_{H,0})$, called the \emph{transfer map}.
\index{transfer map}
We want to give explicit formulas for this morphism. We will
start with a particular case. Let, as before, $\Phi\in W(\Q_p^{ur}/\Q_p)$
be the geometric Frobenius.
Let $\eta_{simple}:{}^L\H\fl{}^L\G$ be the
unramified morphism extending $\eta_0$ and such that $\eta(\Phi)$ is
equal to $(1,\Phi)$ is $p$ splits in $E$, and to $((1,A),\Phi)$ if $p$ is
inert in $E$, where $A$ is defined in proposition
\ref{prop:prolongement_eta_0}.
\index{$\eta_{simple}$}
Let $b_0=b_{\eta_{simple}}:
\Hecke(\G(\Q_p),\K_0)\fl\Hecke(\H(\Q_p),\K_{H,0})$.

Then, if $p$ is inert in $E$, $b_0$ corresponds by the Satake
isomorphisms to the morphism
\[\C[X_*(\Se_G)]^{\Omega_G(\Q_p)}\fl\C[X_*(\Se_H)]^{\Omega_H(\Q_p)}\]
defined by
\[\begin{array}{rcl}X' & \fle & X'_1X'_2 \mbox{ if }n\mbox{ is even}\\
X & \fle & X_1X_2\mbox{ if }n\mbox{ is odd} \\
X_i & \fle & \left\{\begin{array}{ll}X_{1,i} & \mbox{ if }1\leq i\leq q_1 \\
X_{2,i-q_1} & \mbox{ if }q_1+1\leq i\leq q_2 \end{array}\right.\end{array}\]
If $p$ splits in $E$, $b_0$ corresponds by the Satake isomorphisms
to the morphism
\[\C[X_*(\T_G)]^{\Omega_G(\Q_p)}\fl\C[X_*(\T_H)]^{\Omega_H(\Q_p)}\]
defined by
\[\begin{array}{rcl}X & \fle & X \\
X_i & \fle & \left\{\begin{array}{ll}X_{1,i} & \mbox{ if }1\leq i\leq n_1 \\
X_{2,i-n_1} & \mbox{ if }n_1+1\leq i\leq n_2 \end{array}\right.\end{array}\]

Now let $\eta:{}^L\H\fl{}^L\G$ be any unramified $L$-morphism extending
$\eta_0$. Then $\eta=c\eta_{simple}$, where $c:W_{\Q_p}\fl 
Z(\widehat{\H})$ is a $1$-cocycle. Write $\chi_\eta$ for the (unramified)
quasi-character of $\H(\Q_p)$ corresponding to the class of $c$ in 
$\Ho^1(W_{\Q_p},Z(\widehat{\H}))$. Then $b_\eta:\Hecke(\G(\Q_p),\K_0)\fl
\Hecke(\H(\Q_p),\K_{H,0})$ is given by the following formula :
for every $f\in\Hecke(\G(\Q_p),\K_0)$, $b_\eta(f)=\chi_\eta b_0(f)$.

Following Kottwitz (\cite{K-SVLR} p 181), we use this to define $b_\eta$
even if $\eta$ is not unramified. So let $\eta:{}^L\H\fl{}^L\G$ be a
(not necessarily unramified) $L$-morphism extending $\eta_0$. Define
a quasi-character $\chi_\eta$ of $\H(\Q_p)$ as before ($\chi_\eta$ can
be ramified), and \emph{define} $b_\eta:\Hecke(\G(\Q_p),\K_0)\fl
C_c^\infty(\H(\Q_p))$ by the following formula : for every $f\in
\Hecke(\G(\Q_p),\K_0)$, $b_\eta(f)=\chi_\eta b_0(f)$.

\subsection*{The twisted transfer map}

Keep the notations of the previous subsection. Fix an unramified extension
$L$ of $\Q_p$ and write as before $\K_L=\G(\Of_L)$, $d=[L:\Q_p]$ and, if
$L\supset E_\wp$, $a=[L:E_\wp]$; we will use the
same conventions as before when writing the Satake isomorphism for
$\Hecke(\G(L),\G(\Of_L))$ (ie the indeterminates will be $Z$ and the
$Z_i$). Let $\eta:{}^L\H\fl
{}^L\G$ be an unramified $L$-morphism extending $\eta_0$.

Remember the definition of the twisted endoscopic datum associated to
$(\H,s,\eta)$ and to the field extension
$L/\Q_p$ (cf \cite{K-SVLR} p 179-180).
Let $\Phi\in W(\Q_p^{ur}/\Q_p)$ be the geometric Frobenius,
$R=R_{L/\Q_p}\G_L$ and $\theta$ be the automorphism of $R$ corresponding to
$\Phi$. Then
\[\widehat{R}=(\widehat{\G})^d,\]
where the $i$-th factor corresponds to the image of $\Phi^{d-i}$ in 
$\Gal(L/\Q_p)$. The group $W_{\Q_p}$ acts on $\widehat{R}$ via its quotient
$W(\Q_p^{ur}/\Q_p)$, and $\Phi$ acts on $\widehat{R}$ by
\[\Phi(g_1,\dots,g_d)=\widehat{\theta}(\Phi(g_1),\dots,\Phi(g_d))
=(\Phi(g_2),\dots,\Phi(g_d),\Phi(g_1)).\]
In particular, the diagonal embedding $\widehat{\G}\fl\widehat{R}$ is
$W_{\Q_p}$-equivariant, so it extends in an obvious way to a
$L$-morphism ${}^L\G\fl{}^L R$; let $\eta':{}^L\H\fl{}^L R$ denote the
composition of $\eta:{}^L\H\fl{}^L\G$ and of this morphism.
Let $t_1,\dots,t_d\in Z(\widehat{\H})^{\Gal(\overline{\Q}_p/\Q_p)}$ be such
that $t_1\dots t_d=s$, and write $t=(t_1,\dots,t_d)\in\widehat{R}$. Define
a morphism $\widetilde{\eta}:\widehat{\H}\rtimes W(\Q_p^{ur}/\Q_p)\fl
\widehat{R}\rtimes W(\Q_p^{ur}/\Q_p)$ by :
\begin{bulletlist}
\item $\widetilde{\eta}_{|\widehat{\H}}$ is the composition of
$\eta_0:\widehat{\H}\fl\widehat{\G}$ and of the diagonal embedding
$\widehat{\G}\fl\widehat{R}=(\widehat{\G})^d$,
\item $\widetilde{\eta}((1,\Phi))=(t,1)\eta'(1,\Phi)$.

\end{bulletlist}
Then the $\widehat{R}$-conjugacy class of $\widetilde{\eta}$ does not depend
on the choice of $t_1,\dots,t_d$, and $(\H,t,\widetilde{\eta})$ is a twisted
endoscopic datum for $(R,\theta)$. The map
\[b_{\til{\eta}}:\Hecke(\G(L),\K_L)\fl\Hecke(\H(\Q_p),\K_{H,0})\]
induced by $\til{\eta}$ is called the \emph{twisted transfer map} (or the
\emph{unstable base change map}).
\index{twisted transfer map}
\index{unstable base change map}

Assume first that $\eta=\eta_{simple}$, and write $\til{b}_0$ for
$b_{\til{\eta}}$. If $\G$ does not split over $L$,
then the twisted transfer map $\til{b}_0$ corresponds by the Satake
isomorphisms to the morphism induced by
\[\begin{array}{rcl}\C[Z^{\pm 1}]\otimes\C[Z_i^{\pm 1}]
& \fl & \C[X^{\pm 1}]\otimes\C[X_{1,1}^{\pm 1},
\dots,X_{1,n_1}^{\pm 1}]\otimes\C[X_{2,1}^{\pm 1},\dots,X_{2,n_2}^{\pm 1}] \\
Z & \fle & X^d \\
Z_i & \fle & \left\{\begin{array}{ll}X_{1,i}^d & \mbox{ if }1\leq i\leq n_1 \\
-X_{2,i-n_1}^d & \mbox{ if }n_1+1\leq i\leq n\end{array}\right.
\end{array}\]
If $\G$ splits over $L$ (so $L\supset E_\wp$), then the twisted transfer
map $\til{b}_0$ corresponds by the Satake isomorphisms to the morphism
induced by
\[\begin{array}{rcl}\C[Z^{\pm 1}]\otimes\C[Z_i^{\pm 1}]
& \fl & \C[X^{\pm 1}]\otimes\C[X_{1,1}^{\pm 1},
\dots,X_{1,n_1}^{\pm 1}]\otimes\C[X_{2,1}^{\pm 1},\dots,X_{2,n_2}^{\pm 1}] \\
Z & \fle & X^a \\
Z_i & \fle & \left\{\begin{array}{ll}X_{1,i}^a & \mbox{ if }1\leq i\leq n_1 \\
-X_{2,i-n_1}^a & \mbox{ if }n_1+1\leq i\leq n\end{array}\right.
\end{array}\]

Assume now that $\eta$ is any unramified extension of $\eta_0$, and define
an unramified quasi-character $\chi_\eta$ of $\H(\Q_p)$ as in the
previous subsection. Then, for every $f\in\Hecke(\G(L),\K_L)$,
$b_{\til{\eta}}(f)=\chi_\eta\til{b}_0(f)$.

As in the previous subsection, we can use this to define $b_{\til{\eta}}$
for a possibly ramified $\eta$ (this is just \cite{K-SVLR} p 181). Let
$\eta:{}^L\H\fl{}^L\G$ be any $L$-morphism extending $\eta_0$, and attach
to it a (possibly ramified) quasi-character $\chi_\eta$ of $\H(\Q_p)$.
Define $b_{\til{\eta}}:\Hecke(\G(L),\K_L)\fl C_c^\infty(\H(\Q_p))$ by the
following formula : for every $f\in\Hecke(\G(L),\K_L)$,
$b_{\til{\eta}}(f)=\chi_\eta\til{b}_0(f)$.

\section{Twisted transfer map and constant terms}
\label{partie_en_p3}

In this section, we consider, to simplify notation, the situation of the last
subsection of \ref{partie_en_p2},
but all results extend in an obvious way to the groups
$\G(\U(J_1)\times\dots\times\U(J_r))$ of the beginning of
\ref{partie_en_p2}. Assume that $\G$ splits over $L$ (ie that $L$ contains
$E_\wp$), and
fix a $L$-morphism $\eta:{}^L\H\fl{}^L\G$ extending
$\eta_0$; we do not assume that $\eta$ is unramified.

Let $\M$ be a cuspidal standard Levi subgroup of $\G$, and let
$r\in\{1,\dots,q\}$ be such that $\M=\M_{\{1,\dots,r\}}\simeq
(R_{E/\Q}\Gr_m)^r\times\GU^*(m)$, with $m=n-2r$.
Let $(\M',s_M,\eta_{M,0})$ be an
element of $\Ell_\G(\M)$ (cf \ref{groupes4}) whose image in
$\Ell(\G)$ is $(\H,s,\eta_0)$. Assume that $s_M=s_{A,m_1,m_2}$,
with $A\subset\{1,\dots,r\}$ and $m_1+m_2=m$, and where notations are as
in lemma \ref{lemme:Levi_endoscopiques_GU}; define $r_1$ and $r_2$ as in
this lemma. There is a conjugacy class of Levi subgroups of $\H$ associated
to $(\M',s_M,\eta_{M,0})$; let $\M_H$ be the standard Levi subgroup in
that class. Then
$\M_H=\H \cap (\M_{H,1}\times\M_{H,2})$, with
\[\M_{H,i}=\GU^*(n_i)\cap\left(\begin{array}{ccccccc}
R_{E/\Q}\Gr_m & & 0 & & & & 0 \\
& \ddots & & & & & \\
0 & & R_{E/\Q}\Gr_m & & & & \\
& & & \GU^*(m_i) & & & \\
& & & & R_{E/\Q}\Gr_m & & 0 \\
& & & & & \ddots & \\
0 & & & & 0 & & R_{E/\Q}\Gr_m\end{array}\right),\]
where the diagonal blocks are of size $r_i,m_i,r_i$.
On the other hand, $\M_H=\M_{H,l}\times\M_{H,h}$, where
$\M_{H,l}=(R_{E/\Q}\Gr_m)^{r_1}\times (R_{E/\Q}\Gr_m)^{r_2}$
is the linear part of $\M_H$ and
$\M_{H,h}=\G(\U^*(m_1)\times\U^*(m_2))$
is the Hermitian part.
Similarly, $\M=\M_l\times\M_h$, where $\M_l=\M_{H,l}=(R_{E/\Q}\Gr_m)^r$
is the linear part of $\M$ and $\M_h=\GU^*(m)$ is the Hermitian part.
The morphism $\eta:{}^L\H\fl{}^L\G$ determines a $L$-morphism
$\eta_M:{}^L\M_H={}^L\M'\fl{}^L\M$ extending $\eta_{M,0}$, unique up to
$\widehat{\M}$-conjugacy; $\eta_M$ is unramified if $\eta$ is unramified.

As in lemma \ref{lemme:Levi_endoscopiques_GU}, identify $\widehat{\M}$ with the
Levi subgroup
\[\C^\times\times\left(\begin{array}{ccc}\begin{array}{ccc}* & & 0 \\
& \ddots & \\ 0 & & *\end{array} & & 0 \\
 & \GL_m(\C) & \\
0 & & \begin{array}{ccc}* & & 0 \\ & \ddots & \\ 0 & & *\end{array}\end{array}
\right)\]
(blocks of size $r,m,r$) of $\widehat{\G}$.
Identify $\widehat{\M}_H$ to a Levi subgroup of $\widehat{H}$
in a similar way.
Let $s'_M$ be the element of
$Z(\widehat{\M}_H)\subset\widehat{\M}\simeq\C^\times\times(\C^\times)^r\times
\GL_m(\C)\times(\C^\times)^r$ equal to $(1,(1,\dots,1),s_{M_h},(1,\dots,1))$,
where $s_{M_h}$ is the image of $s_M$ by the projection
$\widehat{\M}\fl\GL_m(\C)$ (in the notation of lemma
\ref{lemme:Levi_endoscopiques_GU}, $s'_M=s_{\varnothing,m_1,m_2}$).
Then $(\M_H,s'_M,\eta_{M,0})$ is an elliptic endoscopic triple for $\M$,
isomorphic to $(\M',s_M,\eta_{M,0})$ as an endoscopic $\M$-triple
(but not as an endoscopic $\G$-triple).

Write $b_{s_M},b_{s'_M}:\Hecke(\M(L),\M(\Of_L))\fl C_c^\infty(\M_H(\Q_p))$
for the twisted transfer maps associated to $(\M',s_M,\eta_M)$
and $(\M_H,s'_M,\eta_M)$, and $f\fle f_\M$ (resp $f\fle f_{\M_H}$) for
the constant term map $\Hecke(\G(L),\K_L)\fl\Hecke(\M(L),\M(\Of_L))$ 
or $\Hecke(\G(\Q_p),\K_0)\fl\Hecke(\M(\Q_p),\M(\Z_p))$ (resp.
$\Hecke(\H(\Q_p),\K_{H,0})\fl\Hecke(\M_H(\Q_p),\M_H(\Z_p))$). Then it is
easy to see from the definitions that, if $\eta$ is unramified, then, for
every $f\in\Hecke(\G(L),\K_L)$, $b_{s_M}(f_\M)=(b_{\til{\eta}}(f))_{\M_H}$.
There is
a similar formula for a general $\eta$ : Let $\chi_\eta$ be the quasi-character
of $\H(\Q_p)$ associated to $\eta$ as in the last two subsections of
\ref{partie_en_p2}, and write $b_{s_M,0}$ for the twisted transfer map
defined by $s_M$ in the case $\eta=\eta_{simple}$ (and $\til{b}_0=
b_{\til{\eta}_{simple}}$, as before).
Then, for every $f\in\Hecke(\G(L),\K_L)$, $b_{s_M}(f_\M)=\chi_{\eta|\M_H(\Q_p)}
b_{s_M,0}(f_\M)=\chi_{\eta|\M_H(\Q_p)}(\til{b}_0(f))_{\M_H}$.

Later, we will use the twisted transfer map $b_{s'_M}$ and not
$b_{s_M}$, so we need to compare it to $b_{s_M}$ (ot to $b_{\til{\eta}}$),
at least on certain elements of $\Hecke(\M(L),\M(\Of_L))$. 
First we give explicit formulas for it in the case $\eta=\eta_{simple}$.
Write $\Omega_M(L)=W(\T_G(L),\M(L))$ and $\Omega_{M_H}(\Q_p)=
W(\Se_H(\Q_p),\M_H(\Q_p))$. Then we get Satake isomorphisms
$\Hecke(\M(L),\M(\Of_L))\simeq\C[X_*(\T_G)]^{\Omega_M(L)}$ and
$\Hecke(\M_H(\Q_p),\M_H(\Z_p))\simeq\C[X_*(\Se_H)]^{\Omega_{M_H}(\Q_p)}$,
and, if $\eta=\eta_{simple}$,
then the twisted transfer map $b_{s'_M}$ is induced by the morphism :
\[\begin{array}{rcl}\C[Z^{\pm 1}]\otimes\C[Z_1^{\pm},\dots,Z_n^{\pm 1}]
& \fl & \C[X^{\pm 1}]\otimes\C[X_{1,1}^{\pm 1},
\dots,X_{1,n_1}^{\pm 1}]\otimes\C[X_{2,1}^{\pm 1},\dots,X_{2,
n_2}^{\pm 1}] \\
Z & \fle & X^a \\
Z_{i_k} & \fle & X_{1,k}^a \\
Z_{n+1-i_k} & \fle & X_{1,n_1+1-k}^a \\
Z_{j_l} & \fle & X_{2,l}^a \\
Z_{n+1-j_l} & \fle & X_{2,n_2+1-l}^a \\
Z_i & \fle & \left\{\begin{array}{ll}X_{1,i-r_2}^a & \mbox{ if }r+1\leq i
\leq r+m_1 \\
-X_{2,i-(r_1+m_1)}^a & \mbox{ if }r+m_1+1\leq i\leq r+m\end{array}\right.
\end{array}\]
where we wrote $\{1,\dots,r\}-A=\{i_1,\dots,i_{r_1}\}$ and
$A=\{j_1,\dots,j_{r_2}\}$ with $i_1<\dots<i_{r_1}$ and $j_1<\dots<j_{r_2}$.

Let $\alpha\in\Nat$ such that $n-q\leq\alpha\leq n$. Write
$\mu=\mu_\alpha$, where $\mu_\alpha:\Gr_{m,L}\fl\GU^*(n_i)_L$ is the
cocharacter defined in \ref{notation:mu_p}.
This cocharacter factors through $\M_L$, and we denote by $\mu_M$ the
cocharacter of $\M_L$ that it induces.
Set
\[\phi=\frac{\ungras_{\K_L\mu(\varpi_L^{-1})\K_L}}{vol(\K_L)}\in\Hecke(\G(L),
\K_L)\]
and
\[\phi^\M=\frac{\ungras_{\M(\Of_L)\mu_M(\varpi_L^{-1})\M(\Of_L)}}{vol(
\M(\Of_L))}\in\Hecke(\M(L),\M(\Of_L)).\]
Note that, if $\alpha\leq n-r$, then $\phi^\M$ is the product of a function in
$\Hecke(\M_h(L),\M_h(\Of_L))$ and of the unit element of
$\Hecke(\M_l(L),\M_l(\Of_L))$ (because the image of $\mu_M$ is included in
$\M_{h,L}$ in that case).

The Satake transform of $\phi$ has been calculated in proposition
\ref{prop:transfo_Satake_CB}; it is equal to
\[p^{d\alpha(n-\alpha)/2}Z^{-1}\sum_{{I\subset\{1,\dots,n\}}\atop{|I|=\alpha}}
\prod_{i\in I}Z_i^{-1}.\]
Identify $\Hecke(\G(L),\K_L)$ to a subalgebra of $\Hecke(\M(L),
\M(\Of_L))$ using the constant term morphism (via the Satake isomorphisms,
this corresponds to the obvious inclusion $\C[X_*(\T_G)]^{\Omega_G(L)}
\subset\C[X_*(\T_G)]^{\Omega_M(L)}$).
If $\alpha\geq n-r+1$, then the Satake transform of $\phi^M$ is simply
$(ZZ_1\dots Z_\alpha)^{-1}$. If $\alpha\leq n-r$, then, by proposition
\ref{prop:transfo_Satake_CB} (applied to $\M_h$), the Satake transform
of $\phi^M$ is
\[p^{d(\alpha-r)(n-\alpha-r)/2}(ZZ_1\dots Z_r)^{-1}\sum_{{I\subset\{r+1,
\dots,n-r\}}\atop{|I|=\alpha-r}}\prod_{i\in I}Z_i^{-1}.\]

Let $f^{\H}=b_{\til{\eta}}(\phi)\in C_c^\infty(\H(\Q_p))$, $f^{\M_H}=
b_{s'_M}(\phi^{\M})\in C_c^\infty(\M_H(\Q_p))$ and
$\psi^{\M_H}=b_{s'_M}(\phi_\M)\in C_c^\infty(\M_H(\Q_p))$.
By the definition of $b_{\til{\eta}}$, there exists a quasi-character
$\chi_\eta$ of $\H(\Q_p)$ such that $\chi_\eta^{-1}f^\H\in\Hecke(\H(\Q_p),
\K_{H,0})$; write $f_{\M_H}^{\H}\in C_c^\infty(\H(\Q_p))$ for
$\chi_{\eta|\M_H(\Q_p)}(\chi_\eta^{-1}f^\H)_{\M_H}$.
(Of course, because we used $b_{s'_M}$ and not $b_{s_M}$ to define
$\psi^{\M_H}$, the functions $\psi^{\M_H}$ and $f^\H_{\M_H}$ are
different in general.)

Let $\Omega^{M_H}$ be the subgroup of $\Omega_H(\Q_p)\subset\Sgoth_{n_1}
\times\Sgoth_{n_2}$ generated by the transpositions $((j,n_1+1-j),1)$,
$1\leq j\leq r_1$, and $(1,(j,n_2+1-j))$, $1\leq j\leq r_2$. Then
$\Omega^{M_H}\simeq\{\pm 1\}^{r_1}\times\{\pm 1\}^{r_2}$ (actually,
$\Omega^{M_H}$ is even a subgroup of the relative Weyl group over $\Q$,
$W(\Se_H(\Q),\H(\Q))$).
\index{$\Omega^{M_H}$}

\begin{proposition}\label{prop:identite_en_p} Let $\gamma_H\in\M_H(\Q)$ be
such that $O_{\gamma_H}(f_{\M_H}^{\H})\not=0$. Write $\gamma_H=\gamma_l
\gamma_h$, with $\gamma_h\in\M_{H,h}(\Q)$ and $\gamma_l=((\lambda_{1,1},\dots,
\lambda_{1,r_1}),(\lambda_{2,1},\dots,\lambda_{2,r_2}))\in\M_{H,l}(\Q)=
(E^\times)^{r_1}\times(E^\times)^{r_2}$, and let
\[N_{s_M}(\gamma_H)=\frac{1}{2a}\sum_{i=1}^{r_2}val_p(|
\lambda_{2,i}\overline{\lambda}_{2,i}|_p)\]
(where $val_p$ is the $p$-adic valuation).
\index{valp@$val_p$\quad $p$-adic valuation}
\index{NsM@$N_{s_M}$}
Then $|c(\gamma_H)|_p=p^{d}$, $\frac{1}{2a}val_p(|\lambda_{i,j}\overline{
\lambda}_{i,j}|_p)$ is an integer for every $i,j$ (in particular,
$N_{s_M}(\gamma_H)\in\Z$), and one and only one of the following two assertions
is true :
\begin{itemize}
\item[(A)] There exists $\omega\in\Omega^{M_H}$, uniquely determined by
$\gamma_H$, such that $\omega(\gamma_H)\in\M_{H,l}(\Z_p)\M_{H,h}(\Q_p)$.
\item[(B)] There exists $i\in\{1,2\}$ and $j\in\{1,\dots,r_i\}$ such that
$\frac{1}{2a}val_p(|\lambda_{i,j}\overline{\lambda}_{i,j}|_p)$ is odd.

\end{itemize}

Besides, case (A) can occur only if $\alpha\leq n-r$.

Choose an element $\gamma\in\M(\Q_p)$ coming from $\gamma_H$ (such a
$\gamma$ always exists because $\M$ is quasi-split over $\Q_p$, cf
\cite{K-RCC}). Then
\[O_{\gamma_H}(f_{\M_H}^\H)=<\mu,s><\mu,s'_M>\varepsilon_{s_M}(\gamma)O_{
\gamma_H}(\psi^{\M_H}),\]
where $\varepsilon_{s_M}(\gamma)=(-1)^{N_{s_M}(\gamma)}$.
\index{$\varepsilon_{s_M}$}
If moreover $\gamma_l\in\M_{H,l}(\Z_p)$ (this can happen only if
$\alpha\leq n-r$), then
\[O_{\gamma_H}(f_{\M_H}^\H)=<\mu,s><\mu,s'_M>\delta_{P(\Q_p)}^{1/2}(\gamma)O_
{\gamma_H}(f^{\M_H}),\]
where $\Pa$ is the standard parabolic subgroup of $\G$ with Levi subgroup
$\M$.

\end{proposition}

Let $\gamma_0$ be the component of $\gamma$ in $\M_h(\Q_p)$. For every
$\delta\in\M(L)$ $\sigma$-semi-simple, let $\delta_h$ be the component of
$\delta$ in $\M_h(L)$ and define $\alpha_p(\gamma,\delta)$ and
$\alpha_p(\gamma_0,\delta_h)$ as in as \cite{K-SVLR} p 180 (cf also
subsections \ref{app:alpha_reg} and \ref{app:alpha} of the appendix).

From the definition of $s'_M$, it is
clear that
\[<\alpha_p(\gamma,\delta),s'_M>=<\alpha_p(\gamma_0,\delta_h),s'_M>=<\alpha_p(
\gamma_0,\delta_h),s_M>.\]
After applying corollary \ref{lemme:ss_non_regulier}, we get :\footnote{
Corollary \ref{lemme:ss_non_regulier} applies only to a particular choice of
$\eta$, ie to $\eta=\eta_{simple}$. But it is explained on page 181 of
\cite{K-SVLR} (after formula (7.3)) why this is enough to prove the next
corollary for any choice of $\eta$.}

\begin{corollary}\label{cor:identite_en_p} There is an equality
\[SO_{\gamma_H}(f^\H_{\M_H})=<\mu,s><\mu,s'_M>\varepsilon_{s_M}(\gamma)\sum
_\delta<\alpha_p(\gamma_0,\delta_h),s_M>\Delta_{\M_H,s_M,p}^\M(\gamma_H,\gamma)
e(\delta)TO_\delta(\phi_\M).\]
If moreover $\gamma\in\M_l(\Z_p)\M_h(\Q_p)$, then
\[SO_{\gamma_H}(f^\H_{\M_H})=<\mu,s><\mu,s'_M>\delta_{\Pa(\Q_p)}^{1/2}(\gamma)
\sum_\delta<\alpha_p(\gamma_0,\delta_h),s_M>\Delta_{\M_H,s_M,p}^\M(\gamma_H,
\gamma)e(\delta)TO_\delta(\phi^\M).\]
The two sums above are taken over the set of $\sigma$-conjugacy classes
$\delta$ in $\M(L)$ such that $\gamma$ and $N\delta$ (defined in
\ref{points_fixes6}, after theorem \ref{th:points_fixes_Kottwitz}) are
$\M(\overline{\Q}_p)$-conjugate; for every such $\delta$, we write
$e(\delta)=e(R_{\delta\theta})$, where $R_{\delta\theta}$ is the
$\theta$-centralizer of $\delta$ in $R$ (denoted by $I(p)$ in
\ref{points_fixes6}) and $e$ is the sign of \cite{K-SC}.

\end{corollary}

\begin{proofp} It is obvious from the definition of the twisted transfer
maps that it suffices to prove the proposition for $\eta=\eta_{simple}$.
So we assume this in all the proof.

It is easy to see that the element $\omega\in\Omega^{M_H}$ in (A)
is necessarily unique (if we already know that $|c(\gamma_H)|_p=p^d$).
This comes from the fact that, for every
$\gamma_H\in\M_{H,l}(\Z_p)\M_{H,h}(\Q_p)$ such that $|c(\gamma_H)|_p\not=1$
and for every $\omega\in\Omega^{M_H}$, $\omega(\gamma_H)\not=\gamma_H$.

We know that the Satake transform of $\phi$ is
\[p^{d\alpha(n-\alpha)/2}Z^{-1}\sum_{{I\subset\{1,\dots,n\}}\atop{|I|=\alpha}}
\prod_{i\in I}Z_i^{-1}.\]
For every $I\subset\{1,\dots,n\}$, write
\[n(I)=|I\cap\{n_1+1,\dots,n\}|,\]
and
\[a_I=X^{-a}\prod_{i\in I\cap\{1,\dots,n_1\}}X_{1,i}^{-a}
\prod_{i\in I\cap\{n_1+1,\dots,n\}}X_{2,i-n_1}^{-a}.\]
Then the Satake transform of $f^{\H}_{M_H}$ (that is equal to the Satake
transform of $f^{\H}$) is  
\[S:=p^{d\alpha(n-\alpha)/2}\sum_{{I\subset\{1,\dots,n\}}\atop{|I|=\alpha}}
(-1)^{n(I)}a_I.\]
As $S$ is the product of $X^{-a}$ and of a polynomial in the $X_{i,j}^{-1}$, 
if $O_{\gamma_H}(f_{\M_H}^\H)\not=0$, then $|c(\gamma_H)|_p$ must be equal
to $p^d$.

Let $A_{l,1}=\{1,\dots,r_1\}\cup\{n_1+1-r_1,\dots,n_1\}$,
$A_{l,2}=\{n_1+1,\dots,n_1+r_2\}\cup\{n+1-r_2,\dots,n\}$,
$A_l=A_{l,1}\cup A_{l,2}$, $A_{h,1}=\{r_1+1,\dots,n_1-r_1\}$,
$A_{h,2}=\{n_1+r_2+1,\dots,n-r_2\}$ and $A_h=A_{h,1}\cup A_{h,2}$.
For every $I_l\subset A_l$, write $n_l(I_l)=|I_l\cap A_{l,2}|$ and
\[b_{I_l}=\prod_{j\in I_l\cap A_{l,1}}X_{1,j}^{-a}
\prod_{j\in I_l\cap A_{l,2}}X_{2,j-n_1}^{-a}.\]
For every $I_h\subset A_h$,
write $n_h(I_h)=|I_h\cap A_{h,2}|$, and
\[c_{I_h}=X^{-a}\prod_{j\in I_h\cap A_{h,1}}X_{1,j}^{-a}
\prod_{j\in I_h\cap A_{h,2}}X_{2,j-n_1}^{-a}.\]
As
$a_I=b_{I\cap A_l}c_{I\cap A_h}$ and $n(I)=n_l(I\cap A_l)+n_h(I\cap A_h)$
for every $I\subset\{1,\dots,n\}$,
\[S=p^{d\alpha(n-\alpha)/2}\sum_{k=0}^{\alpha}\left(\sum_{I_l\subset A_l,
|I_l|=k}(-1)^{n_l(I_l)}b_{I_l}\right)
\left(\sum_{I_h\subset A_h,|I_h|=\alpha-k}(-1)^{n_h(I_h)}c_{I_h}\right).\]
For every $k\in\{0,\dots,\alpha\}$, the polynomial $\sum\limits_{I_h\subset
A_h,|I_h|=\alpha-k}(-1)^{n_h(I_h)}c_{I_h}$ is the Satake transform of a
function in $\Hecke(\M_{H,h}(\Q_p),\M_{H,h}(\Z_p))$, that will be denoted by
$\psi_{h,k}$. For every $I_l\subset A_l$, the monomial
$(-1)^{n_l(I_l)}b_{I_l}$ is the Satake transform of a function in
$\Hecke(\M_{H,l}(\Q_p),\M_{H,l}(\Z_p))$, that will be denoted by
$\psi_{I_l}$.
Then
\[f_{\M_H}^\H=p^{d\alpha(n-\alpha)/2}\sum_{k=0}^{\alpha}\psi_{h,k}\left(
\sum_{I_l\subset A_l,|I_l|=k}\psi_{I_l}\right).\]
As $O_{\gamma_H}(f_{\M_H}^{\H})\not=0$, there exist $k\in\{0,\dots,\alpha\}$
and $I_l\subset A_l$ such that $|I_l|=k$ and
$O_{\gamma_H}(\psi_{I_l}\psi_{h,k})\not=0$. 
Write $\gamma_H=\gamma_l\gamma_h$, with $\gamma_h\in\M_{H,h}(\Q_p)$
and $\gamma_l=((\lambda_{1,1},\dots,\lambda_{1,r_1}),(\lambda_{2,1},\dots,
\lambda_{2,r_2}))\in\M_{H,l}(\Q_p)$. Then $O_{\gamma_H}(\psi_{I_l}\psi_{h,k})
=O_{\gamma_l}(\psi_{I_l})O_{\gamma_h}(\psi_{h,k})$. We have $O_{\gamma_l}
(\psi_{I_l})\not=0$ if and only if $\gamma_l$ is in the product of
$\M_{H,l}(\Z_p)$ and of the image of $p$ by the cocharacter corresponding
to the monomial $(X_{1,1}\dots X_{1,r_1}X_{2,1}\dots X_{2,r_2})^{a}
b_{I_l}$, and this implies that
\begin{bulletlist}
\item for every $j\in\{1,\dots,r_1\}$,
\[|\lambda_{1,j}\overline{\lambda}_{1,j}|_p=\left\{\begin{array}{ll}1 & 
\mbox{ if }j\in I_l\mbox{ and }n_1+1-j\not\in I_l \\
p^{-2a} & \mbox{ if }j,n_1+1-j\in I_l\mbox{ or }j,n_1+1-j\not\in I_l \\
p^{-4a} & \mbox{ if }j\not\in I_l\mbox{ and }n_1+1-j\in I_l\end{array}\right.\]
\item for every $j\in\{1,\dots,r_2\}$,
\[|\lambda_{2,j}\overline{\lambda}_{2,j}|_p=\left\{\begin{array}{ll}1 & 
\mbox{ if }n_1+j\in I_l\mbox{ and }n+1-j\not\in I_l \\
p^{-2a} & \mbox{ if }n_1+j,n+1-j\in I_l\mbox{ or }n_1+j,n+1-j\not\in I_l \\
p^{-4a} & \mbox{ if }n_1+j\not\in I_l\mbox{ and }n+1-j\in I_l\end{array}\right.
\]

\end{bulletlist}
This implies in particular that $\frac{1}{2a}val_p(|\lambda_{i,j}\overline{
\lambda}_{i,j}|_p)\in\Z$ for every $i,j$.
On the other hand, $c(\gamma_H)=c(\gamma_h)$, so $|c(\gamma_h)|_p=p^{d}$.

There are three cases to consider :
\begin{itemize}
\item[(1)] Assume that, for every $j\in\{1,\dots,r_1\}$,
either $j\in I_l$ or $n_1+1-j\in I_l$, and that, for every
$j\in\{1,\dots,r_2\}$, either $n_1+j\in I_l$ or $n+1-j\in I_l$.
Let $\omega=((\omega_{1,1},\dots,\omega_{1,r_1}),
(\omega_{2,1},\dots,\omega_{2,r_2}))\in\Omega^{M_h}$ be such that,
for every $j\in\{1,\dots,r_1\}$, $\omega_{1,j}=1$ if $j\not\in I_l$ and
$\omega_{1,j}=-1$ if $j\in I_l$ and, for every $j\in\{1,\dots,r_2\}$,
$\omega_{2,j}=1$ if $n_1+j\not\in I_l$ and $\omega_{2,j}=-1$ if
$n_1+j\in I_l$. It is easy to see that $\omega(\gamma)\in\M_{H,l}(\Z_p)
\M_{H,h}(\Q_p)$, and it is clear that
$\frac{1}{2a}val_p(|\lambda_{i,j}\overline{\lambda}_{i,j}|_p)$ is even for
every $i,j$.\newline
On the other hand, $k=|I_l|=r_1+r_2=r$, so, for $\psi_{h,k}=\psi_{h,r}$ to
be non-zero, we must have $\alpha-r\leq n-2r$, i.e., $\alpha\leq n-r$.
\item[(2)] Assume that there exists $j\in\{1,\dots,r_1\}$ such that
$j,n_1+1-j\in I_l$ or $j,n_1+1-j\not\in I_l$. Then
$|\lambda_{1,j}\overline{\lambda}_{1,j}|_p=p^{-2a}$, 
hence $\frac{1}{2a}val_p(|\lambda_{1,j}\overline{\lambda}_{1,j}|_p)=-1$.
\item[(3)] Assume that there exists $j\in\{1,\dots,r_2\}$ such that
$n_1+j,n+1-j_2\in I_l$ or $n_1+j,n+1-j\not\in I_l$. As in case (2), we see
that $\frac{1}{2a}val_p(|\lambda_{2,j}\overline{\lambda}_{2,j}|_p)=-1$.

\end{itemize}

We now show the last two statements of the proposition. First note that
$<\mu,s><\mu,s_M>=(-1)^{r_2}$.
For every $I\subset\{1,\dots,n\}$, write
\[m(I)=|I\cap\{n_1+r_2+1,\dots,n_1+r_2+m_2\}|.\]
The Satake transform of $\phi_\M$ is equal to the Satake transform of $\phi$,
so the Satake transform of $\psi^{\M_H}$ is
\[p^{d\alpha(n-\alpha)/2}\sum_{{I\subset\{1,\dots,n\}}\atop{|I|=\alpha}}
(-1)^{m(I)}a_I.\]
Hence
\[\psi^{\M_H}=p^{d\alpha(n-\alpha)/2}\sum_{k=0}^{\alpha}\psi_{h,k}
\left(\sum_{I_l\subset A_l,|I_l|=k}(-1)^{n_l(I_l)}\psi_{I_l}\right).\]
To show the first equality of the proposition, it is therefore enough to
see that, for every $I_l\subset A_l$ such that
$O_{\gamma_l}(\psi_{I_l})\not=0$,
$n_l(I_l)=r_2+\varepsilon_{s_M}(\gamma)$ modulo $2$; but this is an easy
consequence of the non-vanishing condition for $O_{\gamma_l}(\psi_{I_l})$
that we wrote above.

Assume that $\gamma_l\in\M_{H,l}(\Z_p)$. Then the only 
$I_l\subset A_l$ such that $O_{\gamma_l}(\psi_{I_l})\not=0$ is
$I_l=\{1,\dots,r_1\}\cup\{n_1+1,\dots,n_1+r_2\}$, and $|I_l|=r$ and
$n_l(I_l)=r_2$. We have already seen that $\psi_{h,r}=0$ unless
$\alpha\leq n-r$, so we may assume this. Then
$O_{\gamma_H}(f_{\M_H}^\H)=O_{\gamma_H}(\psi')$, where $\psi'\in\Hecke(\M_H
(\Q_p),\M_H(\Z_p))$ is the function with Satake transform
\[p^{d\alpha(n-\alpha)/2}(-1)^{r_2}\psi_{h,r}\prod_{j=1}^{r_1}X_{1,j}^{-a}
\prod_{j=1}^{r_2}X_{2,j}^{-a}.\]
Applying the calculation of the twisted transfer morphism to
$\M_h$ and $\M_{H,h}$ instead of $\G$ and $\H$,
we find that the Satake transform of $f^{\M_H}$ is
\[p^{d(\alpha-r)(n-\alpha-r)/2}\psi_{h,r}\prod_{j=1}^{r_1}X_{1,j}^{-a}
\prod_{j=1}^{r_2}X_{2,j}^{-a}.\]
So, to finish the proof of the proposition, it is enough to show that
\[\delta_{P(\Q_p)}^{-1/2}(\gamma)=p^{d(\alpha-r)(n-\alpha-r)/2-d\alpha
(n-\alpha)/2}.\]
As $\gamma$ comes from $\gamma_H$, $c(\gamma)=c(\gamma_H)$.
On the other hand, $\gamma_l\in\M_l(\Z_p)\M_h(\Q_p)$, so
\[\delta_{P(\Q_p)}(\gamma_l)=1.\]
As the image of $\gamma_h\in\GU^*(m)(\Q_p)$ in $\G(\Q_p)$ is
\[\gamma_h=\left(\begin{array}{ccc}c(\gamma_h)I_r & 0 & 0 \\
0 & \gamma_h & 0 \\
0 & 0 & I_r\end{array}\right),\]
we get
\[\delta_{P(\Q_p)}(\gamma)=\delta_{P(\Q_p)}(\gamma_h)=|c(\gamma_h)|_p^{r(r+m)}=
|c(\gamma)|_p^{r(r+m)}=p^{dr(n-r)}.\]
To conclude, it suffices to notice that
\[\alpha(n-\alpha)-(\alpha-r)(n-\alpha-r)=r(n-r).\]

\end{proofp}

\begin{remark}\label{rq:signe_en_p} From the proof above, it is easy to
see that the set of $(i,j)$ such that
$\frac{1}{2a}val_p(|\lambda_{i,j}\overline{\lambda}_{i,j}|_p)$ is odd
has an even number of elements. In particular, the sign
$\varepsilon_{s_M}(\gamma)$ does not change if $N_{s_M}(\gamma_H)$ is
replaced by $\frac{1}{2a}\sum\limits
_{i=1}^{r_1}val_p(|\lambda_{1,i}\overline{\lambda}_{1,i}|_p)$.

\end{remark}

\chapter{The geometric side of the stable trace formula}
\label{FT_stable_geometrique}

\section{Normalization of the Haar measures}
\label{FT_stable_geometrique1}
\index{normalization of Haar measures}

We use the following rules to normalize the Haar measures :
\begin{itemize}
\item[(1)] In the situation of theorem
\ref{th:points_fixes_Kottwitz}, use the normalizations of this theorem.
\item[(2)] Let $\G$ be a connected reductive group over $\Q$. We always
take Haar measures on $\G(\Af)$ such that the volumes of open compact
subgroups are rational numbers. Let $p$ be a prime number such that $\G$
is unramified over $\Q_p$, et let $L$ be a finite unramified extension of
$\Q_p$; then we use the Haar measure on $\G(L)$ such that the volume
of hyperspecial maximal compact subgroups is $1$.
If a Haar measure $dg_f$ on $\G(\Af)$ is fixed, then we use the Haar
measure $dg_\infty$ on $\G(\R)$ such that $dg_fdg_\infty$ is the Tamagawa
measure on $\G(\Ade)$ (cf \cite{O}).
\item[(3)] (cf \cite{K-STF:EST} 5.2)
Let $F$ be a local field of characteristic $0$, $\G$ be a connected reductive
group on $F$ and $\gamma\in\G(F)$ be semi-simple. 
Write $I=\G_\gamma:=\Cent_\G(\gamma)^0$,
and choose Haar measures on $\G(F)$ and $I(F)$.
If $\gamma'\in\G(F)$ is stably conjugate to $\gamma$, then
$I':=\G_{\gamma'}$ is an inner form of $I$, so the measure on $I(F)$
gives a measure on $I'(F)$. When we take the stable orbital integral at
$\gamma$ of a function in $C_c^\infty(\G(F))$, we use these measures on the
centralizers of elements in the stable conjugacy class of
$\gamma$.
\item[(4)] Let $F$ be a local field of characteristic $0$,
$\G$ be a connected reductive group over $F$ and $(\H,s,\eta_0)$ be an
endoscopic triple for $\G$. Let $\gamma_H\in\H(F)$ be semi-simple and
$(\G,\H)$-regular. Assume that there exists an image $\gamma\in
\G(F)$ of $\gamma_H$. Then $I:=\G_\gamma$ is an inner form of
$I_H:=\H_{\gamma_H}$ (\cite{K-STF:EST} 3.1). We always choose Haar
measures on $I(F)$ and $I_H(F)$ that correspond to each other.
\item[(5)] Let $\G$ be a connected reductive group over $\Q$ as in
\ref{points_fixes6}, and let $(\gamma_0;\gamma,\delta)$ be a triple
satisfying conditions (C) of \ref{points_fixes6} and such that the invariant
$\alpha(\gamma_0;\gamma,\delta)$ is trivial. We associate to
$(\gamma_0;\gamma,\delta)$ a group $I$ (connected and reductive over $\Q$)
as in \ref{points_fixes6}. In particular, $I_\R$ is an inner form of
$I(\infty):=\G_{\R,\gamma_0}$. If we already chose a Haar measure on
$I(\R)$ (for example using rule (2), if we have a Haar measure on $I(\Af)$),
then we take the corresponding Haar measure on $I(\infty)(\R)$.

\end{itemize}

\section{Normalization of the transfer factors}
\label{FT_stable_geometrique2}
\index{normalization of transfer factors}

The properties of transfer factors that we will use here are
stated in \cite{K-STF:EST}.
Note that transfer factors have been defined in all generality (for
ordinary endoscopy) by Langlands and Shelstad, cf
\cite{LS1} et \cite{LS2}. The formula of \cite{K-STF:EST} 5.6 is proved in
\cite{LS1} 4.2, and conjecture 5.3 of \cite{K-STF:EST} is proved in
proposition 1 (section 3) of \cite{K-TN}.
\index{transfer factors}

Let $\G$ be one of the unitary groups of \ref{groupes1}, and let
$(\H,s,\eta_0)$ be an elliptic endoscopic triple for $\G$.
Choose a $L$-morphism $\eta:{}^L\H\fl{}^L\G$ extending $\eta_0$.
The local transfer factors associated to $\eta$ are defined only up
to a scalar.

At the infinite place, normalize the transfer factor as in
\cite{K-SVLR} \S7 p 184-185 (this is recalled in
\ref{serie_discrete3}), using the morphism $j$ of
\ref{serie_discrete3} and the Borel subgroup of \ref{rq:le_bon_Borel}.

Let $p$ be a prime number unramified in $E$ (so $\G$ and $\H$ are
unramified over $\Q_p$). Normalize the transfer factor at $p$ as in
\cite{K-SVLR} \S7 p 180-181. If $\eta$ is unramified at $p$, then this
normalization is the one given by the $\Z_p$-structures on $\G$ and $\H$
(it has been defined by Hales in \cite{H} II 7, cf also \cite{Wa3} 4.6).

Choose the transfer factors at other places such that condition 6.10 (b)
of \cite{K-STF:EST} is satisfied. We write $\Delta_{\H,v}^\G$ for the transfer
factors normalized in this way.
\index{$\Delta_{\H,v}^\G$}

Let $\M$ be a cuspidal standard Levi subgroup of $\G$, let $(\M',s_M,
\eta_{M,0})\in\Ell_\G(\M)$, and let $(\H,s,\eta_0)$ be its image in
$\Ell(\G)$. As in \ref{serie_discrete3} and \ref{partie_en_p3}, associate
to $(\M',s_M,\eta_{M,0})$ a cuspidal standard Levi subgroup $\M_H\simeq\M'$
of $\H$ and a $L$-morphism $\eta_M:{}^L\M_H={}^L\M'\fl{}^L\M$ extending
$\eta_{M,0}$ (in \ref{serie_discrete3} and \ref{partie_en_p3}, we took
$\G=\GU(p,q)$, but the general case is similar).
We want to define a normalization of the transfer factors for
$\eta_M$ associated to this data.

At the infinite place, normalize the transfer factor as in
\ref{serie_discrete3}, for the Borel subgroup of $\M$ related to the Borel
subgroup of $\G$ fixed above as in \ref{serie_discrete3}.

If $v$ is a finite place of $\Q$, choose the transfer factor at $v$ that
satisfies the following condition :
\[\Delta_v(\gamma_H,\gamma)_{\M_H}^\M=|D_{\M_H}^\H(\gamma_H)|_v^{1/2}
|D_\M^\G(\gamma)|_v^{-1/2}\Delta_v(\gamma_H,\gamma)_\H^\G,\]
for every $\gamma_H\in\M_H(\Q_v)$ semi-simple $\G$-regular and every
image $\gamma\in\M(\Q_v)$ of $\gamma_H$ (cf \cite{K-NP} lemma 7.5).

We write $\Delta_{\M_H,s_M,v}^\M$ for the transfer factors normalized in
this way.
\index{$\Delta_{\M_H,s_M,v}^\M$}
Note that, if $p$ is unramified in $E$, then
$\Delta_{\M_H,s_M,p}^\M$ depends only on the image of $(\M',s_M,
\eta_{M,0})$ in $\Ell(\M)$ (because it is simply the transfer factor with
the normalization of \cite{K-SVLR} p 180-181, ie, if $\eta_M$ is unramified
at $p$, with the normalization given by the $\Z_p$-structures on
$\M$ and $\M_H$). However, the transfer factors $\Delta_{\M_H,s_M,v}^\M$
at other places may depend on $(\M',s_M,\eta_{M,0})\in\Ell_\G(\M)$, and not
only on its image in $\Ell(\M)$.

\section{Fundamental lemma and transfer}
\label{FT_stable_geometrique3}

We state here the forms of the fundamental lemma and of the transfer
conjecture that we will need in chapter \ref{stabilisation_K}. The local
and adelic stable orbital integrals are defined in
\cite{K-STF:EST} 5.2 et 9.2.
\index{stable orbital integral}
\index{SO@$SO_\gamma$\quad stable orbital integral}

Let $\G$ be one of the groups of \ref{groupes1}, $(\H,s,\eta_0)$ be an elliptic
endoscopic triple of $\G$ and $\eta:{}^L\H\fl{}^L\G$ be a $L$-morphism
extending $\eta_0$.
The transfer conjecture is stated in \cite{K-STF:EST} 5.4 and 5.5. It says
that, for every place $v$ of $\Q$, for every function $f\in C_c^\infty
(\G(\Q_v))$, there exists a function $f^\H\in C_c^\infty(\H(\Q_v))$ such that,
if $\gamma_H\in\H(\Q_v)$ is semi-simple $(\G,\H)$-regular, then
\[SO_{\gamma_H}(f^\H)=\sum_\gamma\Delta_v(\gamma_H,\gamma)e(\G_\gamma)O_\gamma
(f),\]
where the sum is taken over the set of conjugacy classes $\gamma$ in
$\G(\Q_v$) that are images of $\gamma_H$ (so, if $\gamma_H$ has no image
in $\G(\Q_v)$, we want $SO_{\gamma_H}(f^\H)=0$), $\G_\gamma=\Cent_\G(\gamma)^0$
and $e$ is the sign of \cite{K-SC}.
We say that the function $f^\H$ is a transfer of $f$.
\index{transfer conjecture (endoscopic case)}
\index{transfer (endoscopic case)}

The fundamental lemma says that, if $v$ is a finite place where
$\G$ and $\H$ are unramified, if $\eta$ is unramified at $v$ and if
$b:\Hecke(\G(\Q_v),\G(\Z_v))\fl\Hecke(\H(\Q_v),\H(\Z_v))$ is the morphism
induced by $\eta$ (this morphism is made
explicit in \ref{partie_en_p2}), then, for every
$f\in\Hecke(\G(\Q_v),\G(\Z_v))$, $b(f)$ is a transfer of $f$.
\index{fundamental lemma (endoscopic case)}

If $v=\infty$, the transfer conjecture was proved by Shelstad (cf
\cite{Sh}).

For unitary groups, the fundamental lemma and the transfer conjecture were
proved by Laumon-Ngo, Waldspurger and Hales (cf \cite{LN}, \cite{Wa1},
\cite{Wa2} and \cite{Ha}).

We will also need another fundamental lemma. Let $p$ be a finite place where
$\G$ and $\H$ are unramified, $\wp$ be the place of $E$ above $p$ determined
by the fixed embedding $E\fl\overline{\Q}_p$, $j\in\Nat^*$ and $L$ be the
unramified extension of degree $j$ of $E_\wp$ in $\overline{\Q}_p$. Assume
that $\eta$ is unramified at $p$; then it defines a morphism
$b:\Hecke(\G(L),\G(\Of_L))\fl\Hecke(\H(\Q_p),\H(\Z_p))$ (cf \ref{partie_en_p2}
for the definition of $b$ and its description).
The fundamental lemma corresponding to this situation says that,
for every $\phi\in\Hecke(\G(L),\G(\Of_L))$ and every
$\gamma_H\in \H(\Q_p)$ semi-simple and $(\G,\H)$-regular,
\renewcommand\theequation{$*$}
\begin{equation}SO_{\gamma_H}(b(\phi))=\sum_\delta<\alpha_p(\gamma_0;\delta),s>
\Delta_p(\gamma_H,\gamma_0)e(\G_{\delta\sigma})TO_\delta(\phi)
\end{equation}
where the sum is taken over the set of $\sigma$-conjugacy classes $\delta$
in $\G(L)$ such that $N\delta$ is $\G(\overline{\Q}_p)$-conjugate to an
image $\gamma_0\in\G(\Q_p)$ of $\gamma_H$, $\G_{\delta\sigma}$ is
the $\sigma$-centralizer of $\delta$ in $R_{L/\Q_p}\G_L$
and $\alpha_p(\gamma_0;\delta)$
is defined in \cite{K-SVLR} \S7 p 180 (see \cite{K-SVLR} \S7 p 180-181
for more details).
\index{fundamental lemma (twisted case)}
This conjecture (modulo a
calculation of transfer factors) is proved in \cite{Wa3} when $\phi$ is
the unit element of the Hecke algebra $\Hecke(\G(L),\G(\Of_L))$.
The reduction of the general case to this case is done in chapter
\ref{lemme_fondamental}, and the necessary
transfer factor calculation is done in the appendix
by Kottwitz.

\section{A result of Kottwitz}
\label{FT_stable_geometrique4}

We recall here a theorem of Kottwitz about the geometric side of
the stable trace formula for a function that is stable cuspidal at infinity.
The reference for this result is \cite{K-NP}.
\index{stable trace formula (geometric side)}

Let $\G$ be a connected reductive algebraic group over $\Q$. Assume that
$\G$ is cuspidal (cf definition \ref{def:groupe_cuspidal}) and
that the derived group of $\G$ is simply connected. Let $\K_\infty$ be
a maximal compact subgroup of $\G(\R)$.
Let $\G^*$ be a quasi-split inner form of $\G$ over $\Q$,
\index{G@$\G^*$\quad quasi-split inner form of $\G$}
$\overline{\G}$ be an inner form of $\G$ over $\R$ such that
$\overline{\G}/\A_{G,\R}$ is $\R$-anisotropic and $\T_e$ be a maximal
elliptic torus of $\G_\R$. Write
\[\overline{v}(\G)=e(\overline{\G})\vol(\overline{\G}(\R)/\A_G(\R)^0)\]
($e(\overline{\G})$ is the sign associated to $\overline{\G}$ in
\cite{K-SC}), and
\index{vG@$\overline{v}(\G)$}
\index{eG@$e(\G)$}
\[k(\G)=|Im(\H^1(\R,\T_e\cap\G_{der})\fl\H^1(\R,\T_e))|.\]
\index{kG@$k(\G)$}
For every Levi subgroup $\M$ of $\G$, set
\[n_\M^\G=|(\Nor_\G(\M)/\M)(\Q)|.\]
Let $\nu$ be a quasi-character of $\A_G(\R)^0$. Let $\Pi_{temp}(\G(\R),\nu)$
(resp. $\Pi_{disc}(\G(\R),\nu)$) be the subset of $\pi$ in $\Pi_{temp}
(\G(\R))$ (resp. $\Pi_{disc}(\G(\R))$) such that the restriction to
$\A_G(\R)^0$ of the central character of $\pi$ is equal to $\nu$.
\index{$\Pi_{temp}(\G(\R),\nu)$}
\index{$\Pi_{disc}(\G(\R),\nu)$}
Let $C_c^\infty(\G(\R),\nu^{-1})$ be the the set of functions
$f_\infty:\G(\R)\fl\C$ smooth, with compact support modulo
$\A_G(\R)^0$ and such that, for every $(z,g)\in\A_G(\R)^0
\times \G(\R)$, $f_\infty(zg)=\nu^{-1}(z)f_\infty(g)$.
\index{Cc@$C_c^\infty(\G(\R),\nu^{-1})$ }

We say that
$f_\infty\in C_c^\infty(\G(\R),\nu^{-1})$ is \emph{stable cuspidal} if
$f_\infty$ is left and right $\K_\infty$-finite and if the function
\[\Pi_{temp}(\G(\R),\nu)\fl\C,\pi\fle Tr(\pi(f_\infty))\]
vanishes outside $\Pi_{disc}(\G(\R))$ and is constant on the $L$-packets of
$\Pi_{disc}(\G(\R),\nu)$.
\index{stable cuspidal function}

Let $f_\infty\in C_c^\infty(\G(\R),\nu^{-1})$. For every $L$-packet $\Pi$ of
$\Pi_{disc}(\G(\R),\nu)$, write
$Tr(\Pi(f_\infty))=\sum\limits_{\pi\in\Pi}Tr(\pi(f_\infty))$ and
$\Theta_\Pi=\sum\limits_{\pi\in\Pi}\Theta_\pi$.
For every cuspidal Levi subgroup $\M$ of $\G$, define a function
$S\Phi_\M(.,f_\infty)=S\Phi_M^G(.,f_\infty)$ on $\M(\R)$ by the formula :
\[S\Phi_\M(\gamma,f_\infty)=(-1)^{\dim(\A_M/\A_G)}k(\M)k(\G)^{-1}\overline{v}
(\M_\gamma)^{-1}\sum_\Pi\Phi_\M(\gamma^{-1},\Theta_\Pi)Tr(\Pi(f_\infty)),\]
where the sum is taken over the set of $L$-packets $\Pi$ in
$\Pi_{disc}(\G(\R),\nu)$ and $\M_\gamma=\Cent_\M(\gamma)$. Of course,
$S\Phi_\M(\gamma,f_\infty)=0$ unless $\gamma$ is semi-simple and elliptic
in $\M(\R)$. If $\M$ is a Levi subgroup of $\G$ that is not cuspidal, set
$S\Phi_M^G=0$.
\index{SPhi@$S\Phi_M^G$}

Let $f:\G(\Ade)\fl\C$. Assume that $f=f^\infty f_\infty$, with $f^\infty\in
C_c^\infty(\G(\Af))$ and $f_\infty\in C_c^\infty(\G(\R),\nu^{-1})$.
For every Levi subgroup $\M$ of $\G$, set
\[ST_M^G(f)=\tau(\M)\sum_\gamma SO_\gamma(f_\M^\infty)
S\Phi_\M(\gamma,f_\infty),\]
\index{STMG@$ST_M^G$}
where the sum is taken over the set of stable conjugacy classes $\gamma$
in $\M(\Q)$ that are semi-simple and elliptic in $\M(\R)$, and
$f^\infty_\M$ is the constant term of $f^\infty$ at $\M$ (the constant term
depends on the choice of a parabolic subgroup of $\G$ with Levi subgroup
$\M$, but its integral orbitals do not).
Set
\[ST^G(f)=\sum_\M (n_M^G)^{-1}ST_M^G(f),\]
where the sum is taken over the set of $\G(\Q)$-conjugacy classes $\M$
of Levi subgroups of $\G$.
\index{STG@$ST^G$\quad stable trace formule on $\G$}

Let $T^G$ be the distribution of Arthur's invariant trace formula.
\index{TG@$T^G$\quad invariant trace formule on $\G$}
For every $(\H,s,\eta_0)\in\Ell(\G)$ (cf \ref{groupes4}), fix
a $L$-morphism $\eta:{}^L\H\fl{}^L\G$ extending $\eta_0$, and let
\[\iota(\G,\H)=\tau(\G)\tau(\H)^{-1}|\Lambda(\H,s,\eta_0)|^{-1}.\]
Kottwitz proved the following theorem in \cite{K-NP} (theorem 5.1) :
\index{$\iota(\G,\H)$}

\begin{theorem}\label{th:FT_stable_geometrique}
Let $f=f^\infty f_\infty$ be as above. Assume that $f_\infty$ is stable
cuspidal and that, for every $(\H,s,\eta_0)\in\Ell(\G)$, there exists
a transfer $f^\H$ of $f$. Then :
\[T^G(f)=\sum_{(\H,s,\eta_0)\in\Ell(\G)}\iota(\G,\H)ST^H(f^{\H}).\]

\end{theorem}

We calculate $k(\G)$ for $\G$ a unitary group.

\begin{lemma}\label{lemme:calcul_k} Let $p_1,\dots,p_r,q_1,\dots,q_r
\in\Nat$ such that $p_i+q_i\geq 1$ for $1\leq i\leq r$; write 
$n_i=p_i+q_i$, $n=n_1+\dots+n_r$
and $\G=\G(\U(p_1,q_1)\times\dots\times\U(p_r,q_r))$. Then
\[k(\G)=2^{n-r-1}\]
if all the $n_i$ are even, and
\[k(\G)=2^{n-r}\]
otherwise.

\end{lemma}

In particular, $k(R_{E/\Q}\Gr_m)=k(\GU(1))=1$.

\begin{proof} Write $\Gamma(\infty)=\Gal(\C/\R)$.
In \ref{serie_discrete1}, we defined an elliptic maximal torus
$\T_e$ of $\G$ and an isomorphism $\T_e\iso\G(\U(1)^n)$. Tate-Nakayama duality
induces an isomorphism between the dual of
$Im(\Ho^1(\R,T_e\cap\G_{der})\fl\Ho^1(\R,\T_e))$ and
$Im(\pi_0(\widehat{\T}_e^{\Gamma(\infty)})\fl\pi_0((\widehat{\T}_e/Z(\widehat
{\G})^{\Gamma(\infty)})))$ (cf \cite{K-STF:CTT} 7.9). Moreover, there is an
exact sequence
\[(X_*(\widehat{\T}_e/Z(\widehat{\G})))^{\Gamma(\infty)}\fl\pi_0(Z(\widehat{\G}
)^{\Gamma(\infty)})\fl\pi_0(\widehat{\T}_e^{\Gamma(\infty)})\fl\pi_0((\widehat
{\T}_e/Z(\widehat{\G}))^{\Gamma(\infty)})\]
(cf \cite{K-STF:CTT} 2.3), and $(X_*(\widehat{\T}_e/Z(\widehat{\G})))^{\Gamma(
\infty)}=0$ because $\T_e$ is elliptic, hence
\[k(\G)=|\pi_0(\widehat{\T}_e^{\Gamma(\infty)})||\pi_0(Z(\widehat{\G})^
{\Gamma(\infty)})|^{-1}.\]
Of course, $\widehat{\T}_e^{\Gamma(\infty)}=\widehat{\T}_e^{\Gal(E/\Q)}$ et
$Z(\widehat{\G})^{\Gamma(\infty)}=Z(\widehat{\G})^{\Gal(E/\Q)}$. We already
calculated these groups in (i) of lemma \ref{lemme:Tamagawa}. This implies
the result.

\end{proof}

\begin{remark}\label{rq:k_tau} Note that
$k(\G)\tau(\G)=2^{n-1}$.

\end{remark}

\chapter{Stabilization of the fixed point formula}
\label{stabilisation_K}

To simplify the notations, we suppose in this chapter that the group
$\G$ is $\GU(p,q)$, but all the results generalize in an obvious
way to the groups $\G(\U(p_1,q_1)\times\dots\times\U(p_r,q_r))$.

\section{Preliminary simplifications}
\label{stabilisation_prelim}

We first rewrite the fixed point formula using
proposition \ref{prop:calcul_Phi_M}. 

Notations are as in chapters \ref{points_fixes} (especially section
\ref{points_fixes7}), \ref{groupes} and \ref{serie_discrete}. Fix 
$p,q\in\Nat$ such that $n:=p+q\geq 1$, and let $\G=\GU(p,q)$. We may
assume, and we do, that $p\geq q$.
Let $V$ be an irreducible algebraic representations of $\G_\C$ and
$\varphi:W_\R\fl{}^L\G$ be an elliptic Langlands parameter corresponding
to $V^*$ as in proposition \ref{prop:calcul_Phi_M}.
Let $K\subset\C$ be a number field such that $V$ is defined over $K$.
Set
\[\Theta=(-1)^{q(\G)}S\Theta_{\varphi}.\]
Fix $g\in\G(\Af)$, $\K,\K'\subset\G(\Af)$,
$j\in\Nat^*$, prime numbers $p$ and $\ell$ and a place $\lambda$ of $K$
above $\ell$ as in \ref{points_fixes4} and \ref{points_fixes5}. We get
a cohomological correspondence
\[\overline{u}_j:(\Phi^{j}\overline{T}_g)^*IC^{\K}V\fl\overline{T}_1^!IC^{\K}V.
\]

\begin{proposition}\label{prop:FPF_simplifiee}
Write $\G_0=\M_\varnothing=\Pa_\varnothing=\G$ and $\Le_\varnothing=\{1\}$.
For every $s\in\{0,\dots,q\}$, set
\begin{flushleft}$\displaystyle{\Tr_s=(-1)^sm_s(n_{M_S}^G)^{-1}\chi(\Le_S)
\sum_{\gamma_L\in\Le_S(\Q)}
\sum_{(\gamma_0;\gamma,\delta)\in C'_{\G_s,j}}c(
\gamma_0;\gamma,\delta)O_{\gamma_L}(\ungras_{\Le_S(\Z_p)})
O_{\gamma_L\gamma}(f^{\infty,p}_{\M_S})}$\end{flushleft}
\begin{flushright}$\displaystyle{\delta_{P_S(\Q_p)}^{1/2}(\gamma_0)TO_\delta(
\phi_j^{\G_s})\Phi_{\M_S}^{\G}((\gamma_L\gamma_0)^{-1},\Theta),}$
\end{flushright}
where $S=\{1,2,\dots,s\}$ and where $m_s=1$ if $s<n/2$, and $m_{n/2}=
|\X_{n/2}|=2$ if $n$ is even.
Then, if $j$ is big enough :
\[\Tr(\overline{u}_j,R\Gamma(M^{\K}(\G,\X)_\Fi,(IC^{\K}V)_{\Fi}))=
\sum_{s=0}^q \Tr_s.\]
If $g=1$ and $\K=\K'$, then the above formula is true for every
$j\in\Nat^*$.

\end{proposition}

\begin{proof} Let $m\in\Z$ be the weight of $V$ as a representation
of $\G$ (cf \ref{points_fixes3}).
For every $r\in\{1,\dots,q\}$, set $t_r=r(r-n)$. By
proposition \ref{prop:IC_est_pondere}, there is a canonical isomorphism
\[IC^{\K}V\simeq W^{\geq t_1+1,\dots,\geq t_q+1}V.\]
Write
\[\Tr=\Tr(\overline{u}_j,R\Gamma(M^{\K}(\G,\X)_\Fi,(IC^{\K}V)_{\Fi})).\]
If $j$ is big enough, then, by theorem \ref{th:points_fixes_moi},
\[\Tr=\Tr_G+\sum_\Pa \Tr_P,\]
where the sum is taken over the set of standard parabolic subgroups of $\G$.
Set $\Tr'_0=\Tr_G$ and, for every $s\in\{1,\dots,q\}$,
\[\Tr'_s=\sum_{{S'\subset\{1,\dots,s\}}\atop{S'\ni s}}\Tr_{P_{S'}}.\]
We want to show that $\Tr'_s=\Tr_s$.
For $s=0$, this comes from the fact that, for every semi-simple
$\gamma_0\in\G(\Q)$ that is elliptic in $\G(\R)$,
\[\Tr(\gamma_0,V)=\Theta(\gamma_0^{-1})=\Phi_{\G}^{\G}(\gamma_0^{-1},\Theta).\]

Let $s\in\{1,\dots,q\}$; write $S=\{1,\dots,s\}$. Let $S'\subset S$ be
such that $s\in S'$. Then, up to $\Le_{S'}(\Q)$-conjugacy, the only
cuspidal Levi subgroup of $\Le_{S'}$ is $\Le_S=(R_{E/\Q}\Gr_m)^s$. Hence
\begin{flushleft}$\displaystyle{\Tr_{P_{S'}}=(-1)^{\dim(\A_{L_{S}}/
\A_{L_{S'}})}m_{P_{S'}}
(n_{L_S}^{L_{S'}})^{-1}\chi(\Le_S)\sum_{\gamma_L\in\Le_S(\Q)}|D_{\Le_S}^{\Le_
{S'}}(\gamma_L)|^{1/2}\sum_{(\gamma_0;\gamma,\delta)\in C'_{\G_s,j}}
c(\gamma_0;\gamma,\delta)}$\end{flushleft}
\begin{flushright}$\displaystyle{O_{\gamma_L\gamma}((f^{\infty,p})_{\M_S})
O_{\gamma_L}(\ungras_{\Le_S(\Z_p)})\delta_{P_{S'}(\Q_p)}^{1/2}(\gamma_0)
TO_\delta(\phi_j^{\G_s})
\delta_{P_{S'}(\R)}^{1/2}(\gamma_L\gamma_0)
L_{S'}(\gamma_L\gamma_0),}$
\end{flushright}
where
\[L_{S'}(\gamma_L\gamma_0)=\Tr(\gamma_L\gamma_0,R\Gamma(Lie(\N_{S'}),V)_{> t_r
+m,r\in S'}).\]
As $\gamma_0$ is in $\G_s(\Q)$,
\[\delta_{P_{S'}(\Q_p)}(\gamma_0)=\delta_{P_S(\Q_p)}(\gamma_0).\]
Hence : 
\begin{flushleft}$\displaystyle{\Tr'_s=
m_s\chi(\Le_S)\sum_{\gamma_L\in\Le_S(\Q)}
\sum_{(\gamma_0;\gamma,\delta)\in C'_{\G_s,j}}c(\gamma_0;\gamma,\delta)
O_{\gamma_L\gamma}((f^{\infty,p})_{\M_S})O_{\gamma_L}(\ungras_{\Le_S(\Z_p)})
}$\end{flushleft}
\begin{flushright}$\displaystyle{\delta_{P_S(\Q_p)}^{1/2}(\gamma_0)TO_\delta(
\phi_j^{\G_s})\sum_{{S'\subset S}\atop{S'\ni s}}(-1)^{\dim(\A_{L_S}/\A_{L_{S'}}
)}(n_{L_S}^{L_{S'}})^{-1}|D_{\Le_S}^{\Le_{S'}}(\gamma_L)|^{1/2}
\delta_{P_{S'}(\R)}^{1/2}(\gamma_L\gamma_0)L_{S'}(\gamma_L\gamma_0).}$
\end{flushright}

Consider the action of the group $\Sgoth_s$ on $\M_S$ defined in
\ref{serie_discrete4} (so $\Sgoth_s$ acts on $\Le_S=(R_{E/\Q}\Gr_m)^s$ by
permuting the factors, and it acts trivially on $\G_s$).
Then
\[c(\gamma_0;\gamma,\delta)O_{\gamma_L\gamma}((f^{\infty,p})_{\M_S})
O_{\gamma_L}(\ungras_{\Le_S(\Z_p)})\delta_{P_S(\Q_p)}^{1/2}(\gamma_0)
TO_\delta(\phi_j^{\G_s})\]
is invariant by the action of $\Sgoth_s$. Let $\gamma_M\in\M_S(\Q)$ 
be semi-simple and elliptic in $\M_S(\R)$. Write $\gamma_M=\gamma_L
\gamma_0$, with $\gamma_L\in\Le_S(\Q)$ and $\gamma_0\in\G_s(\Q)$. Note that,
for every $S'\subset S$ such that $s\in S'$, $\dim(\A_{L_{S'}}/\A_
{L_S})=\dim(\A_{M_{S'}}/\A_{M_S})$ and
$D_{\M_S}^{\M_{S'}}(\gamma_M)=D_{\Le_{S}}^{\Le_{S'}}(\gamma_L)$.
If $\gamma_M$ satisfies the condition of part (ii) of
proposition \ref{prop:calcul_Phi_M}, then, by this proposition,
\begin{flushleft}$\displaystyle{\sum_{\gamma'\in\Sgoth_s.\gamma_M}\sum_{{S'
\subset S}\atop{S'\ni s}}(-1)^{\dim(\A_{L_S}/\A_{L_{S'}})}(n_{L_S}^{L_{S'}})
^{-1}|D_{\M_S}^{\M_{S'}}(\gamma')|^{1/2}\delta_{P_{S'}(\R)}^{1/2}(\gamma')
L_{S'}(\gamma')}$\end{flushleft}
\begin{flushright}$\displaystyle{=(-1)^s(n_{M_S}^G)^{-1}\sum_{\gamma'\in
\Sgoth_s.\gamma_M}\Phi_{\M_S}^{\G}({\gamma'}^{-1},\Theta),}$\end{flushright}
because the function $\Phi_{\M_S}^{\G}(.,\Theta)$ is invariant by the action
of $\Sgoth_s$, and $n_{M_S}^G=2^ss!=2^s|\Sgoth_s|$. Moreover, also by
proposition \ref{prop:calcul_Phi_M}, for every
$\gamma_L\in\Le_S(\R)$ and $(\gamma_0;\gamma,\delta)\in C_{\G_s,j}-
C'_{\G_s,j}$, $\Phi_{\M_S}^\G((\gamma_L\gamma_0)^{-1},\Theta)=0$ (because
$c(\gamma_L\gamma_0)=c(\gamma_0)<0$).

To finish the proof of the proposition, it is enough to show that, if
$j$ is big enough, then, for every $\gamma_L\in\Le_S(\Q)$ and
$(\gamma_0;\gamma,\delta)\in C'_{\G_s,j}$ such that
$O_{\gamma_L}(\ungras_{\Le_S(\Z_p)})
TO_\delta(\phi_j^{\G_s})O_{\gamma_L\gamma}((f^{\infty,p})_{\M_S})\not=0$,
the element $\gamma_L\gamma_0$ of $\M_S(\Q)$ satisfies the condition of
part (ii) of proposition \ref{prop:calcul_Phi_M}.

Let $\Sigma$ be the set of $(\gamma_L,\gamma_0)\in\M_S(\Q)=(E^\times)^s
\times\G_s(\Q)$ such that there exists $(\gamma,\delta)\in
\G_s(\Af^p)\times\G_s(L)$ with $(\gamma_0;\gamma,\delta)\in C'_{\G_s,j}$
and $O_{\gamma_L}(\ungras_{\Le_S(\Z_p)})TO_\delta(\phi_j^{\G_s})O_{\gamma_L
\gamma}((f^{\infty,p})_{\M_S})\not=0$. By remark
\ref{rq:support_terme_constant}, the function $\gamma_M\fle O_{\gamma_M}(
(f^{\infty,p})_{\M_S})$ on $\M_S(\Af^p)$ has compact support modulo
conjugacy. So there exist $C_1,C_2\in\R^{+*}$ such that, for every
$(\gamma_L=(\lambda_1,\dots,\lambda_s),\gamma_0)\in\Sigma$,
\[|c(\gamma_0)|_{\Af^p}\geq C_1\]
\[\sup_{1\leq r\leq s}|\lambda_r|_{\Af^p}^2\geq C_2.\]
On the other hand, if
$(\gamma_L=(\lambda_1,\dots,\lambda_s),\gamma_0)\in\Sigma$, then
$|\lambda_r|_{\Q_p}=1$ for $1\leq r\leq s$ (because $\gamma_L\in\Le_S(\Z_p)$),
and there exists $\delta\in\G_s(L)$ such that
$TO_{\delta}(\phi_j^{\G_s})\not=0$ and that $N\delta$ is
$\G_s(\overline{\Q}_p)$-conjugate to $\gamma_0$; this implies that
\[|c(\gamma_0)|_{\Q_p}=|c(\delta)|_L=p^d\geq p^j,\]
because $d=j$ or $2j$ .
(If $TO_\delta(\phi_j^{\G_s})\not=0$, then $\delta$ is $\sigma$-conjugate
to an element $\delta'$ of $\G(\Of_L)\mu_{\G_s}(\varpi_L^{-1})\G(\Of_L)$,
so $|c(\delta)|_L=|c(\delta')|_L=p^d$.)
Finally, if $(\gamma_L=(\lambda_1,\dots,\lambda_s),\gamma_0)\in\Sigma$, then
\[|c(\gamma_0)|_\infty\sup_{1\leq r\leq s}|\lambda_r|_\infty^2=
|c(\gamma_0)|_{\Q_p}^{-1}|c(\gamma_0)|_{\Af^p}^{-1}\sup_{1\leq r\leq s}
|\lambda_r|^{-2}_{\Af^p}\leq p^{-j}C_1^{-1}C_2^{-1}.\]
Moreover, if $(\gamma_0;\gamma,\delta)\in C'_{\G_s,j}$, then $c(\gamma_0)>0$.
Hence, if $j$ is such that $p^{j}C_1C_2\geq 1$, then all the elements of
$\Sigma$ satisfy the condition of part (ii) of
proposition \ref{prop:calcul_Phi_M}.

Assume that $g=1$ and $\K=\K'$. Then theorem
\ref{th:points_fixes_moi} is true for every $j\in\Nat^*$. Moreover, by
remark \ref{rq:support_terme_constant}, the support of the function
$\gamma_M\fle O_{\gamma_M}((f^{\infty,p})_{\M_S})$ is contained in the union
of the conjugates of a finite union of open compact subgroups of
$\M_S(\Af^p)$, so we may take $C_1=C_2=1$, and every $j\in\Nat^*$ satisfies
$p^{j}C_1C_2\geq 1$.

\end{proof}

\section{Stabilization of the elliptic part, after Kottwitz}
\label{stabilisation_K1}

In \cite{K-SVLR}, Kottwitz stabilized the elliptic part of the fixed
point formula (ie, with the notations used here, the term
$\Tr_0$). We will recall his result in this section, and, in the next section,
apply his method to the terms $\Tr_s$, $s\in\{1,\dots,q\}$.

For every $(\H,s,\eta_0)\in\Ell(\G)$, fix a $L$-morphism $\eta:{}^L\H\fl
{}^L\G$ extending $\eta_0$ (in this section, we make the $L$-groups with
$W_\Q$).

\begin{theorem}(\cite{K-SVLR} 7.2)\label{th:stab_partie_elliptique} There
is an equality
\[\Tr_0=\sum_{(\H,s,\eta_0)\in\Ell(\G)}\iota(\G,\H)\tau(\H)\sum_{\gamma_H}
SO_{\gamma_H}(f_\H^{(j)}),\]
where the second sum is taken over the set of
semi-simple stable conjugacy classes in $\H(\Q)$ that are elliptic
in $\H(\R)$.

\end{theorem}

We have to explain the notation.
Let $(\H,s,\eta_0)\in\Ell(\G)$.
As in \ref{FT_stable_geometrique4}, write
\[\iota(\G,\H)=\tau(\G)\tau(\H)^{-1}|\Lambda(\H,s,\eta_0)|^{-1}.\]
The function $f_\H^{(j)}$ is a function in $C^\infty(\H(\Ade))$ such that
$f_\H^{(j)}=f_\H^{\infty,p}f_{\H,p}^{(j)}f_{\H,\infty}$, with $f_\H^{\infty,p}
\in
C_c^\infty(\H(\Af^p))$, $f_{\H,p}^{(j)}\in C_c^\infty(\H(\Q_p))$ and
$f_{\H,\infty}
\in C^\infty(\H(\R))$.
\index{fH@$f_\H^{(j)}$}

The first function $f_\H^{\infty,p}$ is simply a transfer of
$f^{\infty,p}\in C_c^\infty(\G(\Af^p))$.

Use the notations of the last subsection of \ref{partie_en_p2},
and set $\eta_p=\eta_{|\widehat{\H}\rtimes W_{\Q_p}}$. Then
$f_{\H,p}^{(j)}\in C_c^\infty(\H(\Q_p))$ is equal to the twisted transfer
$b_{\til{\eta}_p}(\phi_j^\G)$.

For every elliptic Langlands parameter
$\varphi_H:W_\R\fl\widehat{\H}\rtimes W_\R$, set
\[f_{\varphi_H}=d(\H)^{-1}\sum_{\pi\in\Pi(\varphi_H)}f_\pi\]
(with the notations of \ref{serie_discrete1}).
\index{fvarphi@$f_\varphi$}

Let $\B$ be the standard Borel subgroup of $\G_\C$ (ie the subgroup of upper
triangular matrices). It determines as in
\ref{serie_discrete3} a subset $\Omega_*\subset\Omega_G$  and a
bijection $\Phi_H(\varphi)\iso\Omega_*$, 
$\varphi_H\fle\omega_*(\varphi_H)$. Take
\[f_{\H,\infty}=<\mu_\G,s>(-1)^{q(\G)}\sum_{\varphi_H\in\Phi_H
(\varphi)}\det(\omega_*(\varphi_H))f_{\varphi_H},\]
where $\mu_\G$ is the cocharacter of $\G_\C$ associated to the Shimura datum
as in \ref{groupes1}.
\index{$\mu_\G$}
Note that, as suggested by the notation, $f_{\H,p}^{(j)}$ is the only part
of $f_\H^{(j)}$ that depends on $j$.

\begin{remark}\label{rq:G_H_regulier} In theorem 7.2 of
\cite{K-SVLR}, the second sum is taken over the set of semi-simple elliptic
stable conjugacy classes that are
$(\G,\H)$-regular. But proposition \ref{prop:transfert_caracteres}
implies that $SO_{\gamma_H}(f_{\H,\infty})=0$ if $\gamma_H\in\H(\R)$
is not $(\G,\H)$-regular (cf remark \ref{rq:M_M_H_regulier}).

\end{remark}

\section{Stabilization of the other terms}
\label{stabilisation_K2}

The stabilization process that we follow here is mainly due to Kottwitz,
and explained (in a more general situation) in \cite{K-NP}. The
differences between the stabilization of the trace formula (in \cite{K-NP})
and the stabilization of the fixed point formula considered here
are concentrated at the places $p$ and $\infty$. In particular,
the vanishing of part of the contribution of the linear part of
Levi subgroups is particular to the stabilization of the fixed point
formula.

We will use freely the definitions and notations of \ref{groupes4}.

\begin{theorem}\label{th:adieu_partie_lineaire}
\begin{itemize}
\item[(i)] Let $r\in\{1,\dots,
q\}$. Write $\M=\M_{\{1,\dots,r\}}$. Then :
\[\Tr_r=(n_M^G)^{-1}\sum_{(\M',s_M,\eta_{M,0})\in\Ell_\G
(\M)}\tau(\G)\tau(\H)^{-1}|\Lambda_\G(\M',s_M,\eta_{M,0})|^{-1}
ST_{M'}^H(f_\H^{(j)}),\]
where, for every $(\M',s_M,\eta_{M,0})\in\Ell_\G(\M)$, $(\H,s,\eta_0)$
is the corresponding element of $\Ell(\G)$.
\item[(ii)] Write $\G^*=\GU^*(n)$. Let $r\in\Nat$ such that $r\leq n/2$.
Denote by $\M^*$ the standard Levi subgroup of $\G^*$ that corresponds to
$\{1,\dots,r\}$ (as in section \ref{groupes2}). If $r\geq q+1$, then
\[\sum_{(\M',s_M,\eta_{M,0})\in\Ell_{\G^*}
(\M^*)}\tau(\H)^{-1}|\Lambda_{\G^*}(\M',s_M,\eta_{M,0})|^{-1}
ST_{M'}^H(f_\H^{(j)})=0.\]

\end{itemize}
\end{theorem}

\begin{corollary}\label{cor:stab_FT_IC} If $j$ is big enough, then
\[\Tr(\overline{u}_j,R\Gamma(M^{\K}(\G,\X)^*_\Fi,IC^{\K}V_\Fi))=
\sum_{(\H,s,\eta_0)\in\Ell(\G)}\iota(\G,\H)ST^{H}(f_{\H}^{(j)}).\]
If $g=1$ and $\K=\K'$, then this formula is true for every
$j\in\Nat^*$.

\end{corollary}

\begin{remark} Let $\Hecke_K=\Hecke(\G(\Af),\K)$.
Define an object $W_\lambda$ of the
Grothendieck group of representations of $\Hecke_K\times\Gal(\overline{\Q}/F)$
in a finite dimensional $K_\lambda$-vector space by
\[W_\lambda=\sum_{i\geq 0}(-1)^i[\Ho^i(M^{\K}(\G,\X)^*_{\overline{\Q}},IC^{\K}
V_{\overline{\Q}})].\]
Then, for every $j\in\Nat^*$,
\[\Tr(\overline{u}_j,R\Gamma(M^{\K}(\G,\X)^*_\Fi,IC^{\K}V_\Fi))=
\Tr(\Phi_\wp^j h,W_\lambda),\]
where $\Phi_\wp\in\Gal(\overline{\Q}/F)$ is a lifting of the geometric
Frobenius at $\wp$ (the fixed place of $F$ above $p$) and $h=
\vol(\K)^{-1}\ungras_{\K g\K}$.

So the corollary implies that : For every $f^\infty\in\Hecke_\K$ such that
$f^\infty=f^{\infty,p}\ungras_{\G(\Z_p)}$ and for every $j$ big enough (in a
way that may depend on $f^\infty$),
\[\Tr(\Phi_\wp^jf^\infty,W_\lambda)=\sum_{(\H,s,\eta_0)\in\Ell(\G)}
\iota(\G,\H)ST^{H}(f_{\H}^{(j)}),\]
where, for every $(\H,s,\eta_0)\in\Ell(\G)$,
$f_{\H,p}^{(j)}$ and $f_{\H,\infty}$ are defined as before, and
$f_\H^{\infty,p}$ is a transfer of $f^{\infty,p}$.

\end{remark}

\begin{proof}
The corollary is an immediate consequence of
theorem \ref{th:stab_partie_elliptique}, of the above theorem, of
proposition \ref{prop:FPF_simplifiee} and of lemma
\ref{lemme:Levi_endoscopiques}.

\end{proof}

\begin{prooft} As $j$ is fixed, we omit the subscripts ``$(j)$'' in this proof.

We prove (i). As $\M$ is a proper Levi subgroup of $\G$,
$|\Lambda_\G(\M',s_M,\eta_{M,0})|=1$ for every
$(\M',s_M,\eta_{M,0})\in\Ell_\G(\M)$ (lemma
\ref{lemme:Levi_endoscopiques_GU}).
Write, with the notations of the theorem,
\[\Tr'_M=\sum_{(\M',s_M,\eta_{M,0})\in\Ell_\G(\M)}\tau(\H)^
{-1}ST_{M'}^H(f_\H).\]
By the definition of $ST_{M'}^H$ in \ref{FT_stable_geometrique4},
\[\Tr'_M=\sum_{(\M',s_M,\eta_{M,0})\in\Ell_\G(\M)}\tau(\H)^{-1}
\tau(\M')\sum_{\gamma'}SO_{\gamma'}((f^\infty_\H)_{\M'})S\Phi_{M'}^H(\gamma',
f_{\H,\infty}),\]
where the second sum is taken over the set of semi-simple stable conjugacy
classes of $\M'(\Q)$ that are elliptic in $\M'(\R)$.
By proposition \ref{prop:transfert_caracteres}
(and remark \ref{rq:M_M_H_regulier}), the terms indexed by a stable conjugacy
class $\gamma'$ that is not $(\M,\M')$-regular all vanish.
By lemma \ref{lemme:param_groupes_endoscopiques},
\[\Tr'_M=\sum_{\gamma_M}\sum_{\kappa\in\Kgoth_\G(I/\Q)_e}\tau(\M')\tau(\H)^{-1}
\psi(\gamma_M,\kappa),\]
where :
\begin{bulletlist}
\item The first sum is taken over the set of semi-simple stable conjugacy
classes $\gamma_M$ in $\M(\Q)$ that are elliptic in $\M(\R)$.
\item $I=\M_{\gamma_M}$, and $\Kgoth_\G(I/\Q)$ is defined above lemma
\ref{lemme:param_groupes_endoscopiques}.
\item  Let $\gamma_M$ be as above and $\kappa\in\Kgoth_\G(I/\Q)$. Let
$(\M',s_M,\eta_{M,0},\gamma')$ be an endoscopic $\G$-quadruple associated to
$\kappa$ by lemma
\ref{lemme:param_groupes_endoscopiques}. The subset
$\Kgoth_\G(I/\Q)_e$ of $\Kgoth_\G(I/\Q)$ is the set of $\kappa$ such that
$(\M',s_M,\eta_{M,0})\in\Ell_\G(\M)$. If $\kappa\in\Kgoth_\G(I/\Q)_e$, set
\[\psi(\gamma_M,\kappa)=SO_{\gamma'}
((f_\H^\infty)_{\M'})S\Phi_{M'}^H(\gamma',f_{\H,\infty}),\]
where $(\H,s,\eta_0)$ is as before the image of $(\M',s_M,\eta_{M,0})$ in
$\Ell(\G)$.
\index{KGIQe@$\Kgoth_\G(I/\Q)_e$}

\end{bulletlist}

Fix $\gamma_M$, $\kappa\in\Kgoth_\G(I/\Q)_e$ and $(\M',s_M,\eta_{M,0},
\gamma')$ as above. Let $(\H,s,\eta_0)$ be the element of $\Ell(\G)$
associated to $(\M',s_M,\eta_{M,0})$, and define $s'_M$ as in
\ref{partie_en_p3}.
Let $\gamma_0$ be the component of $\gamma_M$ in the Hermitian part
$\G_r(\Q)$ of $\M(\Q)$.
We want to calculate $\psi(\gamma_M,\kappa)$.
This number is the product of three terms :
\begin{itemize}
\item[(a)] Outside of $p$ and $\infty$ : by (ii) of lemma
\ref{lemme:transfert_terme_constant},
\[SO_{\gamma'}((f^{p,\infty}_\H)_{\M'})=\sum_{\gamma}\Delta_{\M',s_M}^
{\M,\infty,p}(\gamma',\gamma)e(\gamma)O_\gamma(f^{\infty,p}_\M),\]
where the sum is taken over the set of semi-simple conjugacy classes
$\gamma=(\gamma_v)_{v\not=p,\infty}$ in $\M(\Af^p)$ such that $\gamma_v$ is
stably conjugate to $\gamma_M$ for every $v$, and $e(\gamma)=\prod\limits_
{v\not=p,\infty}e(\M_{\Q_v,\gamma_v})$.
By \cite{K-STF:EST} 5.6, this sum is equal to
\[\Delta_{\M',s_M}^{\M,\infty,p}(\gamma',\gamma_M)\sum_\gamma<\alpha(\gamma_M,
\gamma),\kappa>e(\gamma)O_\gamma(f^{\infty,p}_\M),\]
where the sum is over the same set as before and
$\alpha(\gamma_M,\gamma)$ is the invariant denoted by $inv(\gamma_M,\gamma)$
in \cite{K-STF:EST} (the article \cite{K-STF:EST} is only stating a conjecture,
but this conjecture has been proved since, see
\ref{FT_stable_geometrique3} for explanations).
Moreover, as the linear part of $\M$ is isomorphic to $(R_{E/\Q}\Gr_m)^r$,
we may replace $\alpha(\gamma_M,\gamma)$ with $\alpha(\gamma_0,\gamma_h)$,
where $\gamma_h$ is the component of $\gamma$ in $\G_r(\Af^p)$.

\item[(b)] At $p$ : By corollary \ref{cor:identite_en_p} (and with the
notations of \ref{partie_en_p3}), 
$SO_{\gamma'}((f_{\H,p})_{\M'})$ is equal to
\[<\mu,s><\mu,s'_M>\varepsilon_{s_M}(\gamma_M)
\sum_\delta<\alpha_p(\gamma_0,\delta_h),s_M>\Delta_
{\M',s_M,p}^{\M}(\gamma',\gamma_M)e(\delta)TO_\delta((\phi_j^\G)_\M),\]
and, if $\gamma_M\in\Le_r(\Z_p)\G_r(\Q_p)$, to
\[<\mu,s><\mu,s'_M>\delta_{\Pa_r(\Q_p)}^{1/2}(\gamma_M)\sum_
\delta<\alpha_p(\gamma_0,\delta_h),s_M>\Delta_{\M',s_M,p}^{\M}(\gamma',
\gamma_M)e(\delta)TO_\delta(\ungras_{\Le_r(\Of_L)}\times\phi_j^{\G_r}).\]
Both sums are taken over the set of $\sigma$-conjugacy classes $\delta$ in
$\M(L)$ such that $\gamma_M\in\Norme\delta$, and
$\delta_h$ is the component of $\delta$ in $\G_r(L)$.

\item[(c)] At $\infty$ : By the definitions and proposition
\ref{prop:transfert_caracteres}, $S\Phi_{\M'}^\H(\gamma',f_{H,\infty})$ is
equal to
\[<\mu,s>(-1)^{\dim(\A_M/\A_G)}k(\M')k(\H)^{-1}
\overline{v}(I)^{-1}\Delta_{\M',s_M,\infty}^{\M}(\gamma',\gamma_M)
\Phi_\M^\G(\gamma_M^{-1},\Theta)\]
(note that $\A_{M'}\simeq\A_M$ and $\A_H\simeq\A_G$ because the endoscopic
data $(\H,s,\eta_0)$ and $(\M',s_M,\eta_{M,0})$ are elliptic,
and that $\overline{v}(I)=\overline{v}(\M'_{\gamma'})$ because $I$ is
an inner form of $\M'_{\gamma'}$ by \cite{K-STF:EST} 3.1).

\end{itemize}

Finally, we find that $\psi(\gamma_M,\kappa)$ is equal to
\begin{flushleft}$\displaystyle{(-1)^{\dim(\A_M/\A_G)}
\varepsilon_{s_M}(\gamma_M)k(\M')k(\H)^{-1}}$\end{flushleft}
\begin{flushright}$\displaystyle{\sum_{(\gamma,\delta)}<\alpha(\gamma_0;\gamma
_h,\delta_h),\kappa>e(\gamma)e(\delta)\overline{v}(I)^{-1}O_\gamma(f^{\infty,p}
_\M)TO_\delta((\phi_j^\G)_\M)\Phi_\M^\G(\gamma_M^{-1},\Theta),}
$\end{flushright}
where the sum is taken over the set of equivalence classes of
$(\gamma,\delta)\in\M(\Af^p)\times\M(L)$ such that $(\gamma_M;\gamma,\delta)$
satisfies conditions (C) of \ref{points_fixes6} and,
if $(\gamma_h,\delta_h)$
is the component of $(\gamma,\delta)$ in $\G_r(\Af^p)\times\G_r(L)$, then
$\alpha(\gamma_0;\gamma_h,\delta_h)\in\Kgoth(I/\Q)^D$ is the invariant
associated to $(\gamma_0;\gamma_h,\delta_h)$ by Kottwitz in
\cite{K-SVLR} \S2 (it is easy to check that $(\gamma_0;\gamma_h,\delta_h)$ also
satisfies conditions (C) of \ref{points_fixes6}).
We say that $(\gamma_1,\delta_1)$ and
$(\gamma_2,\delta_2)$ are equivalent if $\gamma_1$ and $\gamma_2$ are
$\M(\Af^p)$-conjugate and $\delta_1$ and $\delta_2$ are $\sigma$-conjugate in
$\M(L)$. In particular,
$\psi(\gamma_M,\kappa)$ is the product of a term depending only on the
image of $\kappa$ in $\Kgoth_\M(I/\Q)$
and of $k(\M')k(\H)^{-1}\varepsilon_{s_M}(\gamma_M)$.

Moreover, if $\gamma_M\in\Le_r(\Z_p)\G_r(\Q_p)$, then
$\psi(\gamma_M,\kappa)$ is equal to
\begin{flushleft}$\displaystyle{(-1)^{\dim(\A_M/\A_G)}k(\M')k(\H)
^{-1}\delta^{1/2}_{\Pa_r(\Q_p)}(\gamma_M)}$\end{flushleft}
\begin{flushright}$\displaystyle{\sum_{(\gamma,\delta)}<\alpha(\gamma_0;\gamma
_h,\delta_h),\kappa>e(\gamma)e(\delta)\overline{v}(I)^{-1}O_\gamma(f^{\infty,p}
_\M)TO_\delta(\ungras_{\Le_r(\Of_L)}\times\phi_j^{\G_r})\Phi_\M^\G(\gamma_M
^{-1},\Theta),}$\end{flushright}
where the sum is taken over the same set.
Note that, as $\alpha(\gamma_0;\gamma_h,\delta_h)\in\Kgoth_\M(I/\Q)^D$,
the number $<\alpha(\gamma_0;\gamma_h,\delta_h),\kappa>$ depends only on
the image of $\kappa$ in $\Kgoth_\M(I/\Q)$.

Let $\Sigma_L$ be the set of $(\lambda_1,\dots,\lambda_r)\in\Le_r
(\Q)=(E^\times)^r$ such that :
\begin{bulletlist}
\item for every $i\in\{1,\dots,r\}$, $|\lambda_i\overline{\lambda}
_i|_p\in p^{2a\Z}$;
\item there exists $i\in\{1,\dots,r\}$ such that $\frac{1}{2a}val_p(|
\lambda_i\overline{\lambda}_i|_p)$ is odd.

\end{bulletlist}
Remember that we defined in \ref{partie_en_p3} a subgroup $\Omega^{M}$
of the group of automorphisms of $\M$ (it was called $\Omega^{M_H}$
in \ref{partie_en_p3}); if we write
$\M=(R_{E/\Q}\Gr_m)^r\times\G_r$, then $\Omega^{M}$ is the group
generated by the involutions
\[((\lambda_1,\dots,\lambda_r),g)\fle((\lambda_1,\dots,\lambda_{i-1},
\lambda_i^{-1}c(g)^{-1},\lambda_{i+1},\dots,\lambda_r),g),\]
with $1\leq i\leq r$. The order of the group $\Omega^{M}$ is $2^r$.
\index{$\Omega^{M}$}
On the other hand, there is an action of $\Sgoth_r$ on $\M$, given by
the formula
\[(\sigma,((\lambda_1,\dots,\lambda_r),g))\fle((\lambda_{\sigma^{-1}(1)},\dots,
\lambda_{\sigma^{-1}(r)}),g).\]
The subset $\Sigma_L\G_r(\Q)$ of $\M_r(\Q)$ is of course stable by $\Sgoth_r$.
By proposition \ref{prop:identite_en_p}, if $\gamma_M\in\M(\Q)$ is such that
$\psi(\gamma_M,\kappa)\not=0$, then either there exists
$\omega\in\Omega^{M}$ (uniquely determined) such that
$\omega(\gamma_M)\in\Le_r(\Z_p)\G_r(\Q)$, or
$\gamma_M\in\Sigma_L\G_r(\Q)$. By (i) of lemma
\ref{lemme:adieu_partie_lineaire}, if $\gamma_M\in\Sigma_L\G_r(\Q)$, then,
for every $\kappa_M\in\Kgoth_\M(I/\Q)$,
\[\sum_{\kappa\mapsto\kappa_M}\sum_{\gamma'_M\in\Sgoth_r\gamma_M}\varepsilon_
{s_M}(\gamma_M')\tau(\M')k(\M')\tau(\H)^{-1}k(\H)^{-1}=0.\]
Hence
\[\sum_{\gamma_M}\sum_\kappa\tau(\M')\tau(\H)^{-1}\psi(\gamma_M,\kappa)=0,\]
where the first sum is taken over the set of semi-simple stable conjugacy
classes $\gamma_M$ in $\M(\Q)$ that are elliptic in $\M(\R)$ and such that
there is no $\omega\in\Omega^{M}$ such that $\omega(\gamma_M)\in
\Le_r(\Z_p)\G_r(\Q)$.
As $\Omega^{M}$ is a subgroup of $\Nor_\G(\M)(\Q)$, the first expression
above for $\psi(\gamma_M,\kappa)$ shows that
$\psi(\omega(\gamma_M),
\kappa)=\psi(\gamma_M,\kappa)$ for every $\omega\in\Omega^{M}$.
So we get :
\begin{flushleft}$\displaystyle{\Tr'_M=(-1)^{\dim(\A_M/\A_G)}2^r\sum_{\gamma_M}
\delta_{\Pa_r(\Q_p)}^{1/2}(\gamma_M)\sum_{(\gamma,\delta)}e(\gamma)e(\delta)
\overline{v}(I)^{-1}O_\gamma(f^{\infty,p}_\M)TO_\delta(\ungras_{\Le_r(\Of_L)}
\times\phi_j^{\G_r})
}$\end{flushleft}
\begin{flushright}$\displaystyle{\Phi_\M^\G(\gamma_M^{-1},\Theta)
\sum_{\kappa_M\in
\Kgoth_\M(I/\Q)}<\alpha(\gamma_0;\gamma_h,\delta_h),\kappa_M>
\sum_{{\kappa\in\Kgoth_\G(I/\Q)_e}\atop{\kappa\fle\kappa_M}}
\tau(\M')k(\M')\tau(\H)^{-1}k(\H)^{-1},}$\end{flushright}
where the first sum is taken over the set of semi-simple stable conjugacy
classes $\gamma_M$ in $\M(\Q)$ that are elliptic in $\M(\R)$ and the second
sum is taken over the set of equivalence classes $(\gamma,\delta)$ in
$\M(\Af^p)\times\M(L)$ such that $(\gamma_M;\gamma,\delta)$ satisfies
conditions (C) of \ref{points_fixes6}.
By (ii) of lemma \ref{lemme:adieu_partie_lineaire}, for every
$\kappa_M\in\Kgoth_\M(I/\Q)$,
\[\tau(\G)\tau(\M)^{-1}
\sum_{\kappa\mapsto\kappa_M}\tau(\M')k(\M')\tau(\H)^{-1}k(\H)^{-1}=\left\{
\begin{array}{ll}2^{-r} & \mbox{ if }r<n/2 \\
2^{-r+1} & \mbox{ if }r=n/2\end{array}\right..\]
In particular, this sum is independant of $\kappa_M$. So, by the reasoning of
\cite{K-SVLR} \S4 and by the definition of $\overline{v}(I)$ in
\ref{FT_stable_geometrique4} and lemma \ref{lemme:chi_L_S} (and the fact
that $\Le_r$ is a torus), we get
\begin{flushleft}$\displaystyle{\tau(\G)\Tr'_M=(-1)^rm_r\chi(\Le_r)\sum_{
\gamma_l\in\Le_r(\Q)}\sum_{(\gamma_0;\gamma_h,\delta_h)\in C'_{\G_r,j}}
c(\gamma_0;\gamma_h,\delta_h)\delta_{\Pa_r(\Q_p)}^{1/2}(\gamma_0)O_{\gamma_l
\gamma_h}(f^{\infty,p}_\M)}$\end{flushleft}
\begin{flushright}$\displaystyle{O_{\gamma_l}(\ungras_{\Le_r(\Z_p)})TO_{
\delta_h}(\phi_j^{\G_r})\Phi_\M^\G((\gamma_l\gamma_0)^{-1},\Theta),}$
\end{flushright}
where the integer $m_r$ is as in proposition \ref{prop:FPF_simplifiee}
(ie equal to $2$ if $r=n/2$ and to $1$ if $r<n/2$).
This (combined with proposition \ref{prop:FPF_simplifiee}) finishes the proof
of the equality of (i).

We now prove (ii). The proof is very similar to that of (i). Assume that
$r\geq q+1$, and write
\[\Tr'_r=\sum_{(\M',s_M,\eta_{M,0})\in\Ell_{\G^*}
(\M^*)}\tau(\H)^{-1}|\Lambda_{\G^*}(\M',s_M,\eta_{M,0})|^{-1}
ST_{M'}^H(f_\H^{(j)}).\]
Then, as in the prooof of (i), we see that
\[\Tr'_r=\sum_{\gamma_M}\sum_{\kappa\in\Kgoth_{\G^*}(I/\Q)_e}\tau(\M')
\tau(\H)^{-1}\psi(\gamma_M,\kappa),\]
where :
\begin{bulletlist}
\item The first sum is taken over the set of semi-simple stable conjugacy
classes $\gamma_M$ in $\M^*(\Q)$ that are elliptic in $\M^*(\R)$.
\item $I=\M^*_{\gamma_M}$, and $\Kgoth_{\G^*}(I/\Q)$ is defined above lemma
\ref{lemme:param_groupes_endoscopiques}.
\item The subset $\Kgoth_{\G^*}(I/\Q)_e$ of $\Kgoth_{\G^*}(I/\Q)$ is defined as
above. If $\kappa\in\Kgoth_{\G^*}(I/\Q)_e$ and
$(\M',s_M,\eta_{M,0},\gamma')$ is an endoscopic $\G^*$-quadruple associated to
$\kappa$ by lemma \ref{lemme:param_groupes_endoscopiques}, then
\[\psi(\gamma_M,\kappa)=SO_{\gamma'}
((f_\H^\infty)_{\M'})S\Phi_{M'}^H(\gamma',f_{\H,\infty}),\]
where $(\H,s,\eta_0)$ is the image of $(\M',s_M,\eta_{M,0})$ in $\Ell(\G^*)$.

\end{bulletlist}

If there exists a place $v\not=p,\infty$ of $\Q$ such that
$\M^*_{\Q_v}$ does not transfer to $\G_{\Q_v}$ (see lemma
\ref{lemme:transfert_terme_constant2}), then $\Tr'_r=0$ by (ii) of lemma
\ref{lemme:transfert_terme_constant2}. So we may assume that $\M^*_{\Q_v}$
transfers to $\G_{\Q_v}$ for every $v\not=p,\infty$. Then, reasoning
exactly as in the proof of (i) (and using (i) of lemma
\ref{lemme:transfert_terme_constant2}), we see that, for every $\gamma_M$ and
every $\kappa$ as above, $\psi(\gamma_M,\kappa)$ is the product of a term
depending only on the image of $\kappa$ in $\Kgoth_{\M^*}(I/\Q)$
and of $k(\M')k(\H)^{-1}\varepsilon_{s_M}(\gamma_M)$. Now, to show that
$\Tr'_r=0$, we can use (i) of lemma \ref{lemme:adieu_partie_lineaire} as in
the proof of (i).

\end{prooft}

The next two lemmas are results of \cite{K-NP}.

\begin{lemma}\label{lemme:transfert_terme_constant}(cf \cite{K-NP} 7.10
and lemma 7.6) Fix a place $v$ of $\Q$.
Let $\M$ be a Levi subgroup of $\G_{\Q_v}$, $(\M',s_M,\eta_{M,0})\in\Ell_
{\G_{\Q_v}}(\M)$ and $(\H,s,\eta_0)$ be the image of $(\M',s_M,\eta_{M,0})$
in $\Ell(\G_{\Q_v})$.
As in lemma \ref{lemme:Levi_endoscopiques}, we identify $\M'$ with a
Levi subgroup of $\H$. Choose compatible extensions $\eta:{}^L\H\fl{}^L
\G_{\Q_v}$
and $\eta_M:{}^L\M'\fl{}^L\M$ of $\eta_0$ and $\eta_{M,0}$, and
normalize transfer factors as in \ref{FT_stable_geometrique2}.

\begin{itemize}
\item[(i)] Let $f\in C_c^\infty(\G(\Q_v))$. Then, for every $\gamma\in\M(\Q_v)$
semi-simple and $\G$-regular,
\[SO_\gamma(f_\M)=|D_M^G(\gamma)|_v^{1/2}SO_\gamma(f).\]
(Remember that $D_M^G(\gamma)=\det(1-\Ad(\gamma),Lie(\G)/Lie(\M)$.)
\item[(ii)] Let $f\in C_c^\infty(\G(\Q_v))$, and let $f^\H\in C_c^\infty
(\H(\Q_v))$ be a transfer of $f$. Then, for every $\gamma_H\in\M'(\Q_v)$
semi-simple and $(\M,\M')$-regular,
\[SO_{\gamma_H}((f^\H)_{\M'})=\sum_\gamma\Delta_{\M',s_M}^\M(\gamma_H,\gamma)
e(\M_\gamma)O_\gamma(f_\M),\]
where the sum is taken over the set of semi-simple conjugacy classes
$\gamma$ in $\M(\Q_v)$ that are images of $\gamma_H$.

\end{itemize}

\end{lemma}

\begin{proof} (cf \cite{K-NP})
We show (i). Let $\gamma\in\M(\Q_v)$ be semi-simple
and $\G$-regular.
By the descent formula (\cite{A-ITF1} corollary 8.3),
\[O_{\gamma}(f_\M)=|D_M^G(\gamma)|_v^{1/2}O_\gamma(f).\]
On the other hand, as $\M_\gamma=\G_\gamma$ and as the morphism
$\H^1(\Q_v,\M)\fl\H^1(\Q_v,\G)$ is injective,
\footnote{This is explained in \cite{K-NP} A.1, and is true for any
reductive group over a field of characteristic $0$ : Choose a parabolic
subgroup $\Pa$ of $\G$ with Levi subgroup $\M$. Then the map
$\H^1(\Q_v,\M)\fl\H^1(\Q_v,\Pa)$ is bijective, and the map
$\H^1(\Q_v,\Pa)\fl\H^1(\Q_v,\G)$ is injective (the second map has a trivial
kernel by theorem 4.13(c) of \cite{BT}, and it is not hard to deduce from
this that it is injective).}
the obvious map
$\Ker(\H^1(\Q_v,\M_\gamma)\fl\H^1(\Q_v,\M))\fl\Ker(\H^1(\Q_v,
\G_\gamma)\fl\H^1(\Q_v,\G))$ is a bijection. In other words, there is a natural
bijection from the set of
conjugacy classes in the stable conjugacy class of $\gamma$ in
$\M(\Q_v)$ to the set of conjugacy classes in the stable conjugacy class of
$\gamma$ in $\G(\Q_v)$. This proves the equality of (i).

We show (ii). By lemma 2.4.A of \cite{LS2}, it is enough to show the equality
for $\gamma_H$ regular in $\H$. We may even assume that all the images of
$\gamma_H$ in $\M(\Q_v)$ are regular in $\G$; in that case, all the signs
$e(\M_\gamma)$ in the equality that we are trying to prove are trivial.
Applying (i) to $f^\H$, we get
\[SO_{\gamma_H}((f^\H)_{\M'})=|D_{M'}^H(\gamma_H)|_v^{1/2}
SO_{\gamma_H}(f^\H).\]
By definition of the transfer, this implies that
\[SO_{\gamma_H}((f^\H)_{\M'})=|D_{M'}^H(\gamma_H)|_v^{1/2}
\sum_\gamma\Delta_\H^\G(\gamma_H,\gamma)O_\gamma(f),\]
where the sum is taken over the set of conjugacy classes $\gamma$ in
$\G(\Q_v)$ that are images of $\gamma_H$. Such a conjugacy class has a
non-empty intersection with $\M(\Q_v)$, so the equality of (ii) is a
consequence of the descent formula and of the normalization of the transfer
factors.

\end{proof}

\begin{lemma}\label{lemme:transfert_terme_constant2}(cf \cite{K-NP}
lemma 7.4 and A.2) Write as before $\G^*=\GU^*(n)$.
Fix a place $v$ of $\Q$. Let $\M^*$ be a Levi subgroup of $\G^*_{\Q_v}$.
As in \cite{K-NP} A.2, we say that $\M^*$ \emph{transfers to $\G_{\Q_v}$}
if there exists an
inner twisting $\psi:\G^*\fl\G$ such that the restriction of $\psi$ to
$\A_{M^*}$ is defined over $\Q_v$.
\index{transfer to an inner form (for a Levi)}

\begin{itemize}
\item[(i)] Assume that $\M^*$ transfers to $\G_{\Q_v}$, and let
$\psi:\G^*\fl\G$ be an inner twisting such that $\psi_{|\A_{M^*}}$ is
defined over $\Q_v$. Then $\M:=\psi(\M^*)$ is a Levi subgroup of $\G_{\Q_v}$,
and $\psi_{|\M^*}:\M^*\fl\M$ is an inner twisting.
\item[(ii)] Assume that $\M^*$ does not transfer to $\G_{\Q_v}$.
Let $(\M',s_M,\eta_{M,0})\in\Ell_{\G^*_{\Q_v}}(\M^*)$, and let $(\H,s,\eta_0)$
be the image of $(\M',s_M,\eta_{M,0})$ in $\Ell(\G^*_{\Q_v})=\Ell(\G_{\Q_v})$.
As in lemma \ref{lemme:Levi_endoscopiques}, we identify $\M'$ with a
Levi subgroup of $\H$.
Let $f\in C_c^\infty(\G(\Q_v))$, and let $f^\H\in C_c^\infty
(\H(\Q_v))$ be a transfer of $f$. Then, for every $\gamma_H\in\M'(\Q_v)$
semi-simple and $(\M,\M')$-regular,
\[SO_{\gamma_H}((f^\H)_{\M'})=0.\]

\end{itemize}

\end{lemma}

\begin{proof} Point (i) follows from the fact that $\M=\Cent_{\G_{\Q_v}}
(\psi(\A_{M^*}))$.
We prove (ii). By lemma
2.4.A of \cite{LS2}, we may assume that $\gamma_H$ is regular in $\H$;
by continuity, we may even assume that $\gamma_H$ is $\G^*$-regular.
Let $\T_H$ be the centralizer of $\gamma_H$ in $\H$. It is a maximal torus of
$\M_H$ and $\H$, and it transfers to $\M^*$ and $\G^*$ because $\gamma_H$
is $\G^*$-regular.
By (i) of lemma \ref{lemme:transfert_terme_constant}
and the definition of the transfer, to show that
$SO_{\gamma_H}((f^{\H})_{\M'})=0$, it is enough to show that $\T_H$, seen as
a maximal torus in $\G^*$, does not
transfer to $\G$. Assume that $\T_H$ transfers to $\G$. Then there exists
an inner twisting $\psi:\G^*\fl\G$ such that $\psi_{|\T_H}$ is defined over
$\Q_v$; but $\A_{M^*}\subset\T_H$, so $\psi_{|\A_{M^*}}$ is defined over
$\Q_v$, and this contradicts the fact that $\M^*$ does not transfer to
$\G_{\Q_v}$.

\end{proof}

\begin{lemma}\label{lemme:adieu_partie_lineaire} We use the notations of the
proof of theorem \ref{th:adieu_partie_lineaire}. Let
$\kappa_M\in\Kgoth_M(I/\Q)$.
\begin{itemize}
\item[(i)] If $\gamma_M\in\Sigma_L\G_r(\Q)$, then
\[\sum_{{\kappa\in\Kgoth_\G(I/\Q)_e}\atop{\kappa\fle\kappa_M}}\sum_{\gamma'_M
\in\Sgoth_r\gamma_M}\varepsilon_{s_M}(\gamma'_M)\tau(\M')k(\M')
\tau(\H)^{-1}k(\H)^{-1}=0.\]
\item[(ii)]
\[\tau(\G)\tau(\M)^{-1}\sum_{{\kappa\in\Kgoth_\G(I/\Q)_e}\atop{\kappa\fle
\kappa_M}}\tau(\M')k(\M')\tau(\H)^{-1}k(\H)^{-1}=\left\{\begin{array}{ll}
2^{-r} & \mbox{ if }r<n/2 \\
2^{-r+1} & \mbox{ if }r=n/2\end{array}\right..\]

\end{itemize}

\end{lemma}

\begin{proof} Remember that $\M=\M_{\{1,\dots,r\}}\simeq (R_{E/\Q}\Gr_m)^r
\times\GU(p-r,q-r)$.
By remark \ref{rq:k_tau}, $\tau(\H)k(\H)=2^{n-1}$ for every
$(\H,s,\eta_0)\in\Ell(\G)$ and, for every $(\M',s_M,\eta_{M,0})\in\Ell_\G(\M)$,
\[k(\M')\tau(\M')=k(\M)\tau(\M)=\left\{\begin{array}{ll}2^{n-2r-1} & \mbox
{ if } r<n/2 \\ 1 & \mbox{ if }r=n/2
\end{array}\right..\]
In particular, the term $\tau(\M')k(\M')\tau(\H)^{-1}k(\H)^{-1}$ in the two
sums of the lemma does not depend on $\kappa$; it is equal to
$2^{-2r}$ if $r<n/2$, and to $2^{1-n}=2^{-2r+1}$ if $r=n/2$. Besides, by
lemma \ref{lemme:Tamagawa}, $\tau(\G)\tau(\M)^{-1}$ is equal to
$1$ if $r<n/2$, and to $2$ if $r=n/2$.

We calculate
$\Kgoth_{\kappa_M}:=\{\kappa\in\Kgoth_\G(I/\Q)_e|\kappa\fle\kappa_M\}$.
Write $\Gamma=\Gal(\overline{\Q}/\Q)$, and choose an embedding
$\widehat{\M}\subset\widehat{\G}$ as in lemma
\ref{lemme:Levi_endoscopiques_GU}. Then we get isomorphisms
$Z(\widehat{\G})\simeq\C^\times\times\C^\times$ and $Z(\widehat{\M})\simeq
\C^\times\times(\C^\times)^r\times\C^\times\times(\C^\times)^r$ such that
the embedding $Z(\widehat{\G})\subset Z(\widehat{\M})$ is
$(\lambda,\mu)\fle(\lambda,(\mu,\dots,\mu),\mu,(\mu,\dots,\mu))$ and that
the action of $\Gal(E/\Q)$ on $Z(\widehat{\G})$ and $Z(\widehat{\M})$ is
given by the following formulas
\[\tau(\lambda,\mu)=(\lambda\mu^n,\mu^{-1})\]
\[\tau(\lambda,(\lambda_1,\dots,\lambda_r),\mu,(\lambda'_r,\dots,\lambda'_1))=
(\lambda\mu^{n-2r}\lambda_1\dots\lambda_r\lambda'_1\dots\lambda'_r,
({\lambda'_1}^{-1},\dots,{\lambda'_r}^{-1}),\mu^{-1},(\lambda_1^{-1},\dots,
\lambda_r^{-1}))\]
(remember that $\tau$ is the non-trivial element of $\Gal(E/\Q)$).
This implies that $(Z(\widehat{\M})/Z(\widehat{\G}))^\Gamma\simeq
(\C^\times)^r$ is connected (this is a general fact, cf \cite{K-NP}
A.5). By the exact sequence of \cite{K-STF:CTT} 2.3, the morphism
$\H^1(\Q,Z(\widehat{\G}))\fl\H^1(\Q,Z(\widehat{\M}))$ is injective.
By lemma \ref{lemme:Tamagawa} and
\cite{K-STF:CTT} (4.2.2), $\Ker^1(\Q,Z(\widehat{\G}))=\Ker^1(\Q,Z(\widehat{\M}
))=1$; so the following commutative diagram has exact rows :
\[\xymatrix{1\ar[r] & \Kgoth_\G(I/\Q)\ar[r]\ar[d] & (Z(\widehat{I})/Z(\widehat{
\G}))^\Gamma\ar[r]\ar[d] & \H^1(\Q,Z(\widehat{\G}))\ar[d] \\
1\ar[r] & \Kgoth_\M(I/\Q)\ar[r] & (Z(\widehat{I})/Z(\widehat{\M}))^\Gamma\ar[r]
& \H^1(\Q,Z(\widehat{\M}))}\]
Let $x\in\Kgoth_\M(I/\Q)$. Then $x$, seen as an element of $(Z(\widehat{I})/
Z(\widehat{\M}))^\Gamma$, has a trivial image in $\H^1(\Q,Z(\widehat{\M}))$,
so it is in the image of $Z(\widehat{I})^\Gamma\fl(Z(\widehat{I})/
Z(\widehat{\M}))^\Gamma$. In particular, there exists $y\in
(Z(\widehat{I})/Z(\widehat{\G}))^\Gamma$ that is sent to $x$. As the
map $\H^1(\Q,Z(\widehat{\G}))\fl\H^1(\Q,Z(\widehat{\M}))$ is injective,
$y$ has a trivial image in $\H^1(\Q,Z(\widehat{\G}))$, so $y$ is in
$\Kgoth_\G(I/\Q)$. This proves that the map $\alpha:
\Kgoth_\G(I/\Q)\fl\Kgoth_\M(I/\Q)$ in the diagram above is surjective.
We want to determine its kernel. There is an obvious injection
$\Ker(\alpha)\fl(Z(\widehat{\M})/Z(\widehat{\G}))^\Gamma$. By the injectivity
of $\H^1(\Q,Z(\widehat{\G}))\fl\H^1(\Q,Z(\widehat{\M}))$, the image
of $(Z(\widehat{\M})/Z(\widehat{\G}))^\Gamma$ in $(Z(\widehat{I})/
Z(\widehat{\G}))^\Gamma$ is included in $\Kgoth_\G(I/\Q)$; this implies
that $\Ker(\alpha)=(Z(\widehat{\M})/Z(\widehat{\G}))^\Gamma$. Finally,
we get an exact sequence
\[1\fl (Z(\widehat{\M})/Z(\widehat{\G}))^\Gamma\fl
\Kgoth_\G(I/\Q)\fl\Kgoth_\M(I/\Q)\fl 1\]
(it is the exact sequence of \cite{K-NP} (7.2.1)).

If $\kappa\in\Kgoth(\kappa_M)$ and $\gamma'_M\in\Sgoth_r\gamma_M$, write
$\varepsilon_\kappa(\gamma'_M)$ instead of $\varepsilon_{s_M}(\gamma'_M)$
(this sign depends only on $\kappa$, cf remark \ref{rq:signe_en_p}). 
As $I$ is the centralizer of an elliptic element of $\M(\Q)$, it is easy
to see from lemma \ref{lemme:Levi_endoscopiques_GU} that
$\Kgoth(\kappa_M)$ is non-empty and that we can choose
$\kappa_0\in\Kgoth(\kappa_M)$ such that $\varepsilon_{\kappa_0}(\gamma'_M)=1$
for every $\gamma'_M\in\M(\Q)$. Fix such a $\kappa_0$. For every
$A\subset\{1,\dots,r\}$, let $s_A$ be the image in
$Z(\widehat{\M})/Z(\widehat{\G})$ of the element
$(1,(s_1,\dots,s_r),1,(s_r,\dots,s_1))$ of $Z(\widehat{\M})$, where
$s_i=1$ if $i\not\in A$ and $s_i=-1$ if $i\in A$. Lemma
\ref{lemme:Levi_endoscopiques_GU} implies that $\Kgoth(\kappa_M)=
\{\kappa_0+s_A,A\subset\{1,\dots,r\}\}$. If $r<n/2$, then the $s_A$ are
pairwise distinct, so $|\Kgoth(\kappa_M)|=2^r$. If $r=n/2$, then
$s_A=s_{A'}$ if and only if $\{1,\dots,r\}=A\sqcup A'$, so
$|\Kgoth(\kappa_M)|=2^{r-1}$. This finishes the proof of (ii).

We now prove (i). Let $\gamma_M\in\Sigma_L\G_r(\Q)$. We want to show that
\[\sum_{\kappa\in\Kgoth(\kappa_M)}\sum_{\gamma'_M\in\Sgoth_r\gamma_M}
\varepsilon_{s_M}(\gamma'_M)=0.\]
Write $\gamma_M=((\lambda_1,\dots,\lambda_r),g)\in (E^\times)^r\times
\G_r(\Q)$, and let $B$ the set of $i\in\{1,\dots,r\}$ such that
$\frac{1}{2a}val_p(|\lambda_i\overline{\lambda}_i|_p)$ is odd.

It is easy to see from the definition of $\varepsilon_{s_M}$ that, for
every $\sigma\in\Sgoth_r$ and $A\subset\{1,\dots,r\}$,
$\varepsilon_{\kappa_0+s_A}(\sigma(\gamma_M))=
(-1)^{|A\cap\sigma(B)|}$. (If $r=n/2$, this sign is the same for $A$ and
$\{1,\dots,r\}-A$ because $|B|$ is even by remark \ref{rq:signe_en_p}.)
So it is enough to show that, for every $\sigma\in\Sgoth_r$,
$\sum\limits_{A\subset\{1,
\dots,r\}}(-1)^{|A\cap\sigma(B)|}=0$. But
\[\sum_{A\subset\{1,\dots,r\}}(-1)^{|A\cap\sigma(B)|}=2^{n-|B|}\sum_{A\subset
\sigma(B)}(-1)^{|A|},\]
and this is equal to $0$ because $B$ is non-empty by the hypothesis on
$\gamma_M$.

\end{proof}

\begin{lemma}\label{lemme:chi_L_S} Let $s\in\{1,\dots,q\}$. Write
$S=\{1,\dots,s\}$. Then
\[\chi(\Le_S)=\vol(\A_{L_S}(\R)^0\sous\Le_S(\R))^{-1}\]
(where $\A_{L_S}$ is the maximal split subtorus of $\Le_S$, ie
$\Gr_m^s$).

\end{lemma}

\begin{proof} By \cite{GKM} 7.10 and the fact that $\Le_S/\A_{L_S}$ is
$\R$-anisotropic, we get
\[\chi(\Le_S)=(-1)^{q(\Le_S)}\tau(\Le_S)\vol(\A_{L_S}(\R)^0\sous
\Le_S(\R))^{-1}d(\Le_S).\]
As $\Le_S$ is a torus, $q(\Le_S)=0$ and $d(\Le_S)=1$. Moreover, by
lemma \ref{lemme:Tamagawa}, $\tau(\Le_S)=1$.

\end{proof}

\chapter{Applications}
\label{applications}

This chapter contains a few applications of corollary \ref{cor:stab_FT_IC}.
First we show how corollary \ref{cor:stab_FT_IC} implies a variant of
theorem \ref{th:FT_stable_geometrique} for the unitary groups of
\ref{groupes1}. The only reason we do this is to make the other applications
in this chapter logically independant of the unpublished \cite{K-NP}
(this independance is of course only formal, as the whole stabilization of the
fixed point formula in this book was inspired by \cite{K-NP}). Then we gave
applications to the calculation of the (Hecke) isotypical components
of the intersection cohomology and to the Ramanujan-Petersson conjecture.

\section{Stable trace formula}
\label{applications0}

The simplest way to apply corollary \ref{cor:stab_FT_IC} is to use
theorem \ref{th:FT_stable_geometrique} (ie the main result of \cite{K-NP})
to calculate the terms $ST^{H}(f_{\H})$. In this section, we show how
to avoid this reference to the unpublished \cite{K-NP}.

Notations are as in chapter \ref{stabilisation_K}, but with $\G$ any of
the unitary groups defined in \ref{groupes1}. As in
\ref{stabilisation_K1}, fix, for every $(\H,s,\eta_0)\in\Ell(\G)$, a
$L$-morphism $\eta:{}^L\H\fl{}^L\G$ extending $\eta_0$.

\begin{definition}\label{def:f_H} Let $f_\infty\in C^\infty(\G(\R))$. Suppose
that $f_\infty=\sum\limits_\varphi c_\varphi f_\varphi$, where the sum
is taken over the set of equivalence classes of elliptic Langlands
parameters $\varphi:W_\R\fl\widehat{\G}\rtimes W_\R$ and the
$c_\varphi$ are complex numbers that are almost all $0$
($f_\varphi$ is defined in the beginning of \ref{stabilisation_K1}).
Then, for every $(\H,s,\eta_0)\in\Ell(\G)$, set
\[(f_\infty)_\H=<\mu_\G,s>\sum_\varphi c_\varphi\sum_{\varphi_H\in\Phi_H
(\varphi)}\det(\omega_*(\varphi_H))f_{\varphi_H},\]
\index{finfH@$(f_\infty)_\H$}where
the bijection $\Phi_H(\varphi)\iso\Omega_*$, $\varphi_H\fle
\omega_*(\varphi_H)$ is as in \ref{stabilisation_K1} determined by the
standard Borel subgroup of $\G_\C$.

\end{definition}

\begin{remark} By the trace Paley-Wiener theorem of Clozel and Delorme
(\cite{CD}, cf the beginning of section 3 of \cite{A-L2}),
if $f_\infty\in C^\infty(\G(\R))$ is stable cuspidal, then $f_\infty$ satifies
the condition of definition \ref{def:f_H}, so $(f_\infty)_\H$ is defined.

\end{remark}

\begin{definition}\label{def:iota_G_H} Let $a_1,b_1,\dots,a_r,b_r\in\Nat$
be such that $a_i\geq b_i$ for all $i$ and $\G=\G(\U(a_1,b_1)\times\dots
\times\U(a_r,b_r))$. Write $n_i=a_i+b_i$; then
the quasi-split inner form of $\G$ is $\G(\U^*(n_1)\times\dots\times
\U^*(n_r))$. Fix $n_1^+,n_1^-,\dots,n_r^+,n_r^-\in\Nat$ such that
$n_i=n_i^++n_i^-$ for every $i\in\{1,\dots,r\}$ and that $n_1^-+\dots+n_r^-$
is even. Let $(\H,s,\eta_0)$ be the elliptic endoscopic triple for $\G$
associated to these integers as in proposition
\ref{prop:groupes_endoscopiques}. For every $i\in\{1,\dots,r\}$,
if $I_i\subset\{1,\dots,n_i\}$, set
$n_i(I_i)=|I\cap\{n_i^++1,\dots,n_i\}|$. We define a rational number $\iota_
{\G,\H}$ by
\[\iota_{\G,\H}=\iota(\G,\H)|\pi_0(\X)|^{-1}
\sum_{{I_1\subset\{1,\dots,n_1\}}\atop{|I_1|=a_1}}
\dots\sum_{{I_r\subset\{1,\dots,n_r\}}\atop{|I_r|=a_r}}(-1)^{
n_1(I_1)+\dots+n_r(I_r)},\]
where $\X$ is the symmetric space appearing in the Shimura data of
section \ref{groupes1} for $\G$.
\index{$\iota_{\G,\H}$}

\end{definition}

\begin{proposition}\label{prop:FT_stable_geometrique}
Let $f=f^\infty f_\infty$ be as in theorem
\ref{th:FT_stable_geometrique}.  Assume that $f_\infty$ is stable
cuspidal and that, for every $(\H,s,\eta_0)\in\Ell(\G)$, there exists
a transfer $(f^\infty)^\H$ of $f^\infty$. Then :
\[T^G(f)=\sum_{(\H,s,\eta_0)\in\Ell}\iota_{\G,\H}
ST^H((f^\infty)^\H(f_\infty)_\H).\]

\end{proposition}

\begin{remark} It is not very hard to see that proposition
\ref{prop:FT_stable_geometrique} is a consequence of theorem
\ref{th:FT_stable_geometrique}. The goal here is to prove it directly.

\end{remark}

To prove this proposition, we first need an extension of corollary
\ref{cor:stab_FT_IC} (proposition \ref{prop:pour_tout_m} below).

Fix a prime number $p$ where $\G$ is unramified.
Remember that we defined, for every $m\in\Nat^*$, a function $\phi^\G_m$
on $\G(L)$, where $L$ is an unramified extension of $\Q_p$; if
$(\H,s,\eta_0)\in\Ell(\G)$, write $f_{\H,p}^{(m)}$ for the function in
$C_c^\infty(\H(\Q_p))$ obtained by twisted transfer from $\phi_m^\G$
as in \ref{partie_en_p3}. In the proof of proposition
\ref{prop:identite_en_p}, we have calculated the Satake transform $S$ of
$f_{\H,p}^{(m)}$, or more precisely of
$\chi_{\eta_p}^{-1}f_{\H,p}^{(m)}$, where $\chi_{\eta_p}$ is the
quasi-character of $\H(\Q_p)$ associated to $\eta_{|\widehat{\H}\rtimes
W_{\Q_p}}:\widehat{\H}\rtimes
W_{\Q_p}\fl\widehat{\G}\rtimes W_{\Q_p}$ as in the last two subsections of
\ref{partie_en_p2}.
Notice that the expression for the Satake transform $S$
makes sense for any $m\in\Z$.

\begin{definition}\label{def:phi_m_G} If $m\in\Z$, we define
$f_{\H,p}^{(m)}\in C_c^\infty(\H(\Q_p))$ in the following way :
$\chi_{\eta_p}^{-1}f_{\H,p}^{(m)}\in\Hecke(\H(\Q_p),\H(\Z_p))$, and its
Satake transform is given by the polynomial $S$ in the proof of
proposition \ref{prop:identite_en_p} (where of course the integer $a$ in
the definition of $S$ is replaced by $m$).
\index{fHpm@$f_{\H,p}^{(m)}$}

\end{definition}

Fix $f^{\infty,p}\in C_c^\infty(\G(\Af^p))$ and an irreducible algebraic
representation $V$ of $\G_\C$.
For every $(\H,s,\eta_0)\in\Ell(\G)$ and $m\in\Z$, let $f_\H^{(m)}=
f_\H^{p,\infty}f_{\H,p}^{(m)}f_{\H,\infty}\in C^\infty(\H(\Ade))$, where
$f_\H^{p,\infty}$ and $f_{\H,\infty}$ are as in \ref{stabilisation_K1}. 

\begin{proposition}\label{prop:pour_tout_m} Assume that $p$ is inert in $E$.
\footnote{This is not really necessary but makes the proof slightly simpler.}
Then, with the notations of \ref{stabilisation_K2}, for every $m\in\Z$,
\[\Tr(\Phi_\wp^mf^\infty,W_\lambda)=
\sum_{(\H,s,\eta_0)\in\Ell(\G)}\iota(\G,\H)ST^{H}(f_{\H}^{(m)}),\]
where $f^\infty=f^{\infty,p}\ungras_{\G(\Z_p)}$.

\end{proposition}

Notice that, for $m>>0$, this is simply corollary \ref{cor:stab_FT_IC}
(cf the remark following this corollary).

\begin{proof} Fix an (arbitrary) embedding $\iota:K_\lambda\subset\C$ and write
$W=\iota_*(W_\lambda)$. Then $W$ is a virtual complex representation of
$\Hecke_\K\times\Gal(\overline{\Q}/F)$.
As the actions of $\Phi_\wp$ and $f^\infty$ on
$W$ commute, there exist a finite set $I_0$ and families of complex numbers
$(c_i)_{i\in I_0}$ and $(\alpha_i)_{i\in I_0}$
such that, for every $m\in\Z$,
\[\Tr(\Phi_\wp^mf^\infty,W)=
\sum_{i\in I_0}c_i\alpha_i^m.\]

We now want to find a similar expression for the right hand side of
the equality of the proposition. Remember
from the definitions in \ref{FT_stable_geometrique4} that
\[ST^H(f_\H^{(m)})=\sum_{\M_H}(n_{M_H}^H)^{-1}\tau(\M_H)\sum_{\gamma_H}
SO_{\gamma_H}((f_\H^{\infty,p})_{\M_H})SO_{\gamma_H}((f_{\H,p}^{(m)})_{\M_H})
S\Phi_{\M_H}^\H(\gamma_H,f_{H,\infty}),\]
where the first sum is taken over the set of conjugacy classes of cuspidal
Levi subgroups $\M_H$ of $\H$ and the second sum over the set of semi-simple
stable conjugacy classes $\gamma_H\in\M_H(\Q)$ that are elliptic in
$\M_H(\R)$. Note that the first sum is finite. In the second sum, all but
finitely many terms are zero, but the set of $\gamma_H$ such that the term
associated to $\gamma_H$ is non-zero may depend on $m$.

Fix $(\H,s,\eta_0)\in\Ell(\G)$ and a cuspidal Levi subgroup $\M_H$ of $\H$.
By the Howe conjecture, proved by Clozel in \cite{Cl-HC}, the space of
linear forms $\Hecke(\M_H(\Q_p),\M_H(\Z_p))\fl\C$ generated by the
elements $h\fle SO_{\gamma_H}(h)$, for $\gamma_H\in\M_H(\Q_p)$ semi-simple
elliptic, is finite-dimensional. As $p$ is inert in $E$, any semi-simple
$\gamma_H\in\M_H(\Q)$ that is elliptic in $\M_H(\R)$ is also elliptic in
$\M_H(\Q_p)$.
\footnote{Such a $\gamma_H$ is elliptic in $\M_H(\R)$ (resp. $\M_H(\Q_p)$)
if and only if it is contained in no Levi subgroup of $\M_H(\R)$ (resp.
$\M_H(\Q_p)$). But the Levi subgroups of $\M_H(\R)$ and $\M_H(\Q_p)$ are all
defined over $\Q$.}
So we find that the space $D$ of linear forms on $\Hecke(\M_H(\Q_p),
\M_H(\Z_p))$ generated by the $h\fle SO_{\gamma_H}(h)$, for $\gamma_H\in
\M_H(\Q)$ semi-simple and elliptic in $\M_H(\R)$, is finite-dimensional.

On the other hand, by Kazhdan's density theorem (\cite{Ka} theorem 0), every
distribution $h\fle SO_{\gamma_H}(h)$ on $\Hecke(\M_H(\Q_p),\M_H(\Z_p))$ is
a finite linear combination of distributions of the type $h\fle\Tr(\pi(h))$,
for $\pi$ a smooth irreducible representation of $\M_H(\Z_p)$ (that we may
assume to be unramified). So the space $D$ is generated by a finite number
of distributions of that type. 

Using the form of the Satake transform of $\chi_{\eta_p}^{-1}f_{\H,p}^{(m)}$,
it is easy to see that
this implies that there exist a finite set $I_{M_H,H}$ and a family of
complex numbers $(\beta_{M_H,H,i})_{i\in I_{M_H,H}}$ such that, for
every $\gamma_H\in\M_H(\Q)$ semi-simple and elliptic in $\M_H(\R)$,
there exists a family of complex numbers $(d_{M_H,H,i}(\gamma_H))_{i\in
I_{M_H,H}}$ with
\[SO_{\gamma_H}((f_\H^{\infty,p})_{\M_H})SO_{\gamma_H}((f_{\H,p}^{(m)})_{\M_H})
S\Phi_{\M_H}^\H(\gamma_H,f_{H,\infty})=\sum_{i\in I_{M_H,H}}d_{M_H,H,i}
(\gamma_H)\beta_{M_H,H,i}^m\]
for every $m\in\Z$.

Let $m_0\in\Z$. We want to prove the equality of the proposition for
$m=m_0$. Let $N\in\Nat$ such that the equality of the proposition is true
for $m\geq N$ (such a $N$ exists by corollary \ref{cor:stab_FT_IC}). We may
assume that $m_0\leq N$. Let
\[M=|I_0|+\sum_{(\H,s,\eta_0)\in\Ell(\G)}\sum_{\M_H}|I_{M_H,H}|,\]
where the second sum is taken over a set of conjugacy classes of cuspidal
Levi subgroups of $\H$. 
For every $\H$ and $\M_H$ as before, let $\Gamma_{M_H,H}$ be the set of
semi-simple stable conjugacy classes $\gamma_H\in\M_H(\Q)$ that are elliptic
in $\M_H(\R)$ and such that there exists $m\in\Z$ with $m_0\leq m\leq N+M$ and
\[SO_{\gamma_H}((f_\H^{\infty,p})_{\M_H})SO_{\gamma_H}((f_{\H,p}^{(m)})_{\M_H})
S\Phi_{\M_H}^\H(\gamma_H,f_{H,\infty})\not=0.\]
This set is finite. So, by the above calculations, there exist families
of complex numbers $(d_{M_H,H,i})_{i\in I_{M_H,H}}$, for all $\M_H$ and $\H$
as before, such that, for every $m\in\Z$ with $m_0\leq m\leq N+M-1$,
\[\sum_{(\H,s,\eta_0)\in\Ell(\G)}\iota(\G,\H)ST^{H}(f_{\H}^{(m)})=
\sum_{(\H,s,\eta_0)\in\Ell(\G)}\sum_{\M_H}\sum_{i\in I_{M_H,H}}d_{M_H,H,i}
\beta_{M_H,H,i}^m.\]
All the sums above are finite. So we can reformulate this as : there exist
a finite set $J$ (with $|J|=\sum\limits_H\sum\limits_{M_H}|I_{M_H,H}|$) and
families of complex numbers $(d_j)_{j\in J}$ and $(\beta_j)_{j\in J}$ such
that, if $m_0\leq m\leq N+M-1$, then
\[\sum_{(\H,s,\eta_0)\in\Ell(\G)}\iota(\G,\H)ST^{H}(f_{\H}^{(m)})=\sum_{j\in J}
d_j\beta_j^m.\]

So the result that we want to prove is that the equality
\[\sum_{i\in I_0}c_i\alpha_i^m=\sum_{j\in J}d_j\beta_j^m\]
holds for $m=m_0$. But we know that this equality holds if $N\leq m\leq 
N+M-1$, and $M=|I_0|+|J|$, so this equality holds for all $m\in\Z$.

\end{proof}

\begin{proofpn}{\ref{prop:FT_stable_geometrique}}
We may assume that $f^\infty$ is a product $\bigotimes\limits_p f_p$, with
$f_p=\ungras_{\G(\Z_p)}$ for almost all $p$. Let $\K=\prod\limits_p\K_p$ be
a neat open compact subgroup of $\G(\Af)$ such that $f\in\Hecke(\G(\Af),\K)$.
Fix a prime number $p$ that is inert in $E$ and such that $\G$ is unramified
at $p$, $f_p=\ungras_{\G(\Z_p)}$ and $\K_p=\G(\Z_p)$.
Define a virtual representation $W$
of $\Hecke(\G(\Af),\K)\times\Gal(\overline{\Q}/F)$ as in the proof of
proposition \ref{prop:pour_tout_m}. Then, by formula (3.5) and theorem
6.1 of \cite{A-L2}, and by theorem 7.14.B and paragraph (7.19) of
\cite{GKM} :
\footnote{In the articles \cite{A-L2} and \cite{GKM}, the authors consider
only connected symmetric spaces, i.e., they use $\pi_0(\G(\R))\sous\X$
instead of $\X$ (in the cases considered here, $\G(\R)$ acts transitively on
$\X$, so $\pi_0(\G(\R))\sous\X$ is connected).
When we pass from $\pi_0(\G(\R))\sous\X$ to $\X$, the trace of
Hecke operators is multiplied by $|\pi_0(\X)|$.}
\[\Tr(f^\infty,W)=|\pi_0(\X)|T^G(f).\]
On the other hand, using proposition \ref{prop:pour_tout_m} at the place $p$
and for $m=0$, we find
\[\Tr(f^\infty,W)=\sum_{(\H,s,\eta_0)\in\Ell(\G)}\iota(\G,\H)ST^{H}
(f_{\H}^{(0)}).\]
But it is obvious from the definitions of $f_{\H,p}^{(m)}$ and $\iota_{\G,\H}$
that
\[f_{\H,p}^{(0)}=\frac{\iota_{\G,\H}}{\iota(\G,\H)}|\pi_0(\X)|\chi_{\eta_p}
\ungras_{\H(\Z_p)},\]
and we know that $\chi_{\eta_p}\ungras_{\H(\Z_p)}$ is a transfer of
$f_p=\ungras_{\G(\Z_p)}$
by the fundamental lemma (cf \ref{FT_stable_geometrique3}). This finishes
the proof.

\end{proofpn}

\section{Isotypical components of the intersection cohomology}
\label{applications1}

Notations are still as in \ref{applications0}, and
we assume that $\G=\GU(p,q)$, with $n=p+q$ (for the other unitary groups
of \ref{groupes1}, everything would work the same way , but with more
complicated notations). In particular, $V$ is an irreducible algebraic
representations of $\G$ defined over a number field $K$, $\lambda$ is a
place of $K$ over $\ell$ and $\varphi:W_\R\fl\widehat{\G}\rtimes W_\R$ is an
elliptic Langlands parameter corresponding to the contragredient
$V^*$ of $V$ (as in proposition \ref{prop:calcul_Phi_M}).

Let $\Hecke_K=\Hecke(\G(\Af),\K)$.
\index{HK@$\Hecke_\K$\quad Hecke algebra at level $\K$}
Define, as in \ref{stabilisation_K2}, an object $W_\lambda$ of the
Grothendieck group of representations of $\Hecke_K\times\Gal(\overline{\Q}/F)$
in a finite dimensional $K_\lambda$-vector space by
\[W_\lambda=\sum_{i\geq 0}(-1)^i[\Ho^i(M^{\K}(\G,\X)^*_{\overline{\Q}},IC^{\K}
V_{\overline{\Q}})].\]
\index{Wlambda@$W_\lambda$}
Let $\iota:K_\lambda\fl\C$ be an embedding. Then there is an isotypical
decomposition of $\iota_*(W_\lambda)$ as a $\Hecke_\K$-module :
\[\iota_*(W_\lambda)=\sum_{\pi_f}\iota_*(W_\lambda)(\pi_f)\otimes\pi_f^{\K},\]
where the sum is taken over the set of isomorphism classes of irreducible
admissible representations $\pi_f$ of $\G(\Af)$ such that $\pi_f^{\K}\not=0$
and where the $\iota_*(W_\lambda)(\pi_f)$ are virtual representations of
$\Gal(\overline{\Q}/F)$ in finite dimensional $\C$-vector spaces.
As there is only a finite number of $\pi_f$ such that
$\iota_*(W_\lambda)(\pi_f)\not=0$, we may assume, after replacing
$K_\lambda$ by a finite extension, that there exist virtual representations
$W_\lambda(\pi_f)$ of $\Gal(\overline{\Q}/F)$ in finite dimensional
$K_\lambda$-vector spaces such that $\iota_*(W_\lambda(\pi_f))=
\iota_*(W_\lambda)(\pi_f)$. So we get
\[W_\lambda=\sum_{\pi_f}W_\lambda(\pi_f)\otimes\pi_f^{\K}.\]
\index{Wlambdapif@$W_\lambda(\pi_f)$\quad isotypical component of $W_\lambda$}

\begin{notation}
Let $\H$ be a connected reductive group over $\Q$ and $\xi$ be a
quasi-character of $\A_H(\R)^0$. We write $\Pi(\H(\Ade),\xi)$ for the set
\index{$\Pi(\G(\Ade),\xi)$}
of isomorphism classes of irreducible admissible representations of
$\H(\Ade)$ on which $\A_H(\R)^0$ acts by $\xi$. For every $\pi\in\Pi(\H(\Ade),
\xi)$, let $m_{disc}(\pi)$ be the multiplicity of $\pi$ in the discrete part
of $L^2(\H(\Q)\sous\H(\Ade),\xi)$ (cf \cite{A-L2}, \S2).
\index{mdisc@$m_{disc}(\pi)$}

\end{notation}

Let $\xi_G$ be the quasi-character by which the group $\A_G(\R)^0$ acts on
the contragredient of $V$.

For every $(\H,s,\eta_0)\in\Ell(\G)$, fix a $L$-morphism
$\eta:{}^L\H\fl{}^L\G$ extending $\eta_0$ as in proposition
\ref{prop:prolongement_eta_0}.
Let $\Ell^0(\G)$ be the set of
$(\H,s,\eta_0)\in\Ell(\G)$ such that $\H$ is not an inner form of $\G$.
\index{E0G@$\Ell^0(\G)$}
If $n_1,\dots,n_r\in\Nat^*$ and $\H=\G(\U^*(n_1)\times\dots\times\U^*(n_r))$,
we define in the same way a subset $\Ell^0(\H)$ of $\Ell(\H)$ and
fix, for every $(\H',s,\eta_0)\in\Ell(\H)$, a $L$-morphism
$\eta:{}^L\H'\fl{}^L\H$ extending $\eta_0$ as in proposition
\ref{prop:prolongement_eta_0}.

Let $\F_\G$ be the set of sequence $(e_1,\dots,e_r)$ of variable length
$r\in\Nat^*$, with $e_1=(\H_1,s_1,\eta_1)\in\Ell^0_\G$ and, for every
$i\in\{2,\dots,r\}$, $e_i=(\H_i,s_i,\eta_i)\in\Ell^0_{\H_{i-1}}$.
\index{FG@$\F_\G$}

Let $\es=(e_1,\dots,e_r)\in\F_\G$. Write $e_i=(\H_i,s_i,\eta_i)$ and
$\H_0=\G$.
Set $\ell(\es)=r$, $\H_\es=\H_r$, $\eta_\es=\eta_1\circ\dots\circ\eta_r:
{}^L\H_\es\fl{}^L\G$,
\[\iota(\es)=\iota_{\G,\H_1}\iota_{\H_1,\H_2}\dots\iota_{\H_{r-1},\H_r}\]
and
\[\iota'(\es)=\iota(\G,\H_1)\iota_{\H_1,\H_2}\dots\iota_{\H_{r-1},\H_r}.\]

For every finite set $S$ of prime numbers, write
$\Ade_S=\prod\limits_{p\in S}\Q_p$ and $\Ade_f^S=\prod\limits_{p\not\in S}
{}\!\!'\,\Q_p$; if $\pi'_f=\bigotimes{}'\pi'_p$ is an irreducible admissible
representations of $\G(\Af)$, write $\pi'_S=\bigotimes\limits_{p\in S}\pi'_p$
and ${\pi'}^S=\bigotimes\limits_{p\not\in S}{}\!\!'\,\pi'_p$;
if $\G$ is unramified at every $p\not\in S$, write
$\K^S=\prod\limits_{p\not\in S}\G(\Z_p)$. 
\index{AS@$\Ade_S$}
\index{AfS@$\Ade_f^S$}
\index{$\pi_S$}
\index{$\pi^S$}
If $f^{S}\in C_c^\infty(\G(\Af^S))$ and $f_S\in C_c^\infty(\G(\Ade_S))$,
define functions $(f^S)^{\es}\in C_c^\infty(\H_\es(\Af^S))$ and
$(f_S)^\es\in C_c^\infty(\H_\es(\Ade_S))$
by
\[(f^S)^{\es}=(\dots((f^{S})^{\H_1})^{\H_2}\dots)^{\H_r}\]
\[(f_S)^{\es}=(\dots((f_{S})^{\H_1})^{\H_2}\dots)^{\H_r}.\]

Define a function $f^\es_\infty$ on $\H_\es(\R)$ by
\[f^\es_\infty=(\dots((f_\infty)_{\H_1})_{\H_2}\dots)_{\H_r},\]
where $f_\infty=(-1)^{q(\G)}f_\varphi$ (this function is defined in
\ref{stabilisation_K1}). The function
$f^\es_\infty$ is stable cuspidal by definition.

Let $k\in\{1,\dots,r\}$. Consider the morphism
\[\varphi_k:W_\R\stackrel{j}{\fl}{}^L\H_{k,\R}\stackrel{\eta_\infty}{\fl}{}^L
\H_{k-1,\R}\stackrel{p}{\fl}{}^L(\A_{H_{k-1}})_\R,\]
where $j$ is the obvious inclusion, $\eta_\infty$ is induced by $\eta_k$
and $p$ is the dual of the inclusion $\A_{H_{k-1}}\fl\H_{k-1}$. The morphism
$\varphi_k$ is the Langlands parameter of a quasi-character on
$\A_{H_{k-1}}(\R)$, and we write $\chi_k$ for the restriction of this
quasi-character to $\A_{H_{k-1}}(\R)^0$. As $\A_{H_r}=\dots=\A_{H_1}=\A_G$
(because $(\H_k,s_k,\eta_{k,0})$ is an elliptic endoscopic datum for
$\H_{k-1}$ for every $k\in\{1,\dots,r\}$), we may define a quasi-character
$\xi_\es$ on $\A_{H_\es}(\R)^0$ by the formula
\[\xi_\es=\xi_\G\chi_1^{-1}\dots\chi_r^{-1}.\]
This quasi-character satisfies the following property : if $\varphi_{H_\es}:
W_\R\fl{}^L\H_{\es,\R}$ is a Langlands parameter corresponding to a
$L$-packet of representations of $\H_\es(\R)$ with central character
$\xi_\es$ on $\A_{H_\es}(\R)^0$, then $\eta_\es\circ\varphi_{H_\es}:
W_\R\fl{}^L\G_\R$ corresponds to a $L$-packet of representations of $\G(\R)$
with central character $\xi_\G$ on $\A_G(\R)^0$. (This is the construction of
\cite{K-NP} 5.5).
Write $\Pi_\es=\Pi(\H_\es(\Ade),\xi_\es)$. Let $R_\es(V)$ be the set of
$\pi_\infty\in\Pi(\H_\es(\R))$ such that there exists an elliptic Langlands
parameter $\varphi_\es:W_\R\fl{}^L\H_{\es,\R}$ satisfying the following
properties : $\eta_\es\circ\varphi_\es$ is $\widehat{\G}$-conjugate to 
$\varphi$, and $\Tr(\pi_\infty(f_{\varphi_\es}))\not=0$ (remember that
$f_{\varphi_\es}$ is defined in \ref{stabilisation_K1}).
Then $R_\es(V)$ is finite.

If $p$ is a prime number unramified in $E$,
let $\eta_{\es,p}=\eta_{es|\widehat{\H}_\es\rtimes W_{\Q_p}}$
and write $\eta_{\es,p,simple}:{}^L\H_\es\fl
{}^L\G$ for the $L$-morphism extending $\eta_{1,0}\circ\dots\eta_{r,0}$
and equal to the composition of the analogs of the morphism
$\eta_{simple}$ of the last two subsections of \ref{partie_en_p2}. Write
$\eta_{\es,p}=c\eta_{\es,p,simple}$, where $c:W_{\Q_p}\fl Z(\widehat{\H}_\es)$
is a $1$-cocycle. Let $\chi_{\es,p}=\chi_{\eta_{\es,p}}$
be the quasi-character of $\H_\es(\Q_p)$ corresponding to the class of $c$
in $\Ho^1(W_{\Q_p},Z(\widehat{\H}_\es))$.

Suppose that $(\H_1,s_1,\eta_{1,0})$ is the elliptic endoscopic triple for
$\G$ defined by a pair $(n^+,n^-)\in\Nat^2$ as in proposition 
\ref{prop:groupes_endoscopiques} (so $n=n^++n^-$ and $n^-$ is even).
Write
\[\H_\es=\G(\U^*(n_1^+)\times\dots\times\U^*(n_r^+)\times\U^*(n_1^-)\times\dots
\times\U^*(n_s^-)),\]
where the identification is chosen such that $\eta_2\circ\dots\circ\eta_r$
sends
$\widehat{\U^*(n_1^+)}\times\dots\times\widehat{\U^*(n_r^+)}$ (resp.
$\widehat{\U^*(n_1^-)}\times\dots\times\widehat{\U^*(n_s^-)}$) in
$\widehat{\U^*(n^+)}$ (resp. $\widehat{\U^*(n^-)}$).

If $p_1^+,\dots,p_r^+,p_1^-,\dots,p_s^-\in\Nat$ are such that $1\leq
p_i^+\leq n_i^+$ and $1\leq p_j^-\leq n_j^-$ for every $i\in\{1,\dots,r\}$
and $j\in\{1,\dots,s\}$, write
\[\mu=\mu_{p_1^+,\dots,p_r^+,p_1^-,\dots,p_s^-}=(\mu_{p_1^+},\dots,\mu_{p_r^+},
\mu_{p_1^-},\dots,\mu_{p_s^-}):\Gr_{m,E}\fl\H_{\es,E}\]
(cf \ref{notation:mu_p} for the definition of $\mu_p$), and
\[s(\mu)=(-1)^{p_1^-+\dots+p_s^-}.\]
Let $M_\es$ be the set of cocharacters $\mu_{p_1^+,\dots,p_r^+,p_1^-,\dots,
p_s^-}$ with $p=p_1^++\dots+p_r^++p_1^-+\dots+p_s^-$. For every $\mu\in
M_\es$ and every finite place $\wp$ of $F$ where $\H_\es$ is unramified,
we get a representation $r_{-\mu}$ of ${}^L\H_{\es,F_\wp}$, defined in
\ref{def_r_mu}.

For every irreducible admissible representations $\pi_{\es,f}$ of
$\H_\es(\Af)$, let
\[c_\es(\pi_\es)=\sum_{{\pi_{\es,\infty}\in\Pi(\H_\es(\R)),}\atop{\pi_{\es,f}
\otimes\pi_{\es,\infty}\in\Pi_\es}}
m_{disc}(\pi_{\es,f}\otimes\pi_{\es,\infty})\Tr(\pi_{\es,\infty}
(f^\es_\infty))\]
(as $\Tr(\pi_{\es,\infty}(f^\es_\infty))=0$ unless $\pi_{\es,\infty}\in
R_\es(V)$, this sum has only a finite number of non-zero terms).

Write $\Pi_\G=\Pi(\G(\Ade),\xi_G)$. For every irreducible admissible
representation $\pi_f$ of $\G(\Af)$, let
\[c_\G(\pi_f)=\sum_{{\pi_\infty\in\Pi(\G(\R)),}\atop{\pi_\infty\otimes\pi_f\in
\Pi_\G}}m_{disc}(\pi_f\otimes
\pi_\infty)\Tr(\pi_\infty(f_\infty))\]
(this sum has only a finite number of non-zero terms because there are only
finitely many $\pi_\infty$ in $\Pi(\G(\R))$
such that $\Tr(\pi_\infty(f_\infty))\not=0$).
Remember that there is a cocharacter $\mu_G:\Gr_{m,E}\fl\G_E$ associated to
the Shimura datum (cf \ref{groupes1}); this cocharacter gives a
representation $r_{-\mu_G}$ of ${}^L\G_{F_\wp}$, for every finite place
$\wp$ of $F$ where $\G$ is unramified.

Let $\pi_f=\bigotimes\limits_p{}\!'\pi_p$ be an irreducible admissible
representation of $\G(\Af)$ such that $\pi_f^\K\not=0$, and let $\es\in\F_\G$.
Write $R_\es(\pi_f)$ for the set of equivalence classes of irreducible
admissible representations $\pi_{\es,f}=\bigotimes\limits_p{}\!'
\pi_{\es,p}$ of $\H_\es(\Af)$ such that, for almost every prime number $p$
where $\pi_f$ and $\pi_{\es,f}$ are unramified, the morphism
$\eta_\es:{}^L\H_\es\fl{}^L\G$ sends a Langlands parameter of $\pi_{\es,p}$
to a Langlands parameter of $\pi_p$.

Let $p$ be a prime number. Remember that we fixed embeddings
$F\subset\overline{\Q}\subset\overline{\Q}_p$, that determine a place $\wp$
of $F$ above $p$ and a morphism $\Gal(\overline{\Q}_p/F_\wp)\fl
\Gal(\overline{\Q}/F)$. Let $\Phi_\wp\in\Gal(\overline{\Q}_p/F_\wp)$ be
a lift of the geometric Frobenius, and use the same notation for its
image in $\Gal(\overline{\Q}/F)$.
\index{$\Phi_\wp$\quad geometric Frobenius at $\wp$}
If $\H$ is a reductive unramified group
over $\Q_p$ and $\pi_p$ is an unramified representation of $\H(\Q_p)$, denote
by $\varphi_{\pi_p}:W_{\Q_p}\fl{}^L\H_{\Q_p}$ a Langlands parameter of $\pi_p$.

\begin{theorem}\label{th:calcul_L_moche} Let $\pi_f$ be an irreducible
admissible representation of $\G(\Af)$ such that $\pi_f^\K\not=0$. Then
there exists a function $f^\infty\in C_c^\infty(\G(\Af))$ such that, for
almost every prime number $p$ and for every $m\in\Z$,
\begin{flushleft}$\displaystyle{\Tr(\Phi_\wp^m,W_\lambda(\pi_f))=(N\wp)^{md/2}
c_\G(\pi_f)\dim(\pi_f^\K)\Tr(r_{-\mu_G}\circ\varphi_
{\pi_p}(\Phi_\wp^m))}$\end{flushleft}
\begin{center}$\displaystyle{
+(N\wp)^{md/2}\sum_{\es\in\F_\G}(-1)
^{\ell(\es)}\iota(\es)\sum_{\pi_{\es,f}\in R_\es(\pi_f)}c_\es(\pi_{\es,f})
\Tr(\pi_\es((f^\infty)^\es))
}$\end{center}
\begin{flushright}$\displaystyle{
\sum_{\mu\in M_\es}(1-(-1)^{s(\mu)}\frac{\iota'(\es)}{\iota(\es)})
\Tr(r_{-\mu}\circ\varphi_{\pi_{\es,p}\otimes\chi_{\es,p}}(\Phi_\wp^m)),
}$\end{flushright}
where the second sum in the right hand side is taken only over those
$\pi_{\es,f}$ such that $\pi_{\es,p}\otimes\chi_{\es,p}$ is unramified,
$d=\dim(M^{\K}(\G,\X))$ and $N\wp=\#(\Of_{F_\wp}/\wp)$.

\end{theorem}

\begin{remark} The lack of control over the set of ``good'' prime numbers
in the theorem above comes from the fact that we do not have a strong
multiplicity one theorem for $\G$ (and not from a lack of information about
the integral models of Shimura varieties). If $\pi_f$ extends to an
automorphic representation of $\G(\Ade)$ whose base change to $\G(\Ade_E)$
is cuspidal (cf section \ref{GL_n_applications5}), then it is possible
to do better by using corollary \ref{cor:CB}.

\end{remark}

\begin{proof}
It is enough to prove the equality of the theorem for
$m$ big enough (where the meaning of ``big enough'' can depend on
$p$).

Let $R'$ be the set of isomorphism classes of irreducible admissible
representations $\pi'_f$ of $\G(\Af)$ satisfying the following
properties :
\begin{bulletlist}
\item $\pi'_f\not\simeq\pi_f$;
\item $(\pi'_f)^\K\not=0$;
\item $W_\lambda(\pi'_f)\not=0$ or $c_\G(\pi'_f)\not=0$.

\end{bulletlist}
Then $R'$ is finite, so there exists $h\in\Hecke_\K=\Hecke(\G(\Af),\K)$
such that $\Tr(\pi_f(h))=\Tr(\pi_f(\ungras_\K))$
and $\Tr(\pi'_f(h))=0$ for every $\pi'_f\in R'$.

Let $T$ be a finite set of prime numbers such that all the representations
in $R'$ are unramified outside of $T$, that $\G$ is unramified at every
$p\not\in T$, that $\K=\K_T\K^T$ with $\K_T\subset\G(\Ade_T)$ and that
$h=h_T\ungras_{\K^T}$ with $h_T\in\Hecke(\G(\Ade_T),\K_T)$. Then, for every
function $g^T$ in $\Hecke(\G(\Ade_f^T),\K^T)$,
$\Tr(\pi_f(h_Tg^T))=\Tr(\pi_T(\ungras_{\K_T}))\Tr(\pi^T(g^T))$ and
$\Tr(\pi'_f(h_Tg^T))=0$ if $\pi'_f\in R'$.

For every $\es\in\F_\G$, let $R'_\es$ be the set of isomorphism classes of
irreducible admissible representations $\rho_f$ of $\H_\es(\Af)$ such that
$\rho_f\not\in R_\es(\pi_f)$, $\Tr(\rho_f(h^\es))\not=0$ and
$c_\es(\rho_f)\not=0$.
As $\F_\G$ is finite anf $R'_\es$ is finite for every $\es\in\F_\G$, there
exists $g^T\in\Hecke(\G(\Ade_f^T),\K^T)$ such that :
\begin{bulletlist}
\item $\Tr(\pi^T(g^T))=1$;
\item for every $\es\in\F_\G$ and $\rho_f\in R'_\es$, if $k^T$ is the
function on $\H_\es(\Ade_f^T)$ obtained from $g^T$ by the base change
morphism associated to $\eta_\es$, then $\Tr(\rho^T(k^T))=0$.
\end{bulletlist}

Let $S\supset T$ be a finite set of prime numbers such that $g^T=g_{S-T}
\ungras_{\K^S}$, with $g_{S-T}$ a function on $\G(\Ade_{S-T})$. Set
\[f^\infty=h_Tg^T.\]
Let $p\not\in S$ be a prime number big enough for corollary
\ref{cor:stab_FT_IC} to be true. 
Then $f^\infty=f^{\infty,p}\ungras_{\G(\Z_p)}$, and there are functions
$(f^{\infty,p})^{\es}$ and $(f^\infty)^\es$ defined as above,
for every $\es\in\F_\G$.
\quash{
Moreover, after adding a finite number of primes to $S$, we may assume that
$\Tr(\pi_{\es,f}(f^\infty)^\es))=\Tr(\pi_{\es,f}(\ungras_\K)^\es)$ for every
$\pi_{\es,f}\in R_\es(\pi_f)$ such that  $c_\es(\pi_{\es,f})\not=0$.
}

Let $m\in\Z$. Consider the following functions :
\[f^{(m)}=f^{\infty,p}f^{(m)}_pf_\infty\in C_c^{\infty}(\G(\Af^p))C_c^{\infty}
(\G(\Q_p))C^\infty(\G(\R))\]
and
\[f_\H^{(m)}=(f^{\infty,p})^\H f_{\H,p}^{(m)}f_{\H,\infty}\in C_c^{\infty}(\H(
\Af^p))C_c^{\infty}(\H(\Q_p))C^\infty(\H(\R))\]
for every $(\H,s,\eta_0)\in\Ell(\G)$, where :
\begin{bulletlist}
\item $f_p^{(m)}\in\Hecke(\G(\Q_p),\G(\Z_p))$ is the function obtained by
base change from the function $\phi_m^{\G}$ of theorem
\ref{th:points_fixes_Kottwitz};
\item  $f_{\H,p}^{(m)}\in C_c^\infty(\H(\Q_p))$ is the function
obtained by twisted transfer from $\phi_m^{\G}$.

\end{bulletlist}
Then, by corollary \ref{cor:stab_FT_IC} and the choice of $f$,
for $m$ big enough,
\[\Tr(\Phi_\wp^mf,W_\lambda(\pi_f))=
\Tr(\Phi_\wp^mf,W_\lambda)=\sum_{(\H,s,\eta_0)\in
\Ell(\G)}\iota(\G,\H)ST^H(f_{\H}^{(m)}).\]
By proposition \ref{prop:FT_stable_geometrique} and the fact that
$f_\H^{(m)}$ is simply a transfer of $f^{(m)}$ if $\H=\G^*$ (the quasi-split
inner form of $\G$), we get :
\[\Tr(\Phi_\wp^mf,W_\lambda(\pi_f))=
T^G(f^{(m)})+\sum_{\es\in\F_\G}(-1)^{\ell(\es)}
\iota(\es)T^{H_\es}(f^{(m),\es})+\sum_{\es\in\F_\G}(-1)^{\ell(\es)-1}\iota'
(\es)T^{H_\es}(f_\es^{(m)}),\]
where, for every $\es=((\H_1,s_1,\eta_1),\dots,(\H_r,s_r,\eta_r))\in\F_\G$,
we wrote $\H_\es=\H_r$,
\[f^{(m),\es}=(f^{\infty,p})^{\es}(f_p^{(m)})^\es f^\es_\infty\]
and
\[f_\es^{(m)}=(f^{\infty,p})^{\es}f_{\es,p}^{(m)}f^\es_\infty,\]
with
\[f_{\es,p}^{(m)}=(\dots(f_{\H_1,p}^{(m)})^{\H_2}\dots)^{\H_r}.\]

By the calculation of \cite{A-L2} p 267-268 :
\[T^G(f^{(m)})=\sum_{\rho\in\Pi_\G}m_{disc}(\rho)\Tr(\rho(f^{(m)}))\]
\[T^{H_\es}(f^{(m),\es})=\sum_{\rho\in\Pi_\es}m_{disc}(\rho)\Tr(\rho(f^{(m),
\es}))\]
\[T^{H_\es}(f^{(m)}_\es)=\sum_{\rho\in\Pi_\es}m_{disc}(\rho)\Tr(\rho(f^{(m)}_
\es)).\]

Let $\es=((\H_1,s_1,\eta_1),\dots,(\H_r,s_r,\eta_r))\in\F_\G$ and
$\rho=\rho^{\infty,p}\otimes\rho_p\otimes\rho_\infty=\rho_f\otimes\rho_\infty
\in\Pi_\es$. Then
\[\Tr(\rho(f^{(m),\es}))=\Tr(\rho^{\infty,p}((f^{\infty,p})^{\es})\Tr(\rho_p(f^
{(m)}_p)^\es)\Tr(\rho_\infty(f^\es_\infty)).\]
As $\chi_{\es,p}^{-1}(f_p^{(m)})^\es\in\Hecke(\H_r(\Q_p),\H_r(\Z_p))$,
the trace above is $0$ unless $\rho_p\otimes\chi_{\es,p}$ is unramified.
So
\[\Tr(\rho(f^{(m),\es}))=\Tr(\rho_f((f^\infty)^{\es})\Tr(\rho_p(f_p^
{(m)})^\es)\Tr(\rho_\infty(f^\es_\infty)),\]
because both sides are zero if unless $\rho_p\otimes\chi_{\es,p}$ is unramified
and, if $\rho_p\otimes\chi_{\es,p}$ is unramified, then
$\Tr(\rho_p(\chi_{\es,p}\ungras_{\H_r(\Z_p)}))=\dim((\rho_p\otimes\chi_{\es,p})
^{\H_r(\Z_p)})=1$ (and, by
the fundamental lemma, we may assume that $(f_p^{(m)})^\es$
is equal to $\chi_{\es,p}\ungras_{\H_\es(\Z_p)}$).
Assume that $\rho_p\otimes\chi_{\es,p}$ is unramified, and let
$\varphi_{\rho_p\otimes\chi_{\es,p}}
:W_{\Q_p}\fl{}^L\H_{r,\Q_p}$ be a Langlands parameter of
$\rho_p\otimes\chi_{\es,p}$. Then, by
proposition \ref{prop:transfo_Satake_CB} and the calculation of the transfer
of a function in the spherical Hecke algebra in \ref{partie_en_p2}, we get
\[\Tr(\rho_p((f_p^{(m)})^\es))=(N\wp)^{md/2}\sum_{\mu\in M_\es}\Tr(
r_{-\mu}\circ\varphi_{\rho_p\otimes\chi_{\es,p}}(\Phi_\wp^m)).\]

Similarly, using the calculation of the twisted transfer in
\ref{partie_en_p2}, we see that $\Tr(\rho(f_\es^{(m)}))$ is equal to
$0$ if $\rho_p\otimes\chi_{\es,p}$ is ramified, and to
\[\Tr(\rho_f((f^\infty)^\es))\Tr(\rho_\infty(f^\es_
\infty))(N\wp)^{md/2}\sum_{\mu\in M_\es}(-1)^{s(\mu)}\Tr(r_{-\mu}\circ\varphi
_{\rho_p\otimes\chi_{\es,p}}(\Phi_\wp^m))\]
if $\rho_p\otimes\chi_{\es,p}$ is unramified.

Moreover, by the choice of $f^\infty$, if $\rho_f\not\in R_\es(\pi_f)$, then :
\[c_\es(\rho_f)\Tr(\rho_f((f^\infty)^\es))=0.\]

A similar (but simpler) calculation gives, for every
$\rho=\rho_f\otimes\rho_\infty\in\Pi_\G$ : $c_\G(\rho_f)\Tr(\rho_f(f^{\infty,p}
f_p^{(m)}))=0$ if
$\rho$ is ramified at $p$ or if $\rho_f\not\simeq\pi_f$, and, if $\rho_f
\simeq\pi_f$ (so $\rho$ is unramified at $p$), then
\[\Tr(\rho(f^{(m)}))=\dim(\pi_f^{\K})\Tr(\rho_\infty(f_\infty))(N\wp)^{md/2}
\Tr(r_{-\mu_\G}\circ\varphi_{\pi_p}(\Phi_\wp^m)).\]

This calculations imply the equality of the theorem.

\end{proof}

\begin{subremarque}\label{rq:pour_tout_m} Take any $f^\infty$ in
$C_c^\infty(\G(\Af))$. Then
the calculations in the proof of the theorem show that for every
prime number $p$ unramified in $E$ and such that $f^\infty=f^{\infty,p}
\ungras_{\G(\Z_p)}$,
and for every $m\in\Z$ :
\begin{flushleft}$\displaystyle{
\sum_{(\H,s,\eta_0)\in
\Ell(\G)}\iota(\G,\H)ST^H(f_{\H}^{(m)})
=(N\wp)^{md/2}
\sum_{\pi_f}c_\G(\pi_f)\Tr(\pi_f(f^\infty))\Tr(r_{-\mu_G}\circ\varphi_
{\pi_p}(\Phi_\wp^m))}$\end{flushleft}
\begin{center}$\displaystyle{
+(N\wp)^{md/2}\sum_{\es\in\F_\G}(-1)
^{\ell(\es)}\iota(\es)\sum_{\pi_{\es,f}}c_\es(\pi_{\es,f})
\Tr(\pi_{\es,f}((f^\infty)^\es))
}$\end{center}
\begin{flushright}$\displaystyle{
\sum_{\mu\in M_\es}(1-(-1)^{s(\mu)}\frac{\iota'(\es)}{\iota(\es)})
\Tr(r_{-\mu}\circ\varphi_{\pi_{\es,p}\otimes\chi_{\es,p}}(\Phi_\wp^m)),
}$\end{flushright}
where the first (resp. third) sum on the right hand side
is taken over the set of isomomorphism classes of irreducible admissible
respresentations $\pi_f$ (resp. $\pi_{\es,f}$) of $\G(\Af)$ (resp. $\H_\es
(\Af)$) such that $\pi_p$ is unramified (resp. $\pi_{\es,p}\otimes\chi_{\es,p}$
is unramified), and the function $f_{\H,p}^{(m)}$ for $m\leq 0$ is defined
in definition \ref{def:phi_m_G}.

This implies that corollary \ref{cor:stab_FT_IC} is true for every $j\in\Z$,
and not just for $j$ big enough (because that corollary can be rewritten
as an equality $\sum\limits_{i\in I}c_i\alpha_i^j=\sum\limits_{k\in K}d_k
\beta_k^j$, where $(c_i)_{i\in I}$, $(\alpha_i)_{i\in I}$, $(d_k)_{k\in K}$ and
$(\beta_k)_{k\in K}$ are finite families of complex numbers).
This is the statement of proposition
\ref{prop:pour_tout_m} (if $p$ is inert), but note that
proposition \ref{prop:pour_tout_m} was used in the proof of this remark.

\end{subremarque}

For every $i\in\Z$, consider the representation $\Ho^i(M^{\K}(\G,\X)^*_
{\overline{\Q}},IC^{\K}V_{\overline{\Q}})$ of $\Hecke_\K\times\Gal(\overline
{\Q}/F)$. After making $K_\lambda$ bigger, we may assume that all the
$\Hecke_\K$-isotypical
components of this representations are defined over $K_\lambda$. Write
$W_\lambda^i$ for the semi-simplifications of this representations, and let
\[W_\lambda^i=\bigoplus_{\pi_f}W_\lambda^i(\pi_f)\otimes\pi_f^\K\]
be their isotypical decompositions as $\Hecke_\K$-modules (so, as before,
the sum is taken over the set of isomorphism classes of irreducible
admissible representations $\pi_f$ of $\G(\Af)$ such that
$\pi_f^\K\not=0$). Of course,
$W_\lambda=\sum\limits_{i\in\Z}(-1)^i[W_\lambda^i]$ and
$W_\lambda(\pi_f)=\sum\limits_{i\in\Z}(-1)^i[W_\lambda^i(\pi_f)]$ for
every $\pi_f$.
\index{Wilambda@$W_\lambda^i$}
\index{Wilambdapif@$W_\lambda^i(\pi_f)$}

Then, just as in Kottwitz's article \cite{K-LAR} (see also
5.2 of Clozel's article \cite{Cl-RGRA}), we get the following characterization
of the representations $\pi_f$ that appear in $W_\lambda$ :

\begin{subremarque}\label{rq:c_G(pi)} Let $\pi_f$ be an irreducible
admissible representation of $\G(\Af)$ such that $\pi_f^\K\not=0$. Then the
following conditions are equivalent :
\begin{itemize}
\item[(1)] $W_\lambda(\pi_f)\not=0$.
\item[(2)] There exists $i\in\Z$ such that $W_\lambda^i(\pi_f)\not=0$.
\item[(3)]  There exists $\pi_\infty\in\Pi(\G(\R))$ and $i\in\Z$ such that
$m_{disc}(\pi_f\otimes\pi_\infty)\not=0$ and $\Ho^i(\ggoth,\K'_\infty;
\pi_\infty\otimes V)\not=0$.

\end{itemize}
The notations used in condition (3) are those of the proof of lemma
\ref{lemme:conj_RP4}.
Moreover, all this conditions are implied by :
\begin{itemize}
\item[(4)] $c_\G(\pi_f)\not=0$.
\end{itemize}
Assume that, for every $\es\in\F_\G$ and for every $\pi_{e,f}\in
R_\es(\pi_f)$, $c_\es(\pi_{\es,f})=0$. Then (1) implies (4).

\end{subremarque}

\begin{proof} It is obvious that (1) implies (2).

By lemma 3.2 of \cite{K-LAR}, there exists a positive integer $N$ such that,
for every $\pi_\infty\in\Pi(\G(\R))$,
\[\Tr(\pi_\infty(f_\infty))=N^{-1}\sum_{i\in\Z}(-1)^i\dim(\Ho^i(\ggoth,
\K'_\infty,\pi_\infty\otimes V)).\]
This shows in particular that (4) implies (3).

By Matsushima's formula (generalized by Borel and Casselman)
and Zucker's conjecture (proved by Looijenga, Saper-Stern, Looijenga-Rapoport),
for every $i\in\Z$, there is an isomorphism of $\C$-vector spaces :
\[\iota(W_\lambda^i(\pi_f))=\bigoplus_{\pi_\infty\in\Pi(\G(\R))}
m_{disc}(\pi_f\otimes\pi_\infty)
\Ho^i(\ggoth,\K'_\infty,\pi_\infty\otimes V)\]
(remember that $\iota:K_\lambda\fl\C$ is an embedding that was fixed at the
begining of this section). This is explained in the proof of lemma
\ref{lemme:conj_RP4}. The equivalence of (2) and (3) follows from this
formula.

We show that (2) implies (1). This
is done just as in section 6 of \cite{K-LAR}. Let $m$ be the weight of $V$
in the sense of \ref{points_fixes3}. Then the local system $\F^{\K}V$
defined by $V$ is pure of weight $-m$ (cf \ref{points_fixes3}), so the
intersection complexe $IC^{\K}V$ is also pure of weight $-m$. Hence, for
every $i\in\Z$, $W_\lambda^i(\pi_f)$ is pure of weight $-m+i$ as a
representation of $\Gal(\overline{\Q}/F)$ (ie it is unramified and pure
of weight $-m+i$ at almost all places of $F$). In particular, $W_\lambda
^i(\pi_f)$ and $W_\lambda^j(\pi_f)$ cannot have isomorphic irreducible
subquotients if $i\not=j$, so there are no cancellations in the sum
$W_\lambda(\pi_f)=\sum\limits_{i\in\Z}(-1)^i[W_\lambda^i(\pi_f)]$. This
show that (2) implies (1).

We now prove the last statement. By the assumption on $\pi_f$ and
theorem \ref{th:calcul_L_moche}, for almost every prime number $p$ and
every $m\in\Z$,
\[\Tr(\Phi_\wp^m,W_\lambda(\pi_f))=(N\wp)^{md/2}
c_\G(\pi_f)\dim(\pi_f^\K)\Tr(r_{-\mu_G}\circ\varphi_
{\pi_p}(\Phi_\wp^m)).\]
Fix $p$ big enough for this equality to be true. If $W_\lambda(\pi_f)\not=0$, 
then there exists $m\in\Z$ such that $\Tr(\Phi_\wp^m,W_\lambda(\pi_f))\not=0$,
so $c_\G(\pi_f)\not=0$.

\end{proof}

\section{Application to the Ramanujan-Petersson conjecture}
\label{applications2}

We keep the notations of \ref{applications1}, but we take here
$\G=\G(\U(p_1,q_1)\times\dots\times\U(p_r,q_r))$, with
$p_1,q_1,\dots,p_r,q_r\in\Nat$ such that, for every $i\in\{1,\dots,r\}$,
$n_i:=p_i+q_i\geq 1$. Assume that, for every
$i\in\{1,\dots,r\}$, if $n_i\geq 2$, then $q_i\geq 1$.
Write $n=n_1+\dots+n_r$ and $d=\dim M^{\K}(\G,\X)$ (so
$d=p_1q_1+\dots+p_rq_r$).Let $\T$ be the diagonal torus of $\G$.

\begin{theorem}\label{th:conj_RP}
Let $\pi_f$ be an irreducible admissible representation of $\G(\Af)$ such
that there exists an irreducible representation $\pi_\infty$ of $\G(\R)$
with $\Tr(\pi_\infty(f_\infty))\not=0$ and $m_{disc}(\pi_f\otimes
\pi_\infty)\not=0$. For every prime number $p$ where $\pi_f$ is unramified, let
\[(z^{(p)},((z_{1,1}^{(p)},\dots,z_{1,n_1}^{(p)}),\dots,(z_{r,1}^{
(p)},\dots,z_{r,n_r}^{(p)})))\in\widehat{\T}^{\Gal(\overline{\Q}_p/
\Q_p)}\]
be the Langlands parameter of $\pi_p$. Assume that $V$ is pure of weight $0$
in the sense of \ref{points_fixes3} (ie that $\Gr_m$, seen as a subgroup of
the center of $\G$, acts trivially on $V$). Then, for every $p$ where
$\G$ is unramified,
\[|z^{(p)}|=|z_{1,1}^{(p)}\dots z_{1,n_1}^{(p)}|=\dots
=|z_{r,1}^{(p)}\dots z_{r,n_r}^{(p)}|=1.\]
Moreover :

\begin{itemize}
\item[(i)]
Assume that the highest weight of $V$ is regular. Then, for $p$ big enough,
for every $i\in\{1,\dots,r\}$ and
$j\in\{1,\dots,n_i\}$, $|z_{i,j}^{(p)}|=1$.

\item[(ii)] Assume that $r=1$, that $W_\lambda(\pi_f)\not=0$ and that,
for every $\es\in\F_\G$ and for every
$\pi_{e,f}\in R_\es(\pi_f)$,
$c_\es(\pi_{\es,f})=0$. Then, for $p$ big enough :
\begin{bulletlist}
\item if $p$ splits in $E$, then, for every $j\in\{1,\dots,n_1\}$,
$\log_p|z_{1,j}^{(p)}|\in\frac{1}{gcd(p_1,q_1)}\Z$;
\item if $p$ is inert in $E$, then, for every $j\in\{1,\dots,n_1\}$,
$\log_p|z_{1,j}^{(p)}|\in\frac{1}{gcd(2,p_1,q_1)}\Z$.

\end{bulletlist}
\end{itemize}
\end{theorem}

\begin{proof} 
Let $\K$ be a neat open compact subgroup of $\G(\Af)$ such that
$\pi_f^\K\not=\{0\}$.

The center $Z$ of $\G$ is isomorphic in an obvious way to $\G(\U(1)^r)$.
As $\pi_f$ is irreducible, $Z(\Af)$ acts on the space of $\pi_f$ by a
character $\chi:Z(\Af)\fl\C^\times$, that is unramified wherever $\pi_f$ is
and trivial on $\K\cap Z(\Af)$. The character $\chi$ is also trivial on
$Z(\Q)$, because there exists a representation $\pi_\infty$ of $\G(\R)$
such that $\pi_\infty\otimes\pi_f$ is a direct factor of
$L^2(\G(\Q)\sous\G(\Ade),1)$ (where $1$ is the trivial character of
$\A_G(\R)^0$, ie the character by which $\A_G(\R)^0$ acts on $V$).
Hence $\chi$ is trivial on $Z(\Q)(\K\cap Z(\Af))$; as
$Z(\Q)(\K\cap Z(\Af))$ is a subgroup of finite index of $Z(\Af)$,
$\chi$ is of finite order. (As $Z(\R)/\A_G(\R)^0$ is compact, this implies
in particular that the central character of $\pi_\infty\otimes\pi_f$ is
unitary.)

Use \ref{groupes3} to identify $\widehat{Z}$ and
$\C^\times\times(\C^\times)^r$.
For every $p$ where $\pi_f$ is unramified, let
$(y^{(p)},(y_1^{(p)},\dots,y_r^{(p)}))\in\widehat{Z}^{\Gal(\overline{\Q}_p/
\Q_p)}$ be the Langlands parameter of $\chi_p$. As $\chi$ is of finite order,
$|y^{(p)}|=|y_1^{(p)}|=\dots=|y_r^{(p)}|=1$.

The morphism $\widehat{\G}=\C^\times\times\GL_{n_1}(\C)\times
\dots\times\GL_{n_r}(\C)\fl\widehat{Z}=\C^\times\times
(\C^\times)^r$, $(z,(g_1,\dots,g_r))\fle (z,\det(g_1),\dots,
\det(g_r))$ is dual to the inclusion $Z\subset\G$. So, for every
$p$ where $\pi_f$ is unramified,
$z^{(p)}=y^{(p)}$ and $z_{i,1}^{(p)}\dots z_{i,n_i}^{(p)}=
y_i^{(p)}$  for every $i\in\{1,\dots,r\}$. This proves the first statement
of the theorem.

We show (i).
Assume that the highest weight of $V$ is regular. Let $R_\infty$
be the set of $\pi_\infty\in\Pi(\G(\R))$ such that
$\pi_\infty\otimes\pi_f\in\Pi_G$ and $Tr(\pi_\infty(f_\infty))\not=0$.
By the proof of lemma 6.2 of \cite{A-L2}, all the representations
$\pi_\infty\in\Pi(\G(\R))$ such that $\Tr(\pi_\infty(f_\infty))\not=0$
are in the discrete series. So $R_\infty$ is contained in the discrete series
$L$-packet associated to the contragredient of $V$.
In particular, the
function $\pi_\infty\fle\Tr(\pi_\infty(f_\infty))$ is constant on $R_\infty$,
so
\[c_\G(\pi_f)=\sum_{\pi_\infty\in R_\infty}m_{disc}(\pi_\infty\otimes\pi_f)
\Tr(\pi_\infty(f_\infty))\not=0.\]
We prove the result by induction on the set of $(n_1,\dots,n_r)\in (\Nat^*)^r$
such that $n_1+\dots+n_r=n$, with the ordering :
$(n'_1,\dots,n'_{r'})<(n_1,\dots,n_r)$ if and only if $r'>r$.

Assume first that, for every $\es\in\F_\G$ and for every $\pi_{\es,f}\in
\Pi_\es(\pi_f)$, $c_\es(\pi_{\es,f})=0$.
Let $p$ be a prime number big enough for theorem \ref{th:calcul_L_moche} to be
true. Then, for every $m\in\Z$ :
\[\Tr(\Phi_\wp^m,W_\lambda(\pi_f))=(N\wp)^{md/2}c_\G(\pi_f)\dim(\pi_f^\K)
\Tr(r_{-\mu_G}\circ\varphi_{\pi_p}(\Phi_\wp^m)).\]
Let $x_1,\dots,x_s\in\C$ be the eigenvalues of $\Phi_\wp$ acting on
$W_\lambda(\pi_f)$, and $a_1,\dots,a_s\in\Z$ be their multiplicities.
For every $m\in\Z$,
\[\Tr(\Phi_\wp^m,W_\lambda(\pi_f))=\sum_{i=1}^sa_ix_i^m.\]
By lemma \ref{lemme:conj_RP4}, the cohomology of $IC^{\K}V$ is concentrated in
degree $d$. As $IC^{\K}V$ is pure of weight $0$ (because $V$ is pure of weight
$0$),
$\log_p|x_1|=\dots\log_p|x_r|=n(\wp)d/2$, where $n(\wp)=\log_p(N\wp)$.

On the other hand, by lemma \ref{lemme:conj_RP1}, for every $m\in\Z$ :
\[\Tr(r_{-\mu_G}\otimes\varphi_{\pi_p}(\Phi_\wp^{2m}))=
(z^{(p)})^{-2m}
\sum_{{J_1\subset\{1,\dots,n_1\}}\atop{|J_1|=p_1}}\dots
\sum_{{J_r\subset\{1,\dots,n_r\}}\atop{|J_r|=p_r}}
\prod_{i=1}^r\prod_{j\in J_i}(z_{i,j}^{(p)})^{-2m
[F_\wp:\Q_p]}.\]
As $|z^{(p)}|=1$, this implies that, for all $J_1\subset\{1,\dots,n_1\},\dots,
J_r\subset\{1,\dots,n_r\}$ such that $|J_1|=p_1,\dots,|J_r|=p_r$,
\[\sum_{i=1}^r\sum_{j\in J_i}\log_p|z_{i,j}^{(p)}|=0.\]
By the first statement of the theorem and lemma
\ref{lemme:conj_RP2}, we get $\log_p|z_{i,j}^{(p)}|=0$, ie $|z_{i,j}^{(p)}|=1$,
for every $i\in\{1,\dots,r\}$ and $j\in\{1,\dots,n_i\}$.

Assume now that there exists $\es\in\F_\G$ and $\pi_{\es,f}\in\Pi_\es(\pi_f)$
such that $c_\es(\pi_{\es,f})\not=0$.
Write (with the notations of \ref{applications1})
\[\H_\es=\G(\U^*(n'_1)\times\dots\times\U^*(n'_{r'})).\]
Of course, $(n'_1,\dots,n'_{r'})<(n_1,\dots,n_r)$.

Let $\T_H$ be the diagonal torus of $\H_\es$. If $p$ is a prime number where
$\pi_{\es,f}$ is unramified, let
\[(t^{(p)},((t_{1,1}^{(p)},\dots,t_{1,n'_1}^{(p)}),
\dots,(t_{r',1}^{(p)},\dots,t_{r',n'_{r'}}^{(p)})))
\in\widehat{\T}_H^{\Gal(\overline{\Q}_p/\Q_p)}\]
be the Langlands parameter of $\pi_{\es,p}$. By the definition of
$R_\es(\pi_f)$ (and the fact that, in proposition
\ref{prop:prolongement_eta_0}, we chose a unitary character $\mu$),
up to a permutation of the $t_{i,j}^{(p)}$, there is an equality
\[(z^{(p)},u_{1,1}^{(p)}z_{1,1}^{(p)},\dots,u_{1,n_1}^{(p)}z_{1,n_1}^{(p)},
\dots,u_{r,1}^{(p)}z_{r,1}^{(p)},\dots,u_{r,n_r}^{(p)}z_{r,n_r}^{(p)})=
(t^{(p)},t_{1,1}^{(p)},\dots,t_{1,n'_1}^{(p)},\dots,t_{r',1}
^{(p)},\dots,t_{r',n'_{r'}}^{(p)})\]
for almost every $p$, where the $u_{i,j}^{(p)}$ are complex numbers with
absolute value $1$.
So it is enough to show that $|t_{i,j}^{(p)}|=1$ for all $i,j$, if $p$ is
big enough.

As $c_\es(\pi_{\es,f})\not=0$, there exists $\pi_{\es,\infty}\in\Pi(
\H_\es(\R))$ and an elliptic Langlands parameter $\varphi_H:W_\R\fl {}^L\H_
{\es,\R}$ such that
$m_{disc}(\pi_{\es,\infty}\otimes\pi_{\es,f})\not=0$,
$\eta_\es\circ\varphi_H$ is $\widehat{\G}$-conjugate to $\varphi$ and
\[\Tr(\pi_{\es,\infty}(f_{\varphi_H}))\not=0,\]
where $f_{\varphi_H}$ is the stable cuspidal function associated to $\varphi_H$
defined at the end of \ref{stabilisation_K1}. 
By lemma \ref{lemme:conj_RP3},
$\varphi_H$ is the Langlands parameter of a $L$-packet of the discrete
series of $\H_\es(\R)$ associated to an irreducible algebraic representation
of $\H_{\es,\C}$ with regular highest weight and pure of weight $0$.
So the representation $\pi_{\es,f}$ of $\H_\es(\Af)$ satisfies all the
conditions of point (i) of the theorem,
and we can apply the induction hypothesis to
finish the proof. 

We show (ii). Without the assumption on the highest weight of $V$,
the complex $IC^{\K}V$ is still pure of weight $0$, but its cohomology
is not necessarily concentrated in degree $d$.
By the hypothesis on $\pi_f$, for $p$ big enough and for every $m\in\Z$,
there is an equality
\[(N\wp)^{md/2}c_\G(\pi_f)\dim(\pi_f^\K)\Tr(r_{-\mu_G}\circ\varphi_{\pi_p}
(\Phi_\wp^m))=\Tr(\Phi_\wp^m,W_\lambda(\pi_f))
=\sum_{i=1}^sa_ix_i^m,\]
where, as in (i), $x_1,\dots,x_s\in\C$ are the eigenvalues of $\Phi_\wp$
acting on $W_\lambda(\pi_f)$ and $a_1,\dots,a_s\in\Z$ are their multiplicities.
In particular, all the $a_i$ have the same sign (the sign of $c_\G(\pi_f)$),
so $W_\lambda(\pi_f)$ is concentrated either in odd degree or in even degree,
and the weights of $W_\lambda(\pi_f)$ are either all even or all odd.
By applying the same
reasoning as above, we find, for $p$ big enough, a linear system :
\[\sum_{j\in J}\log_p|z_{1,j}^{(p)}|=\frac{1}{2}w_J,\quad J\subset
\{1,\dots,n_1\},\quad|J|=p_1,\]
where the $w_J$ are in $\Z$ and all have the same parity.
As $p_1<n_1$ if $n_1\geq 2$, this implies that
$\log_p|z_{1,j}^{(p)}|-\log_p|z_{1,j'}^{(p)}|\in\Z$ for every
$j,j'\in\{1,\dots,n_1\}$. On the other hand, we know that
$\sum\limits_{j=1}^{n_1}\log_p|z_{1,j}^{(p)}|=0$. So, for every
$J\subset\{1,\dots,n_1\}$ such that $|J|=q_1$,
$\sum\limits_{j\in J}\log_p|z_{1,j}^{(p)}|\in\Z$. 

Let $\alpha\in\R$ be such that $\log_p|z_{1,1}^{(p)}|-\alpha\in\Z$. Then
$\log_p|z_{1,j}^{(p)}|-\alpha\in\Z$ for every $j\in\{1,\dots,n_1\}$, so
$p_1\alpha,q_1\alpha\in\Z$, and $gcd(p_1,q_1)\alpha\in\Z$. Assume that
$p$ is inert in $E$. Then the fact that
$(z^{(p)},(z_{1,1}^{(p)},\dots,z_{1,n_1}^{(p)}))$ is
$\Gal(\overline{\Q}_p/\Q_p)$-invariant implies that, for every
$j\in\{1,\dots,n_1\}$, $\log_p|z_{1,j}^{(p)}|+\log_p|z_{1,n_1+1-j}^{(p)}|=0$.
So $2\alpha\in\Z$.

\end{proof}

\begin{lemma}\label{lemme:conj_RP1} Use the notations of theorem
\ref{th:conj_RP} above. Fix a prime number $p$ where $\pi_f$ is unramified
and $m\in\Z$. Then
\[\Tr(r_{-\mu_G}\circ\varphi_{\pi_p}(\Phi_\wp^m))=
(z^{(p)})^{-m}\sum_{J_1}\dots\sum_{J_r}\prod_{i=1}^r
\prod_{j\in J_i}(\pm z_{i,j}^{(p)})^{-m[F_\wp:\Q_p]},\]
where :
\begin{itemize}
\item[(i)] if $F=\Q$, $p$ is inert in $E$ and $m$ is odd, then, for every
$i\in\{1,\dots,r\}$, the $i$-th sum is taken over the set of subsets $J_i$ of
$\{1,\dots,n_i\}$ such that
\[\{1,\dots,n_i\}-J_i=\{n_i+1-j,j\in J_i\};\]
\item[(ii)] in all other cases, the $i$-th sum is taken over the set of
subsets $J_i$ of $\{1,\dots,n_i\}$ such that $|J_i|=p_i$, and all the signs
are equal to $1$.

\end{itemize}
\end{lemma}

\begin{proof} To make notations simpler, we assume that $r=1$. (The proof is
exactly the same in the general case.)
We first determine the representation $r_{-\mu_G}$ of ${}^L\G_F$. As
$r_{-\mu_G}$ is the contragredient of $r_{\mu_G}$, it is enough to
calculate $r_{\mu_G}$. Remember that $\T$ is the diagonal torus of $\G$,
and that $\widehat{\T}=\C^\times\times (\C^\times)^{n_1}\subset
\widehat{\G}=\C^\times\times\GL_{n_1}(\C)$.
The cocharacter $\mu_G$ of $\T$ corresponds to the following character
of $\widehat{\T}$ :
\[(\lambda,(\lambda_i)_{1\leq i\leq n_1})
\fle\lambda\prod_{i=1}^{p_1}\lambda_i.\]
So the space of $r_{\mu_G}$ is $V_\mu=\bigwedge\limits^{p_1}\C^{n_1}$,
where $\GL_{n_1}(\C)$ acts by $\bigwedge\limits^{p_1}$ of the standard
representation, and $\C^\times$ acts by the character $z\fle z$.
Let $(e_1,\dots,e_{n_1})$ be the canonical basis of $\C^{n_1}$. Then the
family $(e_{i_1}\wedge\dots\wedge e_{i_{p_1}})_{1\leq i_1<\dots<i_{p_1}
\leq n_1}$ is a basis of $V_\mu$.
From the definition of $r_{\mu_G}$ (cf lemma \ref{def_r_mu}), it is easy to
see that $W_E$ acts trivially on $V_\mu$ and that, if $F=\Q$ (so $n_1$ is
even and $p_1=n_1/2$), then an element of $W_\Q-W_E$ sends
$e_{i_1}\wedge\dots\wedge e_{i_{p_1}}$, where $1\leq i_1<\dots<i_{p_1}
\leq n_1$, to $\pm e_{j_1}\wedge\dots\wedge e_{j_{p_1}}$, with
$1\leq j_1<\dots<j_{p_1}\leq n_1$ and $\{n_1+1-j_1,\dots,n_1+1-j_{p_1}\}=\{1,
\dots,n_1\}-\{i_1,\dots,i_{p_1}\}$.
By definition of the Langlands parameter, we may assume that
$\varphi_{\pi_p}(\Phi_\wp)=(((z^{(p)})^{[F_\wp:\Q_p]},((z_1^{(p)})^
{[F_\wp:\Q_p]},\dots,(z_{n_1}^{(p)})^{[F_\wp:\Q_p]})),\Phi_\wp)$
(remember that $\Phi_\wp$ is a lift in $\Gal(\overline{\Q}_p/F_\wp)$ of
the geometric Frobenius).
If $F=E$, $p$ is split in $E$ or $m$ is even, then the image of
$\Phi_\wp^m$ in $W_\Q$ is an element of $W_E$, so $\Phi_\wp^m$ acts trivially
on $V_\mu$ (if $p$ is split in $E$, this comes from the fact that the
image of $W_{\Q_p}$ in $W_\Q$ is included in $W_E$).
If $F=\Q$ and $p$ is inert in $E$, then $\Gal(E_p/\Q_p)\iso\Gal(E/\Q)$, and
the image of $\Phi_\wp$ in $\Gal(E_p/\Q_p)$ generates $\Gal(E_p/\Q_p)$, so
$\Phi_\wp^m\not\in W_E$ for $m$ odd.
The formula of the lemma is a consequence of these remarks and of the
explicit description of $r_{\mu_G}$.

\end{proof}

\begin{lemma}\label{lemme:conj_RP2}\begin{itemize}
\item[(i)] Let $n,p\in\Nat$ be such that
$1\leq p\leq max(1,n-1)$.
Then there exist $J_1,\dots,J_n\subset\{1,\dots,n\}$ such that
$|J_1|=\dots=|J_n|=p$ and that the only solution of the system of linear
equations
\[\sum_{j\in J_i}X_j=0,\qquad 1\leq i\leq n,\]
is the zero solution.
\item[(ii)] Let $r\in\Nat^*$, $n_1,\dots,n_r\geq 2$ and $p_1,\dots,p_r\in\Nat$
be such that $1\leq p_i\leq n_i-1$ for $1\leq i\leq r$. For every
$i\in\{1,\dots,r\}$, choose subsets $J_{i,1},\dots,J_{i,n_i}$ of
$\{1,\dots,n_i\}$ of cardinality $p_i$ and satisfying the property of (i).
Then the only solution of the system of linear equations
\[\left\{\begin{array}{ll}\displaystyle{\sum_{j=1}^{n_i}X_{i,j}=0,} &
\qquad 1\leq i\leq r,\\
\displaystyle{\sum_{i=1}^r\sum_{j\in J_{i,k_i}}X_{i,j}=0,} &
\qquad (k_1,\dots,k_r)\in\{1,
\dots,n_1\}\times\dots\times\{1,\dots,n_r\},\end{array}\right.\]
is the zero solution.

\end{itemize}
\end{lemma}

\begin{proof} We show (i) by induction on $n$.
If $n=1$, the result is obvious. Suppose that $n\geq 2$, and let
$p\in\{1,\dots,n-1\}$. Assume first that $p\leq n-2$. Then, by the induction
hypothesis, there exist $J_2,\dots,J_n\subset\{2,\dots,n\}$ of cardinality
$p$ such that the only solution of the system of linear equations
(with unknowns $X_2,\dots,X_n$)
\[\sum_{j\in J_i}X_j=0,\qquad 2\leq i\leq n\]
is the zero solution. Take $J_1=\{1,2,\dots,p\}$. It is clear that
$J_1,\dots,J_n$ satisfy the condition of (i). Assume now that $p=n-1$.
For every $i\in\{1,\dots,n\}$, let $J_i=\{1,\dots,n\}-\{i\}$. To show that
$J_1,\dots,J_n$ satisfy the condition of (i), it is enough to show that
$\det(A-I_n)\not=0$, where $A\in M_n(\Z)$ is the matrix all of whose entries
are equal to $1$. But it is clear that the kernel of $A$ is of dimension
$n-1$ and that $n$ is an eigenvalue of $A$, so $A$ has no eigenvalue
$\lambda\not\in\{0,n\}$. In particular, $\det(A-I_n)\not=0$.

We show (ii) by induction on $r$. The case $r=1$ is obvious, so we assume
that $r\geq 2$. Let (S) be the system of linear equations of (ii).
For $2\leq i\leq r$, fix $k_i\in\{1,\dots,n_i\}$. Then, by the case $r=1$,
the system (S') :
\[\sum_{i=1}^r\sum_{j\in J_{i,k_i}}X_{i,j}=0,\qquad 
k_1\in\{1,\dots,n_1\},\]
has a unique solution in $(X_{1,1},\dots,X_{1,n_1})$, that is equal to the
obvious solution
\[X_{1,1}=\dots=X_{1,n_1}=-\frac{1}{p_1}\sum_{i=2}^r
\sum_{j\in J_{i,k_i}}X_{i,j}.\]
So the system (S') and the equation $\sum\limits_{j=1}^{n_1}X_{1,j}=0$ imply :
\[X_{1,1}=\dots=X_{1,n_1}=\sum_{i=2}^r\sum_{j\in J_{i,k_i}}X_{i,j}=0.\]
To finish the proof, apply the induction hypothesis to the system
analogous to (S) but with $2\leq i\leq r$.

\end{proof}

We now take $\G=\G(\U^*(n_1)\times\dots\times\U^*(n_r))$. 
Then $\widehat{\G}=\C^\times\times\GL_{n_1}(\C)\times\dots\times\GL_{n_r}(\C)$,
with the action of $W_\Q$ described in \ref{groupes3}. Let $\T$ be the
elliptic maximal torus of $\G$ defined in \ref{serie_discrete1}, and
$u_G\in\G(\C)$ be the element defined in \ref{serie_discrete1}, so that
$u_G^{-1}\T u_G$ is the diagonal torus of $\G$. Let
$\B\supset\T$ be the Borel subgroup of $\G_\C$ image by $\Int(u_G)$ of
the group of upper triangular matrices (we identify $\G_\C$ to
$\Gr_{m,\C}\times\GL_{n_1,\C}\times\dots\times\GL_{n_r}(\C)$ as in
\ref{groupes3}). Identify $\T_\C$ to $\Gr_{m,\C}\times\Gr_{m,\C}^{n_1}
\times\dots\times\Gr_{m,\C}^{n_r}$ and $\widehat{\T}$ to $\C^\times\times
(\C^\times)^{n_1}\times\dots\times(\C^\times)^{n_r}$ as in
\ref{serie_discrete1}. 

Let $V$ be an irreducible algebraic representation of $\G_\C$. Let
$\as=(a,(a_{i,j})_{1\leq i\leq r,1\leq j\leq n_i})\in X^*(\T)$ be the highest
weight of $V$ relatively to $(\T,\B)$; the notation means that $\as$ is the
character
\[(z,(z_{i,j})_{1\leq i\leq r,1\leq j\leq n_i})\fle
z^a\prod_{i=1}^r\prod_{j=1}^{n_i}z_{i,j}^{a_{i,j}}.\]
By definition of the highest weight, $a,a_{i,j}\in\Z$ and
$a_{i,1}\geq a_{i,2}\geq\dots\geq a_{i,n_i}$ for every $i\in\{1,\dots,r\}$.
Notice also that the weight of $V$, in the sense of \ref{points_fixes3},
is $2a+\sum\limits_{i=1}^r\sum\limits_{j=1}^{n_i}a_{i,j}$.

Let $(\H,s,\eta_0)$ be the elliptic endoscopic triple for $\G$ associated
to $((n_1^+,n_1^-),\dots,(n_r^+,n_r^-))$ as in proposition
\ref{prop:groupes_endoscopiques}. Then
$\H=\G(\U^*(n_1^+)\times\U^*(n_1^-)\times\dots\times\U^*(n_r^+)\times
\U^*(n_r^-))$, and we define an elliptic maximal torus $\T_H$ of $\H$ and
a Borel subgroup $\B_H\supset\T_H$ of $\H_\C$ in the same way as
$\T$ and $\B$.
Let
\[\Omega_*=\{\omega=(\omega_1,\dots,\omega_r)\in\Sgoth_{n_1}\times\dots\times
\Sgoth_{n_r}|\forall i,\omega^{-1}_{i|\{1,\dots,n_i^+\}}\mbox{ and }\omega^{-1}
_{i|\{n_i^++1,\dots,n_i\}}\mbox{ are increasing}\}.\]
$\Omega_*$ is the set of representatives of $\Omega(\T_H(\C),\H(\C))\sous
\Omega(\T(\C),\G(\C))$ determined by $\B$ and $\B_H$ as in
\ref{serie_discrete3}.

\begin{lemma}\label{lemme:conj_RP3} 
Let $\varphi:W_\R\fl{}^L\G_\R$ be a Langlands parameter of the $L$-packet of
the discrete series of $\G(\R)$ associated to $V$ and $\eta:{}^L\H_\R\fl{}^L
\G_\R$ be a $L$-morphism extending $\eta_0$ as in proposition 
\ref{prop:prolongement_eta_0}. Remember that we wrote $\Phi_H(\varphi)$
for the set of equivalence classes of Langlands parameters
$\varphi_H:W_\R\fl{}^L\H_\R$ such that $\eta\circ\varphi_H$ and $\varphi$
are equivalent.

Then every $\varphi_H\in\Phi_H(\varphi)$ is the
parameter of a $L$-packet of the discrete series of $\H(\R)$ corresponding to
an algebraic representation of $\H_\C$; this algebraic
representation has a regular highest weight if $\as$ is regular, and its
weight in the sense of \ref{points_fixes3} is equal to the weight of $V$.

\end{lemma}

\begin{proof} We may assume that
\[\varphi(\tau)=((1,(\Phi_{n_1}^{-1},\dots,\Phi_{n_r}^{-1})),\tau)\]
and that, for every $z\in\C^\times$,
\[\varphi(z)=((z^a\overline{z}^{a+S},(B_1(z),\dots,B_r(z))),z),\]
where
\[S=\sum_{i=1}^r\sum_{j=1}^{n_i}a_{i,j}\]
and
\[B_i(z)=diag(z^{\frac{n_i-1}{2}+a_{i,1}}\overline{z}^{\frac{1-n_i}{2}
-a_{i,1}},z^{\frac{n_i-3}{2}+a_{i,2}}\overline{z}^{\frac{3-n_i}{2}
-a_{i,2}},\dots,z^{\frac{1-n_i}{2}+a_{i,n_i}}
\overline{z}^{\frac{n_i-1}{2}-a_{i,n_i}}).\]
(Remember that $W_\R=W_\C\sqcup W_\C\tau$, with $W_\C=\C^\times$, $\tau^2=-1$
and, for every $z\in\C^\times$, $\tau z\tau^{-1}=\overline{z}$.)

Let $C$ be the odd integer associated to $\eta$ as below proposition
\ref{prop:prolongement_eta_0}.
Let $\omega=(\omega_1,\dots,\omega_r)\in\Omega_*$, and let $\varphi_H$ be the
element of $\Phi_H(\varphi)$ associated to $\omega$ as in
\ref{serie_discrete3}.
Write, for every $i\in\{1,\dots,r\}$, $j_{i,s}=\omega_i^{-1}(s)$ if
$1\leq s\leq n_i^+$ and $k_{i,t}=\omega_i^{-1}(t+n_i^+)$ if $1\leq t
\leq n_i^-$. Then we may assume that
\[\varphi_H(\tau)=((1,(\Phi_{n_1^+}^{-1},\Phi_{n_1^-}^{-1},\dots,
\Phi_{n_r^+}^{-1},\Phi_{n_r^-}^{-1})),\tau)\]
and that, for every $z\in\C^\times$,
\[\varphi_H(z)=((z^a\overline{z}^{a+S},(B_1^+(z),
B_1^-(z),\dots,B_r^+(z),B_r^-(z))),z),\]
with
\[B_i^+(z)=diag(z^{\frac{n_i^+-1}{2}+a_{i,1}^+}\overline
{z}^{\frac{1-n_i^+}{2}-a_{i,1}^+},z^{\frac{n_i^+-3}{2}+
a_{i,2}^+}\overline{z}^{\frac{3-n_i^+}{2}-a_{i,2}^+},
\dots,z^{\frac{1-n_i^+}{2}+a_{i,n_i^+}^+}\overline{z}^
{\frac{n_i^+-1}{2}-a_{i,n_i^+}^+}),\]
where $a_{i,s}^+=a_{i,j_{i,s}}+s-j_{i,s}+n_i^-(1-C)/2\in\Z$,
and with
\[B_i^-(z)=diag(z^{\frac{n_i^--1}{2}+a_{i,1}^-}\overline
{z}^{\frac{1-n_i^-}{2}-a_{i,1}^-},z^{\frac{n_i^--3}{2}+
a_{i,2}^-}\overline{z}^{\frac{3-n_i^-}{2}-a_{i,2}^-},
\dots,z^{\frac{1-n_i^-}{2}+a_{i,n_i^-}^-}\overline{z}^
{\frac{n_i^--1}{2}-a_{i,n_i^-}^-}),\]
where $a_{i,t}^-=a_{i,k_{i,t}}+t-k_{i,t}+n_i^+(1+C)/2\in\Z$.
Let $i\in\{1,\dots,r\}$. For all $s\in\{1,\dots,n_i^+-1\}$ and
$t\in\{1,\dots,n_1^--1\}$,
\[a^+_{i,s}-a^+_{i,s+1}=(a_{i,j_{i,s}}-a_{i,j_{i,s+1}})+
(j_{i,s+1}-j_{i,s}-1)\]
\[a^-_{i,t}-a^-_{i,t+1}=(a_{i,k_{i,t}}-a_{i,k_{i,t+1}})+
(k_{i,t+1}-k_{i,t}-1),\]
so $a^+_{i,s}\geq a^+_{i,s+1}$ and $a^-_{i,t}\geq a^-_{i,t+1}$, and the
inequalities are strict if $\as$ is regular.
Notice also that
\[S_H:=\sum_{i=1}^r\sum_{s=1}^{n_i^+}a_{i,s}^++\sum_{i=1}^r\sum_{t=1}^{n_i^-}
a_{i,t}^-=S.\]

So $\varphi_H$ is the paramater of the discrete series of $\H(\R)$
associated to the algebraic representation of $\H_\C$ of highest weight
$(a,((a_{i,s}^+)_{1\leq s\leq n_i^+},(a_{i,t}^-)_{1\leq t\leq n_i^-})_{1\leq i
\leq r})$. This representation has a regular highest weight if $\as$ is
regular by the above calculations, and its weight in the sense of
\ref{points_fixes3} is the same as the weight of $V$ because $2a+S_H=2a+S$.

\end{proof}

We use again the notations of the beginning of this section.

\begin{lemma}\label{lemme:conj_RP4} If the highest weight of $V$ is regular,
then, for every neat open compact subgroup $\K$ of $\G(\Af)$, the cohomology
of the complex $IC^{\K}V$ is concentrated in degree $d$.

\end{lemma}

\begin{proof} By Zucker's conjecture (proved by Looijenga \cite{Lo},
Looijenga-Rapoport \cite{LoR} and Saper-Stern \cite{SS}), the intersection
cohomology of $M^{\K}(\G,\X)^*(\C)$ with coefficients in $\F^{\K}V$
is isomorphic to the $L^2$-cohomology of $M^{\K}(\G,\X)(\C)$ with
coefficients in $\F^{\K}V$.
\index{Zucker's conjecture}
By a result of Borel and Casselman (theorem 4.5 of \cite{BC}),
the $\H^q$ of this
$L^2$-cohomology is isomorphic (as a representation of
$C_c^\infty(\K\sous\G(\Af)/\K)$) to
\[\bigoplus_{\pi}m_{disc}(\pi)(\H^q(\mathfrak{g},\K'_\infty;\pi_\infty
\otimes V)\otimes\pi_f^{\K}),\]
where the sum is taken over the set of isomorphism classes of irreducible
admissible representations of $\G(\Ade)$, $\mathfrak{g}=Lie(\G(\R))\otimes\C$
and
$\K'_\infty=\K_\infty\A_G(\R)^0$, with $\K_\infty$ a maximal compact subgroup
of $\G(\R)$. By the proof of lemma 6.2 of \cite{A-L2}, if $\pi_\infty$ is an
irreducible admissible representation of $\G(\R)$ such that
$\H^*(\mathfrak{g},\K'_\infty;\pi_\infty\otimes V)\not=0$, then $\pi_\infty$
is in the discrete series of $\G(\R)$ (this is the only part where we use
the fact that the highest weight of $V$ is regular).
By theorem II.5.3 of \cite{BW}, if $\pi_\infty$ is in the discrete series of
$\G(\R)$, then $\H^q(\mathfrak{g},\K'_\infty;\pi_\infty\otimes V)=0$
for every $q\not=d$.

\end{proof}

\chapter{The twisted trace formula}
\label{GL_n_applications}

\section{Non-connected groups}
\label{GL_n_applications1}

We first recall some definitions from section 1 of \cite{A-LBWOI}.

Let $\til{\G}$ be a reductive group (not necessarily connnected) over a field
$K$. Fix a connected component $\G$ of $\til{\G}$, and assume that $\G$
generates $\til{\G}$ and that $\G(K)\not=\varnothing$.
Let $\G^0$ be the connected component of $\til{\G}$ that contains $1$.

Consider the polynomial
\[\det((t+1)-\Ad(g),Lie(\G^0))=\sum_{k\geq 0}D_k(g)t^k\]
on $\G(K)$. The smallest integer $k$ for which $D_k$ does not vanish
identically is called the \emph{rank} of $\G$; we will denote by $r$.
An element $g$ of $\G(K)$ is called \emph{regular} if $D_r(g)\not=0$.
\index{rank (of a connected component in a non-connected group)}
\index{regular element}

A \emph{parabolic subgroup} of $\til{\G}$ is the normalizer in
$\til{\G}$ of a parabolic subgroup of $\G^0$. A \emph{parabolic subset}
\index{parabolic subset}
of $\G$ is a non-empty subset of $\G$ that is equal to the intersection of $\G$
with a parabolic subgroup of $\til{\G}$. If $\Pa$ is a parabolic subset of
$\G$, write $\til{\Pa}$ for the subgroup of $\til{\G}$ generated by $\Pa$ and
$\Pa^0$ for the intersection $\til{\Pa}\cap\G^0$ (then $\til{\Pa}=\Nor_
{\til{\G}}(\Pa^0)$ and $\Pa=\til{\Pa}\cap\G$). 

Let $\Pa$ be a parabolic subset of $\G$. The \emph{unipotent radical} $\N_P$
\index{unipotent radical of a parabolic subset}
of $\Pa$ is by definition the unipotent radical of $\Pa^0$. A
\emph{Levi component} $\M$ of $\Pa$ is a subset of $\Pa$ that is equal to
$\til{\M}\cap\G$, where $\til{\M}$ is the normalizer in $\til{\G}$ of a Levi
component $\M^0$ of $\Pa^0$. If $\M$ is a Levi component of $\Pa$, then
$\Pa=\M\N_P$. 
\index{Levi component of a parabolic subset}

A \emph{Levi subset} of $\G$ is a Levi component of a parabolic subset of
$\G$. Let $\M$ be a Levi subset of $\G$.
\index{Levi subset}
Let $\til{\M}$ be the subgroup
of $\til{\G}$ generated by $\M$, $\M^0=\G^0\cap\til{\M}$, $\A_M$ be the
maximal split subtorus of the centralizer of $\M$ in $\M^0$ (so $\A_M\subset
\A_{M^0}$), $X^*(\M)$ be the group of characters of $\til{\M}$ that are defined
over $K$, $\agoth_M=\Hom(X^*(\M),\R)$ and
\[n_M^G=|\Nor_{\G^0(\Q)}(\M)/\M^0(\Q)|.\]
\index{AG@$\A_G$\quad (twisted case)}
\index{agothG@$\agoth_G$\quad (twisted case)}
\index{nMG@$n_M^G$\quad (twisted case)}

Fix a minimal parabolic subgroup $\Pa_0$ of $\G^0$ and a Levi subgroup
$\M_0$ of $\Pa_0$. Write $\A_0=\A_{M_0}$ and $\agoth_0=\agoth_{\M_0}$. If $\Pa$
is a parabolic subset of $\G$ such that $\Pa^0\supset\Pa_0$, then $\Pa$ has a
unique Levi component $\M$ such that $\M^0\supset\M_0$; write $\M_P=\M$.
Let $\Phi(\A_{M_P},\Pa)$ be the set of roots of $\A_{M_P}$ in $Lie(\N_P)$.

Let $W_0^G$ be the set of linear automorphisms of $\agoth_0$ induced
by elements of $\G(K)$ that normalize $\A_0$, and $W_0=W_0^{G^0}$. The
group $W_0$ acts on $W_0^G$ on the left and on the right, and both these
actions are simply transitive.
\index{W0G@$W_0^G$}
\index{W0G0@$W_0^{G^0}$}

Here, we will be interested in the case where $\til{\G}=\G^0\rtimes\langle
\theta\rangle$
and $\G=\G^0\rtimes\theta$, where $\G^0$ is a connected reductive group
over $K$ and $\theta$ is an automorphism of finite order of $\G^0$.

In this situation, we say that an element $g\in\G^0(K)$ is
\emph{$\theta$-semi-simple} (resp. \emph{$\theta$-regular}, resp. \emph{
strongly $\theta$-regular}) if $g\theta\in\G(K)$ is semi-simple (resp. regular,
resp. strongly regular) in $\til{\G}$ (an element of $\gamma$ of $\til{\G}(K)$
is called strongly regular if its centralizer is a torus.)
\index{strongly regular}
\index{$\theta$-semi-simple}
\index{$\theta$-regular}
\index{strongly $\theta$-regular}
Let $\G^0_{\theta-reg}$ be the open subset of $\theta$-regular elements in
$\G^0$.
\index{G0thetareg@$\G^0_{\theta-reg}$\quad set of $\theta$-regular elements}
We say that $g_1,g_2\in\G^0(K)$ are \emph{$\theta$-conjugate} if
$g_1\theta,g_2\theta\in\G(K)$ are conjugate under $\til{\G}(K)$.
\index{$\theta$-conjugate}
If
$g\in\G^0(K)$, let $\Cent_{\G^0}(g\theta)$ be the centralizer of $g\theta\in
\G(K)$ in $\G^0$; we call this group the \emph{$\theta$-centralizer} of $g$.
\index{$\theta$-centralizer}
Write $\G^0_{g\theta}$ for the connected component of $1$ in $\Cent_{\G^0}
(g\theta)$.
\index{G0gtheta@$\G^0_{g\theta}$}
Finally, we say that an element $g\in\G^0(K)$ is
\emph{$\theta$-elliptic} if $\A_{\G^0_{g\theta}}=\A_G$.
\index{$\theta$-elliptic}

Assume that there exists a $\theta$-stable minimal parabolic subgroup
$\Pa_0$ of $\G^0$ and a $\theta$-stable Levi subgroup $\M_0$ of $\Pa_0$.
We say that a parabolic subset $\Pa$ of $\G$ is \emph{standard} if
$\til{\Pa}\supset\Pa_0\rtimes\langle\theta\rangle$, and that a Levi subset
$\M$ of $\G$ is \emph{standard} if there exists a standard parabolic subset
$\Pa$ such that $\M=\M_P$ (so $\til{\M}\supset\M_0\rtimes\langle\theta
\rangle$).
\index{standard parabolic subset}
\index{standard Levi subset}
Then the map $\Pa\fle\Pa^0$ is a bijection from the set of standard parabolic
subsets of $\G$ onto the set of $\theta$-stable standard parabolic subgroups of
$\G^0$. If $\Pa$ is a standard parabolic subset of $\G$, then $\til{\Pa}=
\Pa^0\rtimes\langle\theta\rangle$, $\Pa=\Pa^0\theta$, $\til{\M}_P=\M^0_P
\rtimes\langle\theta\rangle$
and $\M_P=\M^0_P\theta$. It is easy to see that the centralizer of $\M_P$ in
$\M^0_P$ is $Z(\M^0_P)^\theta$; so $\A_{M_P}$ is the maximal split subtorus
of $Z(\M_P^0)^\theta$. 

\begin{example}\label{ex:groupes_non_connexes}
Let $\H$ be a connected reductive quasi-split group over $K$ and $E/K$ be
a cyclic extension. Let $\G^0=R_{E/K}\H_E$, $\theta$ be the isomorphism of
$\G^0$ induced by a fixed
generator of $\Gal(E/K)$, $\til{\G}=\G^0\rtimes\langle
\theta\rangle$ and $\G=\G^0\rtimes\theta$.
Fix a Borel subgroup $\B_{H}$ of $\H$ and a Levi subrgoup $\T_H$ of $\B_H$.
Then $\B^0:=R_{E/K}\B_{H,E}$ is a $\theta$-stable Borel subgroup of $\G^0$,
and $\T^0:=R_{E/K}\T_{H,E}$ is a $\theta$-stable maximal torus of $\G^0$.
The standard parabolic subsets of $\G$ are in bijection with the
$\theta$-stable standard parabolic subgroups of $\G^0$, ie with the
standard parabolic subgroups of $\H$. If $\Pa$ corresponds to $\Pa_H$, then
$\Pa_H=\Pa^0\cap\H$, $\Pa^0=R_{E/K}\Pa_{H,K}$ and $\A_{M_P}=\A_{M_{P_H}}$.

Assume that $K$ is local or global. Then we can associate to $\H$ an
endoscopic datum $(\H^*,\mathcal{H},s,\xi)$ for $(\G^0,\theta,1)$ in the sense
of \cite{KS} 2.1. If $r=[E:K]$, then $\widehat{\G}^0\simeq\widehat{\H}^r$,
with $\widehat{\theta}(x_1,\dots,x_r)=(x_2,\dots,x_r,x_1)$. The diagonal
embedding $\widehat{\H}\fl\widehat{\G}^0$ is $W_K$-equivariant, hence
extends in an obvious way to a $L$-morphism $\xi:\mathcal{H}:={}^L\H\fl{}^L
\G^0$. Finally, take $s=1$.

\end{example}

\vspace{1cm}

Assume that we are in the situation of the example above. In \cite{La-CSCB}
2.4, Labesse defines the norm $\Norme\gamma$ of a $\theta$-semi-simple
element $\gamma$ of $\G^0(K)$ (and shows that it exists); $\Norme\gamma$ is
a stable conjugacy class in $\H(K)$ that depends only on the stable
$\theta$-conjugacy class of $\gamma$, and every element of $\Norme\gamma$ is
stably conjugate to $N\gamma:=\gamma\theta(\gamma)\dots\theta^{[E:K]-1}(\gamma)
\in\G^0(K)=\H(E)$. 
\index{Ngamma@$\Norme\gamma$\quad norm of $\gamma$}

If $\M$ is a Levi subset of $\G$, write, for every $\theta$-semi-simple
$\gamma\in\M^0(K)$,
\[D_M^G(\gamma)=\det(1-\Ad(\gamma)\circ\theta,Lie(\G^0)/Lie(\M^0)).\]
\index{DMG@$D_M^G$\quad (twisted case)}
If $\M$ is a standard Levi subset of $\G$ (or, more generally, any Levi subset
of $\G$ such that $\theta\in\M$), set $\M_H=(\M^0)^\theta=\M^0\cap\H$; then
$\M_H$ is a Levi subgroup of $\H$.

\begin{lemma}\label{lemme:egalite_D_M_G}
Let $\M$ be a standard Levi subset of $\G$. Then, for every
$\theta$-semi-simple $\gamma\in\M^0(K)$ :
\[D_M^G(\gamma)=D_{M_H}^{H}(\Norme\gamma).\]

\end{lemma}

\vspace{1cm}

Assume now that $K=\R$ and $E=\C$.
We will recall results of Clozel and Delorme about $\theta$-stable tempered
representations of $\G^0(\R)$.

Remember that an admissible representation $\pi$ of $\G^0(\R)$ is called
$\theta$-stable if $\pi\simeq\pi\circ\theta$.
\index{$\theta$-stable representation}
In that case, there exists an
intertwining operator $A_\pi:\pi\iso\pi\circ\theta$. We say that $A_\pi$ is
\emph{normalized} if $A_\pi^2=1$.
\index{intertwining operator}
\index{normalized intertwining operator}
The data of a normalized intertwining operator on $\pi$ is equivalent to
that of a representation of $\til{\G}(\R)$ extending $\pi$.
If $\pi$ is irreducible and $\theta$-stable,
then, by a Schur's theorem, it always has a normalized intertwining
operator.

For $\xi$ a quasi-character of $\A_G(\R)^0$, let $C_c^\infty(\G^0(\R),\xi)$
be the set of functions $f\in C^\infty(\G^0(\R))$ that have compact support
modulo $\A_G(\R)^0$ and such that $f(zg)=\xi(z)f(g)$
for every $(z,g)\in\A_G(\R)^0\times\G^0(\R)$.

The following theorem is due to Clozel (cf \cite{Cl-CBRT} 4.1, 5.12, 8.4).

\begin{theorem}\label{th:relevement_caracteres}
Let $\pi$ be an irreducible admissible $\theta$-stable representation of
$\G^0(\R)$ and $A_\pi$ be a normalized intertwining operator on $\pi$.
Let $\xi$ be the quasi-character through which $\A_G(\R)^0$ acts on the space
of $\pi$. Then the map
\[C_c^\infty(\G^0(\R),\xi^{-1})\fl\C,\quad f\fle\Tr(\pi(f)A_\pi)\]
extends to a distribution on $\G^0(\R)$ that is invariant under
$\theta$-conjugacy; this distribution is tempered if $\pi$ is tempered.
Call this distribution the \emph{twisted character} of
$\pi$ and denote it by $\Theta_\pi$.
\index{twisted character}

Let $\varphi:W_\R\fl{}^L\H$ be a tempered Langlands parameter; it defines a
$L$-packet $\Pi_H$ of tempered representations of $\H(\R)$. Write
$\Theta_{\Pi_H}=\sum\limits_{\pi_H\in\Pi_H}\Theta_{\pi_H}$, where, for every
$\pi_H\in\Pi_H$, $\Theta_{\pi_H}$ is the character of $\pi_H$.
Then the representation $\pi$ of $\G^0(\R)=\H(\C)$ associated to $\varphi_
{|W_\C}$ is tempered and $\theta$-stable, and, if $A_\pi$ is a normalized
intertwining operator on $\pi$, there exists $\varepsilon\in\{\pm 1\}$
such that, for every $\theta$-regular
$g\in\G^0(\R)$ :
\[\Theta_\pi(g)=\varepsilon\Theta_{\Pi_H}(\Norme g).\]

\end{theorem}

In particular, $\Theta_\pi$ is invariant under stable $\theta$-conjugacy.

\begin{remark}\label{rq:theta_stable}
Let $\pi$ be an irreducible tempered $\theta$-stable representation of
$\G^0(\R)$. If the infinitesimal character of $\pi$ is equal to that of
a finite-dimensional $\theta$-stable representation of $\G^0(\R)$, then
there exists a tempered Langlands parameter $\varphi:W_\R\fl{}^L\H$ such that
$\pi$ is associated to the parameter $\varphi_{|W_\C}$ (cf \cite{J} (5.16)).

\end{remark}

Assume from now on that $\H(\R)$ has a discrete series.
Let $\K'_\infty$ be the set of fixed points of a Cartan involution of
$\G^0(\R)$ that commutes with $\theta$. Write $\ggoth=Lie(\G)(\C)$. For
every admissible $\theta$-stable representation $\rho$ of $\G^0(\R)$, let
\[ep(\theta,\rho):=\sum_{i\geq 0}(-1)^i\Tr(\theta,\H^i(\ggoth,\K'_\infty;
\rho))\]
\index{eptheta@$ep(\theta,.)$\quad twisted Euler-Poincaré characteristic}
be the twisted Euler-Poincare characteristic of $\rho$. It depends on the
choice of a normalized intertwining operator on $\rho$.
An admissible representation of $\G^0(\R)$ is called \emph{$\theta$-discrete}
if it is irreducible tempered $\theta$-stable and is not a subquotient of
a representation induced from an admissible $\theta$-stable representation
of a proper $\theta$-stable Levi subgroup (cf \cite{AC} I.2.3).
\index{$\theta$-discrete representation}

The following theorem is due to Labesse (cf \cite{La-PCTC} proposition 12).

\begin{theorem}\label{th:pseudo_coeff_tordus}
Let $\pi$ be a $\theta$-discrete representation of $\G^0(\R)$, and let $\xi$
be the quasi-character through which $\A_G(\R)^0$ acts on the space of $\pi$.
Assume that $\pi$ is associated to a Langlands parameter $\varphi_\pi:W_\C\fl{}
^L\H$ satisfying $\varphi_\pi=\varphi_{|W_\C}$, where $\varphi:W_\R\fl{}^L\H$
is a Langlands parameter of the $L$-packet of the discrete series of $\H(\R)$
associated to the contragredient of an irreducible algebraic representation
$V$ of $\H$. As in section 3 of \cite{Cl-CBRT}, we associate to $V$ a
$\theta$-stable algebraic representation $W$
of $\G^0$ and a normalized intertwining operator $A_W$ on
$W$.

Then there exists a function $\phi\in C_c^\infty(\G^0(\R),\xi^{-1})$,
$\K'_\infty$-finite on the right and on the left modulo $\A_G(\R)^0$, such
that, for every admissible $\theta$-stable representation $\rho$ of $\G^0(\R)$
that is of finite length and such that $\rho_{|\A_G(\R)^0}=\xi$
and every normalized intertwining operator
$A_\rho$ on $\rho$,
\[\Tr(\rho(\phi)A_\rho)=ep(\theta,\rho\otimes W).\]

\end{theorem}

Such a function $\phi$ is called \emph{twisted pseudo-coefficient} of $\pi$.
This name is justified by the next remark.
\index{twisted pseudo-coefficient}

\begin{remark}\label{rq:ep_tordu} Let $\pi$ and $\phi$ be as in the above
theorem, and let $A_\pi$ be a normalized intertwining operator on $\pi$.
By the proof of proposition 3.6 of \cite{Cl-RGRA} and by theorem 2 (p 217)
of \cite{De},
\[\Tr(\pi(\phi)A_\pi)\not=0\]
and, for every irreducible $\theta$-stable tempered representation $\rho$ of
$\G^0(\R)$ and every normalized intertwining operator $A_\rho$ on $\rho$,
\[\Tr(\rho(\phi)A_\rho)=0\]
if $\rho\not\simeq\pi$.
In particular, the function $\phi$ is cuspidal (the definition of ``cuspidal''
is recalled, for example, at the beginning of section 7 of \cite{A-ITF2}).

\end{remark}

\begin{definition}\label{def:d_k_tordu} 
Let $\T_e$ be a torus of $\G^0_\R$ such that $\T_e(\R)$ is a maximal torus of
$\K'_\infty$. Set
\[d(\G)=|\Ker(\H^1(\R,\T_e)\fl\H^1(\R,\G^0))|\]
\[k(\G)=|Im(\H^1(\R,(\T_e)_{sc})\fl\H^1(\R,\T_e))|,\]
with $(\T_e)_{sc}$ the inverse image of $\T_e$ by the morphism $\G^0_{sc}\fl
\G^0$ (where $\G^0_{sc}\fl\G^0_{der}$ is the simply connected covering of
$\G^0_{der}$).
\footnote{Cf \cite{CL} A.1 for the definition of $d(\G)$; Clozel and
Labesse use a maximal $\R$-anisotropic torus instead of a maximal
$\R$-elliptic torus, but this does not give the correct result if
$\A_G\not=\{1\}$ (cf \cite{K-NP} 1.1 for the case of connected groups).}
\index{dG@$d(\G)$\quad (twisted case)}
\index{kG@$k(\G)$\quad (twisted case)}

\end{definition}

\begin{remark}\label{rq:d_k_tordu}
As $\G^0$ comes from a complex group by restriction of scalars,
$\H^1(\R,\G^0)=\{1\}$, so $d(\G)=|\H^1(\R,\T_e)|$.
For example, if $\H=\G(\U^*(n_1)\times\dots\times\U^*(n_r))$ (cf \ref{groupes1}
for the definition of this group) with $n:=n_1+\dots+n_r\geq 1$, then
$\T_e=\G(\U(1)^n)$, so $d(\G)=2^{n-1}$.

On the other hand, if the derived group of $\G^0$ is simply connected, then
$k(\G)=|Im(\H^1(\R,\T_e\cap\G_{der})\fl\H^1(\R,\T_e))|$.

\end{remark}

Remember that a $\theta$-semi-simple element $g$ of $\G^0(\R)$ is called
$\theta$-elliptic if $\A_{G^0_{g\theta}}=\A_G(=\A_H)$.

\begin{lemma}\label{lemme:theta_elliptique} Let $g\in\G^0(\R)$ be
$\theta$-semi-simple. Then $g$ is $\theta$-elliptic if and only if $\Norme g$
is elliptic. Moreover, if $g$ is not $\theta$-elliptic, then there exists a
proper Levi subset $\M$ of $\G$ such that $g\theta\in\M(\R)$ and $\G^0_
{g\theta}\subset\M^0$.

\end{lemma}

\begin{proof} Let $g\in\G^0(\R)$ be $\theta$-semi-simple. 
As $h:=g\theta(g)$ is $\G^0(\R)$-conjugate to an element of $\H(\R)$
(cf \cite{Cl-CBRT} p 55), we may assume, after replacing $g$ by a
$\theta$-conjugate, that $h\in\H(\R)$. Let $\Le=\G^0_h$. Then $\Le$ is
stable by the morphism (of algebraic groups over $\R$)
$\theta':x\fle\theta(g)\theta(x)\theta(g)^{-1}$, and it is easy to see that
$\theta'_{|\Le}$ is an involution and that $\G^0_{g\theta}=\Le^{\theta'}$
and $\H_h=\Le^\theta$. This implies that $Z(\H_h)=Z(\Le)^\theta$ and that
$Z(\G^0_{g\theta})=Z(\Le)^{\theta'}$. But $\theta(g)\in\Le(\R)$ (we assumed
that $g\theta(g)\in\H(\R)$, so $g$ and $\theta(g)$ commute), so
$\theta_{|Z(\Le)}=\theta'_{|Z(\Le)}$. Hence $\A_{G^0_{g\theta}}=\A_{H_h}$,
and this proves that $g$ is $\theta$-elliptic if and only if $h$ is elliptic.

Suppose that $g$ is not $\theta$-elliptic. Then $h$ is not elliptic,
so there exists a proper Levi subgroup $\M_H$ of $\H$ such that
$\H_h\subset\M_H$. Let $\M^0=R_{\C/\R}\M_{H,\C}$; it is a $\theta$-stable
proper Levi subgroup of $\G^0$. Let $\M$ be the Levi subset of $\G$
associated to $\M^0$. As $\G^0_h\subset\M^0$, $\G^0_{g\theta}\subset
\M^0$; moreover, $g\in\G^0_h(\R)$, so $g\theta\in\M(\R)$.

\end{proof}

Note that, by the above proof, for every $\theta$-semi-simple $g\in\G^0(\R)$,
if $h\in\Norme(g)$, then $\G^0_{g\theta}$ is an inner form of $\H_h$. Later,
when we calculate orbital integrals at $g\theta$ and $h$, we always choose
Haar measures on $\G^0_{g\theta}$ and $\H_h$ that correspond to each other.

Finally, we calculate the twisted orbital integrals of some of the twisted
pseudo-coefficients defined above at $\theta$-semi-simple elements.
To avoid technical complications, assume that $\G^{der}$ is simply connected.
If $\phi\in C_c^\infty(\G^0(\R),\xi^{-1})$ ($\xi$ is a quasi-character on
$\A_G(\R)^0$) and $g\in\G^0(\R)$, the twisted orbital integral of $\phi$ at $g$
(also called orbital integral of $\phi$ at $g\theta$) is by definition
\[O_{g\theta}(\phi)=\int_{\G^0_{g\theta}(\R)\sous\G^0(\R)}\phi(x^{-1}g
\theta(x))dx\]
(of course, it depends of the choice of Haar measures on $\G^0(\R)$ and
$\G^0_{g\theta}(\R)$).
\index{twisted orbital integral}
\index{Ogtheta@$O_{g\theta}$\quad twisted orbital integral}

Let $\G'$ be a reductive connected algebraic group over $\R$. If $\G'$ has an
inner form $\overline{\G}'$ that is anisotropic modulo its center, set
\[v(\G')=(-1)^{q(\G')}\vol(\overline{\G}'(\R)/\A_{G'_\R}(\R)^0)d(\G')
^{-1},\]
where $d(\G')$ is defined in \ref{serie_discrete1}.
\index{vG@$v(\G)$}

\begin{lemma}\label{lemme:O_pseudo_coeff} Let $V$ be an irreducible algebraic
representation of $\H$, $\varphi:W_\R\fl{}^L\H$ be a Langlands parameter of
the $L$-packet of the discrete series of $\H(\R)$ associated to $V^*$, 
$\pi_{V^*}$ be the representation of $\H(\C)=\G^0(\R)$ corresponding to
$\varphi_{|W_\C}$ (so $\pi_{V^*}$ is $\theta$-discrete) and $\phi_{V^*}$
be a twisted pseudo-coefficient of $\pi_{V^*}$.
Let $g\in\G^0(\R)$ be $\theta$-semi-simple. Then
\[O_{g\theta}(\phi_{V^*})=v(\G^0_{g\theta})^{-1}\Theta_{\pi_{V^*}^\vee}(g)\]
if $g$ is $\theta$-elliptic, and
\[O_{g\theta}(\phi_{V^*})=0\]
if $g$ is not $\theta$-elliptic.

\end{lemma}
\index{$\pi_V$}
\index{$\phi_V$}

The proof of this lemma is inspired by the proof of theorem 2.12 of \cite{CCl}.

\begin{proof} To simplify the notation, we will write $\pi=\pi_{V^*}$
and $\phi=\phi_{V^*}$.

If $V$ is the trivial representation, then the twisted orbital
integrals of $\phi$ are calculated in theorem A.1.1 of \cite{CL};
write $\phi_0=\phi$. (Note that Clozel and Labesse choose the Haar measure on
$\G^0_{g\theta}(\R)$ for which $\vol(\overline{\G}'(\R)/\A_{G'}(\R)^0)=1$,
where $\overline{\G}'$ is an inner form of $\G^0_{g\theta}$ that is
anisotropic modulo its center).

Assume that $V$ is any irreducible algebraic representation of $\H$. Let
$W$ be the $\theta$-stable algebraic representation of $\G^0$
associated to $V$ as in theorem \ref{th:pseudo_coeff_tordus}, with the
normalized intertwining operator $A_\rho$ fixed in that theorem.
Write $\phi'=\Theta_{W}\phi_0$. As $\Theta_{W}$ is invariant by
$\theta$-conjugacy, proposition 3.4 of \cite{Cl-CBRT} implies that, for every
$\theta$-semi-simple $g\in\G^0(\R)$, $O_{g\theta}(\phi')=O_{g\theta}(\phi_0)
\Theta_{\pi^\vee}(g)$. So it is enough to show that $\phi$ and $\phi'$
have the same orbital integrals. By theorem 1 of \cite{KRo}, in order to
prove this, it suffices to show that, for every $\theta$-stable tempered
representation $\rho$ of $\G^0(\R)$ and every normalized intertwining
operator $A_\rho$ on $\rho$,
\[\Tr(\rho(\phi)A_\rho)=\Tr(\rho(\phi')A_\rho).\]
Fix such a representation $\rho$. Then it is easy to see that
\[\Tr(\rho(\phi')A_\rho)=\Tr((\rho\otimes W)(\phi_0)(A_{\rho}\otimes
A_{W})).\]
Hence
\[\Tr(\rho(\phi')A_\rho)=ep(\theta,\rho\otimes W)=\Tr(\rho(\phi)
A_\rho).\]

\end{proof}

\begin{subcorollaire}\label{cor:pseudo_coeff} Use the notations of lemma
\ref{lemme:O_pseudo_coeff} above.
\begin{itemize}
\item[(i)] The function $\phi_{V^*}$ is stabilizing (``stabilisante'')
in the sense of \cite{La-CSCB} 3.8.2.
\item[(ii)] Let $f_{V^*}=\frac{1}{|\Pi(\varphi)|}\sum\limits_{\pi_H\in\Pi
(\varphi)}f_{\pi_H}$,
where $\Pi(\varphi)$ is the discrete series $L$-packet of $\H(\R)$
associated to $\varphi:W_\R\fl{}^L\H$ and,
for every representation $\pi_H$ in the discrete series of $\H(\R)$,
$f_{\pi_H}$ is a pseudo-coefficient of $\pi_H$. Then the functions
$\phi_{V^*}$ and $d(\G)f_{V^*}$ are associated in the sense of
\cite{La} 3.2.
\index{fV@$f_V$}

\end{itemize}
\end{subcorollaire}

\begin{proof} The result follows from lemma \ref{lemme:O_pseudo_coeff}
and the proof of theorem A.1.1 of \cite{CL} (and lemma A.1.2 of
\cite{CL}).

\end{proof}

\section{The invariant trace formula}
\label{GL_n_applications2}

\renewcommand{\theequation}{$*$}

Note first that, thanks to the work of Delorme-Mezo (\cite{DeM}) and
Kottwitz-Rogawski (\cite{KRo}), Arthur's invariant trace formula
(see, e.g., \cite{A-ITF2}) is now available for non-connected groups
as well as for connected groups.

In \cite{A-L2}, Arthur gave a simple form of the invariant trace formula
(on a connected group) for a function that it stable cuspidal at
infinity (this notion is defined at the beginning of section 4 of
\cite{A-L2} and recalled in \ref{FT_stable_geometrique4}).
The goal of this section is to give a similar formula for a (very) particular
class of non-connected groups.

Let $\H$ be a connected reductive quasi-split group over $\Q$; fix a Borel
subgroup of $\H$ and a Levi subgroup of this Borel.
Fix an
imaginary quadratic extension $E$ of $\Q$, and take
$\G^0=R_{E/\Q}\H_E$. Assume that the derived group of $\H$ is simply connected
and that $\H$ is cuspidal (in the sense of definition
\ref{def:groupe_cuspidal}).
A Levi subset $\M$ of $\G$ is called cuspidal if it is conjugate
to a standard Levi subset $\M'$ such that $\M'_H$ is cuspidal.
\index{cuspidal Levi subset}

We first define the analogs of the functions
$\Phi_M^G(.,\Theta)$ of \ref{serie_discrete2}. 

\begin{sublemme}\label{pro:Phi_M_tordu} Let $\pi$ be a $\theta$-stable
irreducible tempered representation of $\G^0(\R)$.
Fix a normalized intertwining operator $A_\pi$ on $\pi$. Assume that there
exists a Langlands parameter $\varphi:W_\R\fl{}^L\H$ such that $\pi$ is
associated to $\varphi_{|W_\C}$. Let $\M$ be a standard cuspidal Levi subset
of $\G$, and let $\T_H$ be a maximal torus of $\M_{H,\R}$ that is anisotropic
modulo $\A_{M_H}$. Write $D$ for the set of $\gamma\in\G^0(\R)$ such that
$\gamma\theta(\gamma)\in\T_H(\R)$.
Then the function
\[D\cap\G^0(\R)_{\theta-reg}\fl\C,\quad
\gamma\fle|D_M^G(\gamma)|^{1/2}\Theta_\pi(\gamma)\]
extends to a continuous function $D\fl\C$, that will be denoted by
$\Phi_M^G(.,\Theta_\pi)$.
\index{$\Phi_M^G$\quad (twisted case)}

\end{sublemme}

Extend $\Phi_M^G(.,\Theta_\pi)$ to a function on $\M^0(\R)$ in the following
way : if $\gamma\in\M^0(\R)$ is $\theta$-elliptic (in $\M^0(\R)$), then it is
$\theta$-conjugate to $\gamma'\in D$, and we set $\Phi_M^G(\gamma,
\Theta_\pi)=\Phi_M^G(\gamma',\Theta_\pi)$; otherwise, we set
$\Phi_M^G(\gamma,\Theta_\pi)=0$.
The function $\Phi_M^G(.,\Theta_\pi)$ is clearly invariant by stable
$\theta$-conjugacy. As every Levi subset of $\G$ is $\G^0(\R)$-conjugate to
a standard Levi subset, we can define in the same way a function
$\Phi_M^G(.,\Theta_\pi)$ for any cuspidal Levi subset $\M$.

The lemma above follows from the similar lemma for connected groups
(lemma \ref{lemme:Phi_M}, due to Arthur and Shelstad), from theorem
\ref{th:relevement_caracteres} (due to Clozel) and from lemma
\ref{lemme:egalite_D_M_G}.

Let $\Pi_{\theta-disc}(\G^0(\R))$ be the set of isomorphism classes of
$\theta$-discrete representations of $\G^0(\R)$. For every
$\pi\in\Pi_{\theta-disc}(\G^0(\R))$, fix a normalizing operator
$A_\pi$ on $\pi$.
\index{$\Pi_{\theta-disc}(\G^0(\R))$}

\begin{subdefinition}\label{def:Phi_M_tordu} Let $\M$ be a cuspidal Levi
subset of $\G$. Let $\xi$ be a $\theta$-stable quasi-character of $\A_{G^0}
(\R)^0$. For every function $\phi\in C_c^\infty(\G^0(\R),\xi^{-1})$ that is
left and right $\K_\infty$-finite modulo $\A_{G^0}(\R)^0$ and every
$\gamma\in\M^0(\R)$, write :
\[\Phi_M^G(\gamma,\phi)=(-1)^{\dim(\A_M/\A_G)}v(\M^0_{\gamma\theta})^{-1}
\sum_{\pi\in\Pi_{\theta-disc}(\G^0(\R))}\Phi_M^G(\gamma,\Theta_{\pi^{\vee}})
\Tr(\pi(\phi)A_\pi),\]
and :
\[S\Phi_M^G(\gamma,\phi)=(-1)^{\dim(\A_M/\A_G)}k(\M)k(\G)^{-1}\overline{v}
(\M^0_{\gamma\theta})^{-1}\sum_{\pi\in\Pi_{\theta-disc}(\G^0(\R))}\Phi_M^G
(\gamma,\Theta_{\pi^{\vee}})\Tr(\pi(\phi)A_\pi).\]
\index{SPhiMG@$S\Phi_M^G$\quad (twisted case)}

\end{subdefinition}

The notations $k$ and $\overline{v}$ are defined in
\ref{FT_stable_geometrique4}, and the notation $v$ is that of
\ref{GL_n_applications1}.

Let $\M_0$ be the minimal $\theta$-stable Levi subgroup of
$\G^0$ corresponding to the fixed minimal Levi subgroup of $\H$ 
($\M_0$ is a torus because $\H$ is quasi-split) and $\xi$ be a $\theta$-stable
quasi-character of $\A_{G^0}(\R)^0$.
Define an action of the group $\til{\G}(\Ade)$ on 
$L^2(\G^0(\Q)\sous\G^0(\Ade),\xi)$ in the following way : 
the subgroup $\G^0(\Ade)$ acts in the usual way, and $\theta$ acts by
$\phi\fle\phi\circ\theta$. For every irreducible $\theta$-stable representation
$\pi$ of $\G^0(\Ade)$ such that $m_{disc}(\pi)\not=0$ (ie such that $\pi$ is
a direct factor of $L^2(\G^0(\Q)\sous\G^0(\Ade),\xi)$, seen as a representation
of $\G^0(\Ade)$), fix a normalized intertwining operator $A_\pi$ on $\pi$.
If $\pi$ and $\pi'$ are such that $\pi_\infty\simeq\pi'_\infty$, choose
intertwining operators that are compatible at infinity.
For the $\theta$-stable irreducible admissible representations of
$\G^0(\R)$ that don't appear in this way, use any normalised intertwining
operator.
If $\pi$ is as above, let
$\til{\pi}^+$ (resp. $\til{\pi}^-$) be the representation of $\til{\G}(\Ade)$
defined by $\pi$ and $A_\pi$ (resp. $-A_\pi$), and let
$m_{disc}^+(\pi)$ (resp. $m_{disc}^-(\pi)$) be the multiplicity of
$\til{\pi}^+$ (resp. $\til{\pi}^-$) in
$L^2(\G^0(\Q)\sous\G^0(\Ade),\xi)$.

We write $C_c^\infty(\G^0(\Ade),\xi^{-1})$ for the vector space of functions
$\phi:\G^0(\Ade)\fl\C$ that are finite linear combinations of functions of the
form $\phi^\infty\otimes\phi_\infty$, with $\phi^\infty\in C_c^\infty(\G^0(
\Af))$ and $\phi_\infty\in C_c^\infty(\G^0(\R),\xi^{-1})$.

Let $\M$ be a Levi subset of $\G$. If $\M$ is cuspidal, then, for every
function $\phi=\phi^\infty\otimes\phi_\infty\in C^\infty(\G(\Ade),\xi^{-1})$,
write
\[T_{M,geom}^G(\phi)=\sum_\gamma
\vol(\M^0_{\gamma\theta}(\Q)\A_M(\R)^0\sous\M^0_{\gamma\theta}(\Ade))
O_{\gamma\theta}(\phi^\infty_{M})\Phi_M^G(\gamma,
\phi_\infty),\]
where the sum is taken over the set of $\theta$-conjugacy classes of
$\theta$-semi-simple elements of $\M^0(\Q)$ and
$\phi^\infty_M$ is the constant term of $\phi^\infty$ at $\M$
(defined in exactly the same way as in the case of connected groups).
\index{TMgeom@$T_{M,geom}^G$}

If $\M$ is not cuspidal, set $T_{M,geom}^G=0$.

For every $t\geq 0$, define $\Pi_{disc}(\G,t)$ and the function
$a_{disc}=a_{disc}^G:\Pi_{disc}(\G,t)\fl\C$ as in section 4 of \cite{A-ITF2}
(p 515-517).
\index{$\Pi_{disc}(\G,t)$}
\index{adisc@$a_{disc}^G(\pi)$}

Let $T^G$ be the distribution of the $\theta$-twisted invariant trace formula
on $\G^0(\Ade)$.
\index{TG@$T^G$\quad twisted invariant trace formula}
The following proposition is the analog of theorem 6.1 of \cite{A-L2} (and of
the formula below (3.5) of this article).

\begin{subproposition}\label{prop:FT_invariante_tordue} 
\footnote{I thank Sug Woo Shin for pointing out that I had forgotten terms on
the spectral side of this proposition.}
Let $\phi=\phi^\infty\phi_\infty\in C^\infty_c(\G^0(\Ade),\xi^{-1})$.
Assume that there exists an irreducible algebraic representation
$V$ of $\H$ such that, if $\varphi:W_\R\fl{}^L\H$ is the Langlands
parameter of the discrete series $L$-packet $\Pi_V$ of $\H(\R)$ associated
to $V$ and $\pi_\infty$ is the $\theta$-discrete representation of
$\G^0(\R)$ with Langlands parameter $\varphi_{|W_\C}$, then $\phi_\infty$
is a twisted pseudo-coefficient of $\pi_\infty$. Then
\[T^G(\phi)=\sum_{\M}(n_M^G)^{-1}T_{M,geom}^G(\phi)=\sum_{t\geq 0}\sum_
{\pi\in\Pi_{disc}(\G,t)}a_{disc}(\pi)\Tr(\pi(\phi)A_\pi),\]
where the first sum is taken over the set of
$\G^0(\Q)$-conjugacy classes of Levi subsets $\M$ of $\G$.

\end{subproposition}

\begin{subremarques}\begin{itemize}
\item[(1)] If $\pi$ is a cuspidal $\theta$-stable representation of
$\G^0(\Ade)$, then $a_{disc}(\pi)=m^+_{disc}(\pi)+m^-_{disc}(\pi)$.
(This is an easy consequence of the definition of $a_{disc}$, cf
\cite{A-ITF2} (4.3) and (4.4).)
\item[(2)] The spectral side of the formula of \cite{A-L2} (formula above
(3.5)) is simpler, because only the discrete automorphic representations
of $\G^0(\Ade)$ can contribute. However, in the twisted case, it is not
possible to eliminate the contributions from the discrete spectrum of
proper Levi subsets (because there might be representations of $\M(\Ade)$
that are fixed by a regular element of $W_0^G$ and still have a regular
archimedean infinitesimal character).

\item[(3)] In theorem 3.3 of \cite{A-ITF2}, the sum is taken over all
Levi subsets $\M$ of $\G$ such that $\M^0$ contains $\M_0$, and the
coefficients are $|W_0^M||W_0^G|^{-1}$ instead of $(n_M^G)^{-1}$; it is easy
to see that these are just two ways to write the same thing.

\end{itemize}
\end{subremarques}

\begin{proof} The second formula (ie the spectral side) is just (a) of
theorem 7.1 of \cite{A-ITF2}, because $\phi$ is cuspidal at infinity.

We show the first formula.
We have to compute the value at $\phi$ of the invariant distributions $I_M^G$
of \cite{A-ITF1}.
\index{IMG@$I_M^G$}
As $\phi_\infty$ is cuspidal, we see by using the
splitting formula (proposition 9.1 of \cite{A-ITF1}) as in \cite{A-L2} \S 3 
that it is enough to compute the $I_M^G$ at infinity, ie to prove the analog of
theorem 5.1 of \cite{A-L2}.
Moreover, by corollary 9.2 of \cite{A-ITF1} applied to the set of places
$S=\{\infty\}$, and thanks to the cuspidality of
$\phi_\infty$, we see that the term corresponding to $\M$ is non-zero only if
$\A_M=\A_{M_\R}$.
\footnote{I thank Robert Kottwitz for patiently explaining to me this subtlety
of the trace formula.}

So we may assume that $\A_M=\A_{M_\R}$.
We want to show that $I_M^G(.,\phi_\infty)=0$ if $\M$ is not cuspidal and
that, for every cuspidal Levi subset $\M$ of $\G$ and for every
$\gamma\in\M^0(\R)$ :
\begin{equation}
I_M^G(\gamma,\phi_\infty)=|D^M(\gamma)|^{1/2}\Phi_M^G(\gamma,
\phi_\infty),\end{equation}
where, if $\gamma\theta=(\sigma\theta)u$ is the Jordan decomposition of
$\gamma\theta$, then
\[D^M(\gamma)=\det((1-\Ad(\sigma)\circ\theta),Lie(\M^0)/Lie(\M^0_{\sigma
\theta})).\]
(Note that, if $\M$ is cuspidal, then $\A_M=\A_{M_\R}$.)
This implies in particular that $I_M^G(\gamma,\phi_\infty)=0$ if $\gamma$
is not $\theta$-semi-simple.
The rest of the proof of theorem 6.1 of \cite{A-L2} applies without any
changes to the case of non-connected groups.

The case where $\M$ is not cuspidal is treated in lemma
\ref{lemme:FT_invariante_tordue2}. In the rest of this proof, we assume
that $\M$ is cuspidal.

For connected groups, the analog of formula $(*)$ for a
semi-simple regular $\gamma$ is theorem 6.4 of \cite{A-IOR2} (cf
formula (4.1) of \cite{A-L2}).
Arthur shows in section 5 of \cite{A-L2} that the analog of $(*)$
for any $\gamma$ is a conseqeunce of this case.

We first show that formula $(*)$ for a $\theta$-semi-simple
$\theta$-regular $\gamma$ implies formula $(*)$ for any $\gamma$, by adapting
the reasoning in section 5 of \cite{A-L2}. 
The reasoning in the second half of page 277 of \cite{A-L2} applies to the case
considered here and shows that it is enough to prove $(*)$ for a
$\gamma\in\M^0(\R)$ such that $\M^0_{\gamma\theta}=\G^0_{\gamma\theta}$. 
Lemma \ref{lemme:FT_invariante_tordue1} below is the analog of lemma 5.3 of
\cite{A-L2}. Once this lemma is known, the rest of the reasoning of
\cite{A-L2} applies. This is because Arthur reduces to the semi-simple
regular case by using the results on orbital integrals at unipotent elements
of \cite{A-L2} p 275-277, and we can apply the same results here, because
these orbitals integrals are taken on the connected group
$\M^0_{\sigma\theta}$ (where, as before, $\gamma\theta=(\sigma\theta)u$
is the Jordan decomposition of $\gamma\theta$).

It remains to show formula $(*)$ for a $\theta$-semi-simple $\theta$-regular
$\gamma$. The article \cite{A-IOR2} is written in the setting of connected
groups, but it is easy to check that, now the invariant formula for
non-connected groups is known, all the article until, and including,
corollary 6.3, applies to the general (not necessarily connected) case.
We can write statements analogous to theorem 6.4 and lemma 6.6 of
\cite{A-IOR2}, by making the following changes :
take a $\theta$-discrete representation $\pi_\infty$ of $\G^0(\R)$
(instead of a discrete series representation of $\G(\R)$), and replace
the character of $\pi_\infty$ by the twisted character.

The proof of lemma 6.6 of \cite{A-IOR2} applies to the non-connected case,
if we replace $\Pi_{temp}$ by the set of isomorphism classes of
$\theta$-stable tempered representations, $\Pi_{disc}$ by 
$\Pi_{\theta-disc}$ and ``regular'' by ``$\theta$-regular''.

The proof of theorem 6.4 of \cite{A-IOR2} proceeds by induction on $\M$,
starting from the case $\M=\G$, and uses lemma 6.6 of \cite{A-IOR2} and
three properties of the characters of discrete series representations :
the differential equations that they satisfy, the conditions at the boundary
of the set of regular elements and the growth properties.
For the twisted characters of $\theta$-discrete representations, there
are of course similar differential equations; the bound that we need
(in the third property) is proved by Clozel in theorem 5.1 of \cite{Cl-CBRT};
as for the boundary conditions, they follow from theorem 7.2 of
\cite{Cl-CBRT} (called theorem \ref{th:relevement_caracteres} in this book)
and from the case of connected groups.
Once these results are known, the reduction to the case $\M=\G$ is
the same as in \cite{A-IOR2}. But the case $\M=\G$ is exactly lemma
\ref{lemme:O_pseudo_coeff}.

\end{proof}

\renewcommand{\theequation}{$**$}

\begin{sublemme}\label{lemme:FT_invariante_tordue1} 
Write $\Phi'_M=|D^M(.)|^{-1/2}I_M^G$. 

Let $\M$ be a cuspidal Levi subset of $\G$ and
$\gamma\in\M^0(\R)$ be such that $\G^0_{\gamma\theta}=\M^0_{\gamma\theta}$.
Let $\gamma\theta=(\sigma\theta)u$ be the Jordan decomposition of
$\gamma\theta$. Then there exist stable cuspidal functions
$f_1,\dots,f_n$ on $\M^0_{\sigma\theta}(\R)$ and a neighbourhood $U$ of $1$
in $\M^0_{\sigma\theta}(\R)$, invariant by $\M^0_{\sigma\theta}(\R)$-conjugacy,
such that, for every $\mu\in U$ :
\begin{equation}
\Phi'_M(\mu\sigma,\phi_\infty)=\sum_{i=1}^n\Phi_{\M^0_
{\sigma\theta}}(\mu,f_i).
\end{equation}

\end{sublemme}

\begin{proof} 
If $\sigma$ is not $\theta$-elliptic in $\M(\R)$, then, by lemma
\ref{lemme:theta_elliptique}, there exists a proper Levi subset
$\M_1$ of $\M_\R$ such that $\sigma\theta\in\M_1(\R)$ and 
$\M_{\R,\sigma\theta}^0\subset\M_1^0$. If $\mu\in\M^0_{\sigma\theta}(\R)$ is
small enough, then $\M^0_{(\sigma\theta)\mu}$ is also included in $\M_1^0$.
Applying the descent property (corollary 8.3 of \cite{A-ITF1}) and using the
cuspidality of $\phi_\infty$, we see that $\Phi'_M(\mu\sigma,\phi_\infty)=0$,
so that we can take $f_i=0$.

We may therefore assume that $\sigma$ is $\theta$-elliptic in $\M^0(\R)$.
We may also assume that $\M$ is standard. Let $\T_H$ be a maximal torus in
$\M_{H,\R}$ that is anisotropic modulo $\A_{M_H}$. Then $\sigma$ is 
$\theta$-conjugate to an element $\sigma'$ such that $\sigma'\theta(\sigma')\in
\T_H(\R)$. As $I_M^G$ is invariant by $\theta$-conjugacy, we may assume
that $h:=\sigma\theta(\sigma)\in\T_H(\R)$. As $\G^0_{\sigma\theta}$
is an inner form of $\H_h$ over $\R$, the maximal torus $\T_H$ of
$\H_h$ transfers to a maximal torus $\T$ of $\G^0_{\sigma\theta}$; of
course, $h\in\T(\R)$. Let $U$ be an invariant neighbourhood of $1$ in
$\M^0_{\sigma\theta}(\R)$ small enough so that $\mu\sigma$ is 
$\theta$-regular if $\mu\in U\cap\T(\R)$ is regular in $\M^0_{\sigma
\theta}$. Then, by formula $(*)$ in the proof of proposition
\ref{prop:FT_invariante_tordue} above for a $\theta$-regular element
(the proof of formula $(*)$ in this case does not depend on the lemma),
for every $\mu\in U\cap\T(\R)$ that is regular in $\M^0_{\sigma\theta}$ :
\[\Phi'_M(\mu\sigma,\phi_\infty)=\Phi_M^G(\mu\sigma,\phi_\infty)=(-1)^
{\dim(\A_M/\A_G)}\Phi_M^G(\mu\sigma,\Theta_{\pi_\infty^{\vee}}),\]
and we know that this is equal to
\[\pm|D_{M_H}^{H}(\Norme(\mu\sigma))|^{1/2}\Theta_{\Pi_{H}^{\vee}}(\Norme(
\mu\sigma))\]
(where the sign depends on the choice of normalized intertwining operator
on $\pi_\infty$).
By the proof of lemma 5.3 of \cite{A-L2} and lemma 4.1 of 
\cite{A-L2}, there exists $f_1,\dots,f_n$ stable cuspidal on $\M^0_{\sigma
\theta}(\R)$ such that $(**)$ is satisfied for every
$\mu\in U\cap\T(\R)$ that is regular in $\M^0_{\sigma\theta}$. 

It remains to show that, for this choice of $f_1,\dots,f_n$ and maybe after
making $U$ smaller, formula $(**)$ is true for every $\mu\in U$.
But the end of the proof of lemma 5.3 of \cite{A-L2} applies without any
changes to the non-connected case.

\end{proof}

\begin{sublemme}\label{lemme:FT_invariante_tordue2} Let $\M$ be
a Levi subset of $\G$. Assume that $\A_M=\A_{M_\R}$ and that $\M$ is not
cuspidal. Then $I_M^G(.,\phi_\infty)=0$.

\end{sublemme}

\begin{proof} We may assume that $\M$ is standard. We first show that
$I_M^G(\gamma,\phi_\infty)=0$ if $\gamma$ is $\theta$-regular in
$\M^0$. Let $\gamma\in\M^0(\R)$ be $\theta$-regular in $\M^0$. We may assume
that $\gamma\theta(\gamma)\in\M_H(\R)$. The centralizer $\T_H$ of
$\gamma\theta(\gamma)$ in $\M_H$ is a maximal torus of $\M_{H,\R}$. By the
assumption on $\M$, the torus $\T_H/\A_{M_\R}$ is not anisotropic, so
there exists a Levi subgroup $\M_{1,H}\not=\M_{H,\R}$ of $\M_{H,\R}$ such
that $\T_H\subset\M_{1,H}$ (for example the centralizer of the $\R$-split part
of $\T_H$). Let $\M_1$ be the corresponding Levi subset of $\G_\R$. Then
$\gamma\in\M_1^0(\R)$ and $\M^0_{1,\gamma\theta}=\M^0_{\gamma\theta}$.
By the descent formula (theorem 8.3 of \cite{A-ITF1}) and the cuspidality
of $\phi_\infty$, $I_M^G(\gamma,\phi_\infty)=0$.

We now show the statement of the lemma.
By formula (2.2) of \cite{A-ITF1}, it is enough to prove that,
for every Levi subset $\M'$ of $\G$ containing $\M$ and every
$\gamma\in\M^0(\R)$ such that $\M^0_{\gamma\theta}=\G^0_{\gamma\theta}$,
$I_{M'}^G(\gamma,\phi_\infty)=0$. If $\M'\not=\M$, this follows from the
descent formula (theorem 8.3 of \cite{A-ITF1}) and the cuspidality of
$\phi_\infty$. It remains to show that $I_M^G(\gamma,\phi_\infty)=0$,
if $\gamma\in\M^0(\R)$ is such that $\M^0_{\gamma\theta}=\G^0_{\gamma\theta}$.
Let $\gamma\theta=(\sigma\theta)u$ be the Jordan decomposition of
$\gamma\theta$. By (2.3) of \cite{A-ITF1}, there exists $f\in C_c^\infty
(\M^0(\R))$ and an open neighbourhood $U$ of $1$ in $\M^0_{\sigma\theta}
(\R)$ such that, for every $\mu\in U$, $I_M^G(\mu\sigma,\phi_\infty)=
O_{\mu\sigma\theta}(f)$. Hence, by the beginning of the proof,
$O_{\mu\sigma\theta}(f)=I_M^G(\mu\sigma,\phi_\infty)=0$ if $\mu\in U$ is such
that $\mu\sigma$ is $\theta$-regular. This implies that $O_{\mu\sigma\theta}
(f)=0$ for every $\mu\in U$. On the other hand, after replacing $\gamma$ by a
$\theta$-conjugate, we may assume that $u\in U$. So $I_M^G(\gamma,
\phi_\infty)=0$.

\end{proof}
\section{Stabilization of the invariant trace formula}
\label{GL_n_applications3}

In this section, we stabilize the invariant trace formula of
proposition \ref{prop:FT_invariante_tordue} 
if $\H$ is one of the quasi-split unitary groups of \ref{groupes1}.
Actually, there is nothing to stabilize; the invariant trace formula is
already stable in this case, and we simply show this.

We use the notations and assumptions of \ref{GL_n_applications2} and 
\ref{FT_stable_geometrique4}.

\begin{subproposition}\label{prop:FT_stable_tordue} Assume that $\H$ is one
of the quasi-split unitary groups of \ref{groupes1} and that $E$ is the
imaginary quadratic extension of $\Q$ that was used to define $\H$.
Let $f=\bigotimes\limits_vf_v\in C_c^\infty(\H(\A),\xi_H^{-1})$ and
$\phi=\bigotimes\limits_v\phi_v\in C_c^\infty(\G^0(\A),\xi^{-1})$
(where $\xi_H$ is the restriction of $\xi$ to $\A_H(\R)^0$).
Assume that, for every finite place $v$ of $\Q$, the functions
$f_v$ and $\phi_v$ are associated in the sense of \cite{La-CSCB} 3.2, that
the function $\phi_\infty$ is of the type considered in proposition
\ref{prop:FT_invariante_tordue} and that $f_\infty=\frac{1}
{|\Pi_V|}\sum\limits_{\pi_H\in\Pi_V}f_{\pi_H}$ (cf corollary
\ref{cor:pseudo_coeff}). 
Then there exists a constant $C\in\R^\times$ (depending only on $\H$ and
the choice of normalized intertwining operators on $\theta$-stable automorphic
representations of $\G^0(\Ade)$),
such that, for every Levi subset $\M$ of $\G$,
\[T_{M,geom}^G(\phi)=C\frac{d(\G)}{\tau(\H)}ST_{M_H}^H(f).\]
In particular,
\[T^G(\phi)=C\frac{d(\G)}{\tau(\H)}ST^H(f).\]

\end{subproposition}

\begin{subremarque}
After maybe choosing different normalized intertwining operators on the
$\theta$-stable automorphic representations of $\G^0(\Ade)$, we may
assume that $C$ is positive.

\end{subremarque}

\begin{proof} To see that the equalities for the $T_{M,geom}^G$ imply
the equality for $T^G$, it is enough to notice that the obvious map from the
set of $\G^0(\Q)$-conjugacy classes of Levi subsets of $\M$ to the set
$\H(\Q)$-conjugacy classes of Levi subgroups of $\H$ is a bijection,
and that $n_{M_H}^H=n_M^G$ if $\M_H$ corresponds
to $\M$.

Let $\M$ be a standard cuspidal Levi subset of $\G$. As the morphism
$\H^1(K,\M^0)\fl\H^1(K,\G^0)$ is injective (see the proof of lemma
\ref{lemme:transfert_terme_constant}) and
$\H^1(K,\G^0)=\{1\}$ for every field $K$, the assumption on
$\G^0$ implies that $d(\M_H,\M^0)=1$,
where $d(\M_H,\M^0)$ is defined in \cite{La-CSCB} 1.9.3.
By lemma \ref{lemme:egalite_D_M_G} and the fact that the descent formula
(corollary 8.3 of \cite{A-ITF1}) works just as well for twisted orbital
integrals, the proof of lemma \ref{lemme:transfert_terme_constant} applies
in the case considered here and shows that the functions $\phi_M$ and
$f_{M_H}$ are associated at every finite place. Using this fact and lemma
\ref{lemme:Phi_SPhi}, we may apply the stabilization process of chapter 4 of
\cite{La-CSCB}, on the group $\M\rtimes<\theta>$, to $T_{M,geom}^G(\phi)$.
As the set of places $\{\infty\}$ is $(\M,\M_H)$-essential (by lemma A.2.1 of
\cite{CL}, whose proof adapts immediately to the case of unitary similitude
groups), we get :
\[T_{M,geom}^G(\phi)=C\tau(\M^0)\tau(\M_H)^{-1}d(\M)k(\M_H)^{-1}
k(\H)ST_{M_H}^H(f),\]
with $C\in\R^\times$
(the factor $2^{-\dim(\agoth_G)}$ of \cite{La-CSCB} 4.3.2 does not appear here
because we are taking functions in $C^\infty_c(\G^0(\Ade),\xi^{-1})$ and not
in $C^\infty_c(\G^0(\Ade))$; and the factor $J_Z(\theta)$ does not appear
because, following Arthur, we consider the action of these functions on
$L^2(\G^0(\Q)\sous\G^0(\Ade),\xi)$ and not on
$L^2(\G^0(\Q)\A_{G^0}(\R)^0\sous\G^0(\Ade))$).

To finish the proof, it is enough to check that :
\[\frac{\tau(\M^0)}{\tau(\M_H)}\frac{d(\M)}{k(\M_H)}k(\H)=
d(\G)\frac{\tau(\G^0)}{\tau(\H)}.\]
By (ii) of lemma \ref{lemme:Tamagawa},
$\tau(\G^0)=\tau(\M^0)=1$.
So the equality above follows from remarks \ref{rq:k_tau} and
\ref{rq:d_k_tordu}.

\end{proof}

In the following lemma, we consider the situation of the beginning of
\ref{GL_n_applications2}, so that $\H$ is a cuspidal connected reductive group
over $\Q$, $E$ is an imaginary quadratic extension of $\Q$ and
$\G^0=R_{E/\Q}\H_E$. Fix a $\theta$-stable Borel subgroup of $\G^0$ (or,
equivalently, a Borel subgroup of $\H$).

\begin{sublemme}\label{lemme:Phi_SPhi} Use the notations of \cite{La-CSCB} 2.7.
Let $\phi_\infty$ be as in proposition \ref{prop:FT_invariante_tordue}, 
$\M$ be a standard cuspidal Levi subset of $\G$, and
$\gamma\in\M^0(\R)$ be $\theta$-semi-simple. Set $f_\infty=\frac{1}{|\Pi_V|}
\sum\limits_{\pi_H\in\Pi_V}f_{\pi_H}$ (cf corollary
\ref{cor:pseudo_coeff}). Then there exists a constant $C\in\R^\times$
(that is independent of $\M$ and
positive for a good choice of normalized intertwining operators)
such that :
\[\begin{array}{rcl}\displaystyle{\sum_{[x]\in\mathfrak{D}(I_\gamma,\M^0;\R)}
e(\delta_x)\Phi_M^G(\delta_x,\phi_\infty)} & = & \displaystyle{k(\M)^{-1}k(\G)
d(\M)S\Phi_M^G(\gamma,\phi_\infty)} \\
& = & \displaystyle{Cd(\M)k(\M_H)^{-1}k(\H)S\Phi_{M_H}^H(\Norme
\gamma,f_\infty),}\end{array}\]
and, if $\kappa\in\mathfrak{K}(I_\gamma,\M^0;\R)-\{1\}$,
\[\sum_{[x]\in\mathfrak{D}(I_\gamma,\M^0;\R)}e(\delta_x)<\kappa,\dot{x}>
\Phi_M^G(\delta_x,\phi_\infty)=0.\]

\end{sublemme}

\begin{proof} Once we notice that $\Phi_M^G(\gamma,\Theta_{\pi^
{\vee}})$ is invariant under stable $\theta$-conjugacy, the proof is
exactly the same as in theorem A.1.1 of \cite{CL}. To show the second line
of the first equality, use the definitions, theorem
\ref{th:relevement_caracteres}, remark \ref{rq:ep_tordu},
lemma \ref{lemme:egalite_D_M_G} and the fact that, if
$\gamma\in\G^0(\R)$ and $h\in\Norme\gamma$, then
$\H_h$ is an inner form of $\G^0_{\gamma\theta}$.

\end{proof}

\vspace{1cm}

We finish this section by recalling a few results on the transfer and
the fundamental lemma for base change.
Assume that we are in the situation of example
\ref{ex:groupes_non_connexes}, with $K$ a local field of characteristic $0$.
If two functions $f\in C_c^\infty(\G^0(K))$ and $h\in C_c^\infty(\H(K))$ 
are associated in the sense of \cite{La-CSCB} 3.2, we also say that $h$ is a
transfer of $f$ to $\H$.
\index{transfer (for the base change)}
Labesse proved the following result.

\begin{subtheoreme}\label{th:transfert_CB}
(\cite{La-CSCB} theorem 3.3.1 and proposition 3.5.2)
Let $f\in C_c^\infty(\G^0(K))$. Then there exists a transfer of $f$ to $\H$.

\end{subtheoreme}

Labesse has also proved a result about inverse transfer. We say that an
element $\gamma_H\in\H(K)$ \emph{is a norm} if there exists
$\gamma\in\G^0(K)$ such that $\gamma_H\in\Norme\gamma$.
\index{is a norm}
Assume that $K$ is non-archimedean.

\begin{subproposition}\label{prop:transfert_reciproque_CB} (\cite{La} 
proposition 3.3.2, proposition 3.5.3) Let
$h\in C_c^\infty(\H(K))$ be such that $SO_{\gamma_H}(h)=0$ for every
semi-simple $\gamma_H\in\H(K)$ that is not a norm.
Then there exists $f\in C_c^\infty(\G^0(K))$ such that
$h$ is a transfer of $f$.

\end{subproposition}

So, in order to determine which functions on $\H(K)$ are transfers of functions
on $\G^0(K)$, we need to describe the set of norms on $\H(K)$. To do this,
we use the results of 2.5 of \cite{La}. In the next lemma, assume that
$\H=\G(\U^*(n_1)\times\dots\times\U^*(n_r))$ with
$n_1,\dots,n_r\in\Nat^*$ (notations are as in \ref{groupes2}) and
that $\G^0=R_{E/\Q}\H_E$, where $E$ is the quadratic extension of $\Q$
used to define $\H$. Take $K=\Q_p$, where $p$ is a prime number.

\begin{sublemme}\label{lemme:description_normes}
Let $D_H=\H/\H^{der}$.
Then a semi-simple element of $\H(\Q_p)$ is a norm if and only if its
image in $D_H(\Q_p)$ is a norm. If $p$ splits and is unramified in $E$, or if
$\H=\GU^*(n)$ with $n$ odd, or if $p$ is unramified in $E$ and (at least) one
of the $n_i$ is odd, then every
semi-simple element of $\H(\Q_p)$ is a norm. \footnote
{Note that the first assertion of this lemma is also a consequence of
lemma 4.2.1 of \cite{Haines} (and of proposition 2.5.3 of \cite{La-CSCB}).
}

\end{sublemme}

\begin{proof} Notice that, if $p$ splits and is unramified in $E$,
then $\G^0(\Q_p)\simeq
\H(\Q_p)\times\H(\Q_p)$, and the naive norm map $\G^0(\Q_p)\fl\G^0(\Q_p)$,
$g\fle g\theta(g)$, is actually a surjection from $\G^0(\Q_p)$ to
$\H(\Q_p)$; there is a similar statement for $D_H$. So the lemma is
trivial in that case.

Hence, for the rest of the proof, we assume that there is only one place $\wp$
of $E$ above $p$, ie that $p$ is inert or ramified in $E$ (this is
just to avoid a discussion of cases; the results of Labesse apply of course
just as well in the general case).
As $\H^{der}$ is simply connected, $D_H(\Q_p)=\Ho^0_{ab}(\Q_p,
\H)$, in the notation of \cite{La-CSCB} 1.6. Similarly, for every Levi subgroup
$\M_H$ of $\H_{\Q_p}$, if we set $D_{M_H}=\M_H/\M_H^{der}$, then $D_{M_H}
(\Q_p)=\Ho^0_{ab}(\Q_p,\M_H)$. Let $\gamma$ be a semi-simple element of
$\H(\Q_p)$, and let $\M_H$ be a Levi subgroup of $\H_{\Q_p}$ such that
$\gamma\in\M_H(\Q_p)$ and $\gamma$ is elliptic in $\M_H$. By proposition
2.5.3 of \cite{La-CSCB}, $\gamma$ is a norm if and only if its image in
$D_{M_H}(\Q_p)$ is a norm. So, to prove the first statement of the lemma,
it is enough to show that an element of $D_{M_H}(\Q_p)$ is a norm if
and only if its image by the canonical map $D_{M_H}(\Q_p)\fl D_H(\Q_p)$ is a
norm. As $\H$ does not split over $\Q_p$,
$\M_H$ is $\H(\Q_p)$-conjugate to a Levi
subgroup of $\H$ (defined over $\Q$), so we may assume that $\M_H$ is a
standard Levi subgroup of $\H$ and is defined over $\Q$. By \ref{groupes2},
there exist $s,m_1,\dots,m_r\in\Nat$ such that $n_1+\dots+n_r=m_1+\dots+m_r+
2s$, $n_i-m_i$ is even for every $i$, and $\M_H\simeq (R_{E/\Q}\Gr_m)^s\times
\G(\U^*(m_1)\times\dots\times\U^*(m_r))$. 
The derived group of $\H$ is $\SU^*(n_1)\times\dots\times\SU^*(n_r)$, so
the map $\G(\U^*(n_1)\times\dots\U^*(n_r))\fl\Gr_m\times(R_{E/\Q}\Gr_m)^r$,
$(g_1,\dots,g_r)\fle (c(g_1),\det(g_1),\dots,\det(g_r))$ induces an
isomorphism
\[D_H\iso\{(\lambda,z_1,\dots,z_r)\in\Gr_m\times(R_{E/\Q}\Gr_m)^r|
\forall i,z_i\overline{z}_i=\lambda^{n_i}\}.\]
Similarly, there is an isomorphism
\[D_{M_H}\iso D_{\M_H,l}\times D_{M_H,h},\]
where $D_{M_H,l}=(R_{E/\Q}\Gr_m)^s$ and
\[D_{M_H,h}=\{(\lambda,z_1,\dots,z_r)\in
\Gr_m\times(R_{E/\Q}\Gr_m)^r|\forall i,z_i\overline{z}_i=\lambda^{m_i}\mbox
{ if }m_i>0\mbox{ and }z_i=1\mbox{ if }m_i=0\}.\]
The canonical map $D_{M_H}\fl D_H$ sends the factor $D_{M_H,l}$ to
$1$ and is induced on the factor $D_{M_H,h}$ by the map
$(\lambda,z_1,\dots,z_r)\fle (\lambda,\lambda^{(n_1-m_1)/2}z_1,\dots,\lambda
^{(n_r-m_r)/2}z_r)$. As every element in $D_{M_H,l}(\Q_p)$ is obviously a norm,
it is now clear that an element of $D_{M_H}(\Q_p)$ is a norm if and only
if its image in $D_H(\Q_p)$ is a norm.

Assume that $\H=\GU^*(n)$ with $n$ odd, and write $n=2m+1$, $m\in\Nat$. Then
it is easy to check that the map $\H\fl R_{E/\Q}\Gr_m$, $g\fle\det(g)c(g)
^{-m}$, induces an isomorphism $D_H\iso R_{E/\Q}\Gr_m$. So every element of
$D_H(\Q_p)$ is a norm, and consequently every semi-simple element of $\H(\Q_p)$
is a norm.

Assume that $\H=\G(\U^*(n_1)\times\dots\times\U^*(n_r))$ with $n_1,\dots,n_r
\in\Nat^*$, that $n_1$ is odd and that $p$ is inert and unramified in $E$.
Write $n_1=2m_1+1$, $m_1\in\Nat$. Then the map $\Gr_m\times(R_{E/\Q}\Gr_m)^r$,
$(\lambda,z_1,\dots,z_r)\fle(z_1\lambda^{-m_1},z_2(z_1^{-1}\lambda^{m_1})
^{n_2},z_3(z_1^{-1}\lambda^{m_1})^{n_3},\dots,z_r(z_1^{-1}\lambda^{m_1})
^{n_r})$ (together with the description of $D_H$ given above) induces an
isomorphism
\[D_H\simeq R_{E/\Q}\Gr_m\times\U(1)^{r-1}.\]
It is obvious that every element of $(R_{E/\Q}\Gr_m)(\Q_p)=E_p^\times$ is a
norm, so, to finish the proof of the lemma, it is enough to show that
every element of $\U(1)(\Q_p)$ is a norm. Let $z\in\U(1)(\Q_p)$. Then
$z$ is an element of $E_p^\times$ such that $z\overline{z}=1$, and we
want to show that there exists $y\in E_p^\times$ such that $z=y\overline{y}
^{-1}$. Write $z=ap^k$, with $a\in\Of_{E_p}^\times$ and $k\in\Z$. Then
$z\overline{z}=a\overline{a}p^{2k}=1$, so $k=0$ and $a\overline{a}=1$, and we
want to show that there exists $b\in\Of_{E_p}^\times$ such that $a=b
\overline{b}^{-1}$. By Hensel's lemma, it is enough to check the analog of this
for the reduction modulo $p$ of $a$. As $p$ is inert and unramified in $E$,
$\Of_{E_p}/(p)=\Fi_{p^2}$.
Let $u:\Fi_{p^2}^\times\fl\Fi_{p^2}^\times$
be the group morphism that sends $b$ to $b\overline{b}^{-1}$. Then $\varphi(b)
=b^{1-p}$ for every $b$, so the image of $\varphi$ is of cardinality
$p^2-1/(p-1)=p+1$. But this image is contained in $\U(1)(\Fi_p)$, and
$\U(1)(\Fi_p)=\{a\in\Fi_{p^2}^\times|a^{p+1}=1\}$ is of cardinality $p+1$,
so $\varphi(\Fi_{p^2}^\times)=\U(1)(\Fi_p)$. This finishes the proof.

\end{proof}

The above lemma (together with the result of Labesse about inverse transfer, ie
proposition \ref{prop:transfert_reciproque_CB}, and the fact that the group of
norms in a torus contains an open neighbourhood of $1$)
has the following immediate consequence :

\begin{sublemme}\label{lemme:image_CB}
If $p$ splits and is unramified in $E$,
or if $\H=\GU^*(n)$ with $n$ odd, or
if $p$ is unramified in $E$ and one of the $n_i$ is odd, then
every function in $C_c^\infty(\H(\Q_p))$ is a transfer of a function in
$C_c^\infty(\G^0(\Q_p))$. In general, the set of functions in
$C_c^\infty(\H(\Q_p))$ that are a transfer of a function in
$C_c^\infty(\G^0(\Q_p))$ is a subalgebra of $C_c^\infty(\H(\Q_p))$, and it
contains all the functions with small enough support.

\end{sublemme}

Transfer is explicit if we are in an unramified situation. Assume
that $K$ is non-archimedean, that the group $\H$ is unramified
over $K$ and that the extension $E/K$ is unramified.
Let $\K_G$ and $\K_H$ be hyperspecial maximal compact subgroups of
$\G^0(K)$ and $\H(K)$ such that
$\K_H=\H(K)\cap\K_G$ and $\theta(\K_G)=\K_G$. 
The $L$-morphism $\xi:{}^L\H\fl{}^L\G^0$ defined in example
\ref{ex:groupes_non_connexes} induces an morphism of algebras
$b:\Hecke(\G^0(\K),\K_G)\fl\Hecke(\H(K),\K_H)$, called base change morphism.
The following theorem, known under the name of ``fundamental lemma for
base change'', is due to Kottwitz (for the unit element of
$\Hecke(\G^0(\K),\K_G)$), Clozel and Labesse (for the other elements).

\begin{subtheoreme}\label{th:lemme_fondamental_CB}(\cite{K-BC}, \cite{Cl-LF},
\cite{La}, \cite{La-CSCB} 3.7) Let $f\in\Hecke(\G^0(K),\K_G)$. Then
$b(f)$ is a transfer of $f$ to $\H$.

\end{subtheoreme}
\index{fundamental lemma for base change}

Let us write down explicit formulas for the base change morphism in the case
of unitary groups.
Let $\H=\G(\U^*(n_1)\times\dots\times\U^*(n_r))$, $E$ 
be the imaginary quadratic extension of $\Q$ used to define $\H$ and
$p$ be a prime number that is unramified in $E$.
The groups $\G^0$ and $\H$ have obvious $\Z_p$-models (cf remark
\ref{rq:groupes_sur_Z}), and we take $\K_G=\G^0(\Z_p)$ and
$\K_H=\H(\Z_p)$. Use the notations of chapter \ref{partie_en_p}.

If $p$ is inert in $E$, the base change morphism is calculated in section
\ref{partie_en_p2} (with $L=E_p$ and $\G=\H$).

Assume that $p$ splits in $E$. Then $\G^0_{\Q_p}\simeq\H_{\Q_p}
\times\H_{\Q_p}$, and, for every $g=(g_1,g_2)\in\G^0(\Q_p)=\H(\Q_p)\times
\H(\Q_p)$, $g_1g_2\in\Norme g$. To simplify notations, we assume that
$r=1$. Then there is an isomorphism (defined in \ref{partie_en_p2})
\[\Hecke(\H(\Q_p),\K_H)\simeq\C[X^{\pm 1}]\otimes\C[X_1^{\pm 1},\dots,
X_n^{\pm 1}]^{\Sgoth_n}.\]
So there is an obvious isomorphism
\[\Hecke(\G^0(\Q_p),\K_G)\simeq \C[Z_1^{\pm 1}]\otimes\C[Z_{1,1}^{\pm 1},\dots,
Z_{n,1}^{\pm 1}]^{\Sgoth_n}\otimes\C[Z_2^{\pm 1}]\otimes\C[Z_{1,2}^{\pm 1},
\dots,Z_{n,2}^{\pm 1}]^{\Sgoth_n},\]
and the base change morphism is induced by
\[Z_j\fle X,\qquad Z_{i,j}\fle X_i.\]

In particular, the base change morphism is surjective if $p$ splits in $E$.
If $p$ is inert in $E$, then the image of the base change morphism is given
in remark \ref{rq:image_BC}; in particular, the base
change morphism is surjective if and only if one of the $n_i$ is odd.

The following lemma will be useful in the applications of the next section.
We assume again that $\H$ is any connected unramified group on a
non-archimedean local field $K$ and that the extension $E/K$ is unramified,
and we choose hyperspecial maximal compact subgroups $\K_H$ and 
$\K_G$ of $\H(K)$ and $\G^0(K)$ as before.

\begin{sublemme}\label{lemme:rep_nr_tordues} Let $\pi$ be a $\theta$-stable
admissible irreducible representation of $\G^0(K)$ and $A_\pi$ be a normalized
intertwining operator on $\pi$. Let $\varepsilon$ be an element of
$\{\pm 1\}$ such that $A_\pi$ acts on $\pi^{\K_G}$ by multiplication by
$\varepsilon$ (such a $\varepsilon$ exists because $A_\pi$ stabilizes
$\pi^{\K_G}$ and $\dim\pi^{\K_G}\leq 1$). Then, for every
$f\in\Hecke(\G^0(K),\K_G)$,
\[\Tr(\pi(f)A_\pi)=\varepsilon\Tr(\pi(f)).\]

\end{sublemme}

\section{Applications}
\label{GL_n_applications4}

Notations here are slightly different from the ones used in
\ref{GL_n_applications3}. Let $\H=\G(\U(p_1,q_1)\times\dots\times
\U(p_r,q_r))$ (this group is defined in \ref{groupes1}),
$\H^*$ be a quasi-split inner form of $\H$ (so $\H^*=\G(\U^*(n_1)\times\dots
\times\U^*(n_r))$, where $n_i=p_i+q_i$), $E$ be the imaginary quadratic
extension of $\Q$ that was used in the definition of $\H$.
Set $\G^0=R_{E/\Q}\H^*_E$.
Of course, $\G^0\simeq R_{E/\Q}\H_E$.

If $V$ is an irreducible algebraic representation of $\H$, let
$\phi_V$ be a twisted pseudo-coefficient of the $\theta$-discrete
representation $\pi_V$ of $\G^0(\R)$ associated to $\varphi_{|W_\C}$, where
$\varphi:W_\R\fl{}^L\H$ is a Langlands parameter of the $L$-packet of
the discrete series of $\H(\R)$ associated to $V$.
\index{$\pi_V$}
\index{$\phi_V$}

Let $\Mcal'_\G$ be the set of conjugacy classes of Levi subsets $\M$ of
$\G$ such that, for every $i\in\{1,\dots,r\}$,
$\M^0\cap R_{E/\Q}\GL_{n_i,E}$ is equal to $R_{E/\Q}\GL_{n_i,E}$ or to a
maximal Levi subgroup of $R_{E/\Q}\GL_{n_i,E}$. 
Let $\M\in\Mcal'_\G$. Then there exist non-negative integers $n_1^+,n_1^-,
\dots,n_r^+,n_r^-$ such that, for every $i\in\{1,\dots,r\}$, $n_i=n_i^++n_i^-$
and $\M^0\cap R_{E/\Q}\GL_{n_i,E}=R_{E/\Q}\GL_{n_i^+,E}\times R_{E/\Q}
\GL_{n_i^-,E}$.
Let $\Mcal_\G$ be the set of $\M\in\Mcal'_\G$ such that we can choose the 
$n_i^+,n_i^-$ so that $n_1^-+\dots+n_r^-$ is even.
\index{MG@$\Mcal_G$}
If $\M$ is in $\Mcal_\G$ and the $n_i^+,n_i^-$ are as above, then we may
assume that $n_1^-+\dots+n_r^-$ is even; let $(\H_M,s_{H_M},\eta_{H_M,0})$
be the elliptic endoscopic datum for $\H$ defined by the $n_i^+,n_i^-$ as in
proposition \ref{prop:groupes_endoscopiques}, and $\eta_{H_M}$ be a
$L$-morphism extending $\eta_{H_M,0}$ as in proposition
\ref{prop:prolongement_eta_0}.
This defines a bijection between $\Mcal_\G$ and the set $\Ell_\H$ of
\ref{applications1}.

Let $\M\in\Mcal_\G$.
Let $\xi:{}^L\H={}^L\H^*\fl{}^L\G^0$ be the $L$-morphism defined in
example \ref{ex:groupes_non_connexes}; as $\M^0=R_{E/\Q}
\H_{M,E}$, we get in the same way a $L$-morphism $\xi_M:{}^L\H_M\fl{}^L\M^0$.
Let $\eta_M$ be the morphism
\[\begin{array}{rcl}{}^L\M^0\simeq (\widehat{\H}_M\times\widehat{\H}_M)\rtimes 
W_\Q & \fl & {}^L\G^0\simeq (\widehat{\H}\times\widehat{\H})\rtimes W_\Q \\
((h_1,h_2),w) & \fle & ((\eta_{H_M,1}(h_1,w),\eta_{H_M,1}(h_2,w)),w),
\end{array}\]
where $\eta_{H_M,1}:{}^L\H_M\fl\widehat{\H}$ is the first component of
$\eta_{H_M}$. It is clear that $\eta_M$ is a $L$-morphism that makes the
following diagram commute :
\[\xymatrix{{}^L\H_M\ar[r]^-{\eta_{H_M}}\ar[d]_{\xi_M} & {}^L\H\ar[d]^\xi \\
{}^L\M^0\ar[r]_-{\eta_M} & {}^L\G^0}\]
Note that the embedding $\widehat{\M}^0\fl\widehat{\G}^0$ induced by
$\eta_M$ is $W_\Q$-equivariant. Let $\eta_{M,simple}:{}^L\M^0\fl{}^L\G^0$
be the obvious $L$-morphism extending this embedding (ie the one that is
equal to identity on $W_\Q$). Write $\eta_M=c_M\eta_{M,simple}$, where
$c_M:W_\Q\fl Z(\widehat{\M}^0)$ is a $1$-cocycle, and let $\chi_M$ be
the quasi-character of $\M^0(\Ade)$ associated to the class of $c_M$ in
$\Ho^1(W_\Q,Z(\widehat{\M}^0)$. (In general, $\chi_M$ can be non-trivial.)

Let $S$ be a set of places of $\Q$. Write $\Ade_S=\prod\limits_
{v\in S}\!\!\!{}'\,\Q_v$ and $\Ade^S=\prod\limits_{v\not\in S}\!\!\!{}'\,\Q_v$.
We say that
a function $f_S\in C_c^\infty(\H(\Ade_S))$ satisfies condition (H) if,
for every $\M\in\Mcal_G$, there exists a transfer $f_S^{H_M}$ of $f_S$ to
$\H_M$ and a function $\phi_{S,M}\in C_c^\infty(\M^0(\Ade_S))$ such that
the functions $\phi_{S,M}$ and $f_S^{H_M}$ are associated at every place
in $S$.
\index{condition (H)}

The next lemma gives examples of functions that satisfy condition (H).
For every place $v$ of $\Q$, we say that a semi-simple
element $\gamma\in\H(\Q_v)$ \emph{is a norm} if there exists $g\in\G^0(\Q_v)$
such that $\gamma\in\Norme g$ (this condition makes sense because
$\Norme g$ is a stable conjugacy class in $\H^*(\Q_v)$ and $\H$ is an inner
form of $\H^*$).
\index{is a norm}

\begin{sublemme}\label{lemme:hypothese_H} Let $v$ be a finite place of $\Q$.
\begin{itemize}
\item[(i)] Every function in $C_c^\infty(\H(\Q_v))$ with support in a small
enough neighbourhood of $1$ satisfies condition (H).
\item[(ii)] Assume that $\H$ is quasi-split over $\Q_v$ (but not
necessarily unramified). Then, for every $\phi\in C_c^\infty(\G^0(\Q_v))$
with support in a small enough neighbourhood of $1$,
there exists $f\in C_c^\infty(\H(\Q_v))$
associated to $\phi$ and satisfying condition (H).
\item[(iii)] Assume that $v$ is unramified in $E$ (so
$\H_{\Q_v}=\H^*_{\Q_v}$ is unramified).
Let $\M\in\Mcal_G$.
Then the commutative diagram 
\[\xymatrix{{}^L\H_M\ar[r]^-{\eta_{H_M}}\ar[d]_{\xi_M} & {}^L\H\ar[d]^\xi \\
{}^L\M^0\ar[r]_-{\eta_M} & {}^L\G^0}\]
gives a commutative diagram
\[\xymatrix{\Hecke(\G^0(\Q_v),\G^0(\Z_v))\ar[r]\ar[d] & 
\Hecke(\H(\Q_v),\H(\Z_v))\ar[d] \\
\chi_{M,v}\Hecke(\M^0(\Q_v),\M^0(\Z_v))\ar[r] & 
\chi_{\eta_{H_M},v}\Hecke(\H_M(\Q_v),\H_M(\Z_v))}\]
(where $\chi_{\eta_{H_M},v}$ is defined as in the last two subsections of
\ref{partie_en_p2}),
satisfying the following properties :
\begin{itemize}
\item[-] the upper horizontal arrow is the base change map;
\item[-] the lower horizontal arrow sends a function
$\chi_{M,v}\phi_v$, with $\phi_v\in\Hecke(\M^0(\Q_v),\M^0(\Z_v))$,
to the function $\chi_{\eta_{H_M},v}f_v$, where $f_v\in
\Hecke(\H_M(\Q_v),\H_M(\Z_v))$ is the image of $\phi_v$ by the base change
map;
\item[-] the left vertical arrow sends a function in
$\Hecke(\G^0(\Q_v),\G^0(\Z_v))$ to the product of its constant term at
$\M^0$ and of $\chi_{M,v}$;
\item[-] the right vertical arrow is the transfer map
defined by $\eta_{H_M}$ as in \ref{partie_en_p2}.

In particular, every function in the image of
the base change morphism $\Hecke(\G^0(\Q_v),\G^0(\Z_v))\fl
\Hecke(\H(\Q_v),\H(\Z_v))$ satisfies condition (H).

\item[(iv)] If $v$ is unramified in $E$ and one of the $n_i$ is odd, then
every function in $C_c^\infty(\H(\Q_v))$ satisfies condition (H).

\end{itemize}

\end{itemize}
\end{sublemme}

\begin{proof} Point (iii) is immediate (because the fundamental lemma
is known, cf \ref{FT_stable_geometrique3}).
Point (iv) is a direct consequence of lemma \ref{lemme:image_CB}.

We show (i). For every $\M\in\Mcal_G$, there is a $\H^*(\Q)$-conjugacy
class of embeddings $\H_M\fl\H^*$; fix an embedding in this class.
Identify $\H_E$ and $\H^*_E$ with $\Gr_{m,E}\times\GL_{n_1,E}
\times\dots\times\GL_{n_r,E}$ using the morphism defined in the beginning
of \ref{groupes3}. Let $\M\in\Mcal_G$. There exists an open neighbourhood
$U_M$ of $1$ in $\H_M(\Q_v)$ such that every semi-simple element in $U_M$
is a norm.
\quash{ (because taking the square of an element is an isomorphism of
varieties in the neighbourhood of $1$).}
Choose an open neighbourhood $V_M$ of $1$ in
$\H^*(E\otimes_\Q\Q_v)=\H(E\otimes_\Q\Q_v)$ such that every semi-simple
element of $\H_M(\Q_v)$ that is $\H^*(E\otimes_\Q\Q_v)$-conjugate to an
element of $V_M$ is $\H_M(\Q_v)$-conjugate to an element of $U_M$
(cf lemma \ref{lemme:beurk} below).

Let $V=\bigcap\limits_{\M\in\Mcal_G}V_M$ and $U=V\cap\H(\Q_v)$.
Then $U$ is an open neighbourhood of $1$ in $\H(\Q_v)$. Let $f\in
C_c^\infty(\H(\Q_v))$ with support contained in $U$. We show that $f$
satisfies condition (H). For every $\M\in\Mcal_G$, choose a transfer
$f^{H_M}$ of $f$ to $\H_M$. To show that there exists a function in
$C_c^\infty(\M^0(\Q_v))$ associated to $f^{H_M}$, it is enough, by
proposition 3.3.2 of \cite{La-CSCB}, to show that, for every semi-simple
$\gamma\in\H_M(\Q_v)$, $SO_\gamma(f^{H_M})=0$ if $\gamma$ is not a norm.
Let $\gamma\in\H_M(\Q_v)$ be semi-simple and such that $SO_\gamma(f^{H_M})\not=
0$. Then, by the definition of the transfer, there exists an image $\delta$ of
$\gamma$ in $\H(\Q_v)$ such that $O_\delta(f)\not=0$. In other words,
$\gamma$ is $\H^*(E\otimes_\Q\Q_v)$-conjugate to an element of $U$.
As $U\subset V_M$, this implies that $\gamma$ is $\H_M(\Q_v)$-conjugate
to an element of $U_M$, hence that $\gamma$ is conjugate to a norm, ie
that $\gamma$ is itself a norm.

We show (ii). By (i), it is enough to check that, if $U$ is a neighbourhood
of $1$ in $\H(\Q_v)$, then there exists a neighbourhood $V$ of $1$ in
$\G^0(\Q_v)$ such that every function $\phi\in C_c^\infty(\G^0(\Q_v))$
with support contained in $V$ admits a transfer
$f\in C_c^\infty(\H(\Q_v))$ with support contained in $U$. This follows
from the proof of theorem 3.3.1 of \cite{La-CSCB}.

\end{proof}

\begin{sublemme}\label{lemme:beurk} Let $F$ be a local field of
characteristic $0$, $E$ be a finite extension of $F$ and $\H$ be a connected
reductive group on $F$. Set $\G=R_{E/F}\H_E$. Let $\M$ be a Levi subgroup
of $\G$. Assume that there exists a connected reductive group $\H_M$ on
$F$ such that $\M=R_{E/F}\H_{M,E}$ ($\H_M$ is not necessarily a subgroup
of $\H$). Let $U$ be a neighbourhood of $1$ in $\H_M(F)$. Then there exists
a neighbourhood $V$ of $1$ in $\G(F)$ such that : for every semi-simple
$\gamma\in\H_M(F)$, if $\gamma$ is $\G(F)$-conjugate to an element of
$V$, then $\gamma$ is $\H_M(F)$-conjugate to an element of $U$.

\end{sublemme}

\begin{proof} Let $(\Se_1,\dots,\Se_r)$ be a system of representatives of
the set of $\H_M(F)$-conjugacy classes of maximal tori of $\H_M$
(this set is finite because the characteristic of $F$ is $0$).
For every $i\in\{1,\dots,r\}$, set $\T_i=R_{E/F}\Se_i$ and
$U_i=U\cap\Se_i(F)$, and choose a neighbourhood $W_i$ of $1$ in
$\T_i(F)=\Se_i(E)$ such that $W_i\cap\Se_i(F)\subset U_i$. Let $i\in\{1,
\dots,r\}$. Then $\T_i$ is a maximal torus of $\M$, hence of $\G$, so, by
lemme 3.1.2 of \cite{La-CSCB}, there exists a neighbourhood $V_i$ of $1$
in $\G(F)$ such that, if an element $t\in\T_i(F)$ has a conjugate in $V_i$,
then $t\in W_i$. Set $V=\bigcap\limits_{i=1}^r V_i$.

Let $\gamma\in\H_M(F)$ be semi-simple and $\G(F)$-conjugate to an element
of $V$. As $\gamma$ is semi-simple, there exists a maximal torus of $\H_M$
containing $\gamma$, so we may assume that there exists
$i\in\{1,\dots,r\}$ such that
$\gamma\in\Se_i(F)$. In particular, $\gamma\in\T_i(F)$. As
$\gamma$ is $\G(F)$-conjugate to an element of $V_i$,
$\gamma\in W_i$. But $W_i\cap\Se_i(F)\subset U_i$, so that
$\gamma\in U_i\subset U$.

\end{proof}

We come back to the situation of the beginning of this section.
Fix a prime number $p$ that is unramified in $E$, a neat open compact
subgroup $\K=\K^p\H(\Z_p)$ (with $\K^p\subset\H(\Af^p)$) of $\H(\Af)$,
an irreducible algebraic
representation $V$ of $\H$ and a function
$f^{p,\infty}\in\Hecke(\H(\Af^p),\K^p)$. Assume that
$f^{p,\infty}$ satisfies condition (H).

Let $\M\in\Mcal_\G$, and define, for every $j\in\Z$, a function
$\phi_M^{(j)}=\phi_M^{p,\infty}\phi_{M,p}^{(j)}\phi_{M,\infty}\in C^\infty
(\M^0(\Ade))$, compactly supported modulo $\A_{M^0}(\R)^0$, in the
following way.
Choose $\phi_M^{p,\infty}\in C_c^\infty(\M^0(\Af^p))$ that is associated at
every place to a transfer $(f^{p,\infty})^{H_M}$ of $f^{p,\infty}$ to $\H_M$.
The calculations of \ref{partie_en_p2} and (iii) of lemma
\ref{lemme:hypothese_H} show that the function
$f_{\H_M,p}^{(j)}$ defined in definition \ref{def:phi_m_G} is the
product of $\chi_{\eta_{H_M},p}$ and of a spherical function in
the image of the base change map $\Hecke(\M^0(\Q_p),\M^0
(\Z_p))\fl\Hecke(\H_M(\Q_p),\H_M(\Z_p))$. Take $\phi_{M,p}^{(j)}\in\Hecke
(\M^0(\Q_p),\M^0(\Z_p))$ to be $\chi_{M,p}\phi'$, where $\phi'$ is
any spherical function in the inverse image of
$\chi_{\eta_{H_M},p}^{-1}f_{\H_M,p}^{(j)}$ by the base change map.
\index{$\phi_{M,p}^{(j)}$}
To define $\phi_{M,\infty}$, use the notations introduced before and in lemma
\ref{lemme:conj_RP3}. Lemma \ref{lemme:conj_RP3} gives
an irreducible algebraic representation $V_\omega$ of $\H_M$ for every
$\omega\in\Omega_*\simeq\Phi_H(\varphi)$. Take
\[\phi_{M,\infty}=\sum_{\omega\in\Omega_*}\det(\omega)\phi_{V_\omega^*},\]
where $\det(\omega)$ is defined in remark \ref{rq:facteur_transfert} and
$\phi_{V_\omega^*}$ is defined at the beginning of this section.
\index{$\phi_{M,\infty}$}

Let
\[c_M=(-1)^{q(\H)}\iota(\H,\H_M)C_M\frac{\tau(\H_M)}{\tau(\M^0)}k(\M)^{-1}
<\mu_H,s_{H_M}>\in\Q^\times,\]
where $\mu_H$ is the cocharacter of $\H_E$ determined by the Shimura
datum as in \ref{groupes1} and $C_M\in\R^\times$ is the constant of
proposition \ref{prop:FT_stable_tordue} for $\M$.
\index{cM@$c_M$}

\begin{subtheoreme}\label{th:stab_FT_IC_GL_n} For every $j\in\Z$,
\[\Tr(f^{p,\infty}\Phi_\wp^j,R\Gamma(M^{\K}(\H,\X)^*_{\overline{\Q}},IC^{\K}V_
{\overline{\Q}}))=\sum_{\M\in\Mcal_\G}c_M T^M(\phi_M^{(j)}),\]
where $(\H,\X)$ is the Shimura datum of \ref{groupes1} and
$\Phi_\wp$ is defined in \ref{applications2}.

\end{subtheoreme}

\begin{proof} The theorem is an easy consequence of corollary
\ref{cor:stab_FT_IC} (and remark \ref{rq:pour_tout_m}),
lemma \ref{lemme:conj_RP3} and proposition
\ref{prop:FT_stable_tordue} (see also remark \ref{rq:ensemble_Sigma} for
the choice of $p$).

\end{proof}

It is possible to deduce from theorem
\ref{th:stab_FT_IC_GL_n} and proposition \ref{prop:FT_invariante_tordue}
an expression for the logarithm of the $L$-function (at a good prime
number) of the intersection complex $IC^{\K}V$, if $\K$ is a small enough
open compact subgroup of $\G(\Af)$.

Remember that we defined in \ref{groupes1} a morphism
$\mu_H:\Gr_{m,E}\fl\H_E$. The formula for $\mu_H$ is :
\[\mu_H:\left\{\begin{array}{rcl}\Gr_{m,E} & \fl & \H_E=\Gr_{m,E}\times
\GL_{n_1,E}\times\dots\times \GL_{n_r,E} \\
z & \fle & (z,\left(\begin{array}{cc}zI_{p_1} & 0 \\ 0 & I_{q_1}\end{array}
\right),\dots,\left(\begin{array}{cc}zI_{p_r} & 0 \\ 0 & I_{q_r}\end{array}
\right))\end{array}\right.\]
For every $\M\in\Mcal_G$, let $M_{\H_M}$ be the set of $\H_M(E)$-conjugacy
classes of cocharacters $\mu_{H_M}:\Gr_{m,E}\fl\H_{M,E}$ such that the
cocharacter $\Gr_{m,E}\stackrel{\mu_{H_M}}{\fl}\H_{M,E}\fl\H_E$
is $\H(E)$-conjugate to $\mu_H$.
Let $\M\in\Mcal_G$. Write as before
$\H_{M,E}=\Gr_{m,E}\times\GL_{n_1^+,E}\times\GL_{n_1^-,E}
\times\dots\times\GL_{n_r^+,E}\times\GL_{n_r^-,E}$. Then every element
$\mu_{H_M}$ of $M_{\H_M}$ has a unique representative of the form
\[z\fle (z,\left(\begin{array}{cc}zI_{p_1^+} & 0 \\ 0 & I_{q_1^+}\end{array}
\right),\left(\begin{array}{cc}zI_{p_1^-} & 0 \\ 0 & I_{q_1^-}\end{array}
\right),\dots,\left(\begin{array}{cc}zI_{p_r^+} & 0 \\ 0 & I_{q_r^+}\end{array}
\right),
\left(\begin{array}{cc}zI_{p_r^-} & 0 \\ 0 & I_{q_r^-}\end{array}\right)),\]
with $p_i^++p_i^-=p_i$. Write $s(\mu_{H_M})=p_1^-+\dots+p_r^-$ and
$d(\mu_{H_M})=p_1^+q_1^++p_1^-q_1^-+\dots+p_r^+q_r^++p_r^-q_r^-$. Let 
$d=d(\mu_H)=p_1q_1+\dots+p_rq_r$ ($d$ is the dimension of
$M^{\K}(\H,\X)$, for every open compact subgroup $\K$ of $\H(\Af)$).

Remember that every cocharacter $\mu_{H_M}:\Gr_{m,E}\fl\H_{M,E}$ defines a
representation $r_{-\mu_{H_M}}$ of ${}^L\H_{M,E}$ (cf lemma \ref{def_r_mu}).

We recall the definition of the $L$-function at a good $p$ of the
intersection complex.

\begin{subdefinition}\label{def:fonction_L_IC} Let $p$ a prime number
as in \ref{points_fixes3}, and let $\wp$ be a place of $E$ above $p$.
Set
\[\log L_\wp(s,IC^{\K}V)=\sum_{m\geq 1}\frac{1}{m}(N\wp)^{-ms}
\Tr(\Phi_\wp^{m*},
R\Gamma(M^{\K}(\G,\X)^*_{\overline{\Q}},IC^{\K}V_{\overline{\Q}})),\]
where $\Phi_\wp\in W_{E_\wp}$ is a lift of the geometric Frobenius,
$N\wp=\#(\Of_{E_\wp}/\wp)$, and
$s\in\C$ (the series converges for $Re(s)>>0$).
\index{Lfunction@$L$-function of the intersection complex}

\end{subdefinition}

\begin{subcorollaire}\label{cor:calcul_L_GL_n} Let $\K$ be a small
enough open compact subgroup of $\H(\Af)$. Then there exist
functions $\phi_M\in C^\infty(\M^0(\Ade))$ with compact support modulo
$\A_{M^0}(\R)^0$, for every $\M\in\Mcal_G$, such that, for every prime
number $p$ as in \ref{points_fixes3} (ie such that $p$ is unramified in $E$
and $\K=\K^p\G(\Z_p)$) and for every place $\wp$ of $E$ above $p$,
\begin{flushleft}$\displaystyle{
\log L_\wp(s,IC^{\K}V)=\sum_{\M\in\Mcal_G}c_M\sum_{t\geq 0}
\sum_{\pi_M\in\Pi_{disc}(\M,t)}
a_{disc}^M(\pi_M)
}$\end{flushleft}
\begin{flushright}$\displaystyle{
\Tr(\pi_M(\phi_M)A_{\pi_M})
\sum_{\mu_{H_M}\in M_{\H_M}}(-1)^{s(\mu_{H_M})}
\log L_\wp(s-\frac{d}{2},(\pi_M\otimes\chi_M)_\wp,r_{-\mu_{H_M}}),
}$\end{flushright}
where, for every $\M\in\M_G$ and $\pi_M\in\Pi_{disc}(\M,t)$,
$(\pi_M\otimes\chi_M)_\wp$ is the local component at $\wp$ of $\pi_M\otimes
\chi_M$, seen as a representation of $\H_M(\Ade_E)$.

\end{subcorollaire}

\begin{proof} By lemma \ref{lemme:hypothese_H}, if $\K$ is a small enough
open compact subgroup of $\H(\Af)$, then the function $\ungras_{\K}$
satisfies condition (H). Fix such a $\K$, and assume also that $\K$ is neat.
Let $S$ be a finite set of prime numbers containing the set of prime
numbers that are ramified in $E$ and such that $\K=\K_S\K^S$, with
$\K_S\subset\H(\Ade_S)$ and $\K^S=\prod\limits_{p\not\in S}\H(\Z_p)$.
For every $\M\in\Mcal_G$, choose a transfer $f_S^{H_M}$ of $\ungras_{\K_S}$
to $\H_M$ and a function $\phi_{M,S}\in C_c^\infty(\M^0(\Ade_S))$ associated
to $f_S^{H_M}$, and write $\phi_M^S=\chi_{M|\M^0(\Af^S)}
\ungras_{\K_M^S}$, where $\K_M^S=\prod
\limits_{p\not\in S}\M^0(\Z_p)$, and $\phi_M=\phi_{M,S}\phi_M^S$.

Let $p\not\in S$ and $j\in\Nat^*$. We want to define, for every
$\mu_{H_M}\in M_{\H_M}$, a function $\phi_{\mu_{H_M},p}^{(j)}\in
\Hecke(\M^0(\Q_p),\M^0(\Z_p))$. Remember that we fixed a place $\wp$ of $E$
above $p$. Let $L$ be the unramified extension of $E_\wp$ of degree $j$ in
$\overline{\Q}_p$. For every $\mu_{H_M}\in M_{\H_M}$, let 
$\phi_{\mu_{H_M},p}^{(j)}$ be the product of $(N\wp)^{j(d-d(\mu_{H_M}))/2}$
and of the image of the function $f_{\mu_{H_M},L}$ in $\Hecke(\M^0(L),
\M^0(\Of_L))$ (defined by $\mu_{M,L}$ as in \ref{partie_en_p1}) by the
morphism
\[\Hecke(\H_M(L),\H_M(\Of_L))\fl\Hecke(\H_M(E_\wp),\H_M(\Of_{E_\wp}))\fl\Hecke(
\M^0(\Q_p),\M^0(\Z_p)),\]
where the first arrow is the base change morphism and the second arrow is
\begin{itemize}
\item[-] identity if $p$ is inert in $E$ (so $\M^0(\Q_p)=\H_M(E_\wp)$);
\item[-] the morphism $h\fle (h,\ungras_{\H_M(\Of_{E_{\wp'}})})$ if $p$
splits in $E$ and $\wp'$ is the second place of $E$ above $p$ (so
$\M^0(\Q_p)=\H_M(E_\wp)\times\H_M(E_{\wp'})$).

\end{itemize}

Let $\pi_p$ be an unramified $\theta$-stable representation of $\M^0(\Q_p)$
and $\varphi_{\pi_p}:W_{\Q_p}\fl{}^L\M^0_{\Q_p}$ be a Langlands parameter of
$\pi_p$. 
As $\pi_p$ is $\theta$-stable, we may assume that $\varphi_{\pi_p}$
factors through the image of ${}^L\H_{M,\Q_p}\fl{}^L\M^0_{\Q_p}$.
Let $\varphi_\wp$ be the morphism $W_{E_\wp}\fl{}^L\H_{M,E_\wp}$ deduced from
$\varphi_{\pi_p}$. If $p$ is inert in $E$, then $\H_M(E_\wp)=\M^0(\Q_p)$,
and $\varphi_\wp$ is a Langlands parameter of $\pi_p$, seen as a
representation of $\H_M(E_\wp)$. If $p$ splits in $E$ and $\wp'$ is the
second place of $E$ above $p$, then $\M^0(\Q_p)=\H_M(E_\wp)\times
\H_M(E_{\wp'})$, so $\pi_p=\pi\otimes\pi'$, where $\pi$ (resp. $\pi'$)
is an unramified representation of $\H_M(E_\wp)$ (resp.
$\H_M(E_{\wp'})$). The morphism $\varphi_\wp$ is a Langlands parameter of
$\pi$.
By theorem \ref{th:transformee_de_Satake} and lemma
\ref{lemme:rep_nr_tordues}, if $A_{\pi_p}$ is a normalized intertwining
operator on $\pi_p$, then
\[\Tr(\pi_p(\phi_{\mu_{H_M},p}^{(j)})A_{\pi_p})=(N\wp)^{jd/2}\Tr(r_{-\mu_{H_M}}
\circ\varphi_\wp(\Phi_\wp^j))\Tr(\pi_p(\ungras_{\M^0(\Z_p)})A_{\pi_p}).\]

Set
\[\phi_{M,p}^{(j)'}=\sum_{\mu_{H_M}\in M_{\H_M}}(-1)^{s(\mu_{H_M})}\phi_{
\mu_{H_M},p}^{(j)}.\]
The calculations of \ref{partie_en_p2} and (iii) of lemma
\ref{lemme:hypothese_H}
imply that the function $\chi_{\eta_{H_M},p}^{-1}f_{H_M,p}^{(j)}\in\Hecke
(\H_M(\Q_p),\H_M(\Z_p))$ is the image by the base change map of the function
$\phi_{M,p}^{(j)'}$ (as before, $f_{H_M,p}^{(j)}$ is the function of
definition \ref{def:phi_m_G}). So we can take
$\phi_{M,p}^{(j)}=\chi_{M,p}\phi_{M,p}^{(j)'}$ in theorem
\ref{th:stab_FT_IC_GL_n}, and the corollary follows from this theorem and
from proposition \ref{prop:FT_invariante_tordue}.

\end{proof}

Another application of theorem \ref{th:stab_FT_IC_GL_n} is the next corollary.
In this corollary, $E$ is still an imaginary quadratic extension of $\Q$.
Fix $n\in\Nat^*$, and let $\theta$ be the involution $(\lambda,g)\fle
(\overline{\lambda},\overline{\lambda}{}^t\overline{g}^{-1})$ of $R_{E/\Q}
(\Gr_{m,E}\times\GL_{n,E})$ (where $(\lambda,g)\fle(\overline{\lambda},
\overline{g})$ is the action of the non-trivial element of $\Gal(E/\Q)$). Then
$\theta$ defines an involution of $\C^\times\times\GL_n(\C)=(R_{E/\Q}(\Gr_{m,E}
\times\GL_{n,E}))(\R)$. 
The morphism  $\theta$ is the involution induced by the non-trivial
element of $\Gal(E/\Q)$, if $\Gr_{m,E}\times\GL_{n,E}$ is identified to
$\GU(n)_E$. 
If $\wp$ is a finite unramified place of $E$ and $\pi_\wp$ is an unramified
representation of $\GU(n)(E_\wp)$, let $\log L(s,\pi_\wp)$ (resp.
$\log L(s,\pi_\wp,\bigwedge\limits^2)$) be the logarithm of the
$L$-function of $\pi_\wp$ and of the representation $id_{\C^\times}\otimes
st$ (resp. $id_{\C^\times}\otimes\bigwedge\limits^2st$) of $\widehat{\GU(n)}=
\C^\times\times\GL_n(\C)$, where $st$ is the standard representation of
$\GL_n(\C)$.

\begin{subcorollaire}\label{cor:G_m_GL_n_E} 
Let $\pi$ be a $\theta$-stable cuspidal automorphic representation
of $\Ade_E^\times\times\GL_n(\Ade_E)$ such that $\pi_{\infty}$ is tempered
(where $\infty$ is the unique infinite place of $E$).
Assume that there exists an irreducible algebraic representation $V$
of $\GU(n)$ such that
$ep(\theta,\pi_{\infty}\otimes W)\not=0$, where $W$ is the
$\theta$-stable representation of $\C^\times\times\GL_n(\C)$ associated to
$V$ (cf theorem \ref{th:pseudo_coeff_tordus}).
Let $m$ be the weight of $V$ in the sense of \ref{points_fixes3}
(ie, the relative integer such that the central subgroup $\Gr_m$ of
$\GU(n)$ acts on $V$ by $x\fle x^m$).
Let $S$ be the union of the set of prime numbers that ramify in $E$ and of the
set of prime numbers under finite places of $E$ where $\pi$
is ramified.
Then there exists a number field $K$, a positive integer $N$ and, for
every finite place $\lambda$ of $K$, a continuous
finite-dimensional representation $\sigma_\lambda$ of
$\Gal(\overline{\Q}/E)$ with coefficients in $K_\lambda$, such that :
\begin{itemize}
\item[(i)] The representation $\sigma_\lambda$ is unramified outside of
$S\cup\{\ell\}$, where $\ell$ is the prime number under $\lambda$,
pure of weight $-m+1-n$ if $n$ is not dividible by $4$ and mixed with weights
between $-m+2(2-n)-1$ and $-m+2(2-n)+1$ if $n$ is dividible by $4$.
If $n$ is dividible by $4$ and the highest weight of $V$ is regular,
then $\sigma_\lambda$ is pure of weight $-m+2(2-n)$.
\item[(ii)] For every place $\wp$ of $E$ above a prime number $p\not\in S$,
for every finite place $\lambda\not|p$ of $K$,
\[\log L_\wp(s,\sigma_\lambda)=N\log L(s+\frac{n-1}{2},\pi_{\wp})\]
if $n$ is not dividible by $4$, and
\[\log L_\wp(s,\sigma_\lambda)=N\log L(s+(n-2),\pi_{\wp}
,\bigwedge^2)\]
if $n$ is dividible by $4$ (where $\pi_{\wp}$ is the local component at $\wp$
of $\pi$, seen as a representation of $\GU(n)(\Ade_E)$).

\end{itemize}

\end{subcorollaire}

\begin{proof} We can, without changing the properties of $\pi$, replace
$\theta$ by its product with an inner automorphism of
$R_{E/\Q}(\Gr_{m,E}\times\GL_{n,E})$.
So we may (and will) assume that
$\theta(\lambda,g)=(\overline{\lambda},\overline{\lambda}J_{p_1,q_1}{}^t
\overline{g}^{-1}
J_{p_1,q_1}^{-1})$, where $p_1,q_1\in\Nat^*$ are such that $p_1+q_1=n$ and
$J_{p_1,q_1}\in\GL_n(\Z)$ is the matrix (defined in \ref{groupes1})
of the Hermitian form that gives the group $\GU(p_1,q_1)$. 

Write as before $\H=\GU(p_1,q_1)$ and $\G^0=R_{E/\Q}\H_E$. Then
$\G^0=R_{E/\Q}\Gr_{m,E}\times R_{E/\Q}\GL_{n,E}$, and the involution $\theta$
of $\G^0$ defined above is equal to the involution induced by the non-trivial
element of $\Gal(E/\Q)$.
Assume that the group $\H$ is quasi-split (but not necessarily unramified)
at every finite place of $\Q$.
As $ep(\theta,\pi_\infty\otimes W)\not=0$
and $\pi_\infty$ is tempered, remark \ref{rq:ep_tordu} and theorem
\ref{th:pseudo_coeff_tordus} imply that $\pi_\infty=\pi_{V^*}$, 
where $\pi_{V^*}$ is the $\theta$-discrete representation of
$\G^0(\R)$ associated to $V^*$ as in lemma
\ref{lemme:O_pseudo_coeff}.

Let $\K_S\subset\G^0(\Ade_S)$ be an open
compact subgroup such that $\Tr(\pi_S(\ungras_{\K_S})A_{\pi_S})\not=0$ (where
$A_{\pi_S}$ is any intertwining operator on $\pi_S$).
By lemma
\ref{lemme:hypothese_H}, by taking $\K_S$ small enough, we may assume that
there exists a function $f_S\in C_c^\infty(\H(\Ade_S))$ associated to
$\phi_S:=\ungras_{\K_S}$ and satisfying condition (H). For every
$\M\in\Mcal_G$, fix a transfer $f_S^{H_M}$ of $f_S$ to $\H_M$ and a function
$\phi_{M,S}\in C_c^\infty(\M^0(\Ade_S))$ associated to $f_S^{H_M}$.
If $p\not\in S$, $\M\in\Mcal_G$ and $\phi_p\in\Hecke(\G^0(\Q_p),\G^0(\Z_p))$,
let $b(\phi_p)\in\Hecke(\H^0(\Q_p),\H^0(\Z_p))$, $b(\phi_p)^{H_M}\in\Hecke
(\H_M(\Q_p),\H_M(\Z_p))$ and $\phi_{M,p}\in\Hecke(\M^0(\Q_p),\M^0(\Z_p))$
be the functions obtained from $\phi_p$ by following the arrows of the
commutative diagram of point (iii) of lemma \ref{lemme:hypothese_H}. 
Finally, for every $\M\in\Mcal_G$, define a function $\phi_{M,\infty}
\in C^\infty(\M^0(\R))$ from $V$ as in theorem \ref{th:stab_FT_IC_GL_n}.

Let $\K_{H,S}\subset\H(\Ade_S)$ be an open compact subgroup
small enough for $f_S$ to be bi-invariant under $\K_{H,S}$. Set
$\K_H=\K_{H,S}\prod\limits_{p\not\in S}\H(\Z_p)$; then
$\K_H$ is an open compact subgroup of $\H(\Af)$, and we may assume
(by making $\K_{H,S}$ smaller) that $\K_H$ is neat.
Then the results of \ref{points_fixes7} apply to $\H$, $\K_H$ and
every $p\not\in S$.
In the beginning of \ref{applications1}, we explained how to get a number field
$K$ and, for every finite place $\lambda$ of $K$, a virtual finite-dimensional
$\lambda$-adic representation $W_\lambda$ of
$\Gal(\overline{\Q}/E)\times\Hecke(\H(\Af),\K_H)$
(the cohomology of the complex $IC^{\K_{H}}V$)
such that there is a decomposition
\[W_\lambda=\bigoplus_{\pi_{H,f}}W_\lambda(\pi_{H,f})\otimes\pi_{H,f}^{\K_H},\]
where the direct sum is taken over the set of isomorphism classes of
irreducible admissible representations $\pi_{H,f}$ of $\H(\Af)$ such that
$\pi_{H,f}^{\K_H}\not=0$, and the $W_\lambda(\pi_{H,f})$ are virtual
$\lambda$-adic representations of $\Gal(\overline{\Q}/E)$.

Let $\Pi_H(\pi_f)$ be the set of isomorphism classes of irreducible
admissible representations $\pi_{H,f}$ of $\H(\Af)$ such that
$\pi_{H,f}^{\K_H}\not=0$ (so $\pi_{H,f}$ is unramified outside of
$S$), that $W_\lambda(\pi_{H,f})\not=0$ and that, for every
$p\not\in S$, if $\varphi_{\pi_{H,p}}:W_{\Q_p}\fl{}^L\H_{\Q_p}$ is a
Langlands parameter of $\pi_{H,p}$, then the composition of
$\varphi_{\pi_{H,p}}$ and of the inclusion ${}^L\H_{\Q_p}\fl{}^L\G^0_{\Q_p}$
(defined in example \ref{ex:groupes_non_connexes}) is a Langlands parameter
of $\pi_p$. 

By the multiplicity $1$ theorem of Piatetski-Shapiro,
$m_{disc}(\pi)=1$. So there is a normalized intertwining operator
on $\pi$ such that $m^+_{disc}(\pi)=1$ and $m^-_{disc}(\pi)=0$;
denote this intertwining operator by $A_\pi$.
Let as before $\mu_H:\Gr_{m,E}\fl\H_E$ be the cocharacter defined by the
Shimura datum, $r_{-\mu_H}$ be the representation of ${}^L\H_E$
determined by $-\mu_H$ (cf \ref{partie_en_p1}) and $d=p_1q_1$. Set
$\phi^\infty=\phi_S\prod\limits_{p\not\in S}\ungras_{\G^0(\Z_p)}$ and
$f^\infty=f_S\prod\limits_{p\not\in S}\ungras_{\H(\Z_p)}$.
Let $\wp$ be a finite place of $E$ above a prime number $p\not\in S$.
Let $\pi_\wp$ be the local component at $\wp$ of $\pi$ (seen as a
representation of $\H(\Ade_E)$) and $\varphi_\wp:W_{E_\wp}\fl{}^L\H_{E_\wp}$
be a Langlands parameter of $\pi_\wp$. We are going to
show that, for every finite place
$\lambda\not|p$ of $K$ and for every $j\in\Z$,
\renewcommand\theequation{$*$}
\begin{equation}
c_G(N\wp)^{dj/2}\Tr(r_{-\mu_H}\circ\varphi_\wp(\Phi_\wp^j))\Tr(\pi_f(\phi^
\infty)A_\pi)=\sum_{\pi_{H,f}\in \Pi_H(\pi_f)}\Tr(\pi_{H,f}(f^\infty))\Tr(\Phi_
\wp^j,W_\lambda(\pi_{H,f})),
\end{equation}
where $\Phi_\wp\in W_{E_\wp}$ is a lift of the geometric Frobenius.
It suffices to show this equality for $j>0$.

Let $\lambda$ be a finite place of $K$ such that $\lambda\not|p$.
Let $j\in\Nat^*$. Define a function $\phi_p^{(j)}\in
\Hecke(\G^0(\Q_p),\G^0(\Z_p))$ using $\mu_H$ and $j$, as in the proof of
corollary \ref{cor:calcul_L_GL_n}. We recall the definition. Let $L$ be an
unramified extension of $E_\wp$ of degree $j$. Then $\phi_p^{(j)}$ is the
image of the function $f_{\mu_H,L}\in\Hecke(\H(L),\H(\Of_L))$ determined by
$\mu_H$ as in \ref{partie_en_p1} by the morphism
\[\Hecke(\H(L),\H(\Of_L))\fl\Hecke(\H(E_\wp),\H(\Of_{E_\wp}))\fl
\Hecke(\G^0(\Q_p),\G^0(\Z_p)),\]
where the first arrow is the base change morphism and the second arrow is
\begin{itemize}
\item[-] identity if $p$ is inert in $E$ (so $\G^0(\Q_p)=\H(E_\wp)$);
\item[-] the morphism $h\fle (h,\ungras_{\H(\Of_{E_{\wp'}})})$ if $p$
splits in $E$ and $\wp'$ is the second place of $E$ above $p$ (so
$\G^0(\Q_p)=\H(E_\wp)\times\H(E_{\wp'})$).
\end{itemize}
Then $b(\phi_p^{(j)})$ is the function $f_{H,p}$ defined after theorem
\ref{th:stab_partie_elliptique}. Moreover, by theorem
\ref{th:transformee_de_Satake} and lemma
\ref{lemme:rep_nr_tordues}, if $\pi_p$ is the local component at $p$ of $\pi$
(seen as a representation of $\G^0(\Ade)$) and $A_{\pi_p}$ is a normalized
intertwining operator on $\pi_p$, then :
\[\Tr(\pi_p(\phi_{p}^{(j)})A_{\pi_p})=(N\wp)^{jd/2}\Tr(r_{-\mu_H}
\circ\varphi_\wp(\Phi_\wp^j))\Tr(\pi_p(\ungras_{\G^0(\Z_p)})A_{\pi_p}).\]

Let $\M\in\Mcal_\G$.
Let $R_M$ be the set of $\pi_M\in\Pi_{disc}(\M,t)$, with $t\geq 0$, such that :
\begin{itemize}
\item[(i)] $\pi_M\otimes\chi_M$
is unramified at every finite place $v\not\in S$;
\item[(ii)] $a_{disc}^M(\pi_M)\not=0$;
\item[(iii)] $\Tr(\pi_{M,S}(\phi_{M,S})A_{\pi_M})\not=0$ and
$\Tr(\pi_{M,\infty}(\phi_{M,\infty})A_{\pi_M})\not=0$ (where $A_{\pi_M}$ is
a normalized intertwining operator on $\pi_M$);
\item[(iv)] if $\M=\G$, then $\pi_M\not\simeq\pi$.
\end{itemize}
Then $R_M$ is finite.

Let $R_H$ be the set of isomorphism classes
of irreducible admissible representations $\pi_{H,f}$ of
$\H(\Af)$ such that :
\begin{itemize}
\item[(i)] $\pi_{H,f}^{\K_H}\not=0$;
\item[(ii)] $\pi_{H,f}\not\in\Pi_H(\pi)$;
\item[(iii)] $W_\lambda(\pi_{H,f})\not=0$.

\end{itemize}
Then $R_H$ is also finite.

By the strong multiplicity $1$ theorem of Jacquet-Shalika for $\G^0$ 
(cf theorem 4.4 of \cite{JS}) and corollary
\ref{cor:CB} (cf also remark \ref{rq:CB}), there exists a function
$g^{S\cup\{p\}}\in\Hecke(\G^0(\Af^{S\cup\{p\}}),\K_G^{S\cup\{p\}})$
(where $\K_G^{S\cup\{p\}}=\prod\limits_{v\not\in S\cup\{p\}}\G^0(\Z_v)$)
such that : 
\begin{bulletlist}
\item $\Tr(\pi^{S\cup\{p\}}(g^{S\cup\{p\}})A_\pi)=\Tr(\pi^{S\cup\{p\}}(\phi
^{S\cup\{p\}})A_\pi)=1;$
\item for every $\pi_{H,f}\in R_H$,
$\Tr(\pi_{H,f}^{S\cup\{p\}}(b(g^{S\cup\{p\}})))=0$ (where
$b(g^{S\cup\{p\}})\in\Hecke(\H(\Af^{S\cup\{p\}}),\K_H^{S\cup\{p\}})$ is
the base change of $g^{S\cup\{p\}}$);
\item for every $\pi_{H,f}\in\Pi_H(\pi)$,
$\Tr(\pi_{H,f}^{S\cup\{p\}}(b(g^{S\cup\{p\}})))=\Tr({\pi_{H,f}}^{S\cup\{p\}}
(f^{S\cup\{p\}}))$ (this actually follows from the first condition and
from the fundamental lemma for base change);
\item for every $\M\in\Mcal_\G$, every $\pi_M\in R_M$ and every
normalized intertwining operator $A_{\pi_M}$ on $\pi_M$,
$\Tr(\pi_M^{S\cup\{p\}}(g^{S\cup\{p\}}_M)A_{\pi_M})=0$
(where $g^{S\cup\{p\}}_M$ is the function obtained from $g^{S\cup\{p\}}$
by following the left vertical arrow in the diagram of (iii) of
lemma \ref{lemme:hypothese_H}).

\end{bulletlist}
Then
\[\Tr(\pi_f(\phi_S\phi^{S\cup\{p\}}\phi_p^{(j)})A_\pi)=\Tr(\pi_f(\phi^\infty)
A_\pi)(N\wp)^{dj/2}\Tr(r_{-\mu_H}\circ\varphi_{\wp}(\Phi_\wp^j)),\]
and, by theorem \ref{th:stab_FT_IC_GL_n} (and the fact that
$\Tr(\pi_\infty(\phi_{G,\infty}))=1$),
\[c_G\Tr(\pi_f(\phi_S\phi^{S\cup\{p\}}\phi_p^{(j)})A_\pi)=
\sum_{\pi_{H,f}\in \Pi_H(\pi_f)}
\Tr(\Phi_\wp^j\times f_Sf^{S\cup\{p\}}\ungras_{\H(\Z_p)},W_\lambda
(\pi_{H,f})\otimes\pi_{H,f}^{\K_H}).\]
This proves equality $(*)$.

Remember that we wanted the group $\GU(p_1,q_1)$ to be quasi-split at
every finite place of $\Q$. We use here the calculations of the Galois
cohomology of unitary groups of section 2 of \cite{Cl-RGRA}.
If $n$ is odd, these calculations imply that the group $\GU(p_1,q_1)$ is
quasi-split at every finite place of $\Q$ for any $p_1$ and $q_1$.
Take $p_1=1$ and $q_1=n-1$. Now assume that $n$ is even. If $n/2$ is odd,
take $p_1=1$ and $q_1=n-1$. If $n/2$ is even, take $p_1=2$ and $q_1=n-2$.
We check that, with these choices, $\GU(p_1,q_1)$ is indeed quasi-split at
every finite place of $\Q$. Let $D$ be the discriminant of $E$. Let $q$ be
a prime number. If $q$ does not divide $D$, then $\GU(p_1,q_1)$ is
unramified at $q$ (so, in particular, it is quasi-split). Assume that
$q$ divides $D$. Then the cohomological invariant of $\GU(p_1,q_1)$ at $q$ is
$0$ if $-1$ is a norm in $\Q_q$, and $q_1+n/2\mod 2$ otherwise. But, by the
choice of $q_1$, $q_1+n/2$ is always even, so $\GU(p_1,q_1)$ is quasi-split
at $q$.

In the rest of proof, take $p_1$ and $q_1$ as in the discussion above.
Note that $d=n-1$ if $n$ is not dividible by $4$, and $d=2(n-2)$ if $n$ is
dividible by $4$.

We now apply lemma \ref{lemme:conj_RP1}. As $\H$ splits over $E$, the
representation $r_{-\mu_H}$ of ${}^L\H_E=\widehat{\H}\times W_E$ determined
by the cocharacter $\mu_H$ of $\H_E$ is trivial on $W_E$. Let $st^\vee$ be
the contragredient of the standard representation of $\GL_n(\C)$ and
$\chi$ be the character $z\fle z^{-1}$ of $\C^\times$.
By lemma \ref{lemme:conj_RP1}, the restriction of $r_{-\mu_H}$ to
$\widehat{\H}=\C^\times\times\GL_n(\C)$ is
$\chi\otimes st^\vee$ if $n$ is not dividible by $4$,
and $\chi\otimes\bigwedge\limits^2 st^\vee$
if $n$ is dividible by $4$.

Let $\wp$ be a finite place of $E$ above a prime number
$p\not\in S$, and $\lambda\not|p$ be a finite place of $K$.
Fix a Langlands parameter
$(z,(z_1,\dots,z_n))$ of $\pi_\wp$ in the maximal torus
$\C^\times\times(\C^\times)^n$ of $\widehat{\H}=\C^\times\times\GL_n(\C)$.
By reasoning as in the beginning of the proof of theorem
\ref{th:conj_RP} (or by applying corollary \ref{cor:CB} and theorem
\ref{th:conj_RP}), we see that $\log_{N\wp}|z|\in\frac{1}{2}\Z$.
For every $\pi_{H,f}\in \Pi_H(\pi_f)$, let $a_i$, $i\in I_{\pi_{H,f}}$, be the
eigenvalues of $\Phi_\wp$ acting on $W_\lambda(\pi_{H,f})$, $n_i\in\Z$,
$i\in I_{\pi_{H,f}}$, be their multiplicities, and
\[b_{\pi_{H,f}}=c_G^{-1}\Tr(\pi_f(\phi^\infty)A_\pi)^{-1}\Tr(\pi_{H,f}(f^\infty
))\]
($b_{\pi_{H,f}}$ does not depend on $\wp$). By equality
$(*)$, for every $j\in\Z$ :
\renewcommand\theequation{$**$}
\begin{equation}
(N\wp)^{dj/2}z^{-j}\sum_{{J\subset\{1,\dots,n\}}\atop{|J|=k}}\prod_{l\in J}
z_l^{-j}=\sum_{\pi_{H,f}\in \Pi_H(\pi_f)}
b_{\pi_{H,f}}\sum_{i\in I_{\pi_{H,f}}}n_i a_i^j,
\end{equation}
where $k=1$ if $n$ is not dividible by $4$, and $k=2$ if $n$ is dividible
by $4$.
So there exists a positive integer $N$ (independent from $\wp$)
such that $N_{\pi_{H,f}}:=Nb_{\pi_{H,f}}\in\Z$, for every $\pi_{H,f}\in
\Pi_H(\pi_f)$. Moreover, for every $\pi_{H,f}\in\Pi_H(\pi_f)$ and $i\in I_{\pi_
{H,f}}$, the product $N_{\pi_{H,f}}n_i$ is positive.
In particular,
\[\sigma_\lambda:=\bigoplus_{\pi_{H,f}\in \Pi_H(\pi_f)}N_{\pi_{H,f}}W_\lambda
(\pi_{H,f})^{\vee}
\]
is a real representation of $\Gal(\overline{\Q}/E)$ (and not just a virtual
representation).
Then equality $(**)$ becomes : for every finite place $\wp$ of $E$ above a
prime number $p\not\in S$, if $\lambda\not|p$, then, for every $j\in\Z$,
\[N(N\wp)^{-jd/2}\Tr((id_{\C^\times}\otimes\bigwedge^kst)
(\varphi_\wp(\Phi_\wp^j))=\Tr(\Phi_\wp^j,
\sigma_\lambda),\]
where $k$ is as before equal to $1$ if $n$ is not dividible by $4$, and
to $2$ if $n$ is dividible by $4$.
This is point (ii) of the lemma (as $\GU(n)$ is an inner form of
$\GU(p_1,q_1)$, we can see $\varphi_\wp$ as the Langlands parameter of the
local component at $\wp$ of $\pi$, seen as a representation of
$\GU(n)(\Ade_E)$).

It remains to determine the weight of $\sigma_\lambda$. As the algebraic
representation $V$ of $\GU(n)$ is pure of weight $m$ in the sense of
\ref{points_fixes3},
the complex $IC^{\K}V$ is pure of weight $-m$. Let
\[W_\lambda=\sum_{i=0}^{2d}(-1)^iW_\lambda^i\]
be the decomposition of $W_\lambda$ according to cohomology degree.
For every irreducible admissible representation $\pi_{H,f}$ of
$\H(\Af)$ such that $\pi_{H,f}^{\K_{H,f}}\not=0$, there is a decomposition
\[W_\lambda(\pi_{H,f})=\sum_{i=0}^{2d}(-1)^iW_\lambda^i(\pi_{H,f}),\]
and the representation $W_\lambda^i(\pi_{H,f})$ of $\Gal(\overline{\Q}/E)$
is pure of weight $-m+i-2d$.

Remember that $(z,(z_1,\dots,z_n))$ is the Langlands parameter of $\pi_\wp$.
Assume first that $n$ is not dividible by $4$. Then equality $(**)$ implies
that $\log_{N\wp}|z_i|\in\frac{1}{2}\Z$ for every
$i\in\{1,\dots,n\}$ (because the $a_i$ and $z$ satisfy the same property).
But we know that $-\frac{1+m}{2}<\log_{N\wp}|z_i|<\frac{1-m}{2}$ for every
$i\in\{1,\dots,n\}$ (cf \cite{Cl-MFA} lemma 4.10; note that the conditions
on $\pi_\infty$ imply that $\pi$ is algebraic regular in the sense of
\cite{Cl-MFA}, and that Clozel uses a different normalization of the
Langlands parameter at $\wp$), so $\log_{N\wp}|z_i|=-\frac{m}{2}$
for every $i$.
This implies that,
if $\pi_{H,f}\in\Pi_H(\pi_f)$, then $W^i(\pi_{H,f})=0$ for every $i\not=d$.
Hence $\sigma_\lambda$ is pure of weight $-m-d=-m+1-n$.

Assume now that $n$ is dividible by $4$. Equality $(**)$ implies only that
$\log_{N\wp}|z_i|\in\frac{1}{4}\Z$.  As before, we know that the
$\log_{N\wp}|z_i|$ are in $]-\frac{1+m}{2},\frac{1-m}{2}[$, so $\log_p|z_i|
\in\{-\frac{1+m}{4},-\frac{m}{4},\frac{1-m}{4}\}$
for every $i\in\{1,\dots,n\}$.
Applying $(**)$ again, we see that the only $W^i(\pi_{H,f})$ that can appear
in $\sigma_\lambda$ are those with $d-1\leq i\leq d+1$. This proves the
bounds on the weights of $\sigma_\lambda$. Assume that the highest weight
of $V$ is regular. Then, by
lemma \ref{lemme:conj_RP4}, $W^i=0$ if $i\not=d$, so $\sigma_\lambda$ is
of weight $-m-d=-m+2(2-n)$.

\end{proof}

To formulate the last two corollaries, we will use the following definition
of Clozel :

\begin{subdefinition}\label{def:algebrique_reguliere}(cf \cite{Cl-MFA}
1.2.3 or \cite{Cl-RGRA} 3.1)
Let $\pi$ be a cuspidal automorphic representation of
$\GL_n(\Ade)$ (resp. $\GL_n(\Ade_E)$, where $E$ is an imaginary quadratic
extension of $\Q$). Then $\pi$ is called \emph{algebraic}
\index{algebraic (for an automorphic representation)}
if there exist a Langlands parameter $\varphi:W_\R\fl\GL_n(\C)$ (resp.
$\varphi:W_\C\fl\GL_n(\C)$) of
$\pi_\infty$ and $p_1,\dots,p_n,q_1,\dots,q_n\in\Z$ such that, for every
$z\in W_\C=\C^\times$,
\[\varphi(z)=\left(\begin{array}{ccc}z^{p_1+\frac{n-1}{2}}\overline{z}^{q_1+
\frac{n-1}{2}} & & 0 \\
& \ddots & \\0 & & z^{p_n+\frac{n-1}{2}}\overline{z}^{q_n+\frac{n-1}{2}}
\end{array}\right).\]
We may assume that $p_1\geq\dots\geq p_n$. The representation $\pi$ is
called \emph{regular algebraic} if $p_1>\dots>p_n$.
\index{regular algebraic (for an automorphic representation)}

If $\pi$ is regular algebraic, then there is an algebraic representation
$W$ of $\GL_n$ (resp. $R_{E/\Q}\GL_{n,E}$) associated to $\pi$ as in
\cite{Cl-MFA} 3.5 and \cite{Cl-RGRA} 3.2 : the highest weight of $W$ is
$(p_1,p_2+1,\dots,p_n+(n-1))$ (resp. $((p_1,p_2+1,\dots,p_n+(n-1)),(q_n,
q_{n-1}+1,\dots,q_1+(n-1)))$). We say that $\pi$ is \emph{very regular} if
the highest weight of $W$ is regular.
\index{very regular (for an automorphic representation)}

\end{subdefinition}

We summarize a few results of Clozel about regular algebraic representations
in the next lemma.

\begin{sublemme}\label{lemme:algebrique_reguliere} Let $E$ be an imaginary
quadratic extension of $\Q$. Let $\G^0=\GL_n$ or $R_{E/\Q}\GL_{n,E}$.
\begin{itemize}
\item[(i)] Let $\pi$ be a cuspidal automorphic representation of $\G^0(\Ade)$.
Then the following conditions are equivalent :
\begin{itemize}
\item[(a)] $\pi$ is regular algebraic.
\item[(b)] The infinitesimal character of $\pi_\infty$ is that of an
algebraic representation of $\G^0$.
\item[(c)] There exists an algebraic representation $W$ of $\G^0$
and a character $\varepsilon$ of $\G^0(\R)$ of order $2$ such that
such that $\Ho^*(\ggoth,\K'_\infty;\varepsilon(\pi_\infty\otimes W^*))
\not=0$, where $\ggoth=Lie(\G^0(\C))$ and $\K'_\infty$ is the set of
fixed points of a Cartan involution of $\G^0(\R)$.

\end{itemize}
Moreover, if $\pi$ is regular algebraic, then $\pi_\infty$ is essentially
tempered.

\item[(ii)] Assume that $\G^0=R_{E/\Q}\GL_{n,E}$. Let $\theta$ be the
involution of $\G^0$ defined by $g\fle {}^t\overline{g}^{-1}$. Let $\pi$
be a $\theta$-stable cuspidal automorphic representation of $\G^0(\Ade)$.
Then the following conditions are equivalent :
\begin{itemize}
\item[(a)] $\pi$ is regular algebraic.
\item[(b)] There exists a $\theta$-stable algebraic representation $W$ of
$\G^0$ such that $ep(\theta,\pi_\infty\otimes W^*)\not=0$ (where
$ep(\theta,.)$ is defined before theorem \ref{th:pseudo_coeff_tordus}).
\end{itemize}
Moreover, if $\pi$ is regular algebraic, then $\pi_\infty$ is tempered.

\end{itemize}

\end{sublemme}

\begin{proof}
\begin{itemize}
\item[(i)] The equivalence of (a) and (b) is obvious from the definition,
and (c) implies (b) by Wigner's lemma. The fact that (a) implies (c) is
lemma 3.14 of \cite{Cl-MFA} (cf also proposition 3.5 of \cite{Cl-RGRA} for
the case $\G^0=R_{E/\Q}\GL_{n,E}$).
The last sentence of (i) is lemma 4.19 of \cite{Cl-MFA}.

\item[(ii)] It is obvious (by Wigner's lemma and (i)) that (b) implies (a).
The fact that (a) implies (b) is proved in proposition 3.5 of \cite{Cl-RGRA}.
The last sentence is a remark made at the beginning of 3.2 of
\cite{Cl-RGRA} : if $\pi$ is regular algebraic, then $\pi_\infty$ is
essentially tempered; as it is also $\theta$-stable, it must be tempered.

\end{itemize}
\end{proof}

As in the lemma above, denote by $\theta$ the automorphism $g\fle{}^t
\overline{g}^{-1}$ of $R_{E/\Q}\GL_{n,E}$.

\begin{subcorollaire}\label{cor:GL_n_E} 
Let $\pi$ be a $\theta$-stable cuspidal automorphic representation
of $\GL_n(\Ade_E)$ that is regular algebraic.
Let $S$ be the union of the set of prime numbers that ramify in $E$ and of the
set of prime numbers under finite places of $E$ where $\pi$
is ramified.
Then there exists a number field $K$, a positive integer $N$ and, for
every finite place $\lambda$ of $K$, a continuous
finite-dimensional representation $\sigma_\lambda$ of
$\Gal(\overline{\Q}/E)$ with coefficients in $K_\lambda$, such that :
\begin{itemize}
\item[(i)] The representation $\sigma_\lambda$ is unramified outside of
$S\cup\{\ell\}$, where $\ell$ is the prime number under $\lambda$,
pure of weight $1-n$ if $n$ is not dividible by $4$ and mixed with weights
between $2(2-n)-1$ and $2(2-n)+1$ if $n$ is dividible by $4$.
If $n$ is dividible by $4$ and $\pi$ is very regular,
then $\sigma_\lambda$ is pure of weight $2(2-n)$.
\item[(ii)] For every place $\wp$ of $E$ above a prime number $p\not\in S$,
for every finite place $\lambda\not|p$ of $K$,
\[\log L_\wp(s,\sigma_\lambda)=N\log L(s+\frac{n-1}{2},\pi_{\wp})\]
if $n$ is not dividible by $4$, and
\[\log L_\wp(s,\sigma_\lambda)=N\log L(s+(n-2),\pi_{\wp}
,\bigwedge^2)\]
if $n$ is dividible by $4$ (where $\pi_{\wp}$ is the local component at $\wp$
of $\pi$, seen as a representation of $\U(n)(\Ade_E)$).

\end{itemize}

(The first $L$-function is the one associated to the standard representation of
$\GL_n(\C)=\widehat{\U(n)}$, and the second $L$-function is the one associated
to the exterior square of the standard representation.)

\end{subcorollaire}

\begin{proof} It is enough to show that there exists a character
$\chi:\Ade_E/E^\times\fl\C^\times$ such that $\chi\otimes\pi$ satisfies the
conditions of corollary \ref{cor:G_m_GL_n_E}, with $V$ of weight $0$. This
follows from lemma VI.2.10 of \cite{HT}.

\end{proof}

Using the base change of Arthur and Clozel (\cite{AC}), it is possible
to deduce from the last corollary results about self-dual automorphic
representations of $\GL_n(\Ade)$. Here, we will treat only the case $n$
odd (which is simpler); in the general case, the next corollary would
not hold for all quadratic imaginary extensions $E$.

\begin{subcorollaire}\label{cor:GL_n_Q} Assume that $n$ is odd.
Let $\tau$ be a self-dual cuspidal automorphic representation of $\GL_n(\Ade)$,
and assume that $\tau$ is regular algebraic.
Let $E$ be a quadratic imaginary extension of $\Q$.
Write $S$ for the union of the set of prime numbers that ramify in $E$ and
the set of prime numbers where $\tau$ is ramified.
Then there exist a number field $K$, a positive integer $N$ and,
for every finite place $\lambda$ of $K$, a (continuous finite-dimensional)
representation $\sigma_\lambda$ of $\Gal(\overline{\Q}/E)$ with coefficients
in $K_\lambda$, such that :
\begin{itemize}
\item[(i)] the representations $\sigma_\lambda$ are unramified outside of $S$
and pure of weight $1-n$;
\item[(ii)] for every finite place $\wp$ of $E$ above a prime number
$p\not\in S$, for every finite place $\lambda\not|p$ of $K$, for every
$j\in\Z$,
\[\Tr(\sigma_\lambda(\Phi_\wp^j))=N(N\wp)^{j(n-1)/2}\Tr(\varphi_{\tau_p}(
\Phi_\wp^j)),\]
where $\varphi_{\tau_p}:W_{\Q_p}\fl\GL_n(\C)$ is a Langlands parameter of
$\tau_p$ and $\Phi_\wp\in W_{E_\wp}$ is a lift of the geometric Frobenius.

\end{itemize}
In particular, $\tau$ satisfies the Ramanujan-Petersson conjecture at
every unramified place.

\end{subcorollaire}

\begin{proof} Let $\theta$ be as before
the involution $g\fle{}^t\overline{g}^{-1}$ of
$R_{E/\Q}\GL_{n,E}$.
If $V$ is an irreducible algebraic representation of $\GL_n$,
it defines a $\theta$-discrete representation $\pi_V$ of
$\GL_n(E\otimes_\Q\R)$ as in lemma \ref{lemme:O_pseudo_coeff}.

Let $\pi$ be the automorphic representation of $\GL_n(\Ade_E)$ obtained from
$\tau$ by base change (cf \cite{AC} theorem III.4.2). Because $n$ is odd,
$\pi$ is necessarily cuspidal (this follows from (b) of loc. cit.).
\footnote{I think Sug Woo Shin for pointing out this useful fact to me.}
By the definition of base change, $\pi$ is regular algebraic.
Let $(p_1,\dots,p_n)\in\Z^n$, with $p_1\geq\dots\geq p_n$, be the $n$-uple
of integers associated to $\tau$ as in definition
\ref{def:algebrique_reguliere}.
As $\tau$ is self-dual, $p_i+p_{n+1-i}=1-n$ for every $i\in\{1,\dots,n\}$.
For every $i\in\{1,\dots,n\}$, set $a_i=p_i+i-1$. Then
$\sum\limits_{i=1}^n a_i=0$.
Let $V$ be the irreducible algebraic representation of $\GL_n$ with
highest weight $(a_1,\dots,a_n)$, and let $W^*$ be the
$\theta$-stable algebraic representation of $R_{E/\Q}\GL_{n,E}$ defined by
$V^*$ as in remark \ref{rq:ep_tordu}. (As the notation suggests, $W^*$
is
the contragredient of the irreducible algebraic representation $W$
of $R_{E/\Q}\GL_{n,E}$ associated to $\pi$ as in definition
\ref{def:algebrique_reguliere}.)
By proposition 3.5 of \cite{Cl-RGRA},
$ep(\theta,\pi_\infty\otimes W^*)\not=0$.
As $\pi_\infty$ is tempered, theorem \ref{th:pseudo_coeff_tordus} and
remark \ref{rq:ep_tordu} imply that $\pi_\infty\simeq\pi_V$,
so that $\pi_\infty$ is $\theta$-discrete.

We may therefore apply corollary \ref{cor:GL_n_E} to $\pi$. We get a
family of representations $\sigma_\lambda$ of $\Gal(\overline{\Q}/E)$.
Point (i) follows from (i) of corollary \ref{cor:GL_n_E}.
It remains to check the equality in point (ii).

Let $\H=\U(n)$, $\H'=\GL_n$, $\G^0=R_{E/\Q}\H_E=R_{E/\Q}\H'_E$,
and let $\theta'$ be the involution $g\fle\overline{g}$ of $\G^0$.
As in \ref{groupes3}, let $\Phi_n\in\GL_n(\Z)$ be the matrix with
coefficients : $(\Phi_n)_{i,j}=(-1)^{i-1}\delta_{i,n+1-j}$.
There is an isomorphism $\widehat{\G}^0\simeq\GL_n(\C)\times\GL_n(\C)$
such that :
\begin{itemize}
\item[-] the embedding $\widehat{\H}'=\GL_n(\C)\fl\widehat{\G}^0$ is
$g\fle (g,g)$;
\item[-] the embedding $\widehat{\H}=\GL_n(\C)\fl\widehat{\G}^0$
is $g\fle (g,\Phi_n{}^tg^{-1}\Phi_n^{-1})$;
\item[-] for every $(g,h)\in\widehat{\G}^0$,
$\widehat{\theta}(g,h)=(\Phi_n{}^th^{-1}\Phi_n^{-1},\Phi_n{}^tg^{-1}
\Phi_n^{-1})$ and $\widehat{\theta}'(g,h)=(h,g)$.

\end{itemize}
Let $\T$ be the diagonal torus of $\G^0$. Let $p\not\in S$ be a prime number.
Denote by $x=((y_1,\dots,y_n),(z_1,\dots,z_n))\in\widehat{\T}^{\Gal(\overline
{\Q}_p/\Q_p)}$ the Langlands parameter of $\pi_p$.
As $\pi_p$ is $\theta$-stable and $\theta'$-stable, we may assume that
$\widehat{\theta}(x)=\widehat{\theta}'(x)=x$, ie that
$y_i=z_i=y_{n+1-i}^{-1}$ for every $i\in\{1,\dots,n\}$.
Assume that $p$ is inert in $E$. Then $\H(E_p)\simeq\H'(E_p)$, and the
Langlands parameter of $\pi_p$, seen as a representation of $\H(E_p)$ or
$\H'(E_p)$, is $(y_1^2,\dots,y_n^2)$; on the other hand, the Langlands
parameter of $\tau_p$ is $(y_1,\dots,y_n)$, hence the image of $\Phi_{E_p}$
by $\varphi_{\tau_p}$ is $(y_1^2,\dots,y_n^2)$.
Assume that $p$ splits in $E$, and let $\wp$ and $\wp'$ be the places of $E$
above $p$. Then $\G^0(\Q_p)=\H(E_\wp)\times\H(E_{\wp'})=\H'(E_{\wp})\times
\H'(E_{\wp'})$.
Write $\pi_p=\pi_{\wp}\otimes\pi_{\wp'}=\pi'_{\wp}\otimes\pi'_{\wp'}$,
where $\pi_\wp$ (resp. $\pi_{\wp'}$, resp. $\pi'_{\wp}$, resp. $\pi'_{\wp'}$)
is an unramified representation of $\H(E_\wp)$ (resp. $\H(E_{\wp'})$,
resp. $\H'(E_\wp)$, resp. $\H'(E_{\wp'})$). Then the Langlands parameter of
$\tau_p$, $\pi_\wp$, $\pi_{\wp'}$, $\pi'_\wp$ or $\pi'_{\wp'}$ is
$(y_1,\dots,y_n)$. These calculations show that point (ii) follows from
(ii) of corollary \ref{cor:GL_n_E}.

\end{proof}

\section{A simple case of base change}
\label{GL_n_applications5}

As an application of the techniques in this chapter (and of the knowledge
about automorphic representations of general linear groups), it is possible
to obtain some weak base change results between general unitary
groups and general linear groups. These results are
spelled out in this section.

Use the notations of the beginning of \ref{GL_n_applications4};
in particular, $\H=\G(\U(p_1,q_1)\times\dots\times\U(p_r,q_r))$,
$\G^0=R_{E/\Q}\H_E$ and $\xi:{}^L\H\fl{}^L\G^0$ is the ``diagonal''
$L$-morphism.

As before, if $\Le$ is a connected reductive group over $\Q$, $p$ is a
prime number where $\Le$ is unramified and $\pi_{L,p}$ is an unramified
representation of $\Le(\Q_p)$, we will denote by $\varphi_{\pi_{L,p}}:
W_{\Q_p}\fl\widehat{\Le}\rtimes W_{\Q_p}$ a Langlands parameter of
$\pi_{L,p}$.

\begin{subdefinition}\label{def:correspondance_en_p} Let $\Le$ be a Levi
subgroup of $\G^0$ (it does not have to be the identity component of a
Levi subset of $\G$). Then there is a $W_\Q$-embedding $\widehat{\Le}\fl
\widehat{\G}^0$, unique up to $\widehat{\G}^0$-conjugacy; fix such an
embedding, and let $\eta_L:{}^L\Le\fl{}^L\G^0$ be the obvious $L$-morphism
extending it (ie the $L$-morphism whose restriction to $W_\Q$ is the
identity).

As $\G^0$ is isomorphic to $R_{E/\Q}(\Gr_{m,E}\times\GL_{n_1,E}\times\dots
\times\GL_{n_r,E})$, $\Le$ is a direct product of $R_{E/\Q}\Gr_{m,E}$ and of
groups of the type
$R_{E/\Q}\GL_{m,E}$, $m\in\Nat^*$. Choose an isomorphism
$\Le\simeq R_{E/\Q}(\Gr_{m,E}\times\GL_{m_1,E}\times\dots\times\GL_{m_l,E})$,
and denote by $\theta_L$ the automorphism $(x,g_1,\dots,g_l)\fle
(\overline{x},\overline{x}{}^t\overline{g}_1^{-1},\dots,\overline{x}
{}^t\overline{g}_l^{-1})$ of $\Le$;
the class of $\theta_L$ in the group of outer automorphisms of $\Le$
does not depend on the choices. (Note also that $\theta$ and
$\theta_{G^0}$ are equal up to an inner automorphism, so we can take
$\theta=\theta_{G^0}$.)
\index{$\theta_L$}

Let $\pi_H$ be an irreducible admissible representation of $\H(\Ade)$ and
$\pi_L$ be an irreducible admissible representation of $\Le(\Ade)$. Let
$v$ be a finite place where $\pi_H$ and $\pi_L$ are unramified. We say
that $\pi_H$ and $\pi_L$ \emph{correspond to each other at $v$} if
$\xi\circ\varphi_{\pi_{H,v}}$ and $\eta_L\circ\varphi_{\pi_{L,v}}$ are
$\widehat{\G}^0$-conjugate.
\index{correspond to each other at $v$}

\end{subdefinition}

\begin{subremarque} Let $v$ be a finite place of $\Q$ that is unramified in
$E$, $\Le$ be a Levi subgroup of $\G^0$, $\pi_{L,v}$ be an
unramified representation of $\Le^0(\Q_v)$ and $\pi_{H,v}$ be an unramified
representation of $\H(\Q_v)$.
Then $\xi\circ\varphi_{\pi_{H,v}}$ and $\eta_L\circ\varphi_{\pi_{L,v}}$ are
$\widehat{\G}^0$-conjugate if and only if, for every $\phi_v\in
\Hecke(\G^0(\Q_v),\G^0(\Z_v))$,
\[\Tr(\pi_v(b_\xi(\phi_v)))=\Tr(\pi_{L,v}((\phi_v)_L)),\]
where $(\phi_v)_L$ is the constant term of $\phi_v$ at $\Le$ (and
$b_\xi$ is, as in \ref{partie_en_p2}, the base change map
$\Hecke(\G^0(\Q_v),\G^0(\Z_v))\fl\Hecke(\H(\Q_v),\H(\Z_v))$).

\end{subremarque}

\begin{subcorollaire}\label{cor:CB}
\begin{itemize}
\item[(i)] Let $\pi_H$ be an irreducible admissible representation of
$\H(\Ade)$. Assume that :
\begin{itemize}
\item[-] there exists a neat open compact subgroup $\K_H$ of $\H(\Af)$ such
that $\pi_{H,f}^{\K_H}\not=\{0\}$;
\item[-] there exist an irreducible algebraic representation
$V$ of $\H$ and $i\in\Z$ such that $m_{disc}(\pi_H)\not=0$ and
$\Ho^i(\hgoth,\K'_\infty;\pi_{H,\infty}\otimes V)\not=0$, where
$\hgoth=Lie(\H(R))\otimes\C$ and the notations are those of
remark \ref{rq:c_G(pi)} (or lemma \ref{lemme:conj_RP4}).

\end{itemize}
(In other words, $\pi_H$ is a discrete automorphic
representation of $\H(\Ade)$ and appears in the intersection cohomology of
some Shimura variety associated to $\H$.)

Then there exists a Levi subgroup $\Le$ of $\G^0$, a
cuspidal automorphic representation $\pi_L$ of $\Le(\Ade)$ and an automorphic
character $\chi_L$ of $\Le(\Ade)$ such that $\pi_L\otimes\chi_L^{-1}$ is
$\theta_L$-stable
and that $\pi_H$ and $\pi_L$
correspond to each other at almost every finite place. If $\Le=\G^0$, then
we can take $\chi_L=1$, and $\pi_L$ is regular algebraic.

\item[(ii)] Assume that $\H$ is quasi-split at every finite place of $\Q$.
Let $\pi$ be a $\theta$-stable cuspidal automorphic
representation of $\G^0(\Ade)$. Assume that there exists an irreducible
algebraic representation $V$ of $\H$ such that $\Tr(\pi_\infty(\phi_V))
\not=0$, where $\phi_V\in C^\infty(\G^0(\R))$ is associated to $V$ as in
lemma \ref{lemme:O_pseudo_coeff}.
(In other words, $\pi$ is regular algebraic, cf lemma
\ref{lemme:algebrique_reguliere}.)
Let $S$ be the union of the set of finite places that are ramified in $E$ and
of the set of places where $\pi$ is ramified.
Then there exists an automorphic
representation $\pi_H$ of $\H(\Ade)$, unramified outside of $S$, such that
$\pi$ and $\pi_H$ correspond to each other at every finite place $v\not\in S$.
Moreover :
\begin{itemize}
\item[-] $\pi_H$ satifies the conditions of (i) for $V$ and some $\K_H$;
\item[-] if $\pi'_H$ satisfies the conditions of (i) and is isomorphic to
$\pi_H$ at almost every finite place, then $\pi'_H$ is cuspidal (in particular,
$\pi_H$ is cuspidal);
\item[-] in the notation of
\ref{applications1}, for every $\es\in\F_\H$, $R_\es(\pi_{H,f})=\varnothing$
(ie $\pi_H$ ``does not come from an endoscopic group of $\H$'').
\end{itemize}

\end{itemize}

\end{subcorollaire}

\begin{subremarque}\label{rq:CB}
Note that, using (ii) of the corollary, we can strengthen
(i) a little. We obtain the following statement :
if, in (i), $\H$ is quasi-split at every finite place of $\Q$ and
$\Le=\G^0$, then $\pi_H$ is unramified at $v$
as soon as $\pi:=\pi_L$ is, and $\pi_H$ and $\pi$ correspond to each other
at every finite place of $\Q$ where $\pi$ is unramified.

\end{subremarque}

In the rest of this section, we use the following notations : if $S$ is a set
of places of $\Q$, then $\Z_S=\prod\limits_{v\in S}\Z_v$ and $\Z^S=\prod
\limits_{v\not\in S\cup\{\infty\}}\Z_v$.
\index{ZS@$\Z_S$}
\index{ZS@$\Z^S$}

\begin{proof} We show (i). Let $\pi_H$, $\K_H$ and $V$ be as in (i).
Let $W_\lambda$ be the virtual $\lambda$-adic representation of
$\Hecke(\H(\Af),\K_H)\times\Gal(\overline{\Q}/E)$ defined by $\K_H$ and
$V$ as in \ref{applications1}. After replacing $\K_H$ by a smaller
open compact subgroup of $\H(\Af)$, we may assume that $\ungras_{\K_H}$
satisfies condition (H) of \ref{GL_n_applications4} (cf (i) of lemma
\ref{lemme:hypothese_H}) and that
$\K_H=\prod_{v\not=\infty}\K_{H,v}$. For every set $S$ of finite places of
$\Q$, we will write $\K_{H,S}=\prod\limits_{v\in S}\K_{H,v}$ and $\K_H^S=\prod
\limits_{v\not\in S\cup\{\infty\}}\K_{H,v}$.
Let $S$ be a finite set of finite places of $\Q$ such that,
for every $v\not\in S\cup\{\infty\}$, $\H$ is unramified at $v$ and
$\K_{H,v}=\H(\Z_v)$. Set $f_S=\ungras_{\K_{H,S}}$.

Let $\Pi'(\pi_H)$ be the set of isomorphism classes of irreducible admissible
representations $\pi'_H$ of $\H(\Ade)$ such that :
\begin{itemize}
\item[-] $\pi'_H$ satisfies the conditions of (i) for the same $\K_H$ and $V$
as $\pi_H$;
\item[-] for almost every finite place $v$ of $\Q$ such that $\H$ is unramified
at $v$, $\xi\circ\varphi_{\pi_{H,v}}$ and $\xi\circ\varphi_{\pi'_{H,v}}$ are
$\widehat{\G}^ 0$-conjugate.

\end{itemize}
(We use the notation $\Pi'(\pi_H)$ to avoid confusion with the set $\Pi_H$
of the proof of (ii).)
Then $\Pi'(\pi_H)$ is finite (in fact, the set of $\pi'_H$ that satisfy the
first condition defining $\Pi'(\pi_H)$ is already finite). So there exists
a finite set $T\supset S$ of finite places of $\Q$ and a function
$f_{T-S}\in\Hecke(\H(\Ade_{T-S}),\K_{H,T-S})$ such that :
\begin{itemize}
\item[-] $f_{T-S}$ is in the image of the base change map $b_\xi:
\Hecke(\G^0(\Ade_{T-S}),\G^0(\Z_{T-S}))\fl\Hecke(\H(\Ade_{T-S}),\K_{H,T-S})$
(hence, by (iii) of lemma \ref{lemme:hypothese_H}, it satisfies condition (H)).
\item[-] Let $\pi'_H$ be an irreducible admissible representation of
$\H(\Ade)$ that satisfies the conditions of (i) for $\K_H$ and $V$. Then
$\Tr(\pi'_{H,T-S}(f_{T-S}))=1$ if $\pi'_H\in\Pi'(\pi_H)$, and $0$
otherwise.
\item[-] For every $v\not\in T$ finite, for every $\pi'_H\in\Pi'(\pi_H)$,
$\xi\circ\varphi_{\pi_{H,v}}$ and $\xi\circ\varphi_{\pi'_{H,v}}$ are
$\widehat{\G}^ 0$-conjugate.

\end{itemize}
Write $f_T=f_Sf_{T-S}$. As in \ref{GL_n_applications4}, $f_T$ determines
functions $\phi_{M,T}\in C_c^\infty(\M^0(\Ade_T))$, for $\M\in\Mcal_\G$.

Fix a prime number $p\not\in T$. Then, with the notations of
\ref{applications1}, for every $m\in\Z$ and every $f^{T\cup\{p\}}$
in the image of the
base change map $b_\xi:\Hecke(\G^0(\Af^{T\cup\{p\}}),\G^0(\Z^{T\cup\{p\}}))\fl
\Hecke(\H(\Af^{T\cup\{p\}}),\K_H^{T\cup\{p\}})$,
\[\Tr(\Phi_\wp^mf_Tf^{T\cup\{p\}},W_\lambda)=\Tr(\pi_{H,f}^{\T\cup\{p\}}(
f^{T\cup\{p\}}))\sum_{\pi'_H\in\Pi'(\pi_H)}\dim((\pi_H')^{\K_H})\Tr(\Phi_\wp^m,
W_\lambda(\pi'_H)).\]
Note that the virtual representation $\sum\limits_{\pi'_H\in\Pi'(\pi_H)}
\dim((\pi_H')^{\K_H})W_\lambda(\pi'_H)$ of $\Gal(\overline{\Q}/E)$ is not
trivial. This is proved as the fact that (2) implies (1) in remark
\ref{rq:c_G(pi)} (see this proof for more details) :
Let $w$ be the weight of $V$ in the sense of
\ref{points_fixes3}. Then, for every $\pi'_H\in\Pi'(\pi_H)$,
$W_\lambda(\pi'_H)=\sum\limits_{i\in Z}(-1)^iW_\lambda^i(\pi'_H)$, with
$W_\lambda^i(\pi'_H)$ a true (not virtual) representation of $\Gal
(\overline{\Q}/E)$ that is unramified and of weight $-w+i$ at almost every
place of $E$. Hence there can be no cancellation in the sum
$\sum\limits_{\pi'_H\in\Pi'(\pi_H)}\dim((\pi_H')^{\K_H})W_\lambda(\pi'_H)$.
(And the assumptions on $\pi_H$ imply that $W_\lambda(\pi_H)$ is not trivial,
by remark \ref{rq:c_G(pi)}.)
So there exists an integer $m\in\Z$ such that
\[C:=\sum_{\pi'_H\in\Pi'(\pi_H)}\dim((\pi_H')^{\K_H})\Tr(\Phi_\wp^m,
W_\lambda(\pi'_H))\not=0.\]

For every $\M\in\Mcal_\G$, write $\phi_{M,p}=\phi_{M,p}^{(m)}$, where
$\phi_{M,p}^{(m)}$ is as in \ref{GL_n_applications4}, and define
$\phi_{M,\infty}$ as in loc. cit.
Then, by theorem \ref{th:stab_FT_IC_GL_n} and the calculations above
(and (iii) of lemma \ref{lemme:hypothese_H}),
for every $\phi^{T\cup\{p\}}\in
\Hecke(\G^0(\Af^{T\cup\{p\}}),\G^0(\Z^{T\cup\{p\}}))$,
\[C\Tr(\pi_{H,f}^{T\cup\{p\}}(f^{T\cup\{p\}}))=\sum_{\M\in\Mcal_\G}c_M T^M
(\phi_M),\]
where $f^{T\cup\{p\}}=b_\xi(\phi^{T\cup\{p\}})$ and $\phi_M=\phi_{M,T}
\phi_{M,p}\phi_M^{T\cup\{p\}}\phi_{M,\infty}$,
with $\phi_M^{T\cup\{p\}}$ equal to
the product of $\chi_{M|\G^0(\Af^{T\cup\{p\}})}$ and of the constant term
of $\phi^{T\cup\{p\}}$ at $\M^0$.

By lemma \ref{lemme:CB1} below, the right hand side of this equality, seen
as a linear form $T$ over $\Hecke(\G^0(\Af^{T\cup\{p\}}),\G^0(\Z^{T\cup
\{p\}}))$,
is a finite linear combination of linear maps of the form
$\phi^{T\cup\{p\}}\fle\Tr(\pi_L((\phi^{T\cup\{p\}})_L)),$ with
$\Le$, $\pi_L$ and $(\phi^{T\cup\{p\}})_L$ as in this lemma.
By the strong multiplicity $1$ theorem of Jacquet and Shalika (cf theorem
4.4 of \cite{JS}), there exist a Levi subgroup $\Le$ of $\G^0$, a
cuspidal automorphic representation $\pi_L$ of $\Le(\Ade)$, an
automorphic character $\chi_L$ of $\Le(\Ade)$, a scalar $a\in\C^\times$,
a finite set $\Sigma\supset T\cup\{p\}$ of
finite places of $\Q$ and a function $\phi_{\Sigma-T\cup\{p\}}\in\Hecke(
\G^0(\Ade_{\Sigma-T\cup\{p\}}),\G^0(\Z_{\Sigma-T\cup\{p\}}))$ such that :
$\pi_L\otimes\chi_L^{-1}$ is $\theta_L$-stable and,
for every $\phi^\Sigma\in\Hecke(\G^0(\Af^\Sigma),\G^0(\Z^\Sigma))$,
\[T(\phi_{\Sigma-T\cup\{p\}}\phi^\Sigma)=a\Tr(\pi_L^\Sigma((\phi^\Sigma)_L))\]
($a$ is non-zero because $T$ is non-zero, and that in turn follows from the
fact that $C\not=0$; the existence of $\chi_L$ such that $\pi_L\otimes\chi_L
^{-1}$ is $\theta_L$-stable comes from lemma \ref{lemme:CB1}).

Let $D=C\Tr(\pi_{H,\Sigma-T\cup\{p\}}(b_\xi(\phi_{\Sigma-T\{p\}}))$.
Then we finally find that,
for every $\phi^\Sigma\in\Hecke(\G^0(\Af^\Sigma),\G^0(\Z^\Sigma))$,
\[D\Tr(\pi_{H,f}^\Sigma(b_\xi(\phi^\Sigma)))=a\Tr(\pi_L^\Sigma((\phi^\Sigma)
_L)).\]
In particular, $D\not=0$ (because $a\not=0$).
This equality implies that $\pi_H$ and $\pi_L$ correspond to each other at
every finite place $v\not\in\Sigma$. 
Assume that $\Le=\G^0$. Then it is obvious from the definition of $T$ and from
the definition of $\phi_{G,\infty}$ in \ref{GL_n_applications4} that
the infinitesimal character of $\pi_L$ is equal to that of an algebraic
representation of $\G^0$. By (i) of lemma \ref{lemme:algebrique_reguliere},
$\pi_L$ is regular algebraic.
This finishes the proof of (i).

We show (ii). Assume that $\H$ is quasi-split at every finite place of $v$,
and let $\pi$, $V$ and $S$ be as (ii). Fix an open compact
subgroup $\K=\prod_v\K_v$ of $\G^0(\Af)$ such that $\Tr(\pi_f(\ungras_\K)
A_{\pi_f})\not=0$, that $\K_v=\G^0(\Z_v)$ for $v\not\in S$
and that $\ungras_\K$ has a
transfer to $\H$ satisfying condition (H) (such a $\K$ exists by
lemma \ref{lemme:hypothese_H}).
Set $\K_S=\prod_{v\in S}\K_v$, $\phi_S=\ungras_{\K_S}$,
and let $f_S\in C_c^\infty(\H(\Ade_S))$ be a transfer of $\phi_S$
satisfying condition (H). Choose an open compact subgroup $\K_{H,S}$
of $\H(\Ade_S)$ such that $f_S$ is bi-$\K_{H,S}$-invariant, and set
$\K_H=\K_{H,S}\prod\limits_{v\not\in S\cup\{\infty\}}\H(\Z_v)$.
After making $\K_{H,S}$ smaller,
we may assume that $\K_H$ is neat.
Let $W_\lambda$ be the $\lambda$-adic virtual representation of
$\Hecke(\H(\Af),\K_H)\times\Gal(\overline{\Q}/E)$ defined by
$\K_H$ and $V$ as in \ref{applications1}.

Let $\phi^S\in\Hecke(\G^0(\Af^S),\G^0(\Z^S))$.
By theorem \ref{th:stab_FT_IC_GL_n} for a prime number $p\not\in S$
where $\phi^S$ is $\ungras_{\G^0(\Z_p)}$ and for $j=0$, there exist scalars
$c'_M\in\R$, for $\M\in\Mcal_\G$, such that
\[\Tr(f_Sb_\xi(\phi^S),W_\lambda)=\sum_{\M\in\Mcal_\G}c'_MT^M(\phi_{S,\M}
\phi^S_\M\phi_{M,\infty}),\]
where $\phi_{S,M}\in C_c^\infty(\M^0(\Ade_S))$ and $\phi^S_\M\in
C_c^\infty(\M^0(\Af^S))$ are the functions associated to
$\phi_S$ and $\phi^S$ as in the beginning of \ref{GL_n_applications4}
(so $\phi^S_\M$ is the product of $\chi_{M|M^0(\Af^S)}$ and of the constant
term of $\phi^S$ at $\M^0$) and $\phi_{M,\infty}\in C^\infty(\M^0(\R))$
is obtained from $V$ as in loc. cit. (so $\phi_{G,\infty}$ is the function
$\phi_V$ that appears in the statement of (ii)).
But, as in he proof of proposition
\ref{prop:FT_stable_geometrique}, we see that, for any $\M\in\Mcal_\G$,
the function $f_{H_M,p}^{(0)}$ of definition
\ref{def:phi_m_G} is equal to the product of $\frac{\iota_{\H,\H_M}}
{\iota(\H,\H_M)}$ and of a transfer of $\ungras_{\H(\Z_p)}$ to $\H_M$.
So, for every $\M\in\Mcal_\G$, $c'_M=\frac{\iota_{\H,\H_M}}{\iota(\H,\H_M)}
c_M$; in particular, $c'_M$ does not depend on $\phi^S$, and
$c'_G\not=0$ (because all the signs in the definition of $\iota_{\H,\H}$ are
obviously equal to $1$).

Consider the function $T$ that sends $\phi^S$ to
$\sum\limits_{\M\in\Mcal_\G}c'_MT^M(\phi_{S,\M}\phi^S_\M\phi_{M,\infty})$.
It is a linear form on $\Hecke(\G^0(\Af^S),\G^0(\Z^S))$. 
By lemma \ref{lemme:CB1}, $T$ is a finite linear combination of
characters on $\Hecke(\G^0(\Af^S),\G^0(\Z^S))$ of the form
$\phi^S\fle\Tr(\pi_L^S((\phi^S)_L))$, where $\Le$ is a Levi subgroup of
$\G^0$ and $\pi_L$ is a cuspidal automorphic representation of $\Le(\Ade)$.
By the strong multiplicity $1$ theorem of Jacquet and Shalika, these
characters are pairwaise distinct.
Hence, by the choice of $\phi_S$, the assumption on $\pi_\infty$
and the fact that $c'_G\not=0$,
the coefficient of the character $\phi^S\fle\Tr(\pi^S(\phi^S))$ in $T$
is non-zero.

Let $R_H$ be the set of equivalence classes of irreducible admissible
representations of $\H(\Ade)$ that satisfy the conditions of (i) for
$\K_H$ and $V$. Then $R_H$ is finite. Define an equivalence relation
$\sim$ on $R_H$ in the following way : if $\pi_H,\pi'_H\in R_H$, then
$\pi_H\sim\pi'_H$ if and only if, for every finite place $v\not\in S$ of $\Q$,
$\xi\circ\varphi_{\pi_{H,v}}$ and $\xi\circ\varphi_{\pi'_{H,v}}$ are
$\widehat{\G}^0$-conjugate ($\pi_H$ and $\pi'_H$ are necessarily unramified at
$v$ because we have chosen $\K_H$ to be hyperspecial outside of $S$).
Let $\Pi_H\in R_H/\sim$. Then, if $\phi^S\in\Hecke(\G^0(\Ade^S),\G^0(\Z^S))$,
the value $\Tr(\pi_H^S(b_\xi(\phi^S)))$ is the same for every $\pi_H\in
\Pi_H$; denote it by $\Tr(\Pi_H^S(b_\xi(\phi^S)))$. Let
\[c(\Pi_H)=\sum_{\pi_H\in\Pi_H}\dim(W_\lambda(\pi_H))\Tr(\pi_{H,S}(f_S))\]
(where the $W_\lambda(\pi_H)$ are as in \ref{applications1}).

By the definition of the $W_\lambda(\pi_H)$, $W_\lambda=\sum\limits_{\pi_H\in
R_H}W_\lambda(\pi_H)\otimes\pi_{H,f}^{\K_H}$. Hence, for every
$\phi^S\in\Hecke(\G^0(\Ade^S),\G^0(\Z^S))$,
\[T(\phi^S)=\Tr(f_Sb_\xi(\phi^S),W_\lambda)=\sum_{\Pi_H\in R_H/\sim}
c(\Pi_H)\Tr(\Pi_H^S(b_\xi(\phi^S))).\]
As the characters of $\Hecke(\G^0(\Ade^S),\G^0(\Z^S))$ are linearly
independant, there exists $\Pi_H\in R_H/\sim$ such that, for every
$\phi^S\in\Hecke(\G^0(\Ade^S),\G^0(\Z^S))$, $\Tr(\Pi_H^S(b_\xi(\phi^S)))=
\Tr(\pi^S(\phi^S))$. Let $\pi_H$ be any element of $\Pi_H$. It is
unramified outside of $S$ and corresponds to $\pi$ at every finite
$v\not\in S$.

It remains to prove the last three properties on $\pi_H$. The first one
is true because $\pi_H\in R_H$ by construction.
The second and third ones are easy
consequences of (i) and of the strong multiplicity $1$ theorem for
$\G^0$.

\end{proof}

\begin{sublemme}\label{lemme:CB1}
Let $S$ be a set of finite places of $\Q$ that contains
all the finite places that ramify in $E$ and let $V$ be an irreducible
algebraic representation of $\H$.
Fix functions $\phi_{M,S}\in C_c^\infty(\M^0(\Ade_S))$, $\M\in\Mcal_\G$,
and let $\phi_{M,\infty}$, $\M\in\Mcal_\G$, be the functions associated to
$V$ as in \ref{GL_n_applications4}.
Consider the linear form $T:\Hecke(\G^0(\Af^S),\G^0(\Z^S))\fl\C$
that sends $\phi^S$ to
\[\sum_{\M\in\Mcal_\G}c_M T^M(\phi_M),\]
where the $\phi_M^S$ are obtained from $\phi^S$ as
in \ref{GL_n_applications4}.

Then $T$ is a finite linear combination of linear maps of the form
$\phi^S\fle\Tr(\pi_L^S((\phi^S)_L)),$ where :
\begin{itemize}
\item[(a)] $\Le$ is a Levi subgroup of $\G^0$.
\item[(b)] $\pi_L$ is a cuspidal automorphic representation of $\Le(\Ade)$
such that there exists an automorphic character $\chi_L$ of $\Le(\Ade)$ such
that $\pi_L\otimes\chi_L^{-1}$ is $\theta_L$-stable; if $\Le=\G^0$, then
we can take $\chi_L=1$.
\item[(c)] $(\phi^S)_L$ is the constant term of $\phi^S$
at $\Le$.

\end{itemize}

\end{sublemme}

\begin{proof} If we did not care about condition (b) on $\pi_L$, the lemma
would be an easy consequence of proposition \ref{prop:FT_invariante_tordue},
(iii) of lemma \ref{lemme:hypothese_H} and lemma \ref{lemme:rep_nr_tordues}.
As we do, we must know more precisely what
kind of automorphic representations appear on the spectral side of
$T^M(\phi_M)$. This is the object of lemma \ref{lemme:CB2} below. Once we
know this lemma (for all $\M\in\Mcal_\G$), the proof is straightforward.

\end{proof}

\begin{sublemme}\label{lemme:CB2}
\footnote{Most of this lemma (and of its proof) was worked out during
conversations with Sug Woo Shin.}
Let $V$ be an irreducible algebraic
representation of $\H$, and $\phi_V\in C^\infty(\G^0(\R))$ be the function
associated to $V$ as in \ref{GL_n_applications4}.
Let $t\geq 0$ and let $\pi\in\Pi_{disc}(\G,t)$ such that $a_{disc}(\pi)\Tr(
\pi_\infty(\phi_V)A_{\pi_\infty})\not=0$ (the notations are those of
\ref{GL_n_applications2}).
Choose a Levi subgroup $\Le$ of $\G^0$, a parabolic subgroup $\QP$ of $\G^0$
with Levi subgroup $\Le$ and a cuspidal
automorphic representation $\pi_L$ of $\Le(\Ade)$ such that $\pi$ is a
subquotient of the parabolic induction $Ind_{\QP}^{\G^0}\pi_L$.

Then there exists an automorphic character $\chi_L$ such that
$\pi_L\otimes\chi_L$ is $\theta_L$-stable.

\end{sublemme}

\begin{proof} Remember that, in this section, we take $\theta$ on
$\G^0=R_{E/\Q}(\Gr_{m,E}\times\GL_{n_1,E}\times\dots\times\GL_{n_r,E})$
to be $(x,g_1,\dots,g_r)\fle(\overline{x},{}^t\overline{g}_1^{-1},\dots,
\overline{x}{}^t\overline{g}_r^{-1})$ ($\theta$ is determined only up to
inner automorphisms, and this is a possible choice of $\theta$).
Let $\T_H$ be the diagonal torus of $\H$. Then $\T_H(\C)$ is the diagonal
torus of $\H(\C)\simeq\C^\times\times\GL_{n_1}(\C)\times\dots\times
\GL_{n_r}(\C)$. Let $\T=R_{E/\Q}\T_{H,E}$, let
$\theta'$ be the automorphism of $\H(\C)$ defined by
the same formula as $\theta$, and choose an isomorphism $\G^0(\C)\simeq
\H(\C)\times\H(\C)$ such that $\theta$ corresponds
to the automorphism $(h_1,h_2)\fle (\theta'(h_2),\theta'(h_1))$ of
$\H(\C)\times\H(\C)$ and that $\T(\C)$ is sent to $\T_H(\C)\times\T_H(\C)$.
Let $\tgoth_H=Lie(\T_H)(\C)$ and $\tgoth=Lie(\T)(\C)$. Then $\tgoth=\tgoth_H
\oplus\tgoth_H$, and $\theta$ acts on $\tgoth$ by $(t_1,t_2)\fle(\iota(t_2),
\iota(t_1))$, where $\iota$ is the involution $(t,(t_{i,j})_{1\leq i\leq r,1
\leq j\leq n_i})\fle (t,(t-t_{i,j})_{1\leq i\leq r,1\leq j\leq n_i})$ of
$\tgoth_H=\C\oplus\C^{n_1}\oplus\dots\oplus\C^{n_r}$.
Let $\lambda\in\tgoth_H^*$ be a representative of the infinitesimal character
of $V$, seen as a representation of $\H(\R)$.
Then $\lambda$ is regular, and $(\lambda,\iota(\lambda))\in
\tgoth^*=\tgoth_H^*\oplus\tgoth_H^*$ represents the infinitesimal character
of the $\theta$-stable representation $W$ of $\G^0(\R)$ associated to
$V$ (as in theorem \ref{th:pseudo_coeff_tordus}).
By the definition of $\phi_V$ in theorem \ref{th:pseudo_coeff_tordus}, the
assumption that $\Tr(\pi_\infty(\phi_V)A_{\pi_\infty})\not=0$ implies
that $ep(\theta,\pi_\infty\otimes W)\not=0$, and, by Wigner's lemma,
this implies that the infinitesimal character of $\pi_\infty$ is
$(-\lambda,-\iota(\lambda))$. In particular, the infinitesimal character of
$\pi_\infty$ is regular.

Let $\Le$ be a Levi subgroup of $\G^0$. We may assume that $\Le$ is
standard (in particular, it is stable by $\theta$, and $\theta_{|\Le}=
\theta_L$), and we write
\[\Le=R_{E/\Q}(\Gr_{m,E}\times\prod_{j=1}^r\prod_{k=1}^{l_j}\GL_{m_{jk},E}),\]
where the $m_{jk}$ are non-negative integers such that $n_j=m_{j,1}+\dots
+m_{j,l_j}$ for every $j\in\{1,\dots,r\}$. We use the notations of
\cite{A-ITF2}, in particular of section 4 of this article. Then
$\agoth_L=\R\oplus\bigoplus\limits_{j=1}^r\bigoplus\limits_{k=1}^{l_j}\R$,
and $\theta$ acts by $(t,(t_{ik}))\fle(t,(t-t_{ik}))$.
In particular,
$\agoth_{G^0}=\R\oplus\bigoplus\limits_{j=1}^r\R$, and $\theta$ acts by
$(t,(t_j))\fle(t,(t-t_j))$, so $\agoth_G=\agoth_{G^0}^
{\theta=1}=\R$. The group $W^{G^0}(\agoth_L)$ is equal to the group of linear
automorphisms of $\agoth_L$ that are induced by an element of
$\Nor_{\G^0(\Q)}(\Le)$, so it is embedded in an obvious way in
$\Sgoth_{l_1}\times\dots\times\Sgoth_{l_r}$ (but this embedding is far
from being an equality in general; for example, $W^{G^0}(\agoth_L)$ contains
the factor $\Sgoth_{l_j}$ if and only if $m_{j,1}=\dots=m_{j,l_j}=n_j/l_j$).
The set $W^G(\agoth_L)$ is equal to $W^{G^0}(\agoth_L)\theta_L$.

It will be useful to determine the subset $W^G(\agoth_L)_{reg}$ of
regular elements of
$W^G(\agoth_L)$. Remember (\cite{A-ITF2} p 517) that an element $s$ of
$W^G(\agoth_L)$ is in $W^G(\agoth_L)_{reg}$ if and only if $\det(s-1)
_{\agoth_L/\agoth_G}\not=0$. By the above calculations,
$\agoth_L/\agoth_G=\bigoplus\limits_{j=1}^r\bigoplus\limits_{k=1}^{l_j}\R$,
and $\theta_L$ acts by multiplication by $-1$. Let $s\in W^G(\agoth_L)$, and
write $s=\sigma\theta_L$, with $\sigma=(\sigma_1,\dots,\sigma_r)\in
W^{G^0}(\agoth_L)\subset\Sgoth_{l_1}\times\dots\times\Sgoth_{l_r}$. Then
$\sigma\theta_L$ acts on $\agoth_L/\agoth_G$ by
$(\lambda_{j,k})_{1\leq j\leq r,1\leq k\leq l_j}\fle (-\lambda_{j,\sigma_j
^{-1}(k)})_{1\leq j\leq r,1\leq k\leq l_j}$, so it is regular if and only
if $\sigma$, acting in the obvious way on $\R^{l_1}\times\dots\times\R^{l_r}$,
does not have $-1$ as an eigenvalue. This is equivalent to the fact that, for
every $j\in\{1,\dots,r\}$, there are only cycles of odd length in the
decomposition of $\sigma_j$ as a product of cycles with disjoint supports.

We will need another fact. Let $s=(\sigma_1,\dots,\sigma_r)\theta_L\in
W^G(\agoth_L)$ as
before (we do not assume that $s$ is regular). Let $\pi_{L,\infty}$ be
an irreducible admissible representation of $\Le(\R)$. As $\Le$ is
standard, $\T$ is a maximal torus of $\Le$. Assume that the infinitesimal
character of $\pi_{L,\infty}$ has a representative of the form
$(\mu,\iota(\mu))\in\tgoth^*$,
with $\mu\in\tgoth_H^*$, and that $\pi_{L,\infty}$ is
$s$-stable. Then, for every $j\in\{1,\dots,r\}$, there are only cycles of
length $\leq 2$ in the decomposition of $\sigma_j$ as a product of cycles
with disjoint supports (in other words, $\sigma_j^2=1$ for every $j\in
\{1,\dots,r\}$).
To prove this, we assume (to simplify notations) that
$r=1$; the general case is similar, because it is possible to reason
independently on each factor $\GL_{n_j,E}$. Write $n=n_1$, $\sigma=\sigma_{n_1}
\in\Sgoth_n$  and
$\Le=R_{E/\Q}(\Gr_{m,E}\times\GL_{m_1,E}\times\dots\times\GL_{m_l,1})$, with
$n=m_1+\dots+m_l$. For every $j\in\{1,\dots,l\}$, set $I_j=\{m_1+\dots+m_{j-1}
+1,\dots,m_1+\dots+m_l\}\subset\{1,\dots,n\}$.
Let $(\mu,\iota(\mu))$, with $\mu\in\tgoth_H^*$, be a representative of the
infinitesimal character of $\pi_{L,\infty}$. Write $\mu=(\mu_0,\dots,\mu_n)\in
\tgoth_H^*=\C\oplus\C^n$. As $\pi_{L,\infty}$ is
$\sigma\theta_L$-stable, its infinitesimal character is
$\sigma\theta_L$-stable. This means that there exists $\tau\in\Sgoth_n$
such that :
\begin{itemize}
\item[(a)] for every $j\in\{1,\dots,l\}$, $\tau(I_j)=I_{\sigma(j)}$ (ie $\tau$,
seen as an element of the Weyl group of $\T(\Q)$ in $\G^0(\Q)$, normalizes
$\Le$ and induces $\sigma$ on $\agoth_L$);
\item[(b)] $(-\mu_1,\dots,-\mu_n)=(\mu_{\tau(1)},\dots,\mu_{\tau(n)})$
(ie $(\mu,\iota(\mu))$ is stable by $\tau\theta$, where $\tau$ is again seen
as an element of the Weyl group of $\T(\Q)$ in $\G^0(\Q)$).

\end{itemize}
Assume that, in the decomposition of $\sigma$ as a product of cycles with
disjoint supports, there is a cycle of length $\geq 3$. Then, by (a), there is
also a cycle of length $\geq 3$ in the decomposition of $\tau$ as a product
of cycles with disjoint supports. But, by (b), this contradicts the fact that
$\mu$ is regular.

We now come back to the proof of the lemma. By the definition of the
coefficients $a_{disc}=a^G_{disc}$ in section 4 of \cite{A-ITF2}, the fact that
$a_{disc}(\pi)\not=0$ implies that there exists a Levi subgroup $\Le$
of $\G^0$ (that we may assume to be standard), a discrete automorphic
representation $\pi_L$ of $\Le(\Ade)$ and an element $s$ of
$W^G(\agoth_L)_{reg}$ such that $\pi_L$ is $s$-stable and
$\pi$ is a subquotient of the parabolic induction of $\pi_L$ (where,
for example, we use the standard parabolic subgroup of $\G^0$ having
$\Le$ as Levi subgroup). In particular, the infinitesimal character
of $\pi_{L,\infty}$ is represented by $(-\lambda,-\iota(\lambda))$
(and it is regular).
By the two facts proved above, $s$ is equal to $\theta_L$. That is,
$\pi_L$ is a $\theta_L$-stable discrete automorphic representation of
$\Le^0(\Ade)$. Hence, if we know that the lemma is true for
discrete automorphic representations (and for any Levi subgroup of $\G^0$),
then we now that it is true in general. So we may
assume that $\pi$ is discrete.

From now on, we assume that the automorphic representation $\pi$ of
$\G^0(\Ade)$ is discrete (and of course $\theta$-stable). Let
$\Le$ be a standard Levi subgroup of $\G^0$ and $\pi_L$ be a cuspidal
automorphic representation of $\Le(\Ade)$ such that $\pi$ is a subquotient
of the parabolic induction of $\Le$ (as before, use the standard parabolic
subgroup of $\G^0$ with Levi subgroup $\Le$). We want to show that $\pi_L$
is $\theta_L$-stable. As $\pi$ is $\theta$-stable, there exists $s\in
W^G(\agoth_L)$ such that $\pi_L$ is $s$-stable (note that $s$ does not have
to be regular now). We also know that $\pi$ is discrete; so, by a
result of Moeglin and Waldspurger (the main theorem of \cite{MW}),
there exist $m_1,\dots,m_r,l_1,\dots,m_r\in\Nat$, an automorphic
character $\chi$ of $\Ade_E^\times$ and cuspidal automorphic representations
$\tau_j$ of $\GL_{m_j}(\Ade_E)$, for $1\leq j\leq r$, such that :
\begin{itemize}
\item[-] $n_j=l_jm_j$ for every $j\in\{1,\dots,r\}$;
\item[-] $\Le=R_{E/\Q}(\Gr_{m,E}\times(\GL_{m_1,E})^{l_1}\times\dots\times
(\GL_{m_r,E})^{l_r})$;
\item[-] $\pi_L=\chi\otimes\pi_1\otimes\dots\otimes\pi_r$, where
$\pi_j$ is the cuspidal automorphic representation
$\tau_j|\det|^{(l_j-1)/2}\otimes\tau_j|\det|^{(l_j-3)/2}\otimes\dots\otimes
\tau_j|\det|^{(1-l_j)/2}$
of $\GL_{m_j}(\Ade_E)^{l_j}$.

\end{itemize}
Write $s=\sigma\theta_L$, with $\sigma=(\sigma_1,\dots,\sigma_r)\in
\Sgoth_{l_1}\times\dots\times\Sgoth_{l_r}=W^{G^0}(\agoth_L)$. As
$\pi$ is a subquotient of the parabolic induction of $\pi_L$, the
infinitesimal character of $\pi_{L,\infty}$ is representated by
$(-\lambda,-\iota(\lambda))$, and so, by the fact proved above, $\sigma^2=1$.
For every $j\in\{1,\dots,r\}$, denote by $\theta_{m_j}$ the automorphism
$g\fle{}^t\overline{g}^{-1}$ of $\GL_{m_j,E}$. Then the fact that
$\pi_L$ is $\sigma\theta_L$-stable means that, for every $j\in\{1,\dots,r\}$
and for every $k\in\{1,\dots,l_j\}$,
\[\tau_j|\det|^{(l_j-2k+1)/2}\simeq(\tau_j|\det|^{(l_j-2\sigma_j(k)+1)/2})\circ
\theta_{m_j}=(\tau_j\circ\theta_{m_j})|\det|^{(2\sigma_j(k)-1-l_j)/2},\]
ie that
\[\tau_j\circ\theta_{m_j}\simeq\tau_j|\det|^{l_j+1-k-\sigma_j(k)}.\]
In particular, we see by taking the absolute values of the central characters
in the equality above that, for every $j\in\{1,\dots,r\}$, the function
$k\fle k+\sigma_j(k)$ is constant on $\{1,\dots,l_j\}$. But, if
$j\in\{1,\dots,r\}$, then
$\sum\limits_{k=1}^{l_j}(k+\sigma_j(k))=l_j(l_j+1)$, so $k+\sigma_j(k)=
l_j+1$ for every $k\in\{1,\dots,l_j\}$.
This show that $\tau_j\simeq\tau_j\circ\theta_{m_j}$ for every $j\in\{1,\dots,
r\}$, and it easily implies that $\pi_L$ is $\theta_L$-stable after we twist
it by an automorphic character.

\end{proof}

\chapter{The twisted fundamental lemma}
\label{lemme_fondamental}

The goal of this chapter is to show that, for the special kind of twisted
endoscopic transfer that appears in the stabilization of the fixed point
formula and for the groups considered in this book (and some others),
the twisted fundamental lemma for the whole Hecke algebra follows from
the twisted fundamental lemma for the unit element.
No attempt has been made to prove a general result, and the method is
absolutely not original : it is simply an adaptation of the method
used by Hales in the untwisted case (\cite{Ha}), and this was inspired
by the method used by Clozel in the case of base change (\cite{Cl-LF})
and by the simplification suggested by the referee of the article
\cite{Cl-LF}.

The definitions and facts about twisted groups recalled in
\ref{GL_n_applications1} will be used freely in this chapter.

\section{Notations}
\label{lemme_fondamental1}

We will consider the following situation : Let $F$ be a local
non-archimedean field of characteristic $0$. Fix an algebraic closure
$\overline{F}$ of $F$, let $\Gamma_F=\Gal(\overline{F}/F)$
and denote by $F^{ur}$ the maximal unramified extension of $F$ in
$\overline{F}$.
\index{$\Gamma_F$}
\index{Fur@$F^{ur}$}
Fix a uniformizer $\varpi_F$ of $F$.
\index{$\varpi_F$}
Let $\G$ be a connected reductive unramified group over $F$. Assume that
$\G$ is defined over $\Of_F$ and that $\G(\Of_F)$ is a hyperspecial
maximal compact subgroup of $\G(F)$. For such a group $\G$, write
$\Hecke_G=\Hecke(\G(F),\G(\Of_F)):=C_c^\infty(\G(\Of_F)\sous\G(F)/
\G(\Of_F))$. 
\index{HG@$\Hecke_G$\quad spherical Hecke algebra of $\G$}
Let $(\H,s,\eta_0)$ be an endoscopic triple for $\G$ (in the sense of
\cite{K-STF:CTT} 7.4). Assume that :
\begin{itemize}
\item[-] $\H$ is unramified over $F$, $\H$ is defined over $\Of_F$ and
$\H(\Of_F)$ is a hyperspecial maximal compact subgroup of $\H(F)$;
% \item[-] $s\in Z(\widehat{\H})^{\Gamma_F}$;
\item[-] there exists a $L$-morphism $\eta:{}^L\H\fl{}^L\G$
extending $\eta_0$ and unramified, ie coming by inflation from a
$L$-morphism
$\widehat{\H}\rtimes\Gal(K/F)\fl\widehat{\G}\rtimes\Gal(K/F)$, where $K$
is a finite unramified extension of $F$.

\end{itemize}
Choose a generator $\sigma$ of $W_{F^{ur}/F}$.
Let $E/F$ be a finite unramified extension of $F$ in $\overline{F}$, and
let $d\in\Nat^*$ be the degree of $E/F$. 
Write $R=R_{E/F}\G_E$. Let $\theta$ be the automorphism of $R$ induced by
the image of $\sigma$ in $\Gal(E/F)$.

Kottwitz explained in \cite{K-SVLR} p 179-180 how to get from this 
twisted endoscopic data for $(R,\theta,1)$ (in the sense of \cite{KS} 2.1).
\index{twisted endoscopic data}
\index{endoscopic data}
We recall his construction.
There is an obvious isomorphism $\widehat{R}=\widehat{\G}^d$, and the actions
of $\widehat{\theta}$ and $\sigma$ are given by the formulas
\[\widehat{\theta}(g_1,\dots,g_d)=(g_2,\dots,g_d,g_1)\]
\[\sigma(g_1,\dots,g_d)=(\sigma(g_2),\dots\sigma(g_d),\sigma(g_1)).\]
In particular, the diagonal embedding $\widehat{\G}\fl\widehat{R}$ is
$W_F$-equivariant, hence extends in an obvious way to a $L$-morphism
${}^L\G\fl{^L}R$. Let $\xi':{}^L\H\fl{}^LR$ be the composition of the
morphism $\eta:{}^L\H\fl{}^L\G$ and of this $L$-morphism.
As $F$ is local, we may assume that $s\in Z(\widehat{\H})^{\Gamma_F}$.
Let $t_1,\dots,t_d\in Z(\widehat{\H})^{\Gamma_F}$ be such that
$s=t_1\dots t_d$. Set $t=(t_1,\dots,t_d)\in\widehat{R}$. 
Let $\xi:{}^L\H\fl{}^LR$ be the morphism such that :
\begin{itemize}
\item[-] $\xi_{|\widehat{\H}}$ is the composition of $\eta_0:\widehat{\H}\fl
\widehat{\G}$ and of the diagonal embedding $\widehat{\G}\fl\widehat{R}$;
\item[-] for every $w\in W_F$ that is a pre-image of $\sigma\in W_{F^{ur}/
F}$, $\xi(1,w)=(t,1)\xi'(1,w)$.

\end{itemize}
Then $(\H,{}^L\H,t,\xi)$ are twisted endoscopic data for $(R,\theta,1)$.
Kottwitz shows (\cite{K-SVLR}, p 180) that the equivalence class of
these twisted endoscopic data does not depend on the choice of
$t_1,\dots,t_d$.

The morphism $\xi$ induces a morphism
\[\Hecke_R\fl\Hecke_H,\]
that will be denoted by $b_\xi$ (in \ref{partie_en_p2}, this morphism is
explicitely calculated for the unitary groups of \ref{groupes1}).

Let $\Delta_\xi$ be the transfer factors defined by $\xi$, normalized as in
\cite{Wa3} 4.6.
\index{$\Delta_\xi$\quad twisted transfer factors}
The twisted fundamental lemma for a function $f\in\Hecke_R$ is the following
statement :
\index{twisted fundamental lemma}
for every $\gamma_H\in\H(F)$ semi-simple and strongly $\G$-regular,
\[\Lambda(\gamma_H,f):=SO_{\gamma_H}(b_\xi(f))-\sum_\delta\Delta_\xi(\gamma_H,
\delta)O_{\delta\theta}(f)=0\]
\index{$\Lambda(\gamma_H,f)$}
where the sum is taken over the set of $\theta$-conjugacy classes of
$\theta$-semi-simple $\delta\in R(F)$.
(Remember that a semi-simple $\gamma_H\in\H(F)$ is called
\emph{strongly $\G$-regular} if
it has an image in $\G(F)$ whose centralizer is a torus.)
\index{strongly $\G$-regular}

\begin{subremarque}\label{rq:E_pas_un_corps} There is an obvious variant
of this statement where $E$ is replaced by a finite product of finite
unramified extensions of $F$ such that
$\Aut_F(E)$ is a cyclic group.

\end{subremarque}

\section{Local data}
\label{lemme_fondamental2}

Notations are as in \ref{lemme_fondamental1}.
Fix a Borel subgroup $\B$ (resp. $\B_H$) of $\G$ (resp. $\H$) 
defined over $\Of_F$ and a maximal torus
$\T_G\subset\B_G$ (resp. $\T_H\subset\B_H$) defined over $\Of_F$. 
Let $\T_R=R_{E/F}\T_{G,E}$ and $\B_R=R_{E/F}\B_{G,E}$.
Denote by $I_R$ (resp. $I_H$) the Iwahori subgroup of
$R(F)$ (resp. $\H(F)$) defined by the Borel subgroup $\B_R$
(resp. $\B_H$).

Let $\Pi(R)$ (resp. $\Pi(\H)$) be the set of equivalence classes of
irreducible $\theta$-stable representations $\pi$ of $R(F)$ (resp. of
irreducible representations $\pi_H$ of $\H(F)$) such that
$\pi^{I_R}\not=\{0\}$ (resp. $\pi_H^{I_H}\not=\{0\}$).
\index{$\Pi(R)$}
\index{$\Pi(\H)$}
For every $\pi\in\Pi(R)$, fix a normalized intertwining operator $A_\pi$
on $\pi$ (remember that an \emph{intertwining operator} on $\pi$ is a
$R(F)$-equivariant
isomorphism $\pi\iso\pi\circ\theta$, and that an intertwining operator
$A_\pi$ is called \emph{normalized} if $A_\pi^d=1$).
\index{normalized intertwining operator}
If $\pi$ is unramified, choose the intertwining operator that fixes the
elements of the subspace $\pi^{R(\Of_F)}$. 

The definition of local data used here is the obvious adaptation of the
definition of Hales (\cite{Ha} 4.1).

\begin{subdefinition}
\index{local data}
\emph{Local data} for $R$ and $(\H,{}^L\H,s,\xi)$
are the data of a set $I$ and of two families of complex numbers,
$(a_i^R(\pi))_{i\in I,\pi\in\Pi(R)}$ and $(a_i^H(\pi_H))_{
i\in I,\pi_H\in\Pi(\H)}$, such that, for every $i\in I$, the numbers
$a_i^R(\pi)$ (resp. $a_i^H(\pi_H)$) are zero for almost every
$\pi\in\Pi(R)$ (resp. $\pi_H\in\Pi(\H)$) and that, for every
$f\in\Hecke_R$, the following conditions are equivalent :
\begin{itemize}
\item[(a)] for every $i\in I$, $\sum\limits_{\pi\in\Pi(R)}a_i^R(\pi)
\Tr(\pi(f)A_\pi)=\sum\limits_{\pi_H\in\Pi(\H)}a_i^H(\pi_H)\Tr(\pi_H(b_\xi(f)))
$;
\item[(b)] for every $\gamma_H\in H(F)$ that is semi-simple, elliptic and
strongly $\G$-regular, $\Lambda(\gamma_H,f)=0$.

\end{itemize}

\end{subdefinition}

\begin{subproposition}\label{prop:utilisation_donnees_locales}
Assume that $\G$ is adjoint, that the endoscopic triple
$(\H,s,\eta_0)$ for $\G$ is elliptic, that the center of $\H$ is connected
and that there exist local data for $R$ and $(\H,{}^L\H,s,\xi)$.

Then, for every $f\in\Hecke_R$ and for every $\gamma_H\in\H(F)$ semi-simple
elliptic and strongly $\G$-regular, $\Lambda(\gamma_H,f)=0$.

\end{subproposition}

\begin{subremarque}\label{rq:xi_sur_extension} If $\G$ is adjoint and
$(\H,s,\eta_0)$ is elliptic, then the morphism $\xi:{}^L\H\fl{}^LR$
comes by inflation from a morphism
$\widehat{\H}\rtimes\Gal(K/F)\fl\widehat{R}\rtimes\Gal(K/F)$, where $K$ is
a finite unramified extension of $F$. Let us prove this.
By the definition of $t$, $\xi(1\rtimes\sigma^d)=\xi'(s\rtimes\sigma^d)$
(remember that $\xi'$ is the composition of $\eta:{}^L\H\fl{}^L\G$
and of the ``diagonal embedding'' ${}^L\G\fl{}^LR$). We know that
$s\in Z(\widehat{\H})^{\Gamma_F}Z(\widehat{\G})$. As $\G$ is adjoint,
$Z(\widehat{\G})$ is finite; because the endoscopic triple
$(\H,s,\eta_0)$ is elliptic, $Z(\widehat{\H})^{\Gamma_F}$ is also finite.
The finite subgroup $Z(\widehat{\G})Z(\widehat{\H})^{\Gamma_F}$ of
$Z(\widehat{\H})$ is invariant
by $\sigma^d$, so there exists $k\in\Nat^*$ such that the restriction of
$\sigma^{dk}$ to this subgroup is trivial. Let $s'=s\sigma^{d}(s)\dots
\sigma^{d(k-1)}(s)$. Then $s'$ is fixed by $\sigma^d$, and
$(s\rtimes\sigma^d)^k=s'\rtimes\sigma^{dk}$. As
$s'$ is in the finite group $Z(\widehat{\G})Z(\widehat{\H})^{\Gamma_F}$, there
exists $l\in\Nat^*$ such that ${s'}^l=1$. Then
\[\xi(1\rtimes\sigma^{dkl})=\xi'((s\rtimes\sigma^d)^{kl})=
\xi'((s'\rtimes\sigma^{dk})^l)=\xi'(1\rtimes\sigma^{dkl}).\]
By the assumption on $\eta$, there exists $r\in\Nat^*$ such that
$\eta(1\rtimes\sigma^r)=1\rtimes\sigma^r$. So we get finally :
$\xi(1\rtimes\sigma^{dklr})=1\rtimes\sigma^{dklr}$.

\end{subremarque}

Before proving the proposition, we show a few lemmas.

Remember that an element of $\H(F)$ or $\G(F)$ is called \emph{strongly
compact} if it belongs to a compact subgroup, and \emph{compact} if its
image in the adjoint group of $\H$ resp. $\G$ is strongly compact
(cf \cite{H} \S2).
\index{compact}
\index{strongly compact}
Every semi-simple ellipic element of $\G(F)$ or $\H(F)$ is compact, and
an element that is stably conjugate to a compact element is also compact
(this follows easily from the characterization of compact elements in
\cite{H} \S2).
If the center of $\G$ is anisotropic (eg if $\G$ is adjoint), then an
element of $\G(F)$ is compact if and only if it is strongly compact.

\begin{sublemme}\label{lemme:DL1} Assume that the centers of $\G$ and $\H$ are
anisotropic. Let $A:\Hecke_R\fl\C$ be a linear form.
Assume that : for every $f\in\Hecke_R$, if $\Lambda_{\gamma_H}(f)=0$
for every $\gamma_H\in\H(F)$ semi-simple elliptic and strongly $\G$-regular,
then $A(f)=0$.

Then $A$ is a linear combination of linear forms $\gamma_H\fle
\Lambda(\gamma_H,f)$, with $\gamma_H\in\H(F)$ semi-simple elliptic and
strongly $\G$-regular.

\end{sublemme}

\begin{proof}
Let $U_H\subset\H(F)$ be the set of compact elements of $\H(F)$ and
$U_R$ be the set of $\theta$-semi-simple elements of $R(F)$ whose
norm contains a compact element of $\G(F)$. Then $U_F$ is compact
modulo $\H(F)$-conjugacy and $U_R$ is compact modulo $\theta$-conjugacy
(these notions are defined before theorem 2.8 of \cite{Cl-LF}).
By the twisted version of the Howe conjecture (ie theorem 2.8 of \cite{Cl-LF}),
the vector space of distributions on $\Hecke_R$ generated by the
$f\fle O_{\delta\theta}(f)$, $\delta\in U_R$ and the $f\fle
O_{\gamma_H}(b_\xi(f))$, $\gamma_H\in U_F$, is finite-dimensional.
If $\gamma_H\in\H(F)$ is elliptic semi-simple, then $\gamma_H\in U_H$,
and every image of $\gamma_H$ in $R(F)$ is in $U_R$.
In particular, the vector space generated by the distributions
$\gamma_H\fle\Lambda(\gamma_H,f)$, for $\gamma_H\in\H(F)$ elliptic semi-simple
and strongly $\G$-regular, is finite-dimensional. The lemma follows from
this.

\end{proof}

Let $\Se_H$ be the maximal split subtorus of $\T_H$, $\Se_R$ be the
maximal split subtorus of $\T_R$ and $\Omega_H=\Omega(\Se_H(F),\H(F))$, 
$\Omega_R=\Omega(\Se_R(F),R(F))$ be the relative Weyl groups.
\index{$\Omega_H$}
\index{$\Omega_R$}
Identify $\Hecke_H$ (resp.
$\Hecke_R$) to $\C[\widehat{\Se}_H/\Omega_H]$ (resp. $\C[\widehat{\Se}_R/
\Omega_R]$) by the Satake isomorphism.
If $z\in\widehat{\Se}_H$ (resp. $\widehat{\Se}_R$) and $f\in\Hecke_H$
(resp. $\Hecke_R$), write $f(z)$ for $f(z\Omega_H)$ (resp. $f(z\Omega_R)$).
\index{SH@$\Se_H$}
\index{SR@$\Se_R$}

We recall the definition of the morphism $b_\xi:\Hecke_R\fl\Hecke_H$ induced
by $\xi$ (cf \cite{Bo} sections 6 and 7).
The group $\Omega_H$ (resp. $\Omega_R$) is naturally isomorphic to the
subgroup of $\Gamma_F$-fixed points of the Weyl group
$\Omega(\widehat{\T}_H,\widehat{\H})$
(resp. $\Omega(\widehat{\T}_R,\widehat{R})$). Let $N_H$ (resp $N_R$)
be the inverse image of $\Omega_H$ (resp. $\Omega_R$) in $\Nor_{\widehat{\H}}
(\widehat{\T}_H)$ (resp. $\Nor_{\widehat{R}}(\widehat{\T}_R)$),
let $Y_H=\widehat{\Se}_H$ (resp. $Y_R=\widehat{\Se}_R$)
\index{YR@$Y_R$}
\index{YH@$Y_H$}
and let
$(\widehat{\H}\rtimes\sigma)_{ss}$ (resp. $(\widehat{R}\rtimes\sigma)_{ss}$) be
the set of semi-simple elements of
$\widehat{\H}\rtimes\sigma\subset\widehat{\H}\rtimes
W_{F^{ur}/F}$ (resp. $\widehat{R}\rtimes\sigma\subset\widehat{R}\rtimes W_{
F^{ur}/F}$) (remember that $\sigma$ is a fixed generator of
$W_{F^{ur}/F}$). As $X_*(\Se_H)=X_*(\T_H)^{\Gamma_F}$ (resp. $X_*(\Se_R)=
X_*(\T_R)^{\Gamma_F}$), the group $\Omega_H$ (resp. $\Omega_R$) acts
naturally on $Y_H$ (resp. $Y_R$). Moreover :
\begin{bulletlist}
\item the restriction to $(\widehat{\T}_H^{\Gamma_H})^0$ (resp.
$(\widehat{\T}_R^{\Gamma_H})^0$) of the morphism $\nu:\widehat{\T}_H\fl Y_H$
(resp. $\nu:\widehat{\T}_R\fl Y_R$) dual of the inclusion $\Se_H\subset\T_H$
(resp. $\Se_R\subset\T_R$) is an isogeny (\cite{Bo} 6.3);
\item the map $\widehat{\T}_H\rtimes\sigma\fl Y_H$
(resp. $\widehat{\T}_R\rtimes\sigma\fl Y_R$) that sends $t\rtimes\sigma$
to $\nu(t)$ induces a bijection
\[(\widehat{\T}_H\rtimes\sigma)/\Int N_H\iso Y_H/\Omega_H\]
\[(\mbox{resp.}\qquad(\widehat{\T}_R\rtimes\sigma)/\Int N_R\iso Y_R/\Omega_R)\]
(\cite{Bo}, lemma 6.4);
\item the inclusion induces a bijection
\[(\widehat{\T}_H\rtimes\sigma)/\Int N_H\iso(\widehat{\H}\rtimes\sigma)_{ss}/
\Int\widehat{\H}\]
\[(\mbox{resp.}\qquad(\widehat{\T}_R\rtimes\sigma)/\Int N_R\iso(\widehat{R}
\rtimes\sigma)_{ss}/\Int\widehat{R})\]
(\cite{Bo}, lemma 6.5).

\end{bulletlist}

In particular, we get bijections $\varphi_H:(\widehat{\H}\rtimes
\sigma)_{ss}/\Int\widehat{\H}\iso Y_H/\Omega_H$ and $\varphi_R:(\widehat{R}
\rtimes\sigma)_{ss}/\Int\widehat{R}\iso Y_R/\Omega_R$.
The morphism $\xi:{}^L\H\fl{}^LR$ is unramified, hence it induces a
morphism $(\widehat{\H}\rtimes\sigma)_{ss}/
\Int\widehat{\H}\fl(\widehat{R}\rtimes\sigma)_{ss}/\Int\widehat{R}$, and this
gives a morphism $b_\xi^*:Y_H/\Omega_H\fl Y_R/\Omega_R$. The morphism
$b_\xi:\C[Y_R/\Omega_R]\fl\C[Y_H/\Omega_H]$ is the dual of $b_\xi^*$.

Let $Y_H^u$ (resp. $Y_R^u$) be the maximal compact subgroup of
$Y_H$ (resp. $Y_R$).

\begin{sublemme}\label{lemme:DL2}
The morphism $b_\xi^*:Y_H/\Omega_H\fl
Y_R/\Omega_R$ sends $Y_H^u/\Omega_H$ to $Y_R^u/\Omega_R$ and
$(Y_H-Y_H^u)/\Omega_H$ to $(Y_R-Y_R^u)/\Omega_R$.

\end{sublemme}

\begin{proof}
Let $K$ be an unramified extension of $F$ such that $\H$ and $\G$ split over
$K$; write $r=[K:F]$.
For every
$g\rtimes\sigma\in\widehat{\H}\rtimes\sigma$ or $\widehat{R}\rtimes\sigma$,
write $N(g\rtimes\sigma)=g\sigma(g)\dots\sigma^{r-1}(g)$.

Let $G'$ be the set
of complex points of an algebraic group over $\C$. Copying the definition of
\cite{H} \S2, say that an element $g\in G'$ is
\emph{strongly compact} if there exists a compact subgroup of $G'$ containing
$g$. It is easy to see that this is equivalent to the fact that there exists a
faithful representation $\rho:G'\fl\GL_m(\C)$ such that the eigenvalues
of $\rho(g)$ all have module $1$. So a morphism of algebraic groups over $\C$
sends strongly compact elements to strongly compact elements.

Let $g\in\widehat{\H}$ be such that $g\rtimes\sigma$ is semi-simple.
We show that $\varphi_H(g\rtimes\sigma)\in Y_H^u/\Omega_H$ if and only if
$N(g\rtimes\sigma)$ is strongly compact.
After replacing $g\rtimes\sigma$ by a $\widehat{\H}$-conjugate, we may assume
that $g\in\widehat{\T}_H$. Assume that $N(g\rtimes\sigma)$
is strongly compact. Then $\nu(N(g\rtimes\sigma))=\nu(g)^r\in Y_H$ is
strongly compact, and this implies that $\nu(g)$ is strongly compact,
ie that $\nu(g)\in Y_H^u$. Assume now that $\varphi_H(g\rtimes\sigma)=
\nu(g)\Omega_H\in Y_H^u/\Omega_H$, ie that $\nu(g)\in Y_H^u$. Then
$\nu(N(g\rtimes\sigma))=\nu(g)^r\in Y_H^u$.
Moreover, $\widehat{\T}_H$ is abelian, so $N(g\rtimes\sigma)\in\widehat{
\T}_H^{\Gamma_H}$. As the restriction of $\nu$ to $\widehat{\T}_H^{\Gamma_H}$
is finite, it is easy to see that the fact that $\nu(N(g
\rtimes\sigma))\in Y_H^u$ implies that $N(g\rtimes\sigma)$ is strongly
compact.

Of course, there is a similar statement for $R$.
Hence, to finish the proof, it is enough to show that, for every
$g\rtimes\sigma\in\widehat{\H}\rtimes\sigma$, $N(g\rtimes\sigma)$ is
strongly compact if and only if $N(\xi(g\rtimes\sigma))$ is strongly
compact. Let $\xi_0:\widehat{\H}\fl\widehat{R}$ be the morphism induced by
$\xi$. Write $t'\rtimes\sigma=\xi(1\rtimes\sigma)$. Then, for every
$g\in\widehat{\H}$, $N(\xi(g\rtimes\sigma))=\xi_0(N(g\rtimes\sigma))
N(t'\rtimes\sigma)$.
By remark \ref{rq:xi_sur_extension}, we may assume (after replacing $K$ by
a bigger unramified extension of $F$) that $\xi(1\rtimes\sigma
^r)=1\rtimes\sigma^r$, ie that $N(t'\rtimes\sigma)=1$. Then the statement of
the lemma follows from the injectivity of $\xi_0$.

\end{proof}

\begin{sublemme}\label{lemme:DL3} For every $\delta\in R(F)$ that is 
$\theta$-regular $\theta$-semi-simple and $\theta$-elliptic, for every
$\gamma_H\in\H(F)$ that is regular semi-simple and elliptic, the distributions
$f\fle O_{\delta\theta}(f)$ and $f\fle
O_{\gamma_H}(b_\xi(f))$ on $\Hecke_R$ are tempered.

\end{sublemme}

\begin{proof} Remember that a distribution on $\Hecke_R$ is called tempered
if it extends continuously to the Schwartz space of 
rapidly decreasing bi-$R(\Of_R)$-invariant functions on $R(F)$
(defined, for example, in section 5 of \cite{Cl-LF}).
\index{tempered distribution}
For the first distribution of the lemma, this is proved in lemma 5.2 of
\cite{Cl-LF}. 
Moreover, the distribution $f_H\fle O_{\gamma_H}(f_H)$ on
$\Hecke_H$ is tempered (this is a particular case of lemma 5.2 of
\cite{Cl-LF}).
So, to prove that the second distribution of the lemma is tempered, it is
enough to show that the morphism $b_\xi:\Hecke_R\fl\Hecke_H$ extends to
the Schwartz spaces. To show this last statement, it is enough to show,
by the proof of lemma 5.1 of \cite{Cl-LF}, that $b_\xi^*$ sends
$Y_H^u/\Omega_H$ to $Y_R^u/\Omega_R$. This follows from lemma
\ref{lemme:DL2} above.

\end{proof}

Call a $\theta$-semi-simple element of $R(F)$ \emph{$\theta$-compact}
if its norm is compact.
\index{$\theta$-compact}
Let $\ungras_c$ be the characteristic function of
the set of semi-simple compact elements of $\H(F)$, and
$\ungras_{\theta-c}$ be the characteristic function of the set of
$\theta$-semi-simple $\theta$-compact elements of $R(F)$.
If $\pi_H$ is an irreducible admissible representation of $\H(F)$, $\pi$
is a $\theta$-stable irreducible admissible representation of $R(F)$ and
$A_\pi$ is a normalized intertwining operator on $\pi$,
define the compact trace of $\pi_H$
and the twisted $\theta$-compact trace of $\pi$ by the formulas :
\[\Tr_c(\pi_H(f_H)):=\Tr(\pi_H(\ungras_cf_H)),\qquad f_H\in C_c^\infty(\H(F))\]
\[\Tr_{\theta-c}(\pi(f)A_\pi):=\Tr(\pi(\ungras_{\theta-c}f)A_\pi),\qquad
f\in C_c^\infty(R(F)).\]
\index{Trc@$\Tr_c$\quad compact trace}
\index{Trthetac@$\Tr_{\theta-c}$\quad twisted $\theta$-compact trace}

\begin{sublemme}\label{lemme:DL4} Let $\pi$ be a $\theta$-stable
irreducible admissible representation of $R(F)$, and let $A_\pi$ be
a normalized intertwining operator on $\pi$. Assume that the distribution
$f\fle\Tr_{\theta-c}(\pi(f)A_\pi)$ on $\Hecke_R$ is not identically
zero. Then $\pi\in\Pi(R)$.

\end{sublemme}

\begin{proof} By the corollary to proposition 2.4 of \cite{Cl-LF},
there exists a $\theta$-stable parabolic subgroup $\Pa\supset\B_R$ of $R$
such that $\pi_{\N_P}$ is unramified, where $\N_P$ is the unipotent
radical of $\Pa$ and $\pi_{\N_P}$ is the unnormalized Jacquet module
(ie the module of $\N_P(F)$-coinvariants of $\pi$). So, if $\N$ is
the unipotent radical of $\B_R$, then $\pi_\N$ is unramified. By proposition
2.4 of \cite{Cas}, this implies that $\pi^{I_R}\not=\{0\}$.

\end{proof}

\begin{sublemme}\label{lemme:DL5} Assume that the centers of $\G$ and $\H$
are anisotropic and that the center of $\H$ is connected.
Let $\delta\in R(F)$ be $\theta$-regular $\theta$-semi-simple and
$\theta$-elliptic, and $\gamma_H\in\H(F)$ be regular semi-simple and elliptic.
Then the distribution $f\fle O_{\delta\theta}(f)$ on $\Hecke_R$ is a
linear combination of distributions $f\fle\Tr_{\theta-c}(\pi(f)A_\pi)$,
with $\pi\in\Pi(R)$, and the distribution
$f_H\fle SO_{\gamma_H}(f_H)$ on $\Hecke_H$ is a linear combination of
distributions $f_H\fle\Tr_c(\pi_H(f_H))$, with $\pi_H\in\Pi(\H)$ coming
from an element of $\Pi(\H_{ad})$.

\end{sublemme}

\begin{proof} We show the first assertion of the lemma.
Let $f\in\Hecke_R$ be such that $\Tr_{\theta-c}(\pi(f)A_\pi)=0$ for every
$\pi\in\Pi(R)$; let us show that $O_{\delta\theta}(f)=0$. As
$\delta$ is $\theta$-elliptic, hence $\theta$-compact,
$O_{\delta\theta}(f)=O_{\delta\theta}(\ungras_{\theta-c}f)$. But
$\Tr(\pi(\ungras_{\theta-c}f)A_\pi)=\Tr_{\theta-c}(\pi(f)A_\pi)=0$ for
every $\pi\in\Pi(R)$, so, by the main theorem of \cite{KRo} and lemma
\ref{lemme:DL4}, $O_{\delta\theta}(\ungras_{\theta-c}f)=0$.

On the other hand, by the twisted version of the Howe conjecture (theorem
2.8 of \cite{Cl-LF}), the space generated by the distributions (on $\Hecke_R$)
$f\fle\Tr_{\theta-c}(\pi(f)A_\pi)$, $\pi\in\Pi(R)$, is finite-dimensional.
This implies the first assertion.

We show the second assertion of the lemma. As $Z(\H)$ is anisotropic and
connected, lemma \ref{lemme:quotient_par_le_centre} implies that, for every
$f_H\in\Hecke_H$, $SO_{\gamma_H}(f_H)=SO_{\gamma'_H}(f'_H)$, where $\gamma'_H$
is the image of $\gamma_H$ in $\H_{adj}(F)$ and $f'_H$ is the image of $f_H$
in $\Hecke_{H_{adj}}$ (defined in lemma \ref{lemme:quotient_par_le_centre}).
So the second assertion of the lemma follows from the first, applied to
the group $\H_{adj}$ (with $\theta=1$).

\end{proof}

Identify the group of unramified characters of $\T_R(F)$ to $Y_R$ in the
usual way (cf \cite{Bo} 9.5). For every $z\in Y_R$, let $\psi_z$ be the
unramified character of $\T_R(F)$ corresponding to $z$ and
denote by $I(z)$ the representation of $R(F)$ obtained by (normalized)
parabolic induction from $\psi_z$ :
\[I(z)=Ind_{\B_R}^R(\delta_{B_R}^{1/2}\otimes\psi_z),\]
where, if $\N$ is the unipotent radical of $\B_R$, $\delta_{B_R}(t)=
|\det(\Ad(t),Lie(\N))|_F$ for every $t\in\T_G(F)$
($\delta_{B_R}^{1/2}\otimes\psi_z$ is seen as a character on $\B_R(F)$ via
the projection $\B_R(F)\fl\T_R(F)$).
\index{Iz@$I(z)$}
If $\widehat{\theta}(z)=z$, then $\psi_z=\psi_z\circ\theta$, so $I(z)$ is
$\theta$-stable. In that case, let $A_{I(z)}$ be the linear endomorphism of
the space of $I(z)$ that sends a function $f$ to the function $x\fle
f(\theta(x))$ (remember that the space of $I(z)$ is a space of functions
$R(F)\fl\C$); then $A_{I(z)}$ is a normalized intertwining operator on
$I(z)$.
We will use similar notations for $\H$ (without the intertwining operators,
of course).

\begin{sublemme}\label{lemme:DL6} Assume that $\G$ is adjoint and that
the center of $\H$ is anisotropic
(the second assumption is true, for example, if $\G$ is adjoint and the
endoscopic triple $(\H,s,\eta_0)$ is elliptic).
Let $\pi\in\Pi(R)$. Then there exist a $\theta$-stable $z'\in Y_R-Y_R^u$,
a $\theta$-stable subquotient $\pi'$ of $I(z')$ and a normalized
intertwining operator $A_{\pi'}$ on $\pi'$ such that, for every
$f\in\Hecke_R$, $\Tr_{\theta-c}(\pi(f)A_\pi)=
\Tr_{\theta-c}(\pi'(f)A_{\pi'})$.

Similarly, if $\pi_H\in\Pi(\H)$ comes from a representation in $\Pi(\H_{ad})$,
then there exist $z'_H\in Y_H-Y_H^u$ and a subquotient $\pi'_H$ of $I(z'_H)$
such that, for every $f_H\in\Hecke_H$, $\Tr_c(\pi_H(f_H))=\Tr_c(\pi'_H(f_H))$.

\end{sublemme}

\begin{proof} By proposition 2.6 of \cite{Cas}, there exists $z_0\in Y_R$
such that $\pi$ is a subrepresentation of $I(z_0)$. By examining the proof
of this proposition, we see that $\theta(z_0)\in\Omega_Rz_0$. By lemma 4.7
of \cite{Cl-LF}, there exists a $\theta$-stable $z$ in $\Omega_Rz_0$. Then
$\pi$ is a subquotient of $I(z)$ (because $I(z_0)$ and $I(z)$ have the
same composition factors).
If $z\not\in Y_R^u$, this finishes the proof of the first statement (take
$z'=z$ and $\pi'=\pi$).
Assume that $z\in Y_R^u$. As $\G$ is adjoint,
by a result of Keys (cf \cite{Ke}, in
particular the end of section 3), the representation $I(z)$ is irreducible,
hence $\pi=I(z)$.
Let $z'\in Y_R-Y_R^u$ be $\theta$-stable. The unramified characters
$\chi_z$ and $\chi_{z'}$ corresponding to $z$ and $z'$ are
equal on the set of $\theta$-compact elements of $\T_R(F)$. Hence, by theorem
3 of \cite{vD}, $\Tr_{\theta-c}(\pi(f)A_\pi)
=\Tr_{\theta-c}(I(z')(f)A_{I(z')})$ for every $f\in C_c^\infty(R(F))$.
This finishes the proof of the first statement (take $\pi'=I(z')$).

The same reasoning (without the twist by $\theta$) applies to $\pi_H$, or
rather to the representation of $\H_{ad}(F)$ inducing $\pi_H$; note that,
as the center of $\H$ is anisotropic, $Y_H=Y_{H_{ad}}$. We need the fact
that $\pi_H$ comes from a representation in $\Pi(\H_{ad})$ to apply
Keys's result.

\end{proof}

In the following lemma, $\N$ is the unipotent radical of $\B_R$ and, for
every representation $\pi$ of $R(F)$, $\pi_N$ is the $\T_R(F)$-module
of $\N(F)$-coinvariants of $\pi$.

\begin{sublemme}\label{lemme:DL7} Let $\pi$ be a $\theta$-stable admissible
representation of $R(F)$ of finite length, and let $A_\pi$ be an intertwining
operator on $\pi$.
The semi-simplification of $\delta_B^{-1/2}\otimes\pi_\N$ is a sum of
characters of $\T_R(F)$; let $z_1,\dots,z_n$ be the points of
$Y_R$ corresponding to the $\theta$-stable unramified characters that appear in
that way.
Then the distribution $f\fle\Tr(\pi(f)A_\pi)$ on $\Hecke_R$ is a linear
combination of distributions $f\fle f(z_i)$, $1\leq i\leq n$.
Moreover, if $\pi$ is a subquotient of $I(z)$, with $z\in Y_R$ $\theta$-stable,
then the $z_i$ are all in $\Omega_Rz$.

\end{sublemme}

Of course, there is a similar result for $\H$.

\begin{proof} As $\pi$ and its semi-simplification have the same character,
we may assume that $\pi$ is irreducible. We may also assume that $\pi$ is
unramified (otherwise the result is trivial). By proposition 2.6 of
\cite{Cas}, there exists $z\in Y_R$ such that $\pi$ is a subquotient of
$I(z)$. Reasoning as in the proof of lemma \ref{lemme:DL6} above, we may
assume that $z$ is $\theta$-stable.
By corollary 2.2 of \cite{Cas}, $I(z)^{R(\Of_F)}$ is $1$-dimensional;
hence $I(z)^{R(\Of_F)}=\pi^{R(\Of_F)}$. By the explicit description of a
basis of $I(z)^{R(\Of_F)}$ in \cite{Car} 3.7 and the definition of the
Satake transform (see for example \cite{Car} 4.2), for every
$f\in\Hecke_R$, $\Tr(f,I(z)^{R(\Of_F)})=\Tr(f,\pi^{R(\Of_F)})=f(z)$.
As $\pi^{R(\Of_F)}$ is $1$-dimensional and stable by $A_\pi$, the
restriction of $A_\pi$ to this subspace is the multiplication by a scalar.
So the distribution $f\fle\Tr(\pi(f)A_\pi)$ is equal to a scalar multiple of
the distribution $z\fle f(z)$.
By theorem 3.5 of \cite{Car}, the $z_i$ are all in $\Omega_Rz$. This finishes
the proof of the lemma.

\end{proof}

For every $\lambda\in X^*(Y_R)$, set 
\[f_\lambda=\sum_{\omega\in\Omega_R}\lambda^\omega\in\C[Y_R]^{\Omega_R}\simeq
\C[Y_R/\Omega_R]\simeq\Hecke_R.\]

\begin{sublemme}\label{lemme:DL8} There exists a non-empty open cone $C$ in
$X^*(Y_R)\otimes_\Z\R$ such that
\begin{itemize}
\item[(a)] for every $\theta$-stable $z\in Y_R$, for every $\theta$-stable
subquotient $\pi$ of $I(z)$ and for every intertwining operator
$A_\pi$ on $\pi$, the restriction to $C\cap X^*(Y_R)$ of the function
$\lambda\fle\Tr_{\theta-c}(\pi(f_\lambda)A_\pi)$ on $X^*(Y_R)$
is a linear combination of the functions
$\lambda\fle\lambda(\omega z)$, $\omega\in\Omega_R$.
\end{itemize}

Assume moreover that there exists an admissible embedding $\T_H\fl\G$ with
image $\T_G$
\footnote{It would be enough to assume this over an unramified extension
$K/F$ such that the base change morphism $\Hecke(R(K),R(\Of_K))\fl
\Hecke(R(F),R(\Of_F))$ is surjective.}
and that the center of $\G$ is connected (both assumptions are automatic
if $\G$ is adjoint, cf lemma \ref{lemme:groupes_qui_marchent}).

Then there exists a non-empty open cone $C$ in $X^*(Y_R)\otimes_\Z\R$ that
satisfies condition (a) above and also the following condition :
\begin{itemize}
\item[(b)]  for every $z_H\in Y_H$, for every subquotient $\pi_H$ of
$I(z_H)$, the restriction to $C\cap X^*(Y_R)$ of the function 
$\lambda\fle\Tr_c(\pi_H(b_\xi(f_\lambda)))$ on $X^*(Y_R)$ is a linear
combination of the functions $\lambda\fle\lambda(\omega b_\xi^*
(z_H))$, with $\omega\in\Omega_R$.

\end{itemize}
\end{sublemme}

\begin{proof} We will need some new notations.
If $\Pa\supset\B_R$ be  is a parabolic subgroup of $R$, let $\N_P$ be the
unipotent radical of $\Pa$, $\M_P$ be the Levi subgroup of $\Pa$ that
contains $\T_R$, $\Omega_{M_P}=\Omega(\Se_R(F),\M_P(F))$ be the relative
Weyl group of $\M_P$, $\delta_P$ be the function $\gamma\fle|\det(\Ad(\gamma)
,Lie(\N_P))|_F$ on $\Pa(F)$, $\agoth_{M_P}=\Hom(X^*(\A_{M_P}),\R)$ and
$a_P=\dim(\agoth_{M_P})$.
Assume that $\Pa$ is $\theta$-stable.
Let $\Pa_0$ and $\M_0$ be the parabolic subgroup and the Levi subgroup of
$\G$ corresponding to $\Pa$ and $\M_P$ (cf example
\ref{ex:groupes_non_connexes}). Denote by $H_{M_0}:\M_0(F)\fl\agoth_{M_0}:=
\Hom(X^*(\A_{M_0}),\R)$ the Harish-Chandra morphism (cf \cite{A-TFRG1} p
917), $\hat{\tau}_{P_0}^G:\agoth_T:=\Hom(X^*(\T_G),\R)\fl\{0,1\}$ the
characteristic function of the obtuse Weyl chamber defined by $\Pa_0$
(cf \cite{A-TFRG1} p 936) and $\hat{\chi}_{N_0}=\hat{\tau}
_{P_0}^G\circ H_{M_0}$ (there is a canonical injective morphism
$\agoth_{M_0}\subset\agoth_T$).
Define a function $\hat{\chi}_{N_P,\theta}$ on $\M_P(F)$ by :
$\hat{\chi}_{N_P,\theta}(m)=\hat{\chi}_{N_0}(\Norme m)$ if $m\in
\M_P(F)$ is $\theta$-semi-simple, and $\hat{\chi}_{N_P,\theta}=0$ otherwise. 
If $\pi$ is a $\theta$-stable admissible representation of $R(F)$ of finite
length and $A_\pi$ is an intertwining operator on $\pi$, denote by
$\pi_{\N_P}$ the Jacquet module of $\pi$ (ie the module of $\N_P$-coinvariants
of $\pi$) and by $A_\pi$ the intertwining operator on $\pi_{\N_P}$ induced
by $A_\pi$. If $f\in\Hecke_R$, denote by $f^{(P)}\in\Hecke_{\M_P}$ the constant
term of $f$ at $\Pa$.

Let $\pi$ be a $\theta$-stable admissible representation of $R(F)$ of
finite length and $A_\pi$ be an intertwining operator on $\pi$.
The corollary to proposition 2.4 of \cite{Cl-LF} says that, for every
$f\in\Hecke_R$,
\renewcommand\theequation{$*$}
\begin{equation}
\Tr_{\theta-c}(\pi(f)A_\pi)=\sum_{\Pa}(-1)^{a_P-a_G}\Tr((\delta_P^{-1/2}
\otimes\pi_{\N_P})(\hat{\chi}_{N_P,\theta}f^{(P)})A_\pi),
\end{equation}
where the sum is taken over the set of $\theta$-stable parabolic subgroups
$\Pa$ of $R$ that contain $\B_R$.

Let $N:\agoth_{\T_R}\fl\agoth_{\T_G}$, $\lambda\fle\lambda+\theta(\lambda)+
\dots+\theta^{d-1}(\lambda)$, and identify $X^*(Y_R)\otimes_\Z\R$ to
$\agoth_{\T_R}$.
Let $\lambda\in X^*(Y_R)$ and let $\Pa\supset\B_R$ be a $\theta$-stable
parabolic subgroup of $R$. Then $f_\lambda^{(P)}=\sum\limits_
{\omega\in\Omega_R}\lambda^\omega\in\C[Y_R]^{\Omega_{M_P}}\simeq\Hecke_{M_P}$,
and it follows easily from the definitions that, for every
$\omega\in\Omega_R$,
\[\hat{\chi}_{N_P,\theta}\sum_{\omega'\in\Omega_{M_P}}\lambda^{\omega'\omega}=
\hat{\tau}_{P_0}^G(N(\lambda^\omega))\sum_{\omega'\in\Omega_{M_P}}\lambda^{
\omega'\omega}.\]
From the definition of the functions $\hat{\tau}_{P_0}^G$, it is clear that
there exists a finite union $D\subset\agoth_{\T_G}$ of hyperplanes (containing
the origin) such that, for every parabolic subgroup $\Pa_0$ of $\G$,
$\hat{\tau}_{P_0}^G$ is constant on the connected components of
$\agoth_T-D$ (take for $D$ the union of the kernels of the fundamental
weights of $\T_G$ in $\B_G$).
Then $D':=N^{-1}(D)\subset\agoth_{T_R}$ is a finite union of hyperplanes
and, for every $\theta$-stable parabolic subgroup $\Pa\supset\B_R$ of $R$,
the function $\hat{\tau}_{P_0}^G\circ N$ is constant on the connected
components of $\agoth_{T_R}-D'$. After replacing $D'$ by $\bigcup\limits
_{\omega\in\Omega_R}\omega(D')$, we may assume that, for every connected
component $C$ of $\agoth_{T_R}-D'$, for all $\lambda,\lambda'\in C$,
for every $\theta$-stable parabolic subgroup $\Pa\supset\B_R$ of $R$
and for every $\omega\in\Omega_R$,
\[\hat{\tau}_{P_0}^G\circ N(\lambda^\omega)
=\hat{\tau}_{P_0}^G\circ N((\lambda')^\omega).\]
Let $C$ be a connected component of  $\agoth_{T_R}-D'$.
The calculations above show that there exist subsets
$\Omega'_{M_P}$ of $\Omega_{R}$, indexed by the set of $\theta$-stable
parabolic subgroups $\Pa$ of $R$ containing $\B_R$, such that :
for every $\lambda\in C$, for every $\Pa$,
\[\hat{\chi}_{N_P,\theta}f_\lambda^{(P)}=\sum_{\omega\in\Omega'_{M_P}}
\lambda^\omega.\]

Let $z\in Y_R$ be $\theta$-stable, let $\pi$ be a $\theta$-stable subquotient
of $I(z)$ and let $A_\pi$ be an intertwining operator on $\pi$.
For every $\theta$-stable parabolic subgroup $\Pa\supset\B_R$ of $R$,
\[\delta_{B_R}^{-1/2}\otimes\pi_{\N_{B_R}}=\delta_{B_R\cap M_P}^{-1/2}
\otimes(\delta_P^{-1/2}\otimes\pi_{\N_P})_{\N_{B_R\cap M_P}}.\] 
By formula $(*)$, the calculation of the functions $\hat{\chi}_{N_P,\theta}f_
\lambda^{(P)}$ above and lemma \ref{lemme:DL7} (applied to the representations
$\delta_P^{-1/2}\otimes\pi_{\N_P}$),
the restriction to $C\cap X^*(Y_R)$ of the function
$\lambda\fle\Tr_{\theta-c}(\pi(f_\lambda)A_\pi)$ is a linear combination
of functions $\lambda\fle\lambda(\omega z)$, with
$\omega\in\Omega_R$. Hence $C$ is a cone satisfying
condition (a) of the lemma.

We now show the second statement of the lemma.
After replacing the embeddings of $\widehat{\T}_H$ and $\widehat{\T}_G$
in $\widehat{\H}$ and $\widehat{\G}$ by conjugates, we may assume that
$\xi_0$ induces a $\Gamma_F$-equivariant isomorphism
$\widehat{\T}_H\iso\widehat{\T}_G$. Use this isomorphism to identify
$\T_H$ and $\T_G$.
By the definition of $\xi$, the restriction to $\widehat{\T}_H$ of
$\xi_0:\widehat{\H}\fl\widehat{\G}$ induces a $\Gamma_F$-equivariant
morphism $\widehat{\T}_H\fl\widehat{\T}_R$. Let $t'\rtimes\sigma=\xi(1\rtimes
\sigma)$. Then $t'$ centralizes the image (by the diagonal embedding) of
$\widehat{\T}_G$ in $\widehat{R}=\widehat{\G}^d$; as $\widehat{\G}_{der}$
is simply connected, $t'\in\widehat{\T}_R$. The isomorphism $\T_H\simeq\T_G$
fixed above induces an isomorphism $\agoth_{\T_G}\simeq\agoth_{\T_H}$, and
we can see the morphism $N:\agoth_{T_R}\fl\agoth_{\T_G}$, $\lambda\fle\lambda+
\theta(\lambda)+\dots+\theta^{d-1}(\lambda)$ defined above as a morphism
$\agoth_{\T_R}\fl\agoth_{\T_H}$. We may identify $\agoth_{\T_R}$ and
$\agoth_{\T_H}$ to $X^*(Y_R)\otimes_\Z\R$ and $X^*(Y_H)\otimes_\Z\R$,
and then $N$ sends $X^*(Y_R)$ to $X^*(Y_H)$. 
It is easy to see that
$b_\xi:\C[Y_R/\Omega_R]\fl\C[Y_H/\Omega_H]$ sends $f_\lambda$ to
$|\Omega_H|^{-1}\sum\limits_{\omega\in\Omega_R}\lambda(t')f_{N(\lambda^\omega)
}$, for every $\lambda\in X^*(Y_R)$, where $\lambda(t')$ denotes the value
of $\lambda$ at the image of $t'$ by the obvious morphism
$\widehat{\T}_R\fl Y_R=\widehat{\Se}_R$.
Let $D_H\subset\agoth_{\T_H}$ be the union of the kernels of the
fundamental weights of $\T_H$ in $\B_H$, and let $D'_H$ be the union of the
$\omega(N^{-1}(D_H))$, for $\omega\in\Omega_R$.
Then, for every connected component $C$ of $\agoth_{\T_R}-D'_H$, for all
$\lambda,\lambda'\in C$, for every parabolic subgroup $\Pa_H\supset\B_H$ of
$\H$ and for every $\omega\in\Omega_R$,
\[\hat{\tau}_{P_H}^H\circ N(\lambda^\omega)=\hat{\tau}_{P_H}^H\circ
N((\lambda')^\omega).\]
By the untwisted version of the reasoning above (applied to the calculation
of compact traces of representations of $\Pi(\H)$), such a connected
component $C$ satisfies condition (b).
Hence a connected component of $\agoth_{\T_R}-(D'\cup D'_H)$
satisfies conditions (a) and (b).

\end{proof}

The next lemma will be used in section \ref{lemme_fondamental3}.
It is a vanishing result similar to proposition
3.7.2 of \cite{La-CSCB}.

\begin{sublemme}\label{lemme:resultat_annulation} Assume that there exists
an admissible embedding $\T_H\fl\G$ with image $\T_G$ and that the center of
$\G$ is connected. Let $\gamma_H\in\H(F)$ be semi-simple elliptic and
strongly $\G$-regular. Assume that, for every $\theta$-semi-simple
$\delta\in R(F)$, no element of $\Norme\delta$ is an image of $\gamma_H$
in $\G(F)$. Then, for every $f\in\Hecke_R$, $O_{\gamma_H}(b_\xi(f))=0$.

\end{sublemme}

As the condition on $\gamma_H$ is stable by stable conjugacy, the lemma
implies that, under the same hypothesis on $\gamma_H$, $SO_{\gamma_H}(b_\xi(f))
=0$ for every $f\in\Hecke_R$.

\begin{proof} The proof is an adaptation of the proof of proposition 3.7.2
of \cite{La-CSCB}.
We first reformulate the condition on $\gamma_H$. By proposition 2.5.3 of
\cite{La-CSCB}, a semi-simple elliptic element of $\G(F)$ is a norm if and
only if its image in $\H^0_{ab}(F,\G)$ is a norm. Using the proof of
proposition 1.7.3 of \cite{La-CSCB}, we get canonical isomorphisms
$\Ho^0_{ab}(F,\G)=\Ho^1(\Gamma_F,Z(\widehat{\G}))^D$ and
$\Ho^0_{ab}(F,\H)=\Ho^1(\Gamma_F,Z(\widehat{\H}))^D$ (where ${}^D$ means
Pontryagin dual). As there is a canonical $\Gamma_F$-equivariant
embedding $Z(\widehat{\G})\subset Z(\widehat{\H})$, we get a canonical morphism
$\H^0_{ab}(F,\H)\fl\H^0_{ab}(F,\G)$. The condition of the lemma on
$\gamma_H$ is equivalent to the following condition :
the image of $\gamma_H$ in $\H^0_{ab}(F,\G)$ is not a norm.

Assume that $\gamma_H$ satisfies this condition. Then there exists a
character $\chi$ of $\H^0_{ab}(F,\G)$ that is trivial on the norms and
such that $\chi_H(\gamma_H)\not=1$, where $\chi_H$ is the character of
$\H(F)$ obtained by composing $\chi$ and the morphism $\H(F)\fl\H^0_{ab}(F,\H)
\fl\H^0_{ab}(F,\G)$. By the proof of lemma 3.7.1 of \cite{La-CSCB},
$\chi$ induces a character $\chi_{T_H}$ of $\T_H(F)$. Let us show that,
for every $f\in\Hecke_R$, $b_\xi(f)=\chi_H b_\xi(f)$ (this finishes the
proof of the lemma, because $\chi_H(\gamma'_H)\not=1$ for every
$\gamma'_H\in\H(F)$ that is conjugate to $\gamma_H$).
To do this, we imitate the proof of lemma 3.7.1 of \cite{La-CSCB}
and we show that
$\Tr(\pi(b_\xi(f)))=\Tr(\pi(\chi_H b_\xi(f)))$ for every $f\in\Hecke_R$ and
every unramified representation $\pi$ of $\H(F)$. 
Identify $\T_H$ to $\T_G$ by an admissible embedding, and let
$N:X_*(\Se_R)\fl X_*(\Se_H)$ be the norm morphism as in the proof of lemma
\ref{lemme:DL8} above.
By the proof of the second part of lemma \ref{lemme:DL8},
every function in $\Hecke_H=\C[\widehat{\Se}_H/\Omega_H]=\C[\Se_H]^{\Omega_H}$
that is in the image of $b_\xi$ is a linear combination of elements
$N(\mu)$, with $\mu\in X_*(\Se_R)$. 
If $z\in\widehat{\Se}_H$ and $\chi_z$ is the unramified character of
$\T_H(F)$ corresponding to $z$, let $\pi_z$ be the unramified representation
of $\H(F)$ obtained from $\chi_z$ ($\pi_z$ is the unique unramified
subquotient of $I(z)$, see eg \cite{Car} p 152). Finally, let
$z_0\in\widehat{\Se}_H$ be the element corresponding to the unramified
character $\chi_{T_H}$ of $\T_H(F)$. As $\chi_{T_H}$ is trivial on the norms,
$N(\mu)(z_0)=1$ for every $\mu\in X_*(\Se_R)$. 

Let $f\in\Hecke_R$. By lemma \ref{lemme:DL7}, for every $z\in
\widehat{\Se}_H$, $\Tr(\pi_z(b_\xi(f)))$ is a linear combination of
the $b_\xi(f)(\omega z)$, with $\omega\in\Omega_H$. Hence, by the discussion
above, for every $z\in\widehat{\Se}_H$, $\Tr(\pi_z(b_\xi(f)))=\Tr(\pi_
{zz_0}(b_\xi(f)))$; but $\Tr(\pi_{zz_0}(b_\xi(f)))=\Tr(\pi_z(\chi_H
b_\xi(f)))$, so $\Tr(\pi_z(b_\xi(f)))=\Tr(\pi_z(\chi_H b_\xi(f)))$. 
This implies that $b_\xi(f)=\chi_H b_\xi(f)$. 

\end{proof}

\begin{proofpn}{\ref{prop:utilisation_donnees_locales}}
Note that, by lemma \ref{lemme:groupes_qui_marchent}, there exists an
admissible embedding $\T_H\fl\G$ with image $\T_G$.

Let $(a_i^R(\pi))_{i\in I,\pi\in\Pi(R)}$ and $(a_i^H(\pi_H))_{i\in I,
\pi_H\in\Pi(\H)}$ be the local data. By the definition of local data,
it is enough to show that, for every $i\in I$ and every $f\in\Hecke_R$,
\[\sum_{\pi\in\Pi(R)}a_i^R(\pi)\Tr(\pi(f)A_\pi)=\sum_{\pi_H\in\Pi(\H)}
a_i^H(\pi_H)\Tr(\pi_H(b_\xi(f))).\]
Fix $i\in I$, and let $A$ be the distribution on $\Hecke_R$ defined by
\[A(f)=\sum_{\pi\in\Pi(R)}a_i^R(\pi)\Tr(\pi(f)A_\pi)-\sum_{\pi_H\in\Pi(\H)}
a_i^H(\pi_H)\Tr(\pi_H(b_\xi(f))).\]
We want to show that $A=0$.
The distribution $A$ is a sum of characters of $\Hecke_R$. In other
words, there exist $z_1,\dots,z_n\in Y_R$ such that $A$ is a linear combination
of the distributions $z_i\fle f(z_i)$; we may assume that $z_i$ and $z_j$ are
not $\Omega_R$-conjugate if $i\not=j$.
Write
\[A(f)=\sum_{i=1}^n c_i f(z_i),\]
with $c_1,\dots,c_n\in\C$. 
By the definition of local data and lemma \ref{lemme:DL1}, $A$ is a linear
combination of distributions $f\fle\Lambda(\gamma_H,f)$, with
$\gamma_H\in\H(F)$ elliptic semi-simple and strongly $\G$-regular.
By lemma \ref{lemme:DL3}, the distribution $A$ is tempered. By lemma 5.5 of
\cite{Cl-LF}, we may assume that $z_1,\dots,z_n\in Y_R^u$.
On the other hand, by lemma \ref{lemme:DL5}, the distribution $A$ is a
linear combination of distributions $f\fle\Tr_{\theta-c}(\pi(f)A_\pi)$
and $f\fle\Tr_c(\pi_H(b_\xi(f)))$, with $\pi\in\Pi(R)$ and $\pi_H\in\Pi(\H)$.
By lemma \ref{lemme:DL8}, there exists a non-empty open cone $C$ of
$X^*(Y_R)\otimes_\Z\R$ and $y_1,\dots,y_m\in Y_R$ such that the
restriction to $C\cap X^*(Y_R)$ of the function
$\lambda\fle A(f_\lambda)$ on $X^*(Y_R)$ is a linear combination of
the functions $\lambda\fle\lambda(y_i)$, $1\leq i\leq m$. By the explicit
description of the $y_i$ given in lemma \ref{lemme:DL8} and lemmas
\ref{lemme:DL6} and \ref{lemme:DL2}, we may assume that $y_1,\dots,y_m\in
Y_R-Y_R^u$. Let $d_1,\dots,d_m\in\C$ be such that
\[A(f_\lambda)=\sum_{i=1}^m d_i\lambda(y_i)\]
for every $\lambda\in C\cap X^*(Y_R)$.
Consider the characters $\varphi$ and $\varphi'$ of
$X^*(Y_R)$ defined by $\varphi(\lambda)=\sum\limits_{i=1}^n\sum\limits
_{\omega\in\Omega_R}c_i\lambda(\omega z_i)$ and $\varphi'(\lambda)=\sum\limits
_{i=1}^m d_i\lambda(y_i)$. Then $\varphi(\lambda)=\varphi'(\lambda)=A(f_
\lambda)$ if $\lambda\in C\cap X^*(Y_R)$. As $C\cap X^*(Y_R)$ generates the
group $X^*(Y_R)$, the characters $\varphi$ and $\varphi'$ are equal.
But the family $(\lambda\fle\lambda(z))_{z\in Y_R}$ of characters of
$X^*(Y_R)$ is free and $\{\omega z_i,1\leq i\leq n,\omega\in\Omega_R\}
\cap\{y_1,\dots,y_m\}=\varnothing$ (because the first set is included in
$Y_R^u$ and the second set is included in $Y_R-Y_R^u$), so
$\varphi=\varphi'=0$. By the linear independance of the characters
$\lambda\fle\lambda(\omega z_i)$, this implies that
$c_1=\dots=c_n=0$, hence, finally, that $A=0$.

\end{proofpn}

\section{Construction of local data}
\label{lemme_fondamental3}

The goal of this section is to construct local data. The method is global and
uses the trace formula.
The first thing to do is to show that there exists a global situation
that gives back the situation of \ref{lemme_fondamental1} at one place.

\begin{sublemme}\label{lemme:situation_globale} Let $F$, $E$, $\G$,
$(\H,s,\eta_0)$ and $\eta$ be as in \ref{lemme_fondamental1}. Assume that
there exists a finite unramified extension $K$ of $E$ such that the groups $\G$
and $\H$ split over $K$ and that the morphisms $\eta$ and $\xi$ come from
morphisms $\widehat{\H}\rtimes\Gal(K/F)\fl\widehat{\G}\rtimes\Gal(K/F)$ and
$\widehat{\H}\rtimes\Gal(K/F)\fl\widehat{R}\rtimes\Gal(K/F)$ (if $\G$ is
adjoint, then such a $K$ exists by remark \ref{rq:xi_sur_extension}).
Then, for everey $r\in\Nat^*$, there exist a number field $k_F$, finite Galois
extensions $k_K/k_E/k_F$, a finite set $S_0$ of finite places of $k_F$,
an element $v_0$ of $S_0$,
connected reductive groups $\underline{\G}$ and $\underline{\H}$ over $k_F$
and $L$-morphisms $\underline{\eta}:{}^L\underline{\H}\fl{}^L\underline{\G}$
and $\underline{\xi}:{}^L\underline{\H}\fl\underline{R}$, where $\underline{R}=
R_{k_E/k_F}\underline{\G}_{k_E}$,
such that :
\begin{itemize}
\item[(i)] The groups $\underline{\H}$ and $\underline{\G}$ are quasi-split
over $k_F$ and split over $k_K$.
\item[(ii)] The set $S_0$ has $r$ elements. Let $v\in S_0$. Then the place
$v$ is inert in $k_K$, and there are isomorphisms
$k_{F,v}\simeq F$, $k_{E,v}\simeq E$, $k_{K,v}\simeq K$, $\underline{\G}_
{v}\simeq\G$, $\underline{\H}_{v} \simeq\H$. The group $\underline{\G}$
has elliptic maximal tori $\T$ that stay elliptic over $k_{F,v}$.
Moreover, the obvious morphism
$\Gal(k_{K,v}/k_{F,v})\fl\Gal(k_K/k_F)$ is an isomorphism (in particular,
the extension $k_E/k_F$ is cyclic).
\item[(iii)] $k_F$ is totally imaginary.
\item[(iv)] $(\underline{\H},{}^L\underline{\H},s,\eta)$ are endoscopic data
for $(\underline{\G},1,1)$, and $(\underline{\H},{}^L\underline{\H},
t,\underline{\xi})$ are endoscopic data for $(\underline{R},
\underline{\theta},1)$, where $\underline{\theta}$ is the automorphism of
$\underline{R}$ defined by the generator of $\Gal(k_E/k_F)$ given by
the isomorphism $\Gal(E/F)\simeq\Gal(k_{E,v_0}/k_{F,v_0})\iso\Gal(k_E/k_F)$
of (ii).
\item[(v)] For every $v\in S_0$, $\underline{\eta}_{v}$ corresponds to $\eta:
{}^L\H\fl{}^L\G$ and $\underline{\xi}_v$ corresponds to $\xi:{}^L\H\fl{}^LR$
(by the isomorphisms of (ii)).
For every infinite place $v$ of $k_F$, the morphism $\underline{\eta}_v:
\widehat{\underline{\H}}\times W_\C=\widehat{\H}\times W_\C\fl\widehat{
\underline{\G}}\times W_\C=\widehat{\G}\times W_\C$ is equal to $\eta_0\times
id_{W_\C}$.
\end{itemize}
Moreover :
\begin{itemize}
\item[(vi)] There exist infinitely many places of $k_F$ that split totally in
$k_K$.
\item[(vii)] For every finite set $S$ of places of $k_F$ such that
$S_0\not\subset S$, the group $\underline{\H}(k_F)$ is dense in
$\prod\limits_{v\in S}\underline{\H}(k_{F,v})$. The same statement is true
if $\underline{\H}$ is replaced by $\underline{\G}$, $\underline{R}$ or by
a torus of $\underline{\H}$, $\underline{\G}$ or $\underline{R}$.
% \item[(viii)] On a $\Ker^1(F,\H)=\Ker^1(F,\G)=\Ker^1(F,R)=\Ker^1_{ab}(F,
% \H\sous R)=\{1\}$, et $\Ker^1(F,\T)=\{1\}$ pour tout tore de $\H$, $\G$ ou
% $R$ (cf \cite{K-STF:CTT} 3.4.1 pour la définition de $\Ker^1$ et
% \cite{La-CSCB} 1.8
% pour celle de $\Ker^1_{ab}(F,\H\sous R)$).

\end{itemize}
\end{sublemme}

\begin{proof} If $r=1$, the existence of $k_F$, $k_E$, $k_K$,
$\underline{\G}$, $\underline{\H}$ and $S_0=\{v_0\}$ satisfying (i), (ii) and
(iii) is a consequence of the proof of proposition 11.1 of \cite{Wa1}.
As in \cite{Cl-LF} p 293, we pass from the case $r=1$ to the general case
by replacing $k_F$ by an extension of degree $r$ where $v_0$ splits totally
(this extension is necessarily linearly disjoint from $k_K$, because
$v_0$ is inert in $k_K$).
By the last sentence of (ii), $\eta$ gives a $L$-morphism
$\widehat{\underline{\H}}\rtimes\Gal(k_K/k_F)\fl\widehat{\underline{\G}}
\rtimes\Gal(k_K/k_F)$ and $\xi$ gives a $L$-morphism
$\widehat{\underline{\H}}\rtimes\Gal(k_K/k_F)\fl\widehat{\underline{R}}
\rtimes\Gal(k_K/k_F)$. Take as $\underline{\eta}$ and $\underline{\xi}$
the $L$-morphisms ${}^L\underline{\H}\fl{}^L\underline{\G}$
and ${}^L\underline{\H}\fl{}^L\underline{R}$ that make the obvious
diagrams commute. Then (iv) and (v) are clear.
Point (vi) follows from the \v{C}ebotarev density theorem
(cf \cite{Ne}, chapter VII, theorem 13.4 and in particular corollary 13.6). 
As all the places of $S_0$ are inert in $k_K$, (vii) follows from (b) of
lemma 1 of \cite{KRo}. 
% Enfin, le point (viii) est prouvé dans \cite{Cl-LF}, p 293.

\end{proof}

The main result of this section is the next proposition.

\begin{subproposition}\label{prop:existence_donnees_locales} Assume that
$\G$ is adjoint and that the endoscopic triple $(\H,s,\eta_0)$ is elliptic.
Let $k_F$, $k_E$, etc, be as in lemma \ref{lemme:situation_globale} above,
with $r=3$.
Assume that, for almost every place $v$ of $k_F$, the fundamental lemma
for the unit element of the Hecke algebra is known for $(\underline{R}
_v,\underline{\theta}_v,1)$ and $(\underline{\H}_v,{}^L\underline{\H}_v,t,
\underline{\xi}_v)$ (at almost every place $v$ of $F$, the local situation is
as in remark \ref{rq:E_pas_un_corps}).

Then there exist local data for $R$ and $(\H,{}^L\H,t,\xi)$.

\end{subproposition}

The proposition is proved at the end of this section, after a few lemmas.

We will need a simple form of the trace formula, due originally to
Deligne and Kazhdan (see the article \cite{He} of Henniart, sections
4.8 et 4.9, for the untwisted case, and lemma I.2.5 of the book \cite{AC}
of Arthur and Clozel for the twisted case). The next lemma is the obvious
generalization (to groups that are not necessarily $\GL_n$) of lemma I.2.5
of \cite{AC}, and the proof of this lemma applies without any change
(in condition (3') on page 14 of \cite{AC}, the assumption that the
functions $\phi_{w_i}$ are all coefficients of the \emph{same} supercuspidal
representation is not necessary).

\begin{sublemme}\label{lemme:formule_traces_simple}
\index{simple trace formula of Deligne and Kazhdan}
Let $F$ be a number
field, $E/F$ be a cyclic extension of degree $d$ and $\G$ be a connected
adjoint group over $F$. Set $R=R_{E/F}\G_E$, fix a generator of
$\Gal(E/F)$ and let $\theta$ be the automorphism of $R$ induced by this
generator. Let $\phi\in C_c^\infty(R(F))$. Denote by $r(\phi)$ the
endomorphism of $L^2:=L^2(R(F)\sous R(\Ade_F))$ obtained by making $\phi$
act by right convolution, and let $I_\theta$ be the endomorphism
$f\fle f\circ\theta^{-1}$ of $L^2$.
Assume that :
\begin{itemize}
\item[(0)] $\phi=\bigotimes\limits_v\phi_v$, where the tensor product is
taken over the set of places $v$ of $F$ and
$\phi_v\in C_c^\infty(R(F_v))$ for every $v$; moreover, at almost every
finite place $v$ where $R$ is unramified, $\phi_v$ is the characteristic
function of a hyperspecial maximal compact subgroup of $R(F_v)$.
\item[(1)] There exists a finite place $v$ of $F$ that splits totally in
$E$ and such that, on $R(F_v)\simeq\G(F_v)^d$, $\phi_v=\phi_1\otimes\dots
\otimes\phi_d$, where the $\phi_i\in C_c^\infty(\G(F_v))$ are supercuspidal
functions (in the sense of \cite{He} 4.8).
\index{supercuspidal function}
\item[(2)] There exists a finite place $v$ of $F$ such that the support of
$\phi_v$ is contained in the set of $\theta$-elliptic $\theta$-semi-simple
and strongly $\theta$-regular elements of $R(F_v)$.

\end{itemize}
Then $r(\phi)I_\theta$ sends $L^2$ to the subspace of cuspidal functions
(in particular, the endomorphism $r(\phi)I_\theta$ has a trace), and
\[\Tr(r(\phi)I_\theta)=\sum_\delta\vol(\G_{\delta\theta}(F)\sous\G_{\delta
\theta}(\Ade_E))O_{\delta\theta}(\phi),\]
where the sum is taken over the set of $\theta$-conjugacy classes
of $\theta$-elliptic $\theta$-semi-simple and strongly $\theta$-regular
$\delta\in R(F)$.

\end{sublemme}

\begin{sublemme}\label{lemme:coeff_de_supercuspidales}(\cite{Ha} lemma 5.1)
Let $F$, $\G$ and $\H$ be as in \ref{lemme_fondamental1} (in particular,
$F$ is a non-archimedean local field, $\G$ is an unramified group over $F$ and
$\H$ is an unramified endoscopic group of $\G$).
Assume that the centers of $\G$ and $\H$ are anisotropic.
Let $\T$ be an unramified elliptic maximal torus of $\H$; assume that $\T$
is defined over $\Of_F$ and that $\T(F)\subset\H(\Of_F)$. Let $j:\T\fl\G$
be an admissible embedding defined over $\Of_F$. Set
$N=\Nor_{\G(F)}(j(\T(F)))$, and make $N$ act on $\T(F)$ via $j$.

Then there exist functions $f\in C_c^\infty(\G(F))$ and $f_H\in C_c^\infty
(\H(F))$ such that :
\begin{bulletlist}
\item $f$ and $f_H$ are supercuspidal (in the sense of \cite{He} 4.8;
in particular, a linear combination of coefficients of supercuspidal
representations is a supercuspidal function);
\item the function $\gamma\fle O_\gamma(f)$ (resp. $\gamma_H\fle
SO_{\gamma_H}(f_H)$) on $\G(F)$ (resp. $\H(F)$) is not identically zero, and
its support is contained in the set of semi-simple strongly regular
(resp. strongly $\G$-regular) elements that are conjugate to an element of
$j(\T(F))$ (resp. $\T(F)$);
\item the function $\gamma_H\fle O_{\gamma_H}(f_H)$ on $\T(F)$
is invariant under the action of $N$.

\end{bulletlist}
\end{sublemme}

The next two lemmas will be useful to construct transfers (and inverse
transfers) of certain functions. The first lemma is a particular case of a
theorem of Vignéras (theorem A of \cite{Vi}).

\begin{sublemme}\label{lemme:caracterisation_integrales_orbitales}
Let $F$ be a non-archimedean local field and $\G$ be a connected reductive
group over $F$. Denote by $\G(F)_{ss-reg}$ the set of semi-simple strongly
regular elements of $\G(F)$. Let $\Gamma:\G(F)_{ss-reg}\fl\C$ be a function
that is invariant by $\G(F)$-conjugacy and such that, for every
$\gamma\in\G(F)_{ss-reg}$, the restriction of $\Gamma$ to
$\G_\gamma(F)\cap\G(F)_{ss-reg}$ is locally constant with compact support.
Then there exists $f\in C_c^\infty(\G(F)_{ss-reg})$ such that, for every
$\gamma\in\G(F)_{ss-reg}$, $\Gamma(\gamma)=O_\gamma(f)$.

\end{sublemme}

We will need a twisted variant of this lemma and some consequences of it.
In the next lemma, $F$ is still a non-archimedean local field and $\G$
a connected reductive group over $F$. We also assume that $F$ is of
characteristic $0$.
Let $E$ be a finite étale $F$-algebra
such that $\Aut_F(E)$ is cyclic. Set $R=R_{E/F}\G_E$, fix a generator of
$\Aut_F(E)$ and denote by $\theta$ the automorphism of $R$ induced by this
generator (so the situation is that of example \ref{ex:groupes_non_connexes},
except that $E$ does not have to be a field).
Use the definitions of \ref{GL_n_applications1}.
Let $\delta\in R(F)$ be $\theta$-semi-simple and strongly $\theta$-regular,
and write $\T=R_{\delta\theta}$.
As in \cite{La-CSCB} 1.8, set
\[\Dgoth(\T,R;F)=\Ker(\H^1(F,\T)\fl\H^1(F,R)).\]
\index{DTRF@$\Dgoth(\T,R;F)$}
As $F$ is local and non-archimedean, the pointed set $\Dgoth(\T,\G;F)$ is
canonically isomorphic to an abelian group (cf \cite{La-CSCB} lemma 1.8.3);
so we will see $\Dgoth(\T,R;F)$ as an abelian group.
If $\delta'\in R(F)$ is stably $\theta$-conjugate to $\delta$, then it defines
an element $inv(\delta,\delta')$ of $\Dgoth(\T,R;F)$. The map
$\delta'\fle inv(\delta,\delta')$ induces a bijection from the set of
$\theta$-conjugacy classes in the stable $\theta$-conjugacy class of $\delta$
to the set $\Dgoth(\T,R;F)$ (cf \cite{La-CSCB} 2.3). 
Remember also (cf \cite{A-LBWOI} \S1) that $\T$ is a torus of $R$ and that,
if $(\T(F)\delta)_{reg}$ is the set of strongly $\theta$-regular elements of
$\T(F)\delta$, then the map $u:(\T(F)\delta)_{reg}\times\T(F)\sous R(F)
\fl R(F)$, $(g,x)\fle x^{-1}g\theta(x)$, is finite on its image and open.

Now forget about $\delta$ and
fix a maximal torus $\T$ of $\G$ (that can also be seen as a torus of
$R$ via the obvious embedding $\G\subset R$).
Let $\Omega$ be the set of $\theta$-semi-simple strongly $\theta$-regular
$\delta\in R(F)$ such that there exists $x\in R(\overline{F})$ with
$R_{\delta\theta}=x\T x^{-1}$.
If $\kappa$ is a character of $\Dgoth(\T,R;F)$, $f\in C_c^\infty(\Omega)$ and
$\delta\in\Omega$ is such that $R_{\delta\theta}=\T$, set
\[O_{\delta\theta}^\kappa(f)=\sum_{\delta'}\langle inv(\delta,\delta'),\kappa
\rangle O_{\delta'\theta}(f),\]
\index{Okappa@$O_{\delta\theta}^\kappa$\quad twisted $\kappa$-integral
orbital}where the sum is taken over the set of
$\theta$-conjugacy classes in the stable $\theta$-conjugacy class of
$\delta$. Then $O_{\delta\theta}^\kappa(f)$ is a (twisted) $\kappa$-orbital
integral of $f$. (This definition is a particular case of
\cite{La-CSCB} 2.7).

\begin{sublemme}\label{lemme:caracterisation_kappa_integrales_orbitales}
Let $\T$ and $\Omega$ be as above.
Let $\Gamma:\Omega\fl\C$ be a function that is invariant by
$\theta$-conjugacy and such that, for every $\delta\in\Omega$, the restriction
of $\Gamma$ to $R_{\delta\theta}(F)\delta\cap\Omega$ is locally constant with
compact support. Then there exists $f\in C_c^\infty(\Omega)$ such that,
for every $\delta\in\Omega$, $\Gamma(\delta)=O_{\delta\theta}(f)$. 

Let $\kappa$ be a character of $\Dgoth(\T,R;F)$. Assume that, for every
$\delta\in\Omega$ such that $R_{\delta\theta}=\T$, for every
$\delta'\in\Omega$ that is stably $\theta$-conjugate to $\delta$,
$\Gamma(\delta)=\langle inv(\delta,\delta'),\kappa\rangle\Gamma(\delta')$.
Then there exists $g\in C_c^\infty(\Omega)$ such that
\begin{itemize}
\item[(a)] for every $\delta\in\Omega$ such that $R_{\delta\theta}=\T$,
$\Gamma(\delta)=O_{\delta\theta}^\kappa(g)$;
\item[(b)] for every character $\kappa'$ of $\Dgoth(\T,R;F)$ such that
$\kappa'\not=\kappa$ and for every $\delta\in\Omega$ such that
$R_{\delta\theta}=\T$, $O_{\delta\theta}^{\kappa'}(g)=0$.

\end{itemize}

Moreover, for every character $\kappa$ of $\Dgoth(\T,R;F)$, for every
$\delta\in\Omega$ and every open neighbourhood $\Omega'$ of $\delta$,
there exists a function $\Gamma$ satisfying the conditions above and
such that $\Gamma(\delta)\not=0$ and that $\Gamma(\delta')=0$ if $\delta'$
is not stably $\theta$-conjugate to an element of $\Omega'$.

\end{sublemme}

\begin{proof} Let $\T=\T_1,\dots,\T_n$ be a system of representatives of the
set of $R(F)$-conjugacy classes of tori of $R$ (defined over $F$)
that are equal to a $x\T x^{-1}$, with $x\in R(\overline{F})$.
For every $i\in\{1,\dots,n\}$, denote by $\Omega_i$ the set of
$\delta\in\Omega$ such that $R_{\delta\theta}$ is $R(F)$-conjugate to
$\T_i$. The $\Omega_i$ are pairwise disjoint open subsets of $\Omega$, and
$\Omega=\bigcup\limits_{i=1}^n\Omega_i$. Let $i\in\{1,\dots,n\}$. If
$\delta\in\Omega_i$ is such that $R_{\delta\theta}=\T_i$, let $u_\delta$
be the function $(\T_i(F)\delta)_{reg}\times\T_i(F)\sous R(F)\fl R(F),
(x,y)\fle y^{-1}xy$. Then $\Omega_i$ is the union of the images of the
$u_\delta$, this images are open, and two of these images are either
equal or disjoint. So there exists a finite family $(\delta_{ij})_{j\in J_i}$
of elements of $\Omega_i$ such that $R_{\delta_{ij}\theta}=\T_i$ for every $j$
and that $\Omega_i=\coprod\limits_{j\in J_i}Im(u_{\delta_{ij}})$. Write
$u_{ij}=u_{\delta_{ij}}$, $\Omega_{ij}=Im(u_{ij})$ and $\Gamma_{ij}=\ungras_
{\Omega_{ij}}\Gamma$. Then the functions $\Gamma_{ij}$ are invariant by
$\theta$-conjugacy, and $\Gamma=\sum\limits_{i,j}\Gamma_{ij}$.

For every $i,i'\in\{1,\dots,n\}$, let $A(i,i')$ be the (finite) set of
isomorphisms $\T_i\iso\T_{i'}$ (over $F$) of the form
$\Int(x)$, with $x\in R(\overline{F})$. If $i=i'$, write $A(i)=A(i,i')$.
Let $i\in\{1,\dots,n\}$. For every $j\in J_i$, $i'\in\{1,\dots,n\}$, $j'
\in J_{i'}$ and $a\in A(i,i')$, the map $\{x\in\T_i(F)|x\delta_{ij}\in
\T_i(F)_{reg}\}\fl\Dgoth(\T_i,R;F)$,
$x\fle inv(x\delta_{ij},a(x)\delta_{i'j'})$ is locally constant. For every
$j\in J_i$, the support of $\Gamma_{ij}\circ u_{ij}$ is contained in a set of
the form $\omega\times\T_i(F)\sous R(F)$, with $\omega$ a compact subset of
$(\T_i(F)\delta_{ij})_{reg}$. So it is easy to see that there exist open
compact subsets $\omega_{ik}$, $k\in K_i$, of $\T_i(F)$, and functions
$\Gamma_{ijk}\in C_c^\infty(\Omega_{ij})$, invariant by $\theta$-conjugacy,
such that :
\begin{itemize}
\item[$(1)$] for every $j\in J_i$ and $k\in K_i$, $\omega_{ik}\delta_{ij}
\subset(\T_i(F)\delta_{ij})_{reg}$, and the support of the function
$\Gamma_{ijk}$ is contained in
$u_{ij}(\omega_{ik}\delta_{ij}\times\T_i(F)\sous R(F))$;
\item[(2)] for every $k\in K_i$, the images of the $\omega_{ik}$ by the
elements of $A(i)$ are pairwise disjoint;
\item[$(3)$] for every $j\in J_i$, $\Gamma_{ij}=\sum\limits_{k\in K_i}
\Gamma_{ijk}$.

\end{itemize}
Let $i\in\{1,\dots,n\}$, $j\in J_i$ and $k\in K_i$. By point (2) above,
the restriction of $u_{ij}$ to $\omega_{ik}\delta_{ij}\times
\T_i(F)\sous R(F)$ is injective. 
Let $U_i$ be an open compact subset of volume $1$ of $\T_i(F)\sous\G(F)$.
Denote by $f_{ijk}$ the product of $\Gamma_{ijk}$ and of the characteristic
function of $u_{ij}(\omega_{ik}\delta_{ij}\times U_i)$.
Then $f_{ijk}\in C_c^\infty(\Omega)$ and, for every $\delta\in\Omega$,
$O_{\delta\theta}(f_{ijk})=\Gamma_{ijk}(\delta)$.
So the function $f:=\sum\limits_{ijk}f_{ijk}$ satisfies the condition of
the first statement of the lemma.

Let $\kappa$ be a character of $\Dgoth(\T,R;F)$. Assume that $\Gamma$
satisfies the condition of the second statement of the lemma.
If $\kappa'$ is a character of $\Dgoth(\T,R;F)$ and $\delta\in\Omega$
is such that $R_{\delta\theta}=\T$, then
\[\begin{array}{rcl}\displaystyle{O_{\delta\theta}^{\kappa'}(f)=\sum_{\delta'}
\langle inv(\delta,\delta'),\kappa'\rangle O_{\delta'\theta}(f)} & = &
\displaystyle{\sum_{\delta'}\langle inv(\delta,\delta'),\kappa'\rangle\Gamma
(\delta')}\\
 & = & \displaystyle{\Gamma(\delta)\sum_
{\delta'}\frac{\langle inv(\delta,\delta'),\kappa'\rangle }{\langle inv(\delta,
\delta'),\kappa\rangle },}\end{array}\]
where the sum is taken over the set of $\theta$-conjugacy classes in the stable
$\theta$-conjugacy class of $\delta$.
So we can take $g=|\Dgoth(\T,R;F)|^{-1}f$.

We show the last statement of the lemma. Let $\kappa$ be a character of
$\Dgoth(\T,R;F)$. Choose (arbitrarily) an element $j_0$ of $J_1$, and write
$\delta_1=\delta_{1,j_0}$. For every $i\in\{1,\dots,n\}$, let $J'_i$ be the set
of $j\in J_i$ such that $\delta_{ij}$ is stably $\theta$-conjugate to an
element of $\T_1(F)\delta_1$; after translating (on the left) $\delta_{ij}$
by an element of $\T_i(F)$, we may assume that $\delta_{ij}$ is stably
$\theta$-conjugate to $\delta_1$ for every $j\in J'_i$. For every $i\in\{1,
\dots,n\}$ and $j\in J'_i$, choose $x_{ij}\in R(\overline{F})$ such that
$\delta_{ij}=x_{ij}\delta_1\theta(x_{ij})^{-1}$,
and let $a_{ij}$ be the element of
$A(1,i)$ induced by $\Int(x_{ij})$.
Let $\omega\subset\T_1(F)$ be an open compact subset such that :
\begin{bulletlist}
\item $1\in\omega$;
\item for every $i\in\{1,\dots,n\}$, the images of $\omega$ by the elements
of $A(1,i)$ are pairwise disjoint;
\item for every $i\in\{1,\dots,n\}$, $j\in J'_i$ and $a\in A(1,i)$,
$a(\omega)\delta_{ij}\subset (\T_i(F)\delta_{ij})_{reg}$, and the function
$x\fle \langle inv(x\delta_1,a(x)\delta_{ij}),\kappa\rangle$ is constant on
$\omega$.

\end{bulletlist}

For every $i\in\{1,\dots,n\}$ and $j\in J'_i$, let $\Gamma_{ij}$ be the product
of the characteristic function of $u_{ij}(a_{ij}(\omega)\delta_{ij}\times
\T_i(F)\sous R(F))$ and of $\langle inv(y\delta_1,a_{ij}(y)\delta_{ij}),\kappa
\rangle^{-1}$, where $y$ is any element of $\omega$.
Set $\Gamma=\sum\limits_{i,j}\Gamma_{ij}$. Let $\delta\in\Omega$ be such that
$R_{\delta\theta}=\T$. Then $\Gamma(\delta)=|A(1)|$ if $\delta$ is 
$\theta$-conjugate to an element of $\omega\delta_1$, and $\Gamma(\delta)=0$
otherwise (in particular, $\Gamma$ is not identically zero).
Let $\delta\in R(F)$ be stably $\theta$-conjugate to $\delta_1$. There exists
a unique pair $(i,j)$, with $i\in\{1,\dots,n\}$ and $j\in J'_i$, such that
$\delta$ is $\theta$-conjugate to an element of $\T_i(F)\delta_{ij}$.
If $\delta$ is not stably $\theta$-conjugate to an element of
$\omega\delta_1$, then $\delta$ is not $\theta$-conjugate to an element of
$a_{ij}(\omega)\delta_{ij}$, and $\Gamma(\delta)=0$. Otherwise,
$\Gamma(\delta)=\langle inv(\delta_1,\delta),\kappa\rangle^{-1}|A(1)|=
\langle inv(\delta_1,\delta),\kappa\rangle^{-1}\Gamma(\delta_1)$. 

Let $\delta\in\Omega$. After changing the order of the $\T_i$ and choosing
other representatives for the $\delta_{1j}$, we may assume that the fixed
$\delta_1$ is $\delta$. As it is always possible to replace
$\omega$ by a smaller open compact (containing $1$), this proves the last
statement of the lemma.

\end{proof}

The two lemmas above have the following consequence :

\begin{sublemme}\label{lemme:construction_transfert} Assume that
$F$, $E$, $\G$, $R$ and $\theta$ are as in lemma
\ref{lemme:caracterisation_kappa_integrales_orbitales}.
Let $(\H,s,\eta_0)$ be an endoscopic triple for $\G$. Assume that
this local situation comes from a global situation as in
lemma \ref{lemme:situation_globale}. In particular, $\H$ is
the first element of endoscopic data $(\H,{}^L\H,t,\xi)$ for $(R,\theta,1)$.
Let $\Delta_\xi$ be the transfer factors defined by $\xi$ (with any
normalization). Then :
\begin{itemize}
\item[(i)] Every function $f\in C_c^\infty(R(F))$ with support in the set of
$\theta$-semi-simple strongly $\theta$-regular elements admits a transfer
to $\H$.
\item[(ii)] Let $\T_H$ be a maximal torus of $\H$. Choose an admissible
embedding $j:\T_H\fl\G$, and make $N:=\Nor_{\G(F)}(j(\T_H(F)))$ act on
$\T_H(F)$ via $j$. Let $f_H\in C_c^\infty(\H(F))$ be a function with support
in the set of strongly regular elements that are stably conjugate to an
element of $\T_H(F)$. Assume that the function $\T_H(F)\fl\C$, $\gamma_H\fle
O_{\gamma_H}(f_H)$, is invariant under the action of $N$. Then there
exists $f\in C_c^\infty(R(F))$ such that $f_H$ is a transfer of $f$ to $\H$.

\end{itemize}
\end{sublemme}

The notion of transfer (or of ``matching functions'') in that case of defined
in \cite{KS} 5.5.
\index{transfer (twisted case)}

\begin{proof} To prove (i), define a function $\Gamma_H$ on the set of
semi-simple strongly $\G$-regular elements of $\H(F)$ by
$\Gamma_H(\gamma_H)=\sum\limits_\delta \Delta_\xi(\gamma_H,\delta)O_{\delta
\theta}(f)$, where the sum is taken over the set of $\theta$-conjugacy classes
of $R(F)$, and apply lemma \ref{lemme:caracterisation_integrales_orbitales}
to $\Gamma_H$. To show (ii), construct a function $\Gamma$ on the set
of $\theta$-semi-simple strongly $\theta$-regular elements of $R(F)$ in the
following way : If there does not exist any $\gamma_H\in\H(F)$ such that
$\Delta_\xi(\gamma_H,\delta)\not=0$, set $\Gamma(\delta)=0$; if there exists
$\gamma_H\in\H(F)$ such that $\Delta_\xi(\gamma_H,\delta)\not=0$, set
$\Gamma(\delta)=\Delta_\xi(\gamma_H,\delta)^{-1}SO_{\gamma_H}(f_H)$.
The function $\Gamma$ is well-defined by the assumption on $f_H$ (and lemma
5.1.B of \cite{KS}). So (ii) follows from theorem 5.1.D of \cite{KS}
and from lemma \ref{lemme:caracterisation_kappa_integrales_orbitales}. 

\end{proof}

\begin{subremarque}\label{rq:comparaison_Kgoth}
We need to be able to compare the groups of endoscopic characters
of \cite{KS} and \cite{La-CSCB}. In the situation of the lemma above, but
with $F$ global or local (and allowed to be archimedean), if $\T_R$ is a
$\theta$-stable maximal torus of $R$ coming from a torus $\T$ of $\G$,
Labesse defined groups $\Kgoth(\T,R;F)_1\subset\Kgoth(\T,R;F)$ (\cite{La-CSCB}
1.8) and Kottwitz and Shelstad defined groups $\Kgoth(\T_R,\theta,R)_1\subset
\Kgoth(\T_R,\theta,R)$ (\cite{KS} 6.4; Kottwitz and Shelstad assume that
$F$ is a number field, but it is possible to write the same definitions if
$F$ is local, erasing of course the quotient by $\Ker^1$ in the definition of
$\Kgoth_1$).
\index{KTRF1@$\Kgoth(\T,R;F)_1$}
\index{KTRF@$\Kgoth(\T,R;F)$}
\index{KTRF1@$\Kgoth(\T_R,\theta,R)_1$}
\index{KTRF@$\Kgoth(\T_R,\theta,R)$}
As we are interested in endoscopic data for the triple
$(R,\theta,1)$ (whose third element, that could in general be an element of
$\Ho^1(W_F,Z(\widehat{\G}))$, is trivial here), we must use the group
$\Kgoth(\T_R,\theta,\T)_1$ (cf \cite{KS} 7.1 and 7.2) to parametrize this
data. Labesse showed that the groups $\Kgoth(\T,R,F)$ and
$\Kgoth(\T_R,\theta,R)$ are canonically isomorphic (cf the end of
\cite{La-CSCB} 2.6). Using the techniques of \cite{La-CSCB} 1.7, it is easy
to see that this isomorphism identifies $\Kgoth(\T,R;F)_1$ and $\Kgoth(\T_R,
\theta,R)_1$.

\end{subremarque}

The next lemma explains what happens if $E=F^d$.

\begin{sublemme}\label{lemme:E_egal_F_puissance_d} Let $F$ be a local or
global field, $\G$ be a connected reductive group over $F$ and $d\in\Nat^*$.
Set $R=\G^d$, and let $\theta$ be the automorphism of $R$ that sends
$(g_1,\dots,g_d)$ to $(g_2,\dots,g_d,g_1)$. Then :
\begin{itemize}
\item[(i)] The set of equivalence classes of endoscopic data for $(R,\theta,1)$
is canonically isomorphic to the set of equivalence classes of endoscopic
data for $(\G,1,1)$.
\item[(ii)] Let $\phi\in C_c^\infty(R(F))$. Assume that
$\phi=\phi_1\otimes\dots\otimes\phi_d$, with $\phi_1,\dots,
\phi_d\in C_c^\infty(\G(F))$. Then, for every $\gamma=(\gamma_1,\dots,
\gamma_d)\in R(F)$,
\[O_{\gamma\theta}(\phi)=O_{\gamma_1\dots\gamma_d}(\phi_1*\dots*\phi_d)\]
(provided, of course, that the measures are normalized in compatible ways).

\end{itemize}
\end{sublemme}

\begin{proof} Point (ii) is a particular case of \cite{AC} I.5.
Point (i) is almost obvious. We explain how the isomorphism is constructed.
The dual group of $R$ is $\widehat{R}=\widehat{\G}^d$, with the diagonal
action of $\Gal(\overline{F}/F)$, so the diagonal embedding $\widehat{\G}\fl
\widehat{R}$ extends in an obvious way to a $L$-morphism
$\eta:{}^L\G\fl{}^LR$. If $(\H,\Hcal,s,\xi)$ are endoscopic data for
$(\G,1,1)$, it defines endoscopic data $(\H,\Hcal,\eta(s),\eta\circ\xi)$ for
$(R,\theta,1)$. Conversely, let $(\H,\Hcal,s,\xi)$ be endoscopic data for
$(R,\theta,1)$. Write $\xi(h\rtimes w)=(\xi_1(h\rtimes w),\dots,\xi_d
(h\rtimes w))\rtimes w$, with $h\rtimes w\in\Hcal\simeq\widehat{\H}
\rtimes W_F$, and $s=(s_1,\dots,s_d)$. Let $\xi_G:\Hcal\fl{}^L\G$,
$h\rtimes w\fle\xi_1(h\rtimes w)\rtimes w$.
Then $(\H,\Hcal,s_1\dots s_d,\xi_G)$ are endoscopic data for $(\G,1,1)$.

\end{proof}

The next lemma is the analog of a statement proved in \cite{Ha}, p 20-22.
It is proved exactly in the same way, using the twisted version of the
Paley-Wiener theorem (cf the article \cite{DeM} of Delorme and Mezo)
instead of the untwisted version.
(The statment on the support of functions in $E$ is shown by using the
control over the support of the functions given by theorem 3 of \cite{DeM}.)

\begin{sublemme}\label{lemme:place_infinie} Let $\G$ be a connected reductive
group over $\C$, $(\H,s,\eta_0)$ be an endoscopic triple for $\G$ and $d\in
\Nat^*$. Set $R=\G^d$, and denote by $\theta$ the automorphism of $R$ defined
by $\theta(g_1,\dots,g_d)=(g_2,\dots,g_d,g_1)$.
Let $\eta=\eta_0\times id_{W_\C}:
{}^L\H\fl{}^L\G$ be the obvious extension of $\eta_0$, and $\xi$ be the
composition of $\eta$ and of the obvious (``diagonal'') embedding
${}^L\G\fl{}^LR$. 
Fix maximal compact subgroups $\K_G$ and $\K_H$ of $\G(\C)$ and $\H(\C)$,
write $\K_R=\K_G^d$ and denote by $C_c^\infty(\G(\C),\K_G)$ (resp.
$C_c^\infty(\H(\C),\K_H)$, $C_c^\infty(R(\C),\K_R)$) the space of $C^\infty$
functions with compact support  and $\K_G$-finite (resp. $\K_H$-finite,
$\K_R$-finite) on $\G(\C)$ (resp. $\H(\C)$, $R(\C)$). Let
$\Pi(\H)$ (resp. $\Pi_\theta(R)$) be the set of isomorphism classes of
irreducible unitary representations of $\H(\C)$ (resp. of $\theta$-stable
irreducible unitary representations of $R(\C)$), and $\Pi_{temp}(\H)$
(resp. $\Pi_{\theta-temp}(R)$) be the subset of tempered representations.
For every $\pi\in\Pi_\theta(R)$, choose a normalized intertwining operator
$A_\pi$ on $\pi$.
For every $\pi\in\Pi_{\theta-temp}(R)$, let $\Pi_H(\pi)$ be the set of
$\pi_H\in\Pi_{temp}(\pi)$ whose functorial transfer to $R$ is $\pi$ (so
$\pi_H$ is in $\Pi_H(\pi)$ if and only if there exists a Langlands
parameter $\varphi_H:W_\C\fl{}^L\H$ of $\pi_H$ such that $\xi\circ\varphi_H$
is a Langlands parameter of $\pi$).
Let $N:C_c^\infty(R(\C),\K_R)=C_c^\infty(\G(\C),\K_G)^{\otimes d}\fl
C_c^\infty(\G(\C),\K_G)$ be the morphism of $\C$-algebras such that, for every
$f_1,\dots,f_d\in C_c^\infty(\G(\C),\K_G)$, $N(f_1\otimes\dots\otimes f_d)=
f_1*\dots *f_d$.

Then there exists a vector space $E\subset C_c^\infty(R(\C),\K_R)$ and a
compact subset $C$ of $\H(\C)$ such that :
\begin{itemize}
\item[(i)] There exists $f\in E$ and a transfer $f^H\in C_c^\infty(\H(\C),
\K_H)$ of $N(f)$ to $\H$ such that the stable orbital integrals of $f^H$
are not identically zero on the set of regular semi-simple elliptic elements
of $\H(\C)$.
\item[(ii)] For every $f\in E$ and every transfer $f^H$ of $N(f)$,
$SO_{\gamma_H} (f^H)=0$ if $\gamma_H$ is not conjugate to an element of $C$.
\item[(iii)] Let $(a(\pi))_{\pi\in\Pi_\theta(R)}$ and $(b(\pi_H))_{\pi_H\in
\Pi(\H)}$ be families of complex numbers such that, for every $f\in E$ and
every transfer $f^H\in C_c^\infty(\H(\C),\K_H)$ of $N(f)$ to $\H$,
the sums $A(f):=\sum\limits_{\pi\in\Pi(R)}a(\pi)\Tr(\pi(f)A_\pi)$ and
$A_H(f^H):=\sum\limits_{\pi_H\in\Pi(\H)}b(\pi_H)\Tr(\pi_H(f^H))$ are
absolutely convergent.
Then the following conditions are equivalent :
\begin{itemize}
\item[(A)]  for every $f\in E$ and every transfer $f^H\in C_c^\infty(\H
(\C),\K_H)$ of $N(f)$ to $\H$, $A(f)=A_H(f^H)$;
\item[(B)] for every $f\in E$, for every transfer $f^H\in C_c^\infty(\H
(\C),\K_H)$ of $N(f)$ to $\H$ and for every $\pi\in\Pi_{\theta-temp}(R)$,
$a(\pi)\Tr(\pi(f)A_\pi)=\sum\limits_{\pi_H\in\Pi_H(\pi)}b(\pi_H)\Tr(\pi_H
(f^H))$.

\end{itemize}
\end{itemize}
\end{sublemme}

The next lemma is proved in \cite{Cl-LF}, p 292.

\begin{sublemme}\label{lemme:presque_tout_gamma} Notations are as in
\ref{lemme_fondamental1}. Let $f\in\Hecke_R$. If there exists a dense subset
$D$ of the set of elliptic semi-simple strongly $\G$-regular elements of
$\H(F)$ such that $\Lambda(\gamma_H,f)=0$ for  every $\gamma_H\in D$,
then $\Lambda(\gamma_H,f)=0$ for every elliptic semi-simple strongly
$\G$-regular $\gamma_H\in\H(F)$.

\end{sublemme}

\begin{sublemme}\label{lemme:isomorphisme_donnees_endoscopiques} Let $F$
be a number field, $E$ be a cyclic extension of $F$ and $\G$ be a connected
reductive group over $F$. Set $R=R_{E/F}\G_E$, choose a generator of
$\Gal(E/F)$ and denote by $\theta$ the automorphism of $R$ induced by this
generator. Let $K$ be a finite extension of $E$ such that $\G$ splits over
$K$. Assume that the center of $\G$ is connected and that there exists a
finite place $v$ of $F$, inert in $K$, such that the morphism
$\Gal(K_v/F_v)\fl\Gal(K/F)$ is an isomorphism.
Then localization induces an injective map from the set of equivalence classes
of endoscopic data for $(R,\theta,1)$ to the set of equivalence classes of
endoscopic data for $(R_v,\theta_v,1)$.

% Let $\T$ be a maximal torus of $\G$. Then $\T$ is elliptic (over $F$)
% if $\T_v:=\T_{F_v}$ is elliptic (over $F_v$). Moreover, the localization
% map $\Kgoth(\T,R;F)\fl\Kgoth(\T_v,R_v;F_v)$ of \cite{La-CSCB} p 43
% is injective.

\end{sublemme}

\begin{proof}
Let $(\H,\Hcal,s,\xi),(\H',\Hcal',s',\xi')$ be endoscopic data for
$(R,\theta,1)$ whose localizations are equivalent (as endoscopic data for
$(R_v,\theta_v,1)$). 
By lemma \ref{lemme:E_egal_F_puissance_d}, the endoscopic data for
$(R_K,\theta_K,1)$ defined by $(\H,\Hcal,s,\xi)$ and $(\H',\Hcal',s',\xi')$
are equivalent to endoscopic
data coming from endoscopic data for $(\G_K,1,1)$. As the derived group of
$\widehat{\G}$ is simply connected (because the center of $\G$ is connected)
and $\G_K$ is split, if $(\G',\Gcal',s_G,\xi_G)$ are endoscopic data for
$(\G_K,1,1)$, then $\G'$ is split, so $\Gcal'\simeq\widehat{\G'}\times W_K$,
and we may assume that $\xi_G$ is the product of an embedding
$\widehat{\G'}\fl\widehat{R}$ and of the identity on $W_K$.
Hence, after replacing $(\H,\Hcal,s,\xi)$ and $(\H',\Hcal',s',\xi')$ by
equivalent data, we may assume that $\xi$ and $\xi'$ come from $L$-morphisms
$\widehat{\H}\rtimes\Gal(K/F)\fl\widehat{R}\rtimes\Gal(K/F)$ and
$\widehat{\H}'\rtimes\Gal(K/F)\fl\widehat{R}\rtimes\Gal(K/F)$.
As the data $(\H,\Hcal,s,\xi)$ and $(\H',\Hcal',s',\xi')$ are equivalent at
$v$, we may identify $\widehat{\H}$ and $\widehat{\H}'$ and assume that
$s=s'$. As $\Gal(K_v/F_v)\iso\Gal(K/F)$ (and $\Gal(\overline{F}/K)$ acts
trivially on $\widehat{\H}$ and $\widehat{\H}'$), the isomorphism
$\widehat{\H}=\widehat{\H}'$ extends to an isomorphism $\Hcal\simeq\Hcal'$
that identifies $\xi$ and $\xi'$. So the data $(\H,\Hcal,s,\xi)$ and
$(\H',\Hcal',s',\xi')$ are equivalent, and the first statement of the lemma
is proved.

% Let $\T$ be a maximal torus of $\G$. It is obvious that $\T$ is elliptic
% if $\T_v$ is.
% We want to show the last statement of the
% lemma. It is equivalent to prove that the map
% $\Egoth(\T,R;F_v)\fl\Egoth(\T,R;\Ade_F/F)$ is surjective.
% Reason as in the proof of lemma 1.9.7 of \cite{La-CSCB}.
% Imitating \cite{Cl-LF} p 293, we see that $\Ker^1_{ab}(F,R)=\{1\}$.
% By proposition 1.8.4 of \cite{La-CSCB}, there is a commutative diagram
% with exact rows :
% \[\xymatrix{1\ar[r] & \Egoth(\T,R;F_v)\ar[r]\ar[d] & \Ho^1_{ab}(F_v,\T)\ar[r]
% \ar[d] & \Ho^1_{ab}(F_v,R)\ar[d] \\
% 1\ar[r] & \Egoth(\T,R;\Ade_F/F)\ar[r] & \Ho^1_{ab}(\Ade_F/F,\T)\ar[r] &
% \Ho^1_{ab}(\Ade/F,R)}\]
% By proposition 1.7.3 of \cite{La-CSCB}, there are canonical isomorphisms
% $\Ho^1_{ab}(F_v,\T)\simeq\pi_0(\widehat{\T}^{\Gamma_v})^D$,
% $\Ho^1_{ab}(\Ade_F/F,\T)\simeq\pi_0(\widehat{\T}^\Gamma)^D$, $\Ho^1_{ab}
% (F_v,R)\simeq\pi_0(Z(\widehat{R})^{\Gamma_v})^D$ and $\Ho^1_{ab}(\Ade_F/F,R)
% \simeq\pi_0(Z(\widehat{R})^\Gamma)^D$, where $\Gamma=\Gal(\overline{F}/F)$
% and $\Gamma_v=\Gal(\overline{F}_v/F_v)$. So the middle vertical map is
% surjective (???). As $R$ splits over $K$ and
% $\Gal(K_v/F_v)\iso\Gal(K/F)$, the vertical map on the right
% is an isomorphism, and this implies the desired result.

\end{proof}

\begin{proofpn}{\ref{prop:existence_donnees_locales}} Write $S_0=\{v_0,v_1,
v_2\}$. Identify $k_{F,v_0}$, $k_{E,v_0}$, etc, to $F$, $E$, etc.
We will prove the proposition by applying the twisted trace formula on
$\underline{R}$ to functions whose local component at $v_0$ is a function in
$\Hecke_R$. 

Let $\T_G$ be an elliptic maximal torus of $\underline{\G}$ such that
$\T_{G,k_{F,v_1}}$ is also elliptic. Fix an
admissible embedding $\T_H\fl\underline{\G}$ with image $\T_G$, where $\T_H$
is an elliptic maximal torus of $\underline{\H}$, and let
$\T_R=R_{k_E/k_F}\T_{G,k_E}$. Let $\kappa$ be the
element of $\Kgoth(\T_R,\theta,k_F)_1=\Kgoth(\T_G,\underline{R};k_F)_1$ (cf
remark \ref{rq:comparaison_Kgoth}) associated to the endoscopic data
$(\underline{\H},{}^L\underline{\H},t,\underline{\xi})$ by the map of
\cite{KS} 7.2. Write $\kappa_{v_1}$ for the image of $\kappa$ by the
localization map $\Kgoth(\T_H,\underline{R};k_F)\fl\Kgoth(\T_{H,v_1},
\underline{R}_{v_1};k_{F,v_1})$ (cf \cite{La-CSCB} p 43). Choose a function
$\phi_{v_1}\in C_c^\infty(\underline{R}(k_{F,v_1}))$ that satisfies the
conditions of lemma \ref{lemme:caracterisation_kappa_integrales_orbitales}
(ie such that the support of $\phi_{v_1}$ is contained in the union of the
stable $\theta$-conjugates of $\T_R(k_{F,v_1})$, that the 
$\kappa_{v_1}$-orbital integrals of $\phi_{v_1}$ are not all zero and
that the $\kappa'_{v_1}$-orbital integrals of $\phi_{v_1}$ are all zero
if $\kappa'_{v_1}\not=\kappa_{v_1}$). Let $f^H_{v_1}$ be a transfer of
$\phi_{v_1}$ to $\underline{\H}_{v_1}$. Let
$f^H_{v_2}\in C_c^\infty(\underline{\H}(k_{F,v_2}))$ be a function with
support in the set of semi-simple strongly $\G$-regular elements, whose
orbital integrals are constant on stable conjugacy classes and whose
stable orbital integrals are not all zero (such a function exists by
lemma \ref{lemme:caracterisation_kappa_integrales_orbitales}, applied
with $\theta=1$ and $\kappa=1$). Fix a function $\phi_{v_2}\in C_c^\infty
(\underline{R}(k_{F,v_2}))$ such that $f^H_{v_2}$ is a transfer of $\phi_{v_2}$
(such a function exists by lemma \ref{lemme:construction_transfert}).

Let $v_3$ and $v_4$ be finite places of $k_F$ where all the data are
unramified (ie where the situation is as in remark
\ref{rq:E_pas_un_corps}); assume moreover that $v_3$ splits totally in
$k_E$ (this is possible by (vi) of lemma \ref{lemme:situation_globale}).
Let $f_{v_3}\in C_c^\infty(\underline{\G}(k_{F,v_3}))$ be as in lemma
\ref{lemme:coeff_de_supercuspidales}. Write
$\phi_{v_3}=f_{v_3}\otimes\dots\otimes f_{v_3}\in C_c^\infty(
\underline{R}(k_{F,v_3}))$ (where we identified $\underline{R}(k_{F,v_3})$ to
$\underline{\G}(k_{F,v_3})^d$), and choose a transfer $f^H_{v_3}$ of
$\phi_{v_3}$ (such a transfer exists by lemma
\ref{lemme:construction_transfert}). Let $f^H_{v_4}\in C_c^\infty(\underline
{\H}(k_{F,v_4}))$ be as in lemma \ref{lemme:coeff_de_supercuspidales}.
Choose a function $\phi_{v_4}\in C_c^\infty(\underline{R}(k_{F,v_4})$
such that $f^H_{v_4}$ is a transfer of $\phi_{v_4}$ (such a function exists
by lemma \ref{lemme:construction_transfert}).

Let $S_\infty$ be the set of infinite places of $k_F$ (by (iii) of lemma
\ref{lemme:situation_globale}, they are all complex). Write
$\H_\infty=\prod\limits_{v\in S_\infty}\underline{\H}(k_{F,v})$
and $R_\infty=\prod\limits_{v\in S_\infty}\underline{R}(k_{F,v})$.
Let $E$ be a subspace of $C_c^\infty(R_\infty)$ and $C_\infty$
be a compact subset of $\H_\infty$ that satisfy the conditions of lemma
\ref{lemme:place_infinie}. Let $\phi_{0,\infty}\in E$ and
$f^H_{0,\infty}$ be a transfer of $\phi_{0,\infty}$ such that the stable
orbital integrals of $f^H_{0,\infty}$ on elliptic elements of $\H_\infty$
are not all zero.

Let $D_1$ be the set of semi-simple strongly $\G$-regular elliptic elements
of $\H(F)$ coming from a $\gamma_H\in\underline{\H}(k_F)$ such that
\begin{itemize}
\item[-] there exists $\delta\in\underline{R}(k_F)$ and an image $\gamma$ of
$\gamma_H$ in $\underline{\G}(k_F)$ such that $\gamma\in\Norme\delta$;
\item[-] for every $v\in\{v_1,v_2,v_3,v_4\}$, $SO_{\gamma_H}(f^H_v)\not=0$;
\item[-] $SO_{\gamma_H}(f^H_{0,\infty})\not=0$.

\end{itemize}
Let $D_2$ be the set of semi-simple strongly $\G$-regular elliptic elements of
$\H(F)$ that have no image in $\G(F)$ that is a norm.
By (vii) of lemma \ref{lemme:situation_globale}, $D:=D_1\cup D_2$ is dense in
the set of semi-simple strongly $\G$-regular elliptic elements of $\H(F)$.
By lemma \ref{lemme:presque_tout_gamma}, we may replace the set of semi-simple
strongly $\G$-regular elliptic elements of $\H(F)$ by $D$ in the definition
of local data. By lemma \ref{lemme:resultat_annulation}, we may even replace
$D$ by $D_1$ in this definition.

Let $\gamma_H\in D_1$ (we use the same notation for the element of $\H(F)$ and
for the element of $\underline{\H}(k_F)$ that induces it).
Let $S$ be a finite set of finite places of $k_F$ such that $\{v_0,v_1,v_2,
v_3,v_4\}\subset S$ and that, for every finite place $v\not\in S$ of $k_F$,
all the data are unramified at $v$, $\gamma_H\in\underline{\H}(
\Of_{k_{F,v}})$ and the fundamental lemma for the unit of the Hecke algebra
is known for $(\underline{R}_v,\underline{\theta}_v,1)$ and
$(\underline{\H}_v,{}^L\underline{\H}_v,t,\underline{\xi}_v)$. For every
$v\in S-\{v_0,v_1,v_2,v_3,v_4\}$, choose associated functions
$f^H_v$ and $\phi_v$ such that $SO_{\gamma_H}(f^H_v)\not=0$ (this is possible
by the end of lemma \ref{lemme:caracterisation_kappa_integrales_orbitales}). 
Let $C_0\ni\gamma_H$ be a compact subset of $\H(F)$ that meets all
the conjugacy classes of semi-simple elliptic elements of $\H(F)$
(such a $C_0$ exists because the center of $\H$ is anisotropic).
By proposition 8.2 of \cite{K-STF:EST}, there is only a finite number of
conjugacy classes of
semi-simple elements $\gamma'_H$ of $\underline{\H}(k_F)$ such that
$\gamma'_H\in C_0$, $\gamma'_H\in C_\infty$, $SO_{\gamma'_H}(f^H_v)\not=0$ for
every $v\in S-\{v_0\}$ and $\gamma'_H\in\underline{\H}(\Of_{k_{F,v}})$ for
every finite place $v\not\in S$. By the end of lemma
\ref{lemme:caracterisation_kappa_integrales_orbitales}, after adding a place
in $S$ and fixing well-chosen functions at that place, we may assume
that $\gamma_H$ is, up to conjugacy, the only semi-simple element of
$\underline{\H}(k_F)$ that satisfies the list of properties given above.
For every finite place $v\not\in S$ of $k_F$, take $\phi_v=\ungras_{
\underline{R}(\Of_{k_{F,v}})}$ and $f^H_v=\ungras_{\underline{\H}(\Of_{k_{F,v}}
)}$. 

Let $\phi_{v_0}\in\Hecke_R$ and $f^H_{v_0}=b_\xi(\phi_{v_0})$. 
Fix $\phi_\infty\in E$ and a transfer $f^H_\infty$ of $\phi_\infty$,
and set $\phi=\phi_\infty\otimes\bigotimes\limits_{v\not=\infty}\phi_v$ and
$f^H=f^H_\infty\otimes\bigotimes\limits_{v\not=\infty}f^H_v$. Then lemma
\ref{lemme:formule_traces_simple} applies to $f$ and $\phi$, thanks to the
choice of the functions at $v_3$ and $v_4$. 
As in \ref{FT_stable_geometrique4} and \ref{GL_n_applications2}, 
let $T^{\underline{R}\rtimes\underline{\theta}}$ and $T^{\underline{\H}}$
be the distributions of the $\underline{\theta}$-twisted invariant trace
formula on $\underline{R}$ and of the invariant trace formula on
$\underline{\H}$. 

By lemma \ref{lemme:formule_traces_simple},
$T^{\underline{\H}}(f)$ is equal to the strongly regular elliptic part
of the trace formula for $\underline{\H}$, so we may use the stabilization of
\cite{L3}. By the choice of $f_{v_2}^H$, the only endoscopic group of
$\underline{\H}$ that appears is $\underline{\H}$ itself; so we need
neither the transfer hypothesis nor the fundamental lemma to stabilize
$T^{\underline{\H}}(f)$.
We get
\[T^{\underline{\H}}(f^H)=ST^{\underline{\H}}_{**}(f^H),\]
where $ST^{\underline{\H}}_{**}$ is the distribution denoted by $ST_e^{**}$
in \cite{KS} 7.4 (the strongly $\underline{\G}$-regular elliptic part of
the stable trace formula for $\underline{\H}$). Moreover, by the choice of
$f^H$,
\[ST^{\underline{\H}}_{**}(f^H)=a_{\phi_\infty}SO_{\gamma_H}
(b_\xi(\phi_{v_0})),\]
where $a_{\phi_\infty}$ is the product of $SO_{\gamma_H}(f_\infty)$ and of
a non-zero scalar that does not depend on $\phi_\infty$ and $\phi_{v_0}$.

Similarly, by lemma \ref{lemme:formule_traces_simple},
$T^{\underline{R}\rtimes\underline{\theta}}(\phi)$ is equal to the
strongly $\underline{\theta}$-regular $\underline{\theta}$-elliptic part
of the trace formula for $\underline{R}\rtimes\underline{\theta}$,
se we may apply the stabilization of chapters 6 and 7 of \cite{KS}.
By the choice of $\phi_{v_1}$, the only endoscopic data of
$(\underline{R},\underline{\theta},1)$ that appear are
$(\underline{\H},\Hcal,t,\underline{\xi})$.
(Equality (7.4.1) of \cite{KS} expresses $T^{\underline{R}\rtimes\underline{
\theta}}(\phi)$ as a sum over elliptic endoscopic data of
$(\underline{R},\underline{\theta},1)$. The proof of lemma 7.3.C and
theorem 5.1.D of \cite{KS} show that the global $\kappa$-orbital integrals
that appear in this sum are products of local $\kappa$-orbital integrals.
By the choice of $\phi_{v_1}$, these products of local $\kappa$-orbital
integrals are zero for the endoscopic data that are not equivalent to
$(\underline{\H},\Hcal,t,\underline{\xi})$ at the place $v_1$. By lemma
\ref{lemme:isomorphisme_donnees_endoscopiques}, $(\underline{\H},\Hcal,t,
\underline{\xi})$ are the only endoscopic data satisfying this condition.)
By equality (7.4.1) and the proof of lemma 7.3.C of \cite{KS}, we get
\[T^{\underline{R}\rtimes\underline{\theta}}(\phi)=b_{\phi_\infty}\sum_\delta
\Delta_\xi(\gamma_H,\delta)O_{\delta\underline{\theta}}(\phi_{v_0}),\]
where the sum is taken over the set of $\theta$-conjugacy classes of
$R(F)$ and $b_{\phi_\infty}$ is the product of $SO_{\gamma_H}(f^H_\infty)$
and of a non-zero scalar that does not depend on $\phi_{v_0}$ and
$\phi_\infty$. 

Hence $\Lambda(\gamma_H,\phi_{v_0})=0$ if and only if, for every
$\phi_\infty\in E$, $a_{\phi_\infty}T^{\underline{R}\rtimes\underline{\theta}}
(\phi)-b_{\phi_\infty}T^{\underline{\H}}(f^H)=0$ (the ``only if'' part comes
from the fact that $a_{\phi_\infty}b_{\phi_\infty}\not=0$ for at least one
choice of $\phi_\infty$). By lemmas
\ref{lemme:formule_traces_simple} and \ref{lemme:place_infinie},
this last condition is equivalent to a family of identities of the form
\begin{flushleft}$\displaystyle{\Tr(\pi_\infty(\phi_\infty)A_{\pi_\infty})
\sum_{\pi_0\in\Pi(R)}a(\pi_0)\Tr(\pi_0(\phi_{v_0}A_{\pi_0}))}$
\end{flushleft}
\begin{flushright}$\displaystyle{=\sum_{\pi_{H,\infty}\in\Pi_{H_\infty}
(\pi_\infty)}
\Tr(\pi_{H,\infty})(f^H_\infty)\sum_{\pi_{H,0}\in\Pi(\H)}b(\pi_{H,\infty},
\pi_{H,0})\Tr(\pi_{H,0}(b_\xi(\phi_{v_0}))),}$\end{flushright}
for $\pi_\infty\in\Pi_{\theta-temp}(R_\infty)$ and $\phi_\infty\in E$,
where the notations are as in lemma \ref{lemme:place_infinie}.
By Harish-Chandra's finiteness theorem (cf \cite{BJ} 4.3(i)), the sums that
appear in these equalities have only a finite number of non-zero terms.

Finally, we showed that the identity $\Lambda(\gamma_H,\phi_{v_0})=0$
is equivalent to a family of identities like those that appear in the
definition of local data.
To obtain local data for $R$ and $(\H,{}^L\H,t,\xi)$, we simply have to
repeat this process for all the elements of $D_1$.

\end{proofpn}

\section{Technical lemmas}
\label{lemme_fondamental4}

We use the notations of \ref{lemme_fondamental1}.

Let $\gamma_H\in\H(F)$ be semi-simple and
$\gamma\in\G(F)$ be an image of $\gamma_H$.
Let $\M_H$ be a Levi subgroup of $\H$ such that $\gamma_H\in\M_H(F)$ and
$\M_{H,\gamma_H}=\H_{\gamma_H}$. Langlands and Shelstad (\cite{LS2} \S1,
see also section 7 of \cite{K-NP}) 
associated to such a $\M_H$ a Levi subgroup $\M$ of $\G$ such that
$\gamma\in\M(F)$ and $\M_\gamma=\G_\gamma$, an endoscopic triple
$(\M_H,s_M,\eta_{M,0})$ for $\M$ and a $L$-morphism
$\eta_M:{}^L\M_H\fl{}^L\M$ extending $\eta_{M,0}$ and such that there is a
commutative diagram
\[\xymatrix{{}^L\M_H\ar[r]^-{\eta_M}\ar[d] & {}^L\M\ar[d] \\
{}^L\H\ar[r]^-\eta & {}^L\G}\]
where the left (resp. right) vertical map is in the
canonical $\widehat{\H}$-conjugacy (resp. $\widehat{\G}$-conjugacy) class
of $L$-morphisms ${}^L\M_H\fl{}^L\H$ (resp. ${}^L\M\fl{}^L\G$).
(If $\G$ and $\H$ are as in \ref{groupes3}, this construction is made more
explicit in \ref{partie_en_p3}.)
Write $\M_R=R_{E/F}\M$ (a $\theta$-stable Levi subgroup of $R$), and let
$\theta_M$ be the restriction of $\theta$ to $\M_R$.
As in \ref{lemme_fondamental1}, associate to $(\M_H,s_M,\eta_M)$ endoscopic
data $(\M_H,{}^L\M_H,t_M,\xi_M)$ for $(\M_R,\theta_M,1)$.

The next lemma is a generalization of the beginning of \cite{H} 12, and can
be proved exactly in the same way (because there is a descent formula for
twisted orbital integrals, cf for example corollary 8.3 of \cite{A-ITF1}).
Note that we also need to use lemma 4.2.1 of \cite{Ha} (and the remarks
below it).

\begin{sublemme}\label{lemme:descente} Assume that, for every proper Levi
subgroup $\M_H$ of $\H$, the twisted fundamental lemma is known for
$\M_R$ and $(\M_H,{}^L\M_H,t_M,\xi_M)$ and for all the functions of
$\Hecke_{M_R}$. Then, for every $f\in\Hecke_R$ and every $\gamma_H\in\H(F)$
that is semi-simple and is not elliptic,
\[\Lambda(\gamma_H,f)=0.\]

\end{sublemme}

\begin{sublemme}\label{lemme:car_sur_classes_theta_stables} Let $\chi$
be a character of $R(F)$ such that $\chi=\chi\circ\theta$.
Then $\chi$ is contant on the $\theta$-semi-simple stable $\theta$-conjugacy
classes of $R(F)$.

\end{sublemme}

\begin{proof} If $\theta=1$ (ie if $E=F$), this is lemma 3.2 of
\cite{Ha}. In the general case, the result follows from the case
$\theta=1$ and from lemma 2.4.3 of \cite{La-CSCB}.

\end{proof}

Notations are still as in \ref{lemme_fondamental1}. Let
$\B_G$ be a Borel subgroup of $\G$ and $\T_G$ be a Levi subgroup of
$\B_G$. Assume that the center of $\G$ is connected and that there exists
a maximal torus $\T_H$ of $\H$ and an admissible embedding $\T_H\fl\G$
with image $\T_G$ (if $\G$ is adjoint, this is always the case by lemma
\ref{lemme:groupes_qui_marchent}). Use the same notations as in lemma
\ref{lemme:DL8} (in particular, $\T_R=R_{E/F}\T_{G,E}$).

Let $Z(R)_\theta$ be the group of $\theta$-coinvariants of the center
$Z(R)$ of $R$. There is a canonical injective morphism $N:Z(R)_\theta\fl
Z(\H)$ (cf \cite{KS} 5.1, p 53). Let $Z$ be a subtorus of $Z(R)_\theta$; denote
by $Z_R$ the inverse image of $Z$ in $Z(R)$ and by $Z_H$ the image of $Z$ in
$Z(\H)$. Let $\chi_H$ be an unramified character of $Z_H(F)$. Write
$\chi_R=(\lambda_C^{-1}(\chi_H\circ N))_{|Z_R(F)}$, where $\lambda_C:Z(R)(F)
\fl\C^\times$ is the character defined in \cite{KS} 5.1 p 53; then $\chi_R$
is also unramified, by (i) of lemma \ref{lemme:H_equivariante} below.

Let $\Hecke_{R,\chi_R}$ (resp. $\Hecke_{\H,\chi_H}$) be the algebra of
functions $f:R(F)\fl\C$ (resp. $f:\H(F)\fl\C$) that are right and left
invariant by $R(\Of_F)$ (resp. $\H(\Of_F)$), have compact support modulo
$Z_R(F)$ (resp. $Z_H(F)$), and such that, for every
$(z,x)\in Z_R(F)\times R(F)$ (resp. $Z_H(F)\times\H(F)$),
$f(zx)=\chi_R^{-1}(z)f(x)$ (resp. $f(zx)=\chi_H^{-1}(z)f(x)$). 
The product is the convolution product, that sends $(f,g)$ to
\[f*g:x\fle\int_{Z_R(F)\sous R(F)}f(xy^{-1})g(y)dy\]
\[(\mbox{resp.}\quad f*g:x\fle\int_{Z_H(F)\sous H(F)}f(xy^{-1})g(y)dy).\]
There is a surjective morphism $\nu_R:\Hecke_R\fl\Hecke_{R,\chi_R}$ (resp.
$\nu_H:\Hecke_H\fl\Hecke_{H,\chi_H}$) that sends $f$ to $x\fle\int_{Z_R(F)}
\chi_R^{-1}(z)f(zx)dz$ (resp. $x\fle\int_{Z_R(F)}\chi_H^{-1}(z)f(zx)dz$).
If $\gamma_H\in\H(F)$ and $f\in\Hecke_{H,\chi_H}$, then we can define
$O_{\gamma_H}(f)$ by the usual formula (the integral converges). If
$\delta\in R(F)$ and $f\in\Hecke_{R,\chi_R}$, set
\[O_{\delta\theta}(f)=\int_{Z_R(F)R_{\delta\theta}(F)\sous R(F)}f(x^{-1}\delta
\theta(x))dx\]
($f(x^{-1}\delta\theta(x))$ depends only on the class of $x$ in
$Z_R(F)R_{\delta\theta}(F)\sous R(F)$, because $\chi_R$ is trivial on
elements of the form $z^{-1}\theta(z)$, $z\in Z_R(F)$). We define
$\kappa$-orbital integrals and stable orbital integrals as before for
functions of $\Hecke_{R,\chi_R}$ and $\Hecke_{H,\chi_H}$.

\begin{sublemme}\label{lemme:H_equivariante}\begin{itemize}
\item[(i)] As in the proof of lemma \ref{lemme:DL8}, write
$\xi(1\rtimes\sigma)=t'\rtimes\sigma$; we may assume that
$t'\in\widehat{\T}_R$. As $Z(R)$ is an unramified torus, there is a canonical
surjective morphism $Z(R)(F)\fl X_*(Z(R)_d)$ with kernel the maximal compact
subgroup of $Z(R)(F)$, where $Z(R)_d$ is the maximal split subtorus of
$Z(R)$ (cf \cite{Bo} 9.5). Define a character $\lambda'_C$ on $Z(R)(F)$
in the following way : if $z\in Z(R)(F)$, and if
$\mu\in X_*(Z(R)_d)=X^*(\widehat{Z(R)}_d)$ is the image of $z$ by the above
morphism, $\lambda'_C(z)$ is the value of $\mu$ at the image of
$t'^{-1}\in\widehat{\T}_R$ by the canonical morphism
$\widehat{\T}_R\fl\widehat{Z(R)}\fl\widehat{Z(R)}_d$.\newline
Then $\lambda_C=\lambda'_C$ (in particular, $\lambda_C$ is unramified).
\item[(ii)] There exists a morphism $b_{\xi,\chi_H}:\Hecke_{R,\chi_R}\fl
\Hecke_{H,\chi_H}$ that makes the following diagram commute
\[\xymatrix{\Hecke_R\ar[r]^-{\nu_R}\ar[d]_{b_\xi} & \Hecke_{R,\chi_R}\ar[d]^
{b_{\xi,\chi_H}} \\
\Hecke_H\ar[r]_-{\nu_H} & \Hecke_{H,\chi_H}}\]
Let $\gamma_H\in\H(F)$ be semi-simple and strongly $\G$-regular. 
Use the morphism $b_{\xi,\chi_H}:\Hecke_
{R,\chi_R}\fl\Hecke_{H,\chi_H}$ to define a linear form
$\Lambda_{\chi_H}(\gamma_H,.)$ on $\Hecke_{R,\chi_R}$ that is the analog
of the linear form $\Lambda(\gamma_H,.)$
on $\Hecke_R$ of \ref{lemme_fondamental1}
(use the same formula). Then the following conditions are equivalent :
\begin{itemize}
\item[(a)] for every $z\in Z_H(F)$, for every $f\in\Hecke_R$,
$\Lambda(z\gamma_H,f)=0$;
\item[(b)] for every $z\in Z_H(F)$, for every $f\in\Hecke_{R,\chi_R}$,
$\Lambda_{\chi_H}(z\gamma_H,f)=0$.

\end{itemize}
\end{itemize}
\end{sublemme}

\begin{proof} Point (i) follows from the definitions of $\lambda_C$ (\cite{KS}
5.1) and of the transfer factor $\Delta_{III}$ (\cite{KS} 4.4).

We show (ii). Let $z\in Z(R)(F)$. For every function $f:R(F)\fl\C$, let
$R_zf$ be the function $x\fle f(zx)$. Denote by $\lambda_z$ the image of $z$
by the canonical map $Z(R)(F)\subset\T_R(F)\fl
X_*(Y_R)$. Then, for every $\lambda\in X^*(Y_R)$, $R_z f_\lambda=
f_{\lambda+\lambda_z}$. Moreover, it is easy to see that, for every
$\lambda\in X^*(Y_R)$, $\nu_R^{-1}(\nu_R(f_\lambda))$ is generated by the
functions $\chi_R(z)^{-1}R_zf_\lambda$, $z\in Z_R(F)$. There are obviously
similar statements for $\H$ instead of $R$.

To show the existence of the morphism $b_{\xi,\chi_H}:\Hecke_{R,\chi_R}\fl
\Hecke_{H,\chi_H}$, it is enough to show that, for every
$\lambda\in X^*(Y_R)$, all the elements of $\nu_R^{-1}(\nu_R(f_\lambda))$
have the same image by $\nu_H\circ b_\xi$. Let $\lambda\in X^*(Y_R)$. Let
$z\in Z_R(F)$; denote by $z_H$ the image of $z$ in $Z_H(F)$. It is enough
to show that
$b_\xi(\chi_R^{-1}(z)R_zf_\lambda)=\chi_H(z_H)^{-1}R_{z_H}b_\xi(f_\lambda)$.
By the explicit calculation of $b_\xi f_\lambda$ in the proof of lemma
\ref{lemme:DL8}, $b_\xi(R_zf_\lambda)=\lambda_z(t')R_{z_H}b_\xi(f_
\lambda)$; hence the equality that we are trying to prove follows from (i)
and from the definition of $\chi_R$.

Let $f\in\Hecke_R$ and $\delta\in R(F)$. It is easy to see that
\[O_{\delta\theta}(\nu_R(f))=\int_{Z(F)}\chi_R^{-1}(z)O_{z\delta
\theta}(f)dz=\int_{Z(F)}\chi_R^{-1}(z)O_{\delta\theta}(R_zf)dz\]
($O_{z\delta\theta}(f)$ depends only on the image of $z$ in $Z(F)$,
because the function $\delta\fle O_{\delta\theta}(f)$ is invariant by
$\theta$-conjugacy). Similarly, for every $f\in\Hecke_H$ and
$\gamma_H\in\H(F)$,
\[O_{\gamma_H}(\nu_H(f))=\int_{Z_H(F)}\chi_H^{-1}(z_H)O_{z_H\gamma_H}(f)dz_H=
\int_{Z_H(F)}\chi_H^{-1}(z_H)O_{\gamma_H}(R_{z_H}f)dz_H.\]
Remember (\cite{KS} 5.1) that $\lambda_C$ is such that, for every semi-simple
strongly regular $\gamma_H\in\H(F)$, every $\theta$-semi-simple strongly
$\theta$-regular $\delta\in R(F)$ and every $z\in Z(R)(F)$,
\[\Delta_\xi(z_H\gamma_H,z\delta)=\lambda_C^{-1}(z)\Delta_\xi(\gamma_H,
\delta),\]
where $z_H$ is the image of $z$ in $Z(\H)(F)$. By this fact and the above
formulas for the integral orbitals, it is clear that (a) implies (b).

Let $\gamma_H\in\H(F)$ be semi-simple strongly $\G$-regular.
Assume that (b) is satisfied for $\gamma_H$; we want to show (a).
Let $\lambda\in X^*(Y_R)$.
Denote by ${}^0Z$ (resp. ${}^0Z_R$, ${}^0Z_H$) the maximal compact subgroup of
$Z(F)$ (resp. $Z_R(F)$, $Z_H(F)$). The function $f_\lambda$ is obviously
invariant by translation by ${}^0Z_R$, and, moreover, all the unramified
$\theta$-stable characters of $R(F)$ are constant on its support
(because $f_\lambda$ is a linear combination of characteristic functions
of sets $R(\Of_F)\mu(\varpi_F)R(\Of_F)$, where $\varpi_F$ is a uniformizer
of $F$ and $\mu\in X^*(Y_R)=X_*(\T_R)^{\Gamma_F}$ is such that
$\lambda\mu^{-1}$ is a cocharacter of $R_{der}$). Hence, for
$z\in Z_R(F)$, $R_zf_\lambda$ depends only on the image of $z$ in
${}^0Z_R\sous Z_R(F)$ and, for every $\theta$-semi-simple stable
$\theta$-conjugacy class $C$ of $R(F)$, there exists a unique $z\in {}^0Z_R
\sous Z_R(F)$ such that, for every $z'\in{}^0Z_R\sous Z_R(F)-\{z\}$ and every
$\delta\in C$, $O_{\delta\theta}(R_{z'}f_\lambda)=0$ (use lemma
\ref{lemme:car_sur_classes_theta_stables}). There are similar results for
$\H$ and $b_\xi(f_\lambda)$. 

Let $C$ be the set of $\delta\in R(F)$ such that $\Delta_\xi(\gamma_H,
\delta)\not=0$. Then $C$ is either the empty set or a $\theta$-semi-simple
$\theta$-regular stable $\theta$-conjugacy class.
So, by the reasoning above and the formulas for the orbital integrals of
$\nu_R(f_\lambda)$ and $\nu_H(b_\xi(f_\lambda))$, there exists
$z\in {}^0Z\sous Z(F)$ such that $O_{\delta\theta}(\nu_R(f_\lambda))=\chi_R(z)
^{-1}O_{\delta\theta}(R_zf_\lambda)$, for every $\delta\in C$, and that
$SO_{\gamma_H}(\nu_H(b_\xi(f_\lambda)))=\chi_H(z_H)^{-1}SO_{\gamma_H}(
R_{z_H}b_\xi(f_\lambda))$, where $z_H\in {}^0Z_H\sous Z_H(F)$ is the image of
$z$.
If $z\not=1$, then $\Lambda(\gamma_H,f_\lambda)=0$, because all the orbital
integrals that appear in this expression are zero. If $z=1$, then
$\Lambda(\gamma_H,f_\lambda)=0$ by condition (b) and the properties of
$\lambda_C$.

\end{proof}

We still denote by $F$ a non-archimedean local field (of characteristic $0$)
and by $E$ a finite unramified extension of $F$.
Let $\G$ be a connected unramified group over $F$, defined over $\Of_F$
and such that $\G(\Of_F)$ is a hyperspecial maximal compact subgroup of
$\G(F)$. Set $R=R_{E/F}\G_E$, and let $\theta$ be the automorphism of $R$
induced by a chosen generator of $\Gal(E/F)$. Let $Z_G$ be a subtorus of
$Z(\G)$ defined over $\Of_F$, and let $Z_R=R_{E/F}Z_{G,E}$.
Let $\G'=\G/Z_G$, $R'=R/Z_R=R_{E/F}\G'_E$, $u:R\fl R'$ be the obvious morphism,
$\Hecke'=\Hecke_{R'}$ and $\Hecke$ be the convolution algebra of
functions $R(F)\fl\C$ that are bi-invariant by $R(\Of_F)$, invariant by
$Z_R(F)$ and with compact support modulo $Z_R(F)$ (with the notations of
lemma \ref{lemme:H_equivariante}, $\Hecke=\Hecke_{R,1}$). As $Z_R$ is
connected, we see as in \cite{Cl-LF} 6.1 (p 284) that Lang's theorem (cf for
example theorem 4.4.17 of \cite{Sp}) and Hensel's lemma imply that
$u:R(\Of_F)\fl R'(\Of_F)$ is surjective. So $u$ induces a morphism of
algebras $\varphi:\Hecke\fl\Hecke'$ (for every $f\in\Hecke$ and every
$x\in R'(F)$, $\varphi(f)(x)$ is equal to $0$ if $x\not\in u(R(F))$ and to
$f(u^{-1}(x))$ if $x\in u(R(F))$).
For every $\delta$ in $R(F)$ or $R'(F)$, denote by $C(\delta)$ (resp.
$C_{st}(\delta)$) the $\theta$-conjugacy (resp. stable $\theta$-conjugacy)
class of $\delta$.

\begin{sublemme}\label{lemme:quotient_par_le_centre} Assume that $\theta$
acts trivially on $\Ho^1(F,Z_R)$.

Let $\delta\in R(F)$ be $\theta$-semi-simple; write $\delta'=u(\delta)$.
Then $u(C_{st}(\delta))=C_{st}(\delta')$. So there exists a (necessarily
finite) family $(\delta_i)_{i\in I}$ of elements of $R(F)$ that are stably
$\theta$-conjugate to $\delta$, such that $C(\delta')=
\coprod\limits_{i\in I}u(C(\delta_i))$. Moreover, for every
$f\in\Hecke$,
\[O_{\delta'\theta}(\varphi(f))=\sum_{i\in I}O_{\delta_i\theta}(f).\]

\end{sublemme}

(As always, we use the Haar measures on $R(F)$ and $R'(F)$ such that the
volumes of $R(\Of_F)$ and $R'(\Of_F)$ are equal to $1$.)

\begin{proof} It is clear that $u(C_{st}(\delta))\subset C_{st}(\delta')$.
We show the other inclusion. Let $\gamma'\in R'(F)$ be stably
$\theta$-conjugate to $\delta'$. As $u(R(F))=\Ker(R'(F)\fl\Ho^1(F,Z_R))$
is the intersection of kernels of $\theta$-stable characters of
$R'(F)$, lemma \ref{lemme:car_sur_classes_theta_stables} implies that there
exists $\gamma\in R(F)$ such that $\gamma'=u(\gamma)$. It is easy to see that
$\gamma$ and $\delta$ are stably $\theta$-conjugate.

Fix a family $(\delta_i)_{i\in I}$ as in the statement of the lemma,
and write $\K=R(\Of_F)$, $\K'=R'(\Of_F)$.
We show the equality of orbital integrals. Let $f\in\Hecke$. We may assume
that $f=\ungras_A$, where $A$ is a compact subset of $R(F)$ such that
$A=Z_R(F)\K A\K$. Then $\varphi(f)=\ungras_{u(A)}$, so
\[O_{\delta'\theta}(\varphi(f))=\sum_{\gamma'}\vol(u(A)\cap R'_{\gamma'\theta}
(F))^{-1},\]
where the sum is taken over a set of representatives $\gamma'$ of the
$\K'$-$\theta$-conjugacy classes of elements of $u(A)$ that are
$\theta$-conjugate to $\delta'$ (in $R'(F)$).
There are similar formulas for the twisted orbital integrals of $f$ at the
$\delta_i$ ; for every $i\in I$,
\[O_{\delta_i\theta}(f)=\sum_\gamma\vol((A\cap R_{\gamma\theta}(F))Z_R(F)/
Z_R(F))^{-1},\]
where the sum is taken over a set of representatives $\gamma$ of the
$\K$-$\theta$-conjugacy classes of elements of $A$ that are
$\theta$-conjugate to $\delta_i$ (in $R(F)$). To show the formula of the
lemma, it is therefore enough to notice that, for every $\gamma\in R(F)$,
$u$ induces an isomorphism $(A\cap R_{\gamma\theta}(F))Z_R(F)/Z_R(F)\iso u(A)
\cap R_{u(\gamma)\theta}$.

\end{proof}

We use again the notations of \ref{lemme_fondamental1}. Assume that
$\eta(1\rtimes\sigma)\in \widehat{\G}_{der}\rtimes\sigma$ and that
$s\in\widehat{\G}_{der}$. Then $(\H,{}^L\H,s,\eta)$ defines in an obvious
way endoscopic data $(\H',{}^L\H',s,\eta')$ for $\G':=
\G/Z(\G)^0$ (because $\widehat{\G}'=\widehat{\G}_{der}$). As in
\ref{lemme_fondamental1}, we get from this endoscopic data
$(\H',{}^L\H',t',\xi')$ for $(R',\theta)$, where $R'=R_{E/F}\G'_E$. 

\begin{sublemme}\label{lemme:quotient_par_le_centre_suite} Assume that
$\theta$ acts trivially on $\Ho^1(F,Z(R)^0)$ and that $\G$ and $\H$
satisfy the conditions of lemma \ref{lemme:H_equivariante}.

Then the fundamental lemma is true for $(R,\theta)$ and $(\H,{}^L\H,t,\xi)$
if and only if it is true for $(R',\theta)$ and $(\H',{}^L\H',t',\xi')$. 

\end{sublemme}

\begin{proof} Let $Z_R=Z(R)^0$, and let $Z_H$ be the image of $Z(\G)^0$ in
$Z(\H)$. Then $\H'=\H/Z_H$. It is easy to check that, if $\gamma_H\in\H(F)$
is semi-simple and strongly $\G$-regular, if $\delta\in R(F)$ is
$\theta$-semi-simple and strongly $\theta$-regular, and if $(\gamma_H,\delta)$
is sent to $(\gamma'_H,\delta')\in\H'(F)\times R'(F)$ by the obvious
projection, then
$\Delta_\xi(\gamma_H,\delta)=\Delta_{\xi'}(\gamma'_H,\delta')$.
Apply lemma \ref{lemme:H_equivariante}
with $\chi_H=1$ and $\chi_R=1$ (this is possible because the character
$\lambda_C$ that appears in this lemma is trivial, thank to the assumption
that $\eta(1\rtimes\sigma)\in\widehat{\G}_{der}\rtimes\sigma$). This lemma
shows that we may replace the Hecke algebras of $R$ and $\H$ by the Hecke
algebras of $Z_R(F)$-invariant or $Z_H(F)$-invariant functions.
To finish the proof, apply lemma \ref{lemme:quotient_par_le_centre}.

\end{proof}

The next lemma and its proof were communicated to me by Robert Kottwitz. (Any
mistakes that I may have inserted are my sole responsibility.)

\begin{sublemme}\label{lemme:groupes_qui_marchent} Let $F$ be a
non-archimedean local field of characteristic $0$, $\G$ be an adjoint
quasi-split group over $F$ and $(\H,s,\eta_0)$ be an endoscopic triple
for $\G$.
Fix a Borel subgroup $\B$ (resp. $\B_H$) of $\G$ (resp.
$\H$) and a Levi subgroup $\T_G$ (resp. $\T_H$) of $\B$ (resp. $\B_H$).
Then there exists an admissible embedding $\T_H\fl\G$ with image $\T_G$.

\end{sublemme}

\begin{proof} Write $\Gamma=\Gal(\overline{F}/F)$.
Choose embeddings $\widehat{\T}_G\subset\widehat{\B}\subset
\widehat{\G}$ and $\widehat{\T}_H\subset\widehat{\B}_H\subset\widehat{\H}$ that
are preserved by the action of $\Gamma$ on $\widehat{\G}$ and
$\widehat{\H}$. 

As $F$ is local, we may assume that $s\in Z(\widehat{\H})^\Gamma$.
By the definition of an endoscopic triple, for every $\tau\in\Gamma$, there
exists $g_\tau\in\widehat{\G}$ such that : for every $h\in\widehat{\H}$,
\renewcommand\theequation{$*$}
\begin{equation}
g_\tau\tau(\eta_0(h))g_\tau^{-1}=\eta_0(\tau(h)).
\end{equation}
In particular, the
$\widehat{\G}$-conjugacy class of $\eta_0(s)$ is fixed by the action of
$\Gamma$ on $\widehat{\G}$. By lemma 4.8 of \cite{Cl-LF}, $\eta_0(s)$ is
$\widehat{\G}$-conjugate to an element of $\widehat{\T}_G^\Gamma$. Replacing
$\eta_0$ by a $\widehat{\G}$-conjugate, we may assume that $\eta_0(s)\in
\widehat{\T}_G^\Gamma$. Then
\[\widehat{\T}_G\subset\Cent_{\widehat{\G}}(\eta_0(s))=\Cent_{\widehat{\G}}
(\eta_0(s))^0=\widehat{\H}\]
($\Cent_{\widehat{\G}}(\eta_0(s))$ is connected because $\widehat{\G}$ is
semi-simple and simply connected), so by further conjugating $\eta_0$ by
an element in $\eta_0(\widehat{\H})$ (which does not change $\eta_0(s)$,
since $s\in Z(\widehat{\H})$), we may also assume that $\eta_0(\widehat{\T}_H)
=\widehat{\T}_G$ and $\eta_0(\widehat{\B}_H)=\widehat{\B}\cap\eta_0(\widehat{
\H})$.

Since $\eta_0(s)$ is fixed by $\Gamma$, for every $\tau\in\Gamma$,
$g_\tau\eta_0(s)g_\tau^{-1}=\eta_0(s)$, so that $g_\tau\in\Cent_{\widehat{\G}}
(\eta_0(s))=\eta_0(\widehat{\H})$. Moreover $(*)$,
together with the fact that $\Gamma$ preserves $(\widehat{\B},\widehat{\T}_G)$
and $(\widehat{\B}_H,\widehat{\T}_H)$, implies that $h_\tau:=\eta_0^{-1}
(g_\tau)$ conjugates $(\widehat{\B}_H,\widehat{\T}_H)$ into itself. Therefore
$h_\tau\in\widehat{\T}_H$, and $(*)$ now shows that $\eta_0$ induces a
$\Gamma$-equivariant isomorphism $\widehat{\T}_H\iso\widehat{\T}_G$. Dual to
this is an admissible isomorphism $\T_H\iso\T_G$.

\end{proof}

Let $F$ be a non-archimedean local field of characteristic $0$.
Let $n,n_1,\dots,n_r\in\Nat^*$.
Set $\PGL_n=\GL_n/Z(\GL_n)$. For every quadratic extension $E$ of
$F$, set $\PGU(n,E)=\GU(n,E)/Z(\GU(n,E))$, where $\GU(n,E)$ is the
unitary group defined by the extension $E/F$ and by the
Hermitian form with matrix
\[J_n:=\left(\begin{array}{ccc}0 & & 1 \\ & \begin{turn}{45}\large\ldots
\end{turn}& \\ 1 & & 0\end{array}\right)\in\GL_n(\Z).\]
More generally, set $\Pro(\U(n_1,E)\times\dots\times\U(n_r,E))=
(\GU(n_1,E)\times\dots\times\GU(n_r,E))/Z$, where $Z=R_{E/\Q}\Gr_m$, embedded
diagonally.
Set $\PGSO_n=\GSO(J_n)/Z(\GSO(J_n))$, where $\GSO(J_n)=\GO(J_n)^0$,
and $\PGSp_{2n}=\GSp(J'_n)/Z(\GSp(J'_n))$, where
\[J'_n=\left(\begin{array}{cc}0 & J \\ -J & 0\end{array}\right)\in\GL_{2n}(\Z)
.\]
If $Y^1,\dots,Y^r\in\{\GSO,\GSp\}$, we denote by
$\Pro(Y^1_{n_1}\times\dots\times Y^r_{n_r})$ the quotient of
$Y^1_{n_1}\times\dots\times Y^r_{n_r}$ by $\Gr_m$ embedded diagonally.

\begin{sublemme}\label{lemme:groupes_qui_marchent2}
Let $\G$ be a simple adjoint unramified group over $F$.
\begin{itemize}
\item[(i)] If $\G$ is of type $A$, then there exists a finite unramified
extension $K$ of $F$, a quadratic unramified extension $E$ of $K$ and a
non-negative integer $n$ such that $\G=R_{K/F}\PGL_n$ or $\G=R_{K/F}\PGU(n,E)$.
If $\G=R_{K/F}\PGL_n$, then $\G$ has no non-trivial elliptic endoscopic
groups. If $\G=R_{K/F}\PGU(n,E)$, then the elliptic endoscopic groups of
$\G$ are the $R_{K/F}\Pro(\GU(n_1,E)\times\GU(n_2,E))$, with $n_1,n_2\in\Nat$
such that $n=n_1+n_2$ and that $n_2$ is even.

\item[(ii)] If $\G$ is of type $B$, then there exists a finite unramified
extension $K$ of $F$ and a non-negative integer $n$ such that
$\G=R_{K/F}\PGSO_{2n+1}$. The elliptic endoscopic groups of $\G$ are the
$R_{K/F}\Pro(\GSO_{2n_1+1}\times\GSO_{2n_2+1})$, with $n_1,n_2\in\Nat$ such
that $n=n_1+n_2$.

\item[(iii)] If $\G$ is of type $C$, then there exists a finite unramified
extension $K$ of $F$ and a non-negative integer $n$ such that
$\G=R_{K/F}\PGSp_{2n}$. The elliptic endoscopic groups of $\G$ are the
$\Pro(\GSO_{2n_1}\times\GSp_{2n_2})$, with $n_1,n_2\in\Nat$ such that
$n=n_1+n_2$ and $n_1\not=1$.

\end{itemize}

\end{sublemme}

In particular, if $\G$ is adjoint of type $A$, $B$ or $C$, then the
hypothesis of proposition \ref{prop:utilisation_donnees_locales} on the
center of $\H$ (ie that this center be connected) is satisfied.

\begin{proof} Let $K$ be the smallest extension of $F$ on which $\G$ splits,
and fix a generator $\sigma$ of $\Gal(K/F)$.
Then $\G_K\simeq (\G')^r$, where $r\in\Nat^*$ and $\G'$ is an adjoint
absolutely simple group over $K$. Let $\theta$ be the automorphism
(over $K$) of $(\G')^r$ induced by $\sigma$. If $\G$ is of type $B$ or $C$,
then $\G'$ is also of type $B$ or $C$, so $\G'$ is equal to $\PGSO_n$ or
$\PGSp_{2n}$, and $\G'$ has no non-trivial outer automorphisms
(cf \cite{Di} IV.6 and IV.7), so we may assume that $\theta$ acts by
permuting the factors of $(\G')^r$. As $K$ is the smallest extension on which
$\G$ splits, $\theta$ has to be a $n$-cycle. So
$\G\simeq R_{K/F}\G'$. To compute the elliptic endosocopic triples for $\G$,
we may assume that $K=F$. Then $\G$ is split and has a connected center,
so its endoscopic groups are also split (cf definition 1.8.1 of \cite{Ng}).
From this observation, it is easy to see that the elliptic endoscopic
groups of $\G$ are the ones given in the statement of the lemma.

Assume that $\G$ is of type $A$. Then there exists $n\in\Nat^*$ such that
$\G'=\PGL_{n,K}$, and $\Out(\G')$ is isomorphic to $\Z/2\Z$ (cf \cite{Di}
IV.6). We may assume that $\theta\in (\Z/2\Z)^r\rtimes\Sgoth_r$, where
$(\Z/2\Z)^r$ acts on $(\G')^r$ via the isomorphism $\Z/2\Z\simeq
\Out(\G')$
(and a splitting of $\Aut(\G')\fl\Out(\G')$) and $\Sgoth_r$ acts on
$(\G')^r$ by permuting the factors. Write $\theta=\epsilon\rtimes\tau$,
with $\varepsilon\in(\Z/2\Z)^r$ and $\tau\in\Sgoth_n$. As in the first case,
$\tau$ has to be a $n$-cycle. After conjugating $\tau$ by an element of
$(\Z/2\Z)^r\rtimes\Sgoth_r$, we may assume that $\varepsilon\in
\{(1,\dots,1),(-1,1,\dots,1)\}$ (because $\varepsilon_1\rtimes\tau$ and
$\varepsilon_2\rtimes\tau$ are conjugate if and only if there exists
$\eta\in (\Z/2\Z)^r$ such that $\varepsilon_1\varepsilon_2=\eta\tau(\eta)$,
and the image of the morphism $(\Z/2\Z)^r\fl (\Z/2\Z)^r$, $\eta\fle\eta\tau
(\eta)$ is $\{(e_1,\dots,e_r)\in(\Z/2\Z)^r|e_1\dots e_r=1\}$). If $\theta=(1,
\dots,1)\rtimes\tau$, then $\G\simeq R_{K/F}\PGL_n$, and it is not hard to
see that $\G$ has no non-trivial elliptic endoscopic triples.
Assume that $\theta=(-1,1,\dots,1)\rtimes\tau$. Then $\theta$ is of order
$2r$, so $[K:F]=2r$ and $\G=R_{K'/F}\PGU(n,K)$, where $K'$ is the subfield of
$K$ fixed by $\theta^r(=(-1,\dots,-1)\rtimes 1)$. The calculation of the
elliptic endoscopic triples of $\G$ is done just as in proposition
\ref{prop:groupes_endoscopiques} (with the obvious changes).

\end{proof}

\section{Results}
\label{lemme_fondamental5}

\begin{subproposition}\label{prop:lemme_fondamental_cas_adjoint}
Let $X\in\{A,B,C\}$. Let $F$ be a non-archimedean local field and
$\G$ be an adjoint unramified group over $F$, of type $X$.
Assume that there exists $N\in\Nat^*$ such that, for all
$F'$, $E'$, $\G'$, $R'$ and $(\H',{}^L\H',t',\xi')$ as in
\ref{lemme_fondamental1}, 
the twisted fundamental lemma is true for the unit of the Hecke algebra if
$\G'$ is adjoint of type $X$, $\dim(\G')\leq\dim(\G)$ and the residual
characteristic of $F'$ does not divide $N$.

Then, for every finite unramified extension $E$ of $F$ and for all twisted
endoscopic data $(\H,{}^L\H,t,\xi)$ for $R:=R_{E/F}\G_E$ as in
\ref{lemme_fondamental1}, the twisted fundamental lemma is true for $R$
and $(\H,{}^L\H,t,\xi)$ and for all the functions in the Hecke algebra.

\end{subproposition}

\begin{proof} By lemma \ref{lemme:descente}, lemma
\ref{lemme:groupes_qui_marchent}, lemma
\ref{lemme:groupes_qui_marchent2}
lemma \ref{lemme:situation_globale},
proposition \ref{prop:existence_donnees_locales} and proposition
\ref{prop:utilisation_donnees_locales}, the twisted fundamental lemma for
$\G$ follows from the twisted fundamental lemma for all proper Levi
subgroups of $\G$ (if $\G$ has no elliptic maximal torus, then lemma
\ref{lemme:descente} is enough to see this). But, by the classification of
adjoint unramified groups of type $X$ given in lemma
\ref{lemme:groupes_qui_marchent2}, every proper Levi subgroup of $\G$ is
isomorphic to a group $\G_0\times\G_1\times\dots\times\G_r$, with $\G_1,
\dots,\G_r$ of the form $R_{K/F}\GL_m$, where $K$ is a finite unramified
extension of $F$ and $m\in\Nat^*$, and $\G_0$ adjoint unramified of type $X$
and such that $\dim(\G_0)<\dim(\G)$. If $1\leq i\leq r$, $\G_i$ has no
non-trivial elliptic endoscopic groups, so the twisted fundamental lemma
for $\G_i$ follows from descent (lemma \ref{lemme:descente}) and from the
fundamental lemma for stable base change, that has been proved in the case
of general linear groups by Arthur and Clozel
(\cite{AC}, chapter I, proposition 3.1).
Hence, to prove the proposition, it suffices to reason by induction on the
dimension of $\G$.

\end{proof}

\begin{subcorollaire}\label{cor:LF_regulier}
We use the notations of \ref{lemme_fondamental1}. If $F=\Q_p$, $\G$ is one
of the unitary groups
$\G(\U^*(n_1)\times\dots\times\U^*(n_r))$ of \ref{groupes1}
and the morphism $\eta$ is the morphism $\eta_{simple}$ of
\ref{partie_en_p2}, then the twisted fundamental lemma is true.

\end{subcorollaire}

\begin{proof} As the center of $\G$ is connected, the corollary follows
from proposition \ref{prop:lemme_fondamental_cas_adjoint} above and
from lemma \ref{lemme:quotient_par_le_centre_suite}, so it is enought
to check that the hypotheses of this lemma are satisfied.
The endoscopic triple $(\H,s,\eta_0)$ satisfies the hypotheses of lemma
\ref{lemme:H_equivariante}, by the explicit description of the endoscopic
triples of $\G$ given in proposition \ref{prop:groupes_endoscopiques}.
It is obvious that $s\in\widehat{\G}_{der}$
and $\eta(1\rtimes\sigma)\in\widehat{\G}_{der}\rtimes\sigma$.
Finally, the center of $\G$ is an induced torus, so its first Galois cohomology
group on any extension of $F$ is trivial.

\end{proof}

The result that we really need in this text is formula $(*)$ of
\ref{FT_stable_geometrique3}. We recall this formula. Notations are still as
in \ref{lemme_fondamental1}, with $E$ a field. Let $\Delta_\eta$ be the
transfer factors for the morphism $\eta:{}^L\H\fl {}^L\G$, with the
normalization given by the $\Of_F$-structures on $\H$ and $\G$ (cf \cite{H}
II 7 or \cite{Wa3} 4.6).
If $\delta\in R(F)$ is $\theta$-semi-simple and
$\gamma\in\Norme\delta$, Kottwitz defined in \cite{K-SVLR} \S7 p 180
an element $\alpha_p(\gamma,\delta)$ of $X^*(Z(\G_\gamma)^{\Gamma_F})$
(remember that $\G_\gamma=\Cent_\gamma(\G)^0$).
The result that we want to prove is the following : For every
$\gamma_H\in\H(F)$ semi-simple, for every $f\in\Hecke_R$, let $\gamma$ be
an image of $\gamma_H$ in $\G(F)$ (such a $\gamma$ exists because $\G$ is
quasi-split). Then
\renewcommand\theequation{$*$}
\begin{equation}
SO_{\gamma_H}(b_\xi(f))=\sum_{\delta}<\alpha_p(\gamma,\delta),s>\Delta_\eta(
\gamma_H,\gamma)e(R_{\delta\theta})O_{\delta\theta}(f),
\end{equation}
where the sum is taken over the set of $\theta$-semi-simple
$\theta$-conjugacy classes $\delta$ of $R(F)$ such that
$\gamma\in\Norme\delta$, $R_{\delta\theta}$ is the connected compoenent of $1$
of the centralizer of $\delta\theta$ in $R$ and
$e(R_{\delta\theta})$ is the sign defined by Kottwitz in \cite{K-SC}.

\begin{subcorollaire}\label{lemme:ss_non_regulier} Assume that $F=\Q_p$ and
that $\G$ is one of the unitary groups $\G(\U^*(n_1)\times\dots\times
\U^*(n_r))$ of \ref{groupes1}. Then formula $(*)$ above is true.

\end{subcorollaire}

\begin{proof} If $\gamma_H$ is strongly regular, then formula $(*)$
follows from corollary \ref{cor:LF_regulier} and from corollary
\ref{RKcor.jj} of the appendix.

The reduction from the general case to the case where $\gamma_H$ is
strongly regular is done in section \ref{app:ghreg} of the appendix (see in 
particular proposition \ref{RKprop.ghreg}).

\end{proof}

\begin{subremarque} The last two corollaries are also true
(with the same proof)
for any group $\G$ with connected center and such that all its endoscopic
triples satisfy the conditions of lemma
\ref{lemme:quotient_par_le_centre_suite}. Examples of such groups are
the symplectics groups of \cite{M3} (cf proposition 2.1.1 of \cite{M3}).

\end{subremarque}

\appendix
\chapter{Comparison of two versions of twisted transfer factors}
\centerline{\bf R. Kottwitz}
\vskip .4in
\counterwithin{equation}{subsection}

In order to stabilize the Lefschetz formula for Shimura varieties over
finite fields, one needs to use  twisted transfer factors for cyclic base
change. Now these twisted transfer factors  can be
expressed in terms of standard transfer factors, the ratio between the two
 being given by a Galois cohomological factor involving an
invariant denoted by $\inv(\gamma,\delta)$ in \cite{KS}. However, in the
stabilization of the Lefschetz formula it is more natural to use a different
invariant $\alpha(\gamma,\delta)$. The purpose of this appendix is to relate
the invariants $\inv(\gamma,\delta)$ and $\alpha(\gamma,\delta)$ (see Theorem
\ref{RKthm.main1}), and then to  justify the use made in \cite{K-SVLR} of  transfer
factors 
\[
\Delta_0(\gamma_H,\delta)=\Delta_0(\gamma_H,\gamma)\langle
\alpha(\gamma,\delta),s\rangle ^{-1},
\] 
 first in the case when $\gamma_H$ is strongly 
$G$-regular semisimple (see Corollary \ref{RKcor.jj}) and then in the more general
case in which $\gamma_H$ is  assumed only to be $(G,H)$-regular (see Proposition 
\ref{RKprop.ghreg}, where, however, the derived group of $G$ is assumed to be
simply connected).

 I would like to thank Sophie Morel for her very helpful comments on
a first version of this appendix. 

\section{Comparison of $\Delta_0(\gamma_H,\delta)$ and
$\Delta_0(\gamma_H,\gamma)$}
In the case of cyclic base change the twisted transfer factors 
$\Delta_0(\gamma_H,\delta)$  of \cite{KS}  are closely related to the
standard transfer factors  $\Delta_0(\gamma_H,\gamma)$ of \cite{LS1}. 
This fact, first observed by Shelstad \cite{Sh2} in the case of base change
for $\mathbb C/\mathbb R$, was one of several guiding principles used to
arrive at the general twisted transfer factors defined in \cite{KS}. Thus
there is nothing really new in this section. After reviewing some basic
notions, we prove Proposition \ref{RKprop.e}, which gives the precise
relationship between $\Delta_0(\gamma_H,\delta)$ and
$\Delta_0(\gamma_H,\gamma)$. 

\subsection{Set-up} We consider a finite cyclic extension $E/F$ of local
fields of characteristic zero. We put $d:=[E:F]$ and choose a generator
$\sigma$ of $\Gal(E/F)$. In addition we choose an algebraic closure
$\overline F$ of $F$ that contains $E$. We write $\Gamma$ for the absolute
Galois group 
$\Gal(\overline F/F)$ and
$W_F$ for the absolute Weil group of $F$. There is then a canonical
homomorphism $W_F \to \Gamma$ that will go unnamed. 

We also consider a quasisplit connected reductive group $G$ over $F$. Put
$R_G:=\Res_{E/F}(G_E)$, where $G_E$ is the $E$-group obtained from $G$ by
extension of scalars, and $\Res_{E/F}$ denotes Weil's restriction of
scalars. As usual there is a natural automorphism $\theta$ of $R_G$ inducing
$\sigma$ on $G(E)$ via the canonical identification $R_G(F)=G(E)$. 

\subsection{Description of $\hat R_G$} 
For any $\Gamma$-group $A$ (that is, a group $A$ equipped with an action of
$\Gamma$) we obtain by restriction a $\Gamma_E$-group $A_E$ (with $\Gamma_E$
denoting the subgroup $\Gal(\overline F/E)$ of $\Gamma$), and we write
$I(A)$ for the $\Gamma$-group obtained from $A_E$ by induction from
$\Gamma_E$ to $\Gamma$. 

Then $I(A)$ has the following description in terms of $A$. Let $J$ denote
the set of embeddings of $E$ in $\overline F$ over $F$, with $j_0$ denoting
the inclusion $E \subset \overline F$. The group $\Gamma$ acts on the left
of $J$ by $\tau j:=\tau \circ j$ (for $\tau \in \Gamma$, $j \in J$), and the
group $\Gal(E/F)$ acts on the right of $J$ by
$j\sigma^i:=j\circ \sigma^i$. An element $x \in I(A)$ is then a map
$j\mapsto x_j$ from $J$ to $A$. An element $\tau \in \Gamma$ acts on $x \in
I(A)$ by 
\[
(\tau x)_j:=\tau(x_{\tau^{-1}j}).
\] 
There is a \emph{right} action of $\Gal(E/F)$ on the $\Gamma$-group $I(A)$
given by 
\[
(x\sigma^i)_j:=x_{j\sigma^{-i}}.
\] 

We have $\hat R_G=I(\hat G)$ as $\Gamma$-group. Bearing in mind that for any
automorphisms $\theta_1$, $\theta_2$ of a connected reductive group
 one has the rule
$\widehat{\theta_1\theta_2}=\hat\theta_2\hat\theta_1$,  we see that the
natural left action of $\Gal(E/F)$ on $R_G$ is converted into a \emph{right}
action of
$\Gal(E/F)$ on $\hat R_G$, and hence that the automorphism $\hat\theta$ of
$\hat R_G$ is given by 
\[
(\hat\theta x)_j=x_{j \sigma^{-1}}.
\] 

There is an obvious embedding
\[
A \hookrightarrow I(A)
\]
of $\Gamma$-groups, sending $a \in A$ to the constant map $J \to A$ with
value $a$, and this map identifies $A$ with the group of fixed points of
$\Gal(E/F)$ on $I(A)$. In particular we get 
\begin{equation*}\label{RKdefi}
i:\hat G \simeq (\hat R_G)^{\hat \theta} \hookrightarrow \hat R_G, 
\end{equation*}
which we extend to an embedding 
\[
i:\LG \to \LR_G
\]
by mapping $g\tau$ to $i(g)\tau$ (for $g \in \hat G$, $\tau \in W_F$). Note
that $i(\LG)$ is the group of fixed points of the automorphism $\Ltheta$ of
$\LR_G$ defined by 
\[
\Ltheta(x\tau):=\hat\theta(x)\tau
\]
for $x \in \hat R_G$, $\tau \in W_F$. 

\subsection{Endoscopic groups and twisted endoscopic groups} 
Let $(H,s,\eta)$ be an endoscopic datum for $G$. Thus $s \in Z(\hat
H)^\Gamma$ and $\eta:\LH \to \LG$ is an $L$-homomorphism that restricts to
an isomorphism $\hat H \to (\hat G_{\eta(s)})^{\circ}$. When the derived
group of $G$ is not simply connected, we should actually allow for a
$z$-extension of $H$, as in \cite{LS1} and \cite{KS}, but since this wrinkle
does not perturb the arguments below in any non-trivial way, we prefer to
ignore it. 

Following Shelstad \cite{Sh2}, we now explain how to regard $H$ as a twisted
endoscopic group for $(R_G,\theta)$. Let $\mathcal Z$ denote the centralizer
of $i\eta(\hat H)$ in $\hat R_G$. Since the centralizer of $\eta(\hat H)$ in
$\hat G$ is $\eta(Z(\hat H))$, we see that $\mathcal Z$ is the subgroup of
$\hat R_G$ consisting of all maps $J \to \eta(Z(\hat H))$. Thus, as a
group, $\mathcal Z$ can be identified with $I(Z(\hat H))$. Since $Z(\hat
H)$ is a $\Gamma$-group, so too is $I(Z(\hat H))=\mathcal Z$, but the
embedding
$\mathcal Z \hookrightarrow \hat R_G$ is \emph{not} $\Gamma$-equivariant. The
subgroup $\mathcal Z$ is however stable under $\hat \theta$.

Using $s \in Z(\hat H)^{\Gamma}$, we now define an element $\tilde s \in
\mathcal Z$ by the rule 
\begin{equation}
\tilde s_j :=
\begin{cases}
s &\text{ if $j=j_0$,}\\
1  &\text{ if $j\ne j_0$.}
\end{cases}
\end{equation}
Thus $\tilde s$ maps to $s$ under the norm map $\mathcal Z \to Z(\hat H)$
(given by $x \mapsto \prod_{j \in J}x_j$). It is easy to see that the
composition 
\[
\hat H \xrightarrow{\eta} \hat G \xrightarrow{i} \hat R_G 
\] 
identifies $\hat H$ with the identity component of the $\hat\theta$-centralizer of
$\tilde s$ in
$\hat R_G$. 

\subsection{Allowed embeddings}\label{RKsub.allow} 
We now have part of what is needed to view
$H$ as a twisted endoscopic group for $(R_G,\theta)$, but in addition to
$\tilde s$ we need  suitable
$\tilde\eta:\mathcal H \to \LR_G$.  In the situation of interest in
the next section of this appendix, we may even take $\mathcal H=\LH$, so
this is the only case we will discuss further.

When $\mathcal H =\LH$, in order to get a twisted endoscopic datum
$(H,\tilde s,\tilde \eta)$ for $(R_G,\theta)$, we need for $\tilde\eta:\LH
\to \LR_G$ to be one of Shelstad's \emph{allowed embeddings} \cite{Sh2},
which is to say that $\tilde\eta$, $i\eta$ must have the same restriction to
$\hat H$, and  that $\tilde\eta(\LH)$ must be contained in the group of fixed
points of the automorphism $\Int(\tilde s)\circ \Ltheta$ of $\LR_G$. 

In subsection \ref{RKsub.mae} we will see that, when $E/F$ is an unramified
extension of $p$-adic fields and $\sigma$ is the Frobenius automorphism,
there exists a canonical allowed embedding $\tilde\eta$ determined by
$\tilde s$. In this section, however, we work with an arbitrary allowed
embedding. 

We are going to use $\tilde\eta$ to produce a $1$-cocycle of $W_F$ in
$\mathcal Z \xrightarrow{1-\hat \theta} \mathcal Z$. For this we need to
compare (as in \cite{KS}) $\tilde\eta$ to the $L$-homomorphism 
\[
i\eta:\LH \to \LG \to \LR_G. 
\]
Since $\tilde\eta$ and $i\eta$ agree on $\hat H$, there is a unique
$1$-cocycle $a$ of $W_F$ in $\mathcal Z$ such that 
\[
\tilde\eta(\tau)=a_\tau i\eta(\tau)
\]
for all $\tau \in W_F$. The pair $(a^{-1},\tilde s)$ is a $1$-cocycle of
$W_F$ in $\mathcal Z \xrightarrow{1-\hat\theta}\mathcal Z$.  Here one must not
forget that the $\Gamma$-action on
$\mathcal Z$ comes from viewing it as $I(Z(\hat H))$. In fact the map
$\tilde \eta \mapsto a$ sets up a bijection between allowed embeddings
$\tilde\eta$ and $1$-cocycles $a$ of $W_F$ in $\mathcal Z$ such that
$(a^{-1},\tilde s)$ is a $1$-cocycle in $\mathcal Z
\xrightarrow{1-\hat\theta}\mathcal Z$.

\subsection{Canonical twisted and standard transfer factors}\label{RKsub.tfv} 
We now choose an $F$-splitting
\cite[p.~224]{LS1} for our quasisplit group $G$. This choice determines
canonical transfer factors $\Delta_0(\gamma_H,\gamma)$ 
(see \cite[p.~248]{LS1}). 

Our $F$-splitting of $G$ can also be viewed as a $\sigma$-invariant
$E$-splitting, and therefore gives rise to an $(F,\theta)$-splitting
\cite[p.~61]{KS} of $R_G$, which then determines
canonical twisted transfer factors $\Delta_0(\gamma_H,\delta)$ (see
\cite[p.~62]{KS}). Our goal is to express $\Delta_0(\gamma_H,\delta)$ as the
product of $\Delta_0(\gamma_H,\gamma)$ and a simple Galois cohomological
factor involving an invariant $\inv(\gamma,\delta)$ that we are now going to
discuss. 

It may be useful to recall (though it will play no role in this appendix) that
when $G$ is unramified, and we fix an $\mathcal O$-structure on $G$ for
which $G(\mathcal O)$ is a hyperspecial maximal compact subgroup of $G(F)$, there
is an obvious notion of $\mathcal O$-splitting, namely an $F$-splitting that is
defined over
$\mathcal O$ and reduces modulo the maximal ideal in $\mathcal O$ to a splitting
for the special fiber of $G$. When  such an $\mathcal
O$-splitting is used, and $H$ is also unramified,  the transfer factors
$\Delta_0(\gamma_H,\gamma)$ so obtained are the ones needed for the fundamental
lemma for  the spherical Hecke algebra on $G$ obtained from $G(\mathcal
O)$. In the case that $E/F$ is unramified, the same is true for the twisted
fundamental lemma for the spherical Hecke algebra for $G(E)$ obtained from
$G(\mathcal O_E)$. 

\subsection{Definition of the invariant
$\inv(\gamma,\delta)$}\label{RKsub.digd}  
 We consider a
maximal
$F$-torus $T_H$ of $H$ and an admissible isomorphism $T_H \simeq T$ between
$T_H$ and a maximal $F$-torus $T$ of $G$. We consider $\gamma_H$ in $T_H(F)$
whose image $\gamma$ in $T(F)$ is strongly $G$-regular. The 
standard transfer factor $\Delta_0(\gamma_H,\gamma)$ is then defined.    
We also consider $\delta\in R_G(F)=G(E)$ whose abstract norm \cite[3.2]{KS}
is the stable conjugacy class of $\gamma$. The twisted transfer
factor
$\Delta_0(\gamma_H,\delta)$ is then defined. 

The position of $\delta$ relative to $\gamma$ is measured by 
\[
\inv(\gamma,\delta) \in H^1(F, R_T \xrightarrow{1-\theta}R_T),
\] 
whose definition \cite[p.~63]{KS} we now recall. Our assumption that the
abstract norm of $\delta$ is $\gamma$ does \emph{not} imply that $\delta$ is
stably $\theta$-conjugate to an element in the $F$-points of the
$\theta$-stable maximal $F$-torus $R_T$ of $R_G$. It does however imply that
there exists $g \in R_G(\overline F)$ such that $g(N\delta)g^{-1}=\gamma$,
where $N:R_G \to R_G$ is the $F$-morphism $x \mapsto x\theta(x)\theta^2(x)
\dotsm
\theta^{d-1}(x)$, and $\gamma$ is viewed as an element of $R_G(F)=G(E)$ via
the obvious inclusion $G(F)\subset G(E)$. Put
$\delta':=g\delta\theta(g)^{-1}$ and define a $1$-cocycle $t$ of  $\Gamma$
by 
$t_\tau:=g\tau(g)^{-1}$ (for
$\tau
\in
\Gamma$). Note that the strong regularity of $\gamma$ implies that its
centralizer in $R_G$ is $R_T$. 

\begin{lemma}
The pair $(t^{-1},\delta')$ is a $1$-cocycle of $\Gamma$ in $R_T
\xrightarrow{1-\theta} R_T$. 
\end{lemma}
\begin{proof}
We first check  that $\delta' \in R_T(\overline F)$. Observe that
$N\delta'=\gamma$. Therefore 
\[
\gamma=\theta(\gamma)=\theta(N\delta')=(\delta')^{-1}(N\delta')\delta'=
(\delta')^{-1}\gamma\delta',
\]
which shows that $\delta'$ centralizes $\gamma$ and hence lies in $R_T$. We
note for later use that the $\theta$-centralizer of $\delta'$ in $R_G$ is
$T$, viewed as the subtorus of $\theta$-fixed points in $R_T$. 

Next, a short calculation using the definitions of $\delta'$ and $t_\tau$
shows that 
\begin{equation}\label{RKeqn.cc}
(\delta')^{-1}t_\tau \tau(\delta')=\theta(t_\tau).
\end{equation}  
To see that $t_\tau \in R_T(\overline F)$, we begin by noting that
\eqref{RKeqn.cc} says that $t_\tau$ $\theta$-conjugates $\tau(\delta')$ into
$\delta'$. 
 Now 
\[
N(\tau(\delta'))=\tau(N\delta')=\tau(\gamma)=\gamma=N(\delta'),
\]
showing that the two elements $\tau(\delta')$ and $\delta'$ in $R_T$ have
the same image under the norm homomorphism $N:R_T \to R_T$, and hence that
there exists $u \in R_T(\overline F)$ that $\theta$-conjugates $\delta'$ into
$\tau(\delta')$. Thus $t_\tau u$ lies in the $\theta$-centralizer (namely
$T=R_T^\theta$) of $\delta'$, which implies that $t_\tau$ lies in $R_T(\overline
F)$.   

The $1$-cocycle condition for $(t^{-1},\delta')$ is none other than
\eqref{RKeqn.cc}, and the proof is complete.  
\end{proof}
 
\begin{definition}
We define $\inv(\gamma,\delta)$ to be the class in $H^1(F, R_T
\xrightarrow{1-\theta}R_T)$ of the $1$-cocycle $(t^{-1},\delta')$.
\end{definition} 

\subsection{Main proposition} 
The last thing to do before stating Proposition \ref{RKprop.e} is to relate
$\mathcal Z$ to $\hat R_T$. This is very easy. Since $T_H$ is a maximal
torus in $H$, there is a canonical $\Gamma$-equivariant embedding $Z(\hat H)
\hookrightarrow \hat T_H$. Our admissible isomorphism  $T_H\simeq T$
 yields $\hat T_H \simeq\hat T$, so that we end up with a
$\Gamma$-equivariant embedding $Z(\hat H) \hookrightarrow \hat T$, to
which we may apply our restriction-induction functor $I$, obtaining a
$\Gamma$-equivariant embedding 
\[
k:\mathcal Z \hookrightarrow \hat R_T,
\]  
which is compatible with the $\hat\theta$-actions as well. We then obtain an
induced homomorphism 
\begin{equation}\label{RKeqn.hom} 
H^1(W_F, \mathcal Z \xrightarrow{1-\hat\theta}\mathcal Z) \to H^1(W_F, \hat
R_T
\xrightarrow{1-\hat\theta}\hat R_T).
\end{equation} 

Near the end of subsection \ref{RKsub.allow} we used $\tilde s$,$\tilde \eta$
to produce a $1$-cocycle $(a^{-1},\tilde s)$ of $W_F$ in $\mathcal Z
\xrightarrow{1-\hat\theta}\mathcal Z$, to which we may apply the homomorphism
$k$, obtaining a $1$-cocycle in $ \hat R_T
\xrightarrow{1-\hat\theta}\hat R_T$, which, since $k$ is injective, we may as
well continue to denote simply by $(a^{-1},\tilde s)$. Recall from Appendix
A of \cite{KS} that there is a $\mathbb C^\times$-valued pairing $\langle
\cdot,\cdot \rangle$ between $H^1(F, R_T
\xrightarrow{1-\theta}R_T)$ and $H^1(W_F, \hat R_T
\xrightarrow{1-\hat\theta}\hat R_T)$. Thus it makes sense to form the complex
number 
\[
\langle \inv(\gamma,\delta),(a^{-1},\tilde s)\rangle.
\] 

\begin{proposition}\label{RKprop.e}
There is an equality 
\[
\Delta_0(\gamma_H,\delta)=\Delta_0(\gamma_H,\gamma) 
\langle \inv(\gamma,\delta),(a^{-1},\tilde s)\rangle.
\]
\end{proposition}
\begin{proof}
Since the restricted root system \cite[1.3]{KS} of $R_T$ can be identified
with the root system of $T$, we may use the same $a$-data and $\chi$-data
for $T$ and $R_T$. When this is done, one has 
\begin{align*}
\Delta_I(\gamma_H,\delta)&=\Delta_I(\gamma_H,\gamma)\\
\Delta_{II}(\gamma_H,\delta)&=\Delta_{II}(\gamma_H,\gamma)\\
\Delta_{IV}(\gamma_H,\delta)&=\Delta_{IV}(\gamma_H,\gamma).
\end{align*} 

It remains only to prove that 
\[
\Delta_{III}(\gamma_H,\delta)=\Delta_{III}(\gamma_H,\gamma) 
\langle \inv(\gamma,\delta),(a^{-1},\tilde s)\rangle.
\] 
To do so we must recall how $\Delta_{III}$ is defined. We use (see
\cite{LS1})   the  chosen $\chi$-data to obtain embeddings 
\begin{align*}
\xi_1:\LT \hookrightarrow \LG,\\
\xi_2:\LT \hookrightarrow \LH.
\end{align*} 
Replacing $\xi_1$ by a conjugate under $\hat G$, we may assume that
$\eta\xi_2$ and $\xi_1$ agree on $\hat T$, and then there exists a unique
$1$-cocycle $b$ of $W_F$ in $\hat T$ so that 
\[
(\eta\xi_2)(\tau)=\xi_1(b_\tau \tau)
\] 
for all $\tau \in W_F$. We then have (see \cite[p.~246]{LS1}) 
\[
\Delta_{III}(\gamma_H,\gamma)=\langle \gamma,b\rangle, 
\] 
where $\langle \cdot,\cdot \rangle$ now denotes the Langlands pairing
between $T(F)$ and $H^1(W_F,\hat T)$. 

Similarly we have two embeddings $i\xi_1,\tilde\eta\xi_2:\LT \hookrightarrow
\LR_G$ that agree on 
$\hat T$, and therefore there exists a unique $1$-cocycle $c$ of $W_F$ in
$\hat R_T$ (which arises here because it is the centralizer in $\hat R_G$ of
$(i\xi_1)(\hat T)$) such that 
\[
(\tilde \eta\xi_2)(\tau)=c_\tau\bigl((i\xi_1)( \tau)\bigr) 
\] 
for all $\tau \in W_F$. Then  
$(c^{-1},\tilde s)$ is a $1$-cocycle of $W_F$ in $\hat R_T
\xrightarrow{1-\hat\theta}\hat R_T$, and  (see pages 40 and 63
of \cite{KS}) 
\[
\Delta_{III}(\gamma_H,\delta)=\langle \inv(\gamma,\delta),(c^{-1},\tilde
s) \rangle. 
\] 

It is clear from the definitions that the $1$-cocycles $a$, $b$, $c$ satisfy
the equality $c=ab$, in which we use $k:\mathcal Z \hookrightarrow \hat R_T$
and $\hat T=(\hat R_T)^{\hat \theta} \hookrightarrow \hat R_T$ to view
$a$, $b$ as $1$-cocycles in $\hat R_T$. Therefore 
\[
(c^{-1},\tilde s)=(a^{-1},\tilde s)(b^{-1},1),
\] 
which shows that 
\[
\Delta_{III}(\gamma_H,\delta)=\langle \inv(\gamma,\delta),(a^{-1},\tilde
s) \rangle \ \langle \inv(\gamma,\delta),(b,1)\rangle^{-1}.
\] 
It remains only to observe that $\langle
\inv(\gamma,\delta),(b,1)\rangle^{-1}
=\langle \gamma,b \rangle$, a consequence of the first part of Lemma
\ref{RKlem.cmp}, to be proved next. Here we use the obvious fact that
the image of
$\inv(\gamma,\delta)$ under $H^1(F, R_T
\xrightarrow{1-\theta}R_T) \to T(F)$ is $\gamma$. 
\end{proof} 
\subsection{Compatibility properties for the pairing $\langle \cdot,\cdot
\rangle$} 
In this subsection we consider a homomorphism $f:T \to U$ of $F$-tori. We
follow all the conventions of Appendix A in \cite{KS} concerning $H^1(F,T
\to U)$ and $H^1(W_F,\hat U \to \hat T)$. We denote by $K$ the kernel of $f$
and by $C$ the cokernel of $f$. Of course $C$ is necessarily a torus, and we
now assume that 
$K$ is also a torus. Dual to the exact sequence 
\[
1 \to K \to T \to U \to C \to 1
\] 
is the exact sequence     
\[
1 \to \hat C \to \hat U \to \hat T \to \hat K \to 1,
\] 
which we use to identify $\hat C$ with $\ker\hat f$ and $\hat K$ with $\cok
\hat f$. From \cite[p.~119]{KS} we obtain two long exact sequences, the
relevant portions of which are 
\begin{align*}
H^1(F,K)\xrightarrow{i'}&H^1(F,T\to U) \xrightarrow{j'} C(F)\\
H^1(W_F,\hat C) \xrightarrow{\hat{i}'}&H^1(W_F,\hat U \to \hat T)
\xrightarrow {\hat{j}'}\hat K^\Gamma.
\end{align*}
The following lemma concerns the compatibility of these two exact sequences
with the pairing \cite{KS} between $H^1(F,T\to U)$ and 
$H^1(W_F,\hat U \to \hat T)$. 

\begin{lemma}\label{RKlem.cmp} The pairing $\langle \cdot,\cdot \rangle$
satisfies the following two compatibilities. 
\begin{enumerate}
\item Let $x \in H^1(F,T\to U)$ and $c \in H^1(W_F,\hat C)$. Then 
\begin{equation*}
\langle x,{\hat i'}c \rangle=\langle j'x,c\rangle^{-1},
\end{equation*} 
where the pairing on the right side is the Langlands pairing between $C(F)$
and $H^1(W_F,\hat C)$. 
\item Let $k \in H^1(F,K)$ and $\hat x \in H^1(W_F,\hat U \to \hat T)$. Then 
\begin{equation*}
\langle i'k,\hat x\rangle=\langle k,{\hat{{j}}'}\hat x \rangle,
\end{equation*} 
where the pairing on the right is the Tate-Nakayama pairing between the
groups 
$H^1(F,K)$ and $\hat K^\Gamma$. 
\end{enumerate}
\end{lemma}
\begin{proof}
Using that the pairing in \cite{KS} is functorial in $T \to U$ (apply this
functoriality to $(K\to 1) \to (T \to U)$ and $(T \to U) \to (1 \to C)$), we
reduce the lemma to the case in which one of $T$,$U$ is trivial, which can
then be handled using the compatibilities (A.3.13) and (A.3.14) of
\cite{KS}. 
\end{proof}

\section{Relation between $\inv(\gamma,\delta)$ 
and $\alpha(\gamma,\delta)$}\label{RKsec2}
We retain all the assumptions and notation of the previous section. In
particular we have the invariant $\inv(\gamma,\delta) \in H^1(F, R_T
\xrightarrow{1-\theta}R_T)$. Throughout this section we
 assume that
$E/F$ is an unramified extension of
$p$-adic fields, and that $\sigma$ is the Frobenius automorphism of $E/F$.
In this situation there is another invariant measuring the position
of
$\delta$ relative to $\gamma$. This invariant arose naturally in
\cite{K-SVLR} in the course of stabilizing the
Lefschetz formula for Shimura varieties over finite fields. This second
invariant, denoted
$\alpha(\gamma,\delta)$, lies in the group $B(T)$ introduced in \cite{K-IAS}
and studied further in \cite{K-IAS2}. 

The goal of this section is to compare $\inv(\gamma,\delta)$ and
$\alpha(\gamma,\delta)$ and then to rewrite the ratio of
$\Delta_0(\gamma_H,\delta)$ and $\Delta_0(\gamma_H,\gamma)$ in terms of
$\alpha(\gamma,\delta)$ rather than $\inv(\gamma,\delta)$. Since the two
invariants lie in different groups, the reader may be wondering what it
means to compare them. Note however that $H^1(F,T)$ injects naturally into
both  $H^1(F, R_T
\xrightarrow{1-\theta}R_T)$ and $B(T)$, which suggests that we need a group
$A$ and a commutative diagram of the form 
\begin{equation*}
\begin{CD}
H^1(F,T) @>>> H^1(F, R_T \xrightarrow{1-\theta}R_T) \\
@VVV @VVV \\
B(T) @>>> A
\end{CD}
\end{equation*}
in which the two new arrows are injective. It should seem plausible that $A$
ought to be a group
$B(R_T \xrightarrow{1-\theta} R_T)$ bearing the same relation to $H^1(F, R_T
\xrightarrow{1-\theta}R_T)$ as $B(T)$ does to $H^1(F,T)$. 

Such a group has already been been introduced and studied in sections 9-13
of \cite{K-IAS2}. The rest of this section will lean heavily on those sections
of \cite{K-IAS2}, whose \emph{raison d'\^etre} is precisely this application to
twisted transfer factors.

This section begins with a review of the relevant material from \cite{K-IAS2},
and then recalls the definition of $\alpha(\gamma,\delta)$. Next comes a
theorem comparing $\inv(\gamma,\delta)$ and $\alpha(\gamma,\delta)$.  
The two invariants do
\emph{not} become equal in
$B(R_T
\xrightarrow{1-\theta} R_T)$; the relation between them is more subtle than that,
as we will see in Theorem \ref{RKthm.main1}. Finally, we express the ratio of twisted
to standard transfer factors in terms of
$\alpha(\gamma,\delta)$.

\subsection{Review of $B(T \to U)$} 
Let $L$ denote the completion of the maximal unramified extension $F^{\un}$
of $F$ in $\overline F$. We use $\sigma$ to denote the Frobenius
automorphism of $L/F$. We are already using $\sigma$ to denote the Frobenius
automorphism of $E/F$, but since $E \subset F^{\un} \subset L$ and $\sigma$
on $L$ restricts to $\sigma$ on $E$, this abuse of notation should lead to
no confusion. 

In this subsection $f:T\to U$ will denote any homomorphism of $F$-tori. We
then have the group \cite[12.2]{K-IAS2} 
\[
B(T \to U):=H^1(\langle \sigma \rangle,T(L)\to U(L)). 
\] 
Elements of $B(T\to U)$ can be represented by simplified $1$-cocycles
\cite[12.1]{K-IAS2} $(t,u)$, where $t \in T(L)$, $u \in U(L)$ satisfy the
cocycle condition $f(t)=u^{-1}\sigma(u)$. Simplified $1$-coboundaries are
pairs $(t^{-1}\sigma(t),f(t))$ with $t \in T(L)$. 

In \cite[11.2]{K-IAS2}  a canonical isomorphism
\[
\Hom_{\cont}(B(T \to U),\mathbb C^\times)\simeq H^1(W_F,\hat U \to \hat T) 
\] 
is constructed; here we implicitly used the canonical isomorphism
\cite[12.2]{K-IAS2} between $B(T\to U)$ and $\mathbf B(T \to U)$, but as we
have no further use for $\mathbf B(T \to U)$, we will not review its
definition.  In particular we have a
$\mathbb C^\times$-valued pairing between $B(T \to U)$ and $H^1(W_F,\hat U
\to \hat T)$. Moreover there is a natural injection \cite[9.4]{K-IAS2} 
\begin{equation}\label{RKeqn.inj}
H^1(F,T\to U) \to B(T \to U).
\end{equation}
Our pairing restricts to one between $H^1(F,T \to U)$ and
$H^1(W_F,\hat U \to \hat T)$, and this restricted pairing agrees
\cite[11.1]{K-IAS2} with the one in Appendix A of \cite{KS}. 

Now we come to the material in \cite[\S13]{K-IAS2}, which concerns the case in
which our homomorphism of tori is of the very special form
$R_T\xrightarrow{1-\theta} R_T$ for some $F$-torus $T$. In this case it is
shown that the exact sequence \cite[(13.3.2)]{K-IAS2} 
\[
1 \to B(T) \to B(R_T\xrightarrow{1-\theta} R_T) \to T(F) \to 1
\]
has a canonical splitting, so that there is a canonical direct product
decomposition 
\begin{equation}\label{RKeqn.dp} 
B(R_T\xrightarrow{1-\theta} R_T)=B(T) \times T(F).
\end{equation}
Similarly it is shown that the exact sequence \cite[(13.3.8)]{K-IAS2} 
\[
1 \to H^1(W_F,\hat T) \to H^1(W_F,\hat R_T \xrightarrow{1-\hat\theta} \hat
R_T) \to \hat T^\Gamma \to 1 
\] 
has a canonical splitting, so that there is a canonical direct product
decomposition 
\begin{equation}\label{RKeqn.ddp} 
H^1(W_F,\hat R_T \xrightarrow{1-\hat\theta} \hat
R_T) = \hat T^\Gamma \times H^1(W_F,\hat T).  
\end{equation} 

Let $x \in B(R_T\xrightarrow{1-\theta} R_T)$ and $\hat x \in H^1(W_F,\hat R_T \xrightarrow{1-\hat\theta} \hat
R_T)$. As we have seen, we may then pair $x$ with $\hat x$, obtaining
$\langle x,\hat x \rangle \in \mathbb C^\times$. Using \eqref{RKeqn.dp} and
\eqref{RKeqn.ddp}, we decompose $x$ as $(x_1,x_2) \in B(T) \times T(F)$, and
$\hat x$ as $(\hat x_1,\hat x_2) \in \hat T^\Gamma \times
H^1(W_F,\hat T)$. We also have the pairing $\langle x_1,\hat x_1 \rangle$
coming from the canonical isomorphisms $B(T)=X_*(T)_\Gamma=X^*(\hat
T^\Gamma)$ of \cite{K-IAS,K-IAS2}, as well as the Langlands pairing $\langle
x_2,\hat x_2
\rangle$. 

\begin{lemma}\label{RKlem.pa}
There is an equality 
\begin{equation*}
\langle x,\hat x \rangle=\langle x_1,\hat x_1 \rangle\langle x_2,\hat
x_2 \rangle^{-1}.
\end{equation*}
\end{lemma}
\begin{proof}
This follows from \cite[Prop.~13.4]{K-IAS2} together with the obvious analog of
Lemma \ref{RKlem.cmp} with $H^1(F,T \to U)$ replaced by $B(T \to U)$.
\end{proof}

\subsection{Review of $\alpha(\gamma,\delta)$}  
\label{app:alpha_reg}
Our assumptions on $\gamma$, $\delta$ are the same as in \ref{RKsub.digd}.
However the group $R_G$ will play no role in the definition of
$\alpha(\gamma,\delta)$, so we prefer to view $\delta$ as an element of
$G(E)$ such that
$N\delta=\delta\sigma(\delta)\dotsm\sigma^{d-1}(\delta)$ is conjugate in
$G(\overline F)$ to our strongly regular element $\gamma \in T(F)$. Then,
since
$H^1(F^{\un},T)$ is trivial, there exists $c \in G(F^{\un}) \subset G(L)$ 
such that 
\begin{equation}\label{RKeqn.c}
c\gamma c^{-1}=N\delta.
\end{equation}
Now define $b \in G(F^{\un})$ by $b:=c^{-1}\delta \sigma(c)$. Applying
$\sigma$ to \eqref{RKeqn.c}, we find that $b$ centralizes $\gamma$, hence lies
in $T(F^{\un}) \subset T(L)$. Making a different choice of $c$ replaces $b$
by a $\sigma$-conjugate under $T(F^{\un})$. Thus it makes sense to
define $\alpha(\gamma,\delta) \in B(T)$ as the $\sigma$-conjugacy class of
$b$. 

\subsection{Precise relation between $\inv(\gamma,\delta)$ and
$\alpha(\gamma,\delta)$} 
Now that we have reviewed $\alpha(\gamma,\delta)$, we can prove one of the
main results of this appendix. We denote by $\inv^B(\gamma,\delta)$ the image
of
$\inv(\gamma,\delta)$ under the canonical injection \eqref{RKeqn.inj} 
\[
H^1(F,R_T \xrightarrow{1-\theta} R_T) \hookrightarrow B(R_T
\xrightarrow{1-\theta} R_T).
\] 
\begin{theorem}\label{RKthm.main1} 
Under the canonical isomorphism \[B(R_T \xrightarrow{1-\theta}
R_T)=B(T)\times T(F),\] the element $\inv^B(\gamma,\delta)$ goes over to the
pair $(\alpha(\gamma,\delta)^{-1},\gamma)$. 
\end{theorem}
\begin{proof}
As usual when working with cocycles, one has to make various choices. In
this proof, unless the choices are made carefully,  $\inv^B(\gamma,\delta)$
will differ from
$(\alpha(\gamma,\delta)^{-1},\gamma)$ by a complicated $1$-cocycle in $R_T
\xrightarrow{1-\theta} R_T$  that
one would  then have to recognize as a $1$-coboundary. We will take care that
this does not happen. 

We have already discussed $R_G$, $\inv(\gamma,\delta)$, and
$\alpha(\gamma,\delta)$. In particular we have chosen $c \in G(F^{\un})$ such
that $c\gamma c^{-1}=N\delta$ and used it to form the element
$b=c^{-1}\delta \sigma(c) \in T(L)$ representing $\alpha(\gamma,\delta) \in
B(T)$. In order to define $\inv(\gamma,\delta)$ we need to choose an element
$g \in R_G(\overline F)$ such that $g(N\delta)g^{-1}=\gamma$. The best
choice for $g$ is by no means the most obvious one. The one we choose lies
in $R_G(F^{\un})$ and is given by a certain function $J \to G(F^{\un})$.  

Recall that $J$ is the set of $F$-embeddings of $E$ in $\overline F$, and
that $j_0 \in J$ is the inclusion $E \subset \overline F$. We now identify
$J$ with $\mathbb Z/d\mathbb Z$, with $i \in \mathbb Z/d\mathbb Z$
corresponding to the embedding $e \mapsto \sigma^i e$ of $E$ in $\overline
F$. Thus $R_G(F^{\un})$ becomes identified with the set of functions $i
\mapsto x_i$ from $\mathbb Z/d\mathbb Z$ to $G(F^{\un})$, and the same is
true with $L$ in place of $F^{\un}$. The Galois action of $\sigma$ on $x \in
R_G(F^{\un})$ is then given by 
\[
(\sigma x)_i=\sigma(x_{i-1}), 
\] 
while the effect on $x$ of the automorphism $\theta$ of $R_G$ is given by 
\[
(\theta x)_i =x_{i+1}.
\] 

For $i=0,1,\dots,d-1$ we put
$g_i:=c^{-1}\delta\sigma(\delta)\sigma^2(\delta)\cdots
\sigma^{i-1}(\delta)$. In particular $g_0=c^{-1}$. Then $i \mapsto g_i$ is
the desired element $g \in R_G(F^{\un})$ satisfying
$g(N\delta)g^{-1}=\gamma$. We leave this computation to the reader,
remarking only that $\gamma$ corresponds to the element $i \mapsto \gamma$ in
$R_G(F^{\un})$, while $\delta$ corresponds to $i \mapsto \sigma^i(\delta)$,
so that $N\delta$ corresponds to $i \mapsto
\sigma^i(\delta)\sigma^{i+1}(\delta)\cdots \sigma^{i+d-1}(\delta)$.

Since our chosen $g$ lies in $R_G(F^{\un})$, the $1$-cocycle
$t_\tau=g\tau(g)^{-1} \in R_T(F^{\un})$ is unramified, which is to say that
$t_\tau$ depends only on the restriction of $\tau$ to $F^{\un}$. Thus we get
a well-defined element $u \in R_T(F^{\un})$ by putting $u:=t_\tau^{-1}$ for
any $\tau \in \Gamma$ such that $\tau$ restricts to $\sigma$ on $F^{\un}$.
It is then clear from the definitions that $\inv^B(\gamma,\delta)$ is
represented by the simplified $1$-cocycle $(u,\delta') \in R_T(L)\times
R_T(L)$. Here, as before, $\delta'=g\delta \theta(g)^{-1}$. 

The element $(u,\delta')$ can be written as the product of two simplified
$1$-cocycles $(u',t')$, $(u'',t'')$ in $R_T(L)\times R_T(L)$. Of course
elements in $R_T(L)$ are given by functions $\mathbb Z/d\mathbb Z \to T(L)$.
We take $u'$ to be the constant function with value $b^{-1}$. We take $t'$
to be the identity. We take $u''$ to be the function given by 
\[ 
u''_i=
\begin{cases}
\gamma &\text{if $i=0$ in $\mathbb Z/d\mathbb Z$,}\\
1 &\text{otherwise.}
\end{cases}
\] 
Finally, we take $t''$ to be the function given by 
\[ 
t''_i=
\begin{cases}
\gamma &\text{if $i=-1$ in $\mathbb Z/d\mathbb Z$,}\\
1 &\text{otherwise.}
\end{cases}
\] 
It is straightforward to verify that $u=u'u''$ and $\delta'=t't''$. Since
$t'=1$ and $u'$ is fixed by $\theta$, it is clear that $(u',t')$ is a
$1$-cocycle. So too is $(u'',t'')$, since its product with $(u',t')$ is a
$1$-cocycle. 

Since $b$ represents $\alpha(\gamma,\delta) \in B(T)$, and since the image
of $b^{-1}$ under $T=R_T^\theta \hookrightarrow R_T$ is $u'$, it is clear
that $(u',t')=(b^{-1},1)$  represents the image of $\alpha(\gamma,\delta)^{-1}$
under the canonical injection $B(T)
\hookrightarrow B(R_T\xrightarrow{1-\theta} R_T)$. 

It remains only to verify that $(u'',t'')$ represents  the image of $\gamma$
under the canonical splitting of the natural surjection 
\[
B(R_T\xrightarrow{1-\theta} R_T) \twoheadrightarrow T(F). 
\]
Since this surjection sends $(u'',t'')$ to $t''_0t''_1\dotsm
t''_{d-1}=\gamma$, we just need to check that the class of $(u'',t'')$ lies
in the  subgroup of $B(R_T\xrightarrow{1-\theta} R_T)$
complementary to $B(T)$ that is described in \cite[p.~326]{K-IAS2}. This is
clear, since $(u'',t'')$ has the form $(\sigma(x),x)$ for $x=t''$, and 
every value of $i \mapsto t''_i$ lies in $T(L)^{\langle \sigma
\rangle}=T(F)$. 
 \end{proof} 

\subsection{More about allowed embeddings} \label{RKsub.mae} 
As mentioned before, now that we are taking $E/F$ to be an unramified
extension of $p$-adic fields, and $\sigma$ to be the Frobenius automorphism
of $E/F$, there is a canonical choice of allowed embedding $\tilde\eta:\LH
\to \LR_G$ determined by $\tilde s$. As we have seen, giving $\tilde\eta$ is
the same as giving a $1$-cocycle $a$ of $W_F$ in $\mathcal Z$ such that
$(a^{-1},\tilde s)$ is a $1$-cocycle of $W_F$ in $\mathcal
Z\xrightarrow{1-\hat\theta} \mathcal Z$. 

Before describing the canonical choice for the $1$-cocycle $a$, we need to
recall the exact sequence 
\[
1 \to I \to W_F \to \langle \sigma \rangle \to 1,
\] 
where $I$ denotes the inertia subgroup of $\Gamma$. By an \emph{unramified}
$1$-cocycle of $W_F$ in $\mathcal Z$ we mean one which is inflated from a
$1$-cocycle of $\langle \sigma \rangle$ in $\mathcal Z^I$. Note that giving
a $1$-cocycle of $\langle \sigma \rangle$ in $\mathcal Z^I$ is the same as
giving an element in $\mathcal Z^I$, namely the value of the $1$-cocycle on
the canonical generator
$\sigma$ of
$\langle \sigma \rangle$. 

\begin{lemma} The element $\tilde s$ satisfies the following properties.  
\begin{enumerate}
\item $\tilde s_j \in Z(\hat H)^\Gamma$ for all $j \in J$. 
\item $\tilde s \in \mathcal Z^I$. 
\item $\hat\theta(\tilde s)=\sigma(\tilde s)$. 
\end{enumerate}
\end{lemma}
\begin{proof}
(1) Recall that $\tilde s_{j_0}=s$ and that $\tilde s_j=1$ for $j \ne j_0$.
Since $s \in Z(\hat H)^\Gamma$, we conclude that (1) is true. 

(2) Since $E/F$ is unramified, the inertia group $I$ acts trivially on $J$.
Therefore for $\tau \in I$ we have 
\[
(\tau\tilde s)_j=\tau(\tilde s_j)=\tilde s_j,
\] 
showing that $\tau$ fixes $\tilde s$. 

(3) Again using that all values of $\tilde s$  are fixed by $\Gamma$, we
compute that 
\[
(\sigma(\tilde s))_j=\sigma(\tilde
s_{\sigma^{-1}j})=\tilde s_{\sigma^{-1}j}=\tilde
s_{j\sigma^{-1}}=(\hat\theta(\tilde s))_j,
\] 
showing that $\sigma(\tilde s)=\hat \theta(\tilde s)$. 
\end{proof}

\begin{corollary}
Let $a$ be the unramified $1$-cocycle of $W_F$ in $\mathcal Z$ sending
$\sigma$ to $\tilde s$. Then $(a^{-1},\tilde s)$ is a $1$-cocycle of $W_F$
in $\mathcal Z \xrightarrow{1-\hat\theta} \mathcal Z$. 
\end{corollary}
\begin{proof}
It follows from the second part of the lemma that $a$ is a valid unramified
$1$-cocycle, and it follows from the third part of the lemma that
$(a^{-1},\tilde s)$ satisfies the $1$-cocycle condition. 
\end{proof}

Combining this Corollary with Theorem \ref{RKthm.main1}, we obtain 
\begin{theorem}\label{RKthm.main2}
There is an equality 
\[
\langle \inv(\gamma,\delta),(a^{-1},\tilde s)\rangle=\langle
\alpha(\gamma,\delta),s \rangle^{-1}, 
\] 
the pairing on the right being the usual one between $B(T)$ and
$\hat T^\Gamma$.
\end{theorem}
\begin{proof}
Using simplified $1$-cocycles of $W_F$ in $\mathcal Z
\xrightarrow{1-\hat\theta} \mathcal Z$, the $1$-cocycle $(a^{-1},\tilde s)$
becomes $(\tilde s^{-1},\tilde s)$, which is of the form $(d^{-1},d)$ for
$d=\tilde s$. Moreover, $\tilde s_j \in Z(\hat H)^\Gamma=(Z(\hat
H)^I)^{\langle
\sigma \rangle}$ for all $j \in J$. It follows from the discussion on pages
327, 328, 331 of \cite{K-IAS2} that $(\tilde s^{-1},\tilde s)$  represents a
class lying in the canonical subgroup of $H^1(W_F, \hat R_T \xrightarrow
{1-\hat\theta} \hat R_T)$ complementary to $H^1(W_F,\hat T)$. It then follows
from Theorem \ref{RKthm.main1}, Lemma \ref{RKlem.pa} and the previous corollary
that 
\begin{equation*} 
\langle \inv(\gamma,\delta),(a^{-1},\tilde s) \rangle =\langle
\alpha(\gamma,\delta)^{-1},s \rangle \langle \gamma,1 \rangle^{-1} =\langle
\alpha(\gamma,\delta),s \rangle^{-1}.
\end{equation*} 
We used that the image of $(a^{-1},\tilde s)$ under $H^1(W_F,\mathcal Z
\xrightarrow{1-\hat\theta} \mathcal Z) \to Z(\hat H)^\Gamma$ is $s$, which
boils down to the  fact that the product of the $d$ values of $\tilde s$
is equal to $s$. 
\end{proof} 

\begin{corollary}\label{RKcor.jj}
When we use the allowed embedding $\tilde\eta$ determined by the special
$1$-cocycle $(a^{-1},\tilde s)$ described above, the twisted transfer factor
$\Delta_0(\gamma_H,\delta)$ is related to the standard transfer factor
$\Delta_0(\gamma_H,\gamma)$ by the equality 
\[
\Delta_0(\gamma_H,\delta)=\Delta_0(\gamma_H,\gamma)\langle
\alpha(\gamma,\delta),s \rangle^{-1}.
\]
\end{corollary}
\begin{proof}
Use the previous theorem together with Proposition
\ref{RKprop.e}. 
\end{proof}

Corollary \ref{RKcor.jj}  justifies the use of $\langle
\alpha(\gamma_0;\delta),s\rangle\Delta_p(\gamma_H,\gamma_0)$ as twisted
transfer factors in \cite[(7.2)]{K-SVLR}, at least for strongly $G$-regular
$\gamma_H$. Under the additional assumption that the derived group of $G$ is
simply connected, the next section will treat all $(G,H)$-regular $\gamma_H$.  
That $\langle
\alpha(\gamma_0;\delta),s\rangle$ (rather than its inverse) appears in
\cite{K-SVLR} is not a mistake; it is due to the fact that the normalization of
transfer factors, both standard and twisted, used in \cite{K-SVLR} is opposite
(see \cite[p.~178]{K-SVLR}) to the one used in \cite{LS1,KS}. However there are
some minor mistakes in the last two lines of page 179 of \cite{K-SVLR}: each of
the five times that $\eta$ appears it should be replaced by $\tilde\eta$, and
the symbols $=t\rtimes \sigma$ near the end of the next to last line should
all be deleted. 

\section{Matching for $(G,H)$-regular elements} 
\label{app:ghreg}
In this section $G$, $F \subset E \subset L$, $\sigma$, $H$ are as in section
\ref{RKsec2}. However we will now consider transfer factors and matching of orbital
integrals for all $(G,H)$-regular semisimple $\gamma_H \in H(F)$. For simplicity
we assume that the derived group of $G$ is simply connected, as this ensures the
connectedness of the centralizer $G_\gamma$ of any semisimple  $\gamma$ in $G$. 

\subsection{Image of the stable norm map} We begin by recalling two facts about
the stable norm map, which we will use to prove a lemma  
needed later when we prove vanishing of certain stable orbital integrals for
non-norms. 

Let $D$ denote the quotient of $G$ by its derived group (which we have assumed to
be simply connected). 
\begin{proposition}[Labesse]
Let $\gamma$ be an elliptic semisimple element in $G(F)$. Then $\gamma$ is a
stable norm from $G(E)$ if and only if the image of $\gamma$ in $D(F)$ is a norm
from $D(E)$. 
\end{proposition}
\begin{proof}
This is a special case of Proposition 2.5.3 in \cite{La-CSCB}. Of course the
implication $\implies$ is obvious and is true even when $\gamma$ is not elliptic. 
\end{proof}

\begin{proposition}[Haines]
Let $M$ be a Levi subgroup of $G$ and let $\gamma$ be a semisimple element in
$M(F)$ such that $G_\gamma \subset M$. Then $\gamma$ is a stable norm from $G(E)$
if and only if it is a stable norm from $M(E)$.  
\end{proposition}
\begin{proof}
This is part of Lemma 4.2.1 in \cite{Haines}. 
\end{proof}

These two results have the following easy consequence.

\begin{lemma}\label{RKlem.nim}
Let $\gamma$ be a semisimple element in $G(F)$ that is not a stable norm from
$G(E)$. Then there exists a neighborhood $V$ of $\gamma$ in $G(F)$ such that no
semisimple element in $V$ is a stable norm from $G(E)$. 
\end{lemma}
\begin{proof}
Let $A$ be the split component of the center of $G_\gamma$. The centralizer $M$ of
$A$ in $G$ is then a Levi subgroup of $G$ containing $G_\gamma$. Note that
$\gamma$ is elliptic in $M(F)$. The property of having a simply connected derived
group is inherited by $M$, and we write $D_M$ for the quotient of $M$ by its
derived group. 

Since $\gamma$ is not a stable norm from $G(E)$, it is certainly not a stable norm
from $M(E)$. By Labesse's result the image $\overline\gamma$ of $\gamma$ in
$D_M(F)$ is not a norm from $D_M(E)$. Since the image of the norm homomorphism
$D_M(E) \to D_M(F)$ is an open subgroup of $D_M(F)$, there is an open neighborhood
of
$\overline\gamma$ in
$D_M(F)$  consisting entirely of non-norms.  Certainly any semisimple element of
$M(F)$ in the preimage
$V_1$ of this neighborhood is not a stable norm from $M(E)$. 

Consider the regular function $m\mapsto \det(1-\Ad(m);\Lie(G)/\Lie(M))$ on $M$. 
Let $M'$ be the Zariski open subset of $M$ where this regular function does not
vanish.  Equivalently $M'$ is the set of points $m \in M$ whose centralizer in
$\Lie(G)$ is contained in $\Lie(M)$, or, in other words, whose connected
centralizer in $G$ is contained in $M$. In particular $\gamma$ belongs to $M'(F)$,
so that $M'(F)$ is another open neighborhood of $\gamma$. Applying Haines' result,
we see that no semisimple element in the open neighborhood
$V_2:=V_1 \cap M'(F)$ of $\gamma$ in $M(F)$ is a stable norm from $G(E)$. 

Finally, consider the morphism $G\times M' \to G$ sending $(g,m')$ to $gm'g^{-1}$.
It is a submersion, so the image $V$ of $G(F)\times V_2$ provides the desired open
neighborhood $V$ of
$\gamma$ in $G(F)$. 
\end{proof}

\subsection{Review of $\alpha(\gamma,\delta)$ in the general case} 
\label{app:alpha}
Let $\gamma$ be a semisimple element in $G(F)$ and put $I:=G_\gamma$, a connected
reductive $F$-group. Suppose that $\gamma$ is the stable norm of some 
$\theta$-semisimple $\delta \in G(E)$, and let $J$ denote the $\theta$-centralizer
 $\{x \in R_G:x^{-1}\delta\theta(x)=\delta\}$ of $\delta$, another
connected reductive $F$-group. 

There exists $c \in G(L)$ such that 
\begin{equation}\label{RKeqn.defc} 
c\gamma c^{-1}=N\delta,
\end{equation} 
where, as before, $N\delta=\delta\sigma(\delta)\dotsm\sigma^{d-1}(\delta) \in
G(E)$. 
Now define $b \in G(L)$ by $b:=c^{-1}\delta \sigma(c)$. Applying
$\sigma$ to \eqref{RKeqn.defc}, we find that $b$ centralizes $\gamma$, hence lies
in $I(L)$. Making a different choice of $c$ replaces $b$
by a $\sigma$-conjugate under $I(L)$. Thus it makes sense to
define $\alpha(\gamma,\delta) \in B(I)$ as the $\sigma$-conjugacy class of
$b$. 

\begin{lemma}
The element $\alpha(\gamma,\delta)$ is basic in $B(I)$. 
\end{lemma}
\begin{proof}
We are free to compute $\alpha(\gamma,\delta)$ using any $c$ satisfying
\eqref{RKeqn.defc}, and therefore we may assume that $c \in G(F^{\un})$. Thus there
exists a positive integer $r$ such that $c$ is fixed by $\sigma^{dr}$. 
Inside the semidirect product $I(L)\rtimes \langle \sigma \rangle$ we then have
$(b\sigma)^{dr}=\gamma^r\sigma^{dr}$, 
and since $\gamma$ is central in $I$, it follows that $b$ is basic \cite{K-IAS} in
$I(L)$. 
\end{proof}

Since $b$ is basic, we may use it \cite{K-IAS,K-IAS2} to twist the Frobenius action
on $I(L)$, obtaining an inner twist $I'$ of $I$ such that $I'(L)=I(L)$ and with
the Frobenius actions $\sigma_{I'}$, $\sigma_I$ on $I'(L)$, $I(L)$ respectively
being related by $\sigma_{I'}(x)=b\sigma_I(x)b^{-1}$ for all $x \in I'(L)=I(L)$. 

Recall that we are writing elements $x \in R_G$ as functions $i \mapsto x_i$ from
$\mathbb Z/d\mathbb Z$ to $G$. There is a homomorphism $p:R_G \to G$ given by
$p(x):=x_0$, but it is only defined over $E$ (not over $F$). The centralizer
$G_{N\delta}$ of $N\delta \in G(E)$ is also defined over $E$, and $p$ restricts to
an $E$-isomorphism $p_J:J \to G_{N\delta}$. Since $c\gamma c^{-1}=N\delta$, the
inner automorphism $\Int(c)$ induces an $L$-isomorphism $I \to G_{N\delta}$.
Therefore $x \mapsto c^{-1}p_J(x)c$ induces an $L$-isomorphism $\psi:J \to I$. 

\begin{lemma}\label{RKlem.btwist}
The $L$-isomorphism $\psi:J \to I$ is an $F$-isomorphism $J \to I'$. In other
words, when we use $b$ to twist the Frobenius action of $\sigma$ on $I$, we obtain
$J$. 
\end{lemma} 
\begin{proof}
Let $x \in J(L)$. We must show that $\psi(\sigma(x))=b \sigma(\psi(x))b^{-1}$. The
left side works out to $c^{-1}\sigma(x)_0 c=c^{-1}\sigma(x_{-1})c$, while the
right side works out to $(c^{-1}\delta
\sigma(c))\sigma(c^{-1}x_0c)(c^{-1}\delta\sigma(c))^{-1}=c^{-1}\delta
\sigma(x_0)\delta^{-1}c$, so we just need to observe that
$\delta\sigma(x_0)\delta^{-1}=\sigma(x_{-1})$, a consequence of the fact that 
$\sigma(x)$ $\theta$-centralizes $\delta$ (apply $p$ to the equality
$\delta\theta(\sigma(x))\delta^{-1}=\sigma(x)$). 
\end{proof}

\subsection{Comparison of $\alpha(\gamma,\delta)$ with $\alpha(\gamma_t,\delta_t)$}
Our next task is to compare $\alpha(\gamma,\delta)$ with $\alpha(\gamma',\delta')$
for suitable $(\gamma',\delta')$ near $(\gamma,\delta)$ with $\gamma'$
regular in $G$. This will be needed in order to understand the behavior of twisted
transfer factors near $(\gamma,\delta)$. As usual in Harish-Chandra's method of
semisimple descent, we obtain suitable $(\gamma',\delta')$ in the following way. 

We retain all the notation of the previous subsection. Choose an elliptic maximal
torus $T$ in $I$. Since $T$ is elliptic, it automatically transfers to the inner
form $J$ of $I$. Let us now see more concretely how this comes about. 

From Lemma \ref{RKlem.btwist} we know that the $L$-isomorphism $\psi:J \to I$ is an
$F$-isomorphism $J \to I'$. Let $i \in I(L)$. Then the $L$-isomorphism $\psi^{-1}
\circ \Int(i):I \to J$ serves to transfer $T$ from $I$ to $J$ if and only if its
restriction to $T$ is defined over $F$. This happens if and only if
$b\sigma(iti^{-1})b^{-1}=i \sigma(t)i^{-1}$ for all $t \in T(L)$, or,
equivalently, if and only if $i^{-1}b\sigma(i) \in T(L)$. Here we used that $T(L)$
is Zariski dense in $T(\bar L)$. 

 Since $T$ is
elliptic in $I$, the image of the map $B(T) \to B(I)$ is the set $B(I)_b$ of basic
elements in $B(I)$ 
\cite[Proposition 5.3]{K-IAS}. Therefore the fiber over $\alpha(\gamma,\delta)$ is
non-empty, which means we may choose $i \in I(L)$ such that
$b_T:=i^{-1}b\sigma(i)
\in T(L)$. As above we then obtain an $F$-embedding $k:T \hookrightarrow J$ (given
by the restriction to $T$ of 
$\psi^{-1}
\circ \Int(i):I \to J$). A standard twisting argument identifies the fiber over
$\alpha(\gamma,\delta)$ with
$\ker[H^1(F,T)\to H^1(F,J)]$, the set that indexes the $J(F)$-conjugacy classes
of embeddings $k':T \hookrightarrow J$ that are stably conjugate to $k$. 
Therefore by varying the choice of $i$, we obtain all the different ways $k':T
\hookrightarrow J$ of transferring $T$ to $J$. We will work with our particular
$i,b_T,k$, but of course everything we do will also apply to the other
choices we could have made.  

Now we are in a position to compute $\alpha(\gamma',\delta')$ for certain suitably
regular $(\gamma',\delta')$ near $(\gamma,\delta)$. Let $t \in T(F)$ and put
$\delta_t:=k(t)\delta \in R_G(F)= G(E)$. Using that $k(t)$ $\theta$-centralizes
$\delta$, we see that the stable norm of $\delta_t$ is represented by 
$\gamma_t:=t^d\gamma
\in T(F)$. Let
$U$ denote the Zariski open subset of $T$ consisting of those $t \in T$ such that
$\gamma_t$ is $G$-regular. For $t \in U(F)$ the centralizer of $\gamma_t$ in $G$
is $T$, and therefore $\alpha(\gamma_t,\delta_t)$ lies in $B(T)$. 

\begin{proposition}\label{RKprop.needyetanotherlabel}
For $t \in U(F)$ the element $\alpha(\gamma_t,\delta_t) \in B(T)$ is represented
by $tb_T \in T(L)$. 
\end{proposition}
\begin{proof}
Recall that $\psi:J \to I$ is given by $\Int(c^{-1})\circ p_J$. Therefore
$k^{-1}:k(T) \to T$ is given by $\Int(i^{-1}c^{-1}) \circ p_J$. Now
$N\delta_t=p(k(t^d))N\delta$. Applying $\Int(i^{-1}c^{-1})$ to both sides of this
equality (and bearing in mind that $i$ centralizes $\gamma$), we find that
$\Int(i^{-1}c^{-1})(N\delta_t)=t^d\gamma$. Therefore $\alpha(\gamma_t,\delta_t)$
is represented by $i^{-1}c^{-1}\delta_t \sigma(c)\sigma(i)$. The identification
$R_G(F)=G(E)$ is induced by $p:R_G \to G$, so that in $G(E)$ we have the equality
$\delta_t=p_J(k(t))\delta$. Therefore our representative for
$\alpha(\gamma_t,\delta_t)$ can be rewritten as $\bigl(
\Int(i^{-1}c^{-1})(p_J(k(t)))\bigr)i^{-1}c^{-1}\delta\sigma(c)\sigma(i)$, 
which simplifies to $tb_T$, as desired. 
\end{proof}

Let $K$ denote the kernel of the homomorphism $T(F) \to B(T)$ that sends $t \in
T(F)$ to the $\sigma$-conjugacy class of $t$ in $T(L)$. It follows easily from
 \cite[\S 7]{K-IAS2} that $K$ is an open subgroup of $T(F)$. The previous
proposition then has the immediate 

\begin{corollary}\label{RKcor.gamt}
For all $t \in U(F)\cap K$ the element $\alpha(\gamma_t,\delta_t) \in B(T)$ maps
to $\alpha(\gamma,\delta)$ under the map $B(T) \to B(I)$ induced by $T \subset I$. 
\end{corollary}

This corollary is exactly what will be needed in the descent argument to come,
through the intermediary of Proposition \ref{RKprop.tdu}. 

\subsection{Twisted transfer factors for $(G,H)$-regular $\gamma_H$} 
Consider a $(G,H)$-regular semisimple element $\gamma_H$ in $H(F)$. The
centralizer
$I_H$ of $\gamma_H$ in $H$ is connected \cite[Lemma 3.2]{K-STF:EST}. Choose an
elliptic maximal torus $T_H$ in $I_H$. In particular $T_H$ is a maximal torus in
$H$ containing $\gamma_H$. Choose an admissible embedding $T_H \hookrightarrow G$. 
We write $\gamma$, $T$ for the images under this embedding of $\gamma_H$, $T_H$
respectively. Then \cite{K-STF:EST} the centralizer $I$ of $\gamma$ in $G$ is an
inner twist of $I_H$. Of course $T$ is an elliptic maximal torus in $I$, and our
chosen isomorphism $T_H\cong T$ exhibits $T$ as the transfer of $T_H$ to the inner
twist $I$ of $I_H$. 

We need a twisted transfer factor $\Delta_0(\gamma_H,\delta)$ for any
$\theta$-semisimple $\delta \in G(E)$ whose stable norm is $\gamma$. These were
not defined in \cite{KS}, but in the current context, that of cyclic base change
for unramified $E/F$, with $\sigma$ being the Frobenius automorphism and the
derived group of $G$ being simply connected, they were defined in \cite{K-SVLR} by
the formula 
\begin{equation}
\Delta_0(\gamma_H,\delta):=\Delta_0(\gamma_H,\gamma)\langle
\alpha(\gamma,\delta),s\rangle ^{-1},
\end{equation}
with $\Delta_0(\gamma_H,\gamma)$ defined as in \cite[2.4]{LS2}. (See the comment
following Corollary \ref{RKcor.jj} concerning the opposite normalization of transfer
factors used in
\cite{K-SVLR}.) The pairing occurring in this formula is between $B(I)_b$ and
$Z(\hat I)^\Gamma$, and comes from the canonical isomorphism \cite{K-IAS,K-IAS2}
$B(I)_b\simeq X^*(Z(\hat I)^\Gamma)$. In forming this pairing, we view $s \in
Z(\hat H)^\Gamma$ as an element of $Z(\hat I)^\Gamma$ via \[
Z(\hat H)^\Gamma \subset Z(\hat I_H)^\Gamma=Z(\hat I)^\Gamma.
\] 
By Corollary \ref{RKcor.jj} this definition of $\Delta_0(\gamma_H,\delta)$ agrees
with the one in \cite{KS} when $\gamma_H$ is $G$-regular. 

We now apply the work we did in the previous subsection to $(\gamma,\delta)$. With
notation as in that subsection we can now formulate 

\begin{proposition}\label{RKprop.tdu} There is an open neighborhood of $1$ in $T(F)$
such that 
\[
\Delta_0(t^d\gamma_H,k(t)\delta)=\Delta_0(\gamma_H,\delta)
\]
for all $t$ in this neighborhood for which $t^d\gamma_H$ is $G$-regular. 
In writing  $t^d\gamma_H$ we are viewing $t$ as an element 
in $T_H(F)$  via our chosen isomorphism $T_H\simeq T$. 
\end{proposition} 
\begin{proof}
This follows from Corollary \ref{RKcor.gamt} and the fact that $t \mapsto
\Delta_0(t^d\gamma_H,t^d\gamma)$ is defined and constant near $t=1$ (see
\cite[2.4]{LS2}, where $\gamma_H$, $\gamma$ are said to be equisingular). 
\end{proof} 

\subsection{Matching of orbital integrals for $(G,H)$-regular $\gamma_H$} 
We continue with $\gamma_H$, $\gamma$, $I_H$, $I$ as in the previous subsection. 
For $f^H \in C^\infty_c(H(F))$ we consider the stable orbital integral 
\[
SO_{\gamma_H}(f^H)=\sum_{\gamma_H'}e(I_{\gamma_H'})O_{\gamma_H'}(f^H),
\] 
where the sum is taken over conjugacy classes of $\gamma_H' \in H(F)$ that are
stably conjugate to $\gamma_H$, and $I_{\gamma_H'}$ denotes the
(connected) centralizer of
$\gamma_H'$ in $H$.  For
$f
\in C^\infty_c(G(E))$ we consider the endoscopic linear combination of twisted
orbital integrals
\[
TO_{\gamma_H}(f)=\sum_\delta e(J_\delta) \Delta_0(\gamma_H,\delta)TO_\delta(f)
\] 
determined by $\gamma_H$. Thus the sum is taken over twisted conjugacy classes of
$\delta \in G(E)$ whose stable norm is $\gamma$, and $J_\delta$ is the twisted
centralizer of $\delta$. When $\gamma$ is not a stable norm from $G(E)$, we have 
$TO_{\gamma_H}(f)=0$, since the sum occurring in its definition is then empty. 

\begin{proposition}\label{RKprop.ghreg} 
Suppose that 
\[
SO_{\gamma_H}(f^H)=TO_{\gamma_H}(f)
\] 
for all $G$-regular semisimple $\gamma_H$ in $H(F)$. Then the same equality holds
for all $(G,H)$-regular semisimple $\gamma_H$ in $H(F)$. 
\end{proposition}

\begin{proof}
We will just sketch the proof since it is essentially the same as that of
Proposition 2 in \cite[p.~640]{K-TN}, as well as those of Lemma 2.4.A in \cite{LS2}
and Proposition 7.2 in \cite{Cl-LF}. 

Fix $(G,H)$-regular semisimple $\gamma_H \in H(F)$. Introduce $T_H$, $T$  
as in the previous subsection. Assume for the moment that $\gamma$ is a stable
norm.  Looking at the degree $0$ part of the germs about $1$ of the functions 
$t \mapsto SO_{t^d\gamma_H}(f^H)$ and $t \mapsto TO_{t^d\gamma_H}(f)$ on $T(F)$, 
we conclude from Proposition \ref{RKprop.tdu} that 
\begin{equation*}
m\sum_{\gamma_H'}(-1)^{q(I_{\gamma_H'})}O_{\gamma_H'}(f^H)=m\sum_\delta
(-1)^{q(J_\delta)}\Delta_0(\gamma_H,\delta)TO_\delta(f),
\end{equation*}
where $q$  assigns to a connected reductive $F$-group the
$F$-rank of its derived group, and $m$ is the common value of the cardinalities of
all the sets $\ker[H^1(F,T_H) \to H^1(F,I_{\gamma_H'})]$ and $\ker[H^1(F,T) \to
H^1(F,J_\delta)]$. 
 Of course
we used sensible Haar measures and Rogawski's formula for the Shalika germ
corresponding to the identity element, just as in the previously cited proofs. We
also used that, when $t^d\gamma$ is $G$-regular, the  $\theta$-conjugacy classes
having stable norm
$t^d\gamma$ are represented by elements of the form $k'(t)\delta$, with $\delta$
again varying through twisted conjugacy classes of
$\delta \in G(E)$ whose stable norm is $\gamma$, and  (for fixed
$\delta$ with stable norm $\gamma$) $k'$ varying through a set of representatives
for the stable conjugacy classes of embeddings $k':T \hookrightarrow J_\delta$ of
the kind appearing in the discussion leading up to Proposition
\ref{RKprop.needyetanotherlabel}.  Dividing both sides of our equality by
$m(-1)^{q(I_0)}$, where $I_0$ is a common quasisplit inner form of all the groups 
$I_{\gamma_H'}$ and $J_\delta$, we obtain  $SO_{\gamma_H}(f^H)=TO_{\gamma_H}(f)$,
as desired. 

When $\gamma$ is not a stable norm from $G(E)$, we must show that
$SO_{\gamma_H}(f^H)=0$. Looking at the degree $0$ part of the stable Shalika
germ expansion for the maximal torus $T_H$ in 
$H$ (see the sentence just before Proposition 1 in \cite[p.~639]{K-TN}), we see
that  it is enough to show that
$SO_{t_H}(f^H)=0$ for all $G$-regular  $t_H \in T_H(F)$ near
$\gamma_H$.  For this it is enough to show that elements $t_H$  near
$\gamma_H$, when viewed in $T(F)$, are not stable norms from $G(E)$, and this
follows from Lemma
\ref{RKlem.nim}.
\end{proof}

\subsection{A correction to \cite{K-TN}}
In the course of looking through section 2 of \cite{K-TN} I noticed an error in
the definition of the Euler-Poincar\'e function $f_{EP}$. The sign character
$\sgn_\sigma$ occurring in the definition of $f_{EP}$ should be defined as follows:
$\sgn_\sigma(g)$ is $1$ if $g$ preserves the orientation of the polysimplex
$\sigma$, and it is $-1$ if $g$ reverses that orientation. When $\sigma$ is a
simplex, $\sgn_\sigma(g)$ is just the sign of the permutation induced by $g$ on
the vertices of $\sigma$. When writing \cite{K-TN} I carelessly assumed that the
same is true for polysimplices, but this is in fact not the case even for the
product of two copies of a $1$-simplex. Then there is a reflection (obviously
orientation reversing) that induces a permutation with cycle structure
$(12)(34)$ (obviously an even permutation) on the four vertices of the square. This
situation actually arises for the Euler-Poincar\'e function on the group
$PGL_2\times PGL_2$.

With this corrected definition of the sign character, the formula
$\sgn_\tau(\gamma)=(-1)^{\dim(\tau)-\dim(\tau(\gamma))}$ used in the proof of
Theorem 2 of \cite{K-TN} becomes correct and so  no change is needed in
that proof. 

%\end{document}

\backmatter

\printindex

\end{document}